\documentclass[reqno,12pt]{amsart}
\usepackage{imakeidx}
\usepackage[normalem]{ulem}
\usepackage[dvipsnames]{xcolor}
\usepackage{amsmath}
\usepackage{amsthm}
\usepackage[margin=3cm]{geometry}
\usepackage{hyperref}
\usepackage{amsfonts}
\usepackage{amssymb}
\usepackage{enumerate}
\usepackage{amstext}
\usepackage{amsbsy}
\usepackage{amsopn}
\usepackage{amscd}
\usepackage{mathtools}
\usepackage[pdftex]{xcolor}
\usepackage{amsxtra}
\usepackage{upref}
\usepackage{graphicx}
\usepackage{epstopdf}

\usepackage[OT2,OT1]{fontenc}
\newcommand\cyr{%
\renewcommand\rmdefault{wncyr}%
\renewcommand\sfdefault{wncyss}%
\renewcommand\encodingdefault{OT2}%
\normalfont
\selectfont}
\DeclareTextFontCommand{\textcyr}{\cyr}

\DeclareFontFamily{OML}{rsfs}{\skewchar\font'177}
\DeclareFontShape{OML}{rsfs}{m}{n}{ <5> <6> rsfs5 <7> <8> <9> rsfs7
  <10> <10.95> <12> <14.4> <17.28> <20.74> <24.88> rsfs10 }{}
\DeclareMathAlphabet{\mathfs}{OML}{rsfs}{m}{n}

\newtheorem{thm}{Theorem}[section]
\newtheorem{lem}[thm]{Lemma}
\newtheorem{sublem}{Sublemma}
\newtheorem*{sublem*}{Sublemma}
\newtheorem{prop}[thm]{Proposition}
\newtheorem{defi}[thm]{Definition}
\newtheorem{cor}[thm]{Corollary}

\newtheorem*{theorem*}{Theorem}

\newtheorem*{lemma*}{Lemma}

\newtheorem*{example*}{Example}

\newtheorem{remark}[thm]{Remark}
\numberwithin{equation}{section}

\renewcommand{\mod}{\mbox{$\,\mathrm{mod}\,$}}

\def\EXP{{\mathbb{E}}}

\renewcommand{\epsilon}{\varepsilon}
\def\wt{\widetilde}

\def\text#1{\textrm{#1}}

\def\emptyset{\varnothing}

\def\vf{\varphi}

\def\L{\Lambda}

\def \R{\mathbb R}
\def \N{{\mathbb N}}
\def\E{\mathbb E}
\def \Z{\mathbb Z}
\def \T{\mathbb T}
\def \H{\mathbb H}

\def\ov{\overline}
\def\un{\underline}

\def\HD{{\mathrm{HD}}}

\def\Q{\mathbb Q}

\def\DS{\displaystyle}

\def\cF{{\mathcal{F}}}
\def\cP{{\mathcal{P}}}

\def\wh{\widehat}
\def\({\biggl(}
\def\){\biggr)}

\def\<{\bold\langle}
\def\>{\bold\rangle}

\def\eps{{\varepsilon}}
\def\Prob{{\mathbb{P}}}
\def\reals{{\mathbb{R}}}

\def\cH{{\mathcal{H}}}
\def\cI{{\mathcal{I}}}

\def\cZ{{\mathcal{Z}}}

\def\tvf{\widetilde{\varphi}}

\def\zc{{\mathfrak{z}}}
\def\Rect{{\mathsf{R}}}
\def\Frob{{\mathsf{Frob}}}
\def\Gen{{\mathsf{Gen}}}
\def\Gain{{\mathsf{Gain}}}

\def\tomega{{\tilde{\omega}}}
\def\BT{{\mathsf{BT}}}

\DeclareMathOperator\const{const}
\DeclareMathOperator\diam{diam}
\DeclareMathOperator\dist{dist}
\DeclareMathOperator\dom{dom}

\DeclareMathOperator\sgn{sgn}

\DeclareMathOperator\Span{span}

\DeclareMathOperator\SL{SL}

\DeclareMathOperator\esssup{ess\,sup}

\def\tpsi{\widetilde{\psi}}
\def\PP{\mathsf{P}}
\def\HH{\mathsf{H}}
\def\sk{\mathsf k}

\title[The Staircase and Irregularities in Uniform Distribution]
{Generic Points on the Staircase and Irregularities in Uniform Distribution}
\author{Dmitry Dolgopyat and Omri Sarig }
\address{Department of Mathematics\\ University of Maryland at College Park, College Park, MD 20740, USA}
\email{dolgop@umd.edu}
\address{Faculty of Mathematics and Computer Science\\ The Weizmann Institute of Science\\ POB 26, Rehovot, Israel}
\email{omsarig@gmail.com}

\def \qed{\hfill$\square$}

\def\tilde{\widetilde}
\def\bar{\overline}

\def\St{{\mathsf{St}}}

\def\diam{{\rm diam}}
\def\dist{{\rm dist}}
\def\eps{{\varepsilon}}

\def\Cov{{\rm Cov}}

\def\Prob{{\mathbb{P}}}

\def\Var{{\rm Var}}

\def\EXP{{\mathbb{E}}}

\def\bbI{\mathbb{I}}

\def\naturals{\mathbb{N}}

\def\reals{\mathbb{R}}

\def\brc{{\bar c}}

\def\brn{{\bar n}}

\def\brxi{{\bar \xi}}

\def\brOmega{{\bar\Omega}}

\def\cC{\mathcal{C}}

\def\cI{\mathcal{I}}

\def\cF{\mathcal{F}}

\def\cH{\mathcal{H}}

\def\cN{\mathcal{N}}

\def\cP{\mathcal{P}}

\def\cU{\mathcal{U}}

\def\cZ{\mathcal{Z}}

\def\fI{\mathfrak{I}}

\def\fl{\mathfrak{l}}

\def\fS{\mathfrak{S}}

\def\hC{{\hat C}}

\def\ho{{\hat o}}

\def\hxi{{\hat\xi}}

\def\hsigma{{\hat\sigma}}

\def\hOmega{{\hat\Omega}}

\def\sC{\mathsf{C}}

\def\sK{\mathsf{K}}

\def\tC{{\tilde C}}

\def\tS{{\tilde S}}

\def\tomega{{\tilde{\omega}}}

\def\tOmega{{\tilde{\Omega}}}

\makeatother

\def\hat{\widehat}

\def\beq{\begin{equation}}
\def\eeq{\end{equation}}

\def\R{\mathbb R}
\def\E{\mathbb E}
\def\N{\mathbb N}
\def\Z{\mathbb Z}

\makeindex
\begin{document}
\maketitle
\begin{abstract}
We study the ergodic integrals for linear flows on the ``infinite staircase," a translation surface with infinite area, and deduce detailed information on the bias properties of local discrepancy functions for irrational numbers of bounded type.  An important new tool in our analysis is the local limit theory for inhomogeneous Markov chains 
developed by the authors in {\em Lecture Notes in Mathematics,} volume {\bf 2331} (2023).
\end{abstract}
\tableofcontents

\newpage

\vspace{3cm}\ \\

\noindent
{\bf Acknowledgements}\\

\medskip
\noindent
 The authors wish to thank Jon Aaronson, Yeor Hafouta, Hitoshi Nakada, and Barak Weiss for useful comments. Part of this research was carried out during a visit of the two authors to the University of Z\"{u}rich, and the authors would like to thank Corinna Ulcigrai her hospitality  during this visit. 
 D.D. thanks Weizmann Institute for excellent working conditions during his visits.
 This work was partially supported by the US-Israel Binational Science Foundation grant 201610. 
 D.D acknowledges NSF grants DMS 2246983 and DMS 2554526 as well as Simons sabbatical fellowship.
 O.S. acknowledges the support of the Minerva grant 714824 and the Israel Science Foundation grants 264/22 and 925/26.

\newpage
\part{Introduction  and Statement of Main Results}
\section{Introduction}

Most of this work is a study of the ergodic integrals of linear flows on the infinite-area translation surface depicted below.

\begin{figure}[h!]
  \centering
\fbox{\includegraphics[width=8cm]{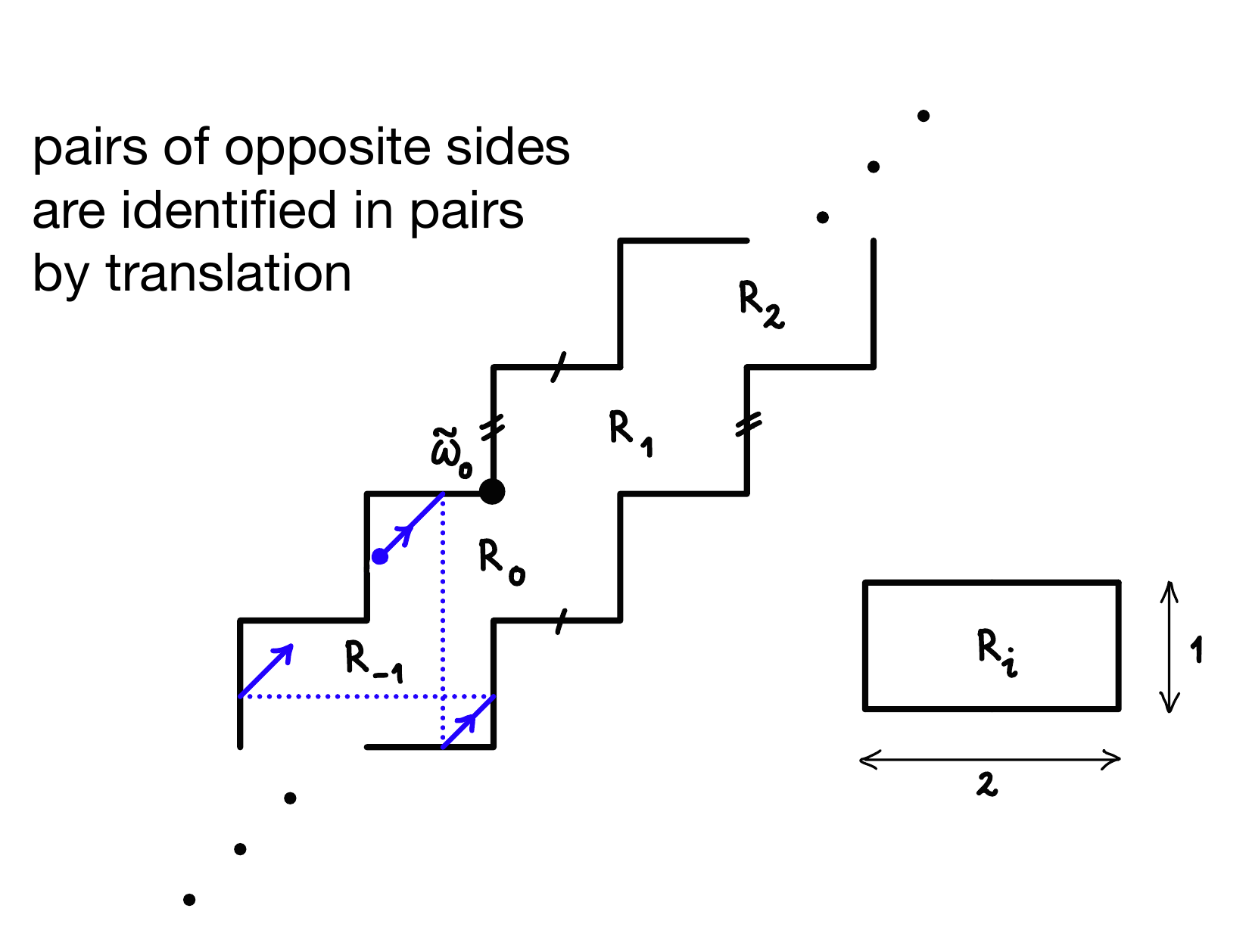}}
  \caption{The linear flow on the infinite staircase }  
\label{Figure-Flow}
\end{figure}

\noindent 
The motivation and the most important applications arise from  the theory of irregularities in uniform distribution of $n\alpha\mod 1$  for $\alpha$ irrational, and this is where we  start.

\subsection{The Error Term in Weyl's Equidistribution Theorem}\label{ss.Background} Henceforth,  $|I|$ denotes the length of an interval $I$;  $\#A$ is the cardinality of a set $A$; and  $\{x\}$ denotes the fractional part of a real number $x$, i.e. $\{x\}\in [0,1)$ and $x-\{x\}\in\Z$.  

Weyl's  equidistribution theorem \index{Weyl's equidistribution theorem}\index{Uniform distribution} states that for every irrational  $\alpha$, and for every sub-interval $I\subset [0,1]$, $$\DS \frac{1}{N}\#\{0\leq j\leq N-1: 
\{j\alpha\}\in I\}\xrightarrow[N\to\infty]{}|I|.$$ 
We are interested in the error term. Define for this purpose the {\em local discrepancy}\index{Discrepancy}\index{Discrepancy!local}
$$
D_N(\alpha,I):=N\times\text{error term}=\#\{0\leq j\leq N-1: \{j\alpha\}\in I\}-N|I|.
$$
Weyl's theorem says that $D_N(\alpha,I)=o(N)$,  but subsequent works show that $D_N(\alpha,I)$ may have unbounded oscillations. Specifically, for every  $\alpha\in\R\setminus\Q$, for a.e. ${\ell}\in (0,1)$,  $\limsup\limits_{N\to\infty} D_N(\alpha,[0,\ell))=\infty$,  $\liminf\limits_{N\to\infty} D_N(\alpha,[0,\ell))=-\infty$ \cite{Oren}, and  
$$
\limsup\limits_{N\to\infty}\frac{|D_{N}(\alpha,[0,\ell))|}{\log N}\geq \frac{1}{400}.
$$
See  \cite{Tijdeman-Wagner} and references therein. 

The nature of the oscillations  depends on the number theoretic properties of 
$\alpha$ and  $\partial I$, and there is much interest in understanding this dependence. 
It is illuminating to fix $I$  and compare the asymptotic behavior of $D_N(\alpha,I)$ for different angles.\footnote{One could also study the  {\em (global) discrepancy\index{Discrepancy!global}} $\Delta_N(\alpha):=\sup\{D_N(\alpha,I):\text{ subintervals }I\}$ \cite{Kuipers-Niederreiter}. 
This work is mostly concerned with  
bounds for  $D_N(\alpha,I)$. These imply lower bounds for $\Delta_N(\alpha)$.}
 In this paper we fix $I=[0,\tfrac{1}{2})$ (the ``symmetric case"), and analyze 
\begin{equation}\label{e.D-N-def}
D_N(\alpha):=D_N(\alpha,[0,\tfrac{1}{2}))=\#\{0\leq j\leq N-1: \{j\alpha\}\in [0,\tfrac{1}{2})\}-\frac{1}{2}N. 
\end{equation}

A surprising and relatively recent discovery is that the error term in Weyl's theorem is not always  spread evenly around zero:\index{Uniform distribution!irregularities in}

\vskip1mm
(1) {\bf Sign Bias (``Heaviness''):}\index{Heaviness}\index{Bias!sign} For some $\alpha$,  $D_N(\alpha)$ is always non-negative. Boshernitzan \& Ralston characterized these  $\alpha$ in \cite{Boshernitzan-Ralston}\index{Boshernitzan-Ralston theorem}.\footnote{
For the existence of  $\alpha$ and $I$ so that $D_N(\alpha,I)$ is bounded on one side but not on the other, see  S\'os \cite{Sos-1974} and Dupain \cite{Dupain} (Veech mentions in \cite{Veech-1977} an unpublished proof by Nelson Markley). The proofs in \cite{Sos-1974} and \cite{Dupain} are quite technical;  for a conceptual proof, see Peres \cite{Peres}.
}

\vskip1mm
(2)  {\bf Global Bias:}\index{Bias!global}\index{Beck's theorem}\index{Temporal CLT}\index{Central limit theorem!Beck's theorem} Beck showed in \cite{Beck-Giant-Leap-1,Beck-Giant-Leap-2} that for every {\em quadratic}\index{Quadratic irrationals!Beck's theorem} irrational $\alpha$, there are constants $A\in\R$ and $B\neq 0$ such that 
$$
\frac{1}{n}\#\left\{1\leq N\leq n: \frac{D_N(\alpha)-A\log n}{B\sqrt{\log n}}\in [a,b]\right\}\xrightarrow[n\to\infty]{}\frac{1}{\sqrt{2\pi}}\int_{a}^{b} e^{-t^2/2}dt.
$$ 
Beck showed that for some $\alpha$, $A\neq 0$. In this case, the histogram of  $(D_1(\alpha),\ldots,D_n(\alpha))$ is given approximately by a Gaussian bell-shaped curve, {\em  centered  away from zero}.

\vskip1mm
(3)  {\bf Local Bias:}\index{Bias!local} \index{Equidistribution} By \eqref{e.D-N-def}, all the values of $D_N(\alpha)$ must belong to $\frac{1}{2}\Z$.
We say that $D_N(\alpha)$ is {\em equidistributed in $\frac{1}{2}\Z$},  if any $\xi_1,\xi_2\in\frac{1}{2}\Z$  appear in $(D_1(\alpha),\ldots,D_n(\alpha))$ with  asymptotically equal frequencies: 
$$
\forall\xi_1,\xi_2\in\frac{1}{2}\Z,\ \ \lim_{n\to\infty}\frac{\#\{1\leq N\leq n: D_N(\alpha)=\xi_1 \}}{\#\{1\leq N\leq n: D_N(\alpha)=\xi_2 \}}=1.
$$
In this work we  show that for many irrationals $\alpha$, $D_N(\alpha)$ is {\em not} equidistributed in $\frac{1}{2}\Z$. We will even construct for each $t$ some  $\alpha$, for which the limit  is  $\exp[t(\xi_1-\xi_2)]$, indicating exponential local bias.
\vskip2mm

The following table describes the bias properties of $D_N(\alpha)$ for a few  irrationals (for proofs, see \cite{Boshernitzan-Ralston}, \cite{Beck-Giant-Leap-1}, and Theorem \ref{t.char-of-equidist-seq} and  Corollary \ref{c.A_n-zero-crit} below):  \index{Bias!examples of numbers with}

\begin{table}[h]
\centering
\caption{Examples of Irrationals with Bias}
\label{table1}
\begin{tabular}{|c|c|c|c|c|c|}
\hline
           & $\sqrt{2}$ & $\sqrt{3}$ & $\sqrt{5}$ & $\sqrt{7}$ & $\alpha_0$\\
 \hline
{\sc Sign Bias} & \text{yes}        &   \text{no}       &    \text{yes}    &    \text{no}   &  \text{no}   \\
 \hline 
{\sc Global Bias} & \text{yes} &   \text{no}     &    \text{yes}    &    \text{yes}    & \text{yes} \\
 \hline
{\sc Local Bias} & \text{yes} &  \text{no}   &   \text{yes}    &    \text{yes}    & \text{no} \\
 \hline
 \end{tabular}
\end{table}
\noindent
%
%
\noindent
The number  $\alpha_0$  is constructed in   Corollary \ref{c.GB-LB}. It has bounded type but is not (and by  Corollary \ref{c.GB-LB} cannot be) a quadratic irrational.  No other combinations of ``yes" and ``no" are possible for irrationals of bounded type (Corollary~\ref{c.GB-LB}). 
 Notice the striking difference between $\sqrt{2}$ and $\sqrt{3}$.

%

\subsection{Number-Theoretic Results} \label{ss.main-results}

It is useful to consider the problem of equidistribution in $\frac{1}{2}\Z$ from a  more general perspective.  Let $\rho:\frac{1}{2}\Z\times\frac{1}{2}\Z\to (0,\infty)$ be a function 
such that $\rho(\xi_1,\xi_2)\rho(\xi_2,\xi_3)=\rho(\xi_1,\xi_3)$ for all $\xi_1,\xi_2,\xi_3\in\frac{1}{2}\Z$. 

\begin{defi}
 A sequence $\{D_N\}_{N\geq 1}$  in $\frac{1}{2}\Z$ is called {\em  $\rho$-equidistributed in}\index{Equidistribution!$\rho$-equidistribution} $\frac{1}{2}\Z$ if 
\begin{equation}\label{e.rho-tiltedness}
\forall\xi_1,\xi_2\in\frac{1}{2}\Z,\ \ \lim_{n\to\infty}\frac{\#\{1\leq N\leq n: D_N(\alpha)=\xi_1 \}}{\#\{1\leq N\leq n: D_N(\alpha)=\xi_2 \}}=\rho(\xi_1,\xi_2).
\end{equation}
\begin{enumerate}[(1)]
\item If $\rho\equiv 1$, we say that $\{D_N\}_{N\geq 1}$ is  {\em equidistributed}.

\item If $\rho\not\equiv 1$, we say that $\{D_N\}_{N\geq 1}$ is  {\em tilted, with profile $\rho$}. \index{Tilt} \index{Tilt!tilt profile} \index{Equidistribution!tilt}

\item If $\rho(\xi_1,\xi_2)=\exp[t(\xi_1-\xi_2)]$ for $t\neq 0$, we say that 
$\{D_N\}_{N\geq 1}$ is  {\em exponentially tilted, with  tilt parameter $t$}. \index{Tilt!exponential} \index{Tilt!tilt parameter}
\end{enumerate}
\end{defi}
\noindent
We would like to know which of these possibilities are realized by $D_N(\alpha)$.

To state our results we need some well-known  facts on continued fractions \cite{Khintchine-Continued-Fractions}. An infinite continued fraction\index{Continued fraction expansion} is a limit $\DS [a_0;a_1,\ldots]:=\lim_{n\to\infty}[a_0;a_1,\ldots,a_n]$, where 
$$
[a_0;a_1,\ldots,a_n]:=a_0+\cfrac{1}{a_1+\cfrac{1}{\cdots+\cfrac{1}{a_n}}}\ .
$$
Any irrational  $\beta$ has a unique representation
$
\beta\!\!=\!\![a_0;a_1,a_2,\ldots], a_0\in\Z$, $a_1,a_2,\ldots\in\N.
$
This is  the {\em regular continued fraction expansion} (regular CFE).\index{Continued fraction expansion!regular}\index{CFE|see {continued fraction expansion}}

 If $\sup a_i<\infty$, we call  $\beta$  an {\em irrational of bounded type}. \index{Bounded type}
Quadratic irrationals have eventually periodic regular CFE,\index{Quadratic irrationals!regular CFE of} and are therefore always of bounded type. It is known that  $\beta$ is of bounded type iff it is {\em badly approximable by rationals}:\index{Badly approximable numbers}\index{Bounded type!and bad approximation by rationals}
\begin{equation}\label{e.badly-approx}
\exists c>0\text{ such that }|\beta-\frac{p}{q}|>cq^{-2}\text{ for all  fractions }\frac{p}{q}. 
\end{equation}

In addition to the regular CFE, \index{Continued fraction expansion!even}
every irrational  $\beta$ also has a unique {\em even continued fraction expansion} (even CFE) 
\begin{equation}\label{e.even-cfe}
\beta=[2m_0;2m_1,2m_2,\ldots], m_0\in\Z, m_1,m_2,\ldots\in\Z\setminus\{0\}.
\end{equation}
 Quadratic irrationals have eventually periodic even CFE, \index{Even continued fraction expansion!of quadratic irrationals}\index{Quadratic irrationals!even CFE of} and  an irrational $\beta$ is of bounded type iff $\sup |m_i|<\infty$ and \eqref{e.even-cfe} does not contain arbitrarily long words of the form $(2,-2,\ldots,\pm 2)$ \cite{Kraaikamp-Lopes}.\index{Even continued fraction expansion!of numbers of bounded type} For more  information on even CFE, see Appendix \ref{AppEven}. 
 
\subsubsection{Equidistribution} The main result of this work is the following: \index{Equidistribution!necessary \& sufficient condition}
\begin{thm}\label{t.char-of-equidist-seq}
Let $\alpha$ be an irrational of bounded type,  and let \eqref{e.even-cfe} be the {\em even} continued fraction expansion of $\beta:=2\alpha-1$. Then $D_N(\alpha)$ is equidistributed in $\frac{1}{2}\Z$ iff 
\begin{equation}\label{e.sgn-criterion}
\lim_{n\to\infty}\frac{1}{n}\sum_{j=1}^n [\sgn(m_{2j})-\sgn(m_{2j-1})]=0.
\end{equation}
\end{thm}
\noindent
For example, if  $\alpha=\sqrt{3}$, then $\beta=2\sqrt{3}-1=[2;2,6,2,6,\ldots]$, the limit is zero, and $D_N(\sqrt{3})$ is equidistributed. But if  $\alpha=\sqrt{7}$, then  $\beta=[4;4,-2,4,10,4,-2,4,10,\ldots]$, \index{Even continued fraction expansion!examples} the limit is not zero, and  $D_N(\sqrt{7})$ is not equidistributed. 

\subsubsection{Tilt}\label{ss.tilt}  Let $\T:=\R/\Z$, $I=[0,\frac{1}{2})$ and $x+I:=\{\{x+y\}:y\in I\}$.  
What are the possible tilt profiles for 
  $D_N(\alpha,x+I)$ and $D_N(\alpha)=D_N(\alpha,I)$? We show that:

\begin{thm}\label{t.tilt-and-exp-tilt}
Suppose $\alpha$ is an irrational of bounded type and $x\in\T$.
If  the sequence $D_N(\alpha,x+I)$ is  tilted, then it is necessarily exponentially tilted. \index{Tilt!implies exponential tilt}
\end{thm}

\begin{thm}\label{t.exp-tilted-exist}
For every $\eta\in\R$, there exists an irrational $\alpha$ of bounded type so that $\{D_N(\alpha)\}$ is exponentially tilted with tilt parameter 
 $\eta$. 
\end{thm}
\index{Tilt!existence of numbers of bounded type with}

\begin{remark}\normalfont
 Nakada's work \cite{Nakada-Keio}\index{Nakada theorem} implies directly that for every irrational $\alpha$ and for every $\eta\in\R$, there exist $x$ such that  $D_N(\alpha,x+I)$ is exponentially tilted with tilt parameter $\eta$, see \S\ref{s.dynamical-reduction}.  Theorem \ref{t.exp-tilted-exist} says that for some $\alpha$, we can take $x=0$.
\end{remark}

\noindent
{\bf Open Problem:} Do there exist irrationals  so that $D_N(\alpha)$ is  $\rho$-tilted but not exponentially tilted? (By Theorem \ref{t.tilt-and-exp-tilt} such $\alpha$ are  of unbounded type.) Which  $\rho$ are possible?

\medskip
Next we ask: How many $x$ are there so that $D(\alpha,x+I)$ is tilted? The work of  Conze \& Keane \cite{Conze-Keane} and Aaronson \& Keane \cite{Aaronson-Keane}  implies directly that for each $\alpha\not\in\Q$, for a.e. $x\in\T$, $D_N(\alpha,x+I)$  is equidistributed (see
\S\ref{s.dynamical-reduction}).\index{Aaronson-Keane theorem}\index{Conze-Keane theorem} Thus,
$$
\text{$
\mathsf T(\alpha):=\{x:D_N(\alpha,x+I)\text{ is tilted}\}\text{ has zero Lebesgue measure for each $\alpha\not\in\Q$}.
$}
$$
Let $\mathrm{HD}(\cdot)$ denote the Hausdorff dimension, and set \index{Tilt!measure, category, and Hausdorff dimension}
$$
\mathsf{T}_{\eta}(\alpha):=\{x\in\T: \text{$D_N(\alpha,x+I)$ is exponentially tilted with tilt parameter $\eta$}\}.
$$
\begin{thm}\label{t.exponential-tilt-HD}
Suppose $\alpha$ is an irrational of bounded type. There exist positive constants $c_1,c_2$ such that for all $\eta$ small enough, 
    $$
    1-c_1 \eta^2\leq \mathrm{HD}(\mathsf{T}_\eta(\alpha))\leq 1-c_2 \eta^2.
    $$
In particular, the set of $x\in\T$ such that $D_N(\alpha,x+I)$  is exponentially tilted has full Hausdorff dimension (and zero Lebesgue measure).
\end{thm}

\begin{thm}\label{t.no-equidist-no-tilt}
The set of $x\in\T$ such that $D_N(\alpha,x+I)$ is not  tilted and not equidistributed is a residual set with full Hausdorff dimension and zero Lebesgue measure.
\end{thm}

\subsubsection{Beck's Theorem for Irrationals of Bounded Type}\label{ss.Beck-Section-Statements}
Beck's theorem on global bias has been extended to  other sequences of low complexity \cite{Levin-2016}, \cite{Dolgopyat-Sarig-TLT},\cite{Aaro-Brom-Chand-IJM}, \cite{Dolgopyat-Sarig-Horo}, \cite{Aaro-Brom-Chand-JMD}, \cite{Paquette-Son}, \cite{Levin-2020}, and to  more general irrationals:\index{Beck's theorem!extensions to bounded type}\index{Bromberg-Ulcigrai theorem}\index{Temporal CLT}

\begin{thm}[Bromberg \& Ulcigrai \cite{Bromberg-Ulcigrai}]\label{e.BU-thm}
Suppose $\alpha$ is an irrational of bounded type and ${\ell}\in (0,1)$ is rational, or more generally a real number in $(0,1)$ such that 
$$
\exists C>0\text{ s.t. }|q\alpha-\ell-p|>C/{|q|}\text{ for all }p\in\Z, q\in\Z\setminus\{0\}\text{ s.t. }\gcd(p,q)=1.
$$
Then there exist $A_n\in\R$ and $B_n>0$ such that for all $a<b$,
$$
\frac{1}{n}\#\left\{1\leq N\leq n: \frac{D_N(\alpha,[0,\ell])-A_n}{B_n}\in [a,b]\right\}\xrightarrow[n\to\infty]{}\frac{1}{\sqrt{2\pi}}\int_a^b e^{-t^2/2}dt.
$$
\end{thm}
\noindent
By contrast, Fr\"uwirth \& Hauke \index{Beck's theorem!Failure for a.e. $\alpha$}\index{Fr\"uwirth-Hauke theorem} showed in \cite{Fruhwirth-Hauke} that the set of $\alpha$ for which such a statement holds for a given $\ell$ has zero Lebesgue measure. See also \cite{Dolgopyat-Sarig-AntiBeck1,Dolgopyat-Sarig-AntiBeck2}.

It is likely that the methods of \cite{Bromberg-Ulcigrai} can be pushed to obtain upper bounds on the growth of $A_n$ and $B_n$. But the authors of that work communicated to us that they cannot estimate $B_n$ from below. Our  tools  allow to do this, in the special case $\ell=\frac{1}{2}$. 
In what follows, 
$
a_n\asymp b_n\text{ means that }\exists C,N>1\text{ s.t. }\forall n>N, a_n/b_n\in [C^{-1},C].
$

\begin{thm}\label{t.A_n-B_n-estimates}
Let $\alpha$ be an irrational of bounded type. Then there exist $A_n=O(\log n)$ and $B_n\asymp\sqrt{\log n}$ such that 
\begin{equation}\label{e.Beck-A_n-B_n-TDLT}
\frac{1}{n}\#\left\{1\leq N\leq n: \frac{D_N(\alpha)-A_n}{B_n}\in [a,b]\right\}\xrightarrow[n\to\infty]{}\frac{1}{\sqrt{2\pi}}\int_a^b e^{-t^2/2}dt.
\end{equation}
\end{thm}
\noindent 
We remark that Beck's original theorem applies to all intervals $[0,\ell]$ with $\ell\in \Q\cap (0,1)$. 
It is likely that Theorem \ref{t.A_n-B_n-estimates}  could also be  extended to such intervals,  along the lines of our proof of Beck's theorem in \cite{Dolgopyat-Sarig-TLT}. However,  we did not pursue this direction.

The behavior of $A_n$ is of interest, because $A_n$ measures the global bias in the histogram of $(D_1(\alpha),\ldots,D_N(\alpha))$.  
The following result
gives more information on $A_n$, when $\ell=\frac{1}{2}$. (Our methods allow to obtain  explicit formulas when  $\frac{1}{2}\neq \ell\in \Q$,
but the resulting formul\ae\, would be significantly more complicated):
\begin{thm}\label{t.A_n}
Let $\alpha$ be an irrational of bounded type, and expand $\beta=2\alpha-1$ into its {\em even} continued fraction expansion $\beta=[2m_0;2m_1,2m_2,\ldots]$. Then one can choose the sequence $A_n$ in Theorem \ref{t.A_n-B_n-estimates} to be of the form 
$
\DS A_n:=\frac{1}{4}\sum_{j=1}^{\nu_n}[\sgn(m_{2j})-\sgn(m_{2j-1})]$, where $\nu_n$ is a non-decreasing sequence with bounded jumps, and $\nu_n\asymp\log n.$
\end{thm}

\index{Bias!global}
\begin{cor}\label{c.A_n-zero-crit} 
Suppose $\alpha$ has bounded type, and $\beta=2\alpha-1$ has even continued fraction expansion $[2m_0;2m_1,2m_2,\ldots]$. The following are equivalent:
\begin{enumerate}[(1)]
\item {\em No Global Bias:} Theorem \ref{t.A_n-B_n-estimates} holds with $A_n=0$;
\item  $\DS \frac{1}{\sqrt{n}} \sum_{j=1}^{n}[\sgn(m_{2j})-\sgn(m_{2j-1})]\xrightarrow[n\to\infty]{}0$.
\end{enumerate}
\end{cor}
\begin{proof}
If \eqref{e.Beck-A_n-B_n-TDLT} holds with two sequences $A_n,A_n'$, then  $\DS \frac{|A_n-A_n'|}{\sqrt{\log n}}\asymp \frac{|A_n-A_n'|}{B_n}\to 0$. Thus \eqref{e.Beck-A_n-B_n-TDLT} holds with $A_n=0$ iff 
$\DS 
\frac{1}{\sqrt{\log n}}\sum_{j=1}^{\nu_n}[\sgn(m_{2j})-\sgn(m_{2j-1})]\to 0$, with $\nu_n$ as in Theorem \ref{t.A_n}. 
This is equivalent to $
\DS \frac{1}{\sqrt{\nu_n}}\sum_{j=1}^{\nu_n}[\sgn(m_{2j})-\sgn(m_{2j-1})]\to 0
$, 
since $\nu_n\asymp \log n$, and to $
\DS  \frac{1}{\sqrt{n}}\sum_{j=1}^{n}[\sgn(m_{2j})-\sgn(m_{2j-1})]\to 0
$, 
since $0\leq \nu_{n+1}-\nu_n\leq \const$.
\end{proof}

We can now compare the following properties:
\begin{enumerate}
\item[(GB)] {\em Global Bias:} \eqref{e.Beck-A_n-B_n-TDLT} holds with $A_n\in\R,B_n\to\infty$ such that  $A_n/B_n\not\to 0$;
\item[(LB)] {\em Local Bias:} $D_N(\alpha)$ is not  equidistributed in $\frac{1}{2}\Z$.\index{Bias!global}\index{Bias!local}\index{Quadratic irrationals!equivalence of local and global bias}
\end{enumerate}
\begin{cor}\label{c.GB-LB} 
For quadratic irrationals, (LB)$\Leftrightarrow$(GB). For irrationals of bounded type, (LB)$\Rightarrow$(GB) but (GB)$\not\Rightarrow$(LB).
\end{cor}
\begin{proof}   By the previous results,   (LB) 
$\DS \Leftrightarrow\frac{1}{{n}} \sum_{j=1}^{n}[\sgn(m_{2j})-\sgn(m_{2j-1})]\not\to 0$,  and 
(GB) $\DS\Leftrightarrow\frac{1}{\sqrt{n}} \sum_{j=1}^{n}[\sgn(m_{2j})-\sgn(m_{2j-1})]\not\to 0$.  Thus, (LB)$\Rightarrow$(GB). 

To see that  (GB)$\not\Rightarrow$(LB), take $\alpha:=\frac{1}{2}(\beta+1)$ where 
$
\beta:=[0;2m_1,2m_2,\ldots,] 
$, $m_i$ is a bounded sequence in $\Z\setminus\{0,1\}$, and 
$\DS \sum_{j=1}^{n}[\sgn(m_{2j})-\sgn(m_{2j-1})]\sim n^{2/3}$.  By Lemma~\ref{l.BT-bounded-pm-two},  $\beta$ has bounded type, and therefore $\beta$ is  badly approximable by rationals. Thus $\alpha$ is badly approximable by rationals, and $\alpha$ has bounded type. By Theorem \ref{t.char-of-equidist-seq},  $D_N(\alpha)$ is equidistributed, and (LB) fails. But by  Corollary \ref{c.A_n-zero-crit}, (GB) holds. 

Now suppose that $\alpha$ is a quadratic irrational. Then $\beta$  is a quadratic irrational,   the {\em even} CFE of $\beta$ is eventually periodic \cite{Kraaikamp-Lopes},\index{Kraaikamp-Lopes theorem}\index{Even continued fraction expansion!of quadratic irrationals} and $$\DS \sum_{j=1}^{n}[\sgn(m_{2j})-\sgn(m_{2j-1})]= \const n+O(1).$$ If the constant is zero, (GB) and (LB) both fail. If the constant is non-zero, (GB) and (LB) both hold. So for quadratic irrationals,  (GB)$\Leftrightarrow$(LB).
\end{proof}

\subsubsection{A ``Local" Beck Theorem}
Beck's theorem describes  $\Prob[D_N(\alpha)\in I_n]$ as $n\to\infty$,  where $N$ is chosen randomly uniformly from $\{1,\ldots,n\}$, and $I_n=A_n+B_n [a,b]$.  Notice that $|I_n|\to\infty$. The following result considers $\Prob[D_N(\alpha)=\xi_n]$ (whence also $I_n$ constant). Such results are called in the probabilistic literature ``local limit theorems."\index{Local Beck theorem}\index{Temporal LLT}

\begin{thm}\label{t.LLT-for-D_N}
For every  irrational  $\alpha$ of bounded type there are constants $A_n=O(\log n)$ as in Theorem \ref{t.A_n}, and $B_n\asymp\sqrt{\log n}$ as follows: If $\xi_n\in\frac{1}{2}\Z$, $t\in\R$, and $\DS \frac{\xi_n-A_n}{B_n}\xrightarrow[n\to\infty]{}t$, then 
$
\DS\#\left\{1\leq N\leq n: D_N(\alpha)=\xi_n\right\}=[1+o(1)]
\frac{ne^{-t^2/2}}{\sqrt{8\pi B_n^2}}\text{ as }n\to\infty.
$
\end{thm}

\begin{cor}\label{c.LLT-for-D_N}
Suppose $\alpha$ is an irrational of bounded type, and $\beta:=2\alpha-1$ has even CFE $[2m_0;2m_1,\ldots]$. If  $\DS \frac{1}{\sqrt{n}}\sum_{j=0}^{n-1}[\sgn(m_{2j})-\sgn(m_{2j-1})]\to 0$, then
$$
\forall \xi\in\frac{1}{2}\Z,\quad
  \#\left\{1\leq N\leq n: D_N(\alpha)=\xi\right\}\asymp n/\sqrt{\log n}.
 $$
 If, in addition, $\alpha$ is a quadratic irrational, then there is a positive constant such that
 $$
\forall \xi\in\frac{1}{2}\Z,\quad
  \#\left\{1\leq N\leq n: D_N(\alpha)=\xi\right\}\sim\const n/\sqrt{\log n}.
 $$
\end{cor}
\begin{proof}
The first statement follows directly from Theorem \ref{t.LLT-for-D_N}; the second is because in the case of quadratic irrationals, $B_n/\sqrt{\log n}\to $positive constant, by Beck's theorem \cite{Beck-Giant-Leap-1,Beck-Giant-Leap-2} (for a short dynamical proof, see \cite{Avila-Dolgopyat-Duryev-Sarig}).
\end{proof}

\noindent
The corollary applies, e.g., to $\alpha=\sqrt{3}$. 

\newpage

\subsection{Reduction to a Problem in Dynamics}\label{s.dynamical-reduction}

\subsubsection{Cylinder Maps}\label{ss.Cylinder-Map} \index{Cylinder map} The following reduction is explained in Klaus Schmidt's paper \cite{Schmidt}. The {\em cylinder map}  is the transformation $T_\alpha:\T\times\Z\to \T\times\Z$ defined by 
$$
T_\alpha(x,\xi)=\bigl(\{x+\alpha\}, \xi+\tau(x)\big),\ \tau(x):=2\cdot 1_{[0,\frac{1}{2})}(\{x\})-1, \text{ $1_A$=indicator function of $A$}.
$$
Let $T^N_\alpha:=T_\alpha\circ\cdots\circ T_\alpha$ ($N$ times), then 
 $T^N_\alpha(x,\xi)=\bigl(\{x+N
\alpha\},\xi+\tau_N(x)\bigr)$, where 
$$
\tau_N(x):=\sum_{j=0}^{N-1}\tau(x+j\alpha)\equiv 2\#\{0\leq j\leq N-1: \{x+j\alpha\}\in [0,\tfrac{1}{2}) \}-N
=  2D_{N}(\alpha,-x+I).
$$ 
Replacing $x$ by $\{-x\}$, we obtain   
\begin{align}
& D_N(\alpha, x+I)=\frac{1}{2}\times\text{ the ordinate of }  T^N_\alpha(-x,0);\\
&\#\{1\leq N\leq n: D_N(\alpha,x+I)=\xi\}=\sum_{j=1}^N (1_{\T\times
\{2\xi\}}\circ T^j_\alpha)(-x,0);\label{e.Cyl-Id-2}\\
&\#\{1\leq N\leq n: D_N(\alpha)=\xi\}=\sum_{j=1}^N (1_{\T\times
\{2\xi\}}\circ T^j_\alpha)(0,0).\label{e.Cyl-Id-3}
\end{align} 
This reduces the study of $D_n(\alpha,x+I)$ and $D_n(\alpha)$ to the study of  $T_\alpha$-orbits. 

For example, Conze \& Keane and Aaronson \& Keane showed in \cite{Conze-Keane}, \cite{Aaronson-Keane} \index{Aaronson-Keane theorem}\index{Conze-Keane theorem}\index{Cylinder map!ergodicity} that for all irrational numbers $\alpha$, the Haar measure $\nu_0$ on $\T\times\Z$ is $T_\alpha$-invariant, ergodic and conservative. Applying  the ratio ergodic theorem, we deduce from \eqref{e.Cyl-Id-2} that
 $$
\lim_{n\to\infty} \frac{\#\{1\leq N\leq n: D_N(\alpha,x+I)=\xi_1\}}{\#\{1\leq N\leq n: D_N(\alpha,x+I)=\xi_2\}}=\frac{\nu_0(\T\times\{2\xi_1\})}{\nu_0(\T\times
\{2\xi_2\})}=1,  \quad \nu_0\text{-a.e.}
$$
Thus  $D_N(\alpha,x+I)$ is equidistributed for Lebesgue a.e. $x\in\T$.
Next, Nakada proved in \cite{Nakada-Keio}\index{Nakada theorem} that for every $t\neq 0$ there is a unique 
$T_\alpha$--ergodic invariant conservative measure $\nu_t$ such that $\nu_t(\T\times\{0\})=1$, and 
$\nu_t\circ Q=e^t \nu_t$, where $Q(x,\xi)=(x,\xi+1)$. Applying  the ratio ergodic theorem, we deduce from \eqref{e.Cyl-Id-2} that for $\nu_t$-a.e. $(-x,\xi)\in\T\times\Z$, 
$$
\forall \xi_1, \xi_2\in \frac{1}{2} \mathbb{Z},\ \lim_{n\to\infty} \frac{\#\{1\leq N\leq n: D_N(\alpha,x+I)=\xi_1\}}{\#\{1\leq N\leq n: D_N(\alpha,x+I)=\xi_2\}}=\frac{\nu_t(\T\times\{2\xi_1\})}{\nu_t(\T\times
\{2\xi_2\})}=e^{2t(\xi_1-\xi_2)}, 
$$
whence $D_N(\alpha,x+I)$ is exponentially tilted for $\nu_t$-a.e. $(-x,0)$. 

This gives us $x'$ such that $D_N(\alpha,x'+I)$ is exponentially tilted. 
But we cannot obtain other profiles of tilt this way, because of the following result:\index{Aaronson-Nakada-Solomyak-Sarig theorem} \index{Cylinder map!ergodic invariant Radon measures}\index{Invariant Radon measures!for the cylinder map}
\begin{thm}[Aaronson, Nakada, Sarig \& Solomyak \cite{ANSS}]\label{t.ANSS} If $\alpha\not\in\mathbb Q$, then  every ergodic $T_\alpha$-invariant Borel measure on $\T\times\Z$ which assigns finite measure to every compact set  is proportional to  $\nu_t$ for some $t\in\R$.
\end{thm}
\noindent
Also, this method 
does not give any information on $D_N(\alpha)$: The behavior of $D_N(\alpha)$ depends on the behavior of the $T_\alpha$-orbit of $(0,0)$, but the ratio ergodic theorem is an almost sure convergence result, and it says nothing about specific initial conditions.

\subsubsection{Generic Points and the Dynamical Formulation of the Problem}\label{ss.generic}
Let $X$ be a 
separable complete metric space such that, except for finitely many {\em singularities}, each point has a compact neighborhood (e.g. $\T\times\Z$ or $\St$ defined below). Let $$X^\ast:=X\setminus\{\text{singularities}\}.$$ 
A Borel measure on $X^\ast$ is called {\em locally finite}, or a {\em Radon measure}, if  $\mu(K)<\infty$ for every compact $K\subset X$.\index{Radon measures}\index{Locally finite measures|see {Radon measure}} Such measures define positive linear functionals on $$C_c(X):=\{f:X\to\R:\text{continuous, with compact support inside $X^\ast$}\}.$$

Suppose $T:X\to X$ is a continuous map,  respectively $\varphi^t:X\to X$ is a continuous semi-flow, and let 
 $\mu$ be a locally finite  ergodic invariant Borel measure on $X^\ast$.
\begin{defi}\label{d.generic}
A point $x_0\in X$ is called {\em $\mu$-generic}\index{Generic points} for  $T$ (resp. $\varphi$), if for any pair of non-negative $h_1,h_2\in C_c(X)$ such that $\mu(h_2)\neq 0$, 
$$
\frac{\sum_{j=0}^{N-1} h_1(T^j x_0)}{\sum_{j=0}^{N-1} h_2(T^j x_0)}\xrightarrow[N\to\infty]{}\frac{\mu(h_1)}{\mu(h_2)},\ 
\left(\text{resp. }
\frac{\int_0^{T} h_1(\varphi^t x_0)dt}{\int_{0}^{T} h_2(\varphi^t x_0)dt}\xrightarrow[T\to\infty]{}\frac{\mu(h_1)}{\mu(h_2)}
\right).
$$
\end{defi}

\smallskip
\noindent
By \eqref{e.Cyl-Id-2}, if $(0,0)$ is $\nu_0$-generic, then $D_N(\alpha)$ is equidistributed, and if $(0,0)$ is $\nu_t$-generic, then $D_N(\alpha)$ is exponentially tilted with tilt parameter $t$. The principal aim  of this work is  to characterize the generic points of $T_\alpha$.

\subsubsection{Linear Flows on the Infinite Staircase}\label{s.infinite-staircase} \index{Constant suspension}  Actually, it is more convenient to characterize the generic points for the constant suspension of $T_\alpha$. Recall that the {\em constant suspension flow} of a continuous map $T\!:\!X\!\to\! X$ is the continuous semi-flow  
$$
\varphi^t:X\times \T\to X\times \T,\text{  given by }\vf^t(z,s)=(T^{\lfloor t\rfloor}(z), \{s+t\})\ \ \ (t\in\R^+).
$$
$\Sigma=X\times\{0\}$ is a Poincar\'e section with a roof function $r\equiv 1$, and the associated Poincar\'e map is conjugate to $T$. Any flow with a Poincar\'e section like that is conjugate to $\varphi$. 
Clearly, $x_0$ is $T$-generic for $\mu$ iff $(x_0,s)$ if $\varphi$-generic for $\mu\times $Lebesgue for some (any) $s$. Thus, the problems of finding the generic points of $\varphi$ and $T$ are equivalent. 

Nevertheless, the constant suspension of $T_\alpha$ is easier to study than $T_\alpha$,  because of the following discovery of
Hooper, Hubert \& Weiss \cite{Hooper-Hubert-Weiss}: $T_\alpha$  is conjugate  to {\em the linear flow on the infinite staircase}. This flow,  which we proceed to define, has many explicit symmetries, which facilitate its study.

\begin{figure}
  \centering
  \fbox{\includegraphics[width=10cm,angle=270]{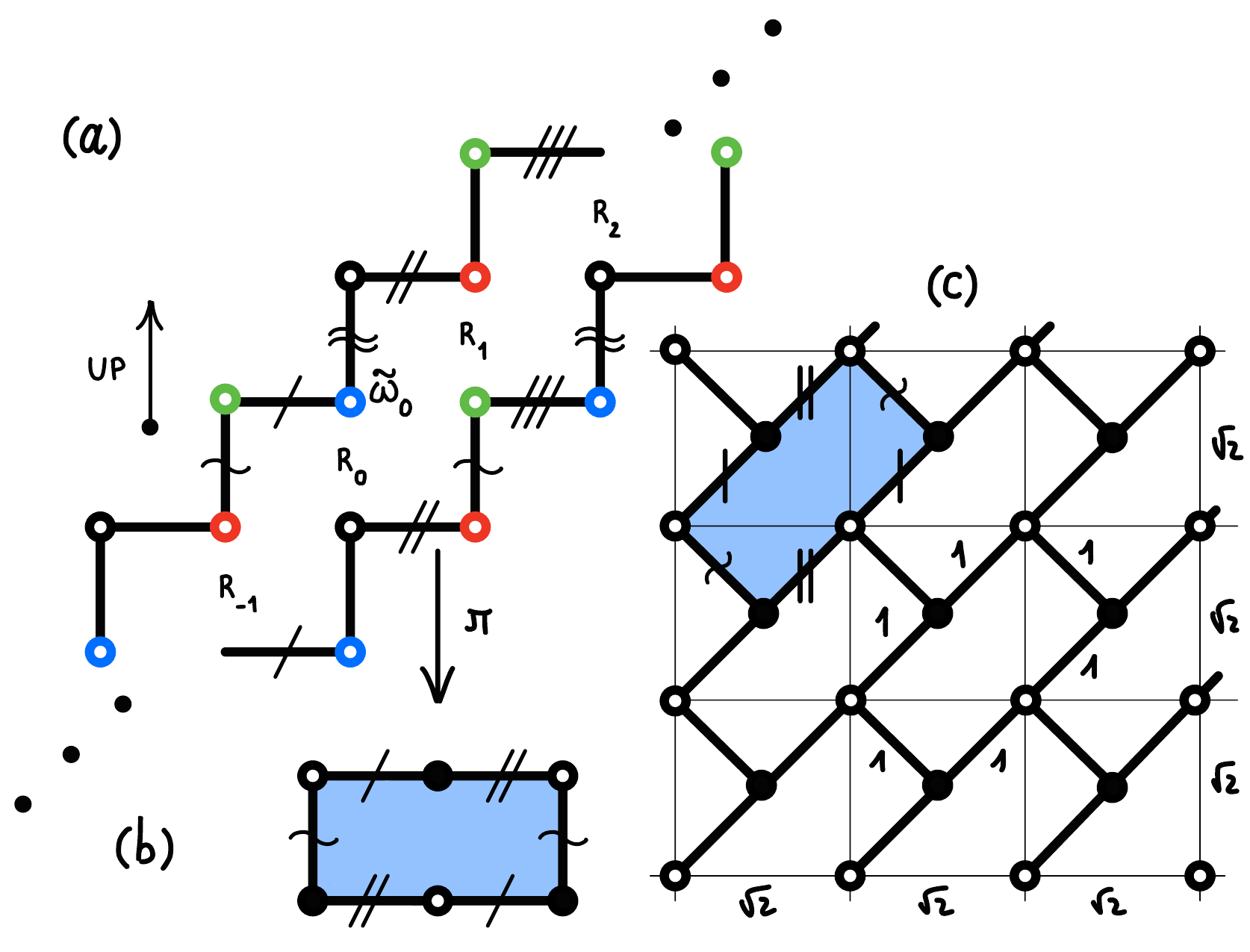}}\\
  \caption{(a) $\mathsf{St}$ and its singularities;  (b) $\mathsf{St}_0$; (c) 
  $\mathsf{St}_0\cong\R^2/\sqrt{2}\Z^2$.}\label{Figure-Staircase}
\end{figure}

The {\em infinite staircase}\index{Infinite staircase} is the identification space $\St$ of the  infinite polygon in 
Figure~\ref{Figure-Staircase}(a), with the natural metric. In this figure, each horizontal rectangle is a 2$\times$1 rectangle with the long side parallel to the $x$ axis, and the short side parallel to the $y$-axis, and pairs of  opposite sides are identified  by translations. The corners fall into four different equivalence classes, each representing a singularity with infinite angle.  Each singularity is the meeting point of infinitely many horizontal rectangles. Consequently,  balls centered at singularities are not compact. \index{Infinite staircase!singularities} \index{Horizontal rectangles}

Let $\St^\ast:=\St\setminus\{\text{singularities}\}$. 
$\St^\ast$ has an obvious atlas of charts which differ by  Euclidean translations. Since translations preserve directions, we can identify the tangent space at all non-singular points with $\R^2$, in the natural way so that ${0\choose 1}$ is direction ``up,"  see Figure \ref{Figure-Staircase}(a).

$\St^\ast$ is naturally partitioned to horizontal rectangles which contain their top side   but
not the corners or the bottom side. We call  one of these rectangles $\Rect_0$, and denote the  others by    $\Rect_k$ $(k\in\Z)$, so that   $\Rect_{k+1}$ is  immediately above $\Rect_k$. 
 Any  non-singular point  $p$ belongs to a unique rectangle. The {\em canonical $\Z$-coordinate} of $p$ is \index{Z-coordinate @ $\Z$-coordinate!canonical}
\begin{equation}\label{e.z-coord-canonical}
\zc(p):=k\text{ s.t. }p\in\Rect_k. 
\end{equation}

Let  $\vec{v}={v_1\choose v_2}$ be a non-zero vector. 
The {\em linear flow with velocity $\vec{v}$} is the semi-flow $\varphi_{\vec{v}}^t:\St\to\St$ which moves a point \index{Linear flow}
$\tilde\omega$  at  velocity $\vec{v}$  subject to identifications. We stipulate  that a line which begins at  or hits a singularity  emerges from the left-most congruent singularity on the sides of $\Rect_0$. If $v_2/v_1$ is irrational, this   happens at most once. 
Let 
$$
\phi^t:=\varphi_{\vec{v}}^t,\text{ where }\vec{v}:={\beta\choose 1}\text{ and }\beta:=2\alpha-1. 
$$
\begin{thm}[Hooper, Hubert \& Weiss \cite{Hooper-Hubert-Weiss}]\label{t.HHW}\index{Hooper-Hubert-Weiss theorem}\index{Infinite staircase!and the cylinder map}\index{Cylinder map!and the infinite staircase}\index{Linear flow!Poincar\'e section for}
$\DS \Sigma:=\bigcup_{k\in\Z}($the top side of $\Rect_k)$  is a Poincar\'e section for $\phi$. The roof function is constant, equal to $1$. The Poincar\'e map is conjugate to the cylinder map $T_\alpha$, via the conjugacy  
$$
\omega:\T\times\Z\to\Sigma,\ \omega(x,k):=\left(\begin{array}{l}
\text{the point on the top side of $\Rect_k$, $y$ distance units right}\\
\text{of the middle, where $y\in [-1,1)$, $y=2x\mod 2$.}
\end{array}\right).
$$ 
\end{thm}

The $T_\alpha$-initial conditions $(0,0)$ and $(-x,0)$  correspond after identifications  to the $\phi$-initial conditions  $\tomega_0:=\omega(0,0)$ and $\tomega_{-x}:=\omega(-x,0)$ in $\St$.  
Notice that if $\tomega$ is on the top side of $\Rect_{2\xi}$, then $\phi^t(\tomega)\in\Rect_{2\xi}$ for all 
$t\in (-1,0]$. Thus by  \eqref{e.Cyl-Id-2}--\eqref{e.Cyl-Id-3},\index{Discrepancy!Infinite staircase and}
\begin{align}
&\#\{1\leq N\leq n:D_N(\alpha)=\xi\}=\int_{-1}^{n-1}1_{\Rect_{2\xi}}(\phi^t(\tomega_0))dt; \label{e.Lin-Flow-1}\\
&\#\{1\leq N\leq n:D_N(\alpha,x+I)=\xi\}=\int_{-1}^{n-1}1_{\Rect_{2\xi}}(\phi^t(\tomega_{-x}))dt. \label{e.Lin-Flow-2} 
\end{align}
We conclude that  the problems we discussed in \S\S \ref{ss.Background}--\ref{ss.main-results} can be recast as problems on the Birkhoff integrals of $\phi$, for the specific initial conditions $\tomega_0$ and $\tomega_{-x}$ .

Let $D:\St\to\St$ denote the map which maps $\Rect_k$ to $\Rect_{k+1}$ by translation (``climb up one stair'').  Then we have the following consequence of Theorems \ref{t.ANSS} and  \ref{t.HHW}:\index{Linear flow!ergodic invariant Radon measures}\index{Hooper-Hubert-Weiss theorem}\index{Invariant Radon measures!for the linear flow on the staircase}

\begin{cor}[Hooper, Hubert \& Weiss \cite{Hooper-Hubert-Weiss}]\label{c.HHW}\index{Invariant Radon measures!$\mu_\xi$}
\label{c.Maharam-Measures} For every $\xi\in\R$ there is a unique locally finite ergodic $\phi$-invariant measure $\mu_{\xi}$ such that $\mu_\xi\circ D=e^{\xi}\mu_\xi$, and $\mu_\xi(\Rect_0)=1$. For $\xi=0$,  $\mu_0=\frac{1}{2}\times$area measure. Any locally finite $\phi$-invariant and ergodic measure is proportional to  $\mu_\xi$ for some $\xi$. 
\end{cor}
\noindent 
To see this, note that the  ergodic invariant measures of the constant suspension of $T$ are the products of the ergodic invariant measures of $T$ and Haar's measure on $\T$.

\subsubsection{The Infinite Staircase as a Regular $\Z$-Cover}\label{ss.Z-d-cover}
 Let $\St^\ast_0$ denote the identification space of the rectangle in Figure  \ref{Figure-Staircase}(b). Let  $\pi:\St^\ast\to\St^\ast_0$ be the map which projects each \index{Infinite staircase!$\Z$-cover of a torus}
$\Rect_k$ to  $\St^\ast_0$ by translation. 
This is a local homeomorphism, with fibres $\pi^{-1}(p):=\{D^n(\wt{p}):n\in\Z\}$ for some (any) $\wt{p}\in\pi^{-1}(p)$. Thus $\pi:\St^\ast\to\St^\ast_0$ is a regular $\Z$-cover.

Looking at Figure \ref{Figure-Staircase}(c), we see that 
$\St_0^\ast$ is isometric to the $\sqrt{2}\times\sqrt{2}$ flat torus with two points (``punctures")\index{Punctures} removed. Let 
$$\St_0:=\R^2/\sqrt{2}\Z^2=\St^\ast_0\cup\{\text{punctures}\}.$$



\subsubsection{The Symmetries of the Infinite Staircase}  
Let $\wt{\psi}\!:\!\St^\ast\!\to\! \St^\ast$ be a differentiable map. Using the canonical identification of the tangent spaces of $\St$ with $\R^2$, we can identify the differential $d\wt{\psi}_\tomega: T_\tomega(\St)\to T_{\tilde \psi(\tomega)}(\St)$  with a matrix, called the {\em derivative}  at $\tomega$.
An {\em affine  automorphism} (or just ``automorphism") is
a homeomorphism $\wt{\psi}:\St\to \St$ s.t. $\wt{\psi}(\St^\ast)=\St^\ast$, $\wt{\psi}$ is differentiable on $\St^\ast$, and $\wt{\psi}$ has constant derivative. \index{Infinite staircase!automorphisms}\index{Automorphisms of the infinite staircase!affine}\index{Affine automorphisms}

If the derivative has positive determinant at all points,  we say that $\wt{\psi}$ is {\em orientation preserving}.\index{Automorphisms of the infinite staircase!orientation preserving}\index{Orientation preserving automorphisms}
A {\em homogeneous automorphism } is an  orientation preserving affine  automorphism $\wt{\psi}$ s.t. $\wt{\psi}\circ D=D\circ\wt{\psi}$, and $\wt{\psi}$ preserves the $D$-orbits of singularities. \index{Automorphisms of the infinite staircase!homogeneity} \index{Homogeneous automorphism} \index{Automorphisms of the infinite staircase!projection to the torus}

Every homogeneous automorphism $\wt{\psi}:\St\to \St$ is the  lift of  an automorphism $\psi$ of the torus $\St_0$ defined by
\begin{equation}\label{Eq-Projected-Auto}
\psi(\omega)=\pi[\wt{\psi}(\wt{\omega})],\text{ for some (any) }\wt{\omega}\in\pi^{-1}(\omega).
\end{equation}
Notice that $\psi$ fixes the punctures of $\St_0$. We call $\psi$ the {\em projection of }
$\wt{\psi}$.\index{Projections of automorphisms} It is easy to see that $\wt{\psi}, \psi$ are area preserving 
\cite[Proposition 2.1]{Avila-Dolgopyat-Duryev-Sarig}.

 The {\em Frobenius function}\label{page-Frobenius} of a homogeneous automorphism $\wt{\psi}$ is  \index{Frobenius function}
$$
\Frob_{\tpsi}:\St^\ast_0\to\Z\ , \ \Frob_{\wt{\psi}}(\omega)=\zc[\wt{\psi}(\wt{\omega})]-\zc(\wt{\omega})\text{ for some (any) }\wt{\omega}\in\pi^{-1}(\omega),
$$
see \eqref{e.z-coord-canonical}. 
The {\em average drift} of $\wt{\psi}$ is \index{Automorphisms of the infinite staircase!average drift}\index{Average drift of a homogeneous automorphism}

\begin{equation}\label{Eq-Drift-Defn}
\delta(\wt{\psi}):=\frac{1}{\mu(\St_0)}\int_{\St_0} 
\Frob_{\wt{\psi}}\, d\mu\ \ (\mu=\text{area measure}).
\end{equation}
We have the identity  $\delta(\wt{\psi}_1\circ\wt{\psi}_2)=\delta(\wt{\psi}_1)+\delta(\wt{\psi}_2)$, see \cite{Avila-Dolgopyat-Duryev-Sarig}. 

Let $
\DS \Gamma(2):=\left\{{\tiny \biggl(\begin{array}{ll}
a & b\\
c & d
\end{array}\biggr)}:a,b,c,d\in\Z,\ ad-bc=1,\  {\tiny \biggl(\begin{array}{ll}
a & b\\
c & d
\end{array}\biggr)}=
{\tiny \biggl(\begin{array}{ll}
1 & 0\\
0 & 1
\end{array}\biggr)}\mod 2
\right\}.
$\index{Automorphisms of the infinite staircase!existence}
\begin{prop}[\cite{Hooper-Hubert-Weiss,Avila-Dolgopyat-Duryev-Sarig}]\label{Prop-Automorphisms}
\begin{enumerate}[(1)]
\item  For every $A\in\Gamma(2)$ there exists a unique homogeneous automorphism $\wt{\psi}$  with $d\wt{\psi}=A$ and $\delta(\wt{\psi})=0$. This automorphism fixes the singularities.

\item  There is also a homogeneous automorphism with drift $\frac{1}{2}$ and derivative ${\tiny \biggl(\begin{array}{cc}
0 & -1\\
1 & 0
\end{array}\biggr)}$.

\item  Every homogeneous automorphism of $\St$ is a product of homeomorphisms from parts (1) and (2). 
\end{enumerate}
\end{prop}


\subsection{Dynamical  Results}\label{s.dynamical-results}
Throughout this section, $\alpha\in
\R\setminus\Q$ is of bounded type,  $\beta:=2\alpha-1$, and $\phi^t:\St\to\St$ is the linear flow on $\St$ with  velocity ${\beta\choose 1}$. The maps $T_{\alpha+m}$ $(m\in\Z)$ are
identical, so 
without loss of generality,  $\alpha\in (0,1)$, and
$\beta\in (-1,1)$.

We expand $\beta$ into the {\em even} CFE
$
\beta=[0; 2m_1,2m_2,\ldots]$,  $m_i\in\Z\setminus\{0\}$. Let $p_n/q_n$ be the associated partial convergents:\index{Partial convergents}\index{Continued fraction expansion!partial convergents $p_n/q_n$}
\begin{equation}
\label{NthConvergent}
p_n/q_n=[0;2m_1,2m_2,\ldots, 2m_n], \text{ and }\mathrm{gcd}(p_n, q_n)=1.
\end{equation}
To fix the signs, we use the standard recursive construction from \cite{Khintchine-Continued-Fractions}: ${p_{-1}\choose q_{-1}}:=
{1\choose 0}$, ${p_{0}\choose q_{0}}:=
{0\choose 1}$, and  $p_k:=2m_k p_{k-1}+p_{k-2}$, $q_k:=2m_k q_{k-1}+q_{k-2}$. Let 
\begin{equation}\label{e.good-matrix}
A_n(\beta):=
\left(\begin{array}{cc} p_{2n-1} & p_{2n}  \\ q_{2n-1} & q_{2n} \end{array}\right)^{-1}\!\!.
\end{equation}
It is known  that $A_n\in\Gamma(2)$
\cite{Kraaikamp-Lopes} (this follows e.g. from identities \eqref{e.A-k} and \eqref{e.A-k-pq} 
in Section \ref{ScReno}).
By Proposition \ref{Prop-Automorphisms} there is a unique $
\wt{\Psi}_n:\St\to\St
$ such that 
\begin{equation}\label{e.Capital-Psi-n}
\wt{\Psi}_n\text{ is the homogeneous automorphism with drift zero and derivative $A_n(\beta)$.}
\end{equation}
Recall that $\zc:\St^\ast\to\Z$ is equal to $k$ on $\Rect_k$, and define
 \begin{equation}\label{e.z-n}
z_n(\tomega):=\begin{cases}
\zc(\wt{\Psi}_n(\tomega)) &\text{ for }\tomega\in \St^\ast,\\
\zc(\wt{\Psi}_n(\phi^t(\tomega))) \text{ for small $t>0$,} & \text{ for $\tomega$ singular.}
\end{cases}
\end{equation}

\subsubsection{Generic Points for $\mu_0$ and Equidistribution}\label{ss.Generic-Zero-Results} 
Recall that $\mu_0=\frac{1}{2}\times$ area measure.\index{Generic points!for $\mu_0$}
\begin{thm}
\label{CrLebGen} 
Suppose $\tomega\in\St$. 
The following are equivalent:
(a) $\tomega$ is $\mu_0$-generic;
(b) $\DS \lim_{n\to \infty} \frac{z_n(\tomega)}{n}=0$;
(c) $\DS \forall k,\ell\in\Z,\ \lim_{T\to \infty} \frac{\int_0^T 1_{\Rect_k} (\phi^t(\tomega)) dt}
{\int_0^T 1_{\Rect_\ell} (\phi^t (\tomega)) dt}=1$.
\end{thm}

\noindent
\begin{cor}
\label{CrCEqui}
$D_N(\alpha,x+I)$ is equidistributed iff $\DS\lim\limits_{N\to\infty}{z_n(\tomega(-x,0))}/{n}=0$.
\end{cor}

\noindent
See \eqref{e.Lin-Flow-2}. In  \S\ref{SSDriftZero}, we  calculate $z_n(\tomega_0)$, and obtain:
\begin{cor}
\label{CrOrigin} 

$\tomega_0:=$the midpoint of the top side of $\Rect_0$ is $\mu_0$-generic iff 
\begin{equation}
\label{SignChanges}
\lim_{n\to\infty} \frac{1}{n} \sum_{k=1}^{n} \left[\sgn(m_{2k-1})-\sgn(m_{2k})\right]=0.
\end{equation}
\end{cor}
\noindent
Theorem \ref{t.char-of-equidist-seq} now follows from \eqref{e.Lin-Flow-1}, see \S\ref{ss.Proof-1} for details.

\begin{thm}
\label{CrLebGen2} 
Suppose $\tomega\in\St$. 
The following are equivalent:
\begin{enumerate}[(a)]
\item[(a)] $\tomega$ is $\mu_0$-generic;
\item[(b)] $\DS \lim_{T\to \infty} \frac{\log \Big( \int_0^T 1_{\St_-} (\phi^t(\tomega)) dt\Big)}
{\log \Big(\int_0^T 1_{\St_+} (\phi^t (\tomega)) dt\Big)}=1$, where  $\St_+:=\bigcup_{k\geq 0}\Rect_k$, $\St_-:=\bigcup_{k<0}\Rect_k$;
\item [(c)] $\DS \lim_{T\to \infty} \left(\log  \int_0^T 1_{\Rect_0} (\phi^t (\tomega)) dt\right)\bigg/\log T=1$.
\end{enumerate}
\end{thm}

\begin{remark}\normalfont
It is natural to wonder if the logarithms in condition (b)
 could be removed. A  related question is whether  $\mu_0$-generic orbits spend asymptotically half of their time in $\St_-$ and $\St_+$.  The answer to both these questions is {\em negative}. See Corollary~\ref{PrBalance} in 
  \S\ref{ss.Limit-Thms-Statements} below. \index{Generic points!for $\mu_0$}
\end{remark}

\subsubsection{Generic Points for $\mu_\xi$ and Tilt} 
\label{SSGP-Tlit}
Let \index{Geometric pressure} 
\begin{equation}
\label{TempPres}
\cP_n(\xi)=\log \left(\int_{\Rect_0} e^{\xi z_n(\tomega)} d\mu_0\right).
\end{equation}
For fixed $n$,  $z_n(\tomega)$ is  bounded, so  $\mathcal P_n(\xi)$ is a finite convex real-analytic  function with asymptotes with slopes $\DS \cP_n'(\pm\infty):=\lim_{\xi \to\pm\infty} \cP_n'(\xi)$. We will see below that for all $n$ large enough,  $\mathcal P_n(\xi)$ is {\em strictly} convex (Lemma \ref{l.strict-convexity-P}). Thus, we may define:
\begin{equation}\label{e.xi-def}
\xi_n(\tomega):=\begin{cases} -\infty & z_n(\tomega)\leq \mathcal P'_n(-\infty)\\
+\infty & z_n(\tomega)\geq \mathcal P'_n(+\infty)\\
\text{the unique $\eta_n$ s.t.  $\cP_n'(\eta_n)=z_n(\tomega)$} & z_n(\tomega)\in (\mathcal P'(-\infty),\mathcal P'_n(+\infty))
\end{cases}
\end{equation}

\index{Generic points!for $\mu_\xi$}

\begin{thm}\label{t.Char-mu-xi-Gen}
\label{ThGen} 
The following are equivalent:
\begin{enumerate}[(a)]
\item[(a)] $\tomega$ is $\mu_\xi$-generic;
\item[(b)] $\DS \lim_{n\to\infty} \xi_n(\tomega)=\xi$;
\item[(c)] $\DS \forall k,\ell\in\Z,\ \lim_{T\to\infty}\frac{\int_0^T 1_{\Rect_k}(\phi^t(\tomega))dt}{\int_0^T 1_{\Rect_\ell}(\phi^t(\tomega))dt}=e^{\xi(k-\ell)}$.
\end{enumerate}
\end{thm}

\medskip
\noindent
{\bf Open Problem:}
Our proof gives a stronger result: There is a constant $R$ such that 
if the limit in (c)  holds for all $k,\ell$ {\em in an interval of length $R$}, then $\tomega$ is generic for some $\mu_\xi.$ Is it enough to assume that the limit in (c) holds just for $k=0,\ell=1$?

\medskip
Analyzing condition (b), we are able to prove the following results: 
\begin{thm}\label{CrWeakBal}
For every $\eta\in\R$ there exists a $\beta$ of bounded type so that 
$\tomega_0$, the midpoint of the top side of $\Rect_0$,  is $\mu_\eta$-generic for the linear flow $\varphi_{\vec{v}}$ with  velocity $\vec{v}={\beta\choose 1}$. 
\end{thm}

Let $\Gen(\xi)$ denote the set of generic points of $\mu_\xi$, and let $\mathrm{HD}$ denote the Hausdorff dimension. \index{Generic points!Hausdorff dimension}\index{Non-generic points} {{Recall that a set is {\em residual}, if it contains a dense $G_\delta$ set.}}\index{Residual set}\index{Hausdorff dimension}

\begin{thm}
\label{CrHDGen}
There are  constants $c_1, c_2>0$ such that for small $\xi$, 
$$ 2-c_1 \xi^2\leq \mathrm{HD}(\mathsf{Gen}(\xi))\leq 2-c_2 \xi^2. $$
\end{thm}
\begin{thm}
\label{CrHDExcept}
The set of points which are not generic for any locally finite ergodic invariant measure of $\phi$ is a 
residual set of Hausdorff
dimension 2 and zero Lebesgue measure.
\end{thm}
\noindent
By \eqref{e.Lin-Flow-2} and the equivalence of 
parts (a) and (c) in Theorem \ref{t.Char-mu-xi-Gen},  $D_N(\alpha,x+I)$ is exponentially tilted with parameter $\xi$ iff $\phi^{-1}(\tomega_{-x})$ is $\mu_\xi$-generic, which is obviously equivalent to the $\mu_\xi$-genericity of $\phi^{t}(\tomega_{-x})$ for some (all) $t\in\R$. Thus 
Theorems \ref{CrWeakBal}, \ref{CrHDGen}, and \ref{CrHDExcept}  directly imply Theorems \ref{t.exp-tilted-exist}, \ref{t.exponential-tilt-HD}, and  \ref{t.no-equidist-no-tilt}.

The following result compares the growth of Birkhoff integrals for generic measures with different parameters: 
\begin{thm}\label{c.diff-xi-diff-asymp}
Suppose $\xi_1>\xi_2\geq 0$ or $\xi_1<\xi_2\leq0$, and let $h$ be a function in $L^1(\mu_{\xi_1})\cap L^1(\mu_{\xi_2})$ such that  $\int h d\mu_{\xi_i}>0$ $(i=1,2)$. Then 
for some $r>0$ independent of $h$, for  $(\mu_{\xi_1}\times\mu_{\xi_2})$-a.e. $(\tomega_1,\tomega_2)$,
$\DS
\frac{\int_0^T h(\phi^t(\tomega_2))dt}{\int_0^T h(\phi^t(\tomega_1))dt}>T^r\text{  for all $T$ large enough.} 
$
\end{thm}

Next we consider a characterization of $\Gen(\xi)$ in terms of the behavior of $z_n(\tomega)$ (compare with Theorem~\ref{CrLebGen}(b)).\index{Birkhoff integrals of the linear flow}
\begin{thm}
\label{CrDriftSeq}
For each $\xi$ there is a sequence $z_n^\xi$ (which depends on $\beta$) such that \index{Generic points!for $\mu_\xi$}
$\tomega$ is $\mu_\xi$-generic  iff $\DS \lim_{n\to \infty} \frac{z_n(\tomega)-z_n^\xi}{n}=0.$ 
\end{thm}

One of the differences between $\mu_0$ and $\mu_\xi$ for $\xi\neq 0$ is that $1_{\St_{\pm}}\not\in L^1(\mu_0)$,  whereas $1_{\St_-}\in L^1(\mu_\xi)$ for $\xi>0$, and $1_{\St_+}\in L^1(\mu_\xi)$ for $\xi<0$ (because $\mu_\xi(\Rect_{k})=e^{\xi k}\mu_\xi(\Rect_0)$).  Let 
 $$\DS \bbI_T(\tomega)\!\!=\!\! \log \left(\int_0^T 1_{[\zc\leq 1]}(\phi^t (\tomega)) dt\right)\text{ and }
\DS \bbI_T^\xi\!\!=\!\!\!\int_{\Rect_0} \bbI_T(\omega) d\mu_\xi(\omega).$$
Note that $ 0\leq \bbI_T(\tomega)\leq \log T$  for all  $\tomega\in \Rect_0$, $T\geq 1$.

\begin{thm}
\label{CrLogOTDetXi}
Suppose $\xi>0$, then: 
\begin{enumerate}[(1)]
\item  
For some $0<c(\xi)<1$,  $c(\xi)\log T\leq \bbI_T^\xi\leq \log T$ for all $T$ large enough;
\item  If $\xi_1>\xi_2>0$, then $\DS\liminf\limits_{T\to\infty}\frac{\bbI_T^{\xi_2}-\bbI_T^{\xi_1}}{\log T}>0$;
\item
{If $\omega$ is generic for $\mu_\xi$ then} 
  $\DS \lim_{T\to\infty}\!\! \frac{\bbI_T(\omega)-\bbI_T^\xi}{\log T}=0$;
\end{enumerate}
A similar result holds for $\xi<0$ and $\xi_1<\xi_2<0$, with $\{\zc\leq 1\}$
replaced by $\{\zc\geq -1\}.$
\end{thm}

%

\subsubsection{Non-Generic Points} 
\label{SSNonGen}
Let $\varphi^t:X\to X$ be a continuous semi-flow on a separable complete metric space $X$ such that every point except for finitely many singularities has a compact neighborhood (e.g. $\St$).\index{Non-generic points} Recall that 
$X^\ast:=X\setminus\{\text{singularities}\}$, and $C_c(X)$ is the space of continuous functions  with compact support inside $X^\ast$. 

\begin{defi}\label{d.historical}
 A locally finite Borel measure \index{Historical measures} $\mu$ on $X^\ast$ is called  {\em historical} for $x_0\in X$, if $\exists T_n\to\infty$ such that for any pair of functions $h_1,h_2\in C_c(X)$ such that $\int h_2 d\mu\neq 0$,
$$
\frac{\int_0^{T_n} h_1(\varphi^t(x_0))dt}{\int_0^{T_n} h_2(\varphi^t(x_0))dt}
\xrightarrow[n\to\infty]{}\frac{\int h_1 d\mu}{\int h_2 d\mu}.
$$
\end{defi}
\noindent
The terminology is motivated by \cite{Ruelle-Historic}. If a point is generic for some measure $\nu$
then $\nu$ is the unique historical measure for that point.
We have the following strengthening
 of Theorem \ref{t.Char-mu-xi-Gen}:

\begin{thm}\label{t.historical-crit} Suppose $\tomega\in\St$ and $\xi_n(\tomega)$ is given by \eqref{e.xi-def}.
\begin{enumerate}[(1)]
\item $\xi\in\R$ is a limit point of $\xi_n(\tomega)$
iff  $\mu_\xi$ is a historical measure of $\tomega$. 
\item If $\infty$ or $-\infty$ are limit points of $\xi_n(\tomega)$, then   there are two finite unions of horizontal rectangles $\mathsf F_1$, $\mathsf F_2$ and  a sequence $ T_n\to\infty$  such that 
$\DS
\frac{\int_0^{T_n}1_{\mathsf F_1}(\varphi^t(\tomega))dt}{\int_0^{T_n}1_{\mathsf F_2}(\varphi^t(\tomega))dt}\xrightarrow[n\to\infty]{}0.
$
\end{enumerate}

\end{thm}

Analyzing the behavior of $\xi_n(z_n(\tomega))$, we are also able to prove:

%

\begin{thm}\label{CrHist}
For every $\tomega\in\St$, every historical measure is proportional  to $\mu_\xi$ for some $\xi$, and the set $\mathcal H(\tomega):=\{\xi\in\R:\mu_\xi\text{ is historical for }\tomega\}$ is either empty, a singleton, or an interval (possibly infinite).
\end{thm}

\subsubsection{Dynamical Beck Theorem and Extensions}\label{ss.Dynamical-Beck-Results}
Let $\mathrm{mes}$ denote Lebesgue's measure on $\R$.  The following ``temporal central limit theorem'' (see \cite{Dolgopyat-Sarig-TLT}) implies the results on Beck's theorem  in \S\ref{ss.Beck-Section-Statements}. 

\index{Temporal CLT}\index{Central limit theorem!temporal}\index{CLT|see{central limit theorem}}
\begin{thm}[Temporal CLT]
\label{ThStairTLT}
There are $\sigma_T$ and  $a_T(\tomega)$ 
such that for every $\tomega$,
\begin{enumerate}[(1)]
\item $\DS
\forall a<b,\ \lim\limits_{T\to\infty}\frac{1}{T} \mathrm{mes}\left\{t\in [0,T]:\frac{\zc(\phi^t \tomega)\!-\!a_T(\tomega)}{\sigma_T}\in [a,b]\right\}=
\frac{1}{\sqrt{2\pi}}\int_a^b e^{-s^2/2}ds
$;
\item 
For some $C$, for all $\tomega$, for all $T>C$, $|a_T(\tomega)|\leq C \log T$ and 
$C^{-1}\leq \frac{\sigma_T}{\sqrt{\log T}}\leq C$;

\item   For some $C$, for all $\tomega$, for all $T>C$, there is $n=n(T)\in\N$ so that   
$C^{-1}<\frac{n}{\log T}<C$, and $a_T(\tomega)=z_n(\tomega)$,     see \eqref{e.z-n}.
\end{enumerate} 
\end{thm}

The next theorem implies the ``local'' version of Beck theorem (Theorem \ref{t.LLT-for-D_N}): \index{Temporal LLT} \index{Local limit theorem!temporal} \index{LLT|see {local limit theorem}}

\begin{thm}[Temporal Local CLT]\label{t.Local-TCLT} 
In the notation of the previous theorem, if $k_T\in\Z$, $s\in\R$, and $\DS\frac{k_T-a_T(\tomega)}{\sigma_T}\xrightarrow[T\to\infty]{}s$, then
$$
\frac{1}{T}\mathrm{mes}\{t\in [0,T]: \phi^t(\tomega)\in\Rect_{k_T}\}=[1+o(1)]\frac{e^{-s^2/2}}{\sqrt{2\pi\sigma_T^2}},\text{ as }T\to\infty.
$$
%
 \end{thm}

\subsubsection{Behavior of Ergodic Integrals}\label{ss.Limit-Thms-Statements}
\label{SSLT}
Recall that $\St_{-}:=\bigcup_{k<0}\Rect_k$ and 
$\St_{+}:=\bigcup_{k\geq 0}\Rect_k$.\index{Birkhoff integrals of the linear flow}

\begin{thm}\label{PrBalance-second half}
For Lebesgue a.e. $\tomega$, for every $\epsilon>0$, for all $T$ large enough,
$
\DS \int_0^T 1_{\St_-}(\phi^t(\tomega))dt\geq \frac{T}{(\log\log T)^{1+\epsilon}}
$. 
Moreover, the power is optimal:
\begin{equation}
\label{Bal3}
 \liminf_{T\to\infty} \frac{\log \left[\int_0^T 1_{\St_-} (\phi^t(\tomega)) dt\right]-\log T}{\log\log\log T}=-1\text{ Lebesgue a.e.}
  \end{equation}
A similar statement holds with $\St_+$ instead of $\St_-$.
\end{thm}
\noindent

By \eqref{Bal3}, for Lebesgue a.e. $\tomega$, for every $\epsilon>0$ there are sequences  $T_k^+, T_k^-\to\infty$ such that 
$
\int_0^{T_k^{\pm}}1_{\St_{\pm}}(\phi^t\tomega)dt=o(T_k^{\pm}) 
$. Necessarily, $
\int_0^{T_k^{\pm}}1_{\St_{\mp}}(\phi^t\tomega)dt\geq T_k^{\pm}-o(T_k^{\pm}) 
$, whence:

\begin{cor}For Lebesgue almost every point, 
\label{PrBalance}
  \begin{equation}
\label{Bal1}
 \liminf_{T\to\infty} \frac{1}{T}\int_0^T 1_{\St_-} (\phi^t\tomega) dt=0,
 \text{ whence }\liminf_{T\to \infty} \frac{\int_0^T 1_{\St_-} (\phi^t\tomega) dt}
{\int_0^T 1_{\St_+} (\phi^t \tomega) dt}=0.
 \end{equation}
 \begin{equation}
\label{Bal2}
\limsup_{T\to\infty} \frac{1}{T}\int_0^T 1_{\St_-} (\phi^t\tomega) dt=1,  \text{ whence }\limsup_{T\to \infty} \frac{\int_0^T 1_{\St_-} (\phi^t\tomega) dt}
{\int_0^T 1_{\St_+} (\phi^t \tomega) dt}=\infty,
 \end{equation}
\end{cor}


 Suppose $\mathcal Z$ is a random variable without atoms. Recall that a parameterized family of random variables $\mathfrak X_T$ {\em converges in distribution} to $\mathcal Z$, if \index{Convergence in distribution}\index{Birkhoff integrals of the linear flow}
 $
 \Prob[\mathfrak X_T\in (a,b)]\xrightarrow[T\to\infty]{}\Prob[\mathcal Z\in (a,b)]$ for every $a<b$.
 We write in this case $\mathfrak X_T\xrightarrow[T\to\infty]{dist}\mathcal Z.$
 We consider such  results for  
$$
\mathfrak X_T=\mathfrak I_T(\tomega,h):=\int_0^T h(\phi^t(\tomega))dt,
$$ for suitable functions $h$ and  {\em random} $\tomega\in\St$. 
Let $\cN$ denote a random variable with the Gaussian distribution with mean zero and variance one,  and let $\cU$ denote a random variable with  uniform distribution on $[0,1].$
\begin{thm}
\label{ThLT-ZD}
There is a function $a_T\asymp T/\sqrt{\log T}$ such that if
$\tomega$ is distributed on the staircase according to a smooth compactly supported probability density with respect to $\mu_0$, then for each
non-negative $h\in L^1(\mu_0)$ such that $\int h d\mu_0=1$, we have
$$
\frac{1}{a_T}\mathfrak I_T(\tomega,h)\xrightarrow[T\to\infty]{dist}\exp\left(-\frac{1}{2}\mathcal N^2\right),
\quad \text{and}\quad \frac{\mathfrak I_T(\tomega,1_{\St_-})}{T}\xrightarrow[T\to\infty]{dist}\mathcal U.
$$
%
%
%
\end{thm}

\begin{thm}
\label{ThLT-NZD}
For each $\xi>0$
there are functions $$a_T(\xi)\asymp\sqrt{\log T} \quad and \quad
(1-c)\log T\leq b_T(\xi)\leq \log T\quad with \quad 0<c<1$$ 
such that for every smooth compactly supported  probability density $\vf$  with respect to $\mu_\xi$, if 
$\tomega$ is distributed on the staircase according to $\vf d\mu_\xi$, and $h\in L^1(\mu_\xi)$ is 
 non-negative with positive integral (e.g. $h=1_{\St_-}$), then
$$ \frac{\log \fI_T(\tomega, h)-b_T(\xi) }{a_T(\xi) }\xrightarrow[T\to\infty]{\text{dist}}  \cN.
$$
The same result holds for $\xi<0$ with  $\St_+$ instead of $\St_-$.\index{Convergence in distribution}\index{Birkhoff integrals of the linear flow}
\end{thm}

\subsection{Historical Context for Our Dynamical Results}
The problem of determining the  behavior of Birkhoff sums for specific rather than for a.e. initial conditions has been considered before, in other scenarios.

\subsubsection{Unique Ergodicity} \index{Unique ergodicity} If a continuous map on a compact metric space has a unique invariant probability measure $\mu$, then the ergodic averages of any continuous function $f$ converge to  $\int f d\mu$, and the convergence is uniform. In particular, every point is $\mu$-generic. This observation, first formalized by Oxtoby \cite{Oxtoby} and Furstenberg \cite{Furstenberg-Strict-Ergodicity}, building on earlier work of Kryloff \& Bogolyubov \cite{Kryloff-Bogoliouboff},  remains the most general way to prove the convergence of ergodic averages at specific points. 

The dynamical systems in this paper are not uniquely ergodic, and they  admit many orbits which are not generic for any  locally finite ergodic invariant measure, finite or infinite (Theorem \ref{CrHDExcept}).

\subsubsection{Generic Points for Finite Measure Preserving Parabolic Dynamical Systems.}
The term {\em ``parabolic dynamical system"} is a loose term for dynamical systems with 
polynomial divergence of nearby orbits.

The crown jewel 
of parabolic dynamics is Ratner's theory for unipotent homogeneous flows \cite{Ratner-1,Ratner-3,Ratner-3}. 
Let $G$ be 
a Lie group, let $\Gamma$ be a lattice in $G$ and let $u(t)$ be a unipotent one-parameter subgroup. The main results of Ratner theory are the {\em measure classification theorem}, 
which says that all finite $u(t)$-invariant measures on $G/\Gamma$ are Haar measures on closed orbits of 
proper subgroups of $G$, and the {\em equidistribution theorem} which states that 
every point is generic for some {\em finite} ergodic invariant measure.\index{Ratner theory} \index{Generic points!unipotent flows}\index{Veech dichotomy}

Much  effort has been directed 
at extending Ratner theory beyond the homogeneous setting, see e.g.
\cite{EMM15, EM18, FL10, FK16, RT06, Thouvenot95}. But  a description of generic points
is available only in a few cases beyond uniquely ergodic systems. Several
recent works show that  Ratner's equidistribution theorem is often false 
for non-homogeneous systems (\cite{BarreiraSchmeling, CM15, Chaika-Khalil-Smillie, Chaika-Smillie-Weiss}).

 One notable positive result in this area is the {\em Veech dichotomy} (\cite{Veech89}). Recall that a translation surface of finite genus is
called a {\em Veech surface} if its $\SL(2,\R)$ orbit is closed in  Teichm\"uller's space. Examples include square tiled surfaces. For Veech surfaces,  any  linear flow at a given direction is either 
periodic or uniquely ergodic. In particular, every point is generic for some invariant  probability measure. See \cite{HubertSchmidt}
for a review of Veech dichotomy and related results.

The infinite staircase is similar to Veech surfaces in that its Veech group (the group of derivatives of affine automorphisms) is a lattice in $\mathrm{SL}(2,\R)$ (Proposition \ref{Prop-Automorphisms}). But since it is the identification space of an infinite rather than a finite polygon,
 it is not covered by the theory above. Indeed, 
   the linear flows on $\St$ in irrational directions 
are neither periodic,  nor uniquely ergodic.

\subsubsection{Infinite Measure Classification for Parabolic Dynamical Systems}
For unipotent flows on infinite volume homogeneous spaces, all ``non-trivial" measures are infinite. It is interesting to  extend Ratner's theory to this case (see \cite{Burger90,BL98,Babillot04,Sarig-Radon, Ledrappier-Sarig-Classification, SarigSchapira, Sarig-Tame,Schapira-Generic,Schapira-HalfGeneric,
Winter15,  Mohammadi-Oh15,Pan18,  OhPan19,LL22, Landesberg-Lindenstrauss,LLLO, Choi-Kim}   and the references therein). 

First we discuss the measure classification problem. 

The simplest non-trivial case is the horocycle
flow of a hyperbolic surface of infinite genus. \index{Babillot-Ledrappier measures}\index{Invariant Radon measures!for horocycle flows}
 In this case the finite ergodic invariant measures are all carried by closed periodic horocycles, but as first discovered by Babillot and Ledrappier \cite{BL98,Babillot04} there could be infinitely many other globally supported singular
  {\em infinite} (locally finite) invariant measures.  
If the surface can be partitioned to a family of  pairs of pants, with boundary components of length bounded above, then  the locally finite ergodic invariant measures all
admit an explicit description in terms of positive (not necessarily $L^2$) minimal eigenfunctions of the Laplace-Beltrami operator  
\cite{Babillot04, Ledrappier-Sarig-Classification, Sarig-Tame}.  As a rule of thumb, the bigger the surface, the more room there is for non-trivial positive eigenfunctions, and this results in additional invariant measures.   

Some of these results have been extended to the case of higher-dimension and/or variable negative curvature \cite{Sarig-Radon,Ledrappier-Fourier,Pollicott-Zd-Covers,Pan18,Landesberg-Lindenstrauss,Landesberg}.
 The locally finite ergodic invariant measures for the dynamical systems in this paper have all been identified  in \cite{ANSS} and \cite{Hooper-Hubert-Weiss}, see Theorem \ref{t.ANSS} and Corollary \ref{c.HHW}. 

\subsubsection{Generic Points for Infinite Measure Preserving Parabolic Systems}
\label{ss.Schapira-Sa} The problem of identifying the generic points of infinite measure preserving unipotent flows  is much more subtle, and it is only  understood in two cases: horocycle flows on hyperbolic surfaces with finite genus \cite{Schapira-Generic,Schapira-HalfGeneric}, and  horocycle flows of  $\Z^d$-covers of compact hyperbolic surfaces \cite{SarigSchapira}. 
The infinite staircase is a also a $\Z$-cover of a finite area surface, and the description of generic points in both cases is similar. 

Consider a $\Z^d$-cover of a compact hyperbolic surface, and enumerate the group of deck transformations by $\mathsf{Deck}=\{D_z:z\in\Z^d\}$ so that $D_z\circ D_w=D_{z+w}$. By \cite{BL98}, 
for each $\xi\in\R^d$ there is a unique locally finite ergodic invariant measure $\mu_{\xi}$ 
such that 
\begin{equation}
\label{DeckInv}
 \mu_{\xi}\circ D_z^{-1}=e^{\<\xi,z\>} \mu_{\xi}\ \ (z\in\Z^d).
 \end{equation}
By \cite{Sarig-Radon}, there are no other locally finite ergodic invariant  measures. 

 Choose some connected pre-compact fundamental domain $\mathfrak F$ for the action of $\mathsf{Deck}$. For all $T$ there is a unique  $z_T(\omega)\in\Z^d$ such that  $g^t \omega\in D_{z_T(\omega)}(\mathfrak{F})$, where $g^t$ is the geodesic flow.  Now let
\begin{equation}
\label{AsymCycle}
 \cZ(\omega)=\lim_{T\to \infty} z_T(\omega)/T
\in\R^d
\end{equation} 
  (if this limit exists). By the work of Babillot \& Ledrappier \cite{Babillot-Ledrappier-Lalley}, for every $\xi\in\R^d$, there is a unique $z(\xi)\in\R^d$ such that 
$\cZ(\omega)=z(\xi)$  $\mu_{\xi}$-a.e.
By \cite{SarigSchapira}, \index{Generic points!for the Babillot-Ledrappier measures}
\begin{equation}\label{e.Drift=Generic}
\text{ $\omega$ is $\mu_{\xi}$-generic iff }\lim_{T\to \infty} \frac{z_T(\omega)}{T}=z(\xi).
\end{equation}
In the special case $\xi=0$, $\mu_{\xi}=$ Lebesgue's measure,  $z(0)=0$, and we deduce  that $\omega$ is generic for Lebesgue's measure iff $z_T(\omega)/T\to 0$.  This characterization is similar to  the characterizations of $\mu_\xi$-generic $\tomega$ in Theorems \ref{CrLebGen} and \ref{CrDriftSeq} in terms of the behavior of $z_n(\tomega)$ as $n\to\infty$.

To see the horocyclic analogue  to   Theorem \ref{t.Char-mu-xi-Gen}(b), we need to recall some additional results from the work of Babillot \& Ledrappier \cite{Babillot-Ledrappier-Lalley}.  
Consider the Babillot-Ledrappier pressure function 
$$\DS \mathcal P_{BL}(\xi)=\lim_{T\to\infty} \frac{1}{T} \log \EXP\left(e^{\<\xi, z_T(\omega)\>} \right) \ \ \ (\xi\in\R^d), $$
where $\EXP$ means the integral with respect to a measure with a smooth
probability density
compactly supported  inside $\mathfrak F$. 
\footnote{
\cite{Babillot-Ledrappier-Lalley} defines the pressure function differently, but it is possible to show using the volume lemmas of \cite{Bowen-Ruelle} that the two definitions coincide.}

By  \cite{Babillot-Ledrappier-Lalley}, $z(\xi)=\nabla \mathcal P_{BL}(\xi).$
Let 
$$\xi_T(\omega):= \text{ the solution to }
 \nabla \mathcal P_{BL}(\xi_T)=\frac{z_T(\omega)}{T},\text{ if it exists.}
$$
Babillot \& Ledrappier proved that $\nabla \mathcal P_{BL}$ is a homeomorphism from $\R^d$ onto an open convex subset of $\R^d$. Thus by
\eqref{e.Drift=Generic}, 
$\omega$ is generic $\Leftrightarrow \nabla \mathcal P_{BL}(\xi_T)\to \nabla\mathcal P_{BL}(\xi)$ $\Leftrightarrow\xi_T\to\xi$. This characterization of genericity is similar to the one in Theorem \ref{t.Char-mu-xi-Gen}(b).

We can use \eqref{e.Drift=Generic} to construct  non-generic points. If  $z_T(\omega)/T$ oscillates without converging, then $\tomega$ is non-generic. It should be possible to show using the theory of non-generic orbits of hyperbolic systems (see e.g. \cite{BarreiraSchmeling}) that the set of such initial conditions is a residual set of full Hausdorff dimension, but this does not seem to have been written anywhere.

Another route to non-genericity is for $z_T(\omega)/T$ to converge to a limit which is not equal to $z(\xi)$ for any $\xi$. By \cite{Babillot-Ledrappier-Lalley}, $\cZ(\omega)=z(\xi)$ for some $\xi$ iff $\cZ(\omega)$ is {\em non-extremal} in the sense that   $\exists \omega_1,\omega_2$ such that  $\cZ(\omega)$ is a non-trivial convex combination of 
$\cZ(\omega_1)$ and $\cZ(\omega_2)$. The question how to construct   $\omega$ with extremal $\cZ(\omega)$ is delicate and interesting, because  it depends on the Riemannian metric of the $\Z^d$-cover and not just on its topology. It is also interesting to understand the orbit closure of the horocycle of such $\omega$. For related work, see  \cite{Landesberg-Farre-Minsky-1,Landesberg-Farre-Minsky-2,DFLM}.

%
%
%

\subsubsection{Linear Flows on the Infinite Staircase}
The principal goal of this paper is to determine the generic points for linear flows  on the infinite staircase, for as many directions as we can. We are able to treat all directions ${\beta\choose 1}$ for $\beta$ of bounded type. 

The special case when $\beta$ is a quadratic irrational was done before in \cite{Avila-Dolgopyat-Duryev-Sarig}   (but the characterization of genericity in terms of the even CFE is new already in this case).

The  difference between the  bounded type case and quadratic irrational case is  that for  quadratic irrationals $\beta$, $\phi$ is {\em renormalized by a single map} in the following sense:  There exist  $0<r<1$ and a homogeneous automorphism  $\wt{\psi}:\St\to\St$ such that 
\begin{equation}\label{e.renorm1}
\wt{\psi}\circ \phi^t=\phi^{r t}\circ\wt{\psi}\ \ (t\in\R).
\end{equation}
But when $\beta$ is not a  quadratic irrational, $\phi$ is not renormalized by a single map, and this accounts  for the new phenomena  in Corollary \ref{c.GB-LB}, and the new difficulties discussed below.\index{Quadratic irrationals!renormalization by a single map}

\subsection{Heuristic Overview of the  Proof and Identification of Main Difficulties}
It is sufficient to prove the  dynamical results in  \S\ref{s.dynamical-results}, because as explained there, they imply all the number-theoretical results in \S\ref{ss.main-results}.

 Let $\phi:\St\to\St$ be the linear flow in direction $\beta\choose 1$, $\beta=2\alpha-1$, $\alpha$ of bounded type.

\medskip
\noindent
{\em The main issue  is to find the asymptotics of 
$
\int_0^T h(\phi^t(\tomega))dt
$
for specific  $\tomega$, as $T\to\infty$, 
for sufficiently many test functions $h:\St\to\R$.} 

\medskip
For simplicity, we consider the   special case $h_0:=1_{\Rect_0}$, since it  contains all essential difficulties. We will also focus on $\tomega:=\tomega_0=$ 
the middle of the top side of $\Rect_0$, the initial condition which determines the behavior of $D_N(\alpha)$, see \eqref{e.Lin-Flow-1}.
In this case,
$$
\frac{1}{T}\int_0^T h_0(\phi^t(\tomega_0))dt=\frac{|A_T\cap\Rect_0|}{|A_T|},\text{  }A_T:=\left(\begin{array}{l}
\text{the line segment in $\St$ starting from $\tomega_0$,}\\
\text{in direction ${\beta\choose 1}$, with length $T\sqrt{1+\beta^2}$,}
\end{array}\right).
$$
Here  $|\cdot|$ denotes the Euclidean length. The problem is to find the asymptotic behavior of  
$|A_T\cap\Rect_0|/|A_T|$, as $T\to\infty$. 

\subsubsection{Renormalization Sequence}\label{ss.renormalization-sequence}
$A_T$ winds around $\St$ in a complicated manner. We will build a sequence of  homogeneous automorphisms which ``simplify it." \index{Renormalization}

If $\wt{\psi}:\St\to\St$ has constant derivative matrix $A$, then $\wt{\psi}$ maps linear segments in direction $\vec{v}$ and  Euclidean length $L$, to linear segments in direction $A\vec{v}$ with Euclidean length $L\|A\vec{v}\|/\|\vec{v}\|$.
We will look for homogeneous automorphisms $\wt{\psi}_T:\St\to\St$ which contract the velocity vector  $\vec{v}={\beta\choose 1}$ so much, that 
$$
A_T^\ast:=\wt{\psi}_T(A_T)\text{ has length in }[c_1,c_2], \text{ with $0<c_1<c_2<1$ independent of $T$.}
$$
Since $|A_T^\ast|\leq c_2$, $A_T^\ast$ is too short to wind, and it remains in a single horizontal rectangle.  Since $|A_T^\ast|\geq c_1$,  $A_T^\ast$ will not degenerate to  a point, when $T\to\infty$. However, the $\Z$-coordinate of the  rectangle which contains $A_T^\ast$ (call it ${z_T^\ast}$) may escape to infinity.

\begin{figure}[h!]
  \centering
  \fbox{\includegraphics[width=6cm,angle=90]{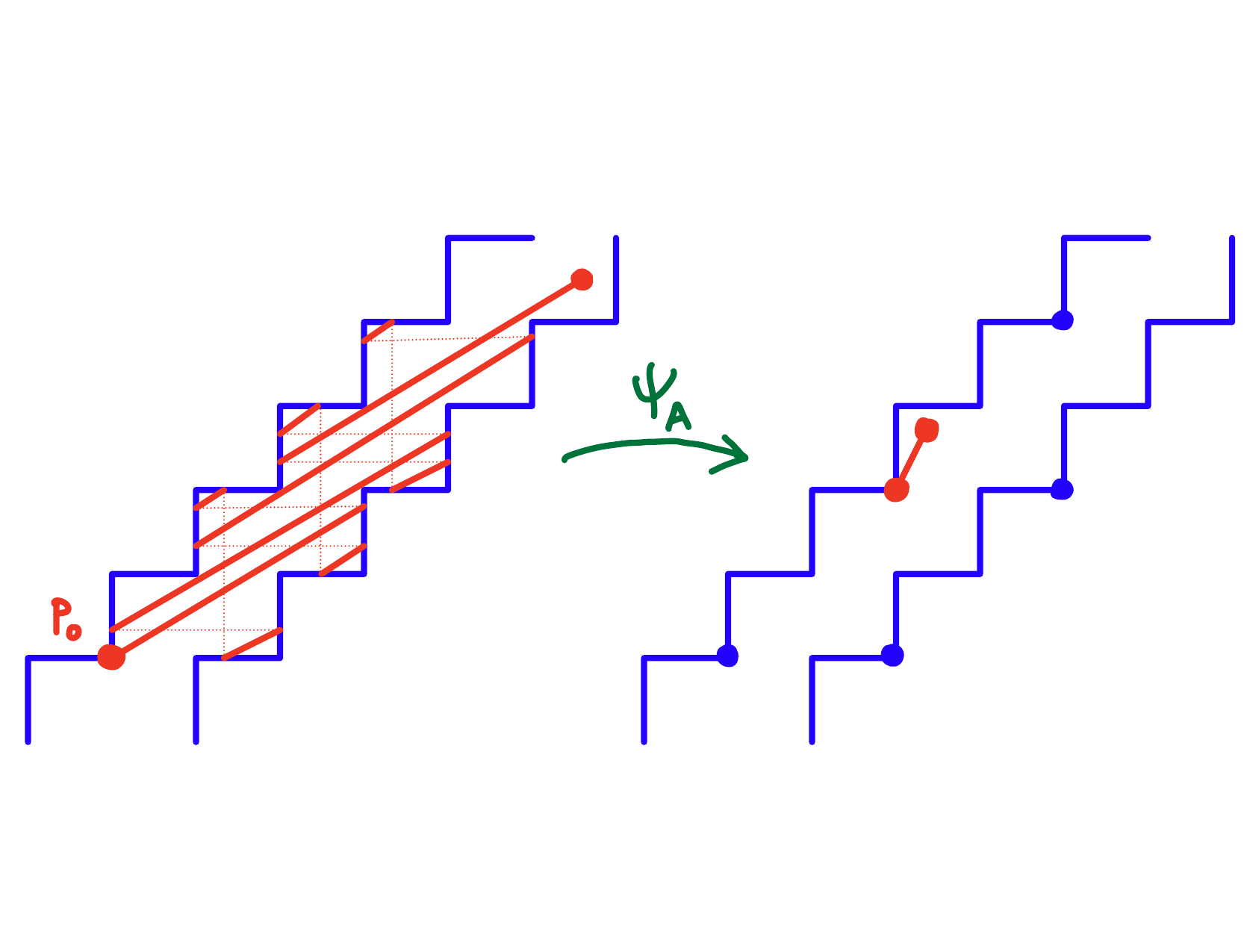}}
  \caption{Simplifying an orbit by a contracting automorphism}\label{Figure-Renormalization}
\end{figure}

Proposition \ref{Prop-Automorphisms} reduces the construction of $\wt{\psi}_T$ to the solution of  the following number theoretic problem: Find  $M_T\in\Gamma(2)$ such that $\|M_T{\beta\choose 1}\|\in [c_1T^{-1},c_2 T^{-1}]$ with $0<c_1<c_2<1$. We will show that   $M_T:=A_{n(T)}(\beta)$ solves this problem, where $A_n(\beta)$ is given by \eqref{e.good-matrix},  and where 
$
n(T)$ is a non-decreasing function with bounded jumps, such that $n(T)\asymp\log T$. Then we will take $\wt{\psi}_T:=\wt{\Psi}_{n(T)}$, with $\wt{\Psi}_n$  as in \eqref{e.Capital-Psi-n}.  With this choice of $\wt{\psi}_T$, the rectangle containing $A_T^\ast$ has $\Z$-coordinate    
$$z_T^\ast:=z_{n(T)}(\tomega_0),
$$
see \eqref{e.z-n}. 
Notice the appearance of the even continued fraction expansion in \eqref{e.good-matrix}.  The bounded type assumption  is needed here, for the lower bound on $\|M_T{\beta\choose 1}\|$, and for  the bounds on  the  jumps and growth of $n(T)$. 

We will sometimes abbreviate, and write $n=n(T)$.

\subsubsection{Probabilistic Formulation}\label{ss.probabilistic-formulation}
Since $\wt{\psi}_T:A_T\to A_T^\ast$ is an affine bijection,  
$$\frac{|A_T\cap\Rect_0|}{|A_T|}\!\!=\!\!
\frac{|\wt{\psi}_T(A_T\cap\Rect_0)|}{|\wt{\psi}_T(A_T)|}\!\!=\!\!\frac{|A_T^\ast\cap\wt{\psi}_T(\Rect_0)|}{|A_T^\ast|}=``\Prob[\wt{\psi}_T^{-1}(\tomega)\in\Rect_0], \text{ as } \tomega\sim\mathcal U(A_T^\ast).\text{''}
$$ 
The expression in quotation marks is the  probability  that $\wt{\psi}_T^{-1}(\tomega)\in\Rect_0$,  when $\tomega$ is selected randomly uniformly from $A_T^\ast$.

Recall the canonical $\Z$-coordinate $\zc$ from \eqref{e.z-coord-canonical}. 
Since $A_T^\ast\subset \Rect_{z_T^\ast}$ and   $z_T^\ast=z_n(\tomega_0)$,  
$
\big(\tomega\in A_T^\ast, \wt{\psi}_T^{-1}(\tomega)\in\Rect_0\big)\Leftrightarrow$  $\zc(\tomega)-\zc(\wt{\Psi}_n^{-1}(\tomega))=z_n(\tomega_0)
$, and therefore  
\begin{equation}\label{e.Int=Prob}
\int_0^{T}h_0(\phi^t(\tomega_0))dt\!\!=\!\!T\cdot \Prob\biggl[\zc(\tomega)
-\zc\bigl(\wt{\Psi}_{n}^{-1}(\tomega)\bigr)\!\!=\!\!z_{n}(\tomega_0)\biggr],\text{ as } \tomega\sim\mathcal U(A_T^\ast), n=n(T).
\end{equation}

\subsubsection{The Random Walk Heuristic}\label{ss.RW-heuristic}
 In Part \ref{p.renormalization}, we will use the even CFE of $\beta$ to find an explicit factorization 
$
\wt{\Psi}_{n}=\wt{\psi}_{n^\ast-1}\circ\cdots\circ \wt{\psi}_0
$ with the following properties:
\begin{enumerate}[(a)]
\item {\em Finiteness:\/} $\{\wt{\psi}_j: j\geq 0\}$ is a finite family of homogeneous automorphisms;
\item {\em Uniform Hyperbolicity:\/} 
$\wt{\psi}_j:(\wt{\psi}_{j-1}\circ\cdots\circ\wt{\psi}_0)(A_T)\to (\wt{\psi}_{j}\circ\cdots\circ\wt{\psi}_0)(A_T)$ is a $\theta$-contraction, with $\theta\in (0,1)$ independent of $j$;
\item {\em Number of Terms:\/} $n^\ast=n^\ast(T)\asymp n(T)$.
\end{enumerate}

Observe that $\zc(\tomega)-\zc\bigl(\wt{\Psi}_{n}^{-1}(\tomega)\bigr)$ in \eqref{e.Int=Prob} is minus the total displacement of the canonical $\Z$-coordinate of the directed path $\tomega^{(n)}\to\tomega^{(n-1)}\to\cdots\to\tomega^{(0)}$, where  
$$
\tomega^{(n)}:=\tomega\ , \ \tomega^{(i-1)}:=\tpsi_{i-1}^{-1}(\tomega^{(i)})\ , \ \tomega^{(0)}
:=(\tpsi_0^{-1}\circ\cdots\circ\tpsi_{n-1}^{-1})(\tomega).
$$
%
%
%
At each step $i=1,\ldots,n^\ast$, the $\Z$-coordinate changes by 
\begin{equation}\label{e.Delta-i}
\Delta_i(\tomega):=\zc(\tomega^{(i-1)})-\zc(\tomega^{(i)})=\zc(\tpsi_{i-1}^{-1}(\tomega^{(i)}))-\zc(\tomega^{(i)}).
\end{equation}
The total displacement is the sum of these changes, therefore by \eqref{e.Int=Prob}, 
\begin{equation}\label{e.Int=RW-Prob}
\int_0^{T}h_0(\phi^t(\tomega_0))dt=T\cdot \Prob\left[\sum_{i=1}^{n^\ast}\Delta_i(\tomega)=-z_n(\tomega_0)\right],\text{ as } \tomega\sim\mathcal U(A_T^\ast).
\end{equation}

The probability in \eqref{e.Int=RW-Prob} bears some formal similarity to the  probability that a ``random walk" on $\Z$ with steps $\Delta_i$ hits the target  $(-z_n(\tomega))$ after $n^\ast$ steps. While false as stated (this is not a random walk), this is a useful heuristic. The maps
$$\wt{\psi}_{j}^{-1}: (\wt{\psi}_{j+1}^{-1}\circ\cdots\circ\wt{\psi}_{n^\ast-1}^{-1})(A_T^\ast)\to (\wt{\psi}_{j}^{-1}\circ\cdots\circ\wt{\psi}_{n^\ast-1}^{-1})(A_T^\ast)$$  are $1/\theta$-expanding, therefore $\tomega_i=(\tpsi_{i+1}^{-1}\circ\cdots\circ\tpsi_{n-1}^{-1})(\tomega)$  are $1/\theta^{(n-i)}$-sensitive to perturbations of $\tomega\in A_T^\ast$.  Such sensitivity generates chaotic behavior, and it is not unreasonable to expect that the distribution of  
$\DS \sum_{i=1}^{n^\ast}\Delta_i$  could be approximated by the distribution of sums of ``classical'' random variables $\DS \sum_{i=1}^{n^\ast}\wh{\Delta}_i$, for which $\DS\Prob[\sum_{i=1}^{n^\ast} \wh{\Delta}_i=-z_n(\tomega_0)]$ can be found using tools from probability theory. This is our plan.

\subsubsection{Identification of the Main Difficulties}\label{s.difficulties}
This plan was carried out successfully in \cite{Avila-Dolgopyat-Duryev-Sarig} in the special case when $\alpha$ is a quadratic irrational. In this case (and in this case only) all the $\wt{\psi}_j$ are equal. In this work  we deal with the general bounded type case,  when the  $\wt{\psi}_j$ are different.

 Suppose for the moment that $\Delta_i(\tomega)$ ($\tomega\sim \mathcal U(A_T^\ast)$) were independent identically distributed (IID)  random variables. Suppose in addition that $z_{n(T)}(\tomega_0)/n^\ast(T)\to 0$ as $T\to\infty$. 
Then the asymptotic behavior of the probability on the right-hand-side of \eqref{e.Int=RW-Prob} can be obtained from the   local limit theorem for moderate deviations for sums of IID random variables. 
Pushing a bit harder, we can also handle the case when $z_{n}(\tomega_0)/n^\ast\not\to 0$, but  $|z_n(\tomega_0)/n^\ast|<\epsilon$ for $\epsilon$ sufficiently small  (this requires the local limit theorem for large deviations). However, even in the case of IID random variables, we run into problems when $|z_n(\tomega_0)/n^\ast|$ is not small, because in this scenario the behavior of the hitting probabilities is no longer universal, and is extremely sensitive to $z_n(\tomega_0)/n^\ast$.\footnote{\label{footnote-3}For example, let $S_n=X_1+\cdots+X_n$, where $X_i$ are IID random variables taking the values $\pm 1$ with equal probabilities. Then $\Prob[S_n=n-t]$ behaves differently as $n\to\infty$,  for $t\in\{-1,0,1,2,\ldots\}$.}

The main difficulty in this work is that the  local limit theorems mentioned above do not apply to the ``random walk'' $\sum\Delta_i$. This is because:
\begin{enumerate}[(1)]
\item $\Delta_i(\tomega')$ are not independent (we randomize once, not at every step);
\item $\Delta_i(\tomega')$  are not identically distributed (because the $\wt{\psi}_j$ are different);
\medskip
\item For some $\alpha$, $z_n(\tomega_0)/n^\ast$ is not sufficiently  small
for the universal large deviations estimates to hold. 
\end{enumerate}
We address these difficulties in the following way.

\medskip
\noindent
{\em Step 1: Symbolic Dynamics.\/} In part I\!I\!I we approximate $\Prob\left[\sum\limits_{i=1}^{n-1}\Delta_i=\zeta_n\right]$ by probabilities of the form 
$
\Prob\left[\sum\limits_{i=1}^{n}f_i(X_i)=\zeta_{n}\bigg|X_{n}=\zeta'_{n}\right],
$ 
where $\{X_n\}$ is a {\em Markov chain}.  

The special case of quadratic irrationals was done in \cite{Avila-Dolgopyat-Duryev-Sarig}. In this case all the $\wt{\psi}_j$ are all equal to the same automorphism $\wt{\psi}:\St\to\St$, this automorphism projects to a hyperbolic toral automorphism of $\St_0$, and the approximating Markov chain arises from the symbolic dynamical description of the normalized area measure, with respect to a suitably chosen Markov partition for $\wt{\psi}$. This Markov chain is stationary, and $f_i$ does not depend on $i$.

But in the general case, the $\wt{\psi}_j$ are different, and  there is no single Markov partition which works for all $\wt{\psi}_j$ at once. Instead we construct a  sequence of partitions $\mathfrak P_i$, with respect to which  $(\psi_0(\tomega),(\psi_1\circ\psi_0)(\tomega),\ldots)$ $(\omega\in\mathcal U(\St_0))$ is distributed like the path of  an {\em inhomogeneous Markov chain}.\index{Markov chain!inhomogeneous}  (An inhomogeneous Markov chain is a Markov chain $X_n$  whose sets of states and transition kernels depend on $n$.) 
This Markov chain leads to the approximation above.  


\medskip
\noindent
{\em Step 2: Local Limit Theorem.} If $\zeta_n/n$ is ``small enough,'' then we can try to use a local limit theorem to find the asymptotic behavior of   
$\DS\Prob\left[\sum_{i=1}^{n}f_i(X_i)=\zeta_n\bigg|X_{n}=x_n\right]$.   
We did this in the special case of quadratic irrationals in \cite{Avila-Dolgopyat-Duryev-Sarig}, when $\{X_i\}$ is  a homogeneous Markov chain and $f_i$ are all equal, using the local limit theorems for homogeneous Markov chains in \cite{Nagaev}, \cite{Guivarch-Hardy}, and \cite{Babillot-Ledrappier-Lalley}. 

But when we began our work on the  general bounded type case, we realized that there are no sufficiently general local limit theorems for  {\em inhomogeneous} Markov chains and {\em non-equal} $f_i$. We had to develop these theorems  ourselves.  This part of the project  took several years to complete, and since it is of independent interest, we decided to publish it separately in the form of a monograph \cite{Dolgopyat-Sarig-Book}. 

To apply the local limit theorems of \cite{Dolgopyat-Sarig-Book}, one needs  to check three conditions: 
\begin{enumerate}[(1)]
\item We need to verify that $f_i$ are uniformly bounded and that the  Markov chain is ``uniformly elliptic," see Definition \ref{d.uniform-ellipticity} and Proposition \ref{Prop-Doeblin}.  
\item We need to show that the variance $\DS \sum_{i=1}^{n}f_i(X_i)$ tends to infinity. 
\item We need to show {\em irreducibility}: There do {\em not} exist an integer $p>1$ and bounded sequences of functions $a_i(X_i)$, $h_i(X_i,X_{i+1})$ such that $\sum \mathrm{Var}(h_i)<\infty$, and \\
$
f_i(X_i)+a_{i+1}(X_{i+1})-a_i(X_i)+h_i(X_i,X_{i+1})\in p\Z\text{ a.s. for all }i.
$
\end{enumerate}
We check these conditions in Section \ref{ScMarkovReg}, using  a combination of geometric and number theoretic tools. The bounded type assumption is essential for the proof.

\medskip
\noindent
{\em Step 3: What to Do When $z_n(\tomega_0)/n$ Is Large?}  
We have already remarked that the local limit theorem does not apply when $z_n(\tomega_0)/n$ is ``not small enough,'' even in the most classical of situations (see footnote \ref{footnote-3}). {\em  But in our case $z_n(\tomega_0)/n$  can be ``too large."}

The first step is to  understand how small is  ``small enough'' for the local limit theorem to hold.
Let $V_n$ denote the variance of $\DS S_n=\sum_{i=1}^n f_i(X_i)$. 
We showed in \cite[Thm 7.26]{Dolgopyat-Sarig-Book} that there is a non-empty interval $(\frak{c}_-,\frak{c}_+)$ so that 
$
\Prob[S_n=\zeta_n]
$
obeys the local limit theorem for large deviations when $\frac{\zeta_n-\E(S_n)}{V_n}\to \zeta\in (\frak{c}_-,\frak{c}_+)$, and does not obey it when  $\frac{\zeta_n-\E(S_n)}{V_n}\to \zeta \not\in [\frak{c}_-,\frak{c}_+]$. (We cannot say anything  in case  $\zeta=\frak{c}_-$ or $\frak{c}_+$.) 

Thus the local limit theorem is inapplicable only when there are  subsequences such that $\frac{\zeta_{n_k}-\E(S_{n_k})}{V_{n_k}}\to \zeta\not\in (\frak{c}_-,\frak{c}_+)$. 
We show that if this happens, then  $\tomega_0$ cannot be generic for any locally finite ergodic invariant measure.

In \cite{Avila-Dolgopyat-Duryev-Sarig}, we proved this in the special case when $\beta$ is a quadratic irrational,  by showing that if $\tomega_0$ is generic for some measure $\mu_\xi$, then $\frac{\zeta_n-\E(S_n)}{V_n}$ must converge to the $\mu_\xi$-almost sure limit of  $\frac{\zeta_n-\E(S_n)}{V_n}$, which  must be inside $(\frak{c}_-,\frak{c}_+)$. The argument  relied on the  approximate exchangeability  of the joint distribution of $\{f_i(X_i)\}_{i\geq 1}$ in case $X_i$ is a homogeneous Markov chain and $f_i$ is independent of $i$, and it does not extend to the case when $X_i$ is inhomogeneous, and $f_i$ depend on $i$.

In this work we are forced to use a different approach.
 We suppose that there is a subsequence such that $\frac{\zeta_{n_k}-\E(S_{n_k})}{V_{n_k}}\to \zeta\not\in (\frak{c}_-,\frak{c}_+)$, and then use ideas from the theory of action minimizing measures to show that 
  there are two finite unions of horizontal rectangles $\mathsf F_1$, $\mathsf F_2$ and $\exists T_n\to\infty$  such that 
$
\frac{\int_0^{T_n}1_{\mathsf F_1}(\phi^t(\tomega_0))dt}
{\int_0^{T_n}1_{\mathsf F_2}(\phi^t(\tomega_0))dt}\xrightarrow[n\to\infty]{}0
$. This rules out genericity for any of the locally finite ergodic  measures. 
\vskip2mm

After all this is done, we have a dichotomy: 
{\em Either} the local limit theorem does not apply, there is a subsequence as above, and our initial condition is not generic for any $\mu_\xi$; 
{\em Or}  the local limit theorem does apply, and we have at our disposal an asymptotic formula for $\int_0^T h(\phi^t(\tomega_0))dt$. 
Our results then follow from this formula, by direct calculations.
  
\subsection{A Remark on Non-Renormalizable Surfaces.}\index{Non-renormalizable surfaces} The infinite staircase is exceptional in that it admits many affine automorphisms.
But this is not a crucial feature of the setup. The reduction to inhomogeneous Markov chains used in this work should also work for linear flows $\phi$ on translation surfaces $X_0$ without affine automorphisms,  provided the following holds: 
There exists a sequence of translation surfaces $X_i$ $(i\geq 1)$ and affine homeomorphisms 
$
X_0\xrightarrow[]{\tpsi_0}X_1\xrightarrow[]{\tpsi_1}X_2\xrightarrow[]{\tpsi_2}\cdots
$
such that

$\bullet$ $\{d\tpsi_i\}$ is pre-compact in $\mathrm{GL}(2,\R)$;

$\bullet$ $A_n:=d\tpsi_{n-1}d\tpsi_{n-2}\cdots d\tpsi_0$
 contracts the direction of $\phi$ at an exponential rate; 
 
$\bullet$ $\{A_n\}_{n\geq 1}$ is ``uniformly hyperbolic" in the sense of Proposition \ref{Prop-Contracting-Automorphisms}; and

{$\bullet$  $X_n$ admit fundamental domains $F_n$ with bounded geometry in the sense
of Lemma~\ref{Lemma-Bounded-Geometry-Gamma(M)}.}

\subsection{Plan of the Monograph.}

Part I of the paper summarizes all our main results. 

Part II provides the renormalization scheme at the base of our approach, see \S\ref{ss.renormalization-sequence}. 
We construct a sequence of automorphisms of $\St$ 
which contract the orbits of $\phi$.

Part III constructs the (inhomogeneous) Markov chain  needed to solve the probabilistic formulation of the problem we presented in \S\ref{ss.probabilistic-formulation},\S\ref{ss.RW-heuristic}. 

Part IV uses these limit theorems to prove the main dynamical results of the paper. 
For the location of proofs of the theorems stated above, see the table  on page \pageref{page-proof-locator}.

Part V contains appendices with background material on even CFE,  denseness of finite line segments  on the torus, and  the results of \cite{Dolgopyat-Sarig-Book} on the local limit theorems for inhomogeneous Markov chains.     Appendix \ref{Appendix-Twist} contains (novel) calculations leading to a formula for $z_n(\tomega_0)$ in terms of the even CFE of $\beta$.

\part{Renormalization}
\label{p.renormalization}
Let $\alpha\in (0,1)$ be an irrational number of bounded type.\index{Canonical matrix product}
In the previous section we related the behavior of $D_n(\alpha)$ 
to the behavior of the linear flow $\phi:\St\to\St$,  with velocity vector ${\beta\choose 1}$,  $\beta:=2\alpha-1\in (-1,1)$.

Suppose  $\wt{\psi}:\St\to\St$ is a homogeneous automorphism with derivative matrix $A$, then $\varphi^t:=\wt{\psi}\circ \phi^t\circ \wt{\psi}^{-1}$ is the linear flow  on $\St$, with velocity vector $A{\beta\choose 1}$. In this part of the paper we implement the heuristic in \S\ref{ss.renormalization-sequence}, and construct a uniformly hyperbolic sequence of homogeneous automorphisms $\wt{\psi}_n$ with derivative matrices $A_n$, so that 
$\|A_n{\beta\choose 1}\|\to 0$. The challenge  is to find integer 
matrices $A_n\in\Gamma(2)$ (necessarily with some big entries) which contract a given vector. The solution comes from the theory of even continued fractions, parts of which  we summarize in Appendix \ref{AppEven}.

\section{Canonical Matrix Products}\label{ScReno}
Let $\beta=[0;2m_1,2m_2,\ldots]$ be the even continued fraction expansion (CFE) of $\beta$.
The {\em canonical matrix products} of $\beta$ are defined by 
\begin{equation}\label{e.A-k}
A_k(\beta):=\left(
       \begin{array}{cc}
         1 & -2m_{2k} \\
         0 & 1 \\
       \end{array}
     \right)
     \left(
       \begin{array}{cc}
         1 & 0 \\
         -2m_{2k-1} & 1 \\
       \end{array}
     \right)\cdots
      \left(
       \begin{array}{cc}
         1 & -2m_{2} \\
         0 & 1 \\
       \end{array}
     \right)
     \left(
       \begin{array}{cc}
         1 & 0 \\
         -2m_{1} & 1 \\
       \end{array}
     \right).
\end{equation}
Since $\Gamma(2):=\{A=\bigl({\tiny\begin{array}{cc} a & b\\
c & d
\end{array}
}\bigr):\det(A)=1, A=\bigl({\tiny\begin{array}{cc} 1 & 0\\
0 & 1
\end{array}
}\bigr)\mod 2 \}$ is a group,  $A_{k}(\beta)\in\Gamma(2)$ for all $k$. In this section we collect some of the properties of $A_k(\beta)$.

\begin{lem}
The following identities hold:
\begin{align}
   A_{k}(\beta)  &=
     \left(
       \begin{array}{cc}
         -2m_{2k} & 1 \\
         1 & 0 \\
       \end{array}
     \right)
     \left(
       \begin{array}{cc}
         -2m_{2k-1} & 1 \\
         1 & 0 \\
       \end{array}
     \right)\cdots
     \left(
       \begin{array}{cc}
         -2m_2 & 1 \\
         1 & 0 \\
       \end{array}
     \right)
     \left(
       \begin{array}{cc}
         -2m_1 & 1 \\
         1 & 0 \\
       \end{array}
     \right); \label{e.A-k-prime}\\
 A_{k}(\beta)&=
\left[\left(
       \begin{array}{cc}
         0 & 1 \\
         1 & 2m_{1}\\
       \end{array}
     \right)
     \left(
       \begin{array}{cc}
         0 & 1 \\
         1 & 2m_{2} \\
       \end{array}
     \right)\cdots
     \left(
       \begin{array}{cc}
         0 & 1 \\
         1 & 2m_{2k-1} \\
       \end{array}
     \right)
     \left(
       \begin{array}{cc}
         0 & 1 \\
         1 & 2m_{2k} \\
       \end{array}\right)\right]^{-1} \label{e.A-k-dagger} ;\\
 A_k(\beta)&= \left(
       \begin{array}{cc}
        p_{2k-1}  & p_{2k}  \\
         q_{2k-1} & q_{2k}  \\
       \end{array}\right)^{{-1}}\!\!\!\!,\text{ where }\frac{p_{j}}{q_{j}}=[0;2m_{1},\ldots,2m_{j}], \text{ and } (p_{j},q_{j})=1.
       \label{e.A-k-pq}
\end{align}
\end{lem}

\begin{proof}
\eqref{e.A-k-prime} is because
$\left(
\begin{array}{cc}
1 & a \\
0 & 1 \\
\end{array}
\right)
\left(
\begin{array}{cc}
1 & 0 \\
b & 1 \\
\end{array}
\right)=\left(
\begin{array}{cc}
ab+1 & a \\
b & 1 \\
\end{array}
\right)=\left(
\begin{array}{cc}
a & 1 \\
1 & 0 \\
\end{array}
\right)
\left(
\begin{array}{cc}
b & 1 \\
1 & 0 \\
\end{array}
\right)
$. 
\eqref{e.A-k-dagger} follows from \eqref{e.A-k-prime} and the identity 
$ \left(
\begin{array}{cc}
0 & 1 \\
1 & a \\
\end{array}
\right)^{-1}\!\!\!=\!\!\left(
\begin{array}{cc}
-a & 1 \\
1 & 0 \\
\end{array}
\right)$.  By  \cite[Ch.~1]{Khintchine-Continued-Fractions},  the reduced form $p_{k}/q_{k}$ of the fraction $[0;2m_1,\ldots,2m_{k}]$ is given recursively by 
$$
\left(
       \begin{array}{cc}
        p_{k}  & p_{k+1}  \\
         q_{k} & q_{k+1}  \\
       \end{array}\right)=\left(
       \begin{array}{cc}
        p_{k-1}  & p_{k}  \\
         q_{k-1} & q_{k}  \\
       \end{array}\right)\left(
       \begin{array}{cc}
        0  & 1  \\
         1 & 2m_{k+1}  \\
       \end{array}\right)\ ;\  \left(
       \begin{array}{cc}
        p_{-1}  & p_{0}  \\
         q_{-1} & q_{0}  \\
       \end{array}\right)=
       \left(
       \begin{array}{cc}
        1  & 0  \\
         0 & 1 \\
       \end{array}\right).
$$
The matrix product inside the square brackets in  \eqref{e.A-k-dagger} satisfies  the same recurrence. So $A_{k}(\beta)=\left(
       \begin{array}{cc}
        p_{2k-1}  & p_{2k}  \\
         q_{2k-1} & q_{2k}  \\
       \end{array}\right)^{-1}$.
\end{proof}

By Khintchine's theorem,  an irrational number has bounded type iff it is badly approximable by rationals in the sense of \eqref{e.badly-approx}.
Since $\alpha$ has bounded type, $\alpha$ is badly approximable. Therefore $\beta=2\alpha-1$ is badly approximable, and therefore  $\beta$ has bounded type. By   Lemma \ref{l.BT-bounded-pm-two},   the following numbers are finite:
\begin{align}
M_1(\beta):=&\text{the maximal absolute value of the digits of the even CFE of $\beta$},\notag\\
M_2(\beta):=&\text{the maximal length of words  of the form $(2,-2,\ldots,\pm 2)$}\notag\\
&\text{ or $(-2,2,\ldots,\pm 2)$ in the even CFE of $\beta$}, \notag\\
M(\beta):=&\max\{M_1(\beta),M_2(\beta)\}+1\ (\text{necessarily, $M(\beta)\geq 3$}).\label{e.M-beta}
\end{align}
\begin{lem}\label{l.disonance}
Suppose $M<\infty$,
 and let $K$ be a compact subset of $(-1,1)$. Then there are constants $C>0, 0<\eta<1$ which only depend on $M$ and $K$ such that if $M(\beta)\leq M$, $\rho\in K$, $k\geq 1$, and $\beta=[0;2m_1,2m_2,\ldots]$ where $m_i\in\Z\setminus\{0\}$, then 
$$
\left\|\left(
\begin{array}{cc}
0 & 1 \\
1 & 2m_1 \\
\end{array}
\right)\left(
\begin{array}{cc}
0 & 1 \\
1 & 2m_2 \\
\end{array}
\right)\cdots
\left(
\begin{array}{cc}
0 & 1 \\
1 & 2m_k \\
\end{array}
\right){\rho\choose 1}\right\|\geq C^{-1}\eta^{-k}.
$$
\end{lem}

We need some preparations to the proof. Henceforth ${\tiny \left(
\begin{array}{cc}
a & b \\
c & d \\
\end{array}
\right)}\cdot z:= \frac{az+b}{cz+d}$. Then ${\tiny \left(
\begin{array}{cc}
a & b \\
c & d \\
\end{array}
\right)}\cdot \left[{\tiny \left(
\begin{array}{cc}
a' & b' \\
c' & d' \\
\end{array}
\right)}\cdot z\right]=\left[{\tiny \left(
\begin{array}{cc}
a & b \\
c & d \\
\end{array}
\right)\!\! \left(
\begin{array}{cc}
a' & b' \\
c' & d' \\
\end{array}
\right)}\right]\cdot z$.  The 
{\em slope} of  $\vec{v}:={x\choose y}\neq {0\choose 0}$ is $\sigma(\vec{v}):=y/x\in [-\infty,\infty]$.  The  {\em  co-slope} of $\vec{v}$ is  $\hsigma(\vec{v}):=1/\sigma(\vec{v})=x/y. $ \index{Slope}\index{Co-slope}
Then 
\begin{equation}\label{e.proj-action}
\hsigma[A\vec{v}]=A\cdot \hsigma(\vec{v})  \text{ for every  } A\in \mathrm{GL}(2,\R).
\end{equation}
\noindent
Let 
$
 v_{2m}:[-1,1]\to \left[\frac{1}{2m+1},\frac{1}{2m-1}\right]\ , \ v_{2m}(t)=\bigl({\tiny
\begin{array}{cc}
0 & 1 \\
1 & 2m \\
\end{array}}
\bigr)\cdot t=\frac{1}{2m+t}.
$
\begin{proof}[Proof of Lemma \ref{l.disonance}] 
Let $B_k\!\!:=\!\!{\tiny \left(
\begin{array}{cc}
0 & 1 \\
1 & 2m_1 \\
\end{array}
\right)\left(
\begin{array}{cc}
0 & 1 \\
1 & 2m_2 \\
\end{array}
\right)\cdots
\left(
\begin{array}{cc}
0 & 1 \\
1 & 2m_k \\
\end{array}
\right)}$. By \eqref{e.proj-action},
$$
\hsigma\left[B_k{\tiny \left(
                   \begin{array}{c}
                     \rho \\
                     1 \\
                   \end{array}
                 \right)
}\right]=(v_{2m_1}\circ\cdots\circ v_{2m_k})(\rho).
$$
By Lemma \ref{l.v-contracts}, there exists $C>0$ and $0<\eta<1$ which only depend on $M$ such that $|(v_{2m_1}\circ\cdots\circ v_{2m_k})'|\leq C\eta^k$ on $[-1,1]$. It follows that for every $\rho,\rho'\in [-1,1]$,
\begin{equation}
\label{DeltaSigma}
\left|\hsigma\left[B_k{\tiny \left(
                   \begin{array}{c}
                     \rho \\
                     1 \\
                   \end{array}
                 \right)
}\right]-\hsigma\left[B_k{\tiny \left(
                   \begin{array}{c}
                     \rho' \\
                     1 \\
                   \end{array}
                 \right)
}\right]\right|\leq C\eta^k.
\end{equation}
Therefore, the angle between $B_k{\rho\choose 1}$ and $B_k{\rho'\choose 1}$ tends to zero uniformly exponentially fast as $k\to\infty$, whenever $\rho,\rho'\in (-1,1)$,

Now suppose $\rho\in K$, where $K$ is a compact subset of $(-1,1)$. Choose $\eps_0$ so small that the closed  $\eps_0$-neighborhood of $K$ is inside $(-1,1)$, and take $\rho':=\rho+\eps_0$ and $\rho'':=\rho-\eps_0$.
Let $Q$ denote the parallelogram generated by the vectors $\frac{1}{2}{\rho'\choose 1}$ and $\frac{1}{2}{\rho''\choose 1}$. Since $\det B_k=1$, the parallelogram $B_k(Q)$ has the same area as $Q$, but the angle between its generating vectors $B_k{\rho'\choose 1}, B_k{\rho''\choose 1}$ tends to zero uniformly exponentially as $k\to\infty$, at a rate which only depends on $M$.

 Necessarily, the length of the  principal diagonal $B_k{1\choose\rho}$ tends to infinity at an exponential rate, which only depends on $M$.
\end{proof}

\begin{lem}\label{l.A-contracts}Suppose $M<\infty$, then there are $C>0$ and $0<\eta<1$ such that if $\beta\in (-1,1)\setminus\Q$ and $M(\beta)\leq M$, then 
$\|A_k(\beta){\beta\choose 1}\|\leq C \eta^k$ for all $k\geq 0$.
\end{lem}

\begin{proof}
Write $\beta=[0;2m_1,2m_2,\ldots,2m_{2k},r_{2k+1}]$, where $r_{2k+1}=[2m_{2k+1};2m_{2k+2},\ldots]$. 
 In the notation of the proof of the previous lemma,
$$
\hsigma\left[{\beta\choose 1}\right]=\beta=(v_{2m_1}\circ\cdots\circ v_{2m_{2k}})
\left(\frac{1}{r_{2k+1}}\right)=\hsigma\left[B_k {1\choose r_{2k+1}}\right]=\hsigma\left[B_k {1/r_{2k+1}\choose 1}\right].
$$
Necessarily, ${\beta\choose 1}$ and $B_k {1/r_{2k+1}\choose 1}$ are collinear, and therefore 
$$
{\beta\choose1}=\pm \frac{\sqrt{1+\beta^2}}{\|B_k{1/r_{2k+1}\choose 1}\|} B_k{1/r_{2k+1}\choose 1}.
$$

By Lemma \ref{l.BT-compact}, there is a compact set $K\subset (-1,1)$ which depends only on $M$ so that $1/r_{2k+1}\in K$ for all $k$.

By Lemma \ref{l.disonance}, there are constants $C>0$ and $0<\eta<1$ which only depend on $M$ so that
$
\left\|B_k{1/r_{2k+1}\choose 1}\right\|\geq C^{-1}\eta^{-k}.
$
So
$
\left\|B_k^{-1}{\beta\choose 1}\right\|\leq C\eta^k \sqrt{1+\beta^2} \sqrt{1+r_{2k+1}^{-2}}\leq 2C\eta^k.
$
By \eqref{e.A-k-dagger},  $B_k^{-1}=A_k$, therefore  $\left\|A_k{\beta\choose 1}\right\|\leq C_1\eta^k,$ where $C_1$ depends only on $M$.
%
\end{proof}

Lemma \ref{l.A-contracts} produces matrices in $\Gamma(2)$ which contract ${\beta\choose 1}$. We will now examine their action on other vectors. We need the following tool \cite{Schweiger-Book}:
\begin{defi}
The {\em Even Gauss Map} $E:(-1,1)\to (-1,1)$, is given by \index{Even Gauss map}
\begin{equation}
\label{EvenGauss}
E(x):=\frac{1}{x}-2m\text{ on }\left[\frac{1}{2m+1},\frac{1}{2m-1}\right)\ \ (m\in\Z).
\end{equation}
\end{defi}
\noindent 
This is designed so that $E([0;2m_1,2m_2,2m_3,\ldots]=[0;2m_2,2m_3,2m_4,\ldots]$.

\begin{lem}\label{l.position}
There exists a constant $0<\eps(M)<1$ as follows.
Suppose $\beta$ and $\gamma$ are irrationals in $(-1,1)$ so that  $M(\beta), M(\gamma)\leq M$, then
\begin{enumerate}[(1)]
\item $
1+\eps(M)<\left|\sigma\left[A_k(\beta){\beta\choose 1}\right]\right|<1/\eps(M)$,
\item $\eps(M)<\left|\sigma\left[A_k(\beta){1\choose -\gamma}\right]\right|<1/(1+\eps(M))$,
\item $
\left|\sin\measuredangle(A_k(\beta){\beta\choose 1},A_k(\beta){1\choose -\gamma})\right|>\eps(M).
$
\end{enumerate}
\end{lem}
\begin{proof}
Define $\hsigma_k^s:=\hsigma\left[A_k(\beta){\beta\choose 1}\right]$. By \eqref{e.A-k-dagger},
\begin{align}
&\hsigma_k^s
     =(v_{2m_1}\circ\cdots\circ v_{2m_{2k}})^{-1}(\beta)=(v_{2m_1}\circ\cdots\circ v_{2m_{2k}})^{-1}(\lim_{n\to\infty} (v_{2m_1}\circ\cdots\circ v_{2m_n})(0))\notag\\
     &=\lim_{n\to\infty} (v_{2m_{2k+1}}\circ\cdots\circ v_{2m_n})(0)=[0;2m_{2k+1},2m_{2k+2},\ldots]=E^{2k}(\beta)=:\beta_k'.\label{e.sigma-s} 
\end{align}
In particular, $M(\hsigma^s_k), M(\frac{1}{\hsigma^s_k}-2m_{2k+1})\leq M$.
By  Lemma \ref{l.BT-compact}, there exists a number  $0<\delta(M)<1$ which only depends on $M$ such that $\delta(M)<|\hsigma_k^s|<1-\delta(M)$ for all $k$. Passing from the co-slope to the slope, we obtain the first part of the lemma. 

Suppose $\gamma=[0;2\ell_1,2\ell_2,\ldots]$ where $\ell_i\in \Z\setminus\{0\}$ and $M(\gamma)\leq M$. Set
$\hsigma^u_k:=\hsigma\left[A_k(\beta){1\choose -\gamma}\right]$. Let $u_{2m}(z):=\bigl({\tiny \begin{array}{cc}
     2m & 1\\ 1 & 0
     \end{array}}\bigr)\cdot z= 2m+\frac{1}{z}$.  By \eqref{e.A-k-prime},
\begin{align}
&\hsigma^u_k
     =(u_{-2m_{2k}}\circ\cdots\circ u_{-2m_1})(-\gamma^{-1})
     =[-2m_{2k};-2m_{2k-1},\ldots,-2m_2,-2m_1,-\gamma^{-1}]\notag\\
     &=[-2m_{2k};-2m_{2k-1},\ldots,-2m_2,-2m_1,-2\ell_1,-2\ell_2,\ldots]\label{e.sigma-u}\\
     &=-[2m_{2k};2m_{2k-1},\ldots,2m_2,2m_1,2\ell_1,2\ell_2,\ldots].\notag
\end{align}
We see that $M(\hsigma_k^u+2m_{2k})\leq 2 M$ (the factor 2 is due to the possible increase in the length of the strings $(2,-2,\cdots,2,-2)$). Therefore there exists a constant $\delta'(M)$ which only depends on $M$ so that
$
\hsigma^u_k\in -2m_{2k}+(-1+\delta',1-\delta').
$
Note that $|2m_{2k}|\leq M-1$, and therefore $|\hsigma^u_k|\leq M$. Next, note that $|2m_{2k}|\geq 2$. If $2m_{2k}\leq -2$, then $\hsigma^u_k\geq 1+\delta'$, and if $2m_{2k}\geq 2$, then $\hsigma^u_k\leq -1-\delta'$. In all cases,
$
1+\delta'\leq |\hsigma_k^u|\leq M.
$
Passing to the slope, we obtain the second part of the  lemma.

The third part of the lemma is an easy consequence of parts one and two.
\end{proof}

\section{Associated Fundamental Domains}\label{ss.Pi(n)} 
$\St_0$ is isometric to $\R^2/G$, where \index{Associated fundamental domains}
 $G$ is the group of translations 
\begin{equation}\label{e.Group-G}
G=\<T,S\>, \ T(x,y)=(x+1,y+1),\  S(x,y)=(x+2, y). 
\end{equation}
Figure  \ref{Figure-Staircase}(c) in \S \ref{s.infinite-staircase}
shows a fundamental domain for $G$.  Here we build  other fundamental domains, with sides parallel to $\wt{E}^s_k:=\Span\{A_k(\beta){\beta\choose 1}\}$ and 
$\wt{E}^u_k:=\Span\{A_k(\beta){1\choose -\beta}\}$. In Part \ref{p.Symbolic-Dynamics}, we will use these domains to build Markov partitions.

\begin{figure}
  \centering
  \fbox{\includegraphics[angle=90,width=12cm]{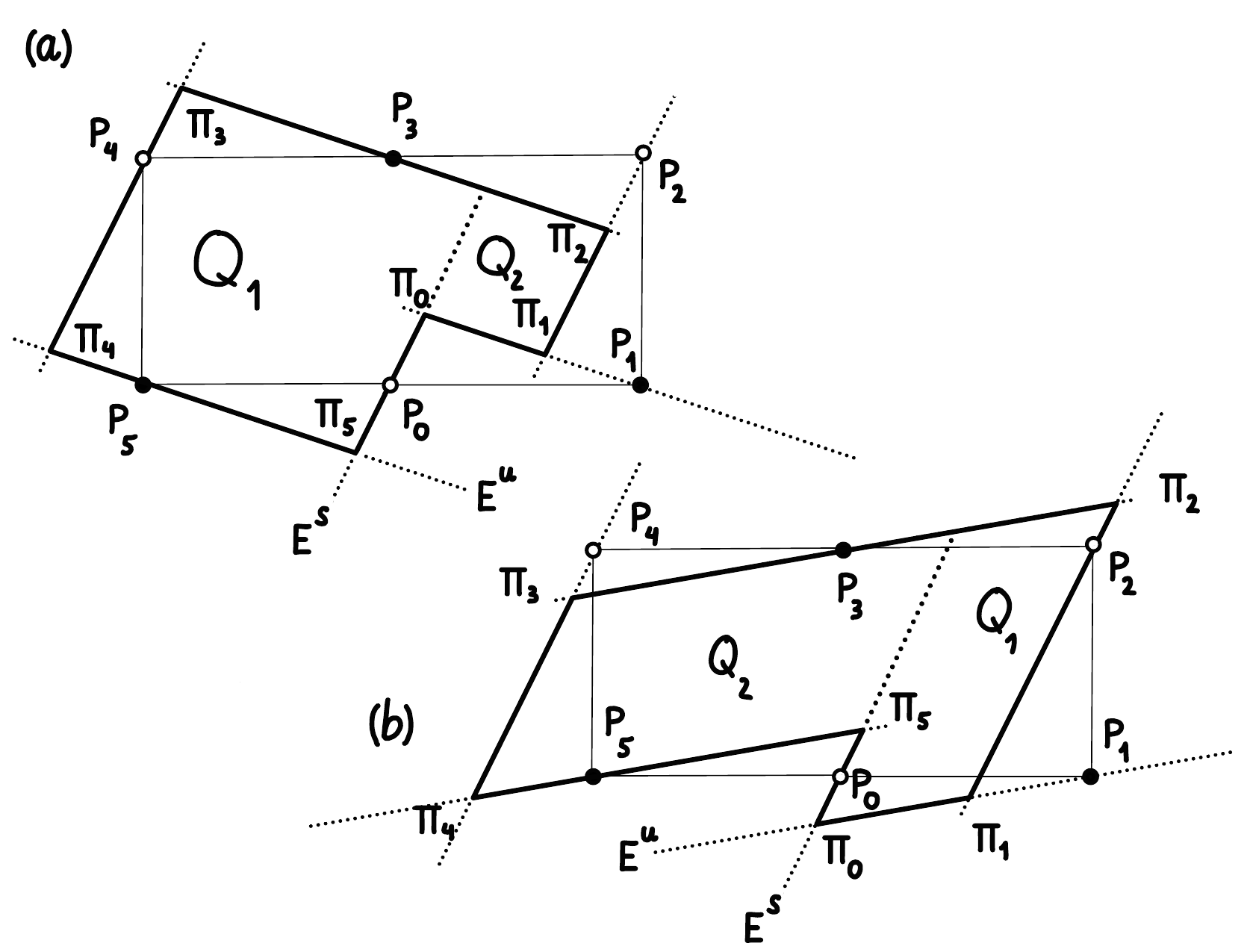}}\\
  \caption{The Polygon $\Pi$ when $\sigma^s>1$ and   (a) the slopes of $E^u, E^s$ have opposite signs; (b) the slopes of $E^u$, $E^s$ have the same sign.}\label{Figure-Pi(n)}
\end{figure}

Suppose $E^s, E^u$ are two linear subspaces in $\R^2$ with respective slopes $\sigma^s, \sigma^u$ such that $|\sigma^s|>1$ and $|\sigma^u|<1$ (e.g. $\wt{E}^s_k, \wt{E}^u_k$, see Lemma \ref{l.position}).  Let 
\begin{align*}
P_0&=(1,0), \ 
P_1=(2,0), \  P_2=(2,1), \  P_3=(1,1), \ P_4=(0,1), \ P_5=(0,0), \ P_6=P_0\\
L_i&:=\begin{cases}
\text{the line through $P_i$ in direction $E^s$, when  $i$ is  even}\\
\text{the line through $P_i$ in direction $E^u$, when  $i$ is odd}
\end{cases} \subset\R^2\ \ (0\leq i\leq 6)\\
\Pi_i&:=\text{ the unique intersection of $L_i$ with $L_{i+1}$ when $0\leq i\leq 5$, and }\Pi_6:=\Pi_0.
\end{align*}
Figure \ref{Figure-Pi(n)} describes the construction when $\sigma^s>1$. 
The case when $\sigma^s<-1$ follows by reflection across the axis  $P_0 P_3$.

\begin{lem}\label{l.jordan-curve}
 The polygonal path $(\Pi_0,\Pi_1,\cdots,\Pi_5,\Pi_0)$ is a Jordan curve, therefore it bounds a  polygon $\Pi(E^u,E^s)$ with vertices $\Pi_i$.
\end{lem}
\begin{proof}
Let  $S_i$ denote open line segment from $\Pi_i$ to $\Pi_{i+1}$. To prove the lemma, we need to show that $\ov{S}_i\cap\ov{S}_j$ are pairwise disjoint when $j\neq i\pm 1\mod 6$, and that $\ov{S}_i\cap\ov{S}_j=\{\Pi_{j}\}$ when $j=i+1\mod 6$. 

By construction, $\Pi_{i},\Pi_{i+1}$ both lie on $L_{i+1}$. Therefore $S_{i}\subset L_{i+1}$. 
If $i=j\mod 2$, then $\ov{S}_i, \ov{S}_{j}$ are disjoint because $L_{i}, L_{j}$ are disjoint parallel lines. 
If $j=i+1\mod 6$ then  $\ov{S}_i\cap\ov{S}_{j}=\{\Pi_{j}\}$, because $L_{i}\cap L_{j}=\{\Pi_j\}$.
It remains to consider the case when  $j-i=3 \mod 6$.

We  call the lines $L_{i}$ with slope  $\sigma^s$ ``almost vertical,'' and lines with slope $\sigma^u$ ``almost horizontal." The ``left,'' ``middle,'' and ``right'' almost vertical lines are 
$L_{4},L_{0},L_{2}$. The   ``top,'' ``middle,'' and ``bottom'' almost horizontal lines are   $(L_{3},L_{1},L_{5})$ when $\sigma^u<0$, and 
$
(L_{3},L_{5},L_{1})$ when  $\sigma^u_n>0$.

\medskip
\noindent
{\em Claim.\/} If  $j-i=3\mod 6$, then  at least one of $\{L_{i+1},L_{j+1}\}$ is a not a middle line.

\medskip
\noindent
{\em Proof of the Claim\/}:  Since $j-i=3\mod 6$,  $i,j$ have different parity. Therefore one of $L_{i+1},L_{j+1}$   is almost vertical, and the other is almost horizontal. Without loss of generality, $L_{i+1}$ is almost vertical, and $L_{j+1}$ is almost horizontal.   If $L_{i+1}=L_0$, then $i=5$,  $j=2$, and $L_{j+1}=L_3$ a top horizontal line. If $L_{i+1}\neq L_0$, then $L_{i+1}$ is a left or right almost vertical line. The claim follows.

\medskip
\noindent
Suppose now that $j-i=3\mod 6$, and assume by contradiction that $\ov{S}_i\cap\ov{S}_j\neq \emptyset$. By the claim,  one of $S_i, S_j$, say $\ov{S}_i$,  lies on a non-middle line $L$. A segment  on a non-middle line $L$ cannot intersect any other segment except at the boundary,  otherwise there would be vertices $\Pi_{i}$ on different sides of  $L$, and $L$ would have to be a middle line. So the  $\ov{S}_i$ intersects $\ov{S}_j$ at $\Pi_i$ or 
$\Pi_{i+1 \mod 6}$. Necessarily, $\ov{S}_j$ intersects $\ov{S}_{i+1 \mod 6}$ or $\ov{S}_{i-1\mod 6}$. But  $j=(i\pm 1)\mod 2$, so as we saw above, this cannot happen.
\end{proof}

\begin{lem}\label{Lemma-Q}
$L_0$ cuts $\Pi$  into two parallelograms $\wt{Q}_1, \wt{Q}_2$, with sides parallel to $E^u, E^s$, such that $\wt{Q}_1$ contains $P_{0}$ in its boundary, and $\wt{Q}_2$  does not contain any $P_{i}$ with $i$ even in its boundary. \index{Associated fundamental domains!partition into parallelograms}
\end{lem}
\noindent
\begin{proof}
This follows easily from the fact that $|\sigma^s|>1$ and $|\sigma^u|<1$.
\end{proof}

\begin{lem}\label{Lemma-Fundamental-Domain}
$\Pi$ is a fundamental domain for the action of $G$ on $\R^2$. Thus, it  has side pairing (perhaps after further division of its sides) which turns it into a model of $\St_0$. 
\end{lem}

\begin{figure}
  \centering
  \fbox{\includegraphics[angle=90,width=12cm]{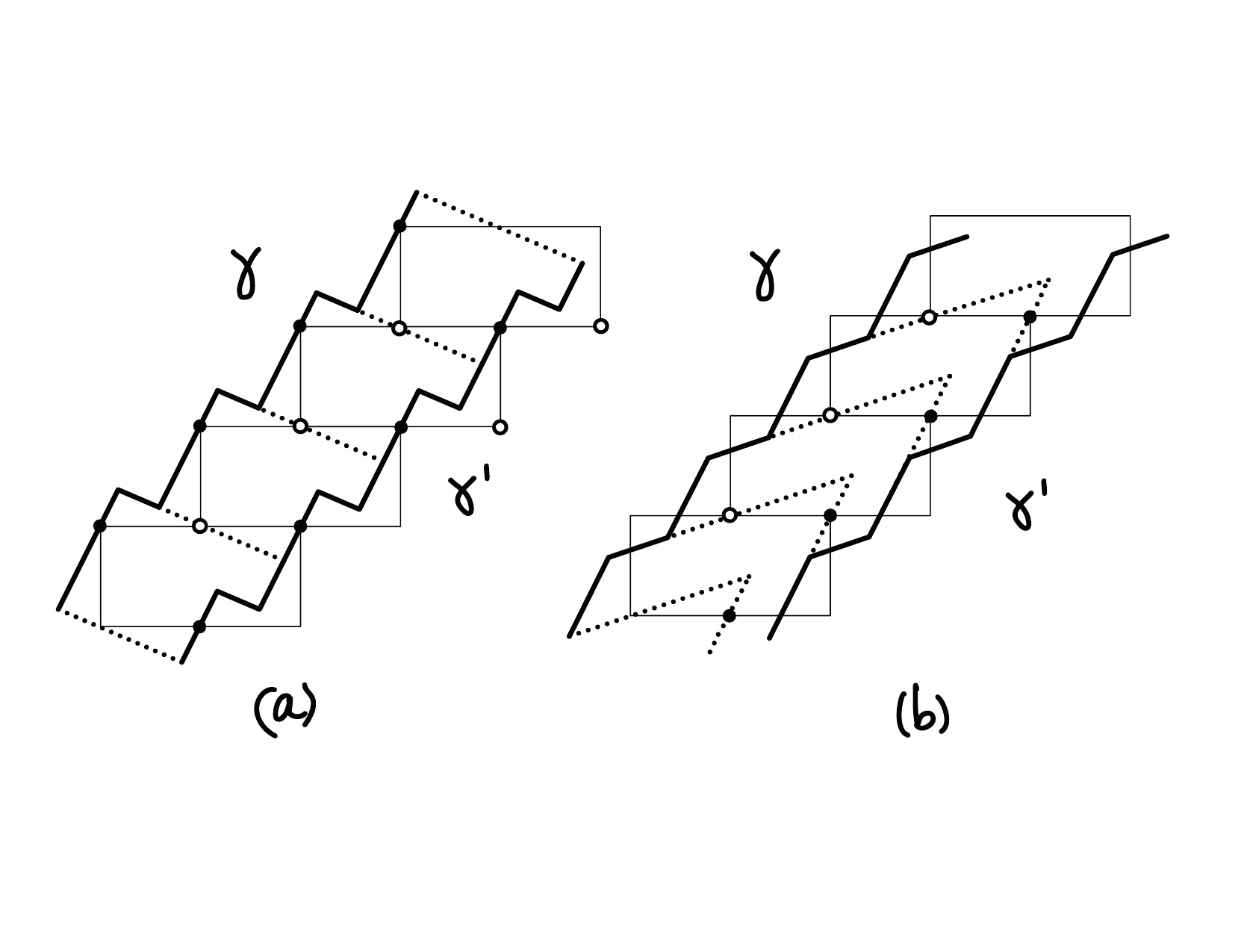}}\\
  \caption{The strip $\Sigma$: (a) Different slope sign; (b) Equal slope sign}\label{Figure-Tesselation}
\end{figure}

\begin{proof} 
Elementary geometric considerations show that  $T^k[\mathrm{int}(\Pi)]$ $(k\in\Z)$ are pairwise disjoint, and that 
$\DS \Sigma:=\bigcup_{k\in\Z}T^k[\Pi]$ is  a strip with two boundary  curves $\gamma,\gamma'$ made from a linear interpolation of the left (resp. right) $E^s$-sides of $\Pi$,  see Figure \ref{Figure-Tesselation}.

We claim that $\gamma'=S(\gamma)$. To see this, consider the polygon  $\widehat{\Pi}$ built as above, but with $P_i$ replaced by $S(P_i)$, and let  $\widehat{\Sigma}:=\bigcup_{k\in\Z}T^k(\widehat{\Pi})$.
Since $P_1=S(P_5)$ and $P_2=S(P_4)$, the  right side $\gamma'$ of $\Sigma$ equals the left side of $\widehat{\Sigma}$. Since $\widehat{\Pi}=S(\Pi)$ and $T\circ S=S\circ T$,  $\wh{\Sigma}=S(\Sigma)$, and the left side of $\wh{\Sigma}$ equals $S(\gamma)$. It follows that $\gamma'=S(\gamma)$.

Since $\gamma'=S(\gamma)$, $S^\ell(\mathrm{int}(\Sigma))$ $(\ell\in\Z)$ are pairwise disjoint, and  
$\DS \bigcup_{\ell\in\Z} S^\ell(\Sigma)=\R^2$. Thus $S^\ell T^k[\mathrm{int}(\Pi)]$ $(k,\ell\in\Z)$ are pairwise disjoint, and $\DS \bigcup_{k,\ell\in\Z}S^\ell T^k(\Pi)=\R^2$. Since $T\circ S=S\circ T$,  $\{T^\ell S^k:k,\ell\in\Z\}=G$, and  $\Pi$ is a fundamental domain for $G$. 

The action of $G$ on $\partial\Pi$ determines a side pairing (perhaps after further subdivision) which makes $\Pi$ congruent to $\St_0$.
\end{proof}

\begin{lem}[Bounded Geometry]\label{Lemma-Bounded-Geometry-Gamma(M)}\index{Bounded geometry}
There exists a constant $\Gamma(M)\!\!>\!\!1$ as follows. Suppose $\beta$ and $\gamma$ are irrationals in $(-1,1)$ so that $M(\beta),M(\gamma)\!\!\leq\!\! M$. Let
$\wt{E}^s_k:=\Span\{A_k(\beta){\beta\choose 1}\}$, and $\wt{E}^u_k:= \Span\{A_k(\beta){1\choose -\gamma}\}$. Then for each $k\geq 0$,
\begin{enumerate}[(1)]
\item $\wt{\Pi}:={\Pi}(\wt{E}^s_k,\wt{E}^u_k)$ is polygon bounded by a Jordan curve;
\item All the angles  of $\wt{\Pi}$ are at least $\Gamma(M)^{-1}$--away from $\pi\Z$;
\item $\Gamma(M)^{-1}\leq \|\wt{\Pi}_i-P_j\|\leq \Gamma(M)$ for all $0\leq i,j\leq 5$;
\item All the sides of $\wt{\Pi}$ have lengths in $[\Gamma(M)^{-1},\Gamma(M)]$;
\item $\wt{Q}_1,\wt{Q}_2$ contain open balls of radius $1/\Gamma(M)$.
\end{enumerate}
\end{lem}
\begin{proof}
Let $\sigma^s_k:=\sigma[A_k(\beta){\beta\choose 1}]$ and  
$\sigma^u_k:=\sigma[A_k(\beta){1\choose -\gamma}]$. 
 By Lemma \ref{l.position}, $|\sigma^s_k|>1$, and $|\sigma^u_k|<1$.
 Part (1) now follows
 from Lemma \ref{l.jordan-curve}. 

The angles of $\wt{\Pi}$ are uniformly bounded away from $\pi\Z$ by a constant which only depends on $M$, because of Lemma \ref{l.position}(3).
This is part (2).

The quantities $\dist(\wt{\Pi}_i,P_j)$ and $\dist (\wt{\Pi}_i,\wt{\Pi}_{i+1})$ can be put in the form $g_{i,j}(\sigma^s_k,\sigma^u_k)$ and $h_i(\sigma^s_k,\sigma^u_k)$, where $g_{ij}, h_i$  are well-defined, positive,  and continuous functions  on
$$
X:=\{(\sigma^s,\sigma^u)\in\R^2: \sigma^u\neq\sigma^s, |\sigma^u|\cdot|\sigma^s|\neq 0\}.
$$
By  Lemma \ref{l.position}, $(\sigma^s_k,\sigma^u_k)$ belong to the compact set 
\begin{multline*}
 C:=\{(\sigma^s,\sigma^u)\in\R^2: |\sigma^s|\in [(1+\epsilon(M))^{-1},\epsilon(M)^{-1}], |\sigma^u|\in [\epsilon(M),(1-\epsilon(M))^{-1}]\}.
\end{multline*}
 On this set,  $g_{ij}$ and $h_i$ are bounded away from zero and infinity by constants which  depend only on $M$.  Parts (3) and (4) follow. 
 
 Part (5) follows from Parts (2) and (4). 
\end{proof}
We call $\Gamma(M)$ the {\em geometric constant}. \index{Geometric constant}

\section{Renormalizing Automorphisms}
\subsection{Denseness Bounds}\label{ss.density-const}
A line on $\mathsf{St}_{0}$ has {\em slope}\index{Slope} $r$, if it has slope $r$ in the horizontal rectangle model for $\mathsf{St}_{0}$ in Figure \ref{Figure-Staircase}(b) 
(\S \ref{s.infinite-staircase}), with  ${0\choose 1}=$ direction ``up").
  Let 
\begin{equation}\label{e.BT}
\BT(M):=\{[a_0; a_1,a_2,\cdots]: a_0\in\Z, a_1,a_2,\ldots\in\N, \sup|a_i|< M\}.
\end{equation}

\begin{lem}
\label{LmDenseWs}
 For every $M>1$ and $0<\epsilon<1$ there exists $T$ as follows: If  $r\in\BT(M)$,  
   then every line in $\St_0$ with slope $r$ and length $T$  is $\epsilon$-dense in $\St_0$.
\end{lem}
\noindent
The proof is standard; we give it for completeness in Appendix \ref{appendix-denseness}

Henceforth we fix $T_0(M)$ so large that any line with slope $r\in\BT(M)$ and length $T_0(M)$ is $1/\Gamma(M)$-dense in $\St_0$. We call $T_0(M)$ the {\em denseness constant}. \index{Denseness constant}

\subsection{Renormalizing Automorphisms}\label{SSRenorm}
Suppose $\alpha\in (0,1)$ has  bounded type, and $\beta:=2\alpha-1$. Fix $M:=M(\beta)+1$ with $M(\beta)$ as in \eqref{e.M-beta}, and let $\Gamma(M)$ and $T_0(M)$ denote the  geometric constant and the denseness constant defined at the end of Sections~\ref{ss.Pi(n)} and \ref{ss.density-const}. Fix a number $\lambda>1$ so large that \index{Renormalizing automorphisms}\index{Renormalization}
\begin{equation}\label{e.choice-of-lambda}
\lambda>\max\{\Gamma(M)^2,\Gamma(M)T_0(M)\}.
\end{equation}

\begin{prop}\label{Prop-Contracting-Automorphisms}
There exists $\epsilon(\alpha), K(\alpha)>0$ and maps  $\wt{\psi}_n:\mathrm{St}\to \mathrm{St}$ as follows:

\begin{enumerate}[(1)]
\item \textbf{\textit{Automorphisms:}} $\wt{\psi}_n$ are homogeneous affine staircase automorphisms with derivatives $A_n\in\Gamma(2)$. They  fix the singularities of the staircase, and have zero drift.
\item \textbf{\textit{Uniform Hyperbolicity:}} For every $n\geq 0$, $\R^2=E^u_n\oplus E^s_n$ where\index{Uniform hyperbolicity}
    \begin{enumerate}[(a)]
    \item $E^s_0=\Span\{{\beta\choose 1}\}$, $E^u_0=\Span\{{1\choose -\beta}\}$, and $A_n(E^t_n)=E^t_{n+1}$ $(t=u,s)$;
    \item $\|A_n|_{E^s_n}\|\leq \lambda^{-1}$ and  $\|A_n|_{E^u_{n}}\|\geq \lambda$;
    \item $|\sin\measuredangle(E^s_n,E^u_n)|>\epsilon(\alpha)$.
    \end{enumerate}
\item \textbf{\textit{Slopes:}} If $\sigma^s_n:=\sigma[A_n\cdots A_1 A_0 {\beta \choose 1}]$ and $\sigma^u_n:=\sigma[A_n\cdots  A_1 A_0{1\choose -\beta}]$, then:
\begin{enumerate}
\item $|\sigma^s_n|>1$ and $|\sigma^u_n|<1$;
\item   $\epsilon(\alpha)<|\sigma^t_n|, |\sigma_n^t-1|<\epsilon(\alpha)^{-1}$ $(t=s,u)$; 
\item   The digits in the regular and even CFE of $\sigma_{n}^{t}$, $t=u,s$, $n\geq 0$ have absolute values less than $K(\alpha)$.
\end{enumerate}
\item \textbf{\textit{Finitely Many Types:}} $\#\{\wt{\psi}_n:n=0,1,2,\ldots\}<\infty$;
\item \textbf{\textit{Bounded Geometry:}} The polygons $\Pi(n):=\Pi(E^s_n,E^u_n)$ constructed in Section \ref{ss.Pi(n)} satisfy all the bounded geometry estimates in Lemma \ref{Lemma-Bounded-Geometry-Gamma(M)}.\index{Bounded geometry}
\end{enumerate}
\end{prop}

\begin{proof}
Let $\eps(\alpha):=\eps(M)$ from Lemma~\ref{l.position}.
Since  $\alpha\in (0,1)$, $\beta\in (-1,1)$. 
Let $\beta:=[0;2m_1,2m_2,2m_3,\ldots]$ be the even CFE of $\beta$. 
By Lemma \ref{l.A-contracts}, there exist $C>1$ and $0<\eta<1$ which only depend on $M$, so that 
$$
\left\|A_{N_0}(\beta'){\beta'\choose 1}\right\|\leq C\eta^{N_0}\left\|{\beta'\choose 1}\right\| \text{for every  $N_0$ and  every $\beta'$ such that $M(\beta')\leq M$}.
$$

We claim that  there is a constant $C'>0$ which only depends on $M$ so that  for every $|\beta''|\leq M$ such that $|\sin\measuredangle ({\beta'\choose 1},{1\choose -\beta''})|>\eps(M)$,
$$
\left\|A_{N_0}(\beta'){ 1\choose -\beta''}\right\|\geq C'\eta^{-N_0}\left\|{ 1\choose -\beta''}\right\|.
$$
To see this note that the area of the parallelogram generated by the vectors
$A_{N_0}(\beta'){ \beta'\choose 1}$ and   $A_{N_0}(\beta'){ 1\choose -\beta''}$ is the same as the area of the parallelogram generated by ${ \beta'\choose 1}$ and   ${1\choose -\beta''}$, whence bounded below by $|\sin\eps(M)|$. Since $\left\|A_{N_0}(\beta'){ \beta'\choose 1}\right\|$ is exponentially small, $\left\|A_{N_0}(\beta'){1\choose -\beta''}\right\|$ is  exponentially large, with bounds depending only on $M$.

We now choose $N_0$ large enough so that $C\eta^{N_0}<\lambda^{-1}$ and $C'\eta^{-N_0}>\lambda$. Then for every $\beta',\beta''$ such that $M(\beta')\leq M$, $|\beta''|\leq M$ and 
$
|\sin({\beta'\choose 1},{1\choose -\beta''})|\geq \eps(M)$,
\begin{equation}\label{e.N_0-choice}
\left\|A_{N_0}(\beta'){\beta'\choose 1}\right\|< \lambda^{-1}\left\|{\beta'\choose 1}\right\|,\text{ and }
\left\|A_{N_0}(\beta'){1\choose -\beta''}\right\|> \lambda\left\|{1\choose -\beta''}\right\|.
\end{equation}
Define ${A}_0\!\!:=\!\!A_{N_0}(\beta)$, $\beta_0',\beta_0''\!\!:=\!\!\beta$,  $E^s_0\!\!:=\!\!\Span\{{\beta_0'\choose 1}\}$ and 
$E^u_0\!\!:=\!\!\Span\{{1\choose -\beta_0''}\}$. Clearly, $M(\beta_0'), M(\beta_0'')\leq M$ and $E^u_0\perp E^s_0$, and therefore by \eqref{e.N_0-choice},
$$
\|A_0|_{E^s_0}\|<\lambda^{-1}\text{ and }\|A_0|_{E^u_0}\|>\lambda.
$$

Let 
$E^s_1:=A_0(E^s_0)$ and $E^u_1:=A_0(E^u_0)$. 
By Lemma \ref{l.position} and since $M(\beta_0'), M(\beta_0'')\leq M$, 
\begin{equation}\label{e.zero-ineq}
|\sin\measuredangle(E^s_1,E^u_1)|>\eps(M).
\end{equation}
Write  $
E^s_1=\Span\{{\beta_1'\choose 1}\}$ and $E^u_1=\Span\{{1\choose -\beta_1''}\}
$, with 
$
\beta_1':=\hsigma\left[A_0{\beta_0'\choose 1}\right]$, and  $-(\beta_1'')^{-1}:=\hsigma\left[A_0{1\choose -\beta_0'}\right]
$.
Co-slope calculations  as  in  \eqref{e.sigma-s} and \eqref{e.sigma-u} show that
\begin{align*}
&\beta_1'=\hsigma\bigl[\bigl(\!{\tiny{\!\begin{array}{c}\beta_1'\\ 1\end{array}}}\!\!\bigr)\bigr]=\hsigma \bigl[A_0\bigl(\!{\tiny{\!\begin{array}{c}\beta_0'\\ 1\end{array}}}\!\!\bigr)\bigr]=A_0\!\cdot\! \beta_0'=E^{2N_0}(\beta_0')=E^{2N_0}(\beta)=[0;2m_{2N_0+1},2m_{2N_0+2},\ldots]\\
&\beta_1'' =-\frac{1}{\hsigma{\tiny{1\choose -\beta_1''} }}=-\frac{1}{\hsigma [A_0{\tiny{1\choose -\beta} }]} =[0;2m_{2N_0},2m_{2N_0-1},\ldots,2m_2, 2m_1, 2m_1,2m_2,2m_3,\ldots ].
\end{align*}
It follows that  $\beta_1', \beta_1''$ are irrationals in $(-1,1)$ such that $M(\beta_1'), M(\beta_1'')\leq M$. 
Let $A_1:=A_{N_0}(\beta_1')=A_{N_0}(E^{2N_0}(\beta))$. By \eqref{e.N_0-choice} and \eqref{e.zero-ineq},
$$
\|A_1|_{E^s_1}\|<\lambda^{-1}\text{ and }\|A_1|_{E^u_1}\|>\lambda.
$$

We now let $E^s_2\!\!:=\!\!A_1(E^s_1)$, $E^u_2\!\!:=\!\!A_1(E^u_1)$. By Lemma \ref{l.position} and since $M(\beta_1'),M(\beta_1'')\!\!\leq\!\! M$, 
$|\sin\measuredangle(E^s_2,E^u_2)|>\eps(M)$. Equations \eqref{e.sigma-s} and \eqref{e.sigma-u} lead as before to
$
E^s_2:=\Span\{{\beta_2'\choose 1}\}$ and $E^u_2:=\Span\{{1\choose -\beta''_2}\}
$, with 
\begin{align*}
&\beta_2'=A_1\cdot \beta_1'=E^{2N_0}(\beta_1')=E^{4N_0}(\beta)=[0;2m_{4N_0+1},2m_{4N_0+2},\ldots]\\
&\beta_2''=-\frac{1}{\hsigma(A_1{\tiny{1\choose -\beta''_1})}} =[0;2m_{4N_0},2m_{4N_0-1},\ldots,2m_2, 2m_1, 2m_1,2m_2,2m_3,\ldots ].
\end{align*}
Again $\beta_2',\beta_2''$ are irrationals in $(-1,1)$ such that $M(\beta_2'),M(\beta_2'')\leq 
M$.  Let $$A_2:=A_{N_0}(\beta_2')=A_{N_0}(E^{4N_0}(\beta)).$$
Then as before,
$\|A_2|_{E^s_2}\|<\lambda^{-1}\text{ and }\|A_2|_{E^u_2}\|>\lambda$,

Continuing in this way,  we obtain Part (2) of the proposition with the matrices
$\DS A_k:=A_{N_0}(E^{2kN_0}(\beta))$, and the spaces  $E^s_k=\Span\{{\beta_k'\choose 1}\}$, $E^u_k=\Span\{{1\choose -\beta_k''}\}$, where 
\begin{align*}
&\beta_k'=[0;2m_{2kN_0+1},2m_{2kN_0+2},\ldots]; \\ 
&\beta_k'' =[0;2m_{2kN_0},2m_{2kN_0-1},\ldots,2m_2, 2m_1, 2m_1,2m_2,2m_3,\ldots ].
\end{align*}
The slopes (not co-slopes) of 
$E^s_k,E^u_k$ are 
\begin{equation}\label{e.E-k} 
\begin{aligned}
& \sigma^s_k\!=\!\frac{1}{\beta_k'}\!=\![2m_{2kN_0+1};2m_{2kN_0+2},\ldots]\\
& \sigma^u_k\!=\!-\beta_k''\!=\! [0;-2m_{2kN_0},-2m_{2kN_0-1},\ldots,\!-2m_2, -2m_1, -2m_1,-2m_2,-2m_3,\ldots ].
\end{aligned}
\end{equation}
Clearly,  $|\sigma^s_k|>1$, $|\sigma^u_k|<1$, and 
$M(\sigma^u_k),M(\sigma^s_k)\leq M$ (see \eqref{e.M-beta}). By Lemma \ref{l.BT-bounded-pm-two}, the partial  quotients of the regular CFE of $\sigma^s_k,\sigma^u_k$ are bounded by 
$M+3$. Necessarily, the distance of 
 $\sigma^s_k, \sigma^u_k$ from  $-1, 0, 1$ is bounded below by a constant which only depends on $M$, whence only on $\alpha$.  This proves  Part (3). 

Let $\wt{\psi}_n$ denote the unique homogeneous automorphism of the staircase with zero drift and derivative matrix $A_n$, see Proposition \ref{Prop-Automorphisms}. This automorphism fixes the singularities of $\St$. 
By construction, $A_n=A_{N_0}([0;2m_{2nN_0+1},2m_{2nN_0+2},2m_{2nN_0+3},\ldots])$. Looking at \eqref{e.A-k} and recalling that $|2m_i|\leq M$ for all $i$, we see that $\{A_0,A_1,\ldots\}$ is finite. Therefore $\{\wt{\psi}_0,\wt{\psi}_1,\ldots\}$ is finite.

We proved parts (1)--(4). Part (5) follows from Lemma \ref{Lemma-Bounded-Geometry-Gamma(M)}, with $\gamma=\beta$. 
\end{proof}

\begin{remark}\label{r.mat-prod-form} \normalfont 
The proof gives  $A_j=A_{N_0}(E^{2jN_0}(\beta))=A_{N_0}([0;2m_{2jN_0+1},2m_{2jN_0+2},\ldots])$. 
By \eqref{e.A-k-prime}, 
$
A_j=\left(
\begin{array}{cc}
-2m_{2(j+1)N_0} & 1 \\
1 & 0
\end{array}
\right)\cdots
\left(
\begin{array}{cc}
-2m_{2jN_0+2} & 1 \\
1 & 0
\end{array}
\right)
 \left(
\begin{array}{cc}
-2m_{2jN_0+1} & 1 \\
1 & 0
\end{array}
\right).
$
Taking the product over $j$ and using \eqref{e.A-k-prime} one more time, we obtain the identities
\begin{equation}\label{e.mat-prod-form}
\tpsi_{n-1}\circ\cdots\circ\tpsi_0=\wt{\Psi}_{nN_0}\ ; \ A_{n-1}\cdots A_2 A_0=	A_{nN_0}(\beta), 
\text{  see \eqref{e.A-k}}.
\end{equation}
Notice that by construction,  we can take the same $N_0$ for all $\beta$ such that $M(\beta)\leq M$. 
\end{remark}

\begin{lem}\label{l.strict-convexity-P}
For all $n$ large enough, $\mathcal P_n(\xi)$ given in \eqref{TempPres} is strictly convex. \index{Geometric pressure!strict convexity}
\end{lem}
\begin{proof}
Recall that $\mathcal P_n(\xi)=\log\left(\int_{\Rect_0}e^{\xi z_n(\tomega)}d\mu_0\right)$, with  $z_n(\tomega)$ as in \eqref{e.z-n}. A direct calculation shows that 
$\DS
\mathcal P''(\xi)=\frac{\int_{\Rect_0} e^{\xi z_n}z_n^2 d\mu_0}{\int_{\Rect_0}e^{\xi z_n}d\mu_0}-\left(\frac{\int_{\Rect_0}e^{\xi z_n}z_nd\mu_0}{\int_{\Rect_0}e^{\xi z_n}d\mu_0}\right)^2.
$
The expression on the right is the variance of the random variable $z_n(\tomega)$, when $\tomega$ is sampled randomly uniformly from the probability measure on $\Rect_0$ with density proportional to $e^{\xi z_n(\tomega)}$. 

Every non-constant random variable has strictly positive variance. Therefore,  to show that $\mathcal P''(\xi)>0$, it is sufficient to show that for all $n$ large enough,  $z_n$ is not equal a.e. to a constant on $\Rect_0$.

Recall that $z_n(\tomega)=\zc(\wt{\Psi}_n(\tomega))$ a.e. on $\Rect_0$. 
{ Fix some linear segment $W\subset\Rect_0$} such that $W$ is parallel to $E^u_0$. Then $|\wt{\Psi}_n(W)|\xrightarrow[n\to\infty]{}\infty$, so for all $n$ sufficiently large, 
$\zc$ is not constant on $\wt{\Psi}_n(W)$. Replacing $W$ by an $(E^u_0,E^s_0)$-parallelogram with $E^u_0$-side $W$, we see that $\zc$ is not almost surely constant on $\wt{\Psi}_n({\Rect_0})$. 
\end{proof}

\subsection{Twists}\label{SSDriftZero}

The renormalizing automorphisms $\wt{\psi}_n$ fix the singularities of the staircase. Thus, if $L_0$ is a ray emanating from a singularity $o$ in direction ${\beta\choose 1}$, then  $L_k:=(\wt{\psi}_{k-1}\circ\cdots\circ \wt{\psi}_0)(L_0)$ is a ray emanating from $o$, in direction $A_{k-1}\cdots A_1 A_0{\beta\choose 1}$.\index{Twist}\index{Automorphisms of the infinite staircase!twist}

However, since each  singularity is the meeting point of infinitely many  horizontal rectangles (Figure \ref{Figure-Staircase}(a)), $L_k$ may start at a different rectangle than $L_0$. Heuristically, what could happen is that each application of $\wt{\psi}_j$ ``twists" the ray into a different rectangle which meets the singularity $o$. \index{Infinite staircase!singularities} \index{Horizontal rectangles}

Our aim in this section is to calculate the $\Z$-coordinate of the horizontal rectangle containing the beginning of $L_k$, and tie its behavior as $k\to\infty$ to the even CFE of $\beta$. Recall the number $N_0$ from \eqref{e.mat-prod-form}.

\begin{thm}
\label{ThDriftOfOrigin}
Let $\beta\in (-1,1)$ be an irrational with even CFE  $\beta:=[0;2m_1,2m_2,\ldots]$, and let $L_0$ denote the ray in direction $\beta\choose 1$,  which begins at 
the singularity $\tomega_0$, in the middle of the top side of rectangle $\mathsf R_0$. Let $\frak{z}_n$ denote the $\Z$-coordinate of the horizontal rectangle which contains the beginning of  
$(\wt{\psi}_{n-1}\circ\cdots\circ\wt{\psi}_0)(L_0)$. Then
$$
\frak{z}_n=\frac{1}{2}\sum_{k=1}^{nN_0} \bigl[\sgn(m_{2k})-\sgn(m_{2k-1})\bigr].
$$
\end{thm}
\noindent
The proof is a direct calculation, recorded for completeness in Appendix \ref{Appendix-Twist}.

\begin{cor}
If the frequency of odd $j$ with a sign change $\sgn(m_{j})\neq \sgn(m_{j+1})$  is zero, then $\frak{z}_n/n\to 0$.
\end{cor}

\noindent
Recall that in our main application $\beta\!\!=\!\!2\alpha\!-\!1$ where $\alpha$ is the angle of rotation on 
$\T\!\!=\!\!\R/\Z$.

\begin{cor}[\cite{Avila-Dolgopyat-Duryev-Sarig}] 
\label{CrSqrt23}

\begin{enumerate}[(1)]
\item
$|\frak{z}_n|\leq nN_0$;
\item If $\beta=2\sqrt{3}-3$ (the direction corresponding to $\alpha=\sqrt{3}\mod 1$), then $\frak{z}_n=0$;
\item  If $\beta=2\sqrt{2}-3$ (the direction corresponding to $\alpha=\sqrt{2}\mod 1$), then $\frak{z}_n=
nN_0$;
\item  
If $\beta=3-2\sqrt{2}$ (the direction corresponding to $\alpha=-\sqrt{2}\mod 1$), then $\frak{z}_n=
-nN_0$.
\end{enumerate}
\end{cor}
\begin{proof} Part (1) is clear. Parts (2)--(4)  follow from the identities
\begin{equation}
\label{Sqrt23ECF}
2\sqrt{3}-3=[0;2,6,2,6,\ldots],\quad
2\sqrt{2}-3=[0;-6,6,-6,6,\ldots], 
\end{equation}
and $-(2\sqrt{2}-3)=[0;6,-6,6,-6,\ldots]$.  (We guessed \eqref{Sqrt23ECF} using the even Gauss map as in the proof of Lemma \ref{l.f-expansion}. The guesses can be confirmed by  solving $t=[0;-6,6,t]$ or $t=[0;2,6,t]$ subject to the constraint $t\in (-1,1)$.)
\end{proof}

\part{Symbolic Dynamics}\label{p.Symbolic-Dynamics}
 In this part of the paper,  we  develop the tools needed for the implementation of  the first step of the heuristic  program outlined  in \S\ref{s.difficulties}.  

\section{Markov Partitions}\label{Section-Polygons} 
Let $\wt{\psi}_n:\St\to\St$ denote the renormalizing automorphisms we constructed in Proposition \ref{Prop-Contracting-Automorphisms}. Each $\wt{\psi}_n:\St\to \St$  projects to a toral automorphism  $\psi_n:\St_0\to\St_0$, given by \eqref{Eq-Projected-Auto}.  

We will construct a sequence of ``Markov partitions" $\mathfrak P_n$ of $\St_0$ with the following property: Given $\omega\in\St_0$, let
\begin{equation}\label{e.X-new}
X_i^{(n)}(\omega):=\text{ the element of $\mathfrak P_i$ which contains }\begin{cases}(\psi_{i}^{-1}\circ\cdots\circ\psi_{n-1}^{-1})(\omega) & i<n\\
 \omega & i=n\\
 (\psi_{i-1}\circ\cdots\circ \psi_n)(\omega) & i>n.
 \end{cases}
\end{equation}
Then the joint distribution of 
$
(X_i^{(n)}(\omega))_{i\geq 0}
$ when  $\omega$ is sampled  randomly uniformly in $\St_0$ is an {\em inhomogeneous  Markov chain}. 
We will also construct ``jump functions" $f_k:\mathfrak P_k\to\R$ so that the random walk
$$
S_n=f_1(X_1^{(n)}(\omega))+\cdots+f_n(X^{(n)}_n(\omega)),\  \omega\text{ sampled randomly uniformly from }\St_0
$$
is a good approximation (in a sense which we will clarify below) of the stochastic process 
$$
\frak{z}\bigl((\wt{\psi}_{0}^{-1}\circ\cdots\wt{\psi}_{n-1}^{-1})(\tomega)\bigr),\  \tomega\text{ sampled randomly uniformly from }\Rect_0.
$$
Later, in Part IV, we will show that the approximation remains useful if we sample $\tomega$ randomly uniformly on $A_T^\ast:=(\wt{\psi}_{n-1}\circ\cdots\circ\wt{\psi}_0)(A_T)$, where $A_T:=\{\phi^t(\tomega_0):0<t<T\}$, and $n=n(T)$ is chosen appropriately. Equation  \eqref{e.Int=Prob} will then allow us to analyze the Birkhoff integrals of $\phi$, by analyzing $S_n$. 

\begin{remark}\normalfont
If $\alpha$ is a quadratic irrational,\index{Quadratic irrationals!renormalization by a single map} then the sequence of matrices in \eqref{e.A-k} is eventually periodic. One can then choose all the  ${\psi}_n$ to be  equal \cite{Avila-Dolgopyat-Duryev-Sarig}, and take all the partitions $\mathfrak P_n$ to be equal to the Adler \& Weiss Markov partition from   \cite{Adler-Weiss-MP}.  
But in the general bounded type case, the sequence $\{\psi_n\}$  cannot be chosen to be  constant, and we  use a modification of  the symbolic coding due to  Arnoux \& Fisher \cite{Arnoux-Fisher}. The difference between our coding and the Arnoux-Fisher coding is in the choice of the fundamental domain for $\R^2/\sqrt{2}\Z^2$. We use the fundamental domains in  Section \ref{ss.Pi(n)}. This leads to a simpler symbolic description of the jump functions $f_k$.\index{Markov partition!Adler-Weiss and Arnoux-Fisher constructions}
\end{remark}

\begin{remark}\normalfont
Sinai and Bowen  provided symbolic coding for a single {\em non-linear} uniformly hyperbolic diffeomorphism in any dimension \cite{Sinai-MP1, Sinai-MP2, Bowen-MP}, which codes the measure of maximal entropy by a Markov chain. This was generalized to sequences of different uniformly hyperbolic diffeomorphisms in 
\cite{GK97, MuentesRibeiro}. 
But  the mixing properties of the 
symbolic systems obtained this way require further analysis to allow the extension of our methods to this context; therefore  we could not use the symbolic codings in these
 references.
\end{remark}

Throughout Part \ref{p.Symbolic-Dynamics}, $\alpha\in (0,1)$ is an irrational number of bounded type, and $\beta:=2\alpha-1$. Let  $M=M(\beta)+1$, where $M(\beta)$ is given in \eqref{e.M-beta}.  Let $\wt{\psi}_n:\St\to\St$ and  $\psi_n:\St_0\to\St_0$ be as above.

Let $E^s_n, E^u_n$ denote the stable and unstable one-dimensional linear spaces defined in Proposition \ref{Prop-Contracting-Automorphisms}, and let 
$
\Pi(n):=\Pi(E^s_n,E^u_n)
$,
see  Section \ref{ss.Pi(n)} and Figure \ref{Figure-Pi(n)}. By Lemma~\ref{Lemma-Fundamental-Domain}, $\Pi(n)$ is a fundamental domain for the action  $G\curvearrowright\R^2$ given in \eqref{e.Group-G}, and 
$$
F_n:=\Pi(n)/G\cong\St_0. 
$$
$G$ induces a side pairing of the sides of $F_n$, perhaps after further subdivision, which turns it into a  model for $\St_0$.\index{Fundamental domain|see {associated fundamental domain}}\index{Associated fundamental domains}

 By Lemma \ref{Lemma-Q}
$
F_n=Q_1(n)\cup Q_2(n),
$ where $Q_i(n)$ are parallelograms with disjoint interiors, and sides parallel to $E^s_n, E^u_n$.   
 $Q_1(n)$  contains the singularity $p_0$ congruent to  $P_0,P_2,P_4$ in its boundary, and $Q_2(n)$ does not.   We call the sides of $Q_i(n)$ parallel to $E^s_n$  {\em stable sides}, \index{Stable side}\index{Unstable side} and those parallel to $E^u_n$ {\em unstable sides}. Let  $\partial^s Q_i(n)$, $\partial^u Q_i(n)$ denote the union of stable, respectively unstable,  sides of $Q_i(n)$.\index{Associated fundamental domains!partition into parallelograms}

Recall that $E^s_n=\Span\{A_{n-1}\cdots A_1 A_0 {\beta\choose 1}\}$ and 
$E^u_n=\Span\{A_{n-1}\cdots A_1 A_0{ 1\choose -\beta}\}$, where $A_i$ are the derivatives of $\wt{\psi}_i$. 
 We define the {\em positive stable direction} on $\Pi(n)$ and $F_n$  to be the direction of the vector 
$A_{n-1}\cdots A_1 A_0{ \beta\choose 1}\in E^s_n$. Similarly, the {\em positive unstable direction} on $\Pi(n)$ is the  direction of the vector $A_{n-1}\cdots A_1 A_0 {1\choose-\beta}\in E^u_n$. \index{Positive stable and unstable directions}

\begin{lem}\label{Lemma-Bottom-Side}
Exactly one stable side of $Q_i(n)$ has the following property: If we move from {\em any} point in the relative interior of this side in the positive unstable direction, then we will initially enter the interior of $Q_i(n)$.
\end{lem}
\begin{proof}
It is clear that such a side exists, the only thing to check is that it is unique.
Otherwise, $\ov{Q_i(n)}$ projects to  a non-trivial closed invariant set for the linear flow in the positive unstable direction. But this flow has dense orbits in $\St_0$, because of the irrationality of the  positive unstable direction.
\end{proof}
\noindent
We will call the stable side given by Lemma \ref{Lemma-Bottom-Side} the {\em bottom stable side}. The other side we will call the {\em top stable side}. These sides may  overlap, but they are not equal. \index{Top (un)stable sides}\index{Bottom (un)stable side}

Let $\Gamma(M)$ and $T_0(M)$ denote the constants from  the end of Sections \ref{ss.Pi(n)} and  \ref{ss.density-const}.  Denote the vertices of $\Pi(n)$ by $\Pi_i(n)$ as in Figure \ref{Figure-Pi(n)}, and let $\Pi_6(n):=\Pi_0(n)$. 
\begin{lem}\label{Lemma-Choice-of-lambda}
For all $0\leq i\leq 5$ and $n\geq 0$,
\begin{align}
&\Gamma(M)^{-2}<{\max\{\|\Pi_i(n)-P_j\|:0\leq i,j\leq 5\}\over
\min\{\|\Pi_i(n+1)-P_j\|:0\leq i,j\leq 5\}}<\Gamma(M)^2 
\label{Eq-Choice-of-lambda}\\
&\begin{array}{l}
\text{Every line in $\mathsf{St}_0$ in direction $E^u_n$ with length $T_0(M)$}\\
\text{must enter $Q_1(n)$ and $Q_2(n)$}.
\end{array}
\label{Eq-Choice-of-lambda-also}
\end{align}
\end{lem}

\begin{proof}
Proposition \ref{Prop-Contracting-Automorphisms} guarantees that $\Pi(n)$ satisfies all the bounded geometry estimates in Lemma \ref{Lemma-Bounded-Geometry-Gamma(M)}. In particular, 
$\Gamma(M)^{-1}\!<\!\|\Pi_i(n)\!-\!P_j\|\!<\!\Gamma(M)$ for all $i,j,n$, whence   \eqref{Eq-Choice-of-lambda}. Next,  $Q_i(n)$ contain open balls of radius $1/\Gamma(M)$, so (\ref{Eq-Choice-of-lambda-also}) follows from Lemma \ref{LmDenseWs} and the choice of $T_0(M)$.
\end{proof}

\begin{lem}\label{Lemma-Invariant-Stable-Sides}
The bottom (resp. top) stable side of $Q_i(n)$ is mapped by $\psi_n$ to the 
relative interior of the  bottom (resp. top) stable side of $Q_1(n+1)$.
\end{lem}

\begin{proof}

 Recall that $\psi_n$ fixes the punctures of $\St_0$:
 $\psi_n(p_0)=p_0$ and $\psi_n(p_1)=p_1$ where $p_0$ is the  congruence  class of $P_0=(1,0)$,   and $p_1$ is the congruence class of  $P_1=(2,0)$.

Let $L_0(n)$ denote the line through $p_0$ in direction $E^s_n$. Given $r$, let $B^s_n(p_0,r)$ denote the symmetric linear segment in $L_0(n)$ with center $p_0$ and radius 
$$
r:=\max\{\|\Pi_i(n)-P_j\|:0\leq i,j\leq 5\}.
$$
By construction, the stable sides of $Q_i(n)$ are contained in $B^s_n(p_0,r)$, see Figure \ref{Figure-Pi(n)}.

Since $\psi_n(p_0)=p_0$, $d\psi_n=A_n$, $A_n(E^s_n)=E^s_{n+1}$, and $\|A_n|_{E^s_n}\|\leq \lambda^{-1}$,  $\psi_n(B^s_n(p_0,r))\subset B^s_{n+1}(p_0,\lambda^{-1}r)$. 
By Lemma \ref{Lemma-Choice-of-lambda} and the choice of $\lambda$ in \eqref{e.choice-of-lambda}, 
$$
\lambda^{-1} r<\frac{1}{\Gamma(M)^2}\max\{\|\Pi_i(n)-P_j\|:0\leq i,j\leq 5\}\leq \min\{\|\Pi_i(n+1)-P_j\|:0\leq i,j\leq 5\}. 
$$
Calling the right-hand-side $r'$, we deduce that 
$\psi_n(B^s_n(p_0,r))\subset B^s_n(p_0,r')$. 

We now observe that $B^s_n(p_0,r')$ is a subset of 
 the  top {\em and} of the bottom stable sides of $Q_1(n+1)$, because each of these segments contains one of $P_0,P_2,P_4$, which all project to $p_0$, and because the distance of these points to the edge is no less than $r'$. 
 
It follows that $\psi_n$ maps the top and bottom stable sides of $Q_i(n)$ into the intersection of the top and bottom sides of $Q_1(n+1)$.
\end{proof}

An {\em $E^t_n$-segment} $(t=u,s)$ is a linear segment parallel to $E^t_n$. An {\em $(E^u_n,E^s_n)$-parallelogram}\index{E$^u,E^s$-parallelogram@$(E^u,E^s)$-Parallelogram} is a parallelogram whose sides are $E^u_n$- and $E^s_n$-segments. \index{E$^u$-segment and E$^s$-segment @ $E^u$-segment and $E^s$-segment} \index{Parallelogram!$(E^u,E^s)$-parallelogram}
An  {\em open $u$-fibre in $F_n$} 
 is a non-empty open $E^u_n$-segment inside $Q_i(n)$,  whose boundary points lie on the stable
  sides of $Q_i(n)$.  Similarly, an {\em open $s$-fibre in $F_n$} is a non-empty open $E^s_n$-segment in $Q_i(n)$, with boundary points on the unstable sides of $Q_i(n)$.  \index{u-Fiber @ $u$-Fiber }\index{s-Fiber @ $s$-Fiber}

   A {\em $u$-fibre in $F_n$}
   is the closure of an open $u$-fibre
   in $F_n$.
An {\em $s$-fiber} in $F_n$ is the closure of an open $s$-fiber in $F_n$. 
 We denote the open  $u$ and $s$-fibres in $F_n$
through $\omega\in F_n$  by 
$W^{u}(\omega, F_{n})$ and $W^s(\omega,F_n)$. More generally, given a $(E^u_n,E^s_n)$-parallelogram $Q$, we let 
$$W^{u}(\omega,Q):=W^{u}(\omega,F_{n})\cap Q\quad \text{and} \quad W^{s}(\omega,Q):=W^{s}(\omega,F_{n})\cap Q.$$

\begin{lem}\label{Lemma-Partition}
For every $n$,
$\psi_n(Q_i(n))=Q_{i,1}^{n+1}\cup\cdots\cup Q_{i, N_i(n+1)}^{n+1}$, where:
\begin{enumerate}[(1)]
\item $Q_{i,j}^{n+1}$ are  $(E^u_{n+1},E^s_{n+1})$-parallelograms.
\item The stable sides of $Q_{i,j}^{n+1}$ are subsets of stable sides of $Q_1(n+1)$ or $Q_2(n+1)$.
\item For  $j<N_i(n+1)$, the top stable side of $Q_{i,j}^{n+1}$ is the bottom stable side of $Q_{i,{j+1}}^{n+1}$.
\item The interiors of $Q_{i,j}^{n+1}$ are pairwise disjoint.
\item $Q_{i,j}^{n+1}$ is a subset of $Q_1(n+1)$ or of $Q_2(n+1)$.
\item For each $\omega \in Q_i(n)$ there are $\varpi_j\in Q_{i,j}^{n+1}$ s.t. 
$\DS\psi_n(W^u(\omega,Q_i(n)))\!=\!\!\!\!\!\bigcup_{j=1}^{N_i(n+1)}\!\!\!\!\! W^u(\varpi_j,Q_{i,j}^{n+1}).$
\item If $Q_{i,j}^n\subset Q_k(n)$, then $\mathrm{int}(Q_{i,j}^n)\cap\psi_n^{-1}(\mathrm{int}(Q_{k,\ell}^{n+1}))\neq\emptyset$ for  $\ell=1,\ldots, N_k(n+1)$.
\item The bottom stable side of $Q_{1,1}^{n+1}$ and the top side of $Q_{1, N_1}^{n+1}$ contain $p_0$. All other stable sides of $Q_{i,j}^{n+1}$ do not.
\item   $N_1(n)$ and $N_2(n)$ are bounded by a constant which only depends only on $M$.
\item The lengths of the stable and unstable sides of $Q_{i,j}^{n+1}$  are uniformly bounded away from zero and infinity by  constants only depend on $M$, but not on $i,j$ and $n$.
\item Each $Q_{i,j}^{n}$ contains a ball of  radius $w>0$, which depends on $M$ but not on $i,j,n$.
\end{enumerate}
\end{lem}

\begin{proof}
Fix two  $u$-fibres $W,W'$ in $Q_i(n)$.
By Lemma \ref{Lemma-Invariant-Stable-Sides}, $\psi_n$ maps $u$-fibres in $F_n$ to the union of $u$-fibres in $F_{n+1}$. Write   
$$
\psi_n(W)=W_1\cup\cdots\cup W_N\ , \quad \psi_n(W')=W_1'\cup\cdots\cup W_{N'}',
$$
where
  $W_i$ and $W_i'$  are  $u$-fibers in $F_n$;   $W_1,\ldots, W_N$ (respectively $W_1',\ldots, W_{N'}'$) are ordered in the positive unstable direction in $ F_{n+1}$; and $\ov{W_1},\ldots, \ov{W_N}$ (respectively $\ov{W_1'},\ldots, \ov{W_{N'}'}$) meet only at their endpoints, if at all. For each $j$, there are $a_j,a_j'\in\{1,2\}$ such that $W_j\subset \ov{Q_{a_j}(n+1)}$ and $W_j'\subset \ov{Q_{a_j'}(n+1)}$. 
We claim that $N=N'$ and $a_j=a_j'$.

By Lemma \ref{Lemma-Invariant-Stable-Sides},  $W_1,W_1'$  begin at $Q_1(n+1)$,  and $W_N,W_{N'}$  end at $Q_1(n+1)$. Therefore $a_1=a_1'$ and $a_N=a_{N'}$. If $a_j=a_j'$ for all $1\leq j\leq \min\{N,N'\}$, then $|W_j|=|W_j'|$ for all $1\leq j\leq \min\{N,N'\}$. Since $|\psi_n(W)|=|\psi_n(W')|$, it must be the case that  $N'=N$ and we are done. 
Assume by contradiction that  there is a (minimal) $1<j<\min\{N,N'\}$ such that  $a_j\neq a_j'$. Then $W,W'$ are separated by a $u$-fiber $V$ such that 
$$
\psi_n(V)=V_1\cup\cdots\cup V_{N''},
$$
where $V_i$ is a $u$-fiber in $Q_{a_i}(n+1)$ between $W_i,W_i'$ for $i=1,\ldots,j-1$, and $V_j$  passes through the unstable boundary of $Q_1(n)$ or $Q_2(n)$. Again, since $\psi_n(V)$ begins and ends inside the open stable sides of $Q_1(n+1)$, $1<j<N''$. 
Looking at Figure \ref{Figure-Pi(n)}  in Section \ref{ss.Pi(n)}
, we see that $V_{j-1}\cup V_j\cup V_{j+1}$ must contain $p_1$ (the projection of $P_1,P_3,P_5$) in its interior.  It follows that  $\psi_n(V)\owns p_1$. But then $p_1\in\mathrm{int}(\psi_n(Q_i(n)))$, whence $p_1\in \mathrm{int}(Q_i(n))$  which is not true. So $a_j=a_j'$ for all $1\leq j\leq \min\{N,N'\}$, whence
$N=N'$ and the claim is proved.

We can now define $N_i(n+1):=N$ and take $Q_{i,j}^{n+1}$ to be the union of $W_j$ where $W$ ranges over all $u$-fibres in $Q_i(n)$.
Since $\psi_n$ is piecewise affine, $Q_{i,j}^{n+1}$ are solid parallelograms, and it is clear from the construction that they satisfy (1)--(7).

For (8), recall that $p_0$ is a fixed point of $\psi_n$ inside the top and bottom stable sides of $Q_1(n)$, and outside the stable sides of $Q_2(n)$. Since
the bottom side of $Q_{1,1}^{n+1}$ and the top side of $Q_{1, N_1}^{n+1}$ are equal to the $\psi_n$--image of the bottom and top stable sides of $Q_1(n)$, they must contain $p_0$.
The same argument shows that the bottom side of $Q_{2,1}^{n+1}$ and the top side of $Q_{2,N_2}^{n+1}$ do {\em not} contain $p_0$.
All other stable sides of $Q_{i,j}^{n+1}$  cannot contain $p_0$,  because they are mapped by $\psi_n^{-1}$ to the interior of $Q_i(n)$, which does not contain $p_0$.

To see part (9), we  recall that  the lengths of the $u$-sides of $Q_{i,j}^{n+1}$ are bounded away from zero and infinity by a constant which depends only on $M$ (Lemma \ref{Lemma-Bounded-Geometry-Gamma(M)}). By construction, the norms of $d\psi_n$ are uniformly bounded (because there are only finitely many such matrices). Now apply the sixth part of the lemma. 

To see part (10), we use the proof of part (9) to bound the $u$-sides, and  note that the $s$-sides of $Q_{i,j}^{n+1}$  are  bounded below, because are  equal to the image of one of the sides of  $\Pi(n)$ by $\psi_n$, and $\{\psi_k:k\geq 0\}$ is a finite set.

Part 11  follows from the uniform bounds on $\measuredangle(E^u_n,E^s_n)$ 
{in Lemma \ref{l.position}.}
\end{proof}

\begin{prop}\label{Prop-Markov}
Let $\mathfrak P_0:=\{Q_1(0), Q_2(0)\}$ and
$$\mathfrak P_{n}:=\{Q_{1,1}^{n},\ldots,Q_{1,N_1(n)}^{n};Q_{2,1}^{n},\ldots,Q_{2,N_2(n)}^{n}\}\ \ (n\geq 0),$$ then the following properties hold:
\begin{enumerate}[(1)]
\item \textit{\textbf{Partition:}} $\{\ov{P}: P\in\mathfrak P_n\}$ covers $\St_0$, and $\{\mathrm{int}(P):P\in\mathfrak P_n\}$ are pairwise disjoint.
\item \textit{\textbf{Product Structure:}} For every $P\in\mathfrak P_{n}$ and $\omega, \varpi \in P$, $W^u(\omega,P)\cap W^s(\varpi,P)$ consists of a single point, which also belongs to $P$.\index{Markov partition}\index{Product structure}
\item \textit{\textbf{Markov Property:}} If $P\in\mathfrak P_n$, $Q\in\mathfrak P_{n+1}$,  \index{Markov property!for a Markov partition}
$\omega \in \mathrm{int}(P)$,  and $\psi_n(\omega)\in \mathrm{int}(Q)$, then
 $\psi_n^{-1}(W^u(\psi_n(\omega),Q))\subset W^u(\omega,P)$, 
 $\psi_n(W^s(\omega ,P))\subset W^s(\psi_n(\omega),Q)$.
\item \textit{\textbf{Bounded Cardinality:}}  $|\mathfrak P_n|\leq C$, where $C$ depends only on $M(\beta)$ in \eqref{e.M-beta}.
\end{enumerate}
\end{prop}
\begin{proof}
Part (1) follows from {Lemma \ref{Lemma-Fundamental-Domain}}
 and the fact that $\psi_n$ is a homeomorphism. Part (2) is because the elements of $\mathfrak P_n$ are parallelograms with sides parallel to $E^s_n, E^u_n$. Part~(4) is a direct consequence of Lemma \ref{Lemma-Partition}(9). Part (3) can be proved as follows.

 Write $P=Q_{i,j}^n$ and $Q=Q_{k,\ell}^{n+1}$.
By Lemma \ref{Lemma-Partition}(5), $Q_{i,j}^n\subset Q_m(n)$ for $m=1$ or $2$.
Necessarily, if $\tomega\in\mathrm{int}(P)$ and $\psi_n(\tomega)\in\mathrm{int}(Q)$, then 
$$\psi_n(\omega)\in \psi_n(\mathrm{int}(Q_m(n)))\cap\mathrm{int}(Q_{k,\ell}^{n+1})\subset \psi_n(\mathrm{int}(Q_m(n)))\cap\psi_n(\mathrm{int}(Q_k(n))).$$ 
So  $\omega\in\mathrm{int}(Q_m(n))\cap \mathrm{int}(Q_k(n))$. It follows that $m=k$ (otherwise the intersection is empty). So $Q_{i,j}^n\subset Q_k(n)$.

By Lemma \ref{Lemma-Partition}(6),  
$\DS \psi_n(W^u(\omega, Q_{i,j}^n))
=\bigcup_{\ell'=1}^{N_k}
W^u(\varpi_{\ell'},Q_{k,\ell'}^{n+1})$ 
for some 
$\varpi_{\ell'} \in Q_{k,\ell'}^{n+1}$. Since $\psi_n(\omega)\in \mathrm{int}(Q_{k,\ell}^{n+1})$ and the interiors of $Q_{k,\ell'}^{n+1}$ are disjoint, it must be the case that $\psi_n(\omega)\in W^u({\varpi_\ell},Q_{k,\ell}^{n+1})$. It follows that  
$W^u(\varpi_{\ell},Q_{k,\ell}^{n+1})=W^u(\psi_n(\omega),Q_{k,\ell}^{n+1})$. So
$$\psi_n(W^u(\omega,Q_{i,j}^n))\supset W^u(\psi_n(\omega),Q_{k,\ell}^{n+1}).$$ 
This is the first inclusion in (3).

For the second inclusion,
we  use $Q_{i,j}^n\!\!\subset\!\! Q_k(n)$ again to see that  
$\psi_n(W^s(\omega ,Q_{i,j}^n))\!\!\subset\!\! \psi_n(W^s(\omega,Q_k(n)))$. This is a line segment through 
$\psi_n(\omega)$, in direction  $E^s_{n+1}$, and with endpoints on the unstable boundary of 
$\DS \psi_n(Q_k(n))\!\!=\!\!\bigcup_{{\ell=1}}^{N_k} Q_{k,\ell}^{n+1}$. Since $\psi_n(\omega)\in\mathrm{int}(Q_{k,\ell}^{n+1})$, the  only such line segment is  $W^s(\psi_n(\omega),Q_{k,\ell}^{n+1})$. 
Therefore $\psi_n(W^s(\omega,Q_{i,j}^n))\subset W^s(\psi_n(\omega),Q_{k,\ell}^{n+1})$.
\end{proof}

We call $\mathfrak P_{n}$ the  {\em Markov partition at time $n$}.

\section{Jump Functions}
\subsection{Modified ${\mathbf{\mathbb Z}}$-Coordinates}\label{Section-Modified-Z-coordinates}
Recall that $D:\St^\ast\to\St^\ast$ is the automorphism which translates horizontal rectangle $\Rect_k$ onto the horizontal rectangle $\Rect_{k+1}$ immediately above it. Every fundamental domain $\wt{F}$  for the action of  $\{D^n:n\in\Z\}$ induces a function, called the {\em $\Z$-coordinate associated to $\wt{F}$}, and  given by \index{Z-coordinate @ $\Z$-coordinate!modified} \index{Modified $\Z$-coordinate}
$$
\zc_{\wt{F}}:\St^\ast\to\Z\ , \ \zc_{\wt{F}}(\tomega)=k\text{ iff }\tomega\in {D^k(\wt{F})}. 
$$
Clearly, 
\begin{equation}\label{e.Z-coord-D}
\zc_{\wt{F}}\circ D=\zc_{\wt{F}}+1.
\end{equation}
We have already encountered the canonical $\Z$-coordinate, defined in \eqref{e.z-coord-canonical}, and corresponding to $\wt{F}=\Rect_0$. 
In this section we construct other $\Z$-coordinates $\wt{\zc}_n:\St^\ast\to\Z$, 
which are better adapted to the action of the renormalizing automorphisms. 

Recall that $\St_0$ is isometric to $\R^2/G$ where $G\curvearrowright\R^2$ is given in \eqref{e.Group-G} (Figure \ref{Figure-Staircase}(c)). 
In Section \ref{ss.Pi(n)} we constructed  a fundamental domain $\Pi(n):=\Pi(E^s_n,E^u_n)$ for $G\curvearrowright\R^2$, with sides parallel to $E^s_n,E^u_n$ 
(Figure \ref{Figure-Pi(n)}  in Section \ref{ss.Pi(n)}). Fix once and for all some point $\varpi\in\St_0$ such that 
$\DS \varpi \not\in \bigcup_{n\geq 0}\partial\Pi(n)/G$ 
(in particular, $\varpi$ is not a puncture), and let $\wt{\varpi}_0$ be the unique lift of $\varpi_0$ to 
$\St^\ast$ s.t. $\zc(\wt\varpi_0)=0$. Let $F_n:=\Pi(n)/G\subset\St_0$, and define
$$
\wt{F}_n\subset \St^\ast=\St\setminus\{\text{singularities}\}
$$
to be the unique  lift of $F_n$ to $\St$ such that 
 $\wt{F}_n$ is pathwise connected;
 $\wt{F}_n$ contains $\wt{\varpi}_0$; and
  $\wt{F}_n$ contains its bottom stable and unstable sides, but not its top stable and unstable sides.
(Here ``top" and ``bottom" are defined using the positive stable and unstable directions, see the discussion before Lemma \ref{Lemma-Bottom-Side}.)
These properties  guarantee that
$\DS \St^\ast=\biguplus_{k\in\Z}D^k(\wt{F}_n)$ (for a picture, see Figure \ref{Figure-Tesselation}
 in Section \ref{ss.Pi(n)}).  Let $\wt{\zc}_n:=\zc_{\wt{F}_n}$, namely:
\begin{equation}
\label{ModifiedZ}
\wt{\zc}_n:\St^\ast\to\Z, \quad \wt{\zc}_n(\omega)=k\Leftrightarrow \omega\in D^k(\wt{F}_n).
\end{equation}

\begin{lem}\label{Lemma-Modified-Approx-Canonical}
$\sup\limits_n\sup\limits_{\omega\in\St^\ast}|\wt{\zc}_n(\omega)-\zc(\omega)|<\infty$.
\end{lem}
\begin{proof}
By Lemma \ref{Lemma-Bounded-Geometry-Gamma(M)},  for some $N$,  
for every $n$, $\DS \wt{F}_n\subset \bigcup_{|n|\leq N}D^n(\Rect_0)$.
\end{proof}

\subsection{Jump Functions}\label{ss.Jump-Functions}
Let $\tomega_n\!\!:=\!\!\tomega$ and $\tomega_k\!\!:=\!\!(\wt{\psi}_k^{-1}\circ\cdots\circ\wt{\psi}_{n-1}^{-1})(\wt{\omega})$ 
$(0\!\!\leq\!\! k\leq n\!-\! 1)$. In \S\ref{ss.RW-heuristic} we related 
the Birkhoff integrals of the linear flow to the ``backward random walk" 
$
\zc(\tomega_n),\zc(\tomega_{n-1}),\ldots,\zc(\tomega_0).
$
By Lemma \ref{Lemma-Modified-Approx-Canonical}, this is equal to  the ``modified" backward random walk 
$
\wt{\zc}_n(\tomega_n),\wt{\zc}_{n-1}(\tomega_{n-1}),\ldots,\wt{\zc}_0(\tomega_0)
$ up to uniformly bounded error. 
The modified random walk has a convenient  symbolic dynamical description in terms of the Markov partitions $\mathfrak P_n$. To describe it, we need the following definition:

\begin{defi}[Jump Functions]\index{Jump functions}\index{Jump functions!symbolic}\index{Frobenius function!jump functions}
Suppose $\omega\in\St_0^\ast$, and choose some (any)  $\wt{\omega}\in\St^\ast$ which projects to $\omega\in\St_0^\ast$. 
Let $\tomega_n:=\tomega$ and $\tomega_k:=(\wt{\psi}_k^{-1}\circ\cdots\circ\wt{\psi}_{n-1}^{-1})(\wt{\omega})$. Set
\begin{equation}\label{Eq-Modified-Frobenius}
\wt{f}^{(n)}_0\equiv 0,\ \ 
\wt{f}_{k}^{(n)}(\omega):=
\wt{\zc}_{k}(\tomega_{k})-\wt{\zc}_{k-1}(\tomega_{k-1})\ (1\leq k\leq n) .
\end{equation}
\end{defi}
\noindent
This is a proper definition: If 
$\wt{\omega},\wt{\omega}'\in\St^\ast$ both project to $\omega\in\St_0^\ast$, then there is some $\ell\in\Z$ such that  
$\wt{\omega}'=D^\ell(\wt{\omega})$; and the right hand side of \eqref{Eq-Modified-Frobenius} is the same for 
$\wt{\omega}$ and $\wt{\omega}'$, because  $\wt{\psi}_i\circ D=D\circ\wt{\psi}_i$ and $\wt{\zc}_j\circ D=\wt{\zc}_j+1$. 

\begin{lem}\label{Lemma-f(X)}
There are functions $f_{k}:\mathfrak P_{k}\to\Z$  such that $f_0\equiv 0$ and  for all $n\geq k$,
$\wt{f}_{k}^{(n)}(\omega)=f_{k}(X^{(n)}_{k}(\omega))\text{ a.e. in $\St_0^\ast$, where $X^{(n)}_{k}$ are given in \eqref{e.X-new}}.
$ 
These functions are uniformly bounded by a constant which only depends on $M(\beta)$.
\end{lem}

\begin{proof}
Fix $\omega\in\St_0$ and let $\omega_n:=\omega$, $\omega_i:=({\psi}_i^{-1}\circ\cdots\circ{\psi}_{n-1}^{-1})({\omega})$ ($0\leq i< n$). 
Discarding a set of measure zero, we may assume that $\omega_i$ do not hit the boundaries of the elements of the Markov partitions.  Fix $0\leq k\leq n-1$, and let $\wt{\omega}_k'$ be the lift of $\omega_k$ to $\St^\ast$ such that 
$\wt{\zc}_k(\wt{\omega}_k')=0$, then
$\wt{\omega}:=(\wt{\psi}_{n-1}\circ \cdots\circ\wt{\psi}_k)(\wt{\omega}_k')$ projects to $\omega$, and 
\begin{equation}\label{Eq-xi-tilde}
\wt{f}^{(n)}_{k+1}(\omega)=\wt{\zc}_{k+1}(\wt{\psi}_k(\wt{\omega}_k')).
\end{equation}

Let $\wt{\gamma}$ denote the linear segment in the positive unstable direction $E^u_k$ from the bottom stable side of  $\wt{F}_k$ to $\wt{\omega}_k$. Then 
 $\wt{\Gamma}:=\wt{\psi}_k(\wt{\gamma})$  terminates at $\wt{\psi}_k(\wt{\omega})$, and starts on the bottom side of $D^{\tau_k}(\wt{F}_{k+1})$, where $\tau_k$, the ``twist"\index{Twist}, is a constant depending only on $k$ (through $\wt{\psi}_k$ and the direction of the  bottom side of $\wt{F}_k$, see  \S\ref{SSDriftZero}). So
\begin{equation}
\wt{\zc}_{k+1}(\wt{\psi}_k(\wt{\omega}_k'))=\wt{\zc}_{k+1}(\text{end}(\wt{\Gamma}))
-\wt{\zc}_{k+1}(\text{beginning}(\wt{\Gamma}))+\tau_k.
\end{equation}

Let $\Gamma$ denote the projection of $\wt{\Gamma}$ to $\St_0$. 
Recall the notation for $\mathfrak P_k$ elements introduced in Lemma \ref{Lemma-Partition}, and
suppose $X^{(n)}_{k+1}(\omega)=Q^{k+1}_{\ell,m}$.  Then $\Gamma$ is the successive concatenation of the  linear segments
\begin{align*}
&\Gamma_1 = \text{ a $u$-fibre in $Q_{\ell,1}^{k+1}$}\\
&\Gamma_2 = \text{ a $u$-fibre in $Q_{\ell,2}^{k+1}$}\\
& \hspace{1cm}\cdots\cdots\cdots\\
&\Gamma_m=\text{ the beginning of a $u$-fibre in $Q_{\ell,m}^{k+1}$ which ends at $\psi_k(\omega_k)$.}
\end{align*}
Let $\Gain_{k+1}({\Gamma}_r):=\wt{\zc}_{k+1}(\text{end}(\wt{\Gamma}_r))-\wt{\zc}_{k+1}(\text{beginning}(\wt{\Gamma}_r))$ for some (any) lift of $\Gamma_r$. 
By the bounded geometry of ${F}_k$ and $F_{k+1}$, $|\Gain_{k+1}|\leq G_0(M)$, where $G_0(M)$ depends only on $M$, and  
\begin{equation}\label{Eq-Sum-of Gains}
\wt{\zc}_{k+1}(\wt{\psi}_k(\wt{\omega}_k))=
\tau_k+\sum_{r=1}^{m-1}\Gain_{k+1}({\Gamma}_r).
\end{equation}
The sum ends at $k=m-1$, because $\Gain({\Gamma}_m)=0$.

 Recall that $\mathfrak P_{k+1}=\{Q^{k+1}_{1,1},\ldots,Q^{k+1}_{1,N_1}; Q^{k+1}_{2,1},\ldots,Q^{k+1}_{2,N_2}\}$, where $N_i=N_i(k+1)$. 

\medskip
\noindent
{\sc Claim:}
If  $\Gamma, \Gamma'$ are two $u$-fibres in $Q_{i,j}^{k+1}$, and 
  $j\neq N_1,N_2$, then 
$\Gain_{k+1}(\Gamma)=\Gain_{k+1}(\Gamma')$.

\medskip
\noindent
{\em Proof of the Claim:} Let $n=k+1$. Two $u$-fibres in $F_{n}$ have the same gain whenever  when we draw them in the plane as $u$-fibres in $\Pi(n)$ and  extend them a bit beyond their endpoints, they enter the same copy $g(\Pi(n))$, $g\in G$. This will happen whenever the endpoints of $\omega,\omega'$ are not separated by the beginning points of a $u$-fibre in $G(\Pi(n))$ which begins or passes through a singularity.

This is the case for all pairs of $u$-fibres in $Q_{i,j}^{k+1}$ whenever $1\leq j\leq N_i-1$.  The endpoints of such fibres cannot be separated by a singularity, because the top stable side of $Q_{i,j}^{k+1}$ does not contain a puncture when $j<N_i$ (if it did, $Q_i(k)$ would contain a puncture in its interior). And they  cannot be separated by the beginning point of a $u$-fibre which passes through a singularity, because $Q_{i,j+1}^{k+1}$ does not contain a puncture in its interior.  The claim follows.

\medskip
The claim guarantees that the quantity
$g_{k+1}(Q_{i,j}^{k+1}):=\Gain_{k+1}(\Gamma)$  for some (all) $u$-fibres $\Gamma$ in $Q_{i,j}^{k+1}$,   
 is properly defined whenever  $j\neq N_1,N_2$ .  
Let $f_{k+1}:\mathfrak P_{k+1}\to\Z$, 
 $$\DS f_{k+1}(Q_{i,j}^{k+1}):=\tau_k+\sum_{r=1}^{j-1}g_{k+1}(Q_{i,r}^{k+1}).$$
By \eqref{Eq-xi-tilde} and \eqref{Eq-Sum-of Gains}, 
$\wt{f}^{(n)}_{k+1}(\omega)=f_{k+1}(X^{(n)}_{k+1}(\omega))$ a.e.

It remains to show that  $f_k$ are uniformly bounded. By definition, 
$|f_{k+1}(Q_{i,j}^{k+1})|\leq |\tau_k|+ j G_0(M)\leq  |\tau_k|+N_i(k+1) G_0(M)$, where 
$N_i(k+1)$ is the maximal number of times $\psi_k(Q_i(k))$ wraps around $F_{k+1}$. By the bounded geometry estimates on $F_k$ and the bounds on $\|d\psi_k\|$, $N_i(k+1)$ can be bounded by a constant depending only on $M(\beta)$. Similarly, $\tau_k$ can be bounded uniformly in terms of $M(\beta)$, because $N_0$ (the constant in \eqref{e.mat-prod-form}) can be bounded uniformly in terms of $M(\beta)$, and therefore  $\{\psi_i:i\geq 0\}$ is contained in a finite set, which depends only on $M(\beta)$.
\end{proof}
\noindent
Notice that $f_k$ are independent of $n$. We call $f_k$ the {\em (symbolic) jump functions}.


\section{An Inhomogeneous Markov Chain}\label{s.Markov-Chain}
\subsection{Cylinders}\label{ss.cylinders} 
The\index{Cylinders} {\em cylinder} $\<a_{n-k},\ldots,\dot{a}_n,\ldots,a_{n+m}\>$ (with {\em entries} $a_i\in \mathfrak P_i$, {\em beginning} $n-k$, {\em center} $n$, and  {\em ending} $n+m$, $k,n,m\in\N\cup\{0\}$, $k\leq n$) is the set
\begin{align}
&\<a_{n-k},\ldots,\dot{a}_n,\ldots,a_{n+m}\>:=[X^{(n)}_{n-k}=a_{n-k},\ldots,X^{(n)}_{n+m}=a_{n+m}] \label{e.cylinder}\\
&:=\{\omega\in\St_0: X^{(n)}_i(\omega)=a_i\ (i=n-k,\ldots, n+m)\}\ \ \ \text{ (see \eqref{e.X-new})}\notag\\
&=\bigcap_{i=n-k}^{n-1} (\psi_{n-1}\circ\cdots\circ\psi_{i})(a_{i})\cap a_n\cap \bigcap_{j=n+1}^{n+m}(\psi_{n}^{-1}\circ\cdots\circ\psi_{j-1}^{-1})(a_j). \notag
\end{align}
 Clearly, 
\begin{equation}\label{e.cylinder-shift}
\psi_{n-k}^{-1}\circ\cdots\circ\psi_{n-1}^{-1}\<a_{n-k},\ldots,\dot{a}_n,\ldots,a_{n+m}\>=\<\dot{a}_{n-k},\ldots,a_{n+m}\>.
\end{equation}

The set of cylinders with beginning $n-k$, center $n$, ending $n+m$ and with non-empty interior is denoted by $\mathfrak C_{k,n,m}$. To describe this set, it is convenient to introduce the following notation for $a\in\mathfrak P_i, b\in \mathfrak P_{i+1}$: 
$a\xrightarrow[]{\psi_i}b\text{ iff }\psi_i(\mathrm{int}(a))\cap\mathrm{int}(b)\neq \emptyset.$
\begin{lem}\label{Lemma-Cylinders}
Suppose $a_i\in\mathfrak P_i$, then  
$\<a_{n-k},\ldots,\dot{a}_n,\ldots,a_{n+m}\>$ has non-empty interior iff $a_i\xrightarrow[]{\psi_i}a_{i+1}$ for all $i$. The interior is empty iff the Lebesgue measure is zero.
\end{lem}
\begin{proof}
By \eqref{e.cylinder-shift}, it is sufficient to prove the lemma for $C=\<\dot{a}_n,a_{n+1},\ldots,a_{n+m}\>$.
 
Clearly, if $\mathrm{int}(C)\neq \emptyset$, then $a_i\xrightarrow[]{\psi_i}a_{i+1}$ for all $i$. 
The converse follows from  the properties of $\mathfrak P_i$ listed in Proposition \ref{Prop-Markov}, as in \cite{Adler-Weiss-MP}. Here are the details.

 If $m=0$, $C=a_n$, and $\mathrm{int}(C)\neq \emptyset$.  Suppose by induction the lemma holds for $m-1$, and assume  $a_n\xrightarrow[]{\psi_n}a_{n+1}\xrightarrow[]{\psi_{n+1}}\cdots \xrightarrow[]{\psi_{n+m-1}}a_{n+m}$. By the induction hypothesis, there exists some point $\omega_{n+1}\in\mathrm{int}(C_{n+1})$, where  $C_{n+1}:=\<\dot{a}_{n+1},\ldots,a_{n+m}\>$. In addition, since $a_{n}\xrightarrow[]{\psi_{n}}a_{n+1}$, there is a point $\omega_n\in\mathrm{int}(a_{n})\cap\psi_n^{-1}(\mathrm{int}(a_{n+1}))$. By the Markov property, 
$$
\psi_{n}^{-1}[W^u(\psi_{n}(\omega_n),a_{n+1})]\subset W^u(\omega_n,a_{n}).
$$
By the product structure property, $\exists\omega^\ast\in W^u(\psi_{n}(\omega_n),a_{n+1})\cap W^s(\omega_{n+1},a_{n+1})$. 

We claim that $\omega\!\!:=\!\!\psi_{n}^{-1}(\omega^\ast)\!\!\in\!\! \mathrm{int}(C)$. First,   
$\omega\!\!\in\!\! \psi_{n}^{-1}[W^u(\psi_{n}(\omega_n),a_{n+1})]$, so by the Markov property, 
$\omega\in W^u(\omega_n,a_{n})\subset \mathrm{int}(a_n)$. 
Second, $\psi_{n}(\omega)\!\!\in\!\! W^s(\omega_{n+1}, a_{n+1})\!\!\subset\!\!\mathrm{int}(a_{n+1})$. Third,
by the Markov property  for every $n<j\leq n+m$,    
$$
(\psi_{j-1}\circ \cdots\circ \psi_{n})(\tomega)\!\!\in\!\! 
(\psi_{j-1}\circ \cdots\circ \psi_{n+1})[W^s(\omega_{n+1}, a_{n+1})]
\!\!\subset\!\! W^s((\psi_{j-1}\circ \cdots\circ \psi_{n+1})(\omega_{n+1}),a_{j})
\!\!\subset\!\! \mathrm{int}(a_{j}).
$$
Thus, $\DS \omega\in \bigcap_{n<j\leq n+m}(\psi_n^{-1}\circ\cdots\circ\psi_{j-1}^{-1})(a_j)=\mathrm{int}(C)$. 
This proves the first part of the lemma. 

The second part is because cylinders are  convex.
\end{proof}

Recall that $\mu_0=\frac{1}{2}\times$ the area measure on $\St_0$ (so $\mu_0(\St_0)=1$).
Every  $a\in\mathfrak P_n$ is a parallelogram with sides parallel to $E^s_n, E^u_n$. Let $\ell^s(a)$ and $\ell^u(a)$ denote the lengths of the stable and unstable sides of $a$, and set 
$$\lambda^s_n:=\|A_n|_{E^s_{n}}\|\text{ and }\lambda^u_n:=\|A_n|_{E^u_{n}}\|.$$
 By Proposition \ref{Prop-Contracting-Automorphisms}(2), $\lambda^s_n\leq \lambda^{-1}<1$, and $\lambda^u_n\geq \lambda>1$. Since
 $\{A_n:n\geq 0\}$ is finite, $\lambda^u_n,\lambda^s_n$ are  bounded away from zero and infinity.

%
%

\begin{lem}\label{Lemma-Cylinder-Area} Let  $C:=\<a_{n-k},\ldots,\dot{a}_n,\ldots,a_{n+m}\>\in\mathfrak C_{k,n,m}$. Then $C$ is an $(E^u_n,E^s_n)$-parallelogram in $\St_0$, and:  \index{Cylinders} \index{Stable side} \index{Unstable side}
\begin{enumerate}[(1)]
\item The stable side of $C$ has length $\ell^s(C)=\ell^s(a_{n-k})\cdot(\lambda^s_{n-k}\cdots \lambda^s_{n-1})$.
\item The unstable side of  $C$  has length $\ell^u(C)=\ell^u(a_{n+m})/(\lambda^u_{n}\cdots \lambda^u_{n+m-1})$.
\item $\mu_0(C)=\frac{1}{2}\ell^s(a_{n-k})\ell^u(a_{n+m})|\sin\measuredangle(E^s_{n},E^u_{n})|(\lambda^s_{n-k}\cdots\lambda^s_{n-1})/(\lambda^u_{n}\cdots \lambda^u_{n+m-1})$. 
\item The area of $C$ is bounded away from zero and infinity by constants which only depend on its length $k+m$.
\end{enumerate}
\end{lem}
\begin{proof}
First we consider  the case $k=0$. Set $\ell:=n+m$, and fix some  $C=\<\dot{a}_n,a_1,\ldots,a_\ell\>$ such that  $\mathrm{int}(C)\neq \emptyset$. By \eqref{e.cylinder}, $\DS C=a_n\cap \bigcap_{j=n+1}^\ell (\psi_n^{-1}\circ\cdots\circ\psi_{j-1}^{-1})(a_j).$
Take $n+1\leq j\leq \ell$, then  $a_j\in\mathfrak P_j$ is an $(E^u_j,E^s_j)$-parallelogram, and $\psi_n^{-1}\circ\cdots\circ\psi_{j-1}^{-1}$ is an affine map which  maps $E^s_j, E^u_j$ to $E^s_n, E^u_n$. Thus $(\psi_n^{-1}\circ\cdots\circ\psi_{j-1}^{-1})(a_j)$ is an $(E^s_n,E^u_n)$-parallelogram. 
Passing to  the intersection over $j$, we see that $C$  is an $(E^s_n,E^u_n)$-parallelogram.

Let us calculate its sides. 
Pick some ${\omega_n}\in\mathrm{int}(C)$, $n+1\leq j\leq \ell$, and let $\omega_j:=(\psi_{j-1}\circ\cdots\circ\psi_n)(\omega_n)$.  By the Markov property, 
$$
(\psi_{j-1}\circ\cdots\circ\psi_n)[W^s(\omega_n,a_n)]\subset W^s(\omega_j,a_j).
$$
Thus
$W^s(\omega_n,a_n)\subset a_n\cap\bigcap_{j=n+1}^{\ell}(\psi_n^{-1}\circ\cdots\circ\psi_{j-1}^{-1})(a_j)=C$. 
$W^s(\omega_n,a_n)$ is an open $s$-fiber in $a_n$. As such, it is an $E^s_n$-segment of length $\ell^s(a_n)$.
It follows that  the stable side of $C$ has length $\ell^s(C)\geq \ell^s(a_n)$. Since $C\subset a_n$,  $\ell^s(C)=\ell^s(a_n)$.

Let $W:=(\psi_n^{-1}\circ\cdots\circ\psi_{\ell-1}^{-1})(W^u(\omega_\ell,a_\ell))$. By the Markov property, for all $n+1\leq j\leq \ell$,   
$$
(\psi_{j-1}\circ\cdots\circ\psi_n)(W)=(\psi_j^{-1}\circ\cdots\circ\psi_{\ell-1}^{-1})(W^u(\omega_\ell,a_\ell))\subset W^u(\omega_j,a_j)\subset a_j.
$$
So $\DS W\subset  \bigcap_{j=n+1}^{\ell}(\psi_n^{-1}\circ\cdots\circ\psi_{j-1}^{-1})(a_j)=C$. 
$W^u(\omega_\ell,a_\ell)$ is an open $u$-fiber in $a_\ell$. As such, it is an $E^u_\ell$-segment of length $\ell^u(a_\ell)$. The derivative of $\psi_n^{-1}\circ\cdots\circ\psi_{\ell-1}^{-1}$ maps $E^u_\ell$ to $E^u_n$, and contracts vectors in $E^u_\ell$ by a factor of $1/(\lambda^u_n\cdots\lambda^u_{\ell-1})$. Therefore $W$ is an  $E^u_n$-segment with length $\ell^u(a_\ell)/(\lambda^u_n\cdots\lambda^u_{\ell-1})$. Since $W\subset C$, $\ell^u(C)\geq \ell^u(a_\ell)/(\lambda^u_n\cdots\lambda^u_{\ell-1})$. On the other hand, $(\psi_{\ell-1}\circ\cdots\circ\psi_n)(C)\subset a_\ell$, so  $\ell^u(C)\lambda^u_n\cdots\lambda^u_{\ell-1}\leq \ell^u(a_\ell)$. Necessarily, 
$
\ell^u(C)=\ell^u(a_\ell)/(\lambda^u_n\cdots\lambda^u_{\ell-1})
$


Now consider a general cylinder $C=\<a_{n-k},\ldots,\dot{a}_n,\ldots,a_{n+m}\>$. By \eqref{e.cylinder-shift}, 
$$
C=(\psi_{n-1}\circ\cdots\circ\psi_{n-k})(C'),\text{ where }C':=\<\dot{a}_{n-k},\ldots,a_{n+m}\>. 
$$ 
The differential of  $\psi_{n-1}\circ\cdots\circ\psi_{n-k}$ maps  $E^t_{n-k}$ to $E^t_n$ $(t=u,s)$. It contracts $E^u_{n-k}$ by  $\lambda^s_{n-k}\cdots\lambda^s_{n-1}$, and  expands $E^s_{n-k}$ by $\lambda^u_{n-k}\cdots\lambda^u_{n-1}$. Thus $C$ is a $(E^u_n,E^s_n)$-parallelogram, and by the discussion above, $\ell^s(C)=\lambda^s_{n-k}\cdots\lambda^s_{n-1}\ell^s(C')=\ell^s(a_{n-k})\lambda^s_{n-k}\cdots\lambda^s_{n-1}$, and 
$\ell^u(C)=\lambda^u_{n-k}\cdots\lambda^u_{n-1}\ell^u(C')=\ell^u(a_{n+m})/(\lambda^u_{n}\cdots\lambda^u_{n+m-1})$. 

Thus $\DS \mu_0(C)=\frac{1}{2}\mathrm{Area}(C)=\frac{1}{2}\ell^s(a_{n-k})\ell^u(a_{n+m})|\sin\measuredangle(E^s_{n},E^u_{n})|
\frac{\lambda^s_{n-k}\cdots\lambda^s_{n-1}}{\lambda^u_{n}\cdots\lambda^u_{n+m-1}}$.

Lemma \ref{Lemma-Partition}(10) gives uniform bounds on $\ell^s(a_{n-k}),\ell^u(a_{n+m})$. Proposition \ref{Prop-Contracting-Automorphisms}(2) gives  uniform bounds on $|\sin\measuredangle(E^s_{n},E^u_{n})|$,  and  Proposition \ref{Prop-Contracting-Automorphisms}(3) gives uniform bounds on $\lambda^s_i, \lambda^u_j$. It follows that $\mu_0(C)\asymp 1$ with bounds depending only on $k+m$.
\end{proof}

\subsection{An Inhomogeneous Markov Chain}\label{ss.IMC} \index{Markov chain!description of  $X_n$} If we choose $\omega$ randomly uniformly in $\St_0$, then the functions $X^{(n)}_i(\omega)$ given by \eqref{e.X-new} become random variables. The following result describes their joint distribution.   Let 
\begin{equation}\label{e.incidence-matrix}
t_{a_i,a_{i+1}}^{(i)}=\begin{cases}  1 & a_i\in\mathfrak P_i, a_{i+1}\in\mathfrak P_{i+1}, \mathrm{int}(a_i)\cap\psi_i^{-1}[\mathrm{int}(a_{i+1})]\neq \emptyset\\
0 & \text{otherwise}.
\end{cases}
\end{equation}
\begin{prop}\label{Prop-Markov-Kernel} 
If $\omega$ is chosen randomly uniformly in 
$\St_0$, then:
\begin{enumerate}[(1)]
\item  $X_i^{(n)}\in\mathfrak P_i$  for all $i\geq 0$.
\item  $\DS \Prob(X_0^{(n)}=a_0,\ldots,X_\ell^{(n)}=a_\ell)=\mu_0(a_0)\prod_{i=0}^{\ell-1} \frac{\ell^u(a_{i+1})t_{a_i,a_{i+1}}^{(i)}}{\lambda_{i}^u\ell^u(a_{i})}.$
\item $(X^{(n)}_i)_{i\geq 0}$ is equal in distribution to an inhomogeneous Markov chain $(X_i)_{i\geq 0}$ with
\begin{enumerate}
\item $i$-th state space $\mathfrak P_i$;
\item initial probability vector $\pi_a:=\Prob[X_0=a]=\mu_0(a)$ $(a\in \mathfrak P_0)$;
\item $i$-th transition kernel $\DS\pi_{i,i+1}(a_i,a_{i+1}):=\Prob[X_{i+1}=a_{i+1}|X_i=a_i]=\frac{\ell^u(a_{i+1})}{\lambda_i^u\ell^u(a_i)}t_{a_i,a_{i+1}}^{(i)}$.
\end{enumerate}
\item $\pi_a$  are positive and uniformly bounded away from zero by constants which only on $M(\beta)$  in \eqref{e.M-beta}, and $\pi_{i,i+1}(a_i,a_{i+1})$ are either zero or uniformly bounded away from zero by constants which depend only on $M(\beta)$.
\end{enumerate}
\end{prop}
\begin{proof} Part (1) holds by definition.

Suppose $\ell\geq n$. 
If we sample $\omega$ randomly uniformly in $\St_0$, then 
\begin{align*}
&\Prob(X_0^{(n)}=a_0,\ldots,X_\ell^{(n)}=a_\ell)=\mu_0[\<a_0,\ldots,\dot{a}_n,\ldots,a_\ell\>]\\
&=\frac{1}{2} \ell^s(a_0)\ell^u(a_\ell)|\sin\measuredangle(E^s_n,E^u_n)|\frac{\lambda^s_{0}\cdots\lambda^s_{n-1}}{\lambda^u_n\cdots \lambda^u_{\ell-1}}\prod_{i=0}^{\ell-1}t_{a_i,a_{i+1}}^{(i)}, \text{ see  Lemmas \ref{Lemma-Cylinders} and \ref{Lemma-Cylinder-Area}}.
\end{align*}
The automorphism  $\psi_{0}^{-1}\circ\cdots\circ\psi_{n-1}^{-1}$ maps the parallelogram with unit sides parallel to $E^s_n,E^u_n$, to a parallelogram with (perpendicular) sides parallel to  $E^s_0,E^u_0$, and with lengths $1/(\lambda^s_0\cdots\lambda^s_{n-1})$, $1/(\lambda^u_0\cdots\lambda^u_{n-1})$. Since the Jacobian is one, the area remains unchanged, whence 
$
|\sin\measuredangle(E^s_n,E^u_n)|=1/(\lambda^s_0\cdots\lambda^s_{n-1}\cdot\lambda^u_0\cdots\lambda^u_{n-1})
$.  Therefore, 
\begin{align*}
&\Prob(X_0^{(n)}=a_0,\ldots,X_\ell^{(n)}=a_\ell)=\frac{1}{2} \frac{\ell^s(a_0)\ell^u(a_\ell)}{\lambda^u_0\cdots \lambda^u_{\ell-1}}\prod_{i=0}^{\ell-1}t_{a_i,a_{i+1}}^{(i)}.
\end{align*}
Since $E^s_0\perp E^u_0$,  $\mu_0(a_0)=\frac{1}{2}\text{Area}(a_0)=\frac{1}{2}\ell^s(a_0)\ell^u(a_0)$, and 
\begin{align*}
&\Prob(X_0^{(n)}=a_0,\ldots,X_\ell^{(n)}=a_\ell)=
\frac{ \mu_0(a_0)\ell^u(a_\ell)}{\lambda^u_0\cdots \lambda^u_{\ell-1}\ell^u(a_0)}\prod_{i=0}^{\ell-1}t_{a_i,a_{i+1}}^{(i)}=\mu(a_0)\prod_{i=0}^{\ell-1} \frac{\ell^u(a_{i+1})}{\lambda^u_i\ell^u(a_i)}t_{a_i,a_{i+1}}^{(i)}.
\end{align*}

This proves part (2) when $\ell\geq n$. To extend this to the case $\ell<n$, we first observe that by Part 6 of Lemma \ref{Lemma-Partition}, 
$\DS \lambda^u_i \ell^u(a_i)=\sum_{a_{i+1}\in\mathfrak P_{i+1}}t_{a_i,a_{i+1}}^{(i)}\ell^u(a_{i+1})$, whence 
$$
\sum_{a_{i+1}\in\mathfrak P_{i+1}}\frac{\ell^u(a_{i+1})}{\lambda^u_i\ell^u(a_i)}t_{a_i,a_{i+1}}^{(i)}=1.
$$ 
Summing the identity 
$\DS \Prob(X_0^{(n)}\!\!=\!\!a_0,\!\ldots\!,X_n^{(n)}\!\!=\!\!a_n)\!\!=\!\! \mu(a_0)\!\prod_{i=0}^{n-1} \frac{\ell^u(a_{i+1})}{\lambda^u_i\ell^u(a_i)}t_{a_i,a_{i+1}}^{(i)}
$ over  $a_n\!\!\in\!\!\mathfrak P_n$, we  obtain  
$
\DS \Prob(X_0^{(n)}\!\!=\!\!a_0,\ldots,X_{n-1}^{(n)}\!\!=\!\!a_{n-1})\!\!=\!\!\mu(a_0)\prod_{i=0}^{n-2} \frac{\ell^u(a_{i+1})}{\lambda^u_i\ell^u(a_i)}t_{a_i,a_{i+1}}^{(i)}.
$
Summing  on $a_{n-1}\in\mathfrak P_{n-1}$ gives  
$\DS
\Prob(X_0^{(n)}=a_0,\ldots,X_{n-2}^{(n)}=a_{n-2})=\mu(a_0)\prod_{i=0}^{n-3} \frac{\ell^u(a_{i+1})}{\lambda^u_i\ell^u(a_i)}t_{a_i,a_{i+1}}^{(i)}.
$ Continuing in this way, we obtain the second part for $\ell<n$.

Part (3) is just a restatement of Part (2) in the terminology of Markov chains. Part (4)
 follows from the bounded geometry estimates in Lemma \ref{Lemma-Partition}. 
\end{proof}

\begin{cor}
The joint distribution of $(X^{(n)}_i)_{i\geq 0}$ does not depend on $n$.
\end{cor}
\noindent
There is no surprise here: By \eqref{e.X-new}, 
\begin{equation}\label{e.X_i-tranformation-law}
X^{(n)}_i(\psi_{n-1}(\omega))=X^{(n-1)}_i(\omega),
\end{equation} 
Since $\psi_{n-1}$ is area preserving, if we sample $\omega$ randomly uniformly in $\St_0$, then $\psi_{n-1}(\omega)$ is distributed randomly uniformly in $\St_0$ too. It follows that $(X^{(n)}_i)_{i\geq 0}$ and $(X^{(n-1)}_i)_{i\geq 0}$ have the same joint distribution.

\medskip
Henceforth we let $\mathsf X:=\{X_i\}_{i\geq 0}$ denote the inhomogeneous Markov chain given by Proposition \ref{Prop-Markov-Kernel}.

\subsection{Uniform Ellipticity and Exponential Mixing}\label{ss.uniform-ellipticity}
Suppose $(Y_n)_{n\geq 0}$ is an inhomogeneous Markov chain  with finite state spaces $\mathfrak S_n$, initial distribution $\mu_0(x)\!\!:=\!\!\Prob[Y_0\!\!=\!\!x]$, 
and transition kernel $\pi_{n,n+1}(x,y)\!\!:=\!\!\Prob[Y_{n+1}=y|Y_n=x]$. 

\begin{defi}\label{d.uniform-ellipticity}\index{Uniform ellipticity}
$(Y_n)_{n\geq 0}$ is   {\em uniformly elliptic}, if there is a positive constant $\epsilon_0$, called an {\em ellipticity constant}, and a representation $\pi_{n,n+1}(y,y')=p_n(y,y')\mu_{n+1}(y')$, where 
\begin{enumerate}[(1)]
\item $\mu_n$ are probability measures on $\mathfrak S_n$;
\item $0\leq p_n(\cdot,\cdot)\leq 1/\epsilon_0$;
\item $\DS \sum_{y\in\mathfrak S_{n+1}} p_n(x,y)p_{n+1}(y,z)>\epsilon_0$ for all $x\in\mathfrak S_n, z\in\mathfrak S_{n+2}$. 
\end{enumerate} 
\end{defi}

\begin{remark}\normalfont
The condition in Definition \ref{d.uniform-ellipticity} is sometimes called the {\em two-step ellipticity 
condition}, because (3) says that one can jump from any state at time $n$ to any state at time $n+2$ with positive probability, in two steps. For a discussion of  other $\gamma$-step ellipticity conditions, see  \S2.4 in \cite{Dolgopyat-Sarig-Book}.
\end{remark}

\begin{prop}\label{Prop-Doeblin}
The Markov chain $(X_n)_{n\geq 0}$ is uniformly elliptic, with ellipticity constant $\epsilon_0$ which only depends on $M(\beta)$ in \eqref{e.M-beta}.  
\end{prop}

\begin{proof}
Let $\mathfrak S_n:=\mathfrak P_n$, and let $\pi_{n,n+1}(\cdot,\cdot)$ denote the transition kernel of $(X_n)_{n\geq 0}$. 
First we claim that the following properties hold:
\begin{enumerate}[(1)]
\item $\mathfrak S_n$ have bounded cardinality;
\item $\exists\epsilon_0'>0$  such that for all $n\geq 1$, $x\in \mathfrak S_n$ and $y\in\mathfrak S_{n+1}$, either $\pi_{n,n+1}(x,y)=0$ or 
$\pi_{n,n+1}(x,y)\geq \epsilon_0'$;
\item For every $n\geq 0$, $x\in\mathfrak S_n$ and $z\in \mathfrak S_{n+2}$, there exists some $y\in\mathfrak S_{n+1}$ for which $\pi_{n,n+1}(x,y)\pi_{n+1,n+2}(y,z)>0$. 
\end{enumerate}

Part (1) is Propositions \ref{Prop-Markov}(4). Part (2) follows Proposition \ref{Prop-Markov-Kernel}(3), by noting that $\ell^u(\cdot)$ are uniformly bounded away from zero and infinity by Lemma \ref{Lemma-Partition}(10), and $\lambda_i^u$ are bounded away from zero and infinity by Proposition \ref{Prop-Contracting-Automorphisms}(4).

To see part (3), it is sufficient to find for every  $a\in \mathfrak P_n$ and $c\in\mathfrak P_{n+2}$ some $b\in\mathfrak P_{n+1}$ such that
$t^{(n)}_{a,b}, t^{(n+1)}_{b,c}\neq 0$. 

Write $a=Q_{\alpha,\beta}^n$ and $c=Q_{\gamma,\delta}^{n+2}$.
By Lemma \ref{Lemma-Partition},    $\mathrm{int}(c)\subset \psi_{n+1}(\mathrm{int}(Q_\gamma(n+1)))$.
Let $W$ be some $u$-fibre in $\mathrm{int}(a)$. By Lemma \ref{Lemma-Bounded-Geometry-Gamma(M)}, 
$\ell(W)\!\!\geq\!\!\Gamma(M)^{-1}.$  By Proposition \ref{Prop-Contracting-Automorphisms}(2)(b) and Equation \eqref{e.choice-of-lambda},
 $\ell(\psi_n(W))\geq T_0(M)$.  By the choice of $T_0(M)$ in \S\ref{ss.density-const}, $\psi_n(W)$ must intersect $\mathrm{int}(Q_\gamma(n+1))$. 
 
 Necessarily, $\psi_n(\mathrm{int}(a))$ intersects the interior of some $b\in\mathfrak P_{n+1}$ inside $Q_\gamma(n+1)$. Thus $t^{(n)}_{a,b}=1$. 
Since $b\subset Q_{\gamma}(n+1)$ and $c=Q_{\gamma,\delta}^{n+2}$, $\psi_{n+1}(b)$ intersects $\mathrm{int}(c)$ (Lemma~\ref{Lemma-Partition}(7)), so $t_{b,c}^{(n+1)}=1$, and we are done.

\smallskip
Conditions (1)--(3) imply uniform ellipticity with $\mu_n:=$uniform measure on $\mathfrak S_n$, 
$p_n(x,y):=\pi_{n,n+1}(x,y)\cdot|\mathfrak S_n|$,
and $\epsilon_0:=(\epsilon_0'/M')^2$,  $M':=\max\{|\mathfrak S_n|:n\geq 0\}$. 
These constants  depend only on the bounded geometry estimates in 
Lemma \ref{Lemma-Bounded-Geometry-Gamma(M)} and on the uniform bounds on $\|d\psi_i\|, \|d\psi_i^{-1}\|$. Therefore, they depend only on $M(\beta)$. 
\end{proof}

Uniform ellipticity implies exponential $\alpha$-, $\rho$-, $\phi$-, and $\psi$-mixing (see Propositions 2.13 and  2.25 in \cite{Dolgopyat-Sarig-Book}). The following consequences suffice for our purposes.
Let $\mathcal F_0^{k}$ denote the $\sigma$-algebra generated by $X_0,\ldots,X_k$, and let $\mathcal F_{\ell}^\infty$ denote the  $\sigma$-algebra generated by $X_{\ell}$, $\ell\geq k$.\label{page-sigma-fields} Proposition \ref{Prop-Doeblin} implies \index{Mixing}\index{Exponential mixing}
 \begin{cor}[Exponential $\rho$-Mixing]\label{Cor-rho-mixing}
There exist constants $C_{mix}>1$ and $\rho\in (0,1)$, which only depend on the ellipticity constant $\epsilon_0$, such that for any $k,n$ and   square integrable functions $h_1, h_2$ such that $h_1$ is $\mathcal F_0^k$-measurable and $h_2$ is $\mathcal F_{k+n}^\infty$-measurable,
$$
\left|\E\bigl(h_1 h_2\bigr)-\E\bigl(h_1\bigr)\E\bigl(h_2\bigr)\right|\leq C_{mix}
\rho^n \sqrt{\mathrm{Var}(h_1)\mathrm{Var}(h_2)}.
$$
\end{cor}
 \begin{cor}[Exponential $\phi$-Mixing]\label{Cor-phi-mixing}
There exist constants $C_{mix}>1$ and $\rho\in (0,1)$, which only depend on the ellipticity constant $\epsilon_0$, such that for any $n>m>0$ and a bounded  function $h_n:\mathfrak P_n\to\R$, 
$
\left|\E\bigl(h_n(X_n)|X_m\bigr)-\E\bigl(h_n(X_n)\bigr)\right|\leq C_{mix}
\rho^{n-m}\|h_n\|_\infty.
$
\end{cor}

\section{A Random Walk Driven by the Markov Chain}
\label{ScMarkovReg}
Let  $\mathsf X=\{X_k\}_{k\geq 0}$ denote the Markov chain from Proposition \ref{Prop-Markov-Kernel}, and let $\mathsf{f}=\{f_k\}_{k\geq 1}$  denote the symbolic jump functions constructed in Lemma \ref{Lemma-f(X)}. The random walk
$$
S_N=f_1(X_1)+\cdots+f_N(X_N).
$$
 plays a central role in the proofs of the main results. 
 
 The purpose of this section is to show that $(\mathsf X,\mathsf f)$ satisfy  the sufficient 
 conditions
  of the local limit theorems in \cite{Dolgopyat-Sarig-Book} (summarized for the convenience of the reader in  Appendix \ref{AppMC}).  There are four  conditions: uniform boundedness of $f_k(X_k)$, uniform ellipticity of $\mathsf X$; $\mathrm{Var}(S_N)\to\infty$; and irreducibility (see below). 
Uniform ellipticity and uniform boundedness  were already verified in \S\S\ref{ss.Jump-Functions} and \ref{ss.uniform-ellipticity}. The proof of the other two conditions uses  the {\em structure constants} of $(\mathsf X,\mathsf f)$, and we begin by recalling their definition from \cite[Chapter 4]{Dolgopyat-Sarig-Book}.

\subsection{Structure Constants}\label{ss.structure-constants} 
Given $n\geq 3$, let $\mathrm{Hex}(n)$ denote the collection of arrays \index{Structure constants}\index{Hexagons}\index{Balance of a hexagon}
$
\left(
x_{n-2} \begin{array}{l}
x_{n-1}\\
y_{n-1}
\end{array}
\begin{array}{l}
x_{n}\\
y_{n}
\end{array} y_{n+1}
\right)
$ with the following {\em admissibility conditions}: $x_i,y_i\in\mathfrak P_i$, and 
\begin{align*}
&\Prob[(X_{n-2},\ldots,X_{n+1})=(x_{n-2};\,x_{n-1},x_{n};\, y_{n+1})|X_{n-2}=x_{n-2}]\neq 0;\\
&\Prob[(X_{n-2},\ldots,X_{n+1})=(x_{n-2};\,y_{n-1},y_{n};\, y_{n+1})|X_{n-2}=x_{n-2}]\neq 0.
\end{align*}
The elements of $\mathrm{Hex}(n)$ are called  {\em hexagons at position $n$}.

The {\em balance} of a hexagon is defined by  
$$
\Gamma\left({
x_{n-2} \begin{array}{l}
x_{n-1}\\
y_{n-1}
\end{array}
\begin{array}{l}
x_{n}\\
y_{n}
\end{array} y_{n+1}}
\right):=[f_{n-1}(x_{n-1})+f_n(x_n)]-[f_{n-1}(y_{n-1})+f_n(y_n)].
$$
We think of this as the difference between the sums of the jump functions along the paths $(x_{n-2},x_{n-1},x_n,y_{n+1})$  and $(x_{n-2},y_{n-1},y_n,y_{n+1})$.

In \cite{Dolgopyat-Sarig-Book}, we introduced a probability measure  $m_{\mathrm{Hex}}^{n}$  on $\mathrm{Hex}(n)$, called the {\em hexagon measure}, which represents the following sampling procedure. Let $\{Y_n\}$ be an  independent copy of our Markov chain  and $\{X_n\}$, then 
\begin{itemize}
\item $(x_{n-2},x_{n-1})$ is sampled from the distribution of $(X_{n-2},X_{n-1})$;
\item $(y_n,y_{n+1})$ is sampled from the distribution of $(Y_n,Y_{n+1})$;
\item $x_{n}$ and $y_{n-1}$ are conditionally independent given the previous choices, and 
\begin{align*}
&\Prob(x_{n}\in E|x_{n-1},y_{n+1})=\Prob(X_{n}\in E|X_{n-1}=x_{n-1}, X_{n+1}=y_{n+1}), \\
&\Prob(y_{n-1}\in E|x_{n-2},y_{n})=\Prob(X_{n-1}\in E|X_{n-2}=x_{n-2},  X_{n}=y_{n}). 
\end{align*}
\end{itemize}

\begin{lem}\label{lemma-Hex-Meas}
 There exists $c_0>0$ so that  $m_{\mathrm{Hex}^n}(H)>c_0$ for every $n$ and  $H\in\mathrm{Hex}(n)$. 
\end{lem}
\begin{proof}
By Proposition \ref{Prop-Markov-Kernel}, the probabilities of any pair $\Prob[X_k=a,X_{k+1}=b]$ and any triplet $\Prob[X_k=a,X_{k+1}=b,X_{k+2}=c]$ are either zero, or are bounded away from zero. By  Proposition \ref{Prop-Doeblin}, the probability of any event  $\Prob[X_{k-1}=a, X_{k+1}=b]$ is bounded away from zero. The lemma follows by direct calculation.
\end{proof}

\begin{defi}\cite{Dolgopyat-Sarig-Book}
The {\em structure constants} of $(\mathsf X,\mathsf f)$  are the numbers
\begin{align}
& u_n:= \bigl(\E_{m_{\mathrm{Hex}}^n}(\Gamma^2)\bigr)^{1/2}\label{e.u},\\
&d_n(\xi):=\bigl(\E_{m_{\mathrm{Hex}}^n}(|e^{i\xi \Gamma}-1|^2)\bigr)^{1/2}\ \ \ \ (\xi\in\R)\label{e.d}.
\end{align}
\end{defi}
Our main result on the structure constants of $(\mathsf X,\mathsf f)$ is the following. 
\begin{thm}\label{Theorem-D_N-Growth} Fix $M>1$. There are  constants $c_\xi>0$ such that \index{Structure constants!growth bounds}  every $\alpha\in\BT(M)$ (see \eqref{e.BT}), the Markov chain and jump functions from \S\ref{ss.Jump-Functions} and \ref{ss.IMC} satisfy the following:
\begin{enumerate}[(1)]
\item For every $N$ large enough,
$
\DS c_0 N\leq \sum_{n=3}^N u_n^2\leq c_0^{-1}N.
$
\item For every $\xi\not\in 2\pi\Z$, for every $N$ large enough,
$
\DS \sum_{n=3}^N d_n^2(\xi)\geq c_\xi N.
$
\end{enumerate}
\end{thm}
\noindent

We stress that this result depends in a crucial way on the assumption that $\alpha$ has bounded type.  
The proof relies on several lemmas. Let $\|x\|:=$ distance of $x$ from $\Z$, and recall that $\tau(x):=2\cdot 1_{[0,\frac{1}{2})}(\{x\})-1$.

\begin{lem}  For every $\eps>0$, there are
$N_{\eps, M},  \delta_{\epsilon,M}>0$ such that for every $\alpha\in\BT(M)$, \begin{equation}\label{e.q}
\exists1\leq q\leq N_{\eps,M}\text{ s.t. } \delta_{\epsilon,M}<\|q\alpha\|<\frac{\eps}{4}\text{ and } \sum_{j=0}^{q-1} \tau(z +j\alpha)=\pm 1 \text{ for all }z\in\T.
\end{equation}
\end{lem}
\begin{proof} Denote the principal denominators of $\alpha$ by $q_k$. 
There are constants $1<M_1<M_2$ which only depend on $M$ so that 
    $
    M_1^k\leq q_k\leq M_2^k\text{ for all }k\geq 1
    $: Write $\alpha=[a_0;a_1,a_2,\ldots]$ and recall from \cite{Khintchine-Continued-Fractions} that $q_{n+1}=a_{n+1}q_n+q_{n-1}$, $q_0=0, q_1=a_1$. Then we can take $M_2:=2M$, and $M_1=\frac{1+\sqrt{5}}{2}$ (the lower bound for Fibonacci's sequence).

Fix $j_\eps\in\N$ so large that $M_1^{j_\eps}>\frac{1}{\sqrt{\eps/4}}$ and let
$
N_{\eps,M}:=M_2^{j_\eps+5}
$,  then 
$q_{j_\eps}\!\!>\!\!2\eps^{-1/2}$, and $q_{j_\eps+4}\!\!<\!\!N_{\eps,M}$.
It is well-known that for at least one 
$j_\eps\!\leq\! j\!\leq\!  j_{\eps}+3$,  $q_j$ is an {\em odd} principal denominator satisfying 
$\DS
\|q_j \alpha\|<\frac{1}{2q_j^2}
$ 
(see e.g. \cite[Lemma 1]{Huveneers}).
Positive integers $q$ such that $\|q\alpha\|<\frac{1}{2q^2}$ are called {\em Legendre denominators}. \index{Legendre denominators}\index{Odd Legendre denominators}
{\em Odd}  Legendre denominators $q$ satisfy
$
\DS\sum_{j=0}^{q-1}\tau(x+j\alpha)=\pm 1\text{ on }\T,
$
see \cite[\S3]{Nakada-Keio} or \cite[Sublemma 1.2]{ANSS}. 
By construction, $q_j< N_{\epsilon, M}$, $\|q_j\alpha\|<\frac{1}{2q_j^2}<\frac{\eps}{4}$, 

By  \cite[Theorems 9,13]{Khintchine-Continued-Fractions},
 $\DS \|q_j\alpha\|>1/(q_j+q_{j+1})>\frac{1}{2}M_2^{-N_{\eps,M}}=:\delta_{\eps,M}
$.
\end{proof}

Recall that $D$ is the deck transformation which increases the $\Z$-coordinate by $+1$. Let $\pi:\St\to\St_0$ be the natural projection, and let $\dist_{\St}$ 
and $\dist_{\St_0}$ denote the distance functions on $\St$ and $\St_0$. 

\begin{lem}\label{Lemma-T} For every $\eps>0$, there exist 
$T_{\eps, M}$ such that for every $\alpha\in\mathrm{BT}(M)$ and  $\upsilon\!\!\in\!\!\St_0$, the following holds: Every line segment in $\St$ of length $T_{\eps,M}$ and direction  ${2\alpha-1 \choose 1}$ contains two points $\wt{\omega},\wt{\varpi}$ such that
$\dist_{\St}(\wt{\varpi},D(\wt{\omega}))<\eps$ and  $\dist_{\St_0}(\pi(\wt{\omega}),\upsilon)<\eps$.
\end{lem}

\begin{proof}
Let $\Sigma$ and $\omega:\T\times\Z\to\Sigma$ be as in Theorem \ref{t.HHW}, and let $\beta:=2\alpha-1$. %
Choose $N_{\epsilon,M}$, $\delta_{\epsilon,M}$ and $1\leq q\leq N_{\epsilon,M}$ as in the previous lemma.

Set $K_M\!\!:=\!\!\sqrt{1+(2M+3)^2}$. 
The linear flow $\phi$ traverses a  distance  $r\!\!:=\!\!\sqrt{1+\beta^2}\leq K_M$ between consecutive visits to the Poincar\'e section $\Sigma$. 
It is sufficient to consider the case when  $\upsilon$ is in the projection of $\Sigma$ to $\St_0$, otherwise we slide it in direction ${\beta\choose 1}$  no more than $K_M$ units of distance to $\Sigma$, and prove the lemma for the resulting point.
 We will show that for such $\upsilon$, the lemma holds with  $T_{\epsilon,M}:=3N_{\epsilon,M}K_M/\delta_M.$ 

Let $L$ be an arbitrary line segment in $\St$ with initial point $p$,  direction ${\beta\choose 1}$, and length bigger than $T_{\epsilon,M}$. 

$\bullet$ Let $\wt{\omega}_0:={\omega}(x_0,n_0)$ denote the first intersection of $L$ with $\Sigma$. This point lies within distance $r$ from $p$, along $L$. 

$\bullet$ The $(k+1)$-th intersection of $L$ with $\Sigma$ is $\wt{\omega}_k:={\omega}(\{x_0+k\alpha\}, n_0+\tau_k(x_0))$, where 
    $
    \DS \tau_k(x)=\sum_{j=0}^{k-1}\tau(x+k\alpha)
    $. It is within distance $(k+1)r$ from $p$, in the direction of $L$. 

$\bullet$ Choose $q$ as in \eqref{e.q}. 
Since $\delta_{\epsilon,M}<\|q\alpha\|<\frac{\epsilon}{4}$,  
$\{\{x_0+k\alpha\}: 0\leq k\leq N_{\epsilon,M}/\delta_{\epsilon,M}\}$ is $\epsilon/4$-dense in $\T$, and  there exists an integer $0\leq k^\ast\leq N_{\eps,M}/\delta_{\epsilon,M}$ such that 
$$
\dist_{\St}(\wt{\omega}_{k^\ast},\wt{\upsilon})<\frac{\epsilon}{2}\text{ for some $\wt{\upsilon}$ such that $\pi(\wt{\upsilon})=\upsilon$.}
$$
\item Let $k^\#\!\!:=\!\!k^\ast+q$ and $\sigma\!\!:=\!\!\tau_q(x_0+k^\ast\alpha)$. By the choice of $q$, 
$\sigma\!\!\in\!\!\{\!-\!1,\!+\!1\}$, and by Theorem~\ref{t.HHW},
$\tomega_{k^\#}=\omega(T_\alpha^q(\{x_0+k^\ast\alpha\},n_{0}+\tau_{k^\ast}(x_0)))
=\omega(\{x_0+k^\ast\alpha+q\alpha\},n_{0}+\tau_{k^\ast}(x_0)+\sigma).
$

Comparing this to $D^{\sigma}(\tomega_{k^\ast})=\omega(\{x_0+k^\ast\alpha\},n_{0}+\tau_{k^\ast}(x_0)+\sigma)$, we find that 
\begin{align*}
&\dist_{\St}(\wt{\omega}_{k^\#},D^{\sigma}(\tomega_{k^\ast}))\leq  2\|q\alpha\|<\frac{\epsilon}{2}.
\end{align*}
 (the factor $2$ is because of the formula for $\omega(\cdot,\cdot)$). By the choice of $q$, $|\sigma|=1$.

We obtained two points, $\wt{\omega}_{k^\ast}$ and $\wt{\omega}_{k^\#}$, which project to the $\eps$-neighborhood of $\upsilon$, such that $\dist_{\St}(\wt{\omega}_{k^\#},D^{\sigma}(\wt{\omega}_{k^\ast}))<\frac{\epsilon}{2}$, where $|\sigma|=1$.
These points are both on $L$, because their distance from  $p$ along $L$ is at most 
$
(k^\#+1)r\leq (2N_{\epsilon,M}/\delta_{\eps,M}+1)K_M<T_{\epsilon,M}
$.   
\end{proof}

\begin{lem}\label{Lemma-D_N-Growth}
For every $\xi\not\in 2\pi\Z$, there exists a positive constant $c_\xi$  so that for all $N$ large enough, 
for all $\alpha\in\BT(M)$,
$
\DS \sum_{n=3}^N d_n^2(\xi)\geq c_{ \xi} N.
$
\end{lem}

\begin{proof}

Let $\mathfrak H(2m,n)$ denote the collection of arrays 
$
\left(
x_n \begin{array}{l}
x_{n+1}\\
y_{n+1}
\end{array}
\cdots
\begin{array}{l}
x_{n+m-1}\\
y_{n+m-1}
\end{array} y_{n+m}
\right)
$ such that  $x_i,y_i\in\mathfrak P_i$ and
$$
\Prob((X_n,\ldots,X_{n+m})=(x_n;\,x_{n+1},\ldots,x_{n+m-1};\, y_{n+m})|X_n=x_n)\neq 0;
$$
$$
\Prob((X_n,\ldots,X_{n+m})=(x_n;\, y_{n+1},\ldots,y_{n+m-1};\,y_{n+m})|X_n=x_n)\neq 0.
$$
The elements of $\mathfrak H(2m,n)$ are called  {\em admissible $2m$-gons} (at position $n+2$).\footnote{This scheme of enumerating the positions is in accordance with the conventions of \cite{Dolgopyat-Sarig-Book}.} 
The {\em balance} of an admissible $2m$-gon is 
$$
\Gamma\left({
x_n \begin{array}{l}
x_{n+1}\\
y_{n+1}
\end{array}
\cdots
\begin{array}{l}
x_{n+m-1}\\
y_{n+m-1}
\end{array} y_{n+m}}
\right):=\sum_{k=1}^{m-1} f_{n+k}(x_{n+k})-\sum_{k=1}^{m-1} f_{n+k}(y_{n+k}).
$$

\medskip
\noindent
{\em Claim 1.\/} For some $m\geq 3$, for every $n\geq 3$, $\Gamma(H)=1$ for some $H\in\mathfrak H(2m,n)$.

\medskip
\noindent
{\em Proof of the Claim.\/} 
Let  $\lambda_1:=1/\max \{\|d\tpsi_k^{-1}\|:k\geq 0\}$ (this is a finite set) and $\lambda_2:=\lambda^{-1}$ with $\lambda$ as in Proposition \ref{Prop-Contracting-Automorphisms}, then 
 $$
 \lambda_2^{-1}\leq \|d\wt{\psi}_k^{-1}|_{E^s_{k+1}}\|\leq \lambda_1^{-1}\text{ for all $k$}.
 $$

 By Lemma~\ref{Lemma-Partition}(11), there is a positive constant $0<w<1$ such that every element of every  $\mathfrak P_k$ contains a ball of radius $w$.
Let $\eps:=w/3$, and  let $M$ be a bound for the digits of the regular CFE of $\alpha$ and $1/(2\alpha-1)$. Choose $T_{\eps,M}$ as in Lemma \ref{Lemma-T}, making sure it is large enough to also satisfy the conclusion of Lemma \ref{LmDenseWs}.

Let   $m:=\lceil \log_{\lambda_2^{-1}}(T_{\eps,M}/w)\rceil$, and fix
 $n\geq 3$.
 Choose some arbitrary  $Q_{n+m}\in\mathfrak P_{n+m}$. By the definition of $w$, $Q_{n+m}$ contains in its interior a linear segment ${L}_{n+m}$ of length $w$ and direction $E^s_{n+m}$. Lift it to some segment $\wt{L}_{n+m}$ in the staircase. 
 Since $L_{n+m}$ does not intersect the boundary of any element of $\mathfrak P_{n+m}$, $\wt{L}_{n+m}$ does not intersect the boundary of $D^k(\wt{F}_{n+m})$ for any $k\in\Z$, and the $\wt{\zc}_{n+m}$-coordinate remains constant on $\wt{L}_{n+m}$.

 Let
$
\wt{L}_n:=\wt{\psi}_n^{-1}\wt{\psi}_{n+1}^{-1}\cdots \wt{\psi}_{n+m-1}^{-1}(\wt{L}_{n+m}).
$
This is a linear segment in direction $E^s_n$, and 
$\DS
|\wt{L}_n|\geq \lambda_2^{-m}w>T_{\eps,M}.
$

Fix some $\upsilon$ at the center of a ball with radius $w$, completely inside the interior of some element of $\mathfrak P_n$. 
By the choice of $T_{\epsilon, M}$,  $\wt{L}_n$ contains two points $\wt{\omega},\wt{\varpi}$ such that 
$$\dist_{\St}(\wt{\varpi},D(\wt{\omega}))<w/3\text{ and }\dist_{\St_0}(\pi(\wt{\omega}),\upsilon)<w/3.$$ 
The line segment $\gamma$ from  $\pi(\wt{\varpi})$ to $\pi(D(\tomega))$ is inside $B(\upsilon,w)$, therefore it does not intersect the boundaries of any of the elements of $\mathfrak P_n$. Therefore, the line segment $\wt{\gamma}$ from $\wt{\varpi}$ to $D(\tomega)$ does not cross the boundary of $D^k(\wt{F}_n)$ for any $k\in\Z$. It follows that 

\begin{equation}\label{e.omega-varpi-1}
 \wt{\zc}_n(\wt{\varpi})=\wt{\zc}_n(D(\tomega))=\wt{\zc}_n(\tomega)+1.
\end{equation}
Similarly, $(\wt{\psi}_{n+m-1}\cdots\wt{\psi}_n)(\wt{\omega}), (\wt{\psi}_{n+m-1}\cdots\wt{\psi}_n)(\wt{\varpi})\in \wt{L}_{n+m}$, so 
\begin{equation}\label{e.omega-varpi-2}
\wt{\zc}_{n+m}[(\wt{\psi}_{n+m-1}\cdots\wt{\psi}_n)(\wt{\omega})]=\wt{\zc}_{n+m}[(\wt{\psi}_{n+m-1}\cdots\wt{\psi}_n)(\wt{\varpi})]
\end{equation}

Let $\omega,\varpi$ denote the projections of $\wt{\omega},\wt{\varpi}$ 
to the punctured torus $\St_0$.
We may assume without loss of generality that the points $\omega_n:=\omega$, $\varpi_n:=\varpi$ and the points
\begin{equation}\label{e.orbit-of-xy}
\omega_{n+j}:=(\psi_{n+j-1}\cdots{\psi}_n)(\omega), \quad
\varpi_{n+j}:=(\psi_{n+j-1}\cdots{\psi}_n)(\varpi) \quad
(j=1,\ldots,m)
\end{equation}
all lie in the interior of the Markov elements of $\mathfrak P_{n+j}$. Otherwise we shift 
$\wt{L}_{n+m}$ in the $E^u_{n+m}$-direction  slightly, and we shift $\omega,\varpi$ 
 in the $E^s_n$ direction slightly.
 
We can now safely define for $0\leq j\leq m$,
\begin{align*}
&\xi_{n+j}:=\mathfrak P_{n+j}-\text{element which contains $\omega_{n+j}$ in its interior};\\
&\eta_{n+j}:=\mathfrak P_{n+j}-\text{element which contains $\varpi_{n+j}$ in its interior}.
\end{align*}
Notice that
$
\xi_n=\eta_n\text{ and }\xi_{n+m}=\eta_{n+m}.
$
The first identity is because
$$
\omega, \varpi \in B(\omega,2w/3)\subset\text{ the interior of the same element of }
\mathfrak P_n.
$$
The second identity is because $\wt{L}_{n+m}$ was chosen to lie completely inside the interior of a Markov element of $\mathfrak P_{n+m}$.

By Lemma \ref{Lemma-Cylinders} and the choice of $\omega_{n+j}, \varpi_{n+j}$, the cylinders $\<\dot{\xi}_n,\ldots,\xi_{n+m-1},\eta_{n+m}\>$ and 
$\<\dot{\xi}_n,\eta_{n+1}\ldots,\eta_{n+m}\>$ have non-empty interior. Therefore, 
$$
\Prob[(X_n,\ldots,X_{n+m})=(\xi_n;\,\xi_{n+1},\ldots,\xi_{n+m-1};\, \eta_{n+m})|X_n=\xi_n]\neq 0;
$$
$$
\Prob[(X_n,\ldots,X_{n+m})=(\xi_n;\, \eta_{n+1},\ldots,\eta_{n+m-1};\,\eta_{n+m})|X_n=\xi_n]\neq 0.
$$
Thus  $\left(
\xi_n \begin{array}{l}
\xi_{n+1}\\
\eta_{n+1}
\end{array}
\cdots
\begin{array}{l}
\xi_{n+m-1}\\
\eta_{n+m-1}
\end{array} \eta_{n+m}
\right)$ is an admissible $2m$-gon in position $n+2$. We have: 
\begin{align*}
&\Gamma\left(
\xi_n \begin{array}{l}
\xi_{n+1}\\
\eta_{n+1}
\end{array}
\cdots
\begin{array}{l}
\xi_{n+m-1}\\
\eta_{n+m-1}
\end{array} \eta_{n+m}
\right)
=\sum_{j=1}^{m-1}f_{n+j}(\xi_{n+j})-
\sum_{j=1}^{m-1}f_{n+j}(\eta_{n+j})\\
&=\sum_{j=1}^{m}f_{n+j}(\xi_{n+j})-
\sum_{j=1}^{m}f_{n+j}(\eta_{n+j}),\text{ because $\xi_{n+m}=\eta_{n+m}$}\\ 
&=\sum_{j=1}^{m}{f}_{n+j}(X^{(n+m)}_{n+j}(\omega_{n+m}))-
\sum_{j=1}^{m}{f}_{n+j}(X^{(n+m)}_{n+j}(\varpi_{n+m})), \text{ see }\eqref{e.X-new}\\
&=
\sum_{j=1}^{m}\wt{f}_{n+j}^{(n+m)}(\omega)-
\sum_{j=1}^{m}\wt{f}_{n+j}^{(n+m)}(\varpi), \text{ see Lemma \ref{Lemma-f(X)}}\\
&=\left[\wt{\zc}_{n+m}(\tomega_{n+m})-
\wt{\zc}_n(\tomega_n)\right]-
\left[\wt{\zc}_{n+m}(\wt{\varpi}_{n+m})-\wt{\zc}_n(\wt{\varpi}_n)\right], \text{ by \eqref{Eq-Modified-Frobenius}}\\
&=[\wt{\zc}_n(\wt{\varpi})-\wt{\zc}_n(\wt{\omega})]+
[\wt{\zc}_{n+m}((\wt{\psi}_{n+m-1}\cdots\wt{\psi}_n)(\wt{\omega}))-
\wt{\zc}_{n+m}((\wt{\psi}_{n+m-1}\cdots\wt{\psi}_n)(\wt{\varpi}))]\\
&=1+0=1,\text{ by \eqref{e.omega-varpi-1} and \eqref{e.omega-varpi-2}}.
\end{align*}

\medskip
\noindent
{\sc Claim 2.\/} \label{page-hexagon-claim}{\em For every $\ell\geq 2$ and $n\geq 0$ there exists $n\leq n'\leq n+m-3$ and an admissible hexagon at position $n'+2$ such that 
$
\DS \Gamma\left(
\xi_{n'} \begin{array}{l}
\xi_{n'+1}\\
\zeta_{n'+1}
\end{array}
\begin{array}{l}
\zeta_{n'+2}\\
\eta_{ n'+2}
\end{array} \eta_{n'+3}
\right)\neq 0\mod\ell
$.
}\\

\begin{figure}
  \centering
  \fbox{\includegraphics[width=7cm,angle=90]{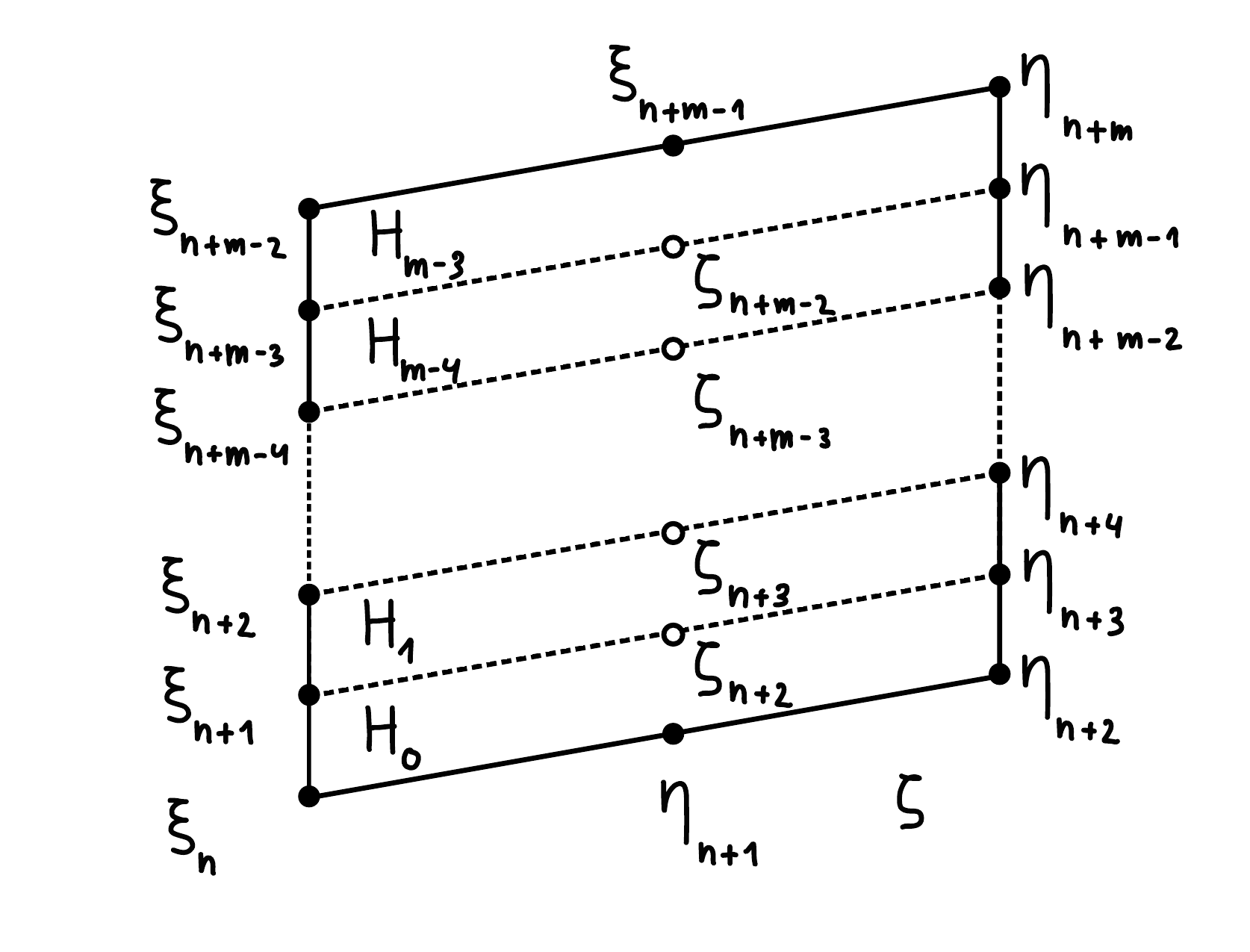}}
  \caption{A decomposition of a $2m$-gon to hexagons}\label{Figure-Hex-Decomposition}
\end{figure}

\noindent
{\em Proof of the Claim.\/} Let $\xi_i,\eta_j$ be in the proof of Claim 1. By uniform ellipticity  (Proposition \ref{Prop-Doeblin}), for every  $j=0,\ldots,m-1$ there is some $\zeta_{n+j+1}\in\mathfrak P_{n+j+1}$ such that
$$
\Prob[(X_{n+j},X_{n+j+1},X_{n+j+2})=(\xi_{n+j},\zeta_{n+j+1},\eta_{n+j+2})|X_{n+j}=\xi_{n+j}, X_{n+j+2}=\eta_{n+j+2}]\neq 0.
$$
Since $\xi_n=\eta_n$ and $\xi_{n+m}=\eta_{n+m}$, we may take
$
\zeta_{n+1}=\eta_{n+1}\text{ and }\zeta_{n+m-1}=\xi_{n+m-1}.
$

Let $H_i:=
\left(
\xi_{n+i} \begin{array}{l}
\xi_{n+i+1}\\
\zeta_{n+i+1}
\end{array}
\begin{array}{l}
\zeta_{n+i+2}\\
\eta_{n+i+2}
\end{array} \eta_{n+i+3}
\right).
$
These are admissible hexagons, and looking at Figure \ref{Figure-Hex-Decomposition}, we see that
$\DS
\Gamma\left(
\xi_n \begin{array}{l}
\xi_{n+1}\\
\eta_{n+1}
\end{array}
\cdots
\begin{array}{l}
\xi_{n+m-1}\\
\eta_{n+m-1}
\end{array} \eta_{n+m}
\right)=\sum_{i=0}^{m-3}\Gamma(H_i),
$
because the terms $f_{n+j}(\zeta_j)$ on  the dashed lines appear twice with opposite signs, and cancel out.
Since the left-hand side equals one,  $\Gamma(H_i)\neq 0\mod\ell$ for some $i$, proving the claim.

\medskip
Suppose $\xi\in 2\pi\Q\setminus 2\pi\Z$. Then  $\xi\in\frac{2\pi}{\ell}\Z$ for some $\ell\in\N$, $\ell\geq 2$. If $\Gamma(H)\neq 0\mod\ell$, then $e^{i\xi\Gamma(H)}\in \{e^{2\pi i k/\ell}:k=1,\ldots,\ell-1\}$, and therefore  $|e^{i\xi\Gamma(H)}-1|>{\brc_\xi}>0$ with $\brc_\xi$ independent of $H$. 
Thus the claim says that for every $n\geq 2$, there are some $n\leq n'\leq n+m-3$  and  $H\in\mathrm{Hex}(n'+2)$ such that $|e^{i\xi \Gamma(H)}-1|\geq \brc_\xi$. By Lemma \ref{lemma-Hex-Meas}, $ d_{n'}^2(\xi)\geq c_0 c_\xi^2$. Thus
$\DS\sum_{n=2}^N d_n^2(\xi)\geq \left(\frac{N}{m}-1\right)c_0c_\xi^2$, and the lemma follows for all $\xi\in 2\pi\Q\setminus 2\pi\Z$.

Suppose $\xi\in 2\pi(\R\setminus\Q)$. There is a constant such that $|\Gamma(H)|\leq M$ for all $H$, therefore if $\Gamma(H)\neq 0$, then  $e^{i\xi\Gamma(H)}\in \{e^{ik\xi}: k\in\Z\cap [-M,M], k\neq 0\}$. By the irrationality of  $\xi/2\pi$, this set does not contain $1$, and since it is finite, $|e^{i\xi\Gamma(H)}-1|>\brc_\xi>0$ with $\brc_\xi$ independent of $H$. Continuing as before, we obtain the lemma for $\xi\in 2\pi(\R\setminus\Q)$.
\end{proof}

\begin{cor}\label{e.U_N-Growth}
Let $u_n:=\E_{m_{\mathrm{Hex}}^n}(\Gamma^2)^{1/2}$. 
There exists constants $c_1,c_2>0$ such that for all $N$ large enough, 
$
\DS c_1N\leq \sum_{n=3}^N u_n^2\leq c_2 N.
$
\end{cor}
\begin{proof}
Recall that $f_k$ are uniformly bounded. Therefore, the balance functions $\Gamma$ are uniformly bounded. So $u_n^2$ are bounded above, and the upper bound follows.

By the  general inequality $|e^{i\theta}-1|^2=4\sin^2\frac{\theta}{2}\leq \theta^2$,  
$ |d_n(\xi)| \leq u_n |\xi|$. Taking $\xi=\frac{2\pi}{\ell}$ with any $\ell\geq 2$, and looking at   Lemma \ref{Lemma-D_N-Growth}, we deduce the linear lower bound. 
\end{proof}

 Together, Lemma \ref{Lemma-D_N-Growth} and Corollary \ref{e.U_N-Growth} prove the Theorem \ref{Theorem-D_N-Growth}. 
 
\subsection{Irreducibility}
\label{SSIrred}

The {\em algebraic range} of $(\mathsf X,\mathsf f)$ is the smallest closed additive subgroup $G_{alg}\subset\R$ for which there are constants $\gamma_k$ such that $f_k(X_k)\in\gamma_k+G_{alg}$ a.s. for all $k$. Such a group exists, see \cite[Lemma 4.15]{Dolgopyat-Sarig-Book}.\index{Algebraic range} 

The {\em essential range}\index{Essential range} of $(\mathsf X,\mathsf f)$ is the smallest closed additive subgroup $G_{ess}\subset\R$ for which there exist  random variables  $h_k(X_k,X_{k+1})$ such that  $f_k(X_k)-h_k(X_k,X_{k+1})\in G_{ess}$ a.s. for all $k$, and so that
  $\DS {\cH}_N:=\sum_{k=1}^N h_k(X_k,X_{k+1})$ is {\em center-tight}, i.e.:\index{Center-tightness}
\begin{equation}\label{e.center-tightness}
\forall\epsilon\, \exists M>0, m_N\in\R \text{ s.t. } \forall N,\  \Prob\left(\left|\cH_N-m_N\right|>M\right)<\epsilon.
\end{equation}
Such a group exists, see \cite[Theorem 4.5]{Dolgopyat-Sarig-Book}. We say that $(\mathsf X,\mathsf f)$ is {\em irreducible}, if its essential range and its algebraic range are equal.\index{Irreducibility} 

\begin{prop}[Irreducibility]\label{Prop-Irred}
$(\mathsf X,\mathsf f)$ is irreducible, and its algebraic and essential ranges are both equal to $\Z$.
\end{prop}
\begin{proof}
The {\em co-range} of $(\mathsf X,\mathsf f)$ is the set\index{Co-range}\index{Essential range!and the co-range}
$
H:=\{\xi\in\R:\sum d_n(\xi)^2<\infty\}.
$ The following facts are proved in \cite[Chapter 4]{Dolgopyat-Sarig-Book}:
\begin{enumerate}[(1)]
\item $H$ is a closed additive group of $\R$, thus $H=\{0\}, t\Z$, or $\R$;
\item If $H=\{0\}$ then the essential range is $\R$; If $H=t\Z$ with $t\neq 0$ then the essential range is $\frac{2\pi}{t}\Z$; And if $H=\R$, then the essential range is $\{0\}$. 
\end{enumerate}
In our case, $H\!\!\supset\!\! 2\pi\Z$ (because $\Gamma$ takes values in $\Z$, so $d_n(2\pi)=0$).  But $H\cap (\R\setminus 2\pi\Z)\!\!=\!\!\emptyset$ by  Theorem \ref{Theorem-D_N-Growth}. Necessarily $H=2\pi\Z$, and the essential range is $\Z$. 
Since $f_k$ take values in the integers,  the algebraic range is contained in $\Z$. However it also contains the essential range, so it  must also equal $\Z$. 
\end{proof}

\subsection{Drift and Variance} Let 
$
S_N=f_1(X_1)+\cdots+f_N(X_N).
$
\begin{prop}\label{Prop-Drift}
There exists a constant $E$ which only depends on $M(\beta)$ such that $|\E(S_N)|\leq E$ for all $N$.\index{Drift of $S_N$}
\end{prop}
\begin{proof}
Let $\mu_0$ denote the normalized area measure on $\St_0$. 
By  Lemma \ref{Lemma-f(X)}
\begin{align*}
&\E(S_N)=\int_{\St_0}\sum_{k=1}^N f_k(X^{(N)}_k(\omega))\, d\mu_0=\int_{\St_0}\sum_{k=1}^N \wt{f}_k^{(N)}(\omega)\, d\mu_0.
\end{align*}
For each $\omega\in\St_0^\ast$, let $\tomega_N\in \St^\ast$ be the unique lift of $\omega$ such that ${\zc}(\tomega_N)=0$, and let 
$\tomega_{k}:=(\tpsi_k^{-1}\circ\cdots\circ\tpsi_{N-1}^{-1})(\tomega_N)$. 
Then
$
\wt{f}^{(N)}_{k+1}(\omega)=\wt{\zc}_{k+1}(\tomega_{k+1})-\wt{\zc}_k(\tomega_k)
$ by definition, whence 
\begin{align*}
&\E(S_N)=\int_{\St_0} [\wt{\zc}_{N}(\tomega_{N})-\wt{\zc}_0(\tomega_0)] d\mu_0(\omega)=
\int_{\St_0} [\zc(\tomega_{N})-\zc(\tomega_0)] d\mu_0(\omega)+E_N,
\end{align*}
where the $\DS |E_N|\leq E:=2\sup_{N}\sup_{\tomega\in\St^\ast}|\wt{\zc}_N(\tomega)-\zc(\tomega)|$. $E$ is finite by Lemma \ref{Lemma-Modified-Approx-Canonical}, and by the bounded geometry of $\Pi(n)$, 
$E$ has 
an upper  bound which depends only on $M(\beta)$. 

Since $\omega$ is uniformly distributed on $\St_0$, and $\pi:\Rect_0\to\St_0$ is measure preserving, $\tomega_N$ is uniformly distributed on $\Rect_0$. So
\begin{align*}
&\E(S_N)=\int_{\Rect_0} [\zc(\tomega_N)-\zc((\tpsi_0^{-1}\circ\cdots\circ\tpsi_{N-1}^{-1})(\tomega_N)] d\mu_0(\tomega_N){\pm E} \\
&=\int_{\Rect_0} [\zc(\tpsi_{N-1}\circ\cdots\circ\tpsi_0)(\tomega_0))-\zc(\tomega_0)] d\mu_0(\tomega_0)
{\pm E} 
\text{ because $\tpsi_i$ are area preserving}\\
&=\delta(\tpsi_{N-1}\circ\cdots\circ\tpsi_0){ \pm E} 
\text{ where $\delta(\cdot)$ is the average drift defined in \eqref{Eq-Drift-Defn}.}
\end{align*}
Now recall that  $\delta(\tpsi_{N-1}\circ\cdots\circ\tpsi_0)=\sum \delta(\tpsi_i)$ (see \cite{Avila-Dolgopyat-Duryev-Sarig}), and  $\delta(\tpsi_i)=0$.
\end{proof}

\begin{prop}\label{Prop-Var-Growth}
There exist positive constants $C_1,C_2$ which only depend on $M(\beta)$, such that for all $N$ large enough, \index{Variance of $S_N$}
$
C_1 N\leq \mathrm{Var}(S_N)\leq C_2 N.
$
\end{prop}
\begin{proof} 
By Theorem \ref{LmVarCycles}
there are positive constants $C_1$ and $C_2$ which only depend on the ellipticity constant and on $\sup_{n\geq 1} \|f_n\|$ such that 
$$
C_1^{-1}\sum_{n=3}^N u_n^2-C_2\leq \mathrm{Var}(S_N)\leq C_1 \sum_{n=3}^N u_n^2+C_2.
$$
The proposition follows from Theorem \ref{Theorem-D_N-Growth}. 
\end{proof}

\section{Random Walks Driven by Renormalizations on the Staircase}\label{s.RW-Driven-By-Renormalizations}
\subsection{Random Walks Driven by Renormalizations}  \label{s.Forward-RW} 
Suppose  $\wt{\phi}_k:\St\to\St$ are  homogeneous automorphisms, and $\Prob$ is a  probability measure on $\St$. The  {\em random walk}  {\em driven by $\wt{\phi}_k$} with {\em initial distribution} $\Prob$ is
$$
\zc[(\wt{\phi}_{n-1}\circ\cdots\circ \wt{\phi}_0)(\tomega)]-\zc(\tomega),\ \text{ $\tomega$ sampled randomly according to $\Prob$}. 
$$  
This work uses three   random walks like that: \index{Random walks driven by renormalizations}

\noindent
$\bullet$  {\bf Forward Random Walk $\boldsymbol{\Upsilon_n}$:} Let $\wt{\phi}_k$ be the homogeneous automorphism $\tvf_k$ with zero drift and derivative \index{Random walks driven by renormalizations!forward ($\Upsilon_N,\Phi_N$)} 
 ${\tiny \left(
       \begin{array}{cc}
         1 & -2m_{2k+2} \\
         0 & 1 \\
       \end{array}
     \right)
     \left(
       \begin{array}{cc}
         1 & 0 \\
         -2m_{2k+1} & 1 \\
       \end{array}
     \right)}$, see \eqref{e.A-k-dagger}. The forward random walk is 
 $\Upsilon_n(\tomega)=\zc[\wt{\Psi}_n(\tomega)]-\zc(\tomega)=\zc[\wt{\Psi}_n(\tomega)]$, where $\tomega$ sampled randomly uniformly from $\Rect_0$, and $\wt{\Psi}_n(\tomega)$ is the homogeneous automorphism with zero drift and derivative matrix $A_n(\beta)$, see \eqref{e.good-matrix}. 
  $\Upsilon_n=+z_n$ a.e., where $z_n$ is given in  \eqref{e.z-n}.  
     See  \S\ref{SSGP-Tlit} for the uses of $\Upsilon_n$ and $z_n$.
    
\smallskip
\noindent
 {\bf $\bullet$ Forward Random Walk $\boldsymbol{\Phi_n}$:} $\Phi_n(\tomega)=\zc[(\tpsi_{n-1}\circ\cdots\circ\tpsi_0)(\tomega)]-\zc(\tomega)$, $\tomega$ sampled randomly uniformly from $\Rect_0$. By \eqref{e.mat-prod-form}, $
\tpsi_k=\tvf_{kN_0+N_0-1}\circ\cdots\circ\tvf_{kN_0}$, so
$\Phi_n=\Upsilon_{nN_0}.
$ \\
$\Phi_n$ is easier to work with than  $\Upsilon_n$, because of the Markov partitions $\mathfrak P_k$.

 \smallskip
 \noindent
 {\bf $\bullet$ Backward Finite Random Walk $\boldsymbol{\Xi_n}$:} \index{Random walks driven by renormalizations!backward $(\Xi_n)$} $\Xi_n^{W}(\tomega)=\zc[(\tpsi_0^{-1}\circ\cdots\circ\tpsi_{n-1}^{-1})(\tomega)]-\zc(\tomega)$, $\tomega$ sampled randomly uniformly from a linear segment $W$. 
We will deviate from the setup of random walks,  and allow $W$ to  depend on $n$, so that $\Xi_n^{W_n}$ are separate random variables, without a joint distribution. In the case which interests us, $W_n$ are  linear segments  in directions $E^s_n$ and lengths bounded away from zero and infinity. This case appears in the study of Birkhoff integrals of the linear flow,   see \eqref{e.Int=Prob}. 

\smallskip
We will relate  $\Upsilon_n, \Phi_n$ and $\Xi_n$ to the simpler  process
$
S_n=f_1(X_1)+\cdots+f_n(X_n)$.

\begin{prop}\label{Prop-Markov-Representation-of-Phi}
There exists a uniformly bounded sequence of functions $\overline{\epsilon}_n:\St^\ast\to\Z$ so that the stochastic processes  \index{Random walks driven by renormalizations!and $S_N$} $\{\Phi_n-\overline{\epsilon}_n\}_{n\geq 1}$ and  $\{S_n\}_{n\geq 1}$ are equal in distribution.
\end{prop}
\begin{proof} Suppose $\tomega\in\Rect_0$,
then $\zc(\tomega)=0$, and 
\begin{align*}
\Phi_n(\tomega)&=\zc[(\tpsi_{n-1}\circ\cdots\circ\tpsi_0)(\tomega)]-\zc(\tomega)
=\wt{\zc}_n[(\tpsi_{n-1}\circ\cdots\circ\tpsi_0)(\tomega)]-\wt{\zc}_0(\tomega)
+\bar{\eps}_n(\tomega),\text{ with }\\
\bar{\eps}_n(\tomega)&:=(\zc-\wt{\zc}_n)((\tpsi_{n-1}\circ\cdots\circ\tpsi_0)(\tomega))-(\zc-\wt{\zc}_0)(\tomega).
\end{align*}

The functions $\bar{\eps}_n$ are uniformly bounded,  by Lemma \ref{Lemma-Modified-Approx-Canonical}. We claim that $(\Phi_n-\bar{\eps})_{n\geq 1}$ is equal in distribution to $S_n$. It is convenient to set  $\tomega_0:=\tomega$, and  
$
\tomega_k:=(\tpsi_{k-1}\circ\cdots\circ\tpsi_0)(\tomega_0)$. Then $\tomega=(\tpsi_{0}^{-1}\circ\cdots\circ\tpsi_{n-1}^{-1})(\tomega_n)$, and   
\begin{align*}
&\Phi_n(\tomega)=\wt{\zc}_n(\tomega_n)-\wt{\zc}_0(\tomega)
+\bar{\eps}_n(\tomega)\\
&=\sum_{k=0}^{n-1}[\wt{\zc}_{k+1}(\tomega_{k+1})-\wt{\zc}_k(\tomega_k)]+\bar{\eps}_n(\tomega)
=\sum_{k=0}^{n-1}f_{k+1}^{(n)}(\pi(\tomega_n))+\bar{\eps}_n(\tomega),\text{ see }\eqref{Eq-Modified-Frobenius}\\
&=\sum_{k=0}^{n-1}f_{k+1}\bigl(X_{k+1}^{(n)}(\pi(\tomega_n))\bigr)+\bar{\eps}_n(\tomega)=\sum_{k=1}^n f_k\big(X_k^{(0)}(\pi(\tomega))\bigr)+\bar{\eps}_n(\tomega),
\end{align*}
by Lemma \ref{Lemma-f(X)} and  \eqref{e.X_i-tranformation-law}. 
The map  $\pi:\Rect_0\to\St_0^\ast$ is area preserving and invertible. So  if $\tomega$ is sampled randomly uniformly in $\Rect_0$, then $\pi(\tomega)$ is distributed randomly uniformly in $\St_0$, and
$
\Phi_n-\bar{\eps}_n
$ is equal in distribution to $\DS \sum_{k=1}^{n} f_k(X_k^{(0)}(\omega'))$, where $\omega'$ is distributed randomly uniformly in $\St_0$. By Proposition \ref{Prop-Markov-Kernel}, $\Phi_n-\bar{\eps}_n=S_n$ in distribution.
\end{proof}

Let $[0; 2m_1,2m_2,\ldots]$ be the even CFE of $\beta:=2\alpha-1$, and let  $N_0$ be the number such that $\tpsi_{n-1}\circ\cdots\circ\tpsi_0=\wt{\Psi}_{nN_0}$ (see 
   Remark \ref{r.mat-prod-form}). Let  $\tvf_k$ be the unique homogeneous automorphism with  zero drift and derivative ${\tiny \left(
       \begin{array}{cc}
         1 & -2m_{2k+2} \\
         0 & 1 \\
       \end{array}
     \right)
     \left(
       \begin{array}{cc}
         1 & 0 \\
         -2m_{2k+1} & 1 \\
       \end{array}
     \right)}$.
      Define
\begin{equation}\label{e.E-cal}
\mathcal E:=
\{\tvf_{n+\ell}\circ\cdots\circ\tvf_n:\quad n\geq 0, \quad 0\leq \ell\leq N_0-1\}.
\end{equation}
Since $m_i$ is bounded,  $\mathcal E$ is finite.

 \begin{lem}\label{Lemma-Psi-and-Phi}
The slopes
 of $d\wt{e}^{-1}{1\choose 0}, d\wt{e}^{-1}{0\choose 1}$ for $\wt{e}\in\mathcal E$ are uniformly bounded away from the slopes of $E^u_n, E^s_n$ for all $n$.
\end{lem}
\begin{proof}
By Proposition \ref{Prop-Contracting-Automorphisms}, the slopes of $E^s_n,E^u_n$ have uniformly bounded  regular CFE.
On the other hand, since $m_i$ are bounded, the slopes of  $d\wt{e}^{-1}{1\choose 0}$ and $d\wt{e}^{-1}{0\choose 1}$ are rational numbers $p/q$ with uniformly bounded denominator, say  $1\leq q\leq Q$. 

The lemma now follows from
 the following well-known fact: For every $M$ there is a $c(M)>0$ such that for every irrational number $\sigma$ with regular continued fraction expansion all of whose digits are less than $M$,  for all relatively prime integers $p,q$, $|\sigma-\frac{p}{q}|\geq c(M)q^{-2}$. See e.g. the proof of Theorem 23 in \cite{Khintchine-Continued-Fractions}.
\end{proof}

Given $n\in\N$ and  $0\leq \ell\leq N_0-1$, let  $\wt{\vf}_{n,\ell}:=\tvf_{nN_0+\ell-1}\circ \cdots\circ\tvf_{nN_0}\in\mathcal E$, then
\begin{equation}\label{e.Psi-and-Phi}
\Upsilon_{nN_0+\ell}=\Phi_{n}+(\zc\circ\wt{\vf}_{n,\ell}-\zc)\circ\wt{\psi}_{n-1}\circ\cdots
    \circ\wt{\psi}_0.
\end{equation}

\begin{lem}\label{l.Phi-Decomp}
For every $N=nN_0+\ell$ and integers  $n\geq 1$ and $0\leq \ell\leq N_0-1$, \index{Random walks driven by renormalizations!and $S_N$}
\begin{align}
&\Upsilon_{N}=\ov{S}_{n}-\Delta\zc_0+\Delta\zc_{n,\ell} \label{e.Phi-Decomp}
\end{align}
where $\Delta\zc_0:=\zc-\wt{\zc}_0$, $\Delta\zc_{n,\ell}:=(\zc\circ \wt{\varphi}_{n,\ell}-\wt{\zc}_n)\circ\wt{\psi}_{n-1}\circ\cdots\circ\wt{\psi}_0$,  $\DS \ov{S}_n(\wt{\omega}):=\sum_{k=1}^n f_k(\ov{X}_k(\wt{\omega}))$,  
\begin{align}
\ov{X}_k(\wt{\omega})& := X^{(0)}_k(\pi(\tomega))=\text{unique element of $\mathfrak P_k$ containing $(\psi_{k-1}\circ\cdots\psi_0)(\pi(\wt{\omega}))$}.\label{e.X-bar}
\end{align} 
If $\tomega$ is sampled randomly uniformly from $\Rect_0$, then $(\ov{X}_k(\tomega))_{k\geq 0}$ is equal in distribution to $(X_k)_{k\geq 0}$, and $(\ov{S}_n(\tomega))_{n\geq 1}$ is equal in distribution to $(S_n)_{n\geq 1}$. The quantities 
$\Delta\zc_0$ and $\Delta\zc_{n,\ell}$ are uniformly bounded. 
\end{lem}
\begin{proof}
See \eqref{e.Psi-and-Phi}, and the proof of Proposition \ref{Prop-Markov-Representation-of-Phi}.
\end{proof}

Let $\ov{\mathcal F}_0^\infty$ denote the $\sigma$-algebra  generated by $\ov{X}_k:\Rect_0\to\mathfrak P_k$ $(k\geq 0)$. Let $\mathcal W_0^s(\Rect_0)$ denote the $\sigma$-algebra of Borel subsets $E\subset \Rect_0$ such that $\DS E=\bigcup_{\tomega\in E}W_0^s(\tomega)$, where 
$$
W^s_0(\tomega):=\pi|_{\mathsf{\Rect}_0}^{-1}\bigl(\text{maximal $E^s_0$-segment in the $\mathfrak P_0$-atom containing $\pi(\tomega)$}\bigr).
$$ 

\begin{lem}\label{l.sigma-algebra}
$\ov{\mathcal F}_0^\infty=\mathcal W^s_0(\Rect_0)$ mod $\mu_0$,  i.e. for every $A\in \ov{\mathcal F}_0^\infty$ there is a $B\in \mathcal W_0^s(\Rect_0)$ such that $\mu_0(A\triangle B)=0$, and vice versa. Moreover, if $\mu_{\Rect_0}:=\mu_0|_{\Rect_0}$,   then 
   for all $g\in L^1(\mu_{\Rect_0})$,  
\begin{equation}\label{e.conditional-expect-id-s}
\E_{\mu_{\Rect_0}}(g|\ov{\mathcal F}_0^\infty)(x)=\frac{1}{|W^s_0(x)|}\int_{W^s_0(x)}g\, ds\ \ \ (\text{integral with respect to arc length}).
\end{equation}
\end{lem}
\begin{proof}
If $\cF', \cF''$ are two $\sigma$-algebras s.t. for each $g\in L^1(\mu_{\Rect_0})$,
$\E(g|\cF')=\E(g|\cF'')$ almost everywhere then for each $A'\in \cF'$ there exists $A''\in \cF''$
such that $\mu_0(A'\triangle A'')=0$, e.g. 
$A''=\{\omega: F=1\}$, where $F$ is an $\cF''$-measurable a.e. version of $ \E(1_{A'}|\cF'')(\omega)$.

Thus, to prove the lemma, it is sufficient to verify \eqref{e.conditional-expect-id-s}. 

By Fubini's theorem, both sides of \eqref{e.conditional-expect-id-s} define bounded linear operators on $L^1$; therefore it is sufficient to prove   \eqref{e.conditional-expect-id-s} on a dense set in $L^1$. We will prove it for $g$ such that $g\circ\pi|_{\Rect_0}^{-1}$ is uniformly continuous on $\St_0$.

Given $\tomega\in\St$, let $C_n(\tomega):=\<\dot{a}_0,\ldots,a_n\>$, where $a_i:=\ov{X}_i(\tomega)$, then  $\ov{\mathcal F}_0^n$ is the discrete $\sigma$-algebra with atoms $\pi^{-1}|_{\Rect_0}(C_n(\tomega))$, and 
$$
\E(g|\ov{\mathcal F}_0^n)(\tomega)=\frac{1}{\mu_0(C_n(\tomega))}\int_{C_n(\tomega)}g\circ\pi|_{\Rect_0}^{-1} d\mu_0.
$$
$C_n(\tomega)$ are $(E^s_0,E^u_0)$-parallelograms with stable side of length $\ell^s(a_0)$, and unstable sides of length tending to zero as $n\to\infty$. These parallelograms decrease to an $E^u_0$-segment of length $\ell^s(a_0)$, containing $\pi(\tomega)$. By length considerations, this segment equals $\pi[W^s_0(\tomega)]$ up to endpoints. By the uniform continuity of $g\circ\pi|_{\Rect_0}^{-1}$, 
$$
\lim_{n\to\infty}\E(g|\ov{\mathcal F}_0^n)(\tomega)=\frac{1}{|\pi(W^s_0(\tomega))|}\int_{\pi(W^s_0(\tomega))}g\circ\pi|_{\Rect_0}^{-1}\, ds=\frac{1}{|W^s_0(\tomega)|}\int_{W^s_0(\tomega)}g\, ds.
$$
By the martingale convergence theorem, the left-hand-side is $\E_{\mu_{\Rect_0}}(g|\ov{\mathcal F}_0^\infty)$.
\end{proof}

Let   $\Xi_n(\tomega):=\zc[(\tpsi_0^{-1}\circ\cdots\circ\tpsi_{n-1}^{-1})(\tomega)]-\zc(\tomega).$ 
\begin{lem}\label{Lemma-Xi-and-S_n}
There exists a uniformly bounded sequence of functions  $\epsilon_n:\St^\ast\to\Z$ such that \index{Random walks driven by renormalizations!and $S_N$}
$
\DS\Xi_n-\epsilon_n=-\sum_{k=1}^{n}f_{k}\circ X^{(n)}_{k}\circ \pi
$, with $X^{(n)}_k$ as in \eqref{e.X-new}.
\end{lem}
\begin{proof}
Let $\epsilon_n:\St^\ast\to\Z$ be the function
\begin{equation}\label{e.epsilon-process}
\epsilon_n:=\zc\circ (\wt{\psi}_0^{-1}\circ\cdots\circ\wt{\psi}_{n-1}^{-1}) -\wt{\zc}_{0}\circ (\wt{\psi}_0^{-1}\circ\cdots\circ\wt{\psi}_{n-1}^{-1})+\wt{\zc}_n-\zc.
\end{equation}
By Lemma \ref{Lemma-Modified-Approx-Canonical}, $\epsilon_n(\cdot)$ is a uniformly bounded sequence of functions.  

Define $\tomega_n:=\tomega$, and $\tomega_k:=(\tpsi_k^{-1}\circ\cdots\circ\tpsi_{n-1}^{-1})(\tomega_n)$. Then, 
\begin{align*}
&\Xi_n(\tomega)-\eps_n(\tomega)=\Xi_n(\tomega_n)-\eps(\tomega_n)=\zc(\tomega_0)-\zc(\tomega_n)-\zc(\tomega_0)+\wt{\zc}_0(\tomega_0)-\wt{\zc}_n(\tomega_n)+\zc(\tomega_n)\\
&=\wt{\zc}_0(\tomega_0)-\wt{\zc}_n(\tomega_n)=-\sum_{k=0}^{n-1} \bigl[\wt{\zc}_{k+1}(\tomega_{k+1})-\wt{\zc}_{k}(\tomega_{k})\bigr]
=-\sum_{k=0}^{n-1}f^{(n)}_{k+1}(\pi(\tomega)), \text{ see \eqref{Eq-Modified-Frobenius}}.
\end{align*}
Now apply Lemma \ref{Lemma-f(X)}.
\end{proof}

 Note that
$\Xi_n^{W_n}(\tomega)-\eps_n(\tomega)$ is {\em not} equal in distribution to $-S_n$, because the random variable  $\Xi_n^{W_n}(\tomega)$ is defined by  sampling $\tomega$ randomly uniformly from the linear segment $W_n$, and not from the rectangle $\Rect_0$. We discuss this sampling procedure below.

\subsection{Conditioning on Pure Unstable Segments}

\label{ss.pure-segments}

Let $E$ be a one-dimensional linear subspace of $\mathbb{R}^2$. 
An {\em $E$-segment} $W$ is a linear segment (in 
$\mathbb{R}^2,\St_0$ or $\St$) parallel to $E$. The  {\em length} of $W$ shall be denoted by $|W|$. The {\em uniform distribution on $W$} is the measure \index{E$^u$-segment and E$^s$-segment @ $E^u$-segment and $E^s$-segment!pure segments!uniform distribution on}\index{Pure $E^u$- and $E^s$-segments!see {E$^u$-segment and E$^s$-segment}} 
$$
\Prob_W:=\frac{1}{|W|}\times\text{ the length measure on $W$.}
$$

Recall that $\wt{F}_n\subset \St$ is the fundamental domain for the action of the deck transformation $D:\St\to\St$, giving rise to the modified $\Z$-coordinate $\wt{\zc}_n$. In particular, $\wt{\zc}_n=\ell$ on $D^\ell(\wt{F}_n)$,  $\pi:D^\ell(\wt{F}_n)\to\St_0$ is bijective, and any cylinder $\<a_{n-k},\ldots,\dot{a}_n\>$ (being a subset of $a_n\in\mathfrak P_n$) has a lift to a parallelogram inside $D^\ell(\wt{F}_n)$.

\begin{defi}\label{d.pure-E-segments}
A {\em pure} $E^s_n$-segment with {\em depth}\index{Depth} $k$ and {\em height}\index{Height} $\ell$ is an $E^s_n$-segment inside $D^\ell(\wt{F}_n)\subset\St$, which projects up to endpoints to a maximal $E^s_n$-segment inside a cylinder $\<a_{n-k},\ldots,\dot{a}_n\>\subset\St_0$.
\end{defi}
\noindent 
 We denote the collection of  pure $E^s_n$-segment with depth $k$ by $\mathfs W_{k,n}$. 
By Lemmas  \ref{Lemma-Cylinder-Area} and  \ref{Lemma-Partition}(10),  
The lengths of  $W\in\mathfs W_{n,k}$ tends to zero  as the depth $k$ tends to $\infty$,   with the following uniform control: 
\begin{align}
\overline{c}_k&:=\sup\{|W|: W\in\bigcup_{n\geq k}\mathfs W_{k,n}\}<\infty,
\text{ and }\overline{c}_k\xrightarrow[k\to\infty]{}0\label{e.c-over},\\
\underline{c}_k&:=\inf\{|W|: W\in\bigcup_{n\geq k}\mathfs W_{k,n}\}>0.\label{e.c-under}
\end{align} 
It is also clear that the $\wt{\zc}_n$-coordinate  is constant along any pure $E^s_n$-segment. 

\begin{defi}\label{d.elementary}
Let $W$ be a pure $E^s_n$-segment. \index{Elementary events} \index{E$^u$-segment and E$^s$-segment @ $E^u$-segment and $E^s$-segment!pure segments!elementary events}
The elements of the $\sigma$-algebra on $W$ generated by  $X^{(n)}_i\circ\pi:W\to\Z$, $i=0,\ldots,n$ are called {\em elementary events} in $W$.  
\end{defi}

The elementary events are the finite disjoint unions of pure  $E^s_n$-segments $W'\subset W$ of depth $n$, and the same height as $W$.
For example,   
$\{\tomega\in W:\Xi_n(\tomega)-\eps(\tomega)=k\}$ is elementary (Lemma \ref{Lemma-Xi-and-S_n}), but $\{\tomega\in W:\Xi_n(\tomega)=k\}$ is not.

The following lemma gives a formula for $\Prob_W(C)$ when $W\in\mathfs W_{n,k}$, and $C$ is an elementary  event in $W$. For simplicity, we restrict to the case when $C$ is a single pure $E^s_n$-segment of depth $n$. If $C\neq\emptyset$, then $C$ and $W$ must have the same height $\ell$, and there is a cylinder $\<a_0,\ldots,\dot{a}_n\>$ such that $\pi(C)$ is a maximal $E^s_n$-segment in $\<a_0,\ldots,\dot{a}_n\>$, and $\pi(W)$  is a maximal $E^s_n$-segment in $\<a_{n-k},\ldots,\dot{a}_n\>$.

\begin{lem}\label{l.prob-of-elementry-event}
Under the above assumptions, \index{Elementary events} \index{E$^u$-segment and E$^s$-segment @ $E^u$-segment and $E^s$-segment!pure segments!elementary events}
\begin{align}\label{e.prob-of-elementary-event}
\Prob_W(C)&=\Prob[X_0=a_0,\ldots,X_{n}=a_{n}|X_{n-k}=a_{n-k},\ldots,X_n=a_n]\\
&=
\Prob[X_0=a_0,\ldots,X_{n-k-1}=a_{n-k-1}|X_{n-k}=a_{n-k}].\notag
\end{align}
\end{lem}

\begin{proof}
Set $Q:=\<a_{n-k},\ldots,\dot{a}_n\>$, and let $\wt{Q}$ denote the lift of $Q$ to $D^\ell(\wt{F}_n)$. 
 By Lemma~\ref{Lemma-Cylinder-Area}, these  are $(E^u_n,E^s_n)$-parallelograms.
 
By assumption, $W$ is a maximal $E^s_n$-segment in $\wt{Q}$. Let
$
{\vartheta}:\wt{Q}\to W
$ denote the {projection along the $E^u_n$-direction}, namely the function which maps $\tomega$ to the unique point in  $W^u(\tomega,\wt{Q})\cap W$, where  $W^u(\tomega,\wt{Q}):=$ the maximal $E^u_n$-segment  through $\tomega$ in $\wt{Q}$. 

%
%
%
%

\medskip
\noindent
{\em Claim 1.\/} Let $\Prob_W, \Prob_{\wt{Q}}$ denote the uniform distributions on $W$ and $\wt{Q}$, then $\Prob_{W}=\Prob_{\wt{Q}}\circ\vartheta^{-1}$. 

\smallskip
\noindent
{\em Proof of the Claim.\/} Let $\ell$ and $m$ denote the Euclidean length and area.  
Suppose $J$ is sub-segment of $W$, then $\vartheta^{-1}(J)$ is a parallelogram with $E^s_n$-side  $J$, and $E^u_n$-side with length $\ell^u(a_n)$. Thus  $m(\vartheta^{-1}(J))=\ell(J)\ell^u(a_n)|\sin\measuredangle(E^u_n,E^s_n)|$. Since $J$ was arbitrary,  $(m\circ\vartheta^{-1})|_{W}=\const.\ell|_{W}$. Normalizing to total mass one, we obtain $\Prob_{\wt{Q}}\circ\vartheta^{-1}=\Prob_{W}$.

\medskip
\noindent
{\em Claim 2.\/} Let  $\wt{P}$ denote the lift of $P:=\<a_0,\ldots,\dot{a}_n\>$ to $D^\ell(\wt{F}_n)$, then  $\vartheta^{-1}(C)=\wt{P}$.

\smallskip
\noindent
{\em Proof of the Claim.\/}    By Lemma \ref{Lemma-Cylinder-Area}, $\wt{P}$ is an $(E^u_n,E^s_n)$-parallelogram, and the $E^u_n$-sides of $\wt{P},\wt{Q}$ have the same length $\ell^u(a_n)$. Since $\wt{P}\subset\wt{Q}$, the maximal $E^u_n$-segments in $\wt{P}$ must  also be maximal $E^u_n$-segments in $\wt{Q}$. 
Thus $\vartheta^{-1}(C)$ equals the  union of all maximal $E^u_n$-segments  in $\wt{P}$ which intersect $C$.  
$C$ is a maximal $E^s_n$-segment of $\wt{P}$, because $\pi(C)$ is a maximal $E^s_n$-segment of $P$,  by assumption. So  
$
\vartheta^{-1}(C)=\wt{P}. 
$

\medskip
By Claims 1 and 2,  $\DS\Prob_W(C)=\Prob_{\wt{Q}}(\vartheta^{-1}(C))=\Prob_{\wt{Q}}(\wt{P})=\frac{\text{Area}(\wt{P})}{\text{Area}(\wt{Q})}=
\frac{\text{Area}(P)}{\text{Area}(Q)}$. By Proposition \ref{Prop-Markov-Kernel},   $\Prob_W(C)={\Prob[X_0=a_0,\ldots,X_n=a_n]}/{\Prob[X_0=a_{n-k},\ldots,X_n=a_n]}$. Since $(X_i)_{i\geq 0}$ is a Markov chain, 
$\Prob_W(C)=\Prob
[X_0=a_0,\ldots,X_{n-k-1}=a_{n-k-1}|X_{n-k}=a_{n-k}]$.
\end{proof}

\begin{cor}
Suppose $W\in\mathfs W_{n,k}$. If $W$ projects to a maximal $E^s_n$-segment in $\<a_{n-k},\ldots,\dot{a}_n\>$,  then for every $n>k$, \index{E$^u$-segment and E$^s$-segment @ $E^u$-segment and $E^s$-segment!pure segments!conditioning on}
\begin{equation}\label{e.Xi-W-local-prob-formula}
\Prob(\Xi_n^W-\eps_n=t)=\Prob\left(S_{n-k-1}=-t-\sum_{j=n-k}^n f_j(a_j)\,\bigg|\,X_{n-k}=a_{n-k}\right).
\end{equation}
\end{cor}

Non-pure $E^s_n$-segments  can be approximated from inside and outside by pure $E^s_n$-segments, as in the following lemma:

\begin{lem}\label{l.approx-by-pure-segments}
For every $L>1$ and $\delta>0$ there are positive integers $k,K>1$ as follows. Suppose  $W_n$ is a sequence of $E^s_n$-segments in $\St$ with lengths in $[L^{-1},L]$.\index{E$^u$-segment and E$^s$-segment @ $E^u$-segment and $E^s$-segment!approximation by pure segments} Then for all $n>K$, there are pure $E^s_n$ segments of depth $k$ $W_{n,1},\ldots,W_{n,m_n+2}$ such that:
\begin{enumerate}[(1)]
\item $W_{n,i}$ have pairwise disjoint relative interiors;
\item $\DS\bigcup_{i=1}^{m_n}W_{n,i}\subset W_n\subset \bigcup_{i=1}^{m_n+2}W_{n,i}$\text{ and each $\overline{W}_{n,i}$ intersects $\ov{W}_n$};
\item $\DS \frac{1}{|W_n|}\sum_{i=1}^{m_n}|W_{n,i}|\Prob_{W_{n,i}}(\cdot)\leq 
\Prob_{W_n}(\cdot)\leq \frac{1}{|W_n|}\sum_{i=1}^{m_n+2}|W_{n,i}|\Prob_{W_{n,i}}(\cdot)
$;
\item $\DS 1-\delta\leq \frac{1}{|W_n|}\sum_{i=1}^{m_n}|W_{n,i}|\leq 1\leq \frac{1}{|W_n|}\sum_{i=1}^{m_n+2}|W_{n,i}|\leq 1+\delta
$;
\item $m_n\leq K$.
\end{enumerate}
\end{lem}
\noindent
The proof follows easily from  \eqref{e.c-over} and \eqref{e.c-under}, and we omit it.

\part{Proofs of the Main Results.}

This part contains the proofs of the main results stated in Part I. To locate the proof of a specific theorem,  see the proof locator on page \pageref{page-proof-locator}.

We assume throughout that $\alpha$ is an irrational of bounded type and $\beta=2\alpha-1$. There is no loss of generality in assuming that $\alpha\in (0,1)$ and that $\beta\in (-1,1)$ has even continued fraction expansion 
$$
\beta=[0;2m_1,2m_2,\ldots]\ \ (m_i\in\Z\setminus\{0\}).
$$
Indeed,  $D_N(\alpha,I)\!\!\equiv\!\! D_N(\{\alpha\},I)$,  
$T_\alpha=T_{\{\alpha\}}$, and the linear flow $\phi_{\vec{v}}$ with velocity 
$\vec{v}={\beta'\choose 1}$,  $\beta'=[2m_0;2m_1,\ldots]$ is conjugated  to $\phi_{\vec{u}}$ with $\vec{u}={\beta\choose 1}$, $\beta:=[0;2m_1,2m_2,\ldots]\in (-1,1)$ by the $\Rect_0$-preserving homogeneous automorphism with derivative ${\tiny{ \biggl(\begin{array}{cc}1 & -2m_0 \\0& 1\end{array}}\biggr)}$.  

Throughout Part IV, $S_n=f_1(X_1)+\cdots+f_N(X_n)$, $V_n=\mathrm{Var}(S_n)$,  $(X_i)_{i\geq 0}$ is the Markov chain from Proposition \ref{Prop-Markov-Kernel}, and $f_i$ are the symbolic jump functions from Lemma \ref{Lemma-f(X)}.
The normalized length measure on a linear segment $W$ is denoted by $\Prob_W$.

\section{Temporal Limit Theorems}
\label{ScCircle}

\subsection{Temporal  Central Limit  Theorems (Theorems \ref{t.A_n-B_n-estimates}, \ref{t.A_n} and \ref{ThStairTLT})}
\label{ScProofBeck}

\begin{lem}\label{l.trivial-prob-fact}
Suppose $Y_n,Y_n'$ are random variables on the same probability spaces $\mathfrak S_n$ such that  $\|Y_n-Y_n'\|_\infty\to 0$ or $\|Y_n/Y_n'-1\|_\infty\to 0$. If $Y_n$ converges in distribution to a limit without atoms, then $Y_n'$ converges in distribution to the same limit. \index{Convergence in distribution}
\end{lem}
\begin{proof}
For every $\epsilon>0$ and $t\neq 0$,  $\Prob[Y_n'<t]$ is eventually between $\Prob[Y_n<t\pm\eps]$.
\end{proof}

Recall the backward random walk $\Xi_n^W$ from Section \ref{s.RW-Driven-By-Renormalizations}.

\begin{lem}\label{Lemma-Pure-Seg-CLT}Let $W_n$ be a sequence of pure $E^s_n$-segments  such that $|W_n|$ is bounded below.
Then for all $a<b$,
$$
\Prob\left[{\Xi_n^{W_n}}/{\sqrt{V_n}}\in (a,b)\right]\xrightarrow[n\to\infty]{}\frac{1}{\sqrt{2\pi}}\int_a^b e^{-t^2/2} dt.
$$
\end{lem}
\begin{proof} 
By \eqref{e.c-over}, the depths of $W_n$ must be bounded above,  so  we can split $\{W_n\}_{n\geq 1}$ to subsequences with eventually  constant depth. Without loss of generality, all $W_n$ have the same  depth $k$. In particular, each $W_n$ projects to a maximal $E^s_n$-segment in some cylinder $\<a_{n-k}^{(n)},\ldots,\dot{a}_n^{(n)}\>$. 
Let  $\DS s_{n,k}=\sum_{j=n-k}^n f_j(a_j^{(n)})=O(1)$. By \eqref{e.Xi-W-local-prob-formula}, 
\begin{equation}\label{e.do}
\Prob_{W_n}\left[\frac{\Xi_n-\eps_n}{\sqrt{V_{n}}}\in (a,b)\right]=
\Prob\left(\frac{S_{n-k-1}+ s_{n,k}}{\sqrt{V_{n}}}\in (-b,-a)\, \bigg|\, X_{n-k}=a_{n-k}^{(n)}
\right).
\end{equation}
 Here $\eps_n$ is the uniformly bounded sequence of functions from Lemma \ref{Lemma-Xi-and-S_n} 

The conditioning is on the last coordinate, not the first coordinate. To deal with this, we reverse time and work with the array 
$$
\begin{array}{lllll}
X_0 &&&&
\text{conditioned on }X_0=a_0^{(k)}\\
X_1 & X_0 &&& \text{conditioned on }X_1=a_1^{(k+1)}\\
 \ \vdots & \ \vdots & \ddots & &\\
 X_{n-k}& X_{n-k-1} & \cdots & X_0 & \text{conditioned on }X_{n-k}=a_{n-k}^{(n)}.
\end{array}
$$
Each row is a finite Markov chain, and one 
can prove exactly as in Proposition \ref{Prop-Doeblin} that all rows are uniformly elliptic with the same ellipticity constant. 

Using the exponential mixing of uniformly elliptic Markov arrays (Proposition \ref{p.mixing-for-arrays}), it is possible to show that
\begin{align*}
&\E'(S_{n-k-1}):= \E(S_{n-k-1}|X_{n-k}=a_{n-k}^{(n)})=\E(S_{n-k-1})+O(1);\\
&V_{n-k-1}':=\Var(S_{n-k-1}|X_{n-k}=a_{n-k}^{(n)})\sim V_{n-k-1};\\
&V_{n-k-1}\sim V_n.
\end{align*}
By Propositions \ref{Prop-Drift} and \ref{Prop-Var-Growth},  $V_{n-k-1}'\to\infty$ and $\E(S_{n-k-1}|X_{n-k}=a_{n-k}^{(n)})=O(1)$.
 
We now apply  Dobrushin's central limit theorem, in its version  for  uniformly elliptic  Markov arrays (Theorem \ref{ThMCCLT}). By Lemma~\ref{l.trivial-prob-fact}, we can replace $A_n$ by minus $\E'(S_{n-k-1})$ and $V_n$ by $V_{n-k-1}'$, and this leads to
 $$
\Prob\left(\frac{S_{n-k}+s_{n,k}}{\sqrt{V_{n}}}\in (-b,-a)\, \bigg|\, X_{n-k}=a_{n-k}^{(n)}
\right)\xrightarrow[n\to\infty]{}\frac{1}{\sqrt{2\pi}}\int_a^b e^{-t^2/2}dt.
$$
By  \eqref{e.do},  $\frac{\Xi_n^{W_n}-\eps_n}{\sqrt{V_n}}$ converges in distribution to the standard normal. Since $s_{n,k}$ and $\eps_n$ are uniformly bounded and $V_n\!\to\!\infty$, the lemma follows from  Lemma \ref{l.trivial-prob-fact}.
\end{proof}

\begin{proof}[\bf Proof of the Temporal CLT for the Linear Flow (Theorem \ref{ThStairTLT})]  Fix 
$\wt{\omega}\in \St$.
Suppose $T\gg 1$, and let 
$
A_T:=\{\phi^t(\wt{\omega}):0<t<T\}.
$

By Proposition \ref{Prop-Contracting-Automorphisms}, 
$
(\tpsi_{n-1}\circ\cdots\circ\tpsi_0)(A_T)
$ is an $E^s_n$-segment of length $\lambda_{n-1}^s\cdots\lambda_0^s T$, where $\lambda^s_i:=\|d\tpsi_i|_{E^s_i}\|\leq \lambda^{-1}<1$.  Since there are only finitely many different $\wt{\psi}_i$, $\lambda^s_i>L^{-1}$ with $L>1$ independent of $i$. So for every $T>1$, there are $n$ such that 
\begin{equation}
\label{DefBrN}
A_T^\ast:=(\wt{\psi}_{n-1}\circ\cdots\circ\wt{\psi}_0)(A_T)\text{ has length in }[1,L].
\end{equation}
 Let $n(T)$ be the minimal number satisfying this condition.
Necessarily,  $n/\log T$ is bounded away from zero and infinity for all $T$ large enough. Let 
$$
\sigma_T:=\sqrt{V_n},\text{ where }V_n:=\mathrm{Var}(S_n), \quad
a_T(\wt{\omega}):=z_{nN_0}(\tomega), \text{ see \eqref{e.z-n} and \eqref{e.mat-prod-form}}.  
$$
$A_T^\ast$ begins inside horizontal rectangle $\Rect_{a_T(\tomega)}$, and  $A_T^\#:=D^{-a_T}(A_T^\ast)$ begins inside $\Rect_0$. Since $A_T^\#$ has bounded length, its $\zc$-coordinate is uniformly bounded on $A_T^\#$. 

By  Proposition \ref{Prop-Var-Growth}, $V_n\asymp n$, therefore $\sigma_T\asymp\sqrt{\log T}$. By Lemmas \ref{Lemma-Modified-Approx-Canonical} and  \ref{Lemma-f(X)},  
$\DS  |a_T(\wt{\omega})|\leq n\sup_k\|f_k\|_\infty+O(1)\leq C\log T$, with  $C$ independent of $\wt{\omega}$.
   

$D^{-a_T}\circ\wt{\psi}_{n-1}\circ\cdots\circ\wt{\psi}_0:A_T\to A_T^\#$ is an affine map, which carries the normalized length measure $\Prob_{A_T}$ on  $A_T$,  to the normalized length measure $\Prob_{A_T^\#}$  on $A_T^\#$. So 
\begin{align*}
&\Prob_{A_T}\left[\frac{\zc-a_T}{\sigma_T}\in (a,b)\right]=
\Prob_{A_T^\#}\left[\frac{\zc\circ\wt{\psi}_{0}^{-1}\circ\cdots\circ\wt{\psi}_{n-1}^{-1}\circ D^{a_T}-a_T}{\sigma_T}\in (a,b)\right]\\
&\overset{!}{=}\Prob_{A_T^\#}\left[\frac{\zc\circ\wt{\psi}_{0}^{-1}\circ\cdots\circ\wt{\psi}_{n-1}^{-1}}{\sigma_T}\in (a,b)\right] =\Prob_{A_T^\#}\left[\frac{\Xi_n+\zc}{\sqrt{V_n}}\in (a,b)\right].
\end{align*}
(The marked identity holds because  $\wt{\psi}_i\circ D=D\circ\wt{\psi}_i$ and $\zc\circ D=\zc+1$).

Fix $\delta>0$ arbitrarily small, and choose pure $E^s_n$-segments $W_{n,1},\ldots,W_{n,m_n+2}$ with depth $k$ as in Lemma \ref{l.approx-by-pure-segments}. Since $W_{n,i}$ have fixed depth, $|W_{n,i}|$ is  bounded away from zero and infinity. Since  $\ov{W}_{n,i}\cap\ov{A}_T^\#\neq \emptyset$ and $\zc$ is uniformly bounded on $A_T^\#$, $\zc$ is uniformly bounded on $W_{n,i}$. Thus by Lemmas \ref{l.trivial-prob-fact} and \ref{Lemma-Pure-Seg-CLT}, 
$$
\lim_{n\to\infty}\Prob_{W_{n,i}}\left[\frac{\Xi_n+\zc}{\sqrt{V_n}}\in (a,b)\right]=\lim_{n\to\infty}\Prob_{W_{n,i}}\left[\frac{\Xi_n}{\sqrt{V_n}}\in (a,b)\right]=\frac{1}{\sqrt{2\pi}}\int_a^b e^{-t^2/2}dt.
$$
Invoking parts (3) and (4) of Lemma \ref{l.approx-by-pure-segments}, we obtain that for all $T$ large enough, 
$$
\frac{1-2\delta}{\sqrt{2\pi}}\int_a^b e^{-t^2/2}dt\leq 
\Prob_{A_T^\#}\left[\frac{\Xi_n+\zc}{\sqrt{V_n}}\in (a,b)\right]\leq  \frac{1+2\delta}{\sqrt{2\pi}}\int_a^b e^{-t^2/2}dt.
$$
Thus  $\Prob_{A_T}\left[\frac{\zc-a_T}{\sigma_T}\in (a,b)\right]=\Prob_{A_T^\#}\left[\frac{\Xi_n+\zc}{\sqrt{V_n}}\in (a,b)\right]\xrightarrow[T\to\infty]{}\frac{1}{\sqrt{2\pi}}\int_a^b e^{-t^2/2}dt$.
\end{proof}

\medskip
\noindent
{\bf Proof of the Temporal CLT for $\boldsymbol{D_N(\alpha)}$ (Theorems \ref{t.A_n-B_n-estimates} and \ref{t.A_n}).}
We  use the notation of the previous proof. 
Let $\tomega_0$ denote the singularity at the middle of the top side of $\Rect_0$,  let $T_n:=n$, and set 
$
 A_n:= \frac{1}{2} a_{T_n}(\tomega_0)\ , \ B_n:= \frac{1}{2}\sigma_{T_n}.
$

Recall that $D_N(\alpha)$ takes values in $\frac{1}{2}\Z$.  By \eqref{e.Lin-Flow-1}, 
\begin{align*}
&\frac{1}{n} \#\{1\leq N\leq n: \frac{D_N(\alpha)-A_n}{B_n}\in (a,b)\}=\sum_{\xi\in \frac{1}{2}\Z,\frac{\xi-A_N}{B_N}\in (a,b) }\!\!\!\!\!\!\!\!\!(1/n) \#\{1\leq N\leq n: D_N(\alpha)=\xi\}\\
&=\sum_{\xi\in\frac{1}{2}\Z,\frac{\xi-A_N}{B_N}\in (a,b) }\frac{1}{T_n}\int_{-1}^{T_n-1}{ 1_{\{2\xi\}}}(\zc(\phi^t(\tomega_0)))dt
=\frac{1}{T_n}\int_{-1}^{T_n-1}1_{(a,b)}\left(\frac{\zc(\phi^t(\tomega_0))/2-A_n}{B_n}\right)dt\\
&=\Prob_{A_{T_n}(\tomega_0)}\left[\frac{\zc-2A_n}{2B_n}\in (a,b)\right]+O(T_n^{-1})
=\Prob_{A_{T_n}(\tomega_0)}\left[\frac{\zc-a_{T_n}}{\sigma_{T_n}}\in (a,b)\right]+o(1).
\end{align*}
Thus,  by Theorem \ref{ThStairTLT}, $\frac{1}{n} \#\{1\leq N\leq n: \frac{D_N(\alpha)-A_n}{B_n}\in (a,b)\}\xrightarrow[n\to\infty]{}\frac{1}{\sqrt{2\pi}}\int_a^b e^{-t^2/2}dt$.  
 Since $|a_T|\leq C\log T$ and $B_T\asymp\sqrt{\log T}$, $A_n=O(\log n)$ and $B_n\asymp\sqrt{\log n}$, giving us Theorem \ref{t.A_n-B_n-estimates}. 
 
 Next,  $A_n=\frac{1}{2}a_{T_n}(\tomega_0)=\frac{1}{2} z_{k_n N_0}(\tomega_0)$, where $k_n=n(T_n)$ is the minimal integer such that $\lambda^s_0\cdots\lambda^s_{k_n-1} T_n\in [1,L]$, see \eqref{DefBrN}.  Thus,  $k_n\asymp \log n$, $k_n$ is non-decreasing, and  $k_{n+1}-k_n$ is bounded. 
Let $\nu_n:=k_nN_0$, then   Theorem \ref{ThDriftOfOrigin} (proved in Appendix \ref{Appendix-Twist}) says that  
$
\DS A_n=\frac{1}{2} z_{\nu_n}(\tomega_0)=\frac{1}{2}\zc_{k_n}
=\frac{1}{4}\sum_{j=1}^{\nu_n}[\sgn(m_{2j})-\sgn(m_{2j-1})].
$
Theorem \ref{t.A_n} follows.\qed

\subsection{Temporal Local  Limit Theorems (Theorems \ref{t.LLT-for-D_N} and \ref{t.Local-TCLT}) }
\label{SSTempLLT}

\begin{lem}\label{Lemma-Pure-Seg-LLT}
 Let $W_n$ be a sequence of pure $E^s_n$-segments  with heights $h_n$, and length bounded below.   Fix $C:=\<\dot{a}_0,\ldots,a_{k-1}\>\neq \emptyset$.  If $k_n\in\Z, t\in\R$, and  $\frac{k_n-h_n}{\sqrt{V_n}}\xrightarrow[n\to\infty]{}t$, then 
$$
\Prob_{W_n}\left\{\tomega:\begin{array}{c}
({\psi}_0^{-1}\circ\cdots\circ{\psi}_{n-1}^{-1})(\pi(\tomega))\in C\\
\wt{\zc}_0[(\wt{\psi}_0^{-1}\circ\cdots\circ\wt{\psi}_{n-1}^{-1})(\tomega)]=k_n
 \end{array}
\right\}=[1+o(1)]\frac{e^{-t^2/2}\mu_0(C)}{\sqrt{2\pi V_n}},\text{ as }n\to\infty.
$$
\end{lem}
\begin{proof}
As in the proof of Lemma \ref{Lemma-Pure-Seg-CLT}, it is sufficient to consider the case when the depth of  $W_n$ is constant, equal to  $d$. 
Then each $W_n$ projects to a maximal $E^s_n$-segment in some cylinder $\<b_{n-d}^{(n)},\ldots,\dot{b}_n\>$. 

Let 
$E_n:=\{\tomega\in W_n: ({\psi}_0^{-1}\circ\cdots\circ{\psi}_{n-1}^{-1})(\pi(\tomega))\in C, \wt{\zc}_0[(\wt{\psi}_0^{-1}\circ\cdots\circ\wt{\psi}_{n-1}^{-1})(\tomega)]=k_n
\}.$
By assumption, $\wt{\zc}_n=h_n$ on $W_n$. By the definition of the jump functions and  \eqref{e.X-new}, on $W_n$, 
$$
\sum_{j=1}^{n}f_j[X^{(n)}_j(\pi(\tomega))]=\wt{\zc}_n(\tomega)-\wt{\zc}_0[(\wt{\psi}_0^{-1}\circ\cdots\circ\wt{\psi}_{n-1}^{-1})(\tomega)]=h_n-\wt{\zc}_0[(\wt{\psi}_0^{-1}\circ\cdots\circ\wt{\psi}_{n-1}^{-1})(\tomega).
$$
It follows that 
$\DS
E_n=\left\{\tomega\in W_n:\begin{array}{l}
X^{(n)}_i(\pi(\tomega))=a_i\text{ for }i=0,\ldots,k-1\\
\DS \sum_{j=1}^{n}f_j[X^{(n)}_j(\pi(\tomega))]=h_n-k_n
\end{array} \right\}.
$
We conclude that  $E_n$ is an elementary event in the sense of Definition \ref{d.elementary}. Thus, by Lemma \ref{l.prob-of-elementry-event} and the Markov property,
\begin{align}
&\Prob_{W_n}(E_n)=\Prob\left[S_n=h_n-k_n,\, X_1=a_1,\ldots,X_{k-1}=a_{k-1}\, \big|\, X_{n-d}=b_{n-d}^{(n)},\ldots,X_n=b_n^{(n)}\right]\notag\\
&=\Prob[X_1=a_1,\ldots,X_{k-1}=a_{k-1}]\times\notag\\
&\hspace{0.7cm}\Prob\left[\sum_{i=k}^{n-d-1} f_i(X_i)=h_n-k_n-\sum_{i=1}^{k-1}f_i(a_i)-\sum_{i=n-d}^n f_i(b_i^{(n)})\, \bigg|\, X_{k-1}=a_{k-1}, X_{n-d}=b_{n-d}^{(n)}\right]\notag\\
&=\mu_0(C)\Prob[S_{n'}'-\E(S_{n'}')=\delta_{n'}|X_0'=a_{k-1}, X_{n'}'=b_{n-d}^{(n)}],\label{e.Prob[S-n-prime]}
\end{align}
 where   $n'=n-d-k$,  $X_i'=X_{i+k-1}$,  $f_j'=f_{j+k-1}$, 
 $\DS S_{n'}'=\sum_{j=1}^{n'}f_j'(X_j')$, and 
$$\delta_{n'}:=h_n-k_n-\sum_{i=1}^{k-1}f_i(a_i)-\sum_{i=n-d}^n f_i(b_i^{(n)})-\left(\E(S_{n-d-1})-\E(S_{k-1})\right)=h_n-k_n+O(1).$$

A calculation using  Corollary  \ref{Cor-rho-mixing} shows that  $V_{n'}':=\Var(S_{n'}')\sim V_n$. 
Thus  
$\frac{\delta_{n'}}{\sqrt{V_{n'}'}}\to -t$. 

Now we apply that the  mixing local limit theorem (Theorem \ref{ThLLTArrays}) to the Markov chain 
$(X_i)_{i\geq k-1}$ 
{started from $X_{k-1}=a_{k-1}$ and conditioned on $X_{n-d}=b_{n-d}.$
The conditioned chain has transition probabilities
\begin{equation}
\label{CondChain}
 \pi_{i, i+1}^{b_{n-d}}(a_i, a_{i+1})=\pi_{i, i+1}(a_i, a_{i+1}) 
\frac{\Prob(X_{n-d}=b_{n-d}|X_{i+1}=a_{i+1})}{\Prob(X_{n-d}=b_{n-d}|X_{i}=a_{i})}
\end{equation}
The conditions of Theorem \ref{ThLLTArrays} are all satisfied.
First, $f_i$ are uniformly bounded by Lemma \ref{Lemma-f(X)}.
Next, uniform ellipticity, irreducibility and the hereditary property for the unconditioned chain 
hold  by Propositions \ref{Prop-Doeblin}, \ref{Prop-Irred} and Lemma \ref{l.hereditary}.
The conditioned chain enjoys the same property because the second factor in \eqref{CondChain}
is uniformly bounded from above and away from zero due to Proposition  \ref{Prop-Markov-Kernel}(4) and 
the uniform ellipticity of the unconditioned chain
(or by $\phi$-mixing given by Corollary \ref{Cor-phi-mixing}). 
Theorem \ref{ThLLTArrays} gives the asymptotic behavior of \eqref{e.Prob[S-n-prime]}, and  the lemma follows.
}
\end{proof}

Recall the definition of the polygon $\wt{F}_0\subset\St$ in \S \ref{Section-Modified-Z-coordinates}. $\wt{F}_0$ is a fundamental domain for the action of the covering group, $\wt{\zc}_0=k$ on $D^k(\wt{F}_0)$, and $\pi:\wt{F}_0\to \St_0$ is continuous and injective. For every non-empty cylinder  $C=\<\dot{a}_0,\ldots,a_{k-1}\>$ and $\ell\in\Z$, let 
\begin{equation}\label{e.curly-C-k}
\begin{aligned}
\wt{C}_\ell&:=\text{the lift of $C$ to $D^\ell(\wt{F}_0)$}=D^\ell[\pi|_{\wt{F}_0}^{-1}(C)];\\
\mathfs C_k&:=\{\wt{C}_\ell: C=\<\dot{a}_0,\ldots,a_{k-1}\>,\ \ell\in\Z,\ \mathrm{int}(\wt{C}_\ell)\cap\mathrm{int}(\Rect_0)\neq \emptyset\}.
\end{aligned}
\end{equation}
By Lemma \ref{Lemma-Modified-Approx-Canonical},  $L_0:=\sup|\zc-\wt{\zc}_0|<\infty$. So $\bigcup\mathfs C_k\subset \bigcup_{|\ell|\leq L_0}\Rect_\ell$, and $\mu_0(\bigcup\mathfs C_k)$ is  bounded.

\begin{lem}\label{Lemma-Constant-Ratio}
There exist $\delta_k\to 0$ as $k\to\infty$
such that for each $\wt{C}_\ell\in\mathfs C_k$ and for every maximal $E^s_0$-segment $W\subset \mathrm{int}(\wt{C}_\ell)$,  
\begin{equation}\label{e.constant-ratio}
\left|\frac{|W\cap{\Rect}_0|}{|W|}-\frac{\mu_0(\wt{C}_\ell\cap{\Rect}_0)}{\mu_0(\wt{C}_\ell)}
\right|<\delta_k.
\end{equation}
\end{lem}
\begin{proof}
Suppose $\wt{C}_\ell\in \mathfs C_k$, where $C=\<\dot{a}_0,\ldots,a_{k-1}\>$. 
By Lemma \ref{Lemma-Cylinder-Area}, $C$ and  $\wt{C}_\ell$ are rectangles  with sides parallel to $E^s_0$ and $E^u_0$; The $E^s_0$-side has length   $\ell^s(a_0)\geq c^{-1}>0$; and the $E^u_0$-side has length $\ell^u(a_{k-1})/\prod_{i=0}^{k-2}\lambda_i^u\leq c \lambda^{-k}$.  Here $c>1$ is some  constant independent of   $k$,$\ell$,  $a_1,\ldots,a_{k-1}$, and coming from Lemma \ref{Lemma-Partition}(10).

Since $E^s_0\parallel {\beta\choose 1}$, the angle between $E^s_0$ and the vertical side of $\Rect_0$ is $\theta:=\tan^{-1}(\beta)$.
By Euclidean geometry, for any maximal $E^s_0$-segments $W,W'$ in $\wt{C}_\ell$ parallel to $E^s_0$, 
$$
\big| |W\cap\Rect_0|-|W'\cap\Rect_0|\big|\leq 2c\lambda^{-k}\tan\theta\text{ and }|W|=|W'|=\ell^s(\wt{C}_\ell)\geq c^{-1}.
$$
So $\DS \left|
\frac{|W\cap\Rect_0|}{|W|}-\frac{|W'\cap\Rect_0|}{\ell^s(\wt{C}_\ell)}
\right|<2c^2\beta\lambda^{-k}.
$ Keeping $W$ fixed and integrating over $W'$ intersecting  the $E^u_0$-side, we obtain
$
\left|\frac{|W|\cap\Rect_0|}{|W|}\ell^u(\wt{C}_\ell)-\frac{\text{Area}(\wt{C}_\ell\cap\Rect_0)}{\ell^s(\wt{C}_\ell)}
\right|<2c^2\beta\lambda^{-k}\ell^u(\wt{C}_\ell).
$
Dividing by $\ell^u(\wt{C}_\ell)$, we obtain the lemma (since $E^u_0\perp E^s_0$).
\end{proof}

\noindent
{\bf Proof of the Temporal Local CLT for Linear Flows (Theorem \ref{t.Local-TCLT}).}
Given integers $k_T$,  let  
$
\mathfrak m_T(\wt{\omega},k_T):=\frac{1}{T}\mathrm{mes}\left\{t\in (0,T): \phi^t(\wt{\omega})\in{\Rect}_{k_T}\right\}
$, where $\text{mes}$ is Lebesgue's measure.  {We need to show that}
  if $\frac{k_T-a_T(\tomega)}{\sigma_T}\to s$, then 
  $
  \sqrt{2\pi\sigma_T^2}\times\mathfrak{m}_T(\wt{\omega},k_T)\to e^{-\frac{1}{2}s^2}.
  $
   
We use the notation of the proof of the temporal CLT in \S\ref{ScProofBeck}.
Write $A_T=A_T'\uplus E$, where $A_T'$ projects to a  union of maximal $E^s_0$-segments in $\mathfrak P_0$-elements, and $E$ (an edge effect) is contained in a union of at most two such pieces. Then 
$$
\mathfrak m_T(\wt{\omega}_0,k_T)=\frac{|A_T'\cap {\Rect}_{k_T}|}{|A_T|}+O(T^{-1}).
$$

Fix  $T$ and $k$ arbitrarily large, and consider the collection  $\{D^{k_T}(\wt{C}_\ell):\wt{C}_\ell\in\mathfs C_k, \ell\in\Z\}$. This is a pairwise disjoint cover of ${\Rect}_{k_T}$, therefore, 
\begin{align*}
\mathfrak m_T(\wt{\omega}_0,k_T)&=
\sum_{\wt{C}_\ell\in\mathfs C_k}\frac{|A_T'\cap D^{k_T}(\wt{C}_\ell)\cap {\Rect}_{k_T}| }{|A_T|}+O(T^{-1}).
\end{align*}
Write $\DS A_T'\cap D^{k_T}(\wt{C}_\ell)=\biguplus_i D^{k_T}(W_i)$, where $W_i$ are maximal $E^s_0$-segments in  $\wt{C}_\ell$, then 
\begin{align*}
&|A_T'\cap D^{k_t}(\wt{C}_\ell)\cap {\Rect}_{k_T}|=\sum_i |D^{k_T}(W_i\cap{\Rect}_0)|\\
&=\sum_i \frac{|D^{k_T}(W_i\cap{\Rect}_0)|}{|D^{k_T}(W_i)|}\times |D^{k_T}(W_i)|=\sum_i \frac{|W_i\cap{\Rect}_0|}{|W_i|}\times |D^{k_T}(W_i)|\\
&=\left(\frac{\mu_0(\wt{C}_\ell\cap{\Rect}_0)}{\mu_0(\wt{C}_\ell)}\pm\delta_k\right)
|A_T'\cap D^{k_T}(\wt{C}_\ell)|, \text{ see Lemma \ref{Lemma-Constant-Ratio}.}
\end{align*}

Dividing by $|A_T|$, and summing over $\wt{C}_\ell\in\mathfs C_k$,  we obtain
\begin{equation}\label{e.8.6}
\mathfrak m_T(\wt{\omega}_0,k_T)=\frac{|A_T'|}{|A_T|}\sum_{\wt{C}_\ell\in\mathfs C_k}
\left(\frac{\mu_0(\wt{C}_\ell\cap{\Rect}_0)}{\mu_0(\wt{C}_\ell)}\pm\delta_k\right)
\left(\frac{|A_T'\cap D^{k_T}(\wt{C}_\ell)|}{|A_T'|}\right)+O(T^{-1}).
\end{equation}
The term $O(T^{-1})$ equals $|E\cap{\Rect}_{k_T}|/|A_T|$, and is uniform in $k$ and $k_T$.

Fix $C=\<\dot{a}_0,\ldots,a_{k-1}\>$ such that $\wt{C}_\ell\in\mathfs C_k$.
In \S \ref{ScProofBeck}, we constructed 
  an integer\\ $n=n(|A_T'|)\asymp \log |A_T'|\asymp\log T$ such that 
$
A_T^\ast:=(\wt{\psi}_{n-1}\circ\cdots\circ\wt{\psi}_0)(A_T')\text{ has length in }[1,L],
$
where $L$ is a global constant. Then
\begin{align*}
\frac{|A_T'\cap D^{k_T}(\wt{C}_\ell)|}{|A_T'|}&=\Prob_{A_T'}\{\wt{\omega}:\wt{\zc}_0(\wt{\omega})=k_T+\ell,\ 
\pi(\wt{\omega})\in \<\dot{a}_0,\ldots,a_{k-1}\>\}\\
&=
\Prob_{A_T^\ast}\left[
\begin{array}{c}
\wt{\zc}_0\circ \wt{\psi}_{0}^{-1}\circ\cdots\circ\wt{\psi}_{n-1}^{-1}=k_T+\ell;\\
 X^{(n)}_i\circ\pi=a_i\ (i=0,\ldots,k-1)
 \end{array} \right].
\end{align*}

We calculate this probability by first approximating $A_T^\ast$ by pure $E^s_n$-segments $W_i$ of arbitrarily large (but constant) depth $d$ as in Lemma \ref{l.approx-by-pure-segments}, and then applying Lemma \ref{Lemma-Pure-Seg-LLT}.  
Since all points in $A_T^\ast$ have $\zc$-coordinate $a_T(\tomega)+O(L)$, each $W_i$ has height $a_T(\tomega)+O(L)$, and  we obtain the following consequence: If $\frac{k_T-a_T}{\sigma_T}\xrightarrow[T\to\infty]{}s$, then 
$$
\frac{|A_T'\cap D^{k_T}(\wt{C}_\ell)|}{|A_T'|}=[1+o(1)]\frac{e^{-\frac{1}{2}s^2}\mu_0(\wt{C}_\ell)}{\sqrt{2\pi\sigma_T^2}}.
$$
Substituting this in \eqref{e.8.6}, and recalling that $|{\mathfs C}_k|<\infty$, we see that 
\begin{align*}
&\sqrt{2\pi\sigma_T^2}\times\mathfrak{m}_T(\wt{\omega}_0,k_T)
=(e^{-\frac{1}{2}s^2}+o(1))\sum_{ \wt{C}_\ell\in\mathfs C_k}[\mu_0(\wt{C}_\ell\cap{\Rect}_0)\pm\delta_k \mu_0(\wt{C}_\ell)]+O\left(\frac{\sigma_T}{T}\right)\\
&=(e^{-\frac{1}{2}s^2}+o(1))\biggl[\mu_0({\Rect}_0)\pm\delta_k 
\mu_0\biggl( \bigcup_{\wt{C}_\ell\in\mathfs{C}_k} \wt{C}_\ell\biggr)\biggr]
+O\left(\frac{\sigma_T}{T}\right).
\end{align*}
Recall that $\mu_0({\Rect}_0)=1$,   $\mu_0(\bigcup\mathfs C_k)$ is uniformly bounded,   $\sigma_T=O(\sqrt{\log T})$, and  $\delta_k$ can be made arbitrarily small by choosing $k$ sufficiently large.  Thus, if $(k_T-a_T)/\sigma_T\to s$, then $\sqrt{2\pi\sigma_T^2}\times\mathfrak{m}_T(\wt{\omega}_0,k_T)\to e^{-\frac{1}{2}s^2}$, as required. 
\qed
\medskip

\noindent
{\bf Proof of the Local Beck Theorem (Theorem \ref{t.LLT-for-D_N}).}
Let $\tomega_0$ denote the midpoint of  the top side of $\Rect_0$. Let $T_n:=n$,  and define as before  $A_n:=\frac{1}{2}a_{T_n}(\tomega_0)$ and $B_n:=\frac{1}{2}\sigma_{T_n}$. We saw above that  $A_n=O(\log n)$ and $B_n\asymp\sqrt{\log n}$.

Suppose $\xi_n\in\frac{1}{2}\Z$, and $\frac{\xi_n-A_n}{B_n}\xrightarrow[n\to\infty]{}t$. Then $\frac{2\xi_n-a_{T_n}(\tomega_0)}{\sigma_{T_n}}\xrightarrow[n\to\infty]{}t$, and
by \eqref{e.Lin-Flow-1}, 
\begin{align*}
\mathfrak{d}_n(\xi_n)&:= \#\{1\leq N\leq n: D_N(\alpha)=\xi_n\}=\mathrm{mes}\{-1<t<T_n-1: \zc(\phi^{t}(\tomega_0))=2\xi_n \}\\
&=T_n\mathfrak{m}_{T_n}(\tomega_0,2\xi_n)+O(1)=[1+o(1)]\frac{T_n e^{-t^2/2}}{\sqrt{2\pi\sigma_{T_n}^2}}+O(1).
\end{align*}
Since  $B_n=\frac{1}{2}\sigma_{T_n}\asymp \sqrt{\log n}$, ${\mathfrak d}_n(\xi_n)=[1+o(1)]\frac{ne^{-t^2/2}}{\sqrt{8\pi B_n^2}}$, as stated in Theorem \ref{t.LLT-for-D_N}.
\qed

\section{Change of Measure and Pressure Functions}
\label{ScPressure}
In the previous section we studied $\mathrm{mes}\{0<t<T:\phi^t(\tomega)\in\Rect_{k_T}\}$ as $T\to\infty$ by analyzing expressions of the form 
$p_n:=\Prob[S_n-\E(S_n)=\delta_n|X_0=a_n, X_1=b_n]$, see  \eqref{e.Prob[S-n-prime]}. 
The temporal limit theorems  assume that  ``the deviations from the mean" $\delta_n$ are {\em local}, i.e. of order $\sqrt{V_n}$.
The theorems describing  $\mu_\xi$-generic points are more complicated, because to handle $\mu_\xi$-a.e. $\tomega$, we must allow  $\delta_n$ such that $\limsup\limits_{n\to\infty}|\delta_n|/\sqrt{V_n}=\infty$. In such cases, the problem of finding $p_n$ becomes a problem in  {large deviations}. 

A classical idea due to H. Cram\'er \index{Cram\'er's change of measure|see {Change of measure}}\index{Change of measure} is that  limit theorems for large deviations can sometimes be reduced to limit theorems for local deviations, by means of a {\em change of measure}: A modification of the underlying probability distribution which tilts  the expectation of $S_n$ in the direction of $\delta_n$ so much, that  the deviation from the new expectation value is bounded --- whence local.

This section prepares the ground for implementing this idea  in our context.

\subsection{Change of Measure}\label{ss.change-of-measure} 
Let $
\pi_{k,k+1}(x,y):=\Prob[X_{k+1}=y|X_k=x]$. 
 $\{X_k\}_{k\geq 0}$ is uniformly elliptic, and  $\{f_k\}_{k\geq 1}$ us uniformly bounded, therefore we can implement the following  general   construction from \cite[\S\S7.3.2--7.3.3]{Dolgopyat-Sarig-Book}: 

\begin{lem}
\label{LmAddPres}
For some $a_n\in\R$, 
there are unique 
 real-analytic  $p_n:\R\to\R$ and unique parameterized families of functions  $h_n(\cdot, \xi):\mathfrak P_n\to\R^+$ $(\xi\in\R)$ such that 
  for all $n\in\
 \N\cup\{0\}$ and  $\xi\in\R$, 
 \begin{equation}
 \sum_x h_n(x, \xi) \Prob(X_n=x)=e^{a_n\xi},\   h_n(\cdot,0)\equiv 1, 
 \end{equation}
\begin{equation}
\label{e.eigenfunction}
\sum_{y\in\mathfrak P_{n+1}} e^{\xi f_{n+1}(y)}\pi_{n,n+1}(x,y) h_{n+1} (y, \xi)=e^{p_n(\xi)} h_n(x, \xi)\ \ (x\in\mathfrak P_n),  
\end{equation}
\begin{equation}\label{e.p-prime-constraint}
\ p_i(0)=0,\ \text{ and } p_1'(0)+\cdots+p_n'(0)=\E(S_n).
\end{equation}
In addition, $h_i(\cdot,\xi)$ is real-analytic in $\xi$, and there are constants  $\sC=\sC(R)$ $(R\!>\!0)$ 
which depend only on $M(\beta)$ such that 
\begin{equation}
\DS \frac{1}{\sC}\leq h_n(\cdot,\xi)\leq \sC, \text{ and }\quad |p_n(\xi)|\leq \sC \text{ for  all
$n\in \naturals$ and $\xi\in [-R, R]$}.
\end{equation}
\end{lem}
\noindent 
Eq. \eqref{e.p-prime-constraint} determines $a_n$. 

\begin{defi}
The {\em change of measure} of $\{X_n\}_{n\geq 0}$ {\em with parameter $\xi$} is the Markov chain $\{\wt{X}^\xi_n\}_{n\geq 0}$ \index{Change of measure}  with   state spaces $\mathfrak P_k$,  initial distribution $\Prob[X_0=\cdot]$, and transition kernel 
\begin{equation}\label{ProbXi}
\Prob[\wt{X}_{n+1}^\xi=y|\wt{X}_n^\xi=x]=\pi_{n,n+1}^\xi(x,y):=\frac{e^{\xi f_{n+1}(y)} h_{n+1}(y, \xi)}{e^{p_n(\xi)} h_n(x,\xi)}\pi_n(x,y).
\end{equation}
\end{defi}
\noindent Eq. \eqref{e.eigenfunction} guarantees that $\DS \sum_{y\in\mathfrak P_{n+1}}\pi_{n,n+1}^\xi(x,y)=1$.
 Henceforth we denote probabilities calculated with respect to  $\{\wt{X}^\xi_n\}_{n\geq 0}$ by $\wt{\Prob}^{\xi}$, and let 
\begin{equation}\label{e.V_n-tilde}
\wt{\E}^\xi(S_n):=\E[f_1(\wt{X}_1^\xi)+\cdots+ f_n(\wt{X}_n^\xi)]\ ,
\ 
\wt{V}_n^\xi:=\mathrm{Var}(f_1(\wt{X}_1^\xi)+\cdots+ f_n(\wt{X}_n^\xi))
\end{equation}
\noindent

We have the following fact from \cite[Lemmas 7.17, 7.21]{Dolgopyat-Sarig-Book}:
\begin{lem}\label{l.V-n-xi}
For every $R>0$ \index{Change of measure!variance estimates}there exists $\mathsf C=\mathsf C(R)$ such that for all for all  $N$ large enough, for all $\xi\in [-R,R]$, $\mathsf C^{-1}\leq \wt{V}^\xi_N/V_N\leq \mathsf C$. For every $\eps>0$ there exists $\xi^\ast(\eps)>0$ such that for all $N$ large enough, for all $|\xi|\leq \xi^\ast$, $e^{-\eps}\leq \wt{V}^\xi_N/V_N\leq e^{\eps}$. 
\end{lem} 

\begin{lem}\label{l.uniform-Doeblin}
Suppose $M(\beta)\leq M$. \index{Uniform ellipticity}\index{Change of measure!uniform ellipticity}
For every $R$, there exists $\eps_0(M,R)$ such that for every $\xi\in [-R,R]$, $\pi_n^\xi(\cdot,\cdot)$ is uniformly elliptic with  ellipticity constant $\eps_0(R)$. 
\end{lem}
\begin{proof}
By Proposition \ref{Prop-Doeblin}, $\pi_{n,n+1}$ is uniformly elliptic with ellipticity constant $\eps_0(M)$. By Lemma \ref{Lemma-f(X)},  $f_n$ is bounded by a constant depending only on $M$, and  by Lemma \ref{LmAddPres}, there is a constant $\mathsf C(M,R)$ for which  $|p_n(\xi)|\leq \mathsf C(M,R)$ and $|\log h(\cdot,\xi)|\leq  \log\mathsf C(M,R)$ when $|\xi|\leq R$  (for the formula of this constant, see the proof of Lemma 7.14 in \cite{Dolgopyat-Sarig-Book}). Now the lemma follows directly from the definition of $\pi_{n,n+1}^\xi$.
\end{proof}

The change of measure method works as follows. Suppose we want to calculate $\Prob[S_n=z_n]$, but the deviations $\delta_n:=z_n-\E(S_n)$ are not local, i.e.  $\limsup|\delta_n|/\sqrt{V_n}=\infty$. 
The (crucial) first step is to find parameters $\xi_n$ so that $\wt{\E}^{\xi_n}(S_n)=z_n+O(1)$. Since the solution $\xi_n$ may depend on $n$, this leads to  an infinite {\em array}  of {\em finite} Markov chains 
$$
\begin{array}{ccccc}
\wt{X}_0^{\xi_0} &  &  & & \\
\wt{X}_0^{\xi_1} & \wt{X}_1^{\xi_1} & &  & \\
\wt{X}_0^{\xi_2} & \wt{X}_1^{\xi_2} & \wt{X}_{2}^{\xi_2} & & \\
\vdots & \vdots &\vdots &\ddots&\\
\wt{X}_0^{\xi_n} & \wt{X}_1^{\xi_n} & \wt{X}_2^{\xi_n} & \cdots\cdots & \wt{X}_{n}^{\xi_n}
\end{array}
$$
The rows in this array are considered separately, without any joint distribution. By Lemmas \ref{l.V-n-xi} and \ref{l.uniform-Doeblin}, 
{\em if $\xi_n$ are  bounded}, then the rows of the array are uniformly elliptic with the same ellipticity constant and  $\wt{V}^{\xi_N}_N\asymp V_N$.  and {\em by the choice of $\xi_n$},  $\wt{\delta}_n:=z_n-\wt{\E}^{\xi_n}(S_n)$ are local (indeed bounded). Thus one can calculate $\wt{\Prob}^{\xi_n}[S_n=z_n]$, using the local limit theorem for  arrays (see Appendix \ref{AppMC}). 
Once this is done, one can  use the explicit form of the change of measure, to translate the asymptotics of $\wt\Prob^{\xi_n}[S_n=z_n]$ into  the asymptotics of $\Prob[S_n=z_n]$. We emphasize:

\medskip
\begin{center}
\textbf{\textit{This method does not work if $\boldsymbol{\xi_n}$  do not exist, or if $\boldsymbol{\xi_n}$  are not bounded.}}
\end{center}

\subsection{Pressure Functions}
These are tools for finding $\xi_n$. Recall the functions $p_n(\xi)$ from Lemma \ref{LmAddPres}. The {\em $n$-th pressure function} is \index{Pressure}
\begin{equation}
\PP_n(\xi):=p_1(\xi)+\cdots+p_n(\xi).\label{AddPres}
\end{equation}

\begin{lem}
\label{LmPressure}
For each $R\!\!>\!\!0$, there is a constant $\sC=\sC(R)$, so that  for all $\xi\in [-R,R]$,
 \begin{equation}
\label{Pres0thDer}
\PP_n(0)=0\text{ and }\left|\PP_n(\xi)-\log\E(e^{\xi S_n})\right|\leq \sC, 
\end{equation} 
 
\begin{equation}
\label{Pres1stDer}
\PP_n'(0)=\E(S_n)\text{ and }\left|\PP_n'(\xi)-\wt{\EXP}^\xi (S_n)\right|\leq \sC, 
\end{equation} 
\begin{equation}
\label{Pres1stDerLM}
 \left|\PP_n'(\xi)-\frac{d}{d\xi} \left[\log \EXP\left(e^{\xi S_n}\right)\right]  \right|\leq \sC, 
\end{equation}
\begin{equation}
\label{Pres2dDer}
\frac{V_n-1}{\sC}  \leq \PP_n''(\xi)\leq \sC[V_n+1] .
\end{equation}
\end{lem}
\noindent

\begin{proof}
Equations \eqref{Pres0thDer}--\eqref{Pres1stDerLM} are the content of parts 1,4 and 5 in Lemma 7.18 in \cite{Dolgopyat-Sarig-Book}.   \eqref{Pres2dDer} follows from Lemma 7.17 in \cite{Dolgopyat-Sarig-Book}, after noting that the function denoted there  by $\ov{P}_n$ differs from  $\PP_n$ by a linear function.
\end{proof}

\begin{cor}
Suppose there exists  $R>0$ such that  $z_n\in (\PP_n'(-R),\PP_n'(R))$, then for all $n$ large enough that $V_n>1$, there exists a unique $\xi_n\in\R$ such that $\PP_n'(\xi_n)=z_n$. These $\xi_n$ satisfy
$
|z_n-\wt{\E}^{\xi_n}(S_n)|\leq \mathsf C(R)\text{ and } |\xi_n|\leq R.
$
\end{cor} 
\begin{proof}
By Lemma \ref{LmAddPres}, $\PP_n(\xi)$ is real-analytic on $\R$. 
By Corollary \ref{Prop-Var-Growth},  $V_n>1$ for all  $n$ large. For such $n$, Equation \eqref{Pres2dDer} says that $\PP_n(\xi)$ is strictly convex on $\R$. It follows that $\PP_n'$ is continuous and strictly increasing.  Since 
$\PP_n'(-R)<z_n<\PP_n({R})$, there exists a unique $\xi_n$ such that $\PP_n'(\xi_n)=z_n$, and $|\xi_n|<R$. By \eqref{Pres1stDer}, $|z_n-\wt{\E}^{\xi_n}(S_n)|\leq \mathsf C(R)$. 
\end{proof} 

\begin{center}
To  find $\xi_n$, we will solve  the equation $\PP_n'(\xi_n)=z_n$.
\end{center}

\subsection{Geometric Pressure Functions}\label{ss.geom-press} \index{Geometric pressure}
One of the difficulties in studying the equation $\PP_n'(\xi_n)=z_n$ is that $\PP_n(\cdot)$ are not given to us explicitly. 
In this section we relate the equation $\PP_n'(\xi)=z$ to  the equation $\cP_{nN_0}'(\eta)=z$, where $N_0$ is the constant  such that
$$
\wt{\Psi}_{nN_0}=\tpsi_{n-1}\circ\cdots\circ\tpsi_0,
$$
and  $\cP_N(\xi)$ is the    {\em geometric pressure function} from \eqref{TempPres}, given  equivalently by
$$
\mathcal P_N(\xi)=\log\E_{\mu_{\Rect_0}}(e^{\xi\Upsilon_N}), 
$$
where $\mu_{\Rect_0}:=\mu_0|_{\Rect_0}=\frac{1}{2}\times\,$the area measure on $\Rect_0$, 
$
\Upsilon_N(\tomega):=\zc[\wt{\Psi}_N(\tomega)]-\zc(\tomega) 
$ and  
  $\wt{\Psi}_N$ is the unique homogeneous automorphism with derivative  $\left(\begin{array}{cc} p_{2N-1} & p_{2N}\\
q_{2N-1} & q_{2N}
\end{array}\right)^{-1}$. Here   $p_N/q_N$ is the reduced form of $[0;2m_1,2m_2,\ldots,2m_N]$ (the truncated even CFE of $\beta$). See \eqref{e.good-matrix} for the choice of signs of $p_n,q_n$.

The following result relates the equations $\PP_n'(\xi)=z_n$, to  the equations  $\mathcal P_{N}'(\xi)=z_n$, for $N$ such that $\lfloor N/N_0\rfloor=n$.  

\begin{prop}\label{p.Pressure-and-Geometric-Pressure}
Fix $z_k\in\R$, and integers $n_k\to\infty$, $0\leq \ell_k\leq N_0-1$,  $N_k=n_k N_0+\ell_k$.\index{Geometric pressure!and pressure} \index{Pressure!and geometric pressure}
\begin{enumerate}[(1)]
\item For all $k$ large enough, there can be at most one $\xi_k$ such that $\mathcal P_{N_k}'(\xi_k)=z_k$, and at most one $\eta_k$ such that $\PP_{n_k}'(\eta_k)=z_k$.
\item There is a bounded sequence of $\xi_k$ such that for all $k$ large enough $\mathcal P_{N_k}'(\xi_k)\!=\!z_k$ iff  there is  bounded sequence of $\eta_k$ such that for all $k$ large enough, $\PP_{n_k}'(\eta_k)=z_k$. Moreover,  $|\xi_k-\eta_k|\to 0$. 
\item If $\mathcal P_{N_k}'(\xi_k)=\PP_{n_k}'(\eta_k)$, then $\xi_k\to\infty$ iff $\eta_k\to\infty$, and $\xi_k\to-\infty$ iff $\eta_k\to-\infty$.
\end{enumerate}
\end{prop}

\subsection{Preparations for the Proof of Proposition \ref{p.Pressure-and-Geometric-Pressure}}

Let $\lambda>1$ be the constant from Proposition \ref{Prop-Contracting-Automorphisms}, then  $\|d\psi_i|_{E^s_i}\|<\lambda^{-1}$ and $\|d\psi_i|_{E^u_i}\|>\lambda$ for all $i$.

\begin{lem}\label{Lemma-is-it-really-this-hard-to-prove?}
There are constants  ${c(R)}, K(R)>0$ and   $C(k,R)$, $N(k,R)$ ($k\geq 1$) as follows. 
For all $k>K(R)$, $n>N(k,R)$ and $N$ such that  $\lfloor N/N_0\rfloor=n$,
$$
\mathcal P_{N}'(\xi)=
e^{\pm c(R)\lambda^{-k}}\bigl(\PP_{n}'(\xi)\pm C(k,R)\bigr), \ \text{ whenever }|\xi|\leq R.
$$ 
\end{lem}

Fix $R>0$ and divide with reminder  $N:=nN_0+\ell$, $0\leq\ell\leq N_0-1$. The proof of the lemma is based on the
following heuristic:
$$
``\mathcal P_N'(\xi)=\frac{d}{d\xi}\log\E_{\mu_{\Rect_0}}(e^{\xi\Upsilon_N})\overset{?}{\approx}\frac{d}{d\xi}\log\E(e^{\xi S_n})\overset{!}{=}\PP_n'(\xi)+O(1)"
$$
Equality $\overset{!}{=}$ is due to \eqref{Pres1stDerLM}, and most of the work is in quantifying the error in $\overset{?}{\approx}$, using the decomposition $\Upsilon_N=\ov{S}_n-\Delta\zc_0+\Delta\zc_{n,\ell}$ from Lemma \ref{l.Phi-Decomp}.

 Let 
$\DS
B_0:=\sup_{n\in\N} \left(\sup_{\Rect_0} \|\Upsilon_N-\ov{S}_n\|\right)+1<\infty. 
$

\begin{sublem}\label{sl.small-xi}
Let $r(k):=1/(2B_0\lambda^k)$,  $C:=\mathsf C(1)+B_0 e^{2B_0}$, { where $\mathsf{C}(R)$ comes from 
Lemma \ref{LmPressure},} and let $K_0:=|\log_\lambda(2B_0)|+1$. Then for all $k>K_0$, for all $|\xi|\leq r(k)$,
$$
\mathcal P_{N}'(\xi)={}e^{\pm\lambda^{-k}}\left(\PP_n'(\xi)\pm C\right).$$
\end{sublem}
\begin{proof} $\DS \mathcal P_N'(\xi)=\frac{\int\, e^{\xi\Upsilon_{N}}\Upsilon_N d\mu_{\Rect_0}}{\int\,e^{\xi\Upsilon_N}d\mu_{\Rect_0}}=\frac{\int\, e^{\xi\Upsilon_{N}}(\ov{S}_n\pm B_0)d\mu_{\Rect_0}}{\int\,e^{\xi\Upsilon_N}d\mu_{\Rect_0}}=
\frac{\int\, e^{\xi\Upsilon_{N}}\ov{S}_n\, d\mu_{\Rect_0} }{\int\, e^{\xi\Upsilon_{N}}d\mu_{\Rect_0}}\pm B_0$
\begin{align}
&=e^{\pm 2B_0|\xi|}\frac{\int\, e^{\xi\ov{S}_n}\ov{S}_n d\mu_{\Rect_0}}{\int\,e^{\xi\ov{S}_n}d\mu_{\Rect_0}}\pm B_0
=e^{\pm 2B_0|\xi|}\frac{d}{d\xi}
\left(\log\E(e^{\xi \ov{S}_n})\right)\pm B_0\notag\\
&=e^{\pm 2B_0|\xi|}(\PP_n'(\xi)\pm\mathsf C(r))\pm B_0, \text{ by \eqref{Pres1stDerLM} and since $\ov{S}_n=S_n$ in distribution.}\label{e.pre-factor}
\end{align}
If   $r(k):=(2B_0\lambda^k)^{-1}$, then $2B_0|\xi|\leq \lambda^{-k}$. If $k>K_0$, then $r(k)<1$, so $\mathsf C(r)$ can be replaced by $\mathsf C(1)$.
\end{proof}

Next we consider the case  $\lambda^{-k}\ll |\xi|\leq R$, when the 
the pre-factor in \eqref{e.pre-factor} is not $1+O(\lambda^{-k})$. This requires a  careful analysis of the decomposition $\Upsilon_N=\ov{S}_n-\Delta\zc_0+\Delta\zc_{n,\ell}$. 

Recall that $\DS \ov{S}_n=\sum_{i=1}^n f_i(\ov{X}_i)$, where 
 $\ov{X}_i:=X^{(0)}_i\circ\pi: \Rect_0\to\mathfrak P_i$, and $X^{(0)}_i(\omega)$ is the $\mathfrak P_i$-element containing $(\psi_{i-1}\circ\cdots\circ\psi_0)(\omega)$.
 
Given non-negative integers $i<j$, let $\ov{\mathfrak F}_{i}^{j}$ denote the $\sigma$-algebra on $\Rect_0$ generated by $\ov{X}_{i},\ldots,\ov{X}_{j}$. For example,
$
\ov{\mathfrak F}_{0}^{k}=\{(\pi|_{\Rect_0})^{-1}(C): C\in \mathfrak P_0^k\}
$,
where $$
\DS \mathfrak P_0^k:=\{\<\dot{a}_0,a_1,\ldots,a_{k}\>:\text{ positive area, } a_i\in\mathfrak P_i,\}.
$$
Let $\ov{\mathcal F}_0^\infty$ denote the $\sigma$-algebra on $\Rect_0$ generated by all the $\ov{X}_i$.  
By Lemma \ref{l.sigma-algebra}, 
\begin{equation}\label{e.conditional-expect-id}
\E_{\mu_{\Rect_0}}(g|\ov{\mathcal F}_0^\infty)(x)=\frac{1}{|W^s_0(x)|}\int_{W^s_0(x)}g\, ds\ \ \ (\text{integral with respect to arc length}).
\end{equation}

\noindent
Recall that $\Delta\zc_0=\zc-\wt{\zc}_0$. 
\begin{sublem}\label{sl.A}
There exist a constant $c_0(R)>0$ and $\ov{\mathcal F}_0^k$-measurable real-valued  functions  $A_{k,\xi}$ on $\Rect_0$ such that for all $\xi\in [-R,R]$ and $k\in\N$,
$$
\E(e^{-\xi\Delta\zc_0}|\ov{\mathcal F}_0^\infty)(x)=e^{\pm c_0(R)\lambda^{-k}}A_{k,\xi}(x).
$$
\end{sublem}
\begin{proof}By identity \eqref{e.conditional-expect-id}, 
it is sufficient to find $c_0(R)>0$ such that for every $\xi\in[-R,R]$ and  $k\in\N$,  if $x,y\in{\Rect}_0$ and  $\pi(x),\pi(y)$ are  in the same $\mathfrak P_0^k$-atom, then 
$$
e^{- c_0(R)\lambda^{-k} }\leq {\frac{1}{|W_0^s(x)|}\int_{W^s_0(x)} e^{-\xi\Delta\zc_0}ds}\bigg/
{\frac{1}{|W_0^s(y)|}\int_{W^s_0(y)} e^{-\xi\Delta\zc_0}ds}\leq e^{c_0(R)\lambda^{-k} }.
$$
 
By Lemma \ref{Lemma-Modified-Approx-Canonical}, ${\Rect}_0$  intersects at most finitely many copies of $\wt{F}_0:=[\wt{\zc}_0=0]$; therefore $\wt{\zc}_0$ attains finitely many different values $z_1,\ldots,z_\nu$ on $\Rect_0$.  
Partition $\Rect_0$ into the pieces $H_1,\ldots,H_\nu$ where $\wt{\zc}_0|_{H_i}=z_i$. Then
$$
\frac{1}{|W^s_0(x)|}\int_{W^s_0(x)} e^{-\xi\Delta\zc_0}ds=\sum_{j=1}^{\nu} e^{-\xi z_j}\frac{|W^s_0(x)\cap H_j|}{|W^s_0(x)|}.
$$
$H_j$ are polygons, bounded by linear segments which are either horizontal, or vertical, or parallel to $E^u_0$ or $E^s_0$.
If $\pi(x),\pi(y)$ are in the same atom  $B\in\mathfrak{P}_0^k$, then $W^s_0(x),W^s_0(y)$ are two $E^s_0$-segments with the same length, and their distance is bounded above by $O(\lambda^{-k})$ (because the width of $B$ in the $E^u_0$ direction is $O(\lambda^{-k})$). It follows that 
$$
\bigl||W^s_0(x)\cap H_j|-|W^s_0(y)\cap H_j|\bigr|=O(\lambda^{-k}) \text{ uniformly in $j$, $k$, $x$, and $y$.}
$$

Since $\pi(x),\pi(y)$ are in the same $\mathfrak P_0^k$-atom, they are in the same $\mathfrak P_0$-atom, and therefore $|W^s_0(x)|=|W^s_0(y)|$. Moreover, this quantity is bounded away from zero and infinity, because it equals $\ell^s(a)$ for some $a\in \mathfrak P_0$. 
It follows that 
$$
\biggl|\frac{|W^s_0(x)\cap H_j|}{|W^s_0(x)|}-\frac{|W^s_0(y)\cap H_j|}{|W^s_0(y)|}\biggr|=O(\lambda^{-k}) \text{ uniformly in $j$, $x$, and $y$,}
$$
whence $|\frac{1}{|W^s_0(x)|}\int_{W^s_0(x)} e^{-\xi\Delta\zc_0}ds-\frac{1}{|W^s_0(y)|}\int_{W^s_0(y)} e^{-\xi\Delta\zc_0}ds|=O(\lambda^{-k})$, uniformly in $x,y$. 
By Lemma \ref{Lemma-Modified-Approx-Canonical},  $\frac{1}{|W^s_0(y)|}\int_{W^s_0(y)} e^{-\xi\Delta\zc_0}ds$ is  bounded away from zero. The  result follows.
\end{proof}

We remind the reader that $\wt{\Psi}_N=\tvf_{n,\ell}\circ\tpsi_{n-1}\circ\cdots\circ\tpsi_0$, with $\tvf_{n,\ell}\in\mathcal E$ (see \eqref{e.Capital-Psi-n} and \eqref{e.E-cal}), and 
$
\Delta\zc_{n,\ell}=(\zc\circ\tvf_{n,\ell}-\wt{\zc}_n)\circ\tpsi_{n-1}\circ\cdots\circ\tpsi_0.
$

 We will build a sequence $m(k)$ and  $\ov{\mathfrak F}_{n-k}^{n+m(k)}$-measurable partitions $\wt{\mathfrak D}_{N,k}$ of $\Rect_0$ with good geometric properties, such that $\Delta\zc_{n,k}$ is constant on most of the atoms of $\wt{\mathfrak D}_{N,k}$. We need the following immediate consequence of Lemma \ref{Lemma-Cylinder-Area}:
\begin{sublem}\label{sl.non-eccentric}
There are constants $M_0$ and $m(k)\geq 0$ ($k\geq 0$) such that every cylinder $\<a_{n-k},\ldots,\dot{a}_n\>$ with $n>k$ can be written as a finite disjoint union of cylinders 
$C'=\<a_{n-k},\ldots,\dot{a}_n,c_{n+1},\ldots,c_{n+m}\>$ so that $m\leq m(k)$
and 
\begin{equation}\label{e.nearly-square}
M_0^{-1}\leq \ell^s(C')/\ell^u(C')\leq M_0.
\end{equation}
\end{sublem} 
\noindent
Cylinders satisfying \eqref{e.nearly-square} are called {\em non-eccentric}. Using Sublemma \ref{sl.non-eccentric}, it is easy to construct, given $n$ and $k$, a partition $
\wt{\mathfrak D}_{k,N}$ of $\Rect_0$, made from sets of the form 
$$
\pi|_{\Rect_0}^{-1}(\psi_0^{-1}\circ\cdots\circ\psi_{n-1}^{-1}(C)),
$$
where $C$ is a non-eccentric cylinder $\<a_{n-k},\ldots,\dot{a}_n,\ldots,a_{n+m}\>$ such that $m\leq m(k)$.

\begin{sublem}\label{s.D-n-k}
There exists $c_1>0$ and a $\wt{\mathfrak D}_{N,k}$-measurable set $D_{N,k}$ such that $\Delta\zc_{n,\ell}$ is constant 
on the interiors of sets in $\wt{\mathfrak D}_{N,k}^\ast:=\{D\in \wt{\mathfrak D}_{N,k}:D\subset D_{N,k}\}$, and  for every $\wt{A}\in\wt{\mathfrak P}_0^n:=\{\pi|_{\Rect_0}^{-1}(C): C\in\mathfrak P_0^n\}$, 
\begin{equation}\label{e.much-better-D_n-small}
\mu_{\Rect_0}(\wt{A}\cap D_{N,k}^c)\leq c_1\lambda^{-k} \mu_{\Rect_0}(\wt{A}).
\end{equation}
\end{sublem}
\noindent 
In particular, $\mu_{\Rect_0}(D_{N,k}^c)\leq c_1\lambda^{-k}$. 

\begin{proof}
Every $D'\in\wt{\mathfrak D}_{N,k}$ has the form 
\begin{equation}\label{e.D-form}
D'=\pi|_{\Rect_0}^{-1}((\psi_{0}^{-1}\circ\cdots\circ\psi_{n-1}^{-1})(C)), \text{ where } 
C=\<a_{n-k},\ldots,\dot{a}_n,\ldots,a_{n+m}\>.
\end{equation}
The 
 discontinuities of  $\wt{\zc}_n$ are located in $\DS \bigcup_{i\in\Z} D^i(\partial\wt{F}_n)$, where 
$\wt{F}_n=[\wt{\zc}_n=0]$. 
Clearly,  
$(\psi_{n-1}\circ\cdots\circ \psi_0)(\mathrm{int}(C))\subset \mathrm{int}(a_n)\subset \mathrm{int}(\pi(\wt{F}_n))$, and 
$$(\tpsi_{n-1}\circ \cdots \circ \tpsi_0)(\mathrm{int}(D'))\cap \bigcup_{i\in\Z} D^i(\partial\wt{F}_n)=\emptyset.$$ 
It follows that  $\wt{\zc}_n\circ\tpsi_{n-1}\circ\cdots\circ\tpsi_0$ is constant on the interior of $D'$. 

Thus, all the discontinuities   of $\Delta\zc_{n,\ell}$ on $\mathfrak D_{N,k}$-atoms are due to jumps in the values  of $\zc\circ\wt{\Psi}_N$. These occur in 
$(\wt{\psi}_0^{-1}\circ\cdots\circ\wt{\psi}_{n-1}^{-1})(\wt{\varphi}_{n,\ell}^{-1}\mathfs S)$, where 
$
{\mathfs S}$ is  the union of boundaries of the  horizontal rectangles.

Recall that $\wt{\varphi}_{n,\ell}\in\mathcal E$, where $\mathcal E$ is a {\em finite} set, independent of $N,\ell$. Let 
$
\DS \wt{\mathfs S}:=\bigcup_{\wt{\varphi}\in\mathcal E}\wt{\varphi}^{-1} \mathfs S,
$
and 
let $\wt{\mathfrak D}_{N,k}^\ast$ denote the collection of sets in $\wt{\mathfrak D}_{N,k}$ whose $\wt{\psi}_{n-1}\circ\cdots\circ\wt{\psi}_0$-images do {\em not} intersect $\wt{\mathfs S}$. Let $D_{N,k}:=$union of  $\wt{\mathfrak D}_{N,k}^\ast$. By the previous paragraph,  $\Delta\zc_{n,\ell}$ is constant on the interior of every element of $\wt{\mathfrak D}_{N,k}^\ast$. It remains to prove \eqref{e.much-better-D_n-small}.

Fix some $\wt{A}\in\wt{\mathfrak P}_0^n$, and some $\wt{D}\in\wt{\mathfrak D}_{N,k}$ which intersects it. Let $A:=\pi(\wt{A})$ and $D:=\pi(\wt{D})$ and write $A=\<\dot{a}_0,\ldots,a_n\>$ and 
$D=(\psi_0^{-1}\circ\cdots\circ\psi_{n-1}^{-1})(C)$, 
with $C\!\!=\!\!\<c_{n-k},\ldots,\dot{c}_n,\ldots,c_{n+m}\>$ non-eccentric, and 
$m\!\!\leq\!\!m(k)$.  Necessarily $a_{n-k}\!\!=\!\!c_{n-k},\ldots,a_n\!\!=\!\! c_n$, and 
$
A\cap D=\<\dot{a}_0,\ldots,a_{n},c_{n+1},\ldots,c_{n+m}\>.
$
Let $C':=\<\dot{a}_n,c_{n+1},\ldots,c_{n+m}\>$. Then 
$$
\mu_{\Rect_0}(\wt{A}\cap\wt{D})=\mu_0(A\cap D)\overset{!}{=}\frac{\mu_0(A)\mu_0(C')}{\mu_0(a_n)}\leq G^{-1}\mu_0(A)\mu_0(C')=G^{-1}\mu_{\Rect_0}(\wt{A})\mu_0(C'),
$$
where $G:=\inf\{\mu_0(a): a\in \mathfrak P_n, n\in\N\}$ (a positive constant,  by Proposition \ref{Prop-Markov-Kernel}). The marked equality follows easily from the identity 
$$\mu_0(A\cap D)=\Prob[X_0=a_0,\ldots,X_n=a_n,X_{n+1}=c_{n+1},\ldots,X_{n+m}=c_{n+m}],
$$
(see Proposition \ref{Prop-Markov-Kernel}), and the rules for calculating probabilities for Markov chains.

$C'$ is a parallelogram with sides parallel to $E^s_n$, $E^u_n$. By Lemma \ref{Lemma-Cylinder-Area}, $\ell^u(C')=\ell^u(C)$, and $\ell^s(C')=\ell^s(a_n)$.  
Recalling that $\mu_0=\frac{1}{2}\times$ Euclidean area measure, we find that  
$\mu_0(C')=\frac{1}{2}\ell^s(a_n)\ell^u(C)|\sin\measuredangle(E^u_n,E^s_n)|\leq \ell^s(a_n)\ell^u(C)$. 

By non-eccentricity and Lemma  \ref{Lemma-Cylinder-Area}, 
\begin{align*}
\mu_0(C')
\leq M_0 \ell^s(a_n) \ell^s(C)
\leq M_0\ell^s(a_n)\ell^s(c_{n-k})\lambda^s_{n-k}\cdots\lambda^s_1.
\end{align*}
Since $\lambda^s_i\leq \lambda^{-1}<1$ and  $\sup\{\ell^s(a): a\in \mathfrak P_n, n\geq 1\}<\infty$, there is global constant $c_1'$ such that $\mu_0(C')\leq c_1'\lambda^{-k}$. Taking $c_1'':=G^{-1}c_1'$, we find that 
$
\mu_{\Rect_0}(\wt{A}\cap\wt{D})\leq c_1''\lambda^{-k}\mu_{\Rect_0}(\wt{A}) 
$.
Summing over all $\wt{D}\in\mathfrak D_{N,k}$ inside $D_{N,k}^c$ such that $\mu_0(\wt{A}\cap\wt{D})\neq 0$, we obtain 
\begin{equation}\label{e.sum-of-sqrts}
\frac{\mu_{\Rect_0}(\wt{A}\cap D_{N,k}^c)}{\mu_{\Rect_0}(\wt{A})}\leq c_1''\lambda^{-k} \#\{\wt{D}\in \mathfrak D_{N,k}\setminus \mathfrak D_{N,k}^\ast: \mu_0(\wt{D}\cap\wt{A})\neq 0\}.
\end{equation}

It remains to bound the number of $\wt{D}\in \mathfrak D_{N,k}\setminus \mathfrak D_{N,k}^\ast$ such that $\mu_0(\wt{D}\cap\wt{A})\neq 0$.
Each $\wt{D}$ like that equals $\pi|_{\Rect_0}^{-1}((\psi_{0}^{-1}\circ\cdots\circ\psi_{n-1}^{-1})(C))$, where $C=\<a_{n-k},\ldots,\dot{a}_n,c_{n+1},\ldots,c_{n+m}\>$, $C$ is non-eccentric, and $C$ intersects $\pi(\wt{\mathfs{S}})$. Moreover different sets $\wt{D}$ are represented by disjoint sets  $C$. We will bound the number of  $C$.

By  Lemma \ref{Lemma-Psi-and-Phi},  $\pi(\wt{\mathfs S})$ is made of a bounded number of line segments, with slopes uniformly bounded {\em away} from the slopes of $E^s_n,E^u_n$. By contrast, $A':=\<a_{n-k},\ldots,\dot{a}_n\>$ is a parallelogram with sides {\em parallel} to $E^s_n,E^u_n$. 
Because of this uniform transversality, 
$\pi(\wt{\mathfs S})\cap A'$ is the union of  linear segments $\Sigma_1,\ldots,\Sigma_\nu$, with  lengths $O(\ell^s(A'))$. Since $\pi(\wt{\mathfs S})$ has finite length and $A'$ is has bounded geometry,  the number $\nu$  is bounded by some global constant $\nu_0$.

Choose one of the $E^u_n$ sides of $A'$, and let  $\mathrm{proj}^u_n$ denote the projection from $A'$ to this side, along the $E^s_n$ direction.  Then $\mathrm{proj}^u_n(A'\cap\pi(\wt{\mathfs S}))$ consists of $\nu$ $E^u_n$-segments $\Sigma^u_1,\ldots,\Sigma^u_\nu$ of length $O(\ell^s(A'))$.  

The sets $C$ we are counting are non-overlapping $(E^u_n,E^s_n)$-sub-parallelograms  of $A'$, which intersect 
$\pi(\wt{\mathfs S})$. Their projections to the $E^u_n$-side of $A_n'$ are non-overlapping $E^u_n$-segments of length $\ell^u(C)$, which intersect $\bigcup\Sigma^u_i$, 

By  Lemma \ref{Lemma-Cylinder-Area},  $\ell^s(C)\!\!=\!\!\ell^s(A')$, and by non-eccentricity, 
$\ell^u(C)\!\!\geq\!\! M_0^{-1}\ell^s(C)\!\!=\!\! M_0^{-1}\ell^s(A')$. 
Length considerations reveal that each $\Sigma_i^u$ can intersect at most
$
\frac{|\Sigma^u_i|}{M_0^{-1}\ell^s(A')}+2
$
$E^u_n$-sides of  sets $C$. So the total number of $C$ is at most
$$
\frac{|\Sigma^u_1|+\cdots+|\Sigma^u_\nu|}{M_0^{-1}\ell^s(A')}+2\nu\leq \frac{O(\ell^s(A'))}{M_0^{-1}\ell^s(A')}+2\nu_0\leq \nu_1,
$$
where $\nu_1$ is a global constant. Thus
$
\#\{\wt{D}\in \mathfrak D_{N,k}\setminus \mathfrak D_{N,k}^\ast: \mu_0(\wt{D}\cap\wt{A})\neq 0\}\leq \nu_1
$.

Now  \eqref{e.sum-of-sqrts} implies \eqref{e.much-better-D_n-small}. 
\end{proof}

\begin{sublem}\label{sl.integral-ineq}
There exists $K_1(R),c_2(R), C_2(R)\!\!>\!\!1$ so that for all $k\geq K_1(R)$, $N\!\!>\!\!N_0 k$,   and $\xi\in [-R,R]\setminus\{0\}$,
 \begin{align}
&\int_{D_{N,k}^c}e^{\xi\Upsilon_N}\, d\mu_{\Rect_0}\leq c_2\lambda^{-k} \int e^{\xi\Upsilon_N}\, d\mu_{\Rect_0};\label{e.D_n-additional-small-denom}\\
&\int_{D_{N,k}^c \cap[\xi\ov{S}_n\geq 0]}e^{\xi\Upsilon_N}\, d\mu_{\Rect_0}\leq c_2\lambda^{-k} \int_{[\xi\ov{S}_n\geq 0]} e^{\xi\Upsilon_N}\, d\mu_{\Rect_0};\label{e.D_n-small-denom}\\
&\bigg|\int\limits_{D_{N,k}^c\cap [\xi\ov{S}_n\geq 0]} e^{\xi\Upsilon_N}\Upsilon_N d\mu_{\Rect_0}\bigg|<c_2\lambda^{-k}\!\!\!\!\!\!\!\!\!\!\!\int\limits_{[\xi\ov{S}_n\geq 0]\cap D_{N,k}}   \!\!\!\!\!\!\!\!\!\!\!\! e^{\xi\Upsilon_N}|\Upsilon_N| d\mu_{\Rect_0}+C_2\int e^{\xi\Upsilon_N}d\mu_{\Rect_0}. \label{e.D_n-small}
\end{align} 
In addition, for every $0<r<R$ there exists $N(r,R)$ such that for all  $n>N(r,R)$,
\begin{equation}\label{e.lalala}
\forall\xi\in[-R,R]\setminus[-r,r],\ \ \ \int_{[\xi\ov{S}_n<0]}e^{\xi\Upsilon_N}|\Upsilon_N|d\mu_{\Rect_0}<\int e^{\xi\Upsilon_N}d\mu_{\Rect_0}.
\end{equation}
\end{sublem}

\begin{proof} $\ov{S}_n$ is constant on atoms of $\wt{\mathfrak P}_0^n$. Therefore, multiplying both sides of \eqref{e.much-better-D_n-small} by the constant value of  $e^{\xi\ov{S}_n}$  on each $\wt{A}\in \wt{\mathfrak P}_0^n$, and then summing over  $\wt{A}$, 
we obtain
$
\int_{D^c_{N,k}} e^{\xi\ov{S}_n}d\mu_{\Rect_0}\leq c_1\lambda^{-k}\int e^{\xi\ov{S}_n}d\mu_{\Rect_0}.
$

Similarly,  multiplying both sides of \eqref{e.much-better-D_n-small} by the constant value of  $e^{\xi\ov{S}_n}$ or $e^{\xi\ov{S}_n}|\ov{S}_n|$ on $\wt{A}$, and then summing over all the $\wt{A}\in\wt{\mathfrak P}_0^n$ such that  $\xi\ov{S}_n\geq 0$ on $\wt{A}$, 
we obtain
\begin{equation}\label{e.better-D_n-small}
\begin{aligned}
&\int_{D_{N,k}^c\cap[\xi\ov{S}_n\geq 0]}e^{\xi\ov{S}_n}\, d\mu_{\Rect_0}\leq c_1\lambda^{-k} 
\int_{[\xi\ov{S}_n\geq 0]}e^{\xi\ov{S}_n} \, d\mu_{\Rect_0};\\
&\left|\int_{D_{N,k}^c\cap[\xi\ov{S}_n\geq 0]}e^{\xi\ov{S}_n}\ov{S}_n\, d\mu_{\Rect_0}\right|\leq c_1
\lambda^{-k} \left|\int_{[\xi\ov{S}_n\geq 0]}e^{\xi\ov{S}_n} \ov{S}_n\, d\mu_{\Rect_0}\right|.
\end{aligned}
\end{equation}
Since 
 $\DS B_0=\sup_n\, \esssup|\Upsilon_N-\ov{S}_n|<\infty$, 
 \eqref{e.D_n-additional-small-denom} and  \eqref{e.D_n-small-denom}  follow (with larger $c$).
 
To obtain  \eqref{e.D_n-small}, suppose for simplicity that $0<\xi<R$, then 
\begin{align}
&\int_{D_{N,k}^c\cap[\xi\ov{S}_n\geq 0]}e^{\xi\Upsilon_N}
|\Upsilon_N|
\, d\mu_{\Rect_0}\leq 
e^{B_0 R}\int_{D_{N,k}^c\cap[\xi\ov{S}_n\geq 0]}
e^{\xi\ov{S}_n}(\ov{S}_n+B_0) d\mu_{\Rect_0}\notag\\ 
&\hspace{0.5cm}\leq e^{B_0 R}\int_{D_{N,k}^c\cap[\xi\ov{S}_n\geq 0]}
e^{\xi\ov{S}_n}\ov{S}_n d\mu_{\Rect_0}+B_0 e^{B_0 R}\int_{D_{N,k}^c\cap [\xi\ov{S}_n\geq 0]}
e^{\xi\ov{S}_n}d\mu_{\Rect_0}\notag \\
&\hspace{0.5cm}\leq c_1e^{B_0 R}\lambda^{-k}\int_{[\xi\ov{S}_n\geq 0]}
e^{\xi\ov{S}_n}\ov{S}_n d\mu_{\Rect_0}+B_0 e^{B_0 R}c_1\lambda^{-k}\int_{[\xi\ov{S}_n\geq 0]}
e^{\xi\ov{S}_n}d\mu_{\Rect_0}\text{ by \eqref{e.better-D_n-small}}\notag\\
&\hspace{0.5cm}\leq c_1e^{2B_0 R}\lambda^{-k}\int_{[\xi\ov{S}_n\geq 0]}
e^{\xi\Upsilon_N}\Upsilon_N d\mu_{\Rect_0}+2B_0 e^{2B_0 R}c_1\lambda^{-k}\int_{[\xi\ov{S}_n\geq 0]}
e^{\xi\Upsilon_N}d\mu_{\Rect_0}\notag\\
&\hspace{0.5cm}\leq c_1e^{2B_0 R}\lambda^{-k}\left[\int_{D_{N,k}^c\cap[\xi\ov{S}_n\geq 0]}
e^{\xi\Upsilon_N}\Upsilon_N d\mu_{\Rect_0}+\int_{D_{N,k}\cap[\xi\ov{S}_n\geq 0]}
e^{\xi\Upsilon_N}\Upsilon_N d\mu_{\Rect_0}\right]\notag \\
&
\hspace{0.5cm}\hspace{8cm}+2B_0 e^{2B_0 R}c_1\lambda^{-k}\int_{[\xi\ov{S}_n\geq 0]}
e^{\xi\Upsilon_N}d\mu_{\Rect_0}\notag\\
&\hspace{0.5cm}\leq c_1e^{2B_0 R}\lambda^{-k}\left[\int_{D_{N,k}^c\cap[\xi\ov{S}_n\geq 0]}
e^{\xi\Upsilon_N} |\Upsilon_N| d\mu_{\Rect_0}+\int_{D_{N,k}\cap[\xi\ov{S}_n\geq 0]}
e^{\xi\Upsilon_N}\Upsilon_n d\mu_{\Rect_0}\right]\notag\\
&\hspace{0.5cm}\hspace{8cm}+2B_0 e^{2B_0 R}c_1\lambda^{-k}\int_{[\xi\ov{S}_n\geq 0]}
e^{\xi\Upsilon_N}d\mu_{\Rect_0}.\notag
\end{align}
%
%
%
%
Subtracting $c_1e^{2B_0 R}\lambda^{-k}\int_{D_{N,k}^c\cap[\xi\ov{S}_n\geq 0]}e^{\xi\Upsilon_N} |\Upsilon_N| d\mu_{\Rect_0}$ from both sides, and dividing by \\$(1-c_1 e^{2B_0 R}\lambda^{-k})$, 
we obtain \eqref{e.D_n-small}, provided that
$
k>K_1(R):=\log_\lambda(2c_1e^{2B_0 R}).
$

We prove \eqref{e.lalala}. Suppose $\xi\in [-R,R]\setminus[-r,r]$. 
The left-hand-side of \eqref{e.lalala} is less than $\DS \int_{[\xi\ov{S}_n<0]}
e^{\xi\ov{S}_n+|\xi|B_0}\left(\frac{|\xi\ov{S}_n|}{r}+B_0\right)d\mu_{\Rect_0}$. 

Since  $e^{t}|t|<1$ and $e^t\leq 1$ for all $ t\leq 0$, this integral is less than 
$$
e^{B_0 R}\left[\frac{1}{r}\int_{[\xi\ov{S}_n<0]}
e^{\xi\ov{S}_n}|\xi\ov{S}_n|d\mu_{\Rect_0}+ B_0\int_{[\xi\ov{S}_n<0]}
e^{\xi\ov{S}_n}d\mu_{\Rect_0}\right]\leq e^{B_0 R}\left(\frac{1}{r}+B_0\right).
$$

Next we  bound the right-hand-side of \eqref{e.lalala} from below. Recall that  $\ov{S}_n=S_n$ in distribution, therefore $\mathrm{Var}(\ov{S}_n)=V_n$. Set
$$\mathcal F_n(\xi):=\frac{1}{V_n}\log\E(e^{\xi\ov{S}_n})
=\frac{1}{V_n}\log\E(e^{\xi{S}_n}).$$ 

Recall that $\E(S_n)$ is bounded, and $V_n\to\infty$.  
By  Proposition \ref{PrFreeEn},
$\mathcal F_n(0)=0$, $\mathcal F_n'(0)=\E(\ov{S}_n)/V_n\to 0$, and  for some $N_1=N_1(R)$, $C_1=C_1(R)$ and all $\xi\in [-R,R]$, if $n>N_1$, then $\mathcal F_n''\geq C_1^{-1}$ on $[-R,R]$. So 
$$
\forall\xi\in [-R,R]\setminus[-r,r],\ \ \ \mathcal F_n(\xi)\geq \frac{\E(\ov{S}_n)}{V_n}\xi+\frac{1}{2C_1}\xi^2\geq \frac{r^2}{2C_1}-\frac{R|\E(S_n)|}{V_n}=\frac{r^2}{2C_1}-o(1).
$$
Thus $\exists N(r,R)$ such that for all $n>N(r,R)$, $\mathcal F_n(\xi)>\frac{r^2}{4C_1}$. It follows that
$$
\int e^{\xi\Upsilon_N}d\mu_{\Rect_0}\geq e^{-B_0}\int e^{\xi\ov{S}_n}d\mu_{\Rect_0}\geq \exp\left(\frac{r^2}{4C_1} V_n-B_0\right).
$$ 
Since $V_n\to\infty$, this is eventually larger than $e^{B_0 R}\left(\frac{1}{r}+B_0\right)$, proving \eqref{e.lalala}.
\end{proof}

Recall the constants $r(k),K_0, K_1(R), N(r,R)$ from Sublemmas \ref{sl.small-xi} and \ref{sl.integral-ineq}, and set 
$
K_3(R):=\max\{K_0,K_1(R)\}$, $N_3(k,R):=N(r(k),R).
$
Recall the function $A_{k,\xi}:\Rect_0\to\R$ from Sublemma \ref{sl.A}, and let $B_{N,k,\xi}:=e^{\xi\Delta\zc_{n,\ell}}1_{D_{n,k}}$, with  $D_{n,k}$ from  Sublemma \ref{s.D-n-k}. Then for all $\xi\in [-R,R]$ and  $N>kN_0$, 
\begin{align*}
&A_{k,\xi}\text{ is $\ov{\mathcal F}_0^k$-measurable, and }  \E_{\mu_{\Rect_0}}(e^{-\xi\Delta\zc_0}|\ov{\mathcal F}_0^\infty)=e^{\pm c_0\lambda^{-k}}A_{k,\xi}\\
& B_{N,k,\xi} \text{ is $\ov{\mathcal F}_{n-k}^{n+m(k)}$-measurable, and $e^{\xi\Delta\zc_{n,\ell}}=B_{N,k,\xi}$ on $D_{n,k}$.}
\end{align*}
Given non-negative integers $i<j$, let $\ov{S}_i^j:=f_i(\ov{X}_i)+\cdots + f_j(\ov{X}_j)$.

\begin{sublem}\label{sl.Num-Denom}
For every $R>0$ there are constants  $c_3(R),C_3(k,R)$ $(k\geq 1)$  so that for all $k>K_3(R)$, $N>N_3(k,R)$, and $\xi\in [-R,R]\setminus[-r(k),r(k)]$, 
$$
\mathcal P_N'(\xi)=e^{\pm c_3(R)\lambda^{-k}}\left(\frac{\mathsf{Num}}{\mathsf{Denom}}\right)\pm C_3(k,R), \text{ where }
\DS\mathsf{Num}:=\int_{\Rect_0} e^{\xi\ov{S}_n}\ov{S}_{k+1}^{n-k-1} A_{k,\xi} B_{N,k,\xi} d\mu_{\Rect_0},$$ and 
$\DS\mathsf{Denom} := \int_{\Rect_0} e^{\xi\ov{S}_n} A_{k,\xi} B_{N,k,\xi} d\mu_{\Rect_0}$.
\end{sublem}
\begin{proof}
$
\DS \mathcal P_N'(\xi)=\frac{d}{d\xi}\log\E_{\mu_{\Rect_0}}[e^{\xi\Upsilon_N}]=\frac{\mathsf N}{\mathsf D}$ where $\mathsf N:= \int e^{\xi\Upsilon_N}\Upsilon_N d\mu_{\Rect_0}$ and $\mathsf D:= \int e^{\xi\Upsilon_N} d\mu_{\Rect_0}$. 
 By \eqref{e.D_n-additional-small-denom},
$$
\mathsf D=  (1\pm c_2(R)\lambda^{-k})\int_{D_{N,k}} e^{\xi\Upsilon_N}d\mu_{\Rect_0}.
$$

By \eqref{e.lalala}, for some constant $\hC_3>1$,     
\begin{align*}
&\mathsf N=\int_{[\xi\ov{S}_n\geq 0]} e^{\xi\Upsilon_N}\Upsilon_N d\mu_{\Rect_0} \pm {\hC_3} 
 \mathsf D\\
 &=\int_{[\xi \ov{S}_n\geq 0]\cap D_{N,k}} e^{\xi\Upsilon_N}\Upsilon_N d\mu_{\Rect_0}+\int_{[\xi\ov{S}_n\geq 0]\cap D_{N,k}^c} e^{\xi\Upsilon_N}\Upsilon_N d\mu_{\Rect_0}\pm {\hC_3}  \mathsf D\notag\\
&=  (1\pm c_2(R)\lambda^{-k})\int_{[\xi\ov{S}_n\geq 0]\cap D_{N,k}} e^{\xi\Upsilon_N}\Upsilon_N d\mu_{\Rect_0} \pm 
{\tC_3} 
\mathsf D, \text{ by \eqref{e.D_n-small}. }                \notag \\
\end{align*}
Since $|\Upsilon_n-\ov{S}_n|\leq B_0$ ,  there are constants $c_3'(R)$ and $C_3'(k,R)$ so that
\begin{equation}
\begin{aligned}
\mathsf D&=e^{\pm c_3'(R)\lambda^{-k}}\int_{D_{N,k}} e^{\xi\Upsilon_N} d\mu_{\Rect_0}\\
\mathsf N&=e^{\pm c_3'(R)\lambda^{-k}}\int_{[\xi\ov{S}_n\geq 0]\cap D_{N,k}} e^{\xi\Upsilon_N}\ov{S}_n\, d\mu_{\Rect_0}\pm C_3'(k,R)\mathsf D. \label{e.num-denom}
\end{aligned}
\end{equation}
We now use the identity $\Upsilon_N=\ov{S}_n-\Delta\zc_0+\Delta\zc_{n,\ell}$ (Lemma \ref{l.Phi-Decomp}):
\begin{align*}
\mathsf D
&=e^{\pm c_3'(R)\lambda^{-k}}\E_{\mu_{\Rect_0}}\left[e^{\xi\Upsilon_N} 1_{D_{N,k}}\right]
=e^{\pm c_3'(R)\lambda^{-k}}\E_{\mu_{\Rect_0}}\left[\E_{\mu_{\Rect_0}}\left(e^{\xi\Upsilon_N} 1_{D_{N,k}}\big|\ov{\mathcal F}_0^\infty\right)\right]\\
&=e^{\pm c_3'(R)\lambda^{-k}}\E_{\mu_{\Rect_0}}\left[\E_{\mu_{\Rect_0}}\left(e^{\xi(\ov{S}_n-\Delta\zc_0+\Delta\zc_{n,\ell})}1_{D_{N,k}}\big|\ov{\mathcal F}_0^\infty\right)\right]\\
&=e^{\pm c_3'(R)\lambda^{-k}}\E_{\mu_{\Rect_0}}\!\!\left[ e^{\xi\ov{S}_n}B_{N,k,\xi}\E_{\mu_{\Rect_0}}\left(e^{-\xi\Delta\zc_0}\big|\ov{\mathcal F}_0^\infty\right)\right]\\
&= e^{\pm c_3'(R)\lambda^{-k}}\E_{\mu_{\Rect_0}}\left[  e^{\xi\ov{S}_n} B_{N,k,\xi}\cdot e^{\pm c_0(R)\lambda^{-k}}A_{k,\xi}\right].  
\end{align*}
Similarly, if $\xi>0$, then
$$\DS\mathsf N= e^{\pm c_3'(R)\lambda^{-k}} \E_{\mu_{\Rect_0}}\left[ 1_{[\xi\ov{S}_n\geq 0]} e^{\xi\ov{S}_n}\ov{S}_n B_{N,k,\xi}\cdot e^{\pm c_0(R)\lambda^{-k}}A_{k,\xi}\right]\pm C_3'(k,R)\mathsf D,$$
and if $\xi<0$ we have the same estimate, but with $e^{\mp c_3'(R)\lambda^{-k}}$ instead of $e^{\pm c_3'(R)\lambda^{-k}}$.

Note that the terms inside expectation have constant sign  on their support. Therefore,  we can pull the $e^{\pm c_0(R)\lambda^{-k}}$ outside, and obtain with $c_3^\ast(R):=c_0(R)+c_3'(R)$, 
\begin{equation}\label{e.D}
\mathsf D=e^{\pm c_3^\ast(R)\lambda^{-k}}\E_{\mu_{\Rect_0}}\left[  e^{\xi\ov{S}_n}A_{k,\xi} B_{N,k,\xi}\right]\equiv e^{\pm c_3^\ast(R)\lambda^{-k}}\mathsf{Denom};
\end{equation}
Similarly, if $\xi>0$, then 
\begin{equation}\label{e.N} 
\mathsf N=e^{\pm c_3^\ast(R)\lambda^{-k}}\int_{[\xi\ov{S}_n\geq 0]} e^{\xi\ov{S}_n}\ov{S}_n A_{k,\xi} B_{N,k,\xi} d\mu_{\Rect_0}\pm C_3'(k,R)\mathsf D.
\end{equation}
If $\xi<0$, replace $e^{\pm c_3^\ast(R)\lambda^{-k}}$ by $e^{\mp c_3^\ast(R)\lambda^{-k}}$. Henceforth, we only consider the case $\xi>0$. 
 
 The same calculations, but with $[\ov{S}_n<0]$ instead of $[\ov{S}_n\geq 0]$ show that 
\begin{align*}
&\int_{[\xi\ov{S}_n< 0]\cap D_{N,k}} e^{\xi\Upsilon_N}\ov{S}_n d\mu_{\Rect_0}=e^{\mp c_3^\ast(R)\lambda^{-k}}\int_{[\xi\ov{S}_n< 0]} e^{\xi\ov{S}_n}\ov{S}_n A_{k,\xi} B_{N,k,\xi} d\mu_{\Rect_0}(x),\text{ and therefore, }\\
&\left|\int_{[\xi\ov{S}_n< 0]} e^{\xi\ov{S}_n}\ov{S}_n A_{k,\xi} B_{N,k,\xi} d\mu_{\Rect_0}(x)\right|=e^{\pm c_3^\ast(R)\lambda^{-k}}\int_{[\xi\ov{S}_n< 0]\cap D_{N,k}} e^{\xi\Upsilon_N}|\ov{S}_n| d\mu_{\Rect_0}\\
&\leq e^{c_3^\ast(R)}\left(\int_{[\xi\ov{S}_n< 0]} e^{\xi\Upsilon_N}|\Upsilon_N| d\mu_{\Rect_0}+B_0 \int_{[\xi\ov{S}_n< 0]} e^{\xi\Upsilon_N}d\mu_{\Rect_0}\right)
\leq e^{c_3^\ast(R)}(B_0+1)\mathsf D, \text{ by \eqref{e.lalala}}. 
\end{align*}
Thus  $\DS \int_{[\xi\ov{S}_n\geq 0]} e^{\xi\ov{S}_n}\ov{S}_n A_{k,\xi} B_{N,k,\xi} d\mu_{\Rect_0}$
\begin{align*}
&= \int_{\Rect_0} e^{\xi\ov{S}_n}\ov{S}_n A_{k,\xi} B_{N,k,\xi} d\mu_{\Rect_0}-\int_{[\xi\ov{S}_n< 0]} e^{\xi\ov{S}_n}\ov{S}_n A_{k,\xi} B_{N,k,\xi} d\mu_{\Rect_0}\\
&= \int_{\Rect_0} e^{\xi\ov{S}_n}\ov{S}_n A_{k,\xi} B_{N,k,\xi} d\mu_{\Rect_0}\pm  e^{c_3^\ast(R)}(B_0+1)\mathsf D\\
&= \int_{\Rect_0} e^{\xi\ov{S}_n}(\ov{S}_{k+1}^{n-k-1}\pm 2k\sup_i\|f_i\|_\infty) A_{k,\xi} B_{N,k,\xi} d\mu_{\Rect_0}\pm  e^{c_3^\ast(R)}(B_0+1)\mathsf D\\
&=\int_{\Rect_0} e^{\xi\ov{S}_n}\ov{S}_{k+1}^{n-k-1} A_{k,\xi} B_{N,k,\xi} d\mu_{\Rect_0}\pm  e^{c_3^\ast(R)}(B_0+1)\mathsf D\\
&\hspace{5cm} \pm 2k\sup_i\|f_i\|_\infty\cdot \int_{\Rect_0} e^{\xi\ov{S}_n} A_{k,\xi} B_{N,k,\xi} d\mu_{\Rect_0}\\
&=\mathsf{Num}\pm e^{2c_3^\ast(R)}(B_0+1)\mathsf{Denom}\pm 2k\sup_i\|f_i\|_\infty \mathsf{Denom},\text{ see \eqref{e.D}}.
\end{align*}
Plugging this into \eqref{e.N}, and collecting constants, we find that for some constants $c^\dagger_3(R)$ and $C_3^\dagger(k,R)$,
$
\mathsf N=e^{\pm c_3^\dagger(R)\lambda^{-k}} \mathsf{Num}\pm C_3^{\dagger}(k,R)\mathsf{Denom}.
$ This and \eqref{e.D} lead to 
$$
\mathcal P_N'(\xi)=\frac{\mathsf N}{\mathsf D}=e^{\pm c_3(R)\lambda^{-k}}\left(\frac{\mathsf{Num}}{\mathsf{Denom}}\right)\pm C_3(k,R),
$$
where $c_3(R)=c_3^\ast(R)+c_3^\dagger(R)$, and $C_3(k,R):=C_3^\dagger(k,R)e^{c_3^\ast(R)}$. 
\end{proof}

\begin{sublem}\label{sl.rho-mixing}
For every $R>0$ there are constants $c_4(R)$  and  $C_4(k,R)$ $(k\geq 1)$ such that for all $k>K_3(R)$, $N>N_3(k,R)$ and $\xi\in [-R,R]\setminus [-r(k),r(k)]$, 
$$
\mathcal P_N'(\xi)=e^{\pm c_4(R)\lambda^{-k}}\bigl(\PP_n'(\xi)\pm C_{4}(k,R)\bigr).
$$
\end{sublem}

\begin{proof}
Recall that $A_{k,\xi}$ is $\ov{\mathcal F}_0^k$-measurable, and $B_{N,k,\xi}$ is $\ov{\mathcal F}_{n-k}^{n+m(k)}$-measurable. Therefore, there are functions 
$\wh{A}_{k,\xi}$ and $\wh{B}_{N,k,\xi}$ such that 
$A_{k,\xi}=\wh{A}_{k,\xi}(\ov{X}_0,\ldots, \ov{X}_k)$, and $B_{N,k,\xi}=\wh{B}_{k,\xi}(\ov{X}_{n-k},\ldots, \ov{X}_{n+m(k)})$. 
 By Proposition \ref{Prop-Markov-Kernel}, the joint distribution of 
$\{\ov{X}_i\}_{i\geq 0}$ with respect to $\mu_{\Rect_0}$ is equal to the joint distribution of the Markov chain $\{X_i\}_{i\geq 0}$. Therefore, 
\begin{align*}
&\mathsf{Num}=\E\left(e^{\xi S_n} S_{k+1}^{n-k-1} \wh{A}_{k,\xi}(X_0,\ldots,X_k) \wh{B}_{N,k,\xi}(X_{n-k},\ldots,X_{n+m(k)})\right),
\end{align*}
 where  $S_{i}^{j}=f_{i}(X_{i})+\cdots+f_{j}(X_{j})$ and $S_n=S_1^n$. 
  
 We  write for short
 $A_{k,\xi}=\wh{A}_{k,\xi}({X}_0,\ldots, {X}_k)$, $B_{N,k,\xi}=\wh{B}_{k,\xi}({X}_{n-k},\ldots, {X}_{n+m(k)})$, and 
$$
\mathsf{Num}=\E\left(e^{\xi S_n} S_{k+1}^{n-k-1} {A}_{k,\xi}{B}_{N,k,\xi}\right).
$$
 Note that integrand has the form $F(X_1,\ldots,X_M)$, with $M:=n+m(k)$. 

We would like to rewrite $\mathsf{Num}$ in terms of the change of measure $\{\wt{X}^\xi_i\}_{i\geq 0}$ from 
\S\ref{ss.change-of-measure}. Recall $h_i$ and $\PP_i$ from Lemma \ref{LmAddPres} and \eqref{AddPres}, then  for any Borel  $F(x_0,\ldots,x_M)$, 
\begin{equation}\label{e.com-F}
\wt{\E}^\xi(F(\wt{X}_0^\xi,\ldots,\wt{X}_M^\xi))=\E\left(e^{\xi{S}_M} \frac{h_{M}({X}_{M},\xi)}
{e^{\PP_M(\xi)}h_{0}({X}_{0},\xi)}F({X}_0,\ldots,{X}_M)\right).
\end{equation}
 This leads to the identity 
$\DS\mathsf{Num}=\wt{\E}^\xi\biggl(e^{-\xi\wt{S}_{n+1}^{M}} \frac{e^{\PP_{M}(\xi)}h_{0}(\wt{X}_{0}^\xi,\xi)}{h_{M}
(\wt{X}_{M}^\xi,\xi)} \wt{A}_{k,\xi} \cdot \wt{S}_{k+1}^{n-k-1}\cdot \wt{B}_{N,k,\xi}
\biggr)
$,
where $\wt{S}_i^j\!=\! f_i(\wt{X}_i^\xi)+\cdots+f_j(\wt{X}_j^\xi)$, $\wt{A}_{k,\xi}\!=\!\wh{A}_{k,\xi}(\wt{X}_1^\xi,\ldots,\wt{X}_M^\xi)$.   
Let ${\mathcal F}_i^j$ denote  the $\sigma$-algebra generated by $\wt{X}^\xi_i,\ldots,\wt{X}^\xi_j$ (or equivalently, by $X_i,\ldots,X_j$). Rearranging terms, we have
\begin{align*}
\mathsf{Num}&=e^{\PP_{M}(\xi)}\sum_{j=k+1}^{n-k-1}\wt{\E}^\xi\biggl( \underset{\mathcal F_0^{k}\text{-measurable}}{\underbrace{h_{0}(\wt{X}_{0}^\xi,\xi) \wt{A}_{k,\xi}}} \times f_j(\wt{X}_j^\xi)
\times
\underset{\mathcal F_{n-k}^{M}\text{-measurable}}{\underbrace{\frac{ \wt{B}_{N,k,\xi}e^{-\xi\wt{S}_{n+1}^{M}}}
{h_{M}(\wt{X}_{M}^\xi,\xi)}}}
\biggr).
\end{align*}

By Lemma \ref{l.uniform-Doeblin}, there is a constant $\eps_0(R)$  such that 
$\{\wt{X}_i^\xi\}_{i\geq 0}$ is uniformly elliptic with the same ellipticity constant $\eps_0(R)$   for every $|\xi|\leq R$.  By Corollary \ref{Cor-rho-mixing},  $\{\wt{X}_i^\xi\}_{i\geq 0}$ are exponentially $\rho$-mixing with the same mixing coefficients $C_{mix}=C_{mix}(R)>0$, $\theta=\theta(R)\in (0,1)$ for all $|\xi|\leq R$. 

Applying exponential $\rho$-mixing twice, first for the pair $\mathcal F_0^k, \mathcal F_j^{\infty}$, and then for the pair $\mathcal F_0^j, \mathcal F_{n-k}^\infty$, and noting that all functions in the expectation are uniformly bounded by a global constant, we obtain that for some constants $c_4'(k,R)$ $(k\geq 1)$, 
\begin{align*}
\mathsf{Num}&=\sum_{j=k+1}^{n-k-1}\wt{\E}^\xi\bigl( h_{0}(\wt{X}_{0}^\xi,\xi) \wt{A}_{k,\xi}\bigr)  \wt{\E}^\xi(f_j(\wt{X}_j^\xi))
\wt{\E}^\xi
\biggl(\frac{e^{\PP_M(\xi)} \wt{B}_{N,k,\xi}e^{-\xi\wt{S}_{n+1}^{M}}}
{h_{M}(\wt{X}_{M}^\xi,\xi)}
\biggr)\pm c_4'(k,R)e^{P_M(\xi)}\\
&=\wt{\E}^\xi(\wt{S}_{k+1}^{n-k-1})\wt{\E}^\xi\bigl( h_{0}(\wt{X}_{0}^\xi,\xi) \wt{A}_{k,\xi}\bigr)\wt{\E}^\xi
\biggl(\frac{ e^{\PP_M(\xi)} \wt{B}_{N,k,\xi}e^{-\xi\wt{S}_{n+1}^{M}}}{h_{M}(\wt{X}_{M}^\xi,\xi)}
\biggr)\pm c_4'(k,R)e^{\PP_M(\xi)}.
\end{align*}
Another application of exponential $\rho$-mixing gives constants $c_4''(k,R),c_4^\ast(k,R)$  for which 
\begin{align}
&\mathsf{Num}{=e^{\PP_M(\xi)}\biggl[\wt{\E}^\xi(\wt{S}_{k+1}^{n-k-1})\wt{\E}^\xi\biggl( \frac{ h_{0}(\wt{X}_{0}^\xi,\xi) \wt{A}_{k,\xi} \wt{B}_{N,k,\xi}e^{-\xi\wt{S}_{n+1}^{M}}}{h_{M}(\wt{X}_{M}^\xi,\xi)}
\biggr)\pm (c_4'' \theta^{n-k}\|\wt{S}_{k+1}^{n-k-1}\|_\infty+ c_4')\biggr]}\notag\\
&=\wt{\E}^\xi(\wt{S}_{k+1}^{n-k-1})\wt{\E}^\xi\biggl( h_{0}(\wt{X}_{0}^\xi,\xi) \wt{A}_{k,\xi}\frac{e^{\PP_M(\xi)} \wt{B}_{N,k,\xi}e^{-\xi\wt{S}_{n+1}^{M}}}{h_{M}(\wt{X}_{M}^\xi,\xi)}
\biggr)\pm c_4^\ast(k,R)e^{\PP_M(\xi)}\notag\\
&=\wt{\E}^\xi(\wt{S}_{k+1}^{n-k-1})\E\left(e^{\xi{S}_n}
A_{k,\xi}B_{N,k,\xi}
\right)\pm c_4^\ast(k,R)e^{\PP_M(\xi)}\quad
\text{ by \eqref{e.com-F}}
\notag\\
&=\wt{\E}^\xi(\wt{S}_{k+1}^{n-k-1})
\times e^{\pm c_3^\ast(R)\lambda^{-k}}\mathsf{Denom}\pm c_4^\ast(k,R)e^{\PP_M(\xi)},\text{ see \eqref{e.D}. }\notag\\
&=(\wt{\E}^\xi(\wt{S}_n)+O(k))\times   e^{\pm c_3^\ast(R)\lambda^{-k}}\mathsf{Denom}\pm c_4^\ast(k,R)e^{\mathsf C(R)m(k)}e^{\PP_n(\xi)},\label{e.curly-N-final}
\end{align}
because $\DS \sup_i\|f_i\|_\infty<\infty$ and $\DS \sup_i \sup_{[-R,R]}p_i<\mathsf C(R)$, see Lemma \ref{LmAddPres}.

By \eqref{Pres1stDer}, there is a constant $\mathsf C(R)$ such that for every $\xi\in [-R,R]$, 
$$
\wt{\E}^\xi(\wt{S}_n)=\PP_n'(\xi)\pm \mathsf C(R).
$$
By \eqref{Pres0thDer} and \eqref{Pres1stDer},  there is a constant $\mathsf C'(R)$ such that 
\begin{align*}
&e^{\PP_n(\xi)}\leq \mathsf C'(R)\E(e^{\xi{S}_n})=\mathsf C'(R)\E_{\mu_{\Rect_0}}(e^{\xi\ov{S}_n})\leq e^{B_0}\mathsf C'(R)\int_{R_0} e^{\xi\Upsilon_N}d\mu_{\Rect_0}=e^{B_0}\mathsf C'(R)\cdot \mathsf{D}\\
&\leq e^{B_0+c_3^\ast(R)}\mathsf C'(R)\cdot \mathsf{Denom}, \text{ see  }\eqref{e.D}.
\end{align*}
Plugging these estimates in \eqref{e.curly-N-final}, and collecting constants, we find that for some constants $c_4^\dagger(R), C_4^\dagger(k,R)$ $(k\geq 1)$, 
$$
\frac{\mathsf{Num}}{\mathsf{Denom}}=e^{\pm c_4^\dagger(R)\lambda^{-k}}\PP_n'(\xi)\pm C_4^\dagger(k,R).
$$

We now add the assumptions that $k>K_3(R)$, $N>N_3(k,R)$, and $\xi\in [-R,R]\setminus[-r(k),r(k)]$, and invoke Sublemma \ref{sl.Num-Denom}, to obtain the result.
\end{proof} 

\begin{proof}[{\bf Proof of Lemma \ref{Lemma-is-it-really-this-hard-to-prove?}}] Suppose 
$k\!\!>\!\!K_3(R)$, then $k\!\!>\!\!K_0(R)$. 
The assertion of the lemma holds for all $\xi\!\!\in\!\! [-r(k),r(k)]$, by Sublemma \ref{sl.small-xi}, and for all 
$\xi\!\!\in\!\! [-R,R]\setminus [-r(k),r(k)]$, by  Sublemma \ref{sl.rho-mixing}.
\end{proof}

\begin{lem}\label{Lemma-Uniform-Convexity}
Let $V_n:=\Var(S_n)$.
For every $R>0$ there are $\underline{\delta}(R),\overline{\delta}(R), N(R)>0$ such that for all $n\geq N(R)$ and $N$ such that $\lfloor N/N_0\rfloor=n$,  \index{Pressure!second derivative and variance}\index{Geometric pressure!second derivative and variance}
$$
\mathcal P_N''(\xi), \PP_n''(\xi)\in [ \underline{\delta}(R)V_n, \overline{\delta}(R)V_n]\text{ for all }\xi\in [-R,R].
$$
\end{lem}
\begin{proof}
The assertion for $\PP_n''$ is shown in Lemma \ref{LmPressure}, see \eqref{Pres2dDer}. We prove the assertion for $\mathcal P_N''$.
Consider the following three  probability measures on $\Rect_0$:
\begin{align*}
&\mu_{\Rect_0}:=\text{ normalized Euclidean area}\ , \ 
\wh{m}_{\xi,N}:=\frac{e^{\xi\Upsilon_N}d\mu_{\Rect_0}}{\int_{\mathsf{\Rect}_0} e^{\xi\Upsilon_N}d\mu_{\Rect_0}}\ , \
\ov{m}_{\xi,n}:=\frac{e^{\xi \ov{S}_n}d\mu_{\Rect_0}}{\int_{\mathsf{\Rect}_0} e^{\xi \ov{S}_n}d\mu_{\Rect_0}}. 
\end{align*}
By Lemma \ref{l.Phi-Decomp},  $|\Upsilon_N-\ov{S}_n|\leq B_0$ with $B_0$ independent of $N$, so ${d\ov{m}_{\xi,N}}/{d\wh{m}_{\xi,N}}=e^{\pm 2B_0}$. 

By differentiation, 
$\mathcal P_N''(\xi)=\mathrm{Var}_{\wh{m}_{\xi,N}}(\Upsilon_N):=\int \Upsilon_N^2 d\wh{m}_{\xi,N}-
\left(\int \Upsilon_N d\wh{m}_{\xi,N}\right)^2.
$
 Using the identity $\mathrm{Var}(X+Y)=\mathrm{Var}(X)+\mathrm{Var}(Y)+2\mathrm{Cov}(X,Y)$ and  the Cauchy-Schwarz inequality $|\mathrm{Cov}(X,Y)|\leq\sqrt{\mathrm{Var}(X)\mathrm{Var}(Y)}$, we obtain 
$$
\mathcal P_N''(\xi)=\mathrm{Var}_{\wh{m}_{\xi,N}}(\ov{S}_n)
+O(1)+O(\sqrt{\mathrm{Var}_{\wh{m}_{\xi,N}}(\ov{S}_n)}).
$$
Thus to prove the lemma, it is sufficient to show that there is a constant $C(R)$ such that for all $n$ large enough, for all  $\xi\in [-R,R]$, 
$$
C(R)^{-1}V_n\leq \mathrm{Var}_{\wh{m}_{\xi,N}}(\ov{S}_n)\leq C(R)V_n.
$$
Let $\wh{\mu}_{\xi,N}:=\E_{\wh{m}_{\xi,N}}(\ov{S}_n):=\int \ov{S}_n d\wh{m}_{\xi,N}$,  $\ov{\mu}_{\xi,n}:=\E_{\ov{m}_{\xi,n}}(\ov{S}_n):=\int \ov{S}_n d\ov{m}_{\xi,n}$. 

\medskip
\noindent
{\em Lower Bound\/:} 
$\mathrm{Var}_{\wh{m}_{\xi,N}}(\ov{S}_n)=\frac{\int e^{\xi\Upsilon_n}(\ov{S}_n-\wh{\mu}_{\xi,N})^2 d\mu_{\Rect_0}  }{\int e^{\xi\Upsilon_n}d\mu_{\Rect_0} }\geq e^{-2B_0|\xi|}\frac{\int e^{\xi\ov{S}_n}(\ov{S}_n-\wh{\mu}_{\xi,N})^2 d\mu_{\Rect_0}  }{\int e^{\xi\ov{S}_n}d\mu_{\Rect_0} }.
$
The expectation of a random variable $Y$ is the constant $\mu$ which minimizes $\E[(Y-\mu)^2]$, and the variance of $Y$ is the minimal  value  thus obtained. 

So  changing $\wh{\mu}_{\xi,N}$ to $\ov{\mu}_{\xi,n}$ in the  integral {\em decreases} the fraction to $\mathrm{Var}_{\ov{m}_{\xi,N}}(\ov{S}_n)$, and 
\begin{align*}
\mathrm{Var}_{\wh{m}_{\xi,N}}(\ov{S}_n)&\geq e^{-2B_0|\xi|}\mathrm{Var}_{\ov{m}_{\xi,N}}(\ov{S}_n)\equiv e^{-2B_0|\xi|}\left(\frac{\int e^{\xi\ov{S}_n}\ov{S}_n^2 d\mu_{\Rect_0}}{\int e^{\xi\ov{S}_n}d\mu_{\Rect_0}}-\left(\frac{\int e^{\xi\ov{S}_n}\ov{S}_n d\mu_{\Rect_0}}{\int e^{\xi\ov{S}_n}d\mu_{\Rect_0}}\right)^2
\right)\\
&\equiv e^{-2B_0|\xi|}\frac{d^2}{d\xi^2}\left(\log \int e^{\xi\ov{S}_n}d\mu_{\Rect_0}\right)=e^{-2B_0|\xi|}\frac{d^2}{d\xi^2}\left(\log \E(e^{\xi{S}_n})\right).
\end{align*}
By Theorem 7.3 in \cite{Dolgopyat-Sarig-Book} {(see also Proposition \ref{PrFreeEn})},
there is a constant $C(R)\!>\!1$ such that for all $n$ large enough,
\begin{equation}\label{e.F-is-unif-convex}
C(R)^{-1}V_n\leq \frac{d^2}{d\xi^2}\left(\log  \E(e^{\xi{S}_n})\right)\leq C(R)V_n\text{ for all }\xi\in [-R,R].
\end{equation}
In particular, for all $n$ large enough, $\mathrm{Var}_{\wh{m}_{\xi,N}}(\ov{S}_n)\geq C(R)^{-1}e^{-2B_0R}V_n$.

\medskip
\noindent
{\em Upper Bound\/:} Adding a constant does not change the variance, so 
\begin{align*}
\mathrm{Var}_{\wh{m}_{\xi,N}}(\ov{S}_n)&=
\mathrm{Var}_{\wh{m}_{\xi,N}}(\ov{S}_n-\ov{\mu}_{\xi,n})=
\E_{\wh{m}_{\xi,N}}[(\ov{S}_n-\ov{\mu}_{\xi,n})^2]-
\E_{\wh{m}_{\xi,N}}[(\ov{S}_n-\ov{\mu}_{\xi,n})]^2\\
&\leq \E_{\wh{m}_{\xi,N}}[(\ov{S}_n-\ov{\mu}_{\xi,n})^2]\leq e^{2B_0|\xi|}\E_{\ov{m}_{\xi,N}}[(\ov{S}_n-\ov{\mu}_{\xi,n})^2],
\text{ because $\frac{d\wh{m}_{\xi,N}}{d\ov{m}_{\xi,n}}\leq e^{2B_0|\xi|}$ }\\
&=e^{2B_0|\xi|}\mathrm{Var}_{\ov{m}_{\xi,N}}(\ov{S}_n)=
e^{2B_0|\xi|}\frac{d^2}{d\xi^2}\left(\log  \E(e^{\xi\ov{S}_n})\right)=
e^{2B_0|\xi|}\frac{d^2}{d\xi^2}\left(\log  \E(e^{\xi{S}_n})\right).
\end{align*}
By \eqref{e.F-is-unif-convex},  for all $n$ large enough, for all $\xi\in[-R,R]$, $\mathrm{Var}_{\wh{m}_{\xi,N}}(\ov{S}_n)\leq C(R)e^{2B_0R}V_n$. 
\end{proof}

\subsection{Proof of Proposition \ref{p.Pressure-and-Geometric-Pressure}}

\begin{lem}\label{Lemma-An-Bn}
Let $A_n(\xi),B_n(\xi)$ be two $C^2$ convex functions on $\R$, satisfying the following conditions, for some sequence $V_n\to\infty$ and $M,\lambda>1$:
\begin{enumerate}[(C1)]
\item $\forall R>0$ $\exists$ $c(R)$, $K(R)$, $C(k,R)$ such that for all $k> K(R)$ and for all $n$ sufficiently large,  $A_n'=e^{\pm c(R)\lambda^{-k}}(B_n'\pm C(k,R))$ on $[-R,R]$;
\item $\forall R>0$ $\exists \underline{\delta}(R),\overline{\delta}(R)>0$ such that  $A_n'', B_n''\in[\underline{\delta}(R)V_n,\overline{\delta}(R) V_n]$ on $[-R,R]$;
\item $\exists M$ such that for all $n$, $|A_n'(0)| , |B_n'(0)|\leq M$. 
\end{enumerate}
Then for any sequence of numbers $z_n\in\R$: 
\begin{enumerate}[(1)]
\item The equation $A_n'(\xi)=z_n$ (resp. $B_n'(\xi)=z_n$) can have at most one solution. 
\item There is a bounded sequence of $\xi_n$ such that $A_n'(\xi_n)=z_n$ for all $n$ large enough, iff there is a bounded sequence of $\eta_n$ such that $B_n'(\eta_n)=z_n$ for all $n$ large enough. If this happens, $|\xi_n-\eta_n|\to 0$.
\item If  $A_n'(\xi_n)=B_n'(\eta_n)$ for all $n$ large enough, then $\xi_n\to +\infty$ iff 
$\eta_n\to +\infty$, and  $\xi_n\to -\infty$ iff  $\eta_n\to -\infty$.  
\end{enumerate}
\end{lem}

\begin{proof}
By (C2), $A_n'', B_n''>0$, therefore $A_n'$ and $B_n'$ are strictly increasing, therefore the equations $A_n'(\xi)=z_n$ and $B_n'(\xi)=z_n$ can have at most one solution. This is part 1.

Suppose $\exists \xi_n\in (-R,R)$ such that $A_n'(\xi_n)=z_n$ for all $n$ large enough, and let $\eps$ be  an arbitrary positive number. 

Set
$\DS \eps_k:=\max_{\sigma=\pm 1} |e^{\sigma c(R+1)\lambda^{-k}}-1|=
e^{c(R+1)\lambda^{-k}}-1$, and choose $k=k(\eps,R)$ so large that $k>K(R+1)$, and 
 \begin{align}\label{e.k-choice}
& e^{-c(R+1)\lambda^{-k}}\underline{\delta}(R+1)>\eps_k(1+R\overline{\delta}(R))
\text{ and } \eps_k<\eps \underline{\delta}(R+2)/(R\overline{\delta}(R)).
 \end{align}
By  (C2) and (C3), 
\begin{equation}\label{e.A_n'-UB}
|A_n'|,|B_n'|\leq M+R\overline{\delta}(R)V_n\text{ on }[-R,R]. 
\end{equation}
Choose $N_0\!\!=\!\!N_0(\eps,k,R)$ so that for all $n\!\!>\!\!N_0$, $V_n\!\!>\!\!M$, and (C2) holds with $R,R+1,R+2$. 
For such $n$, by the mean value theorem, 
for any $t=e^{\pm c(R+1)\lambda^{-k}}$, 
\begin{align*}
&tA_n'(R+1)\geq t\bigl(A_n'(R)+\underline{\delta}(R+1)V_n\bigr)
=A_n'(R)+t\underline{\delta}(R+1)V_n-
(1-t)A_n'(R)\\
&\geq A_n'(R)+V_n\left[t\underline{\delta}(R+1)-\eps_k
\bigl(MV_n^{-1}+
R\overline{\delta}(R)\bigr)\right], \text{ by \eqref{e.A_n'-UB} and since $|1-t|\leq \eps_k$ }\\
&\geq A_n'(R)+V_n\underset{\text{positive, by \eqref{e.k-choice}}}{\underbrace{\left[e^{-c(R+1)\lambda^{-k}}\underline{\delta}(R+1)-\eps_k
\bigl(1+
R\overline{\delta}(R)\bigr)\right]}}>A_n'(R)\text{ }.
\end{align*} 

Thus 
$e^{\pm c(R+1)\lambda^{-k}}A_n'(R+1)\!\!>\!\!A_n'(R)\text{ for all $n\!\!>\!\!N_0$},
$ regardless of the sign of $A_n'(R+1)$ or the value of $t$, provided $t\!\!=\!\!e^{\pm c(R+1)\lambda^{-k}}$. 
Fix $N_1\!\!>\!\!N_0$ such that for all $n\!\!>\!\!N_1$, 
 $\underline{\delta}(R+2)V_n\!\!>\!\!C(k,R+1)$.
By  (C2), and then (C1),
\begin{align*}
&B_n'(R+2)\geq B_n'(R+1)+\underline{\delta}(R+2)V_n\\
&= \left(e^{\pm c(R+1)\lambda^{-k}}A_n'(R+1)\pm C(k,R+1)\right)+\underline{\delta}(R+2)V_n\\
&\geq e^{\pm c(R+1)\lambda^{-k}} A_n'(R+1)
>A_n'(R).
\end{align*}
Thus, for all $n$ large enough, $B'_n(R+2)>A_n'(R)$. 
By symmetry,  $B_n'(-R-2)\leq A_n'(-R)$ (apply the previous argument to the pair $(A_n(-t), B_n(-t))$). So  for all $n$ large,
$$
z_n=A_n'(\xi_n)\in (A_n'(-R),A_n'(R))\subset (B_n'(-R-2), B_n'(R+2)).
$$
By the intermediate value theorem, for such $n$  there exists $|\eta_n|\leq R+2$ such that $B_n'(\eta_n)=z_n$. By the first part of the lemma, these $\eta_n$ are unique.

Next we show that $|\xi_n-\eta_n|\to 0$. By the mean value theorem, either $\xi_n=\eta_n$, or 
 $\DS
 \frac{|B_n'(\xi_n)-B_n'(\eta_n)|}{|\xi_n-\eta_n|}\geq \underline{\delta}(R+2) V_n$. Either way,  $\DS |\xi_n-\eta_n|\leq \frac{|B_n'(\xi_n)-B_n'(\eta_n)|}{\underline{\delta}(R+2)V_n}$. 
 By (C1) and the choice of $\eps_k$, 
 \begin{align*}
 &|B_n'(\xi_n)-B_n'(\eta_n)|=|B_n'(\xi_n)-A_n'(\xi_n)|=|(e^{\pm c(R+1)\lambda^{-k}}A_n'(\xi_n)\pm C(k,R+1))-A_n'(\xi_n)|\\
 &\leq C(k,R+1)+\eps_k|A_n'(\xi_n)|
 \leq   C(k,R+1)+\eps_k(M+R\overline{\delta}(R)V_n).
 \end{align*}
 Thus, 
$\DS \limsup_{n\to\infty}|\xi_n-\eta_n|\leq
\limsup_{n\to\infty}\frac{C(k,R+1)+\eps_k(M+R\overline{\delta}(R)V_n)}
{\underline{\delta}(R+2) V_n}=\eps_k\times \frac{R\overline{\delta}(R)}{\underline{\delta}(R+2)}<\eps
$,
see \eqref{e.k-choice}. Since $\eps$ was arbitrary, $|\xi_n-\eta_n|\to 0$. 

This proves  direction $(\Rightarrow)$ in part (2) of the lemma. Direction $(\Leftarrow)$ follows by symmetry, after noting that   $(A_n,B_n)$ satisfies (C1)--(C3) iff  $(B_n,A_n)$ satisfies (C1)--(C3).

Part (3) of the lemma has a similar proof. Again it is enough to prove the implication $(\Rightarrow)$.  It is also sufficient to consider the case $\xi_n\to+\infty$, otherwise we work with $A_n(-\xi), B_n(-\xi)$. Suppose $A_n'(\xi_n)=B_n'(\eta_n)$. 

Fix $R>1$ arbitrary. Arguing as above but with $B_n$ replacing $A_n$, we construct $k=k(R)$ and $N_0$ so 
that for all $n>N_0$,  $e^{-c(R+2)\lambda^{-k}}B_n'(R+1)>B_n'(R)$.
Fix $N_1$ such that for all  $n>N_1$, $\xi_n>R+2$ and 
 $\underline{\delta}(R+2)V_n>C(k,R+2).$
By the mean value theorem,  $B_n'(R+2)\geq B_n'(R+1)+\underline{\delta}(R+2)V_n$, whence 
\begin{align*}
&B_n'(\eta_n)=A_n'(\xi_n)>A_n'(R+2)>e^{-c(R+2)\lambda^{-k}}(B_n'(R+2)
-C(k,R+2))\\
&>
e^{-c(R+2)\lambda^{-k}}(B_n'(R+1)+\underline{\delta}(R+2)V_n-C(k,R+2))
\!\!>\!\!e^{-c(R+2)\lambda^{-k}}B_n'(R+1)\!\!>\!\!B_n'(R).
\end{align*}
So  $\eta_n>R$ for all $n>\max\{N_0,N_1\}$. Since $R$ was arbitrary, $\eta_n\to +\infty$. 
\end{proof}

\medskip
\noindent
{\bf Proof of Proposition \ref{p.Pressure-and-Geometric-Pressure}.} The proposition follows from Lemma \ref{Lemma-An-Bn}, applied to  $A_k(\xi):=\mathcal P_{N_k}(\xi)$ and $B_k(\xi):=\PP_{n_k}(\xi)$. The lemma is applicable: 

Take $V_k:=\Var(S_k)$, then $V_k\to\infty$ by Proposition \ref{Prop-Var-Growth}. 

Condition (C1) is guaranteed by Lemma \ref{Lemma-is-it-really-this-hard-to-prove?}.

Condition (C2) is the content of  Lemma \ref{Lemma-Uniform-Convexity}. 

Condition (C3)  can be verified as follows. $\PP_n'(0)$ is bounded, because $\PP_n'(0)=\E(S_n)$ by construction, and $\E(S_n)=O(1)$ by Proposition \ref{Prop-Drift}. To see that  $\mathcal P_N'(0)$ is  bounded, we first observe by direct differentiation that  $\mathcal P_N'(0)=\int \Upsilon_N d\mu_{\Rect_0}$. By Lemma \ref{l.Phi-Decomp}, $|\Upsilon_N-\ov{S}_n|$ is uniformly bounded,  and $\ov{S}_n$ is equal to $S_n$ in distribution. Therefore  $\mathcal P_N'(0)=\E(S_n)+O(1)=O(1)$.\qed

\section{Characterizations of Generic Points and Consequences}
\label{ScGen}

We remind the reader that  for every $\xi\in\R$,  $\mu_\xi$ is the unique $\phi$-ergodic invariant locally finite measure such that $\mu_\xi\circ D=e^{\xi}\mu_\xi$, and $\mu_\xi(\Rect_0)=1$. By \cite{Hooper-Hubert-Weiss}, any $\phi$-ergodic invariant Radon measure is proportional to some $\mu_\xi$.
Define as in \S\ref{s.dynamical-results}
\begin{align*}
&z_N(\tomega):=\begin{cases}
\zc(\wt{\Psi}_N(\tomega)) & \text{ for }\omega\in\St^\ast\\
\zc(\wt{\Psi}_N(\phi^t(\tomega)))\text{ for all $t>0$ small} & \text{ for }\tomega\text{ singular}
\end{cases}\ ,\ \wt{\Psi}_N\text{ as in \eqref{e.Capital-Psi-n};}\\
&\mathcal P_N(\xi):=\log\int_{\Rect_0} e^{\xi z_N(\tomega)}d\mu_0(\tomega);\\
&\xi_N(\tomega):=
\begin{cases}
-\infty & z_N(\tomega)\leq \mathcal P_N'(-\infty):=\lim\limits_{\eta\to -\infty}\mathcal P_N'(\eta),\\
+\infty & z_N(\tomega)\geq \mathcal P_N'(+\infty):=\lim\limits_{\eta\to +\infty}\mathcal P_N'(\eta),\\
\text{ the unique $\eta_N$ s.t. $\mathcal P_N'(\eta_N)=z_N(\tomega)$} & 
z_N(\tomega)\in (\mathcal P_N'(-\infty), \mathcal P_N'(+\infty)).
\end{cases}
\end{align*}
There is a  constant $N_0$ such that  $\Psi_{nN_0}=\tpsi_{n-1}\circ\cdots\circ\tpsi_0$, see Remark \ref{r.mat-prod-form}.

\subsection{A Sufficient Condition for Genericity} \label{ScBulk}
In this section we prove:
\begin{prop}
\label{PrBulk}
If $\xi_N(\omega)\xrightarrow[N\to\infty]{}\xi$ and $\xi$ is  finite, then 
$\tomega$ is generic for $\mu_\xi.$ If  $\xi$ is a finite limit point of $\xi_N(\omega)$ then $\mu_\xi$ is a historical measure for $\tomega$ (see Definition \ref{d.historical}). 
\end{prop}

\begin{lem}\label{l.solution-xi}
Let $z_n$ be a sequence of real numbers. If there exists a  bounded sequence of $\xi_n\in\R$ such that $\PP_n'(\xi_n)=z_n+O(1)$, then for all $n$ large enough there exist $\xi_n'\in\R$ such that $\PP_n'(\xi_n')=z_n$. Moreover, $\xi_n'=\xi_n+O(n^{-1})$. 
\end{lem}
\begin{proof}
Suppose $|\xi_n|< R$ and $|\PP_n'(\xi_n)-z_n|\leq K$ for all $n$. By \eqref{Pres2dDer} and Proposition~\ref{Prop-Var-Growth}, there exists a positive constant $c$ such that $\PP_n''>cn$ on $[-2R,2R]$, for all $n$ large enough. Necessarily,  
$
\PP_n'(\xi_n+\frac{2K}{cn})>\PP_n'(\xi_n)+cn\cdot \frac{2K}{cn}\geq z_n-K+2K>z_n$. Similarly,   $
\PP_n'(\xi_n-\frac{2K}{cn})<z_n$. By the intermediate value theorem, there exists $\xi_n'\in [\xi_n-\frac{2K}{cn},\xi_n+\frac{2K}{cn}]$ such that $\PP_n'(\xi_n')=z_n$. Since $P_n'$ is strictly increasing,  $\xi_n'$ is unique.
\end{proof}

Let $ \mathfs G:=\bigcup_{\ell\geq 1}\mathfs G_\ell$, where   
 $\mathfs G_\ell$ is the collection of parallelograms  $\wt{C}\subset \St$ of the form $
\wt{C}:=(D^j\circ \pi|_{\wt{F}_0}^{-1})\<\dot{a}_0,\ldots,a_\ell\>.
$ Each $\wt{C}\in\mathfs G$
is an $(E^u_0,E^s_0)$-parallelogram.  The number $j$, called the {\em height of $\wt{C}$}, is the constant value of $\wt{\zc}_0$ on $\mathrm{int}(\wt{C})$.

Recall the change of measure $\{\wt{X}^\xi_n\}_{n\geq 0}$ from \S\ref{ss.change-of-measure}, and the functions $h_i(\cdot,\xi)$ on $\mathfrak P_i$ from Lemma \ref{LmAddPres}. Define $M_\xi:\mathfs G\to \R$ by 
\begin{equation}
\label{DefMM}
M_\xi(\wt{C}):=e^{\xi j} h_0(a_0,\xi) \Prob[\wt{X}^\xi_0=a_0,\ldots,\wt{X}^\xi_\ell=a_\ell],\text{ for }\wt{C}=(D^j\circ \pi|_{\wt{F}_0}^{-1})\<\dot{a}_0,\ldots,a_\ell\>.
\end{equation}

Note that $M_\xi\circ D=e^{\xi} M_\xi$. Eventually it will transpire that  $M_\xi=\const \mu_\xi|_{\mathfs G}$, but we do not know this a priori.

Recall that a pure $E^s_n$-segment $W$ with depth $k$ and height $\ell$ is an $E^s_n$-segment in $D^\ell(\wt{F}_n)$, which projects to a maximal $E^s_n$-segment in some cylinder $\<b_{n-k},\ldots,\dot{b}_n\>$. The height of $W$ is $z_n(W):=\ell=$ the constant value of  $\wt{\zc}_n$ on $W$. Let $\Prob_{W}$ denote the normalized length measure on $W$, and let $\PP_n$ denote the pressure function from \eqref{AddPres}.

\begin{lem}\label{Lemma-Pure-Seg-Ratio-Ergodic-Thm}
Suppose
$n_s$ is a strictly increasing sequence of natural numbers, and $\xi_{n_s}, z_{n_s}, \xi$ are real numbers such that   $\PP_{n_s}'(\xi_{n_s})=z_{n_s}$ and $\xi_{n_s}\to\xi$. 
 Let $W_{n_s}$ be a sequence of pure $E^s_{n_s}$-segments of bounded depth. If  
 $z_{n_s}(W_{n_s})=z_{n_s}+O(1)$, then   
\begin{equation}\label{e.engine2}
\forall \, \wt{C}_1,\wt{C}_2\in\mathfs G,\ \lim_{k\to\infty}\frac{\Prob_{W_{n_s}}\left\{\wt{\omega}\in W_{n_s}:  
(\wt{\psi}_0^{-1}\circ\cdots\circ\wt{\psi}_{n_s-1}^{-1})
(\wt{\omega})\in\wt{C}_1\right\}}
{\Prob_{W_{n_s}}\left\{\wt{\omega}\in W_{n_s}:  (\wt{\psi}_0^{-1}\circ\cdots\circ\wt{\psi}_{n_s-1}^{-1})
(\wt{\omega})\in\wt{C}_2
\right\}}=\frac{M_\xi(\wt{C}_1)}{M_\xi(\wt{C}_2)}.
\end{equation}
For any other sequence of pure $E^s_{n_s}$-segments $W_{n_s}'$ as above,   
for all $\wt{C}\in\mathfs G$, for all $k$ sufficiently large,
\begin{equation}\label{e.engine2bis}
\frac{\Prob_{W_{n_s}'}\left\{\wt{\omega}'\in W_{n_s}':  (\wt{\psi}_0^{-1}\circ\cdots\circ
\wt{\psi}_{n_s-1}^{-1})(\wt{\omega})\in\wt{C}
\right\}}
{\Prob_{W_{n_s}}\left\{\wt{\omega}'\in W_{n_s}:  
(\wt{\psi}_0^{-1}\circ\cdots\circ\wt{\psi}_{n_s-1}^{-1})(\wt{\omega})\in\wt{C}
\right\}}\asymp 1.
\end{equation}
\end{lem}
\begin{proof} 
We assume for simplicity that  $(n_s)_{s\geq 1}$ is the sequence $(1,2,3,\ldots)$, and write $n_s=n$. The routine modifications needed for the general case are left to the reader.

We may assume  without loss of generality that $z_n(W_n)=z_n$ (otherwise we use Lemma~\ref{l.solution-xi} to change $\xi_n$ to $\xi_n'=\xi_n+O(n^{-1})$ such that $\PP_n'(\xi_n')=z_n(W_n)$). 

There is also no loss of generality in assuming $W_n$ has { fixed}  depth $k$, otherwise we split $\{W_n\}$ into a finite number of subsequences with this property. 
In this case,  $\pi(W_n)$ is a maximal $E^s_n$-segment in a cylinder $Q_n=\<b^{(n)}_{n-k},\ldots,\dot{b}^{(n)}_n\>$. 

Fix some 
 $\wt{C}=(D^{j}\circ\pi|_{\wt{F}_0}^{-1})\<\dot{a}_0,\ldots,a_\ell\>$, and define 
\begin{align*}
p_n(\wt{C})&:=\Prob_{W_n}\left\{\wt{\omega}'\in W_n:  (\wt{\psi}_0^{-1}\circ\cdots\circ\wt{\psi}_{n-1}^{-1})(\wt{\omega}')\in\wt{C}
\right\}.
\end{align*}
 Recall $X^{(n)}_i$ from \eqref{e.X-new}. By Lemma \ref{Lemma-f(X)},  
 $\DS\sum_{i=1}^n f_i(X^{(n)}_i(\pi(\tomega)))
 \!\!=\!\wt{\zc}_n((\tpsi_{n-1}\circ\cdots\tpsi_0)(\tomega))\!-\!\wt{\zc}_0(\tomega)$, so  for every $\tomega\in W_n$,
$$ 
 (\wt{\psi}_0^{-1}\circ\cdots\circ\wt{\psi}_{n-1}^{-1})(\wt{\omega})\in\wt{C}
\Leftrightarrow  \left(\begin{array}{l}
{X}_i^{(n)}(\pi(\tomega))=a_i\ (0\leq i\leq \ell)\\
\DS \sum_{i=1}^n f_i(X^{(n)}_i(\pi(\tomega)))=z_n-j
\end{array}\right).
$$ 
Thus $p_n(\wt{C})$ is the probability of an elementary event in $W_n$, and by Lemma \ref{l.prob-of-elementry-event} 
\begin{align}
&p_n(\wt{C})=\Prob\left[
\begin{array}{c}
S_n=z_n-j\\
{X}_0=a_0,\ldots, X_\ell=a_\ell
\end{array}
\bigg| 
{X}_{n-k}=b_{n-k}^{(n)},\ldots,{X}_{n}=b_{n}^{(n)}\right]\notag\\
&=\Prob\left[
\begin{array}{c}
S_n=z_n-j\\
{X}_0=a_0,\ldots, X_\ell=a_\ell
\end{array}
\ ,\ 
{X}_{n-k}=b_{n-k}^{(n)},\ldots,{X}_{n}=b_{n}^{(n)}\right]\bigg/ \mu_0(Q_n).\label{e.engine3}
\end{align}

It is tempting to try to use  the mixing local limit theorem (Theorem \ref{ThLLTArrays}) to find $p_n(\wt{C})$,  but this result  assumes  that $\frac{z_n-\E(S_n)}{\sqrt{V_n}}$ converges to a finite limit, and  in our case $\frac{z_n-\E(S_n)}{\sqrt{V_n}}$ could tend to infinity.  
We will overcome this difficulty using the change of measure  discussed in  \S\ref{ss.change-of-measure}, and the assumptions  $\PP_n'(\xi_n)=z_n$,  $\xi_n\to\xi$.

A direct calculation shows that 
\begin{equation}\label{e.change-of-measure-id}
\Prob[\wt{X}^{\xi_n}_0=a_0,\ldots,\wt{X}_n^{\xi_n}=a_n]=\frac{e^{\xi \sum_{i=1}^n f_i(a_i)}h_{n}(a_{n},\xi_n)}{e^{\PP_n(\xi_n)}h_0(a_0,\xi_n)}\Prob\left(X_0=a_0,\ldots,X_n=a_n\right).
\end{equation}
After some algebraic manipulations this leads 
to the identity
\begin{align*}
&p_n(\wt{C})\!=\!\frac{e^{\PP_n(\xi_n)} h_0(a_0,\xi_n)}
{e^{\xi_n(z_n-j)} h_n(b^{(n)}_n,\xi_n)\mu_0(Q_n)}\Prob\left[\sum_{i=1}^n f_i(\wt{X}_i^{\xi_n})=z_n-j,\!\begin{array}{l}
\wt{X}_i^{\xi_n}=a_i\ \   (0\leq i\leq \ell)\\  
\wt{X}_i^{\xi_n}=b_i^{(n)}
 (n-k\leq i\leq n)
\end{array}\!
\right]\!.
\end{align*}
$\wt{X}^{\xi_n}_i$ is a  Markov chain, so if 
$\DS w_n\!\!:=\!\!\sum_{i=0}^{\ell-1} f_i(a_i)+\!\!\sum_{i={n-k+1}}^n f_i(b^{(n)}_i)$,\; 
$\DS \wt{S}_{\ell}^{n-k}\!\!:=\!\!\sum_{i=\ell}^{n-k} f_i(\wt{X}_i^{\xi_n})$, then 
\begin{align}
&p_n(\wt{C})=   e^{\xi_n j} h_0(a_0,\xi_n)\Prob[\wt{X}^{\xi_n}_0=a_0,\ldots,\wt{X}^{\xi_n}_\ell=a_\ell]\notag\\
&\hspace{0.25cm}\times
\Prob\left[\begin{array}{l}
\DS\sum_{i=\ell}^{n-k} f_i(\wt{X}_i^{\xi_n})=z_n-w_n-j
\end{array}
\bigg|\begin{array}{l}
\wt{X}_\ell^{\xi_n}=a_\ell\\  
\wt{X}_{n-k}^{\xi_n}=b_{n-k}^{(n)}
\end{array}
\right] \label{e.p(C)}
\\&\hspace{0.25cm} 
\times \frac{\Prob[\wt{X}_{n-k}^{\xi_n}=b^{(n)}_{n-k}|\wt{X}_\ell^{\xi_n}=a_\ell]}{\Prob[\wt{X}_{n-k}^{\xi_n}=b^{(n)}_{n-k}]}\notag
\\
& \hspace{0.25cm} \times\Prob[\wt{X}_{n-k}^{\xi_n}=b^{(n)}_{n-k}]\times\Prob
[\wt{X}_{n-k}^{\xi_n}=b^{(n)}_{n-k},\ldots,\wt{X}_n^{\xi_n}=b^{(n)}_n|\wt{X}_{n-k}^{\xi_n}=b^{(n)}_{n-k}]\times
\frac{e^{\PP_n(\xi_n)-\xi_n z_n} }
{h_n(b^{(n)}_n,\xi_n)\mu_0(Q_n)}.\notag
\end{align}

The first term converges to $M_\xi(\wt{C})$, because $\xi_n\to\xi$ and the transition kernel of $\wt{X}^\xi_i$ given in \eqref{ProbXi} is continuous in $\xi$. The third term tends to one as $n\to\infty$ by $\phi$-mixing (Corollary \ref{Cor-phi-mixing}).
The fourth term does not depend on $\wt{C}$ and will cancel out in  \eqref{e.engine2}. It remains to analyze the second term. 

Let $\DS\wt{V}_n^{\xi_n}:=\mathrm{Var}(\sum_{i=1}^{n} f_i(\wt{X}_i^{\xi_n}))$.

\medskip
\noindent
{\em Claim 1.\/}  
$
\DS\Prob\left[\begin{array}{l}
\DS\sum_{i=\ell}^{n-k} f_i(\wt{X}_i^{\xi_n})=z_n-w_n-j
\end{array}
\bigg|\begin{array}{l}
\wt{X}_\ell^{\xi_n}=a_\ell\\  
\wt{X}_{n-k}^{\xi_n}=b_{n-k}^{(n)}
\end{array}
\right]=\frac{1+o(1)}{\sqrt{2\pi \wt{V}_n^{\xi_n}}}.
$

\noindent
{\em Proof of the Claim.\/}
 By Lemmas \ref{l.uniform-Doeblin} and Corollary \ref{Cor-rho-mixing}, 
$\mathrm{Var}(\wt{S}_\ell^{n-k})=\wt{V}_n^{\xi_n}+O(1)$, and  
by Lemma 7.21 in \cite{Dolgopyat-Sarig-Book}, $\wt{V}_n^{\xi_n}\asymp V_n$. Thus $\mathrm{Var}(\wt{S}_\ell^{n-k})\asymp n$ (Proposition \ref{Prop-Var-Growth}). 

Since $f_i$ are uniformly bounded, $w_n=O(\ell+k)$. Next,  \eqref{Pres1stDer} and the assumption $\PP_n'(\xi_n)=z_n$ lead to 
$$
\E(\wt{S}_\ell^{n-k})=\E(\wt{S}_1^n)+O(
\ell+k)\equiv\wt{\E}^{\xi_n}(S_n)+O(1)=\PP_n'(\xi_n)+O(1)=z_n+O(1).
$$
Thus $\DS \frac{(z_n-w_n-j)-\E(\wt{S}_\ell^{n-k})}{\sqrt{\mathrm{Var}(\wt{S}_\ell^{n-k})}}
=O\left(\frac{\ell+k+|j|}{n}\right)\xrightarrow[n\to\infty]{}0.$
Now we can apply the mixing local limit theorem (Theorem \ref{ThLLTArrays}) to the array

$$
\begin{array}{llllll}
\wt{X}_\ell^{\xi_{\ell+k}} &  &  & & &\text{given $\wt{X}^{\xi_{\ell+k}}_\ell=a_\ell$}\\
\wt{X}_\ell^{\xi_{\ell+k+1}} & \wt{X}_{\ell+1}^{\xi_{\ell+k+1}} & &  & &\text{given $\wt{X}^{\xi_{\ell+k+1}}_\ell=a_\ell$}\\
\wt{X}_\ell^{\xi_{\ell+k+2}} & \wt{X}_{\ell+1}^{\xi_{\ell+k+2}} & \wt{X}_{\ell+2}^{\xi_{\ell+k+2}} & & &
\text{given $\wt{X}^{\xi_{\ell+k+2}}_\ell=a_\ell$}\\
\vdots & \vdots &\vdots &\ddots& & \hspace{1.3cm}\vdots\\
\wt{X}_\ell^{\xi_n} & \wt{X}_{\ell+1}^{\xi_n} & \wt{X}_{\ell+2}^{\xi_n} & \cdots\cdots & \wt{X}_{n-k}^{\xi_n} 
&\text{given $\wt{X}^{\xi_{n}}_\ell=a_\ell$}
\end{array}
$$

The conditions of the theorem hold: Since $\xi_n\to\xi$, there exists some $R>0$ such that $\xi_n\in [-R,R]$ for all $n$.   Uniform ellipticity follows from Lemma \ref{l.uniform-Doeblin}.  
Irreducibility and the hereditary property follow from Lemmas  \ref{LmChange} and \ref{l.hereditary}. Finally, 
$\DS \Prob\left[\wt{X}^{\xi_n}_{n-k}=b^{(n)}_{n-k}\right]$ is bounded below, for the following reasons. 

We saw in the proof of Proposition \ref{Prop-Doeblin}, that there is a constant $\eps_0'>0$ such that for all $i$, $a\in \mathfrak P_i$, and $b\in\mathfrak P_{i+2}$, there exists $c\in\mathfrak P_{i+1}$ such that 
$$
\pi_{i,i+1}(a,c)\pi_{i+1,i+2}(c,b)>(\eps_0')^2.
$$
Looking at \eqref{ProbXi} and Lemma \ref{LmAddPres}, we see that there is a constant $\eps_1>0$ independent of $i,n,a,c$ such that 
$c
\pi_{i,i+1}^{\xi_n}(a_i,a_{i+1})\pi_{i+1,i+2}^{\xi_n}(a_{i+1},a_{i+2})>\eps_1.
$
Thus 
\begin{align}
&\Prob[\wt{X}^{\xi_n}_{n-k}=b^{(n)}_{n-k}]  =
\sum_{a,c}\Prob[\wt{X}^{\xi_n}_{n-k-2}=a, \wt{X}^{\xi_n}_{n-k-1}=c, \wt{X}^{\xi_n}_{n-k}=b_{n-k}^{(n)}]\label{e.state-lower-bound-xi}\\
&=
\sum_{a,c}\Prob[\wt{X}^{\xi_n}_{n-k-2}=a]
\pi_{n-k-2,n-k-1}^{\xi_n}(a,c)\pi_{n-k-1,n-k}^{\xi_n}(c,b_{n-k}^{(n)})\geq \eps_1 \sum_{a}\Prob[\wt{X}^{\xi_n}_{n-k-2}=a]=\eps_1.\notag
\end{align}
This completes the verification of the conditions of {Theorem \ref{t.ThLLT}}. Invoking the theorem, we obtain the claim.

\medskip
In summary,  $p_n(\wt{C})=[1+o(1)]M_\xi(\wt{C})c_n(W_n)$,  where
\begin{align}\label{e.12.6}
& c_n(W_n):=\frac{e^{\PP_n(\xi_n)-\xi_n z_n}}{h_n(b^{(n)}_n,\xi_n)
\sqrt{2\pi \wt{V}_n^{\xi_n}}}\frac{\Prob
[\wt{X}_{n-k}^{\xi_n}=b^{(n)}_{n-k},\ldots,\wt{X}_n^{\xi_n}=b^{(n)}_n]}{\mu_0(Q_n)}
\end{align}
Since $c_n(W_n)$ is independent from $\wt{C}$,  ${p_n(\wt{C}_1)}/p_n(\wt{C}_2)\xrightarrow[n\to\infty]{}{M_\xi(\wt{C}_1)}/{M_\xi(\wt{C}_2)}$.

We proved the first part of the lemma, \eqref{e.engine2}. 

For the second part, 
suppose $W_n'$ is another  sequence of pure $E^s_n$-segments of bounded depth and heights $z_n(W_n')=z_n+O(1)$. There is no loss of generality in assuming that $W_n'$ has fixed depth $k'$. Then $\pi(W_n')$ is a maximal $E^s_n$-segment in $Q_n':=\<c_{n-k'}^{(n)},\ldots,\dot{c}_n^{(n)}\>$. 
To prove \eqref{e.engine2bis}, it is sufficient to show:

\medskip
\noindent
{\em Claim 2.\/} 
 $c_n(W_n')\asymp c_n(W_n)$.  

\medskip
\noindent
{\em Proof of the Claim.\/} 
Let $\xi_n'$ be the solutions to $\PP_n'(\xi_n')=z_n(W_n')$. By Lemma \ref{l.solution-xi}, $|\xi_n'-\xi_n|=o(1/n)$. In particular, $\xi_n',\xi_n\to\xi$. We use the decomposition 
$$
\frac{c_n(W_n')}{c_n(W_n)}=
\underset{\mathrm{I}}{\underbrace{\frac{e^{\PP_n(\xi_n')-\xi_n' z_n}}{e^{\PP_n(\xi_n)-\xi_n z_n}}}}\cdot 
\underset{\mathrm{II}}{\underbrace{
\frac{h_n(b_n^{(n)},\xi_n)\sqrt{2\pi \wt{V}_n^{\xi_n}}}
{h_n(c_n^{(n)},\xi_n')\sqrt{2\pi \wt{V}_n^{\xi_n'}}}}}
\cdot
\underset{\mathrm{III}}{\underbrace{
\frac{\mu_0(Q_n)}{\mu_0(Q_n')}}}
\cdot
\underset{\mathrm{IV}}{\underbrace{
\frac{\Prob
[\wt{X}_{n-k'}^{\xi_n'}=c^{(n)}_{n-k'},\ldots,\wt{X}_n^{\xi_n'}=c^{(n)}_n]}{\Prob
[\wt{X}_{n-k}^{\xi_n}=b^{(n)}_{n-k},\ldots,\wt{X}_n^{\xi_n}=b^{(n)}_n]}}}.
$$

\medskip
$\boldsymbol{\mathrm{I}\asymp 1:}$
Choose some $R>0$ such that $\xi_n,\xi_n'\in [-R,R]$ for all $n$. 
By  \eqref{Pres1stDer}, $|\PP_n'(\xi)|=O(n)$ uniformly on $[-R,R]$, whence:
\begin{itemize}
\item  $z_n=O(n)$, and 
$
|(\xi_n-\xi_n')z_n|=O(n^{-1}n)=O(1).
$
\item  
$|\PP_n(\xi_n)-\PP_n(\xi_n')|=O(n\cdot |\xi_n-\xi_n'|)=O(n n^{-1})=O(1). 
$
\end{itemize}
Thus $|(\xi_n-\xi_n')z_n|+|\PP_n'(\xi_n')-\PP_n(\xi_n)|=O(1)$, and  $\mathrm{I}\asymp 1$. 

\medskip
$\boldsymbol{\mathrm{II}\asymp 1:}$ There are constants such that  $h_n(\cdot,\xi)=\sC(R)^{\pm 1}$ and $\wt{V}^{\xi}_n= \sC(R)^{\pm 1}V_n$ for all $\xi\in [-R,R]$, see Lemmas \ref{LmAddPres} and \ref{l.V-n-xi}.

\medskip
$\boldsymbol{\mathrm{III}\asymp 1:}$ $\mu_0(Q_n),\mu_0(Q_n')\asymp 1$, by Lemma \ref{Lemma-Cylinder-Area} (4).

\medskip
$\boldsymbol{\mathrm{IV}\asymp 1:}$  Choose $R$ such that $\xi_n,\xi_n'\in [-R,R]$ for all $n$, then it is easy to see using Lemma \ref{LmAddPres} that  $\pi_{j,j+1}^\xi(a,b)\asymp \pi_{j,j+1}(a,b)$ for $\xi\in [-R,R]$, uniformly in $j,a,b$.   

We saw above that $\Prob[\wt{X}^{\xi_n}_{n-k}=b^{(n)}_{n-k}]\in (\eps_1,1]$. 
Also $\Prob[X_{n-k}=b^{(n)}_{n-k}]=\mu_0(b^{(n)}_{n-k})\asymp 1$, and $\pi_{j,j+1}^\xi(a,b)\asymp \pi_{j,j+1}(a,b)$ uniformly. So 
\begin{align*}
&\Prob
[\wt{X}_{n-k}^{\xi_n}=b^{(n)}_{n-k},\ldots,\wt{X}_n^{\xi_n}=b^{(n)}_n]=\Prob[\wt{X}^{\xi_n}_{n-k}=b^{(n)}_{n-k}]\prod_{j=n-k}^{n-1}\pi_{j,j+1}^{\xi_n}(b^{(n)}_j,b^{(n)}_{j+1})\\
&\asymp\Prob[X_{n-k}=b^{(n)}_{n-k}]\prod_{j=n-k}^{n-1}\pi_{j,j+1}(b^{(n)}_j,b^{(n)}_{j+1})=\mu_0(Q_n)\asymp 1.
\end{align*}
Similarly, $\Prob
[\wt{X}_{n-k'}^{\xi_n'}=c^{(n)}_{n-k'},\ldots,\wt{X}_n^{\xi_n}=c^{(n)}_n]\asymp \mu_0(Q_n')\asymp 1$.
So $\mathrm{IV}\asymp 1$.

\medskip
It follows that $c_n(W_n')/c_n(W)\asymp 1$.
\end{proof}


Given $\wt{C}\in\mathfs G$, let $\wt{C}^u$ be the bottom $E^u$-side of $\wt{C}$. Given $x\in\wt{C}$,  let $W^s(x,\wt{C})$ be the maximal $E^s_0$-segment in $\wt{C}$ passing through $x$, and let  $\fl|_{W^s(x,\wt{C})}$ denote the length measure on $W^s(x,\wt{C})$. 

A locally finite measure $\nu$ on $\St$ is called  {\em locally flow invariant in $\mathfs G$},\index{Local flow invariance}\index{Locally flow invariant measure} if  for any $\wt{C}\in \mathfs G$, there is a  measure $\zeta_{\wt{C}^u}$ on $\wt{C}^u$ such that 
$$
\nu|_{\wt{C}}=\int_{\wt{C}^u}\fl|_{W^s(x,\wt{C})}d\zeta_{\wt{C}^u}(x).
$$
Since $E^s_0$ is the direction of the flow   $\phi$, $\wt{C}^u$ is a transverse to the flow direction, and $\wt{C}$ is a flow box; therefore any $\phi$-invariant measure is locally flow invariant in $\mathfs G$. The measure $\zeta_{\wt{C}^u}$  describes the flux through $\wt{C}^u$. 

\begin{lem} \label{l.Extending-M-xi}
$M_\xi$ has a unique extension to a locally finite  measure $\nu_\xi$, which is locally flow invariant on $\mathfs G$. \index{Invariant Radon measures!$\mu_\xi$}
\end{lem}
\begin{proof}
Fix  $\wt{C}\in\mathfs G$ of the form $(D^j\circ\pi|_{\wt{F}_0}^{-1})\<\dot{a}_0,\ldots,a_n\>$, and let $\wt{C}^u$ denote the bottom $E^u_0$-side of $\wt{C}$. Given $
W:=\wt{C}^u\cap (D^j\circ\pi|_{\wt{F}_0}^{-1})\<\dot{a}_0,\ldots,a_n; a_{n+1},\ldots,a_{n+\ell}\>
$, let 
$$
\zeta_{\wt{C}^u}(W):=\frac{e^{j}h_0(a_0,\xi)}{\ell^s(a_0)}\Prob[\wt{X}^\xi_0=a_0,\ldots, \wt{X}^\xi_{n+\ell}=a_{n+\ell}].
$$
Each $W$ is a sub-segment of $\wt{C}^u$ with length $O(\lambda^{-(n+\ell)})$, and the collection of  $W$s is a semi-algebra  on $\wt{C}^u$, which --- since $\ell$ is arbitrary --- generates the Borel $\sigma$-algebra on $\wt{C}^u$. By the Carath\'eodory's extension theorem, $\zeta_{\wt{C}^u}$ extends to a Borel measure on $\wt{C}^u$.

Set 
$\DS \nu_{\xi,\wt{C}}:=\int_{\wt{C}^u}\fl|_{W^s(x,\wt{C})}\, d\zeta|_{\wt{C}^u}(x).
$ 
We claim that there is a measure $\nu_\xi$ on $\St$ such that $\nu_{\xi,\wt{C}}=\nu_\xi|_{\wt{C}}$ for all $\wt{C}$. If  
 $\wt{C},\wt{D}\in\mathfs G$ intersect with positive measure, then either $\wt{C}\subset\wt{D}$ or
 $\wt{D}\subset \wt{C}$, and  $j(\wt{D})=j(\wt{C})$, $|D^s|=|C^s|$, and  $a_0(\wt{C})=a_0(\wt{D})$. So
 $\zeta_{\wt{D}^u}|_{\wt{C}\cap\wt{D}}=\zeta_{\wt{C}^u}|_{\wt{C}\cap\wt{D}}$,  and $(\nu_{\xi,{\wt{C}}})|_{\wt{C}\cap \wt{D}}=\nu_{\xi,\wt{D}}|_{\wt{C}\cap\wt{D}}$.  Thus  $\nu_{\xi,\wt{C}}$ $(\wt{C}\in\mathfs G)$ have a common extension $\nu_\xi$. Clearly, $\nu_{\xi}$ is locally finite, locally flow invariant on $\mathfs G$, and  $\nu_\xi(\wt{C})=M_\xi(\wt{C})$ for all $\wt{C}\in\mathfs G$.

Every extension of $M_\xi$ to a measure which is  locally flow invariant on $\mathfs G$ must induce the same flux measure on $\wt{C}^u$, and is therefore  equal to $\nu_\xi$.
\end{proof}

\begin{lem}\label{l.Lemma-nu-xi-is-BL-measure}
Fix $T_k\to\infty$ and $\wt{\omega}\in\St$ whose forward $\phi$ orbit
does not hit a singularity, and suppose there is a locally finite  measure $\nu$ which is locally flow-invariant in $\mathfs G$, and so that for any  pair $(h_1,h_2)$ of indicator functions of sets in $\mathfs G$, 
$\int h_i d\nu\neq 0$, and 
\begin{equation}\label{e.RET}
\frac{\int_{0}^{T_k} h_1(\phi^t(\wt{\omega}))\,dt}{\int_{0}^{T_k} h_2(\phi^t(\wt{\omega}))\,dt}\xrightarrow[T\to\infty]{}\frac{\int h_1 d\nu}{\int h_2 d\nu};
\end{equation}
\begin{equation}\label{e.BL-eqn}
\nu(h_i\circ D)=e^{\xi}\nu(h_i).
\end{equation}
Then   $\nu=\const\mu_\xi$ (whence $\nu$ is $\phi$-ergodic and invariant), and  
\eqref{e.RET}--\eqref{e.BL-eqn} hold for any pair of continuous functions with compact support such that $\int h_2 d\nu>0$. 
\end{lem}
\noindent
\begin{proof} 
Suppose $\wt{C}\in\mathfs G$. By \eqref{e.RET}, the orbit of $\tomega$
enters $\wt{C}$ and $\wt{C}^c$ infinitely many times. Therefore it crosses $\wt{C}$ infinitely many times. Each crossing takes a constant time equal to the length of the $E^s_0$-side of $\wt{C}$. 

Thus, if $h=1_{\wt{C}}$, then 
\begin{equation}\label{e.RET-conservative}
\int_{0}^\infty h(\phi^t(\wt{\omega}))\,dt=\infty.
\end{equation}

By Lemma \ref{Lemma-Cylinder-Area}, the elements of  $\mathfs G_k$ are rectangles   with ``long``  $E^s_0$-side of length $\asymp 1$, and ``short`` $E^u_0$-side with length $O(\lambda^{-k})$. 
Call a piece of $\wt{C}$ cut by maximal $E^u_0$ segments in $\wt{C}$ an {\em $E^u_0$-slice}. Slices of elements of $\mathfs G_k$ are rectangles with $E^u_0$-sides of length $O(\lambda^{-k})$, and $E^s_0$-sides  which could be arbitrarily narrow.
Thus 
$
\mathfs H:=\{\text{$E^u_0$-slices of sets in $\mathfs G$}\}
$ generates the Borel $\sigma$-algebra. 

Let  $\wt{D}$ be an $E^u$-slice of $\wt{C}\!\in\!\mathfs G$. 
 $E^s_0$ is the direction of the flow, and geometric considerations show that  
$
\DS \int_0^{T_k}\! 1_{\wt{D}}(\phi^t(\wt{\omega}))dt=\rho
{\int_0^{T_k} 1_{\wt{C}}(\phi^t(\wt{\omega}))dt}+O(1)$,  where $\DS\rho=\frac{\ell^s(\wt{D})}{\ell^s(\wt{C})}$.
By the local flow invariance, $\DS \rho=\frac{\nu(\wt{D})}{\nu(\wt{C})}$, and by \eqref{e.RET-conservative}, 
$
\DS \frac{\int_{0}^{T_k} 1_{\wt{D}}(\phi^t(\wt{\omega}))\,dt}{\int_{0}^{T_k} 1_{\wt{C}}(\phi^t(\wt{\omega}))\,dt}\xrightarrow[T\to\infty]{}\frac{\nu(\wt{D})}{\nu(\wt{C})}.
$

Repeating this argument for other slices of other elements of $\mathfs G$, and employing \eqref{e.RET} for pairs of indicators of  $\mathfs G$-elements, we deduce that  \eqref{e.RET} holds for all pairs of indicators of {\em slices} of   $\mathfs G$-elements. 

Since one can cover $\St^\ast$ by slices of $\mathfs G$-elements of { uniformly small diameter}, a standard approximation argument gives \eqref{e.RET} and  \eqref{e.RET-conservative} for all pairs of continuous functions with compact support and positive integrals. 

By \eqref{e.RET-conservative},  $\nu$ must be $\phi$-invariant.
The ergodic components of $\nu$ are all locally finite  measures.
By \cite{Hooper-Hubert-Weiss}, every $\phi$-ergodic invariant locally finite measure is proportional to $\mu_\xi$. Thus, by the ergodic decomposition theorem for conservative $\sigma$-finite invariant measures, there is a probability  measure $\eta$ on $\R$ such that 
$
\nu=\int_\R \mu_{\xi'}d\eta(\xi').
$
Recall that $D$ is the deck transformation which increases the $\zc$-coordinate by one.  If $\nu\circ D=e^{\xi}\nu$ on $\mathfs G$, then $\nu\circ D^j=e^{j\xi}\nu$ on $\mathfs G$ for all $j\in\Z$ (because $\mathfs G$ is $D$-invariant). Thus for every $\wt{C}\in\mathfs G$, $e^{j\xi}\nu(\wt{C})=\int_\R e^{j\xi'}\mu_{\xi'}(\wt{C})d\eta(\xi')$ for all $j\in\Z$. Necessarily, $\eta$ is the point mass at $\xi$,  $\nu$ is proportional to $\mu_\xi$, and $\nu$ is ergodic.  
\end{proof}

\begin{proof}[\bf Proof of Proposition \ref{PrBulk}] 
Recall that $\wt{\psi}_{n-1}\circ\cdots\circ\wt{\psi}_0=\wt{\Psi}_{nN_0}$. 
Suppose $\xi\in\R$, $\eta_N\to\xi$, and  $\mathcal P_{N}'(\eta_N)=z_{N}(\tomega)$. Specializing to the subsequence $N_n=nN_0$ and applying  Proposition \ref{p.Pressure-and-Geometric-Pressure}, we obtain finite numbers $\xi_n\to\xi$ such that $\PP_n'(\xi_n)=z_{nN_0}(\tomega)$.

Fix $0<\epsilon<0.01$, $k\geq 1$, and a  constant $L_0$ so large that for every $n$, any pure $E^s_n$-segment of depth $k$ has length less than $\epsilon L_0$. This is possible by \eqref{e.c-over}.

Let $A_T:=\{\phi^t(\wt\omega):0<t<T\}$.  Arguing as in the proof of the temporal CLT in Section \ref{ScProofBeck}, there is a  global constant $c>4$ 
and natural numbers $n=n(T)\asymp \log T$ $(T>0)$, such that
\begin{equation}\label{e.n(T)}
A_T^\ast:=(\wt{\psi}_{n-1}\circ\cdots\circ\wt{\psi}_0)(A_T) 
\text{  has length in $[L_0,cL_0]$.  }
\end{equation}
Since $|A_T^\ast|\leq cL_0$ is bounded and $|\wt{\zc}_n-\zc|$ is uniformly bounded, 
$$
\wt{\zc}_n(\cdot)=z_{nN_0}(\tomega)+O(1)\text{ uniformly on }A_T^\ast.
$$

We denote the uniform distribution on $A_T^\ast$ by $\Prob_{A^\ast_T}$.  Since $\wt{\psi}_{n-1}\circ\cdots\circ\wt{\psi}_0:A_T\to A_T^\ast$ is  affine, for every $\wt{C}\in\mathfs G$, 
\begin{align}
&\frac{1}{T}\int_0^T 1_{\wt{C}}(\phi^t(\wt{\omega}))dt=\frac{|A_T\cap \wt{C}|}{|A_T|}=\frac{|A_T^\ast\cap (\wt{\psi}_{n-1}\circ\cdots\circ\wt{\psi}_0)(\wt{C})|}{|A_T^\ast|}\notag
\\
&=\Prob_{A_T^\ast}\left\{\wt{\omega}'\in A_T^\ast:  (\wt{\psi}_0^{-1}\circ\cdots\circ\wt{\psi}_{n-1}^{-1})(\wt{\omega}')\in\wt{C}
\right\}.\label{e.sirena}
\end{align}

By Lemma \ref{l.approx-by-pure-segments}, we can approximate $A_T^\ast$  by finite unions of pure $E^s_n$-segments $W_{n,i}$ of  depth  $k$, which intersect $A_T^\ast$, and so that 
$$
\DS 
\bigcup_{i=1}^{m_n}W_{n,i}\subset A_T^\ast\subset \bigcup_{i=1}^{m_n+2}W_{n,i}\ ; \  \frac{1}{|A_T^\ast|}\sum_{i=1}^{m_n}|W_{n,i}|\Prob_{W_{n,i}}(\cdot)\leq \Prob_{A_T^\ast}(\cdot)\leq  \frac{1}{|A_T^\ast|}\sum_{i=1}^{m_n+2}|W_{n,i}|\Prob_{W_{n,i}}(\cdot)
.
$$
Since $W_{n,i}\cap A_T^\ast\neq \emptyset$ and $|W_{n,i}|=O(1)$, their  heights satisfy $z_n(W_{n,i})=z_{nN_0}(\tomega)+O(1)$.

Thus, for any pair of sets $\wt{C}_1,\wt{C}_2\in\mathfs G$, 
\begin{align*}
&\frac{\int_0^T 1_{\wt{C}_1}(\phi^t(\wt{\omega}))dt}
{\int_0^T 1_{\wt{C}_2}(\phi^t(\wt{\omega}))dt}=
\frac{\Prob_{A^\ast_T}[(\wt{\psi}_{n-1}\circ\cdots\circ\wt{\psi}_0)(\wt{C}_1)]}{\Prob_{A^\ast_T}[(\wt{\psi}_{n-1}\circ\cdots\circ\wt{\psi}_0)(\wt{C}_2)]}\leq 
\frac{\sum\limits_{i=1}^{m_n+2}|W_{n,i}|\Prob_{W_{n,i}}[(\wt{\psi}_{n-1}\circ\cdots\circ\wt{\psi}_0)(\wt{C}_1)]}{\sum\limits_{i=1}^{m_n}|W_{n,i}|\Prob_{W_{n,i}}[(\wt{\psi}_{n-1}\circ\cdots\circ\wt{\psi}_0)(\wt{C}_2)]}\\
&\leq [1+o(1)]\frac{M_\xi(\wt{C}_1)}{M_\xi(\wt{C}_2)}\left[1+\frac{\sum_{i=m_n+1}^{m_n+2}|W_{n,i}|\Prob_{W_{n,i}}[(\wt{\psi}_{n-1}\circ\cdots\circ\wt{\psi}_0)(\wt{C}_2)]}{\sum_{i=1}^{m_n}|W_{n,i}|\Prob_{W_{n,i}}[(\wt{\psi}_{n-1}\circ\cdots\circ\wt{\psi}_0)(\wt{C}_2)]}
\right],\text{ by \eqref{e.engine2}}\\
&\leq [1+o(1)]\frac{M_\xi(\wt{C}_1)}{M_\xi(\wt{C}_2)}\left[1
+O\left(\frac{\sum_{i=m_n+1}^{m_n+2}|W_{n,i}|}{\sum_{i=1}^{m_n}|W_{n,i}|}\right)
\right],\text{ by \eqref{e.engine2bis}.}
\end{align*}
By the choice of $L_0$,  
$\DS \frac{\sum_{i=m_n+1}^{m_n+2}|W_{n,i}|}
{\sum_{i=1}^{m_n}|W_{n,i}|}\leq \frac{2\epsilon L_0}{L_0-2\epsilon L_0}<\frac{2\epsilon}{1-2\epsilon}
<3\epsilon.$ So 
$$
\limsup_{T\to\infty} \frac{\int_0^T 1_{\wt{C}_1}(\phi^t(\wt{\omega}))dt}
{\int_0^T 1_{\wt{C}_2}(\phi^t(\wt{\omega}))dt}\leq \frac{M_\xi(\wt{C}_1)}{M_\xi(\wt{C}_2)}[1+O(\epsilon)]=\frac{\nu_\xi(\wt{C}_1)}{\nu_\xi(\wt{C}_2)}[1+O(\epsilon)],
$$
with $\nu_\xi$ as in Lemma \ref{l.Extending-M-xi}.
By symmetry, $\DS \liminf_{T\to\infty} \frac{\int_0^T 1_{\wt{C}_1}(\phi^t(\wt{\omega}))dt}
{\int_0^T 1_{\wt{C}_2}(\phi^t(\wt{\omega}))dt}\geq \frac{\nu_\xi(\wt{C}_1)}{\nu_\xi(\wt{C}_2)}[1-O(\epsilon)]$, 
and since $\epsilon$ is arbitrary,
$\DS
\lim_{n\to\infty}\frac{\int_0^T 1_{\wt{C}_1}(\phi^t(\wt{\omega}))dt}
{\int_0^T 1_{\wt{C}_2}(\phi^t(\wt{\omega}))dt}=\frac{\nu_\xi(\wt{C}_1)}{\nu_\xi(\wt{C}_2)},
$ for all $\wt{C}_1,\wt{C}_2\in\mathfs G$. 

This shows that $\nu_\xi$ satisfies condition \eqref{e.RET} of Lemma \ref{l.Lemma-nu-xi-is-BL-measure}. 
Condition \eqref{e.BL-eqn} holds as well, because $\nu_\xi|_{\mathfs G}=M_\xi$ and $M_\xi\circ D=e^\xi M_\xi$ by construction.
By Lemma \ref{l.Lemma-nu-xi-is-BL-measure}, $\wt{\omega}$ is generic for $\mu_\xi$.
This proves the first part of Proposition \ref{PrBulk}. 

\medskip
Suppose there are integers $N_k\to\infty$ and finite numbers  $\xi,\eta_{N_k}\in\R$  such that $\eta_{N_k}\to\xi$ and $\mathcal P_{N_k}'(\eta_{N_k})=z_{N_k}(\tomega)$. Let $N_0$ be the constant from \eqref{e.mat-prod-form}, and set $n_k:=\lfloor{N_k/N_0}\rfloor$. By Proposition \ref{p.Pressure-and-Geometric-Pressure}, there are  $\xi_{n_k}\to\xi$ such that $\PP_{n_k}'(\xi_{n_k})=z_{N_k}(\tomega)$. 

Build the sequence of times $T_k\to\infty$ such that the length of 
$(\wt{\psi}_{n_k-1}\circ\cdots\circ\wt{\psi}_0)(A_{T_{n_k}})$ is between $L_0$ and $cL_0$. Proceeding as before, we  can show that 
$$
\lim_{n\to\infty}\frac{\int_0^{T_k} 1_{\wt{C}_1}(\phi^t(\wt{\omega}))dt}
{\int_0^{T_k} 1_{\wt{C}_2}(\phi^t(\wt{\omega}))dt}=\frac{\nu_\xi(\wt{C}_1)}{\nu_\xi(\wt{C}_2)}
\text{ for all $\wt{C}_1,\wt{C}_2\in\mathfs G$. }
$$
By Lemma \ref{l.Lemma-nu-xi-is-BL-measure}, $\nu_\xi$ is a historical measure of $\tomega$. 
\end{proof}

\subsection{A Necessary Condition for Genericity}
\label{ScEdge}
In this section we prove:
\begin{prop}
\label{PrEdge} If $\xi_N(\wt{\omega})$ does not converge to a finite limit as $N\to\infty$, then 
$\tomega$ is not generic for any $\phi$-ergodic invariant locally finite measure on $\St.$
\end{prop}

We already know from Proposition \ref{PrBulk} that if $\xi_N(\tomega)$ has different finite limit points, then $\tomega$ has several historical measures, and is not generic for any measure. It remains to consider the case when $\xi_N(\tomega)$ has an infinite limit point. 
We need the following lemma:

\begin{lem}\label{l.Pure-Seg-Ratio-Edge-Case}
For every $Z_0>0$ there exists a constant $M_0$ as follows. 
Suppose
 $n_s$ is a strictly increasing sequence of positive integers, and 
\begin{enumerate}[(a)]
\item 
$z_{n_s}$ is a sequence of integers such that for all $R>0$, for all $s$ large enough, either $z_{n_s}=\PP_{n_s}'(\xi_{n_s})$ with $\xi_{n_s}>R$, or  $z_{n_s}\geq \lim\limits_{t\to +\infty}\PP_{n_s}'(t)$;
\item $W_{n_s}, \overline{W}_{n_s}$ are two sequences of pure $E^s_{n_s}$-segments of bounded depths, and with heights $z_{n_s}(W_{n_s})$, $z_{n_s}(\ov{W}_{n_s})$  such that  $|z(W_{n_s})-z_{n_s}|, 
|z(\overline{W}_{n_s})-z_{n_s}|\leq Z_0$;  
\item $\Prob_{W_{n_s}}
[ (\wt{\psi}_{n_s-1}\circ\cdots\circ\wt{\psi}_0)(\wt{F}_0)]$ is positive for all $s$ large, where $\wt{F}_0=[\wt{\zc}_0=0]$.
\end{enumerate}
 Then $\DS \Prob_{\ov{W}_{n_s}}\left[(\wt{\psi}_{n_s-1}\circ\cdots\circ\wt{\psi}_0)
 \left(\bigcup_{\ell=-M_0}^{M_0}D^{-\ell}\wt{F}_0\right)\right]$ is positive for all $s$ large, and 
 \begin{equation}\label{e.P-ratio}
 \lim_{k\to\infty}\frac{\Prob_{W_{n_s}}[
 (\wt{\psi}_{n_s-1}\circ\cdots\circ\wt{\psi}_0)(\wt{F}_0)]}
 {\Prob_{\ov{W}_{n_s}}[(\wt{\psi}_{n_s-1}\circ\cdots\circ\wt{\psi}_0)
 (\bigcup_{\ell=-M_0}^{M_0}D^{\ell}\wt{F}_0)]}=0.
 \end{equation}
The same conclusion holds if (a) is replaced by 
\begin{enumerate}[(a')]
\item
$z_{n_s}$ is a sequence of integers such that for all $R>0$, for all $s$ large enough, either 
$z_{n_s}=\PP_{n_s}'(\xi_{n_s})$ with $\xi_{n_s}<-R$, or  $z_{n_s}\leq \lim\limits_{t\to -\infty}\PP_{n_s}'(t)$.
\end{enumerate}
\end{lem}

We assume for simplicity that  $(n_s)$ is the sequence $(1,2,3,\ldots)$, and write $n_s=n$.  The routine extension to the general case is left to the reader.

We give the proof assuming (a), and comment on the case  (a') at the end of the proof. 
There is no loss of generality in assuming that  $W_n$ have constant depth $k$, and $\ov{W}_n$ have constant depth $\ov{k}$. By the definition of pure $E^s_n$-segments, $\pi(W_n)$ is a maximal $E^s_n$-segment in some cylinder $\<a_{n-k}^{(n)},\ldots,\dot{a}_n^{(n)}\>$, and $\pi(\ov{W}_n)$ is  a  maximal $E^s_n$-segment in some cylinder $\<b_{n-\ov{k}}^{(n)},\ldots,\dot{b}_n^{(n)}\>$. Let 
\begin{equation}\label{e.P-Q-words}
\un{a}^{(n)}:=(a^{(n)}_{n-k},\ldots,a^{(n)}_n), \quad \un{b}^{(n)}:=(b^{(n)}_{n-\bar k},\ldots,b^{(n)}_n).
\end{equation}

Let $\{X_i\}_{i\geq 0}$ denote the  Markov chain from Proposition \ref{Prop-Markov-Kernel}, and let $f_k=f_k(X_k)$ be the functions constructed in Lemma \ref{Lemma-f(X)}. Recall that 
$\DS K:=\sup_{n} \|f_n\|_\infty<\infty.$
Set
\begin{align}
&\Sigma^n:=\biggl\{(x_0,\ldots,x_n)\in\prod_{i=0}^n\mathfrak P_i: \Prob\biggl((X_0,\ldots,X_n)=\un{x}\bigg|(X_{n-k},\ldots,X_n)=\un{a}^{(n)}
\biggr)\neq 0
\biggr\};\notag\\
&\Omega^n:=\{(x_i)_{0\leq i\leq n}\in \Sigma^n:\sum_{i=0}^n f_i(x_i)=\wh{z}_n\}\ , \text{ where }   \wh{z}_n:=z_n(W_n). \label{e.Omega-n-Def}
\end{align}
$\Sigma^n, \Omega^n\neq \emptyset$  for all $n$ sufficiently large, because by  (c), $\wt{F}_0\cap (\wt{\psi}_{n-1}\circ\cdots\circ\wt{\psi}_0)^{-1}(W_n)\neq \emptyset$, and    
 $
\Omega^n\cap \Sigma^n\owns(X^{(n)}_0(\pi(\omega')),\ldots,X^{(n)}_n(\pi(\omega'))),
 $ 
for any $\omega'\in \wt{F}_0\cap (\wt{\psi}_{n-1}\circ\cdots\circ\wt{\psi}_0)^{-1}(W_n)$, see \eqref{e.X-bar}, Proposition \ref{Prop-Markov-Kernel},  and Lemma \ref{Lemma-f(X)}.

\setcounter{sublem}{0}
\begin{sublem}\label{sl.blocks}
For all $L$ large enough, $\{0,\ldots, n\}$ splits into  segments $I_1,\ldots,I_{\nu_n}$, 
listed in increasing order, such that
\begin{enumerate}[(1)]
\item For all $p\neq \nu_n$, $\DS L\leq \mathrm{Var}\biggl(\sum_{n\in I_p}f_n(X_n)\biggr)\leq 2L$; 
\item For all $p$, $\DS 
 \bigg|\Cov\biggl(\sum_{n\in I_p} f_n(X_n), \sum_{m\not\in I_p} f_m(X_m)\biggr)\bigg|\leq 0.01 L$;
\item 
$\DS\sup_n ({\max_{1\leq p\leq \nu_n}} |I_p|/L)<\infty$;  
$\DS\inf_n ({ \min_{1\leq p\leq \nu_n-1}} |I_p|/L)>0$; and  $\nu_n\asymp n$.
 \end{enumerate}
\end{sublem}
\begin{proof}
The proof relies on the  exponential $\rho$-mixing property of $\{X_i\}_{i\geq 1}$ (Corollary \ref{Cor-rho-mixing}):  There exist constants $C_{mix}>1$ and $0<\rho<1$ such that for all $i,j$, 
$$
|\mathrm{Cov}(f_i(X_i),f_j(X_j))|\leq C_{mix}K^2\rho^{|i-j|}.
$$ 

Given integers $0<m<n$, let  $\DS V_{m,n}:=\mathrm{Var}\bigl(\sum_{m< i\leq n}f_i(X_i)\bigr)$ and $V_k:=V_{1,k}$. Then 
$\DS
V_{m,n}\geq V_n-V_m-\sum_{i=1}^{m}
\sum_{j=m+1}^n
C_{mix}K^2\rho^{j-i}\geq V_n-V_m-\const.
$ By  Proposition \ref{Prop-Var-Growth}, $V_n\asymp n$. 
Fixing $m$ and letting $n\to\infty$, we obtain
$V_{m,n}\xrightarrow[n\to\infty]{}\infty$ for each $m$. 

At the same time, 
$\DS
|V_{m,j+1}-V_{m,j}|\leq 2\sum_{i=m}^{j+1} |\mathrm{Cov}(f_{i}(X_{i}),f_{j+1}(X_{j+1}))|\leq \frac{C_{mix}K^2}{1-\rho} 
$. 
So 
\begin{equation}
\label{VmnIneq}
V_{m,n}\leq C_{mix}K^2(1-\rho)^{-1}(n-m).
\end{equation}

Suppose  $L>200C_{mix}K^2(1-\rho)^{-2}$. Then $V_{1,2}\leq L$, and by the previous  estimates there is some $n_1$ such that $L\leq V_{1,n_1}\leq 2L$. We let $I_1:=\{1,2,\ldots,n_1\}$. Similarly $V_{n_1+1,n_1+2}\leq L$, and there is $n_2>n_1$ such that $L\leq V_{n_1+1,n_2}\leq 2L$. Proceeding by induction, we obtain a partition of $\{0,1,\ldots,n\}$ into increasing integer segments $I_1,\ldots,I_{\nu_n}$ which satisfy the first part  of the sublemma. 
 
The second part holds as well,  because  if $I_p=\{a,a+1,\ldots,b\}$, then
\begin{align*}
&\bigg|\mathrm{Cov}\biggl(\sum_{n\in I_p} f_n(X_n),\sum_{n\not\in I_p} f_n(X_n)\biggr)\bigg|\leq C_{mix}K^2\sum_{n=a}^b\biggl(\sum_{m=1}^{a-1}\rho^{n-m}+\sum_{m=b+1}^{\infty}
\rho^{m-n}\biggr)\\
&\leq 2C_{mix}K^2(1-\rho^{-1})^2<0.01 L.
\end{align*} 

 Let $|I_p|$ denote the ``length" (i.e. cardinality) of $I_p$, and let  
$\DS V(I_p):=\mathrm{Var}\bigl(\sum_{j\in I_p}f(X_j)\bigr)$. 
 By the variance growth theorem for uniformly elliptic Markov chains 
 (Theorem \ref{LmVarCycles})
 there are constants $C_1,C_2$ (depending only on $K$ and the ellipticity constant $\eps_0$ from Proposition~\ref{Prop-Doeblin}) such that 
$$
V(I_p)\geq C_1^{-1}\sum_{j\in I_p}u_j^2-C_2,\text{ where $u_j$ are the structure constants in \eqref{e.u}. }
$$
We saw in the proof of Lemma \ref{Lemma-D_N-Growth} that there exists a constant $m_0>0$ (independent of $n$ and $p$) so that for all $j$ there is $j\leq j'\leq j+m_0$ and an admissible hexagon at position $j'+2$ with non-zero balance (see p. \pageref{page-hexagon-claim}). The balance is an integer, so by Lemma \ref{lemma-Hex-Meas}, there is a positive constant $c_0$ such that $u_{j'}^2\geq c_0$.
Since the gaps between $j'$s such that $u_{j'}^2\geq c_0$ are of size less than $m_0$, $\sum_{j\in I_p}u_j^2\geq c_0(|I_p|-m_0)/m_0$, whence 
$$
2L\geq V(I_p)\geq C_1^{-1}c_0\frac{|I_p|-m_0}{m_0}-C_2.
$$
Taking $L$ much larger than $C_2+1$, we obtain $|I_p|\leq \const L$. 

By construction,  for all $p\!\!<\!\!\nu_n$, $V(I_p)\!\!\geq\!\! L$. 
By \eqref{VmnIneq}
$V(I_p)\leq C_{mix}K^2(1-\rho)^{-1}|I_p|$. So for all segments except perhaps the last one, 
$
|I_p|\geq \const L.
$
All constants are independent of $L$ and $n$.

It immediately follows that for fixed $L$, the number of segments is bounded above and below by constants time $n$.
\end{proof}

Let 
$\DS 
\Sigma^n(I_p)\!:=\!\{(x_i)_{i\in I_p}\!\in \!\prod_{i\in I_p}\mathfrak P_i\!:\! \un{x}\in\Sigma^n\}.
$
A block $(x_i)_{i\in I_p}\!\in\! \Sigma^n(I_p)$  is called {\em optimal}, if \index{Optimal block}\index{Sub-optimal block}
\begin{equation}\label{e.def-of-optima}
\sum_{i\in I_p} f_i(x_i)\geq \sum_{i\in I_p} f_i(y_i) \quad \begin{array}{l}\text{for all $(y_i)_{i\in I_p}\in\Sigma^n(I_p)$ such that $y_j=x_j$}\\
\text{at the endpoints $j=\min(I_p), j=\max(I_p)$},
\end{array}
\end{equation}
and {\em sub-optimal}, otherwise. 

\begin{sublem}\label{sl.sub-optimal-exist}
 If $L$ is sufficiently large, then   for each $p<\nu_n$, for each optimal $I_p$  there exists
  a sub-optimal block with the same endpoints as $I_p$.
\end{sublem}
\begin{proof}
 Let $m_0$ be an  upper bound  for the gaps in the set  
 $$\{j\geq 3: \exists\text{  hexagons in $\mathrm{Hex}(j)$ with non-zero balance}\}$$ (see the proof of Claim 2 on page \pageref{page-hexagon-claim}).  By Sublemma \ref{sl.blocks}, the length of all the segments except perhaps the last one is bounded below by a constant times $L$. We choose $L$ so large, that $|I_p|>m_0+100$ for all $p$ except perhaps the last one.

Suppose $I_p=\{a,a+1,\ldots,b\}$, and  $(x_a,\ldots,x_b)$ is optimal.
 Since $|I_p|>m_0+100$, 
 there exists some $\DS a+4\leq i\leq b-3$,
and some hexagon $\left(y_{i-2} \begin{array}{cc} y_{i-1} & y_i\\
z_{i-1} & z_i
\end{array} z_{i+1}\right)\in \mathrm{Hex}(i)$ with non-zero balance. 
By  uniform ellipticity (Proposition \ref{Prop-Doeblin})
$$
\exists b_{i-3}\text{ such that } \Prob[(X_{i-4},X_{i-3},X_{i-2})=(x_{i-4},b_{i-3},y_{i-2})]\neq 0,$$
$$\exists b_{i+2}\text{ such that } \Prob[(X_{i+1},X_{i+2},X_{i+3})=(z_{i+1},b_{i+2},x_{i+3})]\neq 0.$$ 
By the Markov property, the following two blocks belong to $\Sigma(I_p)$:
$$
\un{y}:=(x_a,\ldots,x_{i-4};b_{i-3};\ \ \ y_{i-2},y_{i-1},y_i,z_{i+1}\ \ \ ;  
 b_{i+2}; x_{i+3},\ldots,x_b);
$$
$$
\un{z}:=(x_a,\ldots,x_{i-4} ;b_{i-3};\ \ \ y_{i-2},z_{i-1},z_i,z_{i+1}\  \ \ ; 
 b_{i+2} ; x_{i+3},\ldots,x_b).
$$
Since the balance of $\left(y_{n-2} \begin{array}{cc} y_{n-1} & y_n\\
z_{n-1} & z_n
\end{array} z_{n+1}\right)$ is not zero,  the sums $\DS \sum_{i\in I_p}f_i$ along $\un{y}$ and $\un{z}$ are different. Therefore, at least one of $\un{y},\un{z}$ is sub-optimal. 
\end{proof}

Recall the definition of  $\Omega^n$ 
from \eqref{e.Omega-n-Def}.

\begin{sublem}\label{sl.ManyMax}Given $\un{x}\in\Omega^n$, let  $e_n(\un{x})$ denote the number of segments $I_p$ inside $[0, n-k]$ such that $(x_i)_{i\in I_p}$ is sub-optimal. 
Then 
\begin{equation}
\label{ManyMax}
 \lim_{n\to\infty} \frac{1}{n}\max \{e_n(\un{x}):\un{x}\in\Omega^n\}=0.
\end{equation}
\end{sublem}
\begin{proof}Assume by contradiction that \eqref{ManyMax} is false. Then  for some $\delta>0$,  there are $n_j\to\infty$ and  
$\un{x}^{(j)} \in \Omega^{n_j}$ such that
${e_{n_j}(\un{x}^{(j)})}>\delta{n_j}.$

We replace each sub-optimal $(x_i^{(j)})_{i\in I_p}$ in $[0,n_j-k]$ by an optimal block $(y_i)_{i\in I_p}$ with the same endpoints. The result is a word $\un{y}\in\Sigma^n$ such that for all $j$ large enough,
$$
\sum_{i=0}^{n_j} f_i(y_i)\geq \sum_{i=0}^{n_j} f_i(x_i^{(j)})+e_{n_j}(\un{x}^{(j)})\geq \wh{z}_{n_j}+\delta n_j\equiv z_{n_j}(W_{n_j})+\delta n_j\equiv z_{n_j}+\delta n_j+O(1). 
$$
 Therefore, for all $j$ large enough,
 $\DS \sum_{i=0}^{n_j} f_i(y_i)>z_{n_j}+\frac{1}{2}\delta n_j$.
Thus there exists an integer $m_j\geq \frac{1}{2}\delta n_j$ such that 
$\DS \Prob(S_{n_j}=z_{n_j}+m_j)\!\neq \!0 $. 
Since the non-zero transition probabilities of $\{X_i\}$ are bounded away from zero 
(Proposition \ref{Prop-Markov-Kernel}), there exists  $\varkappa>0$ such that 
$$\DS \Prob(S_{n_j}=z_{n_j}+m_j)\geq e^{-\varkappa n_j}.$$

Fix  $R>6\varkappa/\delta$. The function $F_n(t):=\log \E(e^{t S_n})$ is  convex,  analytic,  $F_n(0)=0$, and  $F_{n_j}(R)\geq \log(e^{R(z_{n_j}+m_j)}\Prob[S_{n_j}=z_{n_j}+m_j])>R(z_{n_j}+m_j)-\varkappa n_j\geq R \bigl(z_{n_j}+\frac{1}{3}\delta n_j\bigr)
 $, whence by convexity, 
$$F_{n_j}'(R)\geq \frac{F_{n_j}(R)-F_{n_j}(0)}{R}\geq z_{n_j}+\frac{1}{3}\delta n_j.$$
Employing \eqref{Pres1stDerLM} on the interval $[-R,R]$, we find that for all $j$ sufficiently large, 
$
\PP_{n_j}'(R)\!\!\geq \!\! z_{n_j}.
$
 By strict convexity, $z_{n_j}\!\!<\!\! \limsup\limits_{t\to\infty}P_{n_j}'(t)$,  whence by condition (a), 
 for all $j$ large enough, there exists $\xi_{n_j}$  with $z_{n_j}\!\!=\!\!\PP_{n_j}'(\xi_{n_j})$.  Now the strict monotonicity of $\PP_{n_j}'$ implies that $\xi_{n_j}\!\!<\!\!R$ contradicting condition (a).
This contradiction proves \eqref{ManyMax}.
\end{proof} 
 
\medskip
Recall the word $\un{b}^{(n)}$ from \eqref{e.P-Q-words}, and set  
\begin{align*}
&\wh{\Sigma}^n:=\biggl\{(x_0,\ldots,x_n)\in\prod_{i=0}^n\mathfrak P_i: \Prob\biggl((X_0,\ldots,X_n)=\un{x}\bigg|(X_{n-\ov{k}},\ldots,X_n)=\un{b}^{(n)}
\biggr)\neq 0
\biggr\}\\
&\hOmega^n(M):=\biggl\{(x_i)_{0\leq i\leq n}\in \wh{\Sigma}^n:\sum_{i=0}^n f_i(x_i)\in [\wh{z}_n-M,\wh{z}_n+M]\biggr\}.
\end{align*}
Let $\mu_n$ be the (atomic) measure on $\wh{\Sigma}^n$ given by ${ \mu_n}(E):=\Prob[(X_0,\ldots,X_n)\in E]$. 
 
 Let $c_1$ be the constant from Sublemma \ref{sl.blocks} such $|I_p|\leq c_1 L$ for all $p$, and recall that $\|f_j\|_\infty\leq K$ for all $j$.  Let $M:=2c_1 K(L+\overline{k}+1)$. 
\begin{sublem}\label{sl.DownSlope} 
$\DS \lim_{n\to\infty} \frac{\Prob((X_0,\ldots,X_n)\in\Omega^n)}{\Prob((X_0,\ldots,X_n)\in\hOmega^n(M))}=0.$
\end{sublem}
\begin{proof}We show that there are positive constants $c,\mathsf K, n_J$ $(J\in\N)$ so that
for each $J\in\N$ and $n>n_J$, there exist a subset $\brOmega^n \subset \Omega^n$
and $J$  maps
$ \DS \Theta_1, \Theta_2, \dots, \Theta_J:\brOmega^n\mapsto  \hOmega^n(M)$ as follows:
\begin{enumerate}[(a)]
\item  $\Prob[(X_0,\ldots,X_n)\in \brOmega^n]\geq (1/2)\Prob[(X_0,\ldots,X_n)\in \Omega^n]$;
\item $\Theta_1(\brOmega^n),\ldots, \Theta_J(\brOmega^n)$ are  pairwise disjoint;
 \item For every $\un{x}\in \brOmega^n$, $\mu_n(\Theta_j(\un{x}))\geq c\mu_n(\un{x})$;
\item Each point in the image of $\Theta_j$ has at most $\mathsf{K}$ 
$\Theta_j$-preimages.
\end{enumerate}

To see why this implies the sublemma, take $n>n_J$ and  use  (d) to split
$\DS
\brOmega^n=\biguplus_{\ell=1}^{\mathsf K}\brOmega^n_{j,\ell}$ so that $\Theta_j$   is injective on each $\brOmega^n_{j,\ell}$. (We allow some of these sets to be empty. No deep measurable selection theorems are needed here,  $\brOmega^n$ is a finite set.) Then:
  \begin{align*}
 &\Prob[(X_0,\ldots,X_n)\in\wh{\Omega}^n(M)]=\mu_n[\wh{\Omega}_n(M)]\geq \sum_{j=1}^J \mu_n[\Theta_j(\brOmega^n)],\text{ by (b)}\\
 &\geq \sum_{j=1}^J \int \frac{1}{\mathsf{K}}\sum_{\ell=1}^{\mathsf K}1_{\Theta_{j}(\brOmega^n_{j,\ell})}\, d\mu_n=\frac{1}{\mathsf{K}}\sum_{j=1}^J \sum_{\ell=1}^{\mathsf K}\mu_n[{\Theta_j(\brOmega_{j,\ell})}]\\
&\geq \frac{1}{\mathsf{K}}\sum_{j=1}^J \sum_{\ell=1}^{\mathsf K}c\mu_n[\brOmega_{j,\ell}],\text{ by (c) and the injectivity of $\Theta_j$ on $\brOmega_{j,\ell}$}\\
&=\frac{1}{\mathsf K} \sum_{j=1}^J c\mu_n(\brOmega^n)=\frac{cJ}{\mathsf K}\mu_n(\brOmega^n)=\frac{cJ}{\mathsf K}\Prob[(X_0,\ldots,X_n)\in \brOmega^n]
\\&
\geq \frac{cJ}{2\mathsf K}\Prob[(X_0,\ldots,X_n)\in \Omega^n] , \text{ by }(a). 
 \end{align*}
Thus $\DS\limsup_{n\to\infty}\frac{\Prob((X_0,\ldots,X_n)\in\Omega^n)}{\Prob((X_0,\ldots,X_n)\in\hOmega^n(M))}\leq \frac{2\mathsf K}{cJ}$ for all $J$, and the limit  must be zero.
 
We now construct $\Theta_1, \dots, \Theta_J.$
Call a {segment $I_j$} {\em typical} if \index{Typical segment}\index{Atypical segment}
$$
\Prob\left(\text{$(X_i)_{i\in I_j}$ is optimal}\big| (X_0,\ldots,X_n)\in \Omega^n \right)\geq c_J:=1-1/(10J),
$$
and {\em atypical} otherwise. 

The frequency of atypical segments tends to zero:
\begin{align*}
&\frac{1}{\nu_n}
\sum_{p=1}^{\nu_n}1_{\{j:I_j\text{ is atypical}\}}(p)\leq \frac{1}{\nu_n}\sum_{p=1}^{\nu_n}\frac{\Prob[(X_i)_{i\in I_p}\text{ is {\em sub}-optimal}|(X_0,\ldots,X_n)\in\Omega^n]}{1-c_J}\\
&=10J\E\left(\frac{1}{\nu_n}\sum_{p=1}^{\nu_n}1_{[(X_i)_{i\in I_p}\text{ is sub-optimal}]}(X_0,\ldots,X_n)\bigg|(X_0,\ldots,X_n)\in\Omega^n\right)\\
&\leq 10J\cdot \frac{n}{\nu_n}\cdot \frac{1}{n}\max_{\Omega^n}(e_n+k)\xrightarrow[n\to\infty]{}0\text{ by \eqref{ManyMax}, and since $\nu_n\asymp n$}. 
\end{align*} 

Necessarily, the frequency of typical segments tends to one. Fix $J\in\N$. Then for all $n$ sufficiently large, the first $J$ typical segments $I_{j_1},\ldots,I_{j_J}$ are all inside $[0,n-k-2]\cap[0,n-\ov{k}-2]$. Let
\begin{align*}
&\brOmega^n:=\{\un{x} \in \Omega^n:  (x_i)_{i\in I_p}\text{ is optimal for every }p=j_1, \dots, j_J\},
\end{align*}
then $\Prob[(X_0,\ldots,X_n)\in \brOmega^n\, |\, (X_0,\ldots,X_n)\in \Omega^n]$
\begin{align*}
&\geq 1-\sum_{p\in \{j_1,\ldots,j_J\}} \Prob[(X_i)_{i\in I_{p}}\text{ is sub-optimal}\, |\, (X_0,\ldots,X_n)\in \Omega^n]\\
&\geq 1-J(1-c_J)\geq 0.9, \text{ (the summands are less than $1-c_J$ because of typicality). }
\end{align*}
So $\Prob[(X_0,\ldots,X_n)\in \brOmega^n]\geq 0.9\Prob[(X_0,\ldots,X_n)\in \Omega^n]$, giving us (a).

Next we construct $\Theta_i$. 
Impose some arbitrary total order on each set $\mathfrak P_j$. 
For each $1\leq i\leq J$, let $\Theta_i$ be the map on $\ov{\Omega}^n$,  which
\vskip1mm 

\noindent
(1) 
replaces the  $j_i$ block of $\un{x}\in \ov{\Omega}^n$ by the
lexicographically smallest sub-optimal segment  with the same endpoints. (Such segments exist by Sublemma \ref{sl.sub-optimal-exist}.) 

\vskip1mm 
\noindent
(2) replaces the suffix $(x_{n-\overline{k}-2},\ldots,x_n)$ by 
$(x_{n-\overline{k}-2},b,\un{b}^{(n)})$, where $b$ is the minimal element of $\mathfrak P_{n-\overline{k}}$ such that 
$\DS 
\Prob[(X_{n-\overline{k}-2},X_{n-\overline{k}-1},X_{n-\overline{k}})\!\!=\!\!(x_{n-\overline{k}-2},
b,{b}^{(n)}_{n-\overline{k}})]\neq 0
$
(such $b$ exists, because of uniform ellipticity.)
\vskip1mm
\noindent
Since $\|f_i\|_\infty\leq K$ and  $\Theta_i$ changes at most $c_1 L+\ov{k}+1$ coordinates, 
 $\Theta_i(\brOmega^n)\subset \hOmega^n(M)$. 

  $\Theta_1(\ov{\Omega}^n),\ldots,\Theta_J(\overline{\Omega}^n)$ are pairwise disjoint, since if $k\neq \ell$, then the $j_{k}$-th block in every element of $\Theta_{k}(\overline{\Omega}^n)$ is sub-optimal by construction, whereas the  $j_{k}$-th block in every element of $\Theta_{\ell}(\overline{\Omega}^n)$ is optimal by  choice of $\brOmega^n$.
  
  The Radon-Nikodym derivatives of  $\Theta_i$ are bounded away from zero by constants independent of $i,J$ and $n$, because $\Theta_i$ changes a bounded number of coordinates whilst respecting the admissibility of transitions, and the (non-zero) transition probabilities of the Markov chain $\{X_i\}$ are uniformly bounded below (Proposition \ref{Prop-Markov-Kernel}). 
  
   The same reasoning, but using the uniform boundedness of $|\mathfrak P_i|$,  allows one to bound the number of preimages of $\Theta_j$.
This completes the proof of (a)--(d). 
\end{proof}

\noindent
{\bf Proof of Lemma \ref{l.Pure-Seg-Ratio-Edge-Case}.}
A calculation similar  to the one leading to \eqref{e.engine3} shows that 
\begin{align*}
&\Prob_{W_n}\bigl[(\wt{\psi}_{n-1}\circ\cdots\circ\wt{\psi}_0)(\wt{F}_0)\bigr]
=\Prob_{W_n}\left\{\wt{\omega}\in W_n: \begin{array}{l}
\sum\limits_{i=1}^n f_i(X_i^{(n)}(\pi(\tomega)))=z(W_n)=\wh{z}_n
\end{array}
\right\}\\
&
=\Prob\left\{S_n=\wh{z}_n,(X_{n-k},\ldots,X_n)=\un{a}^{(n)}
\right\}/\mu_0(\<\un{a}^{(n)}\>)\\
&=\Prob((X_0,\ldots,X_n)\in\Omega^n)/\mu_0(\<\un{a}^{(n)}\>)\asymp \Prob((X_0,\ldots,X_n)\in\Omega^n),
\end{align*}
because $\mu_0(\<\un{a}^{(n)}\>)\asymp 1$, since  $\un{a}^{(n)}$ has constant length.

By our assumptions,  $|z(\ov{W}_n)-\wh{z}_n|\leq 2Z_0$. Let  $M_0:=M+2Z_0$. Then we also have
\begin{align*}
&\Prob_{\ov{W}_n}\left[(\wt{\psi}_{n-1}\circ\cdots\circ\wt{\psi}_0)
\biggl(\bigcup_{\ell=-M_0}^{M_0}
D^\ell \wt{F}_0)\biggr)\right]
=\Prob_{\ov{W}_n}\left\{\begin{array}{c}
|\sum\limits_{i=1}^n f_i(X_i^{(n)}(\pi(\tomega)))-z(\ov{W}_n)|\leq M_0
\end{array}
\right\}\\
&
\geq \Prob\left\{|S_n-\wh{z}_n|\leq M\ , \ (X_{n-\ov{k}},\ldots,X_n)=\un{b}^{(n)}
\right\}\big/\mu_0(\<\un{b}^{(n)}\>)\\
&=\Prob\bigl((X_0,\ldots,X_n)\in\hOmega^n(M)\bigr)/\mu_0(\<\un{b}^{(n)}\>)
\asymp \Prob((X_0,\ldots,X_n)\in\hOmega^n(M)).
\end{align*}
Equation \eqref{e.P-ratio} now follows from Sublemma \ref{sl.DownSlope}.

This completes the proof assuming condition (a). If we replace (a) by (a'), a similar argument works, except that one needs to reverse the inequality in the definition of optimality in \eqref{e.def-of-optima}, and then update the proof of Sublemma \ref{sl.ManyMax} accordingly, working towards the inequalities $F_{n_j}'(-R)<z_{n_j}-\frac{1}{3}\delta n_j$ and $\PP_n'(-R)\leq z_{n_j}$. 
\qed

\begin{figure}
  \centering
  \fbox{\includegraphics[angle=90,width=9cm]{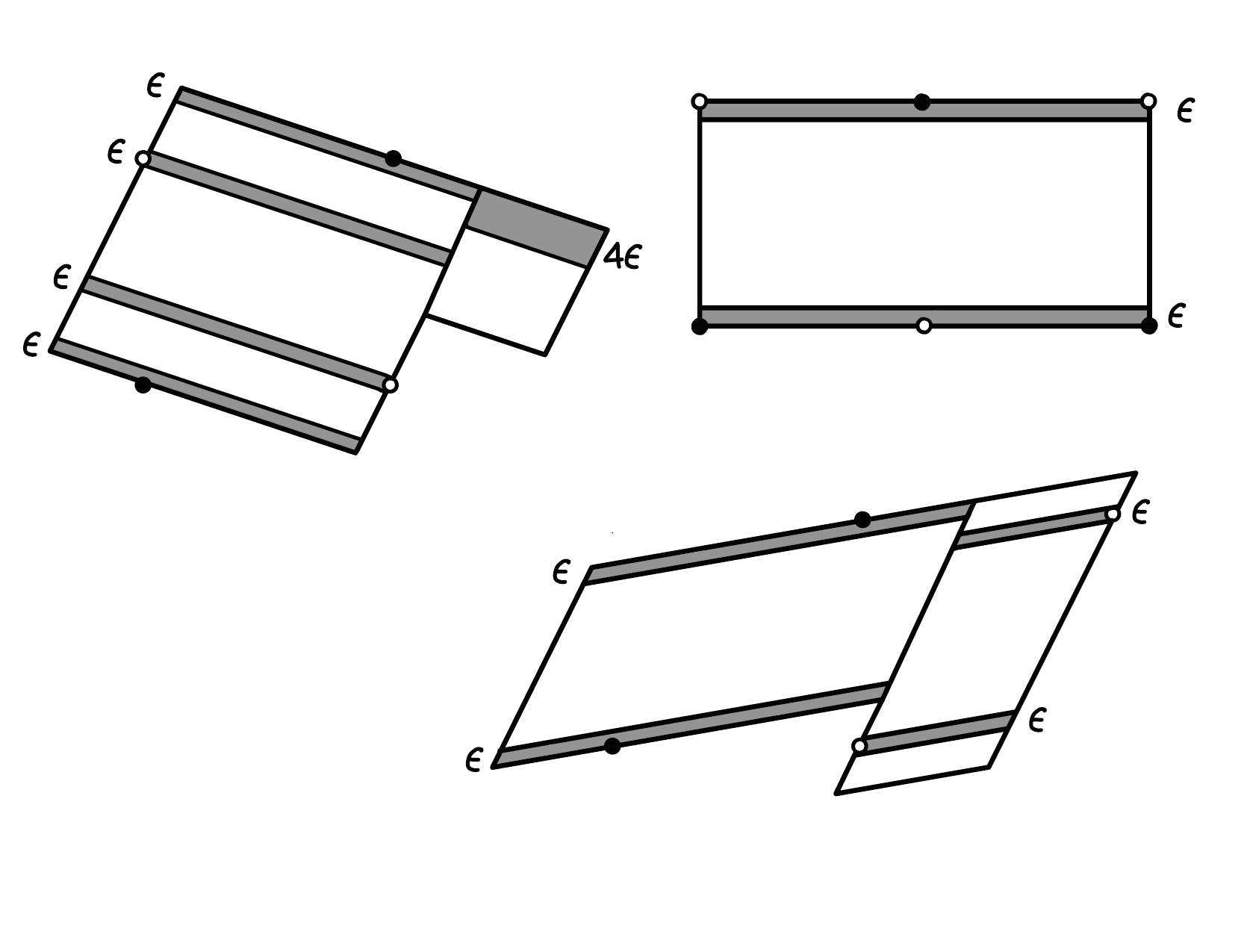}}\\
  \caption{The sets $\mathsf M_{\eps,k}$ for $\Rect_k$ and $D^k(\wt{F}_0)$ are the complement of the shaded strips. The black and hollow points represent the singularities. }\label{Figure-M-sets}
\end{figure}

\begin{lem}\label{l.Generic-Good-For-Rectangles}
If $\wt{\omega}$ is generic for some  $\phi$-ergodic invariant locally finite measure $\mu$, then for any pair of horizontal rectangles $\Rect_i, \Rect_j$ with non-zero measure, 
$$
\DS\frac{\int_0^T 1_{\Rect_i}(\phi^t(\wt{\omega}))dt}{\int_0^T 1_{\Rect_j}(\phi^t(\wt{\omega}))dt}\xrightarrow[T\to\infty]{}\frac{\mu(\Rect_i)}{\mu(\Rect_j)}, \text{ and }\int_0^\infty 1_{\Rect_i}(\phi^t(\wt{\omega}))dt=\infty.
$$
A similar statement holds if $\Rect_i, \Rect_j$  are replaced by $D^i(\wt{F}_0)$, $D^j(\wt{F}_0)$. 
\end{lem}
\begin{proof}
The (small) difficulty is that the indicators of $\Rect_k$ and  $D^k(\wt{F}_0)$  
$(k=i,j)$ cannot be sandwiched between continuous functions with compact support,  because the singularities on their boundaries do not have compact neighborhoods.

Geometric considerations
show that $\int_0^T 1_{\partial\Rect_k}(\phi^t(\tomega'))dt=0$ for all $T$ and $\tomega'$, therefore by the ratio ergodic theorem, $\mu(\partial\Rect_k)=0$ $(k=i,j)$. 
Divide $\Rect_k$  into three horizontal strips: The top and bottom with width $\epsilon\in (0,\frac{1}{3})$, and the middle, with width $1-2\epsilon$. Let $\mathsf M_{\epsilon,k}$ denote the middle strip. Since $\mu(\partial\Rect_k)=0$,  $\mu(\mathsf M_{k,\eps})\xrightarrow[\eps\to 0^+]{}\mu(\Rect_k)$. 

 $\mathsf M_{\epsilon,k}$ are bounded away from the singularities of the staircase, so their indicators can be sandwiched between  continuous functions with compact support, and  nearly the same integrals. 
 Therefore, for generic $\wt{\omega}$, 
$
\DS \lim_{T\to\infty}\frac{\int_0^T 1_{\mathsf M_{\epsilon,i}}(\phi^t(\wt{\omega}))dt}
{\int_0^T 1_{\mathsf M_{\epsilon,j}}(\phi^t(\wt{\omega}))dt}
=\frac{\mu(\mathsf M_{\epsilon,i})}{\mu(\mathsf M_{\epsilon,j})}.
$
Moreover,  the numerator and denominator tend to infinity (see the proof of \eqref{e.RET-conservative}). Clearly, 
$$
\int_0^T 1_{\Rect_k}(\phi^t(\wt{\omega}))dt=(1-2\epsilon)^{-1}\int_0^T 1_{\mathsf M_{\epsilon,k}}(\phi^t(\wt{\omega}))dt+O(1).
$$ Thus, 
$
\DS\lim_{T\to\infty}\frac{\int_0^T 1_{\Rect_i}(\phi^t(\wt{\omega}))dt}
{\int_0^T 1_{\Rect_j}(\phi^t(\wt{\omega}))dt}
=\frac{\mu(\mathsf M_{\epsilon,1})}{\mu(\mathsf M_{\epsilon,2})}
$.
Passing to the limit $\eps\to 0$, gives the first part of lemma for horizontal rectangles. The second part $\int_0^\infty 1_{\Rect_k}(\phi^t(\tomega))dt=\infty$ follows, as in the proof of \eqref{e.RET-conservative}. 

The case $D^k(\wt{F}_0)$ is handled in a similar way, using the sets $\mathsf M_{\epsilon,k}$  in Figure \ref{Figure-M-sets}.
\end{proof}

\begin{proof}[Proof of Proposition \ref{PrEdge}]  Fix $\tomega\in\St$, and  define 
$$
\xi_n:=\begin{cases} -\infty & z_n(\tomega)\leq \lim\limits_{t\to -\infty} \PP'_n(t)\\
+\infty & z_n(\tomega)\geq \lim\limits_{t\to\infty}\PP'_n(t)\\
\text{the unique $\xi_n$ s.t.  $\PP_n'(\xi_n)=z_n(\tomega)$} & z_n(\tomega)\in (\PP'(-\infty),\PP'_n(+\infty)).
\end{cases}
$$
The definition is proper, because $\PP_n'$ is continuous and strictly increasing. 

By Proposition \ref{p.Pressure-and-Geometric-Pressure}, $\xi_N(\tomega)$  (defined in \eqref{e.xi-def}) converges to a finite limit iff $\xi_n$ converges to a finite limit.
Therefore, to prove Proposition \ref{PrEdge} it is sufficient to show that if $\xi_n$ does not converge  to a finite limit, then $\tomega$ is not generic for any locally finite $\phi$-ergodic invariant measure. 
There are three cases
a to consider:
\begin{enumerate}[(1)]
\item  $\{\xi_n\}_{n\geq 1}$ is eventually finite and  bounded, but has  more than one limit point;
 
 \item  $\exists n_k\to\infty$ s.t. $|\xi_{n_k}|<\infty$ and $\xi_{n_k}\to\infty$,  or  $\xi_{n_k}=\infty$ for all $k$.
  
\item  $\exists n_k\to\infty$ s.t. $|\xi_{n_k}|<\infty$ and $\xi_{n_k}\to -\infty$, or  
   $\xi_{n_k}=-\infty$ for all $k$.
 \end{enumerate}

In case 1,   $\xi_N(\tomega)$ has several finite limit points (Proposition \ref{p.Pressure-and-Geometric-Pressure}),  therefore 
 $\wt{\omega}$ has  several different historical measures (Proposition \ref{PrBulk}), and   $\tomega$ cannot be generic. Next we consider cases 2 and 3.

Let $A_T(\tomega):=\{\phi^t(\wt{\omega}):0<t<T\}$.
Fix $0<\eps<1$, $k\geq 1$ and a constant $L_0$ so large that for every $n$, any pure $E^s_n$-segment of depth $\ell$ has length less than $\eps L_0$. We saw in the proof of Proposition \ref{PrBulk}, that there is a  constant $c>4$ such that for every $T$ there is $n\asymp\log T$ such that 
$
A_T^\ast(\tomega):=(\tpsi_{n-1}\circ\cdots\circ\tpsi_0)(A_T)$  has length  $|A_T^\ast|\in [L_0,cL_0]$. 

For fixed $n$, $T\mapsto |(\tpsi_{n-1}\circ\cdots\circ\tpsi_0)(A_T)|$ is continuous and strictly increasing. Therefore, for every $n$, there is $T_n$ such that 
$
|A^\ast_{T_n}|\in [L_0,cL_0], \text{ and } T_n\to\infty.
$
$T_n$ does not depend on $\tomega$. It is completely determined by $c,L_0$, and $\|d\tpsi_i|_{E^s_i}\|$ $(i\geq 0)$.

By Lemma \ref{l.approx-by-pure-segments},  one can approximate $A_{T_n}^\ast$ by finite unions of pure $E^s_n$-segments $W_{n,i}$ of bounded depth which intersect $A_{T_n}^\ast$, so that 
$\DS
\bigcup_{i=1}^{K_n} W_{n,i}\subset A^\ast_{T_n}(\tomega)\subset \bigcup_{i=1}^{K_n+2} W_{n,i}
$, and 
\begin{equation}\label{e.approx}
\frac{1}{|A_{T_n}^\ast|}\sum_{i=1}^{K_n}|W_{n,i}|\Prob_{W_{n,i}}(\cdot)\leq \Prob_{A_{T_n}^\ast}(\cdot)\leq \frac{1}{|A_{T_n}^\ast|}\sum_{i=1}^{K_n+2}|W_{n,i}|\Prob_{W_{n,i}}(\cdot).
\end{equation}
$K_n$ is bounded, and since $|W_{n,i}|$ and $|A_{T_n}^\ast|$ are uniformly bounded, there exists $Z_0$ finite so that  $|z_n(W_{n,i})-z_{nN_0}(\tomega)|\leq Z_0$. (This is certainly true for the canonical $\Z$-coordinates; For the modified $\Z$-coordinates, use Lemma \ref{Lemma-Modified-Approx-Canonical}.)

Suppose first that $\exists k_j\uparrow\infty$ such that  $\Prob_{W_{n_{k_j}, i}}((\wt{\psi}_{n_{k_j}}\circ\cdots\circ\wt{\psi}_0)(\wt{F}_0))=0$ for each $1\leq i\leq K_{n_{k_j}}+2$. By \eqref{e.approx},  
$\Prob_{A_{T_{n_{k_j}}}^\ast}[(\wt{\psi}_{n_{k_j}}\circ\cdots\circ\wt{\psi}_0)(\wt{F}_0)]=0$, whence  $\Prob_{A_{T_{n_{k_j}}}}(\wt{F}_0)=0$, whence  
$
\int_0^{\infty}1_{\wt{F}_0}(\phi^t (\tomega))dt=0
$.  
By Lemma \ref{l.Generic-Good-For-Rectangles}, $\tomega$ is not generic. 

Suppose now that  for every $k$ large enough, there exists some $1\leq i_k\leq N_k+2$ such that   $\Prob_{W_{n_{k_j}, i}}((\wt{\psi}_{n_{k_j}}\circ\cdots\circ\wt{\psi}_0)(\wt{F}_0))\neq 0$. 

Then Lemma \ref{l.Pure-Seg-Ratio-Edge-Case} applies  with condition (a) in case 2 and condition (a') in case (3), and there exists a finite union $\wt{B}:=\bigcup_{\ell=-M_0}^{M_0}D^\ell\wt{F}_0$ such that 
$$
\limsup_{k\to\infty}\frac{\sum\limits_{i=1}^{K_{n_k}+2}|W_{n_k,i}|\Prob_{W_{n_k,i}}((\wt{\psi}_{n_{k_j}}\circ\cdots\circ\wt{\psi}_0)(\wt{F}_0))}
{\Prob_{W_{n_k,i_k}}((\wt{\psi}_{n_{k_j}}\circ\cdots\circ\wt{\psi}_0)(\wt{B}))}=0. 
$$
By \eqref{e.approx}, $\DS\lim_{k\to\infty}\frac{\Prob_{A_{T_{n_k}^\ast}}((\wt{\psi}_{n_{k_j}}\circ\cdots\circ\wt{\psi}_0)(\wt{F}_0))}
{\Prob_{A_{T_{n_k}^\ast}}((\wt{\psi}_{n_{k_j}}\circ\cdots\circ\wt{\psi}_0)(\wt{B}))}=0$, whence 
$$
\lim_{k\to\infty}\frac{\int_0^{T_{n_k}} 1_{\wt{F}_0}(\phi^t(\wt{\omega})) dt}
{\int_0^{T_{n_k}} 1_{\wt{B}}(\phi^t(\wt{\omega}))dt}\equiv\lim_{k\to\infty}\frac{\Prob_{A_{T_{n_k}}}(\wt{F}_0)}{\Prob_{A_{T_{n_k}}}(\wt{B})}
=\lim_{k\to\infty}\frac{\Prob_{A_{T_{n_k}}}((\wt{\psi}_{n_{k_j}}\circ\cdots\circ\wt{\psi}_0)(\wt{F}_0))}
{\Prob_{A_{T_{n_k}}}((\wt{\psi}_{n_{k_j}}\circ\cdots\circ\wt{\psi}_0)(\wt{B}))}
=0.
$$
Again, by Lemma \ref{l.Generic-Good-For-Rectangles}, $\tomega$ cannot be generic.  
In summary, $\tomega$ is not generic for any locally finite $\phi$-ergodic invariant measure in any of the cases (1)--(3). 
\end{proof}

\begin{remark}
\label{RkMesClass}
\normalfont
Corollary \ref{c.Maharam-Measures} states that every locally finite $\phi$-ergodic invariant measure is proportional to a unique measure $\mu_\xi$ such that $\mu_\xi(\Rect_0)=1$ and $\mu_\xi\circ D=e^{\xi}\mu_\xi$. We can now give an alternative proof of this claim, assuming $\alpha$ has bounded type. 

Given $\xi\in\R$, let $\nu_\xi$ denote the unique extension of  $M_\xi$ to a locally finite measure which is locally flow invariant on $\mathfs G$ (Lemma \ref{l.Extending-M-xi}). Looking at the definition of $M_\xi$, we see that $e^{-\xi}M_\xi\circ D=M_\xi$, and $0<M_\xi(\wt{F}_0)<\infty$. Since  $\nu_\xi$ is uniquely determined,  $e^{-\xi}\nu_\xi\circ D=\nu_\xi$. 
Also,  $\nu_\xi(\wt{F}_0)=M_\xi(\wt{F}_0)\in (0,\infty)$, whence $\nu_\xi(\Rect_0)\in (0,\infty)$.

Suppose $\mu$ is a locally finite $\phi$-ergodic invariant measure.  Let $\tomega$ be a $\mu$-generic point. By Proposition \ref{PrEdge}, $\xi_n(\tomega)$ is bounded. The proof of Proposition \ref{PrBulk} shows (without relying on Corollary \ref{c.Maharam-Measures}) that if $\xi$ is a limit point of $\xi_n(\tomega)$, then $\nu_\xi$ is a historical measure for $\tomega$. But by genericity, the only historical measure for $\tomega$ is $\mu$. So  $\xi_n(\tomega)\to\xi$,  $\tomega$ is $\nu_\xi$-generic, and  $\mu$ is proportional to $\nu_\xi$.  Necessarily, $\nu_\xi$ is ergodic and invariant, and we saw above that $\nu_\xi\circ D=e^{\xi}\nu_\xi$ and $0<\nu_\xi(\Rect_0)<\infty$.  Setting $\mu_\xi=\nu_\xi/\nu_\xi(\Rect_0)$, we obtain the result. 
\end{remark}

\subsection{Consequences (Theorems  \ref{t.char-of-equidist-seq},  \ref{t.tilt-and-exp-tilt}, \ref{CrLebGen}, \ref{t.Char-mu-xi-Gen},  \ref{CrDriftSeq}, \ref{t.historical-crit},  \ref{CrHist})}\label{ss.Proof-1}\ \\

\smallskip
\noindent
{\bf Proof of Theorem \ref{t.historical-crit}.} 
It is convenient to start with the second part of the theorem, which says that if $\xi$ is a limit point of $\xi_N(\tomega)$ and  $|\xi|=\infty$, then there are two finite unions of horizontal rectangles $\mathsf F_1,\mathsf F_2$ and a sequence $T_n\to\infty$ such that 
\begin{equation}\label{e.Prop-1.31-statement}
\frac{\int_0^{T_n}1_{\mathsf F_1}(\phi^t(\tomega))dt}{\int_0^{T_n}1_{\mathsf F_2}(\phi^t(\tomega))dt}\xrightarrow[n\to\infty]{}0.
\end{equation}
Suppose $N_k\to\infty$ and $\xi_{N_k}(\tomega)\to\infty$ (this includes the case when $\xi_{N_k}(\tomega)=\infty$ for infinitely many $k$).
Let $n_{k}:=\lfloor N_{k}/N_0\rfloor$. By Proposition \ref{p.Pressure-and-Geometric-Pressure}, for all $k$ large enough either $z_{n_k N_0}(\tomega)\geq \PP_{n_k}'(\infty):=\lim\limits_{t\to\infty}\PP_{n_k}'(t)$, or  
$z_{n_k N_0}(\tomega)\leq \PP_{n_k}'(-\infty):=\lim\limits_{t\to -\infty}\PP_{n_k}'(t)$, 
or  $\PP_{n_k}'(\eta_{n_k})=z_{n_k N_0}(\tomega)$ for some $\eta_{n_k}$ such that $|\eta_{n_k}|\to\infty$.

We claim that a similar statement holds for  $D^j(z_{n_k N_0}(\tomega))$, for any $j\in\Z$. Fix $j\in\Z$.  Either for all $k$ large enough
$z_{n_k N_0}(D^j(\tomega))\geq \PP_{n_k}'(\infty)$ or  
$z_{n_k N_0}(D^j(\tomega))\leq \PP_{n_k}'(-\infty)$,  
or $\exists k_i\uparrow\infty$ and  $\exists\eta_{n_{k_i},j}$ such that $\PP_{n_{k_i}}'(\eta_{n_{k_i},j})=z_{n_{k_i} N_0}(D^j(\tomega))$. 

Necessarily, $|\eta_{n_{k_i},j}|\to\infty$. Otherwise (after passing to a further subsequence) $|\eta_{n_{k_i},j}|\leq R$ for all $i$. 
By Proposition \ref{Prop-Var-Growth} and \eqref{Pres2dDer}, there is a constant  $\mathsf C>0$ such that  $\PP_n''\geq \mathsf C n$ everywhere on $[-(R+1),(R+1)]$, for all $n>\mathsf C$. If  also $n>|j|/\mathsf C$, then by the mean value theorem,
\begin{align*}
\PP_{n_{k_i}}'(R+1)&>\PP_{n_{k_i}}'(R)+\mathsf C n\geq \PP_{n_{k_i}}'(\eta_{n_{k_i},j})+|j|\\
&=z_{n_{k_i}N_0}(D^j(\tomega))+|j|= z_{n_{k_i}N_0}(\tomega)+j+|j|\geq z_{n_{k_i}}(\tomega)\\
\PP_{n_{k_i}}'(-R-1)&<\PP_{n_{k_i}}'(-R)-\mathsf C n\leq \PP_{n_{k_i}}'(\eta_{n_{k_i},j})-|j|\\
&=z_{n_{k_i}N_0}(D^j(\tomega))-|j|= z_{n_{k_i}N_0}(\tomega)+j-|j|\leq z_{n_{k_i}N_0}(\tomega).
\end{align*}
It follows that the equation  $\PP_{n_{k_i}}'(\eta)=z_{n_{k_i}N_0}(\tomega)$ has a solution in $[-(R+1),R+1]$. But the unique solution is $\eta_{n_{k_i}}$, and $|\eta_{n_{k_i}}|\to\infty$.

Applying a diagonal argument, we obtain a subsequence $n_k^\ast\uparrow \infty$ such that for each $j\in\Z$, one of the following alternatives holds:
\begin{itemize}
\item $z_{n_k^\ast N_0}(D^j(\tomega))\geq \PP_{n_k^\ast}'(\infty)$ for all $k$ large enough, or
\item $z_{n_k^\ast N_0}(D^j(\tomega))\leq \PP_{n_k^\ast}'(-\infty)$ for all $k$ large enough, or
\item for all $k$ large enough, $z_{n_{k}^\ast N_0}(D^j(\tomega))=\PP_{n_{k}^\ast}'(\eta_{n_{k}^\ast,j})$ for some 
$\eta_{n_{k}^\ast,j}\in\R$, and $\eta_{n_{k}^\ast,j}$ converges to a limit in $\{-\infty,\infty\}$. 
\end{itemize}

We now argue as in the proof of Proposition \ref{PrEdge}, and construct $T_k\to\infty$ such that 
$|\tpsi_{n_k^\ast-1}\circ\cdots\circ\tpsi_0(A_{T_k})|\in [L_0,cL_0]$ (with $c,L_0$ as in that proof). 
Again, one can use Lemma \ref{l.Pure-Seg-Ratio-Edge-Case} to produce constants $M_j<\infty$ such that 
$$
\lim_{k\to\infty}\frac{\int_0^{T_k} 1_{\wt{F}_0}(\phi^t(D^j(\wt{\omega}))) dt}
{\int_0^{T_k} 1_{\wt{B}_j}(\phi^t(D^j(\wt{\omega})))dt}=0,\text{ where }
 \wt{B}_j:=\bigcup_{|i|\leq M_j}D^i(\wt{F}_0).
$$
Since $\phi^t\circ D=D\circ\phi^t$, $\DS
\lim_{n\to\infty}\frac{\int_0^{T_k} 1_{D^{-j}\wt{F}_0}(\phi^t(\wt{\omega})) dt}
{\int_0^{T_k} 1_{D^{-j}\wt{B}_j}(\phi^t(\wt{\omega}))dt}=0$ for all $j\in\Z$, and it follows that  for all $N$, 
$\DS
\lim_{k\to\infty}\frac{\int_0^{T_k} 1_{\bigcup_{|j|\leq N}D^{-j}\wt{F}_0}(\phi^t(\wt{\omega})) dt}
{\int_0^{T_k} 1_{\bigcup_{|j|\leq N}D^{-j}\wt{B}_j}(\phi^t(\wt{\omega}))dt}=0.
$

By Lemma \ref{Lemma-Modified-Approx-Canonical}, there is a  constant $N$ such that 
$\DS
\bigcup_{|j|\leq N}D^{-j}\wt{F}_0\equiv \big[|\wt{\zc}_0|\leq N\big]\!\!\supset\!\! [\zc\!\!=\!\!0]\!\!=\!\!\Rect_0$. 
Similarly, there is a constant $N'$ such that 
$\DS \bigcup_{|j|\leq N}D^{-j}\wt{B}_j \subset \big[|\zc|\leq N'\big]
$. Then \eqref{e.Prop-1.31-statement} holds with $\mathsf F_1:=\Rect_0$, and 
$\DS \mathsf F_2:=\bigcup_{|i|\leq N'}\Rect_i$.

\medskip
Next we prove the  first part of the theorem, which says that for every $\xi$ finite, $\xi$ is a limit point of $\xi_N(\tomega)$ iff $\mu_\xi$ is a historical measure of $\tomega$. 

Direction ($\Rightarrow$) follows from Proposition \ref{PrBulk}. To see direction $(\Leftarrow)$, suppose $\mu_\xi$ is a historical measure.  
The proof of Lemma \ref{l.Generic-Good-For-Rectangles} can be easily extended to show that for some sequence $T_k\to\infty$, for any pair of horizontal rectangles $\Rect_k, \Rect_\ell$, 
\begin{equation}\label{e.R-k-R-l}
\frac{\int_0^{T_k} 1_{\Rect_k}(\phi^t(\tomega))dt}{\int_0^{T_k} 1_{\Rect_\ell}(\phi^t(\tomega))dt}\xrightarrow[k\to\infty]{}e^{(k-\ell)\xi}. 
\end{equation}
Define integers $n_k:=n(T_k)$, as in the proof of Propositions \ref{ScBulk} and \ref{ScEdge}, so that $A_{T_k}^\ast:=(\tpsi_{n_k-1}\circ\cdots\circ\tpsi_0)(A_{T_k})$ has length in $[L_0,cL_0]$. By the argument above, $\xi_{n_k N_0}(\tomega)$ has no infinite limit points. By \eqref{e.R-k-R-l} and the proof of Proposition \ref{ScBulk}, the only possible finite limit point is $\xi$.\qed

\begin{proof}[Proof of Theorem  \ref{t.Char-mu-xi-Gen}] 
The theorem asserts the equivalence of the  conditions
\begin{enumerate}[(a)]
\item $\tomega$ is $\mu_\xi$-generic;
\item $\xi_n(\tomega)\to\xi$;
\item $\forall k,\ell\in\Z$, $\lim\limits_{T\to\infty}\frac{\int_0^T 1_{\Rect_k}(\phi^t(\tomega))dt}{\int_0^T 1_{\Rect_\ell}(\phi^t(\tomega))dt}=e^{\xi(k-\ell)}$. 
\end{enumerate}
(a)$\Rightarrow$(b): 
By Proposition \ref{PrEdge}, $\xi_n(\tomega)$ must converge to a finite limit. By Proposition~\ref{PrBulk}, this limit must be $\xi$.  

\smallskip
\noindent
(b)$\Rightarrow$(a) due to Proposition \ref{ScBulk}.

\smallskip
\noindent
(b)$\Rightarrow$(c) By Proposition \ref{PrBulk}, $\tomega$ a generic point of  $\mu_\xi$. By Lemma \ref{l.Generic-Good-For-Rectangles}, 
$$\lim\limits_{T\to\infty}\frac{\int_0^T 1_{\Rect_k}(\phi^t(\tomega))dt}{\int_0^T 1_{\Rect_\ell}(\phi^t(\tomega))dt}=\frac{\mu_\xi(\Rect_k)}{\mu_\xi(\Rect_\ell)}=\frac{\mu_\xi(D^{k-\ell}(\Rect_\ell))}{\mu_\xi(\Rect_\ell)}=e^{\xi(k-\ell)}.
$$

\smallskip
\noindent
(c)$\Rightarrow${(b)} By (c), $
\liminf\limits_{T\to\infty}\frac{\int_0^T 1_{\mathsf F_1}(\phi^t(\tomega))dt}{\int_0^T 1_{\mathsf F_2}(\phi^t(\tomega))dt}>0$ for any pair of finite unions of horizontal rectangles $\mathsf F_1, \mathsf F_2$. 
Therefore, by Theorem \ref{t.historical-crit}, all the limit points of $\xi_n(\tomega)$ are finite. 
Suppose $\eta$ is a finite limit point. By Proposition \ref{PrBulk}, $\mu_\eta$ is a historical measure of $\tomega$, and there exists some $T_n\to\infty$ so that $\lim\limits_{T\to\infty}\frac{\int_0^T 1_{\Rect_k}(\phi^t(\tomega))dt}{\int_0^T 1_{\Rect_\ell}(\phi^t(\tomega))dt}=e^{\eta(k-\ell)}$. Necessarily $\eta=\xi$. Thus all the limit points of $\xi_n(\tomega)$ equal $\xi$, a statement which is equivalent to (b). 
\end{proof}

\begin{lem}\label{Lemma-Solution}
For every $\epsilon>0$, there exists $\delta>0$ such that if $|z_n/n|<\delta$, then for all $n$ large enough, there exists $\xi_n\in [-\epsilon,\epsilon]$ such that $\PP_n'(\xi_n)=z_n$.
\end{lem}
\begin{proof}
By Lemma \ref{LmPressure} and Proposition \ref{Prop-Var-Growth} $\PP_n(0)=0$, $\PP_n'(0)=\E(S_n)=O(1)$, $V_n\asymp n$, and for some positive constant $c$, $\PP_n''(\xi)\geq c n$ on $[-1,1]$ for all $n$ large enough. Therefore, for all $0<\eps<1$,  
$$
\PP_n'(\epsilon)\geq \PP_n'(0)+cn\epsilon=cn\epsilon+O(1)\text{ and }\PP_n'(-1)\leq P_n'(0)-cn\epsilon=-cn\epsilon+O(1).
$$
Thus, if $|z_n/n|<\frac{1}{2}c\epsilon$ and $n$ is large enough,  then 
$
\PP_n'(-\epsilon)<z_n<\PP_n'(\epsilon).
$

Since $\PP_n$ is analytic and strictly convex, there exists a unique solution to the equation $\PP_n'(\xi_n)=z_n$, and $|\xi_n|<\epsilon$.
\end{proof}

\smallskip
\noindent
{\bf Proof of Theorem \ref{CrLebGen}.} The theorem asserts the equivalence  of the  conditions
\begin{enumerate}[(a)]
\item $\tomega$ is $\mu_0$-generic;
\item $z_n(\tomega)/n\to 0$;
\item $\DS \forall k,\ell\in\Z,\ \lim_{T\to \infty} \frac{\int_0^T 1_{\Rect_k} (\phi^t(\tomega)) dt}
{\int_0^T 1_{\Rect_\ell} (\phi^t (\tomega)) dt}=1$.
\end{enumerate}
The equivalence (a)$\Leftrightarrow$(c) is a special case of Theorem \ref{t.Char-mu-xi-Gen}, and was proved above. It remains to prove (a)$\Leftrightarrow$(b).  Let $z_n:=z_n(\tomega), \xi_n:=\xi_n(\tomega)$. 

Suppose (b). 
By Lemma \ref{Lemma-Solution},  $\exists \eta_n\to 0$
such that $\PP_{n}'(\eta_n)=z_n$. By Proposition \ref{p.Pressure-and-Geometric-Pressure}, $\xi_N(\tomega)\to 0$, and $\tomega$ is $\mu_0$-generic by Theorem \ref{t.Char-mu-xi-Gen}. 

Suppose (a), then Theorem \ref{t.Char-mu-xi-Gen} and Proposition \ref{p.Pressure-and-Geometric-Pressure} tell us that  $z_n=\PP_n'(\xi_n)$ where $\xi_n\to 0$. By \eqref{Pres1stDer} and 
Proposition~\ref{Prop-Var-Growth}, $\PP_n'(0)=\E(S_n)=O(1)$,  and for all the $n$ such that $\xi_n=0$, $z_n/n=\E(S_n)/n=O(1/n)\to 0$. 
For the  $n$ such that $\xi_n\neq 0$, we have 
$$
\frac{z_n}{n}=\frac{\PP_n'(\xi_n)-\PP_n'(0)}{n}+O(n^{-1})=\frac{1}{n}
\cdot\frac{\PP_n'(\xi_n)-\PP_n'(0)}{\xi_n}\cdot\xi_n+O(n^{-1}).
$$ 
For all large $n$, $|\xi_n|\!\!<\!\!1$.
By the mean value theorem, $\DS \frac{1}{n}\left|\frac{\PP_n'(\xi_n)-\PP_n'(0)}{\xi_n}\right|\leq \frac{1}{n}\sup_{[-1,1]}\PP_n''$,  and this quantity is bounded above  
 (see \eqref{Pres2dDer} and recall that $V_n\asymp n$). Thus again $z_n/n=O(|\xi_n|+n^{-1})\to 0$.  It follows that $z_n/n\to 0$ for the full sequence of $n$. \qed

\smallskip
\noindent
{\bf Proof of Theorem \ref{t.char-of-equidist-seq}.}
Without loss of generality $\alpha\in [0,1)$, and   $\beta=2\alpha-1\in (-1,1)$. Expand $\beta$ into its even continued fraction expansion: $\beta=[0; 2m_1,\ldots]$. Theorem \ref{t.char-of-equidist-seq} states that $D_N(\alpha)$ is equidistributed in $\frac{1}{2}\Z$ iff $\frac{1}{n}\sum_{j=1}^n [\sgn(m_{2j})-\sgn(m_{2j-1})]\xrightarrow[n\to\infty]{}0$. 

Suppose $D_N(\alpha)$ is equidistributed in $\frac{1}{2}\Z$.  Then $D_N(\alpha)$ attains any value in $\frac{1}{2}\Z$  at least once, and therefore 
$$\limsup\limits_{N\to\infty}D_N(\alpha)=\infty\text{ and }
\liminf\limits_{N\to\infty}D_N(\alpha)=-\infty.
$$ Since $|D_N(\alpha)\!-\!D_{N+1}(\alpha)|\!\leq\! \frac{1}{2}$,  $D_N(\alpha)$ attains any value in $\frac{1}{2}\Z$ infinitely many times, and
\begin{equation}\notag
\forall k\in\Z,\ \#\{1\leq N\leq n: D_N(\alpha)=\frac{1}{2}k\}\xrightarrow[n\to\infty]{}\infty.
\end{equation}

Let $\tomega_0$ denote the midpoint of the top side of horizontal rectangle $\Rect_0$. By \eqref{e.Lin-Flow-1},
\begin{equation}\label{e.1.11}
\forall k\in\Z,\ \int_{-1}^{n-1} 1_{\Rect_{k}}(\phi^t(\tomega_0))dt=\#\{1\leq N\leq n: D_N(\alpha)=\frac{1}{2}k\}.
\end{equation}
Necessarily, $\int_{-1}^{n-1} 1_{\Rect_{k}}(\phi^t(\tomega_0))dt\to\infty$, and the equidistribution of $D_N(\alpha)$ translates to the statement that 
$\DS \forall k,\ell\in\Z,\ \lim_{T\to \infty} \frac{\int_0^T 1_{\Rect_k} (\phi^t(\tomega)) dt}
{\int_0^T 1_{\Rect_\ell} (\phi^t (\tomega)) dt}=1$. 

By Theorem \ref{CrLebGen}, $\DS\frac{z_n(\tomega_0)}{n}\xrightarrow[n\to\infty]{}0.$ 
Theorem \ref{ThDriftOfOrigin}, proved in Appendix \ref{Appendix-Twist}, says that 
$
\DS z_n(\tomega_0)=\frac{1}{2}\sum_{k=1}^{nN_0}[\sgn(m_{2k})-\sgn(m_{2k-1})],$  with $N_0$ as in \eqref{e.mat-prod-form}. 
It follows that 
\begin{equation}\label{e.lari}
\DS\frac{1}{n}\sum_{k=1}^{n}[\sgn(m_{2k})-\sgn(m_{2k-1})]\xrightarrow[n\to\infty]{}0.
\end{equation}

Conversely, \eqref{e.lari} implies that $z_n(\tomega_0)/n\to 0$, whence by Theorem \ref{CrLebGen},   $\tomega_0$ is $\mu_0$-generic, and 
$\DS \forall k,\ell\in\Z,\ \lim_{T\to \infty} \frac{\int_0^T 1_{\Rect_k} (\phi^t(\tomega)) dt}
{\int_0^T 1_{\Rect_\ell} (\phi^t(\tomega)) dt}=1.$
The numerator and denominator must tend to infinity, see the argument leading to \eqref{e.RET-conservative}. Therefore, by \eqref{e.1.11},
$$
\forall k,\ell\in\Z, \ \lim_{n\to\infty}\frac{\#\{1\leq N\leq n: D_N(\alpha)=\frac{1}{2}k\}}
{\#\{1\leq N\leq n: D_N(\alpha)=\frac{1}{2}\ell\}}=1,
$$
and $D_N(\alpha)$ is equidistributed in $\frac{1}{2}\Z$.
\qed

\begin{proof}
[{\bf Proof of Theorem \ref{t.tilt-and-exp-tilt}}] The theorem states that if  $D_N(\alpha,x+I)$ is tilted, i.e. 
$$
\forall k,\ell\in\Z,\ \lim_{n\to\infty}\frac{\#\{1\leq N\leq n: D_N(\alpha,x+I)=\frac{1}{2}k\}}{\#\{1\leq N\leq n: D_N(\alpha,x+I)=\frac{1}{2}\ell\}}=\rho(k/2,\ell/2)>0,
$$
then $D_N(\alpha,x+I)$ it is exponentially tilted, i.e.  $\rho(s,t)=\exp[\eta(s-t)]$ for some $\eta\in\R$.

Assume $D_N(\alpha,x\!+\!I)$ is tilted. One shows, as in the previous proof, that 
$
\forall k\in\Z
$, 
$
\ \#\{1\leq N\leq n: D_N(\alpha,x+I)=\frac{1}{2}k\}\xrightarrow[n\to\infty]{}\infty.
$

Let $\tomega=\tomega_{-x}:=\omega(-x,0)\in\St$ be the initial condition of $\phi$,  corresponding to the initial condition $(-x,0)$ for the cylinder map, as in   Theorem \ref{t.HHW}. By \eqref{e.Lin-Flow-2}, 
\begin{align}
&\forall k,\ell\in\Z,\ \int_{-1}^{n-1} 1_{\Rect_{k}}(\phi^t(\tomega))dt=\#\{1\leq N\leq n: D_N(\alpha,x+I)=\frac{1}{2}k\}\xrightarrow[n\to\infty]{}\infty,\notag\\
&\underset{n\in\N}{\lim_{n\to\infty}}\frac{\int_{-1}^{n-1} 1_{\Rect_{k}}(\phi^t(\tomega))dt}{\int_{-1}^{n-1} 1_{\Rect_{\ell}}(\phi^t(\tomega))dt}
={\lim_{n\to\infty}}\frac{\#\{1\leq N\leq n: D_N(\alpha,x+I)=\frac{1}{2}k\}}{\#\{1\leq N\leq n: D_N(\alpha,x+I)=\frac{1}{2}\ell\}}=\rho(k/2,\ell/2)>0,\notag\\
&\hspace{-0.3cm}\text{ and therefore, }\underset{T\in\R}{\lim_{T\to\infty}}\frac{\int_{0}^{T} 1_{\Rect_{k}}(\phi^t(\tomega))dt}{\int_{0}^{T} 1_{\Rect_{\ell}}(\phi^t(\tomega))dt}=\rho(k/2,\ell/2)>0. \label{e.rho-limit}
\end{align}
In particular,   $\DS {\liminf_{T\to\infty}}\frac{\int_{0}^{T} 1_{\mathsf F_1}(\phi^t(\tomega))dt}{\int_{0}^{T} 1_{\mathsf F_2}(\phi^t(\tomega))dt}>0$ whenever $\mathsf F_1,\mathsf F_2$ are finite unions of  horizontal rectangles.
By Theorem \ref{t.historical-crit}, all the limit points of $\xi_n(\tomega)$ must be finite. 

Let $\xi_1,\xi_2$ be two finite limits points of  $\xi_n(\tomega)$. By Proposition \ref{PrBulk}, there are 
{$T_n^j\to\infty$  such that 
$\DS
\frac{\int_{0}^{T_n^j} 1_{\Rect_{k}}(\phi^t(\tomega))dt}{\int_{0}^{T_n^j} 1_{\Rect_{\ell}}(\phi^t(\tomega))dt}\to 
e^{\xi_j(k-\ell)}.$}
By \eqref{e.rho-limit}, $e^{\xi_1(k-\ell)}=\rho(k/2,\ell/2)=e^{\xi_2(k-\ell)}$, and $\xi_1=\xi_2$. Thus $\xi_n(\tomega)$ converges to a finite limit $\xi$.  By Theorem \ref{t.Char-mu-xi-Gen} and \eqref{e.rho-limit}, $\rho(k/2,\ell/2)=e^{\xi(k-\ell)}$ for all $k,\ell\in\Z$.  
\end{proof}

\smallskip
\noindent
{\bf Proof of Theorem \ref{CrDriftSeq}.}
The theorem states the existence of sequences $z_n^\xi\in\Z$ such that $\tomega$ is $\mu_\xi$-generic iff $z_n(\tomega)=z_n^\xi+o(n)$. Equivalently, if $\wt{\omega}'$ is generic for $\mu_\xi$, then
\begin{equation}\label{LmGenPair}
\text{ $\forall\tomega'',\ \wt{\omega}''$ is generic for $\mu_\xi$ iff 
$\frac{z_N(\wt{\omega}'')-z_N(\wt{\omega}')}{N}\xrightarrow[N\to\infty]{}0$.}
\end{equation}
(Take $z_n^\xi:=z_n(\tomega')$.)
 Let $z_n':=z_{nN_0}(\wt{\omega}')$ and $z_n'':=z_{nN_0}(\wt{\omega}'')$. Let $\xi_n',\xi_n''$ be the solutions to the equations $\PP_n'(\xi_n')=z_n'$,  
 $\PP_n'(\xi_n'')=z_n''$.

\medskip
\noindent
{\bf $(\Rightarrow)$:} 
Suppose $\wt{\omega}''$ is generic for $\mu_\xi$. By Theorem \ref{t.Char-mu-xi-Gen} and Proposition \ref{p.Pressure-and-Geometric-Pressure}, $\xi_n', \xi_n''\to \xi$.  In particular, $\xi_n'$ and $\xi_n''$  both eventually lie in the interval $[\xi-1,\xi+1]$. By Lemma~\ref{LmPressure} and Proposition~\ref{Prop-Var-Growth}, there are positive constants $c_1,c_2$ such that for all $n$ large enough,  $c_1<\PP_n''/n<c_2$ on $[\xi-1, \xi+1]$.  
By the mean value theorem,
$$ \frac{|z_n''-z_n'|}{n}
\leq \biggl(\sup_{[\xi-1,\xi+1]}\frac{\PP_n''}{n}\biggr)\cdot|\xi_n''-\xi_n'|\leq c_2|\xi_n''-\xi_n'|. 
$$
Since $\xi_n', \xi_n''\to \xi$, ${|z_n''-z_n'|}/{n}\to 0$. Since $z_N(\tomega')=z_{\lfloor N/N_0\rfloor}'+O(1)$ and $z_N(\tomega'')=z_{\lfloor N/N_0\rfloor}''+O(1)$, $|z_N(\tomega')-z_N(\tomega'')|/N\to 0$. 

\medskip
\noindent
$(\Leftarrow)$: Suppose $|z_N(\tomega')-z_N(\tomega'')|/N\to 0$, then 
 ${|z_n''-z_n'|}/{n}\to 0$. Fix $0<\eps<1$. By Theorem \ref{t.Char-mu-xi-Gen}, Proposition \ref{p.Pressure-and-Geometric-Pressure}, and the  $\mu_\xi$-genericity of $\tomega'$, $\xi_n'\to\xi$. Therefore, for all $n$ large enough, $|\xi_n'-\xi|<\epsilon$. By  assumption,  for all $n$ large enough, $|z_n''-z_n'|<c_1 n \epsilon$. 
By the convexity of $\PP_n$, for such $n$,  
$$
\PP_n'(\xi_n'+\epsilon)\geq \PP_n'(\xi_n')+c_1 n\epsilon=z_n'+c_1 n \epsilon>z_n''; 
$$
$$
\PP_n'(\xi_n'-\epsilon)\leq \PP_n'(\xi_n')-c_1 n\epsilon=z_n'-c_1 n \epsilon<z_n''. 
$$
So $z_n''\in [\PP_n'(\xi_n'-\epsilon),\PP_n'(\xi_n'+\epsilon)]$. Since $P_n$ is analytic and strictly convex on $[\xi-1,\xi+1]$, there is a unique finite solution $\xi_n''$ to the equation $z_n''=\PP_n'(\xi_n'')$, and $|\xi_n''-\xi_n'|<\epsilon$ for all $n$ sufficiently large. Since $\epsilon$ was arbitrary and $\xi_n'\to\xi$, $\xi_n''\to\xi$. By Proposition \ref{p.Pressure-and-Geometric-Pressure}, $\xi_N(\tomega'')\to\xi$, and by Theorem \ref{t.Char-mu-xi-Gen} $\wt{\omega}'$ is generic for $\mu_\xi$.
This concludes the proof of \eqref{LmGenPair}, and the proof of Theorem \ref{CrDriftSeq}.
\qed

\begin{lem}
\label{LmTurtle}\ 
\begin{enumerate}[(1)]
\item If $z_n:=z_{nN_0}(\tomega)$, $\PP_n'(\eta_n)=z_n$, and $\eta_n$ is bounded, then $\DS\lim_{n\to\infty} (\eta_{n+1}-\eta_n)=0.$ 
\item If  $\xi_n(\tomega)$ are bounded, then
$\DS \lim_{n\to\infty} [\xi_{n+1}(\tomega)-\xi_n(\tomega)]=0.$ 
\end{enumerate}
\end{lem}

\begin{proof} 
Suppose that 
$
\eta_n\in [-R,R]\text{  for all $n$}.$  Let $K$ be a uniform bound for the size of the jump functions $f_i$, then $|z_{n+1}-z_n|\leq K$ for all $n$. 
Choose a constant $\mathsf C$ as in \eqref{Pres1stDer} so that $|\PP_n'(\eta)-\wt{\E}^{\eta}(S_n)|\leq \mathsf C$ for all $|\eta|\leq R$. Then for every $\eta\in [-R,R]$, 
$$
|\PP_{n+1}'(\eta)-\PP_{n}'(\eta)|\leq |\wt{\E}^\eta(S_{n+1})-
\wt{\E}^\eta(S_{n})|+2\mathsf C\leq 2\mathsf C+K.
$$
Thus
$\DS \PP_{n+1}'(\eta_n)=\PP_{n}'(\eta_n)+O(1)=z_n+O(1)$, 
and since  
$|z_{n+1}-z_n|=O(1)$,  
$$
\PP_{n+1}'(\eta_{n+1})-\PP_{n+1}'(\eta_n)=z_{n+1}-z_n+O(1)=O(1).
$$
By the mean value theorem, there exists some  $\zeta_n\in [\eta_n, \eta_{n+1}]$ such that  $$\PP_{n+1}'(\eta_{n+1})-\PP_{n+1}'(\eta_n)=\PP_n''(\zeta_n) (\eta_{n+1}-\eta_n).$$
By \eqref{Pres2dDer}, and since $\eta_n\in [-R,R]$, there exists a constant $c$ which only depends on $R$ such that 
 $\PP_n''(\eta_n) \geq c V_n\xrightarrow[n\to\infty]{}\infty$. So $|\eta_{n+1}-\eta_n|=O(1/V_n)$.
 This completes the proof of the first part of the lemma 
 since $V_n\asymp cn$ due to Proposition \ref{Prop-Var-Growth}.
 
The second part now follows from Proposition \ref{p.Pressure-and-Geometric-Pressure}.
\end{proof}

\smallskip
\noindent
{\bf Proof of Theorem \ref{CrHist}.}
The first part of the theorem states that every historical measure $\mu$ is proportional to $\mu_\xi$ for some $\xi\in\R$. Let us prove this. 

Suppose $\mu$ is a historical measure for some $\tomega$. Without loss of generality, $\mu(\wt{F}_0)\neq 0$. Otherwise, we choose $j$ such that  $\mu(D^j(\wt{F}_0))\neq 0$, and replace $\mu$ and $\tomega$ by  $\mu\circ D^j$ and $D^j(\tomega)$, noting that if we can show that $\mu\circ D^j$ is proportional to some $\mu_\xi$, then it will immediately follow that $\mu$ is proportional to $\mu_\xi\circ D^{-j}=e^{-\xi j}\mu_\xi$.

Since $\mu$ is a historical measure, there is a sequence $T_k\to\infty$ such that  for all $h_1,h_2\in C_c(\St^\ast)$ with non-zero $\mu$-integrals,
$
\DS \frac{\int_0^{T_k} h_1(\phi^t(\tomega))dt}{\int_0^{T_k} h_2(\phi^t(\tomega))dt}\xrightarrow[k\to\infty]{}\frac{\int h_1 d\mu}{\int h_2 d\mu}. 
$ 

Let  $A_{T_k}:=\{\phi^t(\tomega):0<t<T_k\}$. Arguing as  in the proof of Proposition \ref{PrBulk}, we can construct $n_k\to\infty$ such that 
$$
A_{T_{k}}^\ast:=(\tpsi_{n_k-1}\circ\cdots\circ\tpsi_0)(A_{T_{n_k}})\text{ has length in $[L_0, cL_0]$} 
$$
(with $c$ and $L_0$ as in that proof). 

Without loss of generality, $\xi_{n_k}(\tomega)\xrightarrow[k\to\infty]{}\xi\in [-\infty,\infty]$, otherwise choose a sequence $k_i\uparrow\infty$  such that $\xi_{n_{k_i}}(\tomega)\to\xi$,  and replace $T_k$ by $T_{k_i}$. 

Suppose $|\xi|=\infty$, then one can show as in  the proof of Proposition \ref{PrEdge} that
there is a set $\DS \mathsf F=\bigcup_{|i|<N}D^i(\wt{F}_0)$ such that 
$
\DS \lim_{k\to\infty} \frac{\int_0^{T_k} 1_{\wt{F}_0}(\phi^t(\tomega))dt}{\int_0^{T_k} 1_{\mathsf F}(\phi^t(\tomega))dt}=0.
$
But by the choice of $T_k$, and arguing as in the proof of Lemma \ref{l.Generic-Good-For-Rectangles}, the limit is $\mu(\wt{F}_0)/\mu(\mathsf F)\neq 0$. 
Thus $|\xi|<\infty$.

Since $|\xi|<\infty$, the proof of Proposition \ref{PrBulk} gives 
$
\DS \frac{\int_0^{T_k} h_1(\phi^t(\tomega))dt}{\int_0^{T_k} h_2(\phi^t(\tomega))dt}\xrightarrow[k\to\infty]{}\frac{\int h_1 d\mu_\xi}{\int h_2 d\mu_\xi}
$, for all $h_1,h_2\in C_c(\St^\ast)$ with non-zero $\mu_\xi$-integral. 
It follows that $\mu$ is proportional to $\mu_\xi$, proving the first part of the theorem.

\medskip
Let $\mathcal H(\tomega)=\{\xi\in\R: \mu_\xi\text{ is a historical measure for $\tomega$}\}$. The second part of the theorem states that $\mathcal H(\tomega)$ is either empty, or a single point, or a (possibly infinite)  interval. It is sufficient to show that if $\xi_1<\xi_2$ belong to  $\mathcal H(\tomega)$, then $[\xi_1,\xi_2]\subset \mathcal H(\tomega)$. 

The proof of the first part of the theorem shows that if $\xi_i\in\mathcal H(\tomega)$, then there exists sequences $n_k(1),n_k(2)\uparrow\infty$ such that $\xi_{n_k(i)}(\tomega)\to\xi_i$. Build by induction a sequence $\ell_k\uparrow\infty$ such that $(\ell_{2k-1})_{k\geq 1}$ is a subsequence of $(n_{k}(1))_{k\geq 1}$, and $(\ell_{2k})_{k\geq 1}$ is a subsequence of $(n_k(2))_{k\geq 1}$. Then 
$
\xi_{\ell_{2k-1}}(\tomega)\to\xi_1\text{ and }\xi_{\ell_{2k}}(\tomega)\to\xi_2
$.

By Lemma \ref{LmTurtle},  
$\Lambda_k\!\!:=\!\!\{\xi_n(\tomega):\ell_{2k-1}\!\!\leq \!\!n\!\!\leq\!\! \ell_{2k}\}\cap [\xi_1,\xi_2]$ is  
$\delta_k$-dense in $[\xi_1,\xi_2]$ with $\delta_k\to 0$. 
Thus, for every $\xi\in [\xi_1,\xi_2]$ there is a subsequence $j_k\to\infty$ such that $\xi_{j_k}(\tomega)\to\xi$. By Theorem \ref{t.historical-crit}, $[\xi_1,\xi_2]\subset\mathcal H(\tomega)$.
\qed

\section{Limit Theorems for Birkhoff Integrals}\label{ScBirkhoff}
\subsection{Preliminary Estimates}
In this section we give estimates for intergals of the form $\int_0^T h(\phi^t(\tomega))dt$. Recall the pressure functions $\PP_n(\xi)$ from Section \ref{ScPressure}. 
The estimates involve the Legendre transforms $\HH_n(\eta)$ of $\PP_n(\xi)/V_n$, defined by \index{Legendre transform} \index{Pressure!Legendre transform of}
$$
\HH_N(\eta):=\xi\eta-\frac{\PP_N(\xi)}{V_N}\text{ for the unique $\xi$ s.t. }\frac{\PP_N'(\xi)}{V_N}=\eta, \text{ when it exists}.
$$

Let $[\ov{a}_N^R, \ov{b}_N^R]:=[\PP_N'(-R)/V_N,\PP_N'(R)/V_N]$. 
The next result is consequence of Lemma \ref{LmPressure} and \cite[Lemmata 7.23--7.25]{Dolgopyat-Sarig-Book}:

\begin{lem}\label{l.rate-function} 
$\HH_N$ is well-defined and analytic on 
$[\ov{a}_N^R, \ov{b}_N^R]$, for all $R>0$. In addition,  
\begin{enumerate}[(1)]
\item $\exists c,R>0$ such that for all $N$ large enough, $
[\ov{a}_N^R, \ov{b}_N^R]\supset [-c,c]$.
\item  $H_N(\frac{\E(S_N)}{V_N})=0$, $H_N'(\frac{\E(S_N)}{V_N})=0$ and  $\forall R>0$, $\exists\rho(R)>1$ such that for all $N$ large enough, $\rho^{-1}\leq \HH_N''\leq \rho$ on $[\ov{a}_N^R, \ov{b}_N^R]$.
\item For every $\epsilon>0$ there exist $\delta>0$ and $N_\epsilon$ such that for all $N>N_\epsilon$ and \\
$\DS \eta\!\!\in\!\! \left[\frac{\E(S_N)}{V_N}-\delta,\frac{\E(S_N)}{V_N}+\delta \right]$,\;\; 
$\DS
\frac{e^{-\epsilon} }{2}\left(\eta-\frac{\E(S_N)}{V_N}\right)^2\!\!\!\leq \HH_N(\eta)\leq \frac{e^{+\epsilon}}{2}\left(\eta-\frac{\E(S_N)}{V_N}\right)^2\!.
$
\item In particular, if  $\DS \frac{z_N}{V_N}\to 0$, then 
$\DS V_N\HH_N\left(\frac{z_N}{V_N}\right)=\frac{1+o(1)}{2}\left(\frac{z_N-\E(S_N)}{\sqrt{V_N}}\right)^2.$
\end{enumerate}
\end{lem}

Recall that $N_0$ is the integer  such that $\Psi_{nN_0}=\tpsi_{n-1}\circ\cdots\circ\tpsi_0$, and  $\xi_{nN_0}(\tomega), z_{nN_0}(\tomega)$ are given in  \eqref{e.z-n} and \eqref{e.xi-def}. Recall also the  integer-valued functions $n(T)$ from the proof of the temporal CLT, defined to be the minimal integers  such that 
\begin{equation}\label{e.n(T)-again}
T\prod_{k=0}^{n(T)-1} \|d\psi_k|_{E^s_k}\|\in [1,L].
\end{equation}
$L$ was chosen just before \eqref{DefBrN} to guarantee the existence 
of  $n(T)$, and $n(T)\asymp \log T$.  

\begin{prop}\label{p.Bikhoff-Sums-Estimate}
Suppose $h\in C_c(\St^\ast)$ is non-negative and  $\int h d\mu_0\neq 0$. Then for every $R>0$ there are $C(R,h), N(R,h)>1$ 
as follows: If $T>N(R,h)$, $n=n(T)$,   $|\xi_{nN_0}(\tomega)|<R$, then \index{Birkhoff integrals of the linear flow}
$$
\frac{1}{T}\int_0^T h(\phi^t(\wt{\omega}))dt=C(R,h)^{\pm 1}\left[
\frac{1}{\sqrt{V_{n}}}\int h d\mu_0\cdot \exp\left({-V_{n}\mathcal \HH_{n}\left(\frac{z_{nN_0}(\tomega)}{V_{n}}\right)}\right)\right].
$$
\end{prop}
\begin{proof}
Every non-negative $h\in C_c(\St^\ast)$ can be sandwiched between two functions with nearly the same integral, equal to  linear combinations  of indicators of sets of the form 
$
\wt{C}:=(D^j\circ\pi|_{\wt{F}_0}^{-1})\<\dot{a}_0,\ldots,a_\ell\>.
$
Therefore, it is enough to prove the lemma for $h=1_{\wt{C}}$.

Throughout this proof, $n=n(T)$, and 
 $z_n:=z_{nN_0}(\tomega)$.  
 
Recall that $\mathcal P_{nN_0}'(\xi_{nN_0}(\tomega))=z_n$. 
 By Proposition \ref{p.Pressure-and-Geometric-Pressure}, for all $R>0$, for all $T$  large enough, for every $\tomega$ such that $|\xi_{nN_0}(\tomega)|\leq R$, there is $\xi_n$ such that 
$$
\PP_n'(\xi_n)=z_n. 
$$ 
Moreover, $\xi_n=\xi_{nN_0}(\tomega)+o(1)$, as $T\to\infty$.  Without loss of generality, $|\xi_n|\leq R$ (otherwise increase $R$ a bit).  Note that in general, $\xi_n\neq \xi_{nN_0}(\tomega)$.

 Let $W_n$ be  a pure $E^s_n$-segments of depth $k$, and height $z_n(W_n)=z_n+O(1)$. By definition, $\pi(W_n)$ is a maximal $E^s_n$-segment in some cylinder 
$Q_n:=\<b^{(n)}_{n-k},\ldots,\dot{b}^{(n)}_n\>$. Let 
$$p_n(\wt{C})
:=\Prob_{W_n}[(\wt{\psi}_{n-1}\circ\cdots\circ\wt{\psi}_0)(\wt{C})].$$
Since $|\xi_n|\leq R$, equation \eqref{e.12.6} gives the following estimate  for $p_n(\wt{C})$:
\begin{align*}
p_n(\wt{C})&=[1+o(1)]\frac{e^{\PP_n(\xi_n)-\xi_nz_n}}{h_n(b^{(n)}_n,\xi_n)\sqrt{2\pi \wt{V}_n^{\xi_n}}}\times M_{\xi_n}(\wt{C})\times
\frac{\Prob[\wt{X}^{\xi_n}_{n-k}=b^{(n)}_{n-k},\ldots,\wt{X}^{\xi_n}_{n}=b^{(n)}_{n}]}{\Prob[X_{n-k}=b^{(n)}_{n-k},\ldots,X_{n}=b^{(n)}_{n}]}
\end{align*}

$\bullet$ 
 $\PP_n(\xi_n)-\xi_n z_n=-V_n\HH_n(\frac{z_n}{V_n})$, because $\PP_n'(\xi_n)=z_n$.   

$\bullet$ $\wt{V}^{\xi_n}_n$ is by definition the variance  of $S_n$ with respect to the change of measure $\{\wt{X}^{\xi_n}_i\}$. By Lemma 7.17 in \cite{Dolgopyat-Sarig-Book}, 
there is  a constant $C_1(R)>0$ such that for all $n$ sufficiently large, for every $\xi\in [-R,R]$,     $\wt{V}_n^{\xi}=C_1(R)^{\pm 1}V_n$. 

$\bullet$ $h_n(\cdot,\xi)=C_2(R)^{\pm 1}$ for some $C_2(R)>1$ and every $\xi\in [-R,R]$  , see Lemma \ref{LmAddPres}. 

$\bullet$ $M_\xi(\wt{C})=e^{\xi j}h_0(a_0,\xi)\Prob[\wt{X}^\xi_0=a_0,\ldots,\wt{X}^\xi_\ell=a_\ell]$ by definition. We saw above that if $|\xi|\leq R$, then  $e^{\xi j}h_0(a_0,\xi)=C_2(R)^{\pm 1} e^{\pm |j| R}$. Looking at the  transition kernel of the change of measure, and recalling that the jump functions are uniformly bounded,  we see that there is a constant $C_2'(R,\ell)>1$ such that for every $\xi\in [-R,R]$, 
$$
\Prob[\wt{X}^\xi_0=a_0,\ldots,\wt{X}^\xi_\ell=a_\ell]=C_2'(R,\ell)^{\pm 1} \Prob[X_0=a_0,\ldots,X_\ell=a_\ell].
$$ 
Thus, $M_\xi(\wt{C})=C_3(R,\ell)^{\pm 1}\mu_0(\wt{C})$ for all $|\xi|\leq R$.

$\bullet$ Similarly, there is a constant  $C_3'(R,k)>1$ such that 
\begin{align*}
&\Prob[\wt{X}^\xi_{n-k}={b}_{n-k}^{(n)},\ldots,\wt{X}^\xi_n={b}_n^{(n)}|\wt{X}^\xi_{n-k}={b}_{n-k}^{(n)}]\\
&\hspace{1cm}=C_2'(R,k)^{\pm 1} \Prob[{X}_{n-k}={b}_{n-k}^{(n)},\ldots,{X}_n={b}_n^{(n)}|{X}_{n-k}={b}_{n-k}^{(n)}].
\end{align*}
Let $\eps_0(R)$ be the common ellipticity constant of $\{\wt{X}_i^\xi\}$ for all $|\xi|\leq R$ (Lemma \ref{l.uniform-Doeblin}). By  
Corollary 2.10 in \cite{Dolgopyat-Sarig-Book}, there is a constant $C_4'(R)>1$ which  depends on $\eps_0(R)$, such that  for all $|\xi|\leq R$ and $k<n$, $\Prob[\wt{X}^\xi_{n-k}={b}_{n-k}^{(n)}]=C_4'(R)^{\pm 1} \Prob[{X}_{n-k}={b}_{n-k}^{(n)}].
$
Letting $C_4(R,k):=C_3'(R,k)C_4'(R)$, we find that 
$\frac{\Prob[\wt{X}^{\xi_n}_{n-k}=b^{(n)}_{n-k},\ldots,\wt{X}^{\xi_n}_{n}=b^{(n)}_{n}]}{\Prob[X_{n-k}=b^{(n)}_{n-k},\ldots,X_{n}=b^{(n)}_{n}]}
=C_4(R,k)^{\pm 1}$. 
 
(It is possible to show that $C_1,C_2,C_3\xrightarrow[R\to 0]{}1$, but we cannot  claim that $C_4\xrightarrow[R\to 0]{}1$. The issue is with $C_4'(R)$: While $\Prob[\wt{X}^\xi_{n-k}={b}_{n-k}^{(n)}]/ \Prob[{X}_{n-k}={b}_{n-k}^{(n)}]\xrightarrow[\xi\to 0]{}1$ for {\em fixed } $n$, we do not know how to show  that the limit holds {\em uniformly} in $n$.)

Let  $C(R,\ell,k):=\prod C_i$, then 
$$
p_n(\wt{C})=C(R,\ell,k)^{\pm 1}\frac{1}{\sqrt{2\pi V_n}}\exp\left({-V_{n}\HH_{n}\left(\frac{z_{n}}{V_{n}}\right)}\right)\mu_0(\wt{C}). 
$$

We now proceed as in the proof of Proposition \ref{PrBulk}. Let $A_T:=\{\phi^t(\tomega):0<t<T\}$. By the choice of  $n=n(T)$, 
$$
A_T^\ast:=(\tpsi_{n-1}\circ\cdots\circ\tpsi_0)(A_T)\text{ has length bounded away from zero and infinity}.
$$
Using Lemma \ref{l.approx-by-pure-segments}, we can approximate $A_T^\ast$ by finite unions of  pure $E^s_n$-segments $W_{n,1},\ldots,W_{n,m_n+2}$ of fixed depth, heights $z_n(W_{n,i})=z_n(\tomega)+O(1)$,  so that 
$$
\frac{1}{|A_T^\ast|}\sum_{i=1}^{m_n}|W_{n,i}|\Prob_{W_{n,i}}(\cdot)\leq \Prob_{A_T^\ast}(\cdot)\leq  \frac{1}{|A_T^\ast|}\sum_{i=1}^{m_n+2}|W_{n,i}|\Prob_{W_{n,i}}(\cdot),\text{ and }
$$
$$
e^{-\epsilon}<\frac{1}{|A_T^\ast|}\sum_{i=1}^{m_n}|W_{n,i}|<\frac{1}{|A_T^\ast|}\sum_{i=1}^{m_n+2}|W_{n,i}|<e^\epsilon.
$$
Recall from \eqref{e.sirena}   that $\frac{1}{T}\int_0^T 1_{\wt{C}}(\phi^t(\tomega))dt=\Prob_{A_T^\ast}(\wt{C})$. It follows that  for all $n$ large enough, if $|\xi_n|\leq R$, then 
$$
\frac{1}{T}\int_0^T 1_{\wt{C}}(\phi^t(\tomega))dt=C(R,\ell,k)^{\pm 1}e^{\pm 2\epsilon}\int 1_{\wt{C}}d\mu_0\cdot \frac{1}{\sqrt{2\pi V_n}}\exp\left({-V_{n}\HH_{n}\left(\frac{z_{n}}{V_{n}}\right)}\right). 
$$
If we keep $k$ and $\epsilon$ fixed (say equal to one), then $C(R,\ell,k)$ only depends on $\ell$, whence only on $\wt{C}$.
This proves the proposition for $h=1_{\wt{C}}$ with $\wt{C}$ as above. To get the statement for non-negative $h\in C_c(\St^\ast)$, we  sandwich $h$ between two finite linear combinations of indicators of $\wt{C}_i$ as above, with integrals close to $\int h d\mu_0$.
\end{proof}

 Artigiani, Castorinni, Ravotti, and Tumarkin  \cite{ACRT}
have obtained a  more precise
expansion for general $\Z^m$-covers of compact translation surfaces, but under the additional assumption that ${\beta\choose 1}$ is an eigenvector of some (single) automorphism of the cover, with eigenvalue of modulus strictly smaller than one. See \cite[Chapter 5]{Tumarkin-PhD} for additional discussion.

By \eqref{e.n(T)-again}, for every $n\in\N$ sufficiently large, there is $T>0$ such that $n(T)=n$. 
\begin{prop}\label{p.Large-Deviations} Suppose $n_j\uparrow\infty$ and choose $T_j\to\infty$ such that $n(T_j)=n_j$. \index{Birkhoff integrals of the linear flow}
\begin{enumerate}[(1)]
\item If $\lim\limits_{j\to\infty} \frac{z_{n_j}(\tomega)}{n_j}>0$, then 
$
\exists\kappa>0
$ s.t. $\frac{1}{T_j}\int_0^{T_j} 1_{[\zc<C]}(\phi^t(\wt{\omega}))dt=O(T_j^{-\kappa})$;

\item If $\lim\limits_{j\to\infty}\frac{z_{n_j}(\tomega)}{n_j}<0$, then 
$
\exists\kappa>0
$ s.t. $\frac{1}{T_j}\int_0^{T_j} 1_{[\zc>C]}(\phi^t(\wt{\omega}))dt
=O(T_j^{-\kappa})$.
\end{enumerate}
\end{prop}
\begin{proof}
Let $z_n:=z_{nN_0}(\tomega)$, and suppose $\lim z_{n_j}/n_j>0$.

Let $W_n$ be a sequence of pure $E^s_n$-segments with   depth $k$, heights $z_n(W_n)=z_n+O(1)$. Then    $\pi(W_n)$ is a maximal $E^s_n$-segment in some cylinder $Q_n=\<b^{(n)}_{n-k},\cdots,\dot{b}^{(n)}_n\>$. Let 
\begin{align*}
q_n(W_n):=\Prob_{W_n}\left\{\wt{\omega}'\in W_n: (\wt{\psi}_0^{-1}\circ\cdots\circ\wt{\psi}_{n-1}^{-1})(\wt{\omega}')\in [\zc<C]
\right\}.
\end{align*}

By Lemma \ref{Lemma-Modified-Approx-Canonical}, $|\wt{\zc}_0-\zc|, |\wt{\zc}_n-\zc|$ are  bounded by a constant independent of $n$. So  there exists $C'<0$ such that  
\begin{align*}
q_n(W_n)\leq \Prob_{W_n}\left\{\wt{\omega}'\in W_n: \wt{\zc}_0\bigl((\wt{\psi}_0^{-1}\circ\cdots\circ\wt{\psi}_{n-1}^{-1})(\wt{\omega}')
\bigr)
-\wt{\zc}_n(\wt{\omega}')<C'-z_n\right\}.
\end{align*}
By Lemma \ref{Lemma-f(X)}, $\DS \wt{\zc}_0\circ\tpsi_0^{-1}\circ\cdots\circ\tpsi_{n-1}^{-1}-\wt{\zc}_n=-\sum_{i=1}^{n} \wt{f}_i=-\sum_{i=1}^n f_i\circ X_i^{(n)}\circ\pi$. Therefore, $q_n(W_n)$ is a probability of an elementary event in $W_n$, and by Lemma \ref{l.prob-of-elementry-event}, 
\begin{align*}
q_n(W_n)\!\!\leq\! \Prob\!\left[-S_n<C'-z_n|X_{n-k}=b^{(n)}_{n-k},\ldots,X_n=b^{(n)}_n\right]
\!\!\leq\!\! \frac{\Prob\left[S_n>z_n-C'\right]}{\Prob[X_{n-k}=b^{(n)}_{n-k},\ldots,X_n=b^{(n)}_n]}.
\end{align*}

The denominator equals $\mu_0(Q_n)$, and is bounded below by a constant which only depends on $k$ (Lemma \ref{Lemma-Cylinder-Area}). Therefore, for some $c(k)>0$ independent of $n$, 
\begin{align*}
q_n(W_n)& 
\leq c(k)\Prob\left[S_n>z_n-C'\right].
\end{align*}
We will analyze these probabilities using Chernoff bounds. 

Recall that $V_n\!\!:=\!\!\Var(S_n)\!\!\asymp \!\!n$, and let  $\mathcal F_n(t)\!\!:=\!\!\frac{1}{V_n}\log \E(e^{t S_n})$. These are analytic functions, and by 
 Proposition \ref{PrRateFn},
there exist $c_0, C_0, n_0\!\!>\!\!0$ as follows. Let $\mathcal I_n$ denote the Legendre transform of $\mathcal F_n$; Then for all $n\!\!>\!\!n_0$, 
 $\mathcal I_n(\cdot)$ is well-defined, finite, and analytic on a neighborhood of 
 $(\frac{\E(S_N)}{V_N}-c_0,\frac{\E(S_N)}{V_N}+c_0)$, and for all 
 $\DS \eta\in \left[\frac{\E(S_N)}{V_N}-c_0,\frac{\E(S_N)}{V_N}+c_0\right]$, 
$$\DS
C_0^{-1}\left(\eta-\frac{\E(S_N)}{V_N}\right)^2\leq \mathcal I_N(\eta)\leq C_0\left(\eta-\frac{\E(S_N)}{V_N}\right)^2.
$$
Since $\lim (z_{n_j}/n_j)>0$ and  $V_n\asymp n$,   for all sufficiently large $j$,  $z_{n_j}/V_{n_j}>\delta>0$. Since $\E(S_n)=O(1)$, we may assume without loss of generality that
$
z_{n_j}-C'-|\E(S_{n_j})|>\delta V_{n_j}.
$
In this case, $q_{n_j}(W_{n_j})\leq c(k)\Prob[S_{n_j}>\E(S_{n_j})+\delta V_{n_j}]$.
By Chernoff's bound, 
\begin{align*}
q_{n_j}(W_{n_j})&\leq c(k) \exp\left[-V_{n_j}\mathcal I_{n_j}\left(\frac{\E(S_{n_j})}{V_{n_j}}+\delta\right)\right]
\leq c(k)\exp[-C_0^{-1}\delta^2 V_{n_j}].
\end{align*}

Let $A_T:=\{\phi^t(\wt{\omega}):0<t<T\}$. Since $n(T_j)=n_j$, 
$
A_{T_j}^\ast:=(\wt{\psi}_{n_j-1}\circ\cdots\circ\wt{\psi}_0)(A_{T_j})
$
has length bounded away from zero and infinity. So  $A_{T_j}^\ast$ can be covered by a bounded  number of pure $E^s_{n_j}$-segments $W_{n_j,1},\ldots,W_{n_j,m_j+2}$ of fixed depth $k$, as in Lemma \ref{l.approx-by-pure-segments}   (the value of $k$ does not matter).
Since  $m_j$ and  $|W_{n_j,i}|$ are uniformly bounded, 
$$
\frac{1}{T_j}\int_0^{T_j} 1_{[\zc<C]}(\phi^t(\wt{\omega}))dt=\Prob_{A_{T_j}}[\zc<C]\leq \sum_{i=1}^{m_j+2}|W_{n_j,i}|q_{n_j}(W_{n_j,i})=O(\exp(-C_0^{-1}\delta^2 V_{n_j})). 
$$
Since $V_{n_j}\asymp n_j\asymp\log T_j$, there is a positive constant $\kappa$ such that $C_0^{-1}\delta^2 V_{n_j}\geq \kappa \log T_j$, and the first part of the proposition follows. 
The second part has a similar proof.
\end{proof}

%

\subsection{Consequences (Theorems \ref{CrLebGen2}, \ref{c.diff-xi-diff-asymp}, \ref{CrLogOTDetXi}, \ref{PrBalance-second half},\ref{ThLT-ZD}, \ref{ThLT-NZD})} We use the notation of the previous section.

\medskip
\noindent
{\bf Proof of Theorem \ref{CrLebGen2}.} The theorem asserts the equivalence of the conditions
\begin{enumerate}[(a)]
\item[(a)] $\tomega$ is $\mu_0$-generic;
\item[(b)] $\DS \lim_{T\to \infty} \frac{\log \Big( \int_0^T 1_{\St_-} (\phi^t\tomega) dt\Big)}
{\log \Big(\int_0^T 1_{\St_+} (\phi^t \tomega) dt\Big)}=1$, where  
$\DS \St_+:=\bigcup_{k\geq 0}\Rect_k$, $\DS \St_-:=\bigcup_{k<0}\Rect_k$;
\item [(c)] $\DS \lim_{T\to \infty} \left(\log  \int_0^T 1_{\Rect_0} (\phi^t \tomega) dt\right)\bigg/\log T=1$. 
\end{enumerate}
Throughout this proof, $z_n:=z_{nN_0}(\tomega)$, and $\xi_n:=\xi_{nN_0}(\tomega)$. 

\medskip
\noindent
{\bf (a)$\boldsymbol{\Rightarrow}$(b):}   Since $\tomega$ is $\mu_0$-generic, $z_n/n\to 0$, and $\xi_n\to 0$ (Theorems \ref{CrLebGen}, \ref{t.Char-mu-xi-Gen}).

It is sufficient to prove (b) along any sequence $T_k\to\infty$.   Let $n_k=n(T_k)$, where $n(\cdot)$ denotes the function 
 defined by \eqref{e.n(T)-again}. Set $z_{n_k}:=z_{n(T_k)}$. 

Suppose first that
$
z_{n_k}\geq 0
$ for all $k$ large.

Then 
\begin{align*}
&\frac{1}{T_k}\int_0^{T_k}\!\!\! 1_{\St_+}(\phi^t(\wt{\omega}))dt
=\Prob_{A_{T_k}}[\zc\geq 0]=\Prob_{A_{T_k}} \!\!\left(\frac{\zc-z_{n_k}}{\sqrt{V_{n_k}}}
\geq -\frac{z_{n_k}}{\sqrt{V_{n_k}}}\right)
\geq \Prob_{A_{T_k}}\!\!\left(\frac{\zc-z_{n_k}}{\sqrt{V_{n_k}}}
\geq 0\right).
\end{align*}
By the temporal CLT, the last probability tends to $\frac{1}{2}$. 
So if $k$ is large enough, then  $\frac{T_k}{3} \leq \int_0^{T_k} 1_{\St_+}(\phi^t(\wt{\omega}))dt \leq T_k$, and 
\begin{equation}\label{e.tocatta}
\log \left(\int_0^{T_k} 1_{\St_+}(\phi^t(\wt{\omega}))dt\right)=\log T_k +O(1).
\end{equation}
Fix some  $h\in C_c(\St^\ast)$ such that $0\leq h\leq 1_{\St_-}$ and $\int h d\mu_0>0$. By 
Proposition \ref{p.Bikhoff-Sums-Estimate},
\begin{align*}
1&\geq \frac{1}{T_k}\int_0^{T_k}1_{\St_-}(\phi^t(\wt{\omega}))dt
\geq \frac{1}{T_k}\int_0^{T_k} h(\phi^t(\wt{\omega}))dt
\geq  \frac{C}{\sqrt{V_{n_k}}}e^{-V_{n_k} \HH_{n_k}(\frac{z_{n_k}}{V_{n_k}})}
\end{align*}
with $C>0$.
By Lemma \ref{l.rate-function}(3) and the estimates 
 $\E(S_{n})=O(1)$, $V_{n}\asymp n$,  $z_n/n\to 0$, 
 \begin{equation}\label{e.HH-o(V_n)}
\HH_{n}\left(\frac{z_{n}}{V_n}\right)\asymp \left(\frac{z_n-\E(S_n)}{n}\right)^2\xrightarrow[n\to\infty]{}0.
\end{equation}
Since $V_{n}\asymp n(T)\asymp \log T$, this leads to the bounds 
$$
T_k\geq \int_0^{T_k} 1_{\St_-}(\phi^t(\wt{\omega}))dt\geq \frac{CT_k}{\sqrt{V_{n_k}}}\exp[o(V_{n_k})]
=\exp\left[\log T_k+o(\log T_k)\right].
$$
Thus $\DS
\log \left(\int_0^{T_k} 1_{\St_-}(\phi^t(\wt{\omega}))dt\right)=
\log T_k+o(\log T_k)
$. Combining this with \eqref{e.tocatta} yields the limit in (b), along any sequence $T_k\to\infty$ for which  
{ $z_{n_k}\geq 0$ for all $k$ large.}

A similar argument gives  $
\log \left(\int_0^T 1_{\St_-}(\phi^t(\wt{\omega}))dt\right)\!\!=\!\!\log T\!+\!O(1)
$, $
\log \left(\int_0^T 1_{\St_+}(\phi^t(\wt{\omega}))dt\right)=
\log T+o(\log T)$, and (b) along any sequence $T_k\to\infty$ such that  $z_{n_k}\leq 0$ for all $k$ large.

 The remaining case is when $z_{n_k}=z_{n(T_k)}$ changes sign infinitely often. Split $(T_k)_{k\geq 1}=(T_k')_{k\geq 1}\cup (T_k'')_{k\geq 1}$  so that $z_{n(T_k')}\geq 0$ and  $z_{n(T_k'')}\leq 0$ for all $k$. By the previous argument, the limit in (b) holds along $(T_k')_{k\geq 1}$ and $(T_k'')_{k\geq 1}$, hence also along $(T_k)_{k\geq 1}$.

\medskip
\noindent
{\bf (b)$\boldsymbol{\Rightarrow}$(a)}: 
Assume by contradiction that  $\exists n_j\to\infty$ such that $\lim (z_{n_j}/n_j)>0$.  Choose $T_j$ so that $n(T_j)=n_j$. 
We saw above that 
$\DS 
 \log\left(\int_0^{T_j}1_{\St_+}(\phi^t(\wt{\omega}))dt\right)=\log T_j+O(1).
$
On the other hand, by Proposition \ref{p.Large-Deviations},  for some $\kappa>0$ and all $j$ large enough, 
$\DS
\frac{1}{T_j} \int_0^{T_j}1_{\St_-}(\phi^t(\wt{\omega}))dt=O(T_j^{-\kappa})
$, whence $\DS  \log \int_0^{T_j}1_{\St_-}(\phi^t(\wt{\omega}))dt\leq (1-\kappa)\log T_k+O(1)$.
But this implies that  
$\liminf\limits_{j\to\infty}
\log\left(\int_0^{T_j}1_{\St_+}(\phi^t(\wt{\omega}))dt\right)
\bigg/\log\left(\int_0^{T_j}1_{\St_-}(\phi^t(\wt{\omega}))dt\right)>1$, in contradiction to (b).  This proves that $\limsup\limits_{n\to\infty}z_n/n\leq 0$. Similarly, one shows that $\liminf\limits_{n\to\infty}z_n/n\geq 0$, whence $z_n/n\to 0$, and $\tomega$ is $\mu_0$-generic, by Theorem \ref{CrLebGen}.

\medskip
\noindent
{\bf (a)$\boldsymbol{\Rightarrow}$(c)}: Fix some continuous function  $0\leq h\leq 1_{\Rect_0}$ with compact support and with positive integral. Since  $\tomega$ is $\mu_0$-generic, $z_n/n\to 0$ and $\xi_n\to 0$ (Theorem \ref{CrLebGen}).
 By Proposition \ref{p.Bikhoff-Sums-Estimate} and \eqref{e.HH-o(V_n)}, for some positive $C$ and all $T$ large enough,
\begin{align*}
1&\geq \frac{1}{T}\int_0^T 1_{\Rect_0}(\phi^t(\wt{\omega}))dt
\geq \frac{1}{T}\int_0^T h(\phi^t(\wt{\omega}))dt\geq  \frac{C}{\sqrt{2\pi V_{n(T)}}}\exp[o(V_{n(T)})].
\end{align*}
Since $V_n\asymp n$ and $n(T)\asymp \log T$,  $\log \int_0^T 1_{\Rect_0}(\phi^t(\wt{\omega}))dt=\log T+o(\log T)$, and (c) follows.

\smallskip
\noindent
{\bf (c)$\boldsymbol{\Rightarrow}$(a):} Fix a subsequence $n_j\uparrow\infty$ such that $\frac{z_{n_j}}{n_j}$ converges to a limit $z\in [-\infty,\infty]$. 
Assume by contradiction that $z>0$. Then Proposition \ref{p.Large-Deviations}   gives us a constant $\kappa>0$ and $T_j\to\infty$ such that 
$\DS
\frac{1}{T_j}\int_0^{T_j}1_{\Rect_0}(\phi^t(\tomega))dt\leq 
\frac{1}{T_j}\int_0^{T_j}1_{[\zc<1]}(\phi^t(\tomega))dt=O(T_j^{-\kappa})
$. But this contradicts (c). Similarly, it is impossible that $z<0$. 
Thus the only possible limit point of $z_n/n$ is zero. By Theorem \ref{CrLebGen}, $\tomega$ is $\mu_0$-generic.\qed

\smallskip
\noindent
{\bf 
Proof of Theorem \ref{c.diff-xi-diff-asymp}.} This theorem follows from the following statement: If  $\xi_1>\xi_2\geq 0$ or $\xi_1<\xi_2\leq 0$,  then there is some $r>0$ such that for any pair of $\mu_{\xi_i}$-generic point  $\tomega_i$, and for every non-identically zero and  non-negative function $h\in C_c(\St)$, for all $T$ large enough, 
$\DS
\frac{\int_0^T h(\phi^t(\tomega_2))dt}{\int_0^T h(\phi^t(\tomega_1))dt}>T^r.
$

Let $n=n(T)$, with $n(T)$ as in \eqref{e.n(T)-again}. By Proposition \ref{Prop-Var-Growth},  $V_n\asymp n\asymp\log T$.

 Let $z_{n}^{(i)}:=z_{nN_0}(\tomega_i)$. 
By Proposition \ref{p.Bikhoff-Sums-Estimate},  for all $T$ large enough, 
\begin{align*}
&\log\left(\frac{\int_0^T h(\phi^t(\tomega_2))dt}{\int_0^T h(\phi^t(\tomega_1))dt}\right)\geq V_n\left[\HH_n\left(\frac{z_{n}^{(1)}}{V_n}\right)-\HH_n\left(\frac{z_{n}^{(2)}}{V_n}\right)\right]+O(1).
\end{align*}
To prove the proposition, it is sufficient to show that $\liminf\limits_{T\to\infty}\left[\HH_n\left(\frac{z_{n}^{(1)}}{V_n}\right)-\HH_n\left(\frac{z_{n}^{(2)}}{V_n}\right)\right]>c$, where $c$  is a positive constant independent of $\tomega_i$. Suppose first that $\xi_1>\xi_2\geq 0$.

Since $\tomega_i$ is $\mu_{\xi_i}$-generic, 
 $z_n^{(i)}=\mathcal P'_{nN_0}(\xi_{nN_0}(\tomega_i))$, where $\xi_{nN_0}(\tomega_i)\to \xi_i$. 
 By Proposition~\ref{p.Pressure-and-Geometric-Pressure}, $\exists\xi_n^{(i)}\to\xi_i$ such that $\PP_n'(\xi_n^{(i)})=z_n^{(i)}$.  
 
Fix some $R\!\!>\!\!0$ such that $\xi_n^{(i)}\!\!\in\!\! [-R/2,R/2]$ for all $n$ sufficiently large. 
By Lemma~\ref{LmPressure}, there is  some $\rho_1>0$ such that 
$
\rho_1^{-1}V_n\leq \PP_n''\leq \rho_1 V_n$ on $[-R,R]$. 
Since $\xi_1>\xi_2\geq 0$, 
\begin{align*}
z_n^{(2)}-\E(S_n)&=\frac{\PP_n'(\xi_2+o(1))-\PP_n'(0)}{\xi_2+o(1)}(\xi_2+o(1))
\geq \rho^{-1}_1(\xi_2+o(1))V_n\geq o(V_n),\\
z_n^{(1)}-z_n^{(2)}&=\PP_n'(\xi_1+o(1))-\PP_n'(\xi_2+o(1))
\geq \rho^{-1}_1(\xi_1-\xi_2+o(1))V_n.
\end{align*}

Since  $\xi_n^{(i)}\in[-R/2,R/2]$ for all $n$ sufficiently large, 
$
\frac{z_n^{(i)}}{V_n}\in [\bar{a}_n^R, \bar{b}_n^R]$ for all $n$ sufficiently large,  
where $[\bar{a}_n^R, \bar{b}_n^R]$ are the intervals in Lemma \ref{l.rate-function}. 
By Lemmas \ref{LmPressure} and  \ref{l.rate-function}, there is  some $\rho_2>0$ such that 
$\rho^{-1}_2\leq \HH_n''\leq \rho_2$  on $[\bar{a}_n^R, \bar{b}_n^R]$. It follows that 

\begin{align*}
&\HH_n\left(\frac{z_{n}^{(1)}}{V_n}\right)-\HH_n\left(\frac{z_{n}^{(2)}}{V_n}\right)=\int_{{z_{n}^{(2)}}/{V_n}}^{z_{n}^{(1)}/{V_n}}\int_{\E(S_n)/V_n}^t \HH_n''(s)dsdt\ (\text{because }\HH_n'\left(\frac{\E(S_n)}{V_n}\right)=0)\\
&\geq \frac{1}{2}\rho_2^{-1}\left[\left(\frac{z_n^{(1)}-\E(S_n)}{V_n}\right)^2-\left(\frac{z_n^{(2)}-\E(S_n)}{V_n}\right)^2\right]\\
&=\frac{1}{2\rho_2}\times\frac{z_n^{(1)}-z_n^{(2)}}{V_n}\times  \left(\frac{z_n^{(1)}-z_n^{(2)}}{V_n}+\frac{2(z_n^{(2)}-\E(S_n))}{V_n}\right)
\end{align*}
\begin{align*}
&\geq \frac{1}{2\rho_2}\times \frac{\xi_1-\xi_2+o(1)}{\rho_1}\times\left(\frac{\xi_1-\xi_2+o(1)}{\rho_1}+\frac{2o(V_n)}{V_n}\right)\\
&>\frac{1}{2\rho_1\rho_2^2}(\xi_1-\xi_2)^2+o(1)>c:=\frac{(\xi_1-\xi_2)^2}{3\rho_1\rho_2^2},\text{ for all $n$ sufficiently large.}
\end{align*}

This proves  result when $\xi_1\!\!>\!\!\xi_2\!\!\geq\!\! 0$.
If  $\xi_1\!\!<\!\!\xi_2\!\!\leq\!\! 0$, we show first that $z_n^{(2)}-\E(S_n)\leq -o(V_n)$ and $z_n^{(1)}-z_n^{(2)}\leq -\rho_1^{-1}(|\xi_1-\xi_2|+o(1))V_n$, and proceed as before. 
\qed

\smallskip
\noindent
{\bf Proof of Theorem \ref{CrLogOTDetXi}.} 
Given $\xi>0$, denote
 $\bbI_T(\tomega)\!\!:=\!\! \log \left(\int_0^T 1_{[\zc\leq 1]}(\phi^t \tomega) dt\right)$,   and $\bbI_T^\xi\!\!:=\!\!\!\int_{\Rect_0} \bbI_T(\omega) d\mu_\xi(\omega).$ The theorem says that 
\begin{enumerate}[(1)] 
\item For every $\mu_\xi$-generic $\tomega$, 
$\DS \lim_{T\to\infty}\!\! \frac{\bbI_T(\omega)-\bbI_T^\xi}{\log T}=0$;
\item If $\xi_1>\xi_2>0$, then $\liminf\limits_{T\to\infty}\frac{\bbI_T^{\xi_2}-\bbI_T^{\xi_1}}{\log T}>0$;
\item  $c(\xi)\log T\leq \bbI_T^\xi\leq \log T$ for all $T$ large enough, with $c(\xi)>0$.
 \end{enumerate}

We begin by showing that part (1) holds for $\mu_\xi$-a.e. $\tomega$ and then prove that it actually holds 
for all $\mu_\xi$-generic points.

Since $\xi>0$, $\DS \mu_\xi[\zc\leq 1]=\sum_{n\leq 1}\mu_{\xi}(D^{-\xi}(\Rect_0))=\sum_{n\geq 1}e^{-n\xi}<\infty$, and $1_{[\zc\leq 1]}\in L^1(\mu_\xi)$.  Fix some $h\in L^1(\mu_\xi)$ as in Proposition \ref{p.Bikhoff-Sums-Estimate}. 
 
For $\mu_\xi$-a.e.  $\tomega$, $\tomega$ is $\mu_\xi$-generic {\em and } the ratio ergodic theorem holds at $\tomega$ for the pair $(h,1_{[\zc\leq 1]})$ (even though $1_{[\zc\leq 1]}\not\in C_c(\St)$). By Theorem \ref{t.Char-mu-xi-Gen},  $\xi_{nN_0}(\tomega)\to\xi$, and by Theorem \ref{CrDriftSeq} there is a sequence of constants $z_n^\xi$ independent of $\tomega$ such that $z_{nN_0}(\tomega)=z_n^\xi+o(n)$. Moreover, we may take $z_n^\xi:=z_{nN_0}(\tomega_0)$ where $\tomega_0$ is  $\mu_\xi$-generic, and $\xi_{nN_0}(\tomega_0)\to\xi$.  
 By the ratio ergodic theorem for $(h,1_{[\zc\leq 1]})$,  for $\mu_\xi$-a.e.  $\tomega$, $\bbI_T(\tomega)=\log\int_0^T h(\phi^t(\tomega))dt+o(1)$, and by Proposition \ref{p.Bikhoff-Sums-Estimate}
\begin{equation}\label{e.brahms-2}
\bbI_T(\tomega)=\log T-V_n\HH_n\left(\frac{z_{nN_0}(\tomega)}{V_n}\right)+o(\log T).
\end{equation}

We would like to replace $z_{nN_0}(\tomega)$ by the constants $z_n^\xi$. Let us estimate the error this introduces.  
By Proposition \ref{p.Pressure-and-Geometric-Pressure}, $\exists\eta_n=\xi+o(1)$ such that $\PP_n'(\eta_n)=z_n^\xi$. Looking at \eqref{Pres1stDer}, \eqref{Pres2dDer}, and recalling that $\E(S_n)=O(1)$ and $V_n\to\infty$, we deduce that 
\begin{equation}\label{e.z-div-V}
\frac{z_n^\xi}{V_n}=\frac{\PP_n'(\eta_n)-\PP_n'(0)}{V_n}=\frac{1}{V_n}\int_0^{\xi+o(1)}\PP_n''(t)dt\asymp 1.
\end{equation}
Recall that $[\bar{a}_n^R,\bar{b}_n^R]\!\!=\!\![\PP_n'(-R)/V_n, \PP_n'(R)/V_n]$.
 Clearly, $\E(S_n)/{V_n}\!\!=\!\!\PP_n'(0)/V_n\!\in\! [\bar{a}_n^R,\bar{b}_n^R]$ for all $R>0$, and $\eta_n\in [\bar{a}_n^R,\bar{b}_n^R]$ for all $R>\xi$ and $n$ large enough. 
Fix $R,n_0>0$ such that   $\xi_{nN_0}(\tomega),\xi_{nN_0}(\tomega_0),\eta_n\in [-R/2,R/2]$ for all $n>n_0$. Then 
for all $n$ large enough,  
$$\frac{\E(S_n)}{V_n}, \frac{z_{nN_0}(\tomega)}{V_n}, \frac{z_n^\xi}{V_n}\in [\ov{a}_n^R,\ov{b}_n^R].$$ 
By Lemma \ref{l.rate-function},
$\HH_n\!\left(\frac{\E(S_n)}{V_n}\right)\!\!=\!0$, 
$\HH_n'\!\left(\frac{\E(S_n)}{V_n}\right)\!\!=\!0$, and $\rho(R)^{-1}\!\!\leq\!\! \HH_n''\!\!\leq\!\! \rho(R)\text{ on }[\bar{a}_n^R,\bar{b}_n^R].
$
So  
\begin{align}
&\biggl|\HH_n(\frac{z_{nN_0}(\tomega)}{V_n})-\HH_n(\frac{z_n^\xi}{V_n})\biggr|
\leq \left|\int_{z_{nN_0}(\tomega)/V_n}^{z_n^\xi/V_n} \int_{\E(S_n)/V_n}^t \HH_n''(s)dsdt \right|\notag\\
&\leq \frac{1}{2}\rho(R) \left(\frac{z_n^\xi-z_{nN_0}(\tomega)}{V_n}\right)^2=\frac{1}{2}\rho(R) \left(\frac{o(n)}{V_n}\right)^2=o(1),\text{ since $V_n\asymp n$}.\label{e.Rod}
\end{align}

Substituting this estimate in \eqref{e.brahms-2}, we obtain 
\begin{equation}
\label{bbIToH}
\bbI_T(\tomega)=\log T-V_n\HH_n\bigg(\frac{z_n^\xi}{V_n}\bigg)+\eps_T(\tomega), 
\end{equation}
where $\eps_T(\tomega)=o(\log T)$ (because $V_n\asymp n\asymp\log T$). In addition,
\begin{equation}
\bbI_T^\xi=\log T-V_n\HH_n\bigg(\frac{z_n^\xi}{V_n}\bigg)+\int_{\Rect_0}\eps_T(\tomega)d\mu_\xi.\label{e.wanderer}
\end{equation}

Thus, $\DS\left|\frac{\bbI_T(\tomega)-\bbI_T^\xi}{\log T}\right|\leq \frac{|\eps_T(\tomega)|}{\log T}+ \int_{\Rect_0}\frac{|\eps_T(\tomega)|}{\log T} d\mu_\xi(\tomega)$.
We already know that $\frac{|\eps_T(\tomega)|}{\log T}\to 0$ $\mu_\xi$-a.e. To see that $\int_{\Rect_0}\frac{|\eps_T(\tomega)|}{\log T} d\mu_\xi(\tomega)\to 0$, we prove  that $\|\eps_T\|_\infty =O(\log T)$ and apply the bounded convergence theorem.
  Here is the proof. Arguing as above, we see that since  $\DS \frac{\E(S_n)}{V_n}=o(1)$ and  $\DS \frac{z_n^\xi}{V_n}\asymp 1$, 
\begin{align}
&\DS\HH_n(z_n^\xi/V_n)= \int_{\E(S_n)/V_n}^{z_n^\xi/V_n} \int_{\E(S_n)/V_n}^t \HH_n''(s)dsdt \asymp \left(\frac{z_n^\xi}{V_n}-\frac{\E(S_n)}{V_n}\right)^2\asymp 1.\label{e.H_n-O(1)}
\end{align}
Since  $V_n\asymp n\asymp\log T$, there exists a positive constant $M$ so large that  
$$\left|\log T-V_n\HH_n\left(\frac{z_n^\xi}{V_n}\right)\right|< M\log T.$$ 
 On $\Rect_0$, $0\leq \bbI_T(\tomega)\leq \log T$. Therefore,   for all $T$ large enough  $\|\eps_T\|_\infty\leq (M+1)\log T$.
 Now we can apply the bounded convergence theorem and complete the proof that  $\left|\frac{\bbI_T(\tomega)-\bbI_T^\xi}{\log T}\right|\to 0$ for $\mu_\xi$-a.e. $\tomega$.

We will now prove that this statement holds for all $\mu_\xi$-{\em generic} $\tomega$.

Taking a compactly supported continuous function $h$ with compact support inside $\{\zc\leq 1\}$ 
so that 
$\bbI_T(\omega)\geq \log\left(\int_0^T h(\phi^t \tomega) dt\right) $,
we get by the foregoing argument that for each $\mu_\xi$-generic $\tomega$
\begin{equation}
\label{bbILower}
\bbI_T(\tomega)\geq \log T-V_n\HH_n\bigg(\frac{z_n^\xi}{V_n}\bigg)+o(\log T).
\end{equation}
On the other hand, invoking Lemma \ref{l.approx-by-pure-segments} and arguing as in the proof of Theorem 
\ref{ThStairTLT}
we obtain that there is a constant $K$ such that
\begin{equation}
\label{BeckLike}
 \bbI_T(\tomega)-\log T\leq \log \left( \Prob(S_n\geq z_{nN_0}(\omega)+K)\right)+O(1). 
 \end{equation}
Let $\mathcal F_n(\xi):=\frac{1}{V_n}\log\E(e^{\xi S_n})$. By Chernoff's bound, for each $\xi\in\R$ we have
\begin{align*} 
&\Prob(S_n\geq z_{nN_0}(\omega)+K)
=\Prob\left( e^{\xi S_n}\geq  e^{\xi (z_{n N_0}+K)}\right)
\leq \exp\left(V_n \left[ \mathcal F_n(\xi)-\xi \frac{z_{n N_0} (\omega)+K}{V_n}\right]\right) \\
&\leq \exp\left(V_n \left[ \frac{P_n(\xi)+O(1)}{V_n}-\xi \frac{z_{n N_0} (\omega)+K}{V_n}\right]\right),
\text{ see \eqref{Pres0thDer}.}
\end{align*}
Taking the supremum over $\xi$ gives
\begin{equation}
\label{Chernoff}
\log \Prob(S_n\geq z_{nN_0}(\omega)+K)\leq -V_n\HH_n \left(\frac{z_{nN_0}(\tomega)+K}{V_n}\right) +O(1).
\end{equation}
The computation leading to \eqref{bbIToH} now shows that the right hand side equals to
$\DS
-V_n\HH_n\left(\frac{z_n^\xi}{V_n}\right)+o(V_n)
$.\footnote{This calculation only requires the assumption that $z_{nN_0}(\tomega)=z_n^\xi+o(1)$, hence it applies to all $\mu_\xi$-generic points.}

Hence combining \eqref{BeckLike}, \eqref{Chernoff} and the estimate $V_n\asymp n\asymp\log T$ gives
\begin{equation}
\label{bbIUpper}
\bbI_T(\tomega)\leq \log T-V_n\HH_n\bigg(\frac{z_n^\xi}{V_n}\bigg)+o(\log T).
\end{equation}
Taking into account \eqref{bbILower} we obtain that for every $\mu_\xi$-generic $\tomega$, 
$$\bbI_T(\tomega)= \log T-V_n\HH_n\bigg(\frac{z_n^\xi}{V_n}\bigg)+o(\log T).$$
Comparing this with \eqref{e.wanderer}, and recalling that $\int_{\Rect_0}\frac{\epsilon_T}{\log T} d\mu_\xi\!\to\! 0$, we conclude that  $\left|\frac{\bbI_T(\tomega)-\bbI_T^\xi}{\log T}\right|\to 0$ for each $\mu_\xi$-generic $\tomega$. This is part (1) of the theorem.

Next we prove (2). Pick some non-negative $h\in C_c(\St)$ such that $\int h d\mu_{\xi_i}\neq 0$ $(i=1,2)$. Since $h,1_{[\zc\leq 1]}\in L^1(\mu_{\xi_i})$,  we can apply the ratio ergodic theorem and deduce that for $\mu_{\xi_i}$-a.e. $\tomega_i$,  $\log\int_0^T h(\phi^t(\tomega_i))dt=\bbI_T(\tomega_i)+o(1)$. Thus, by part (1) proved above, 
$$
\frac{\log\int_0^T h(\phi^t(\tomega_i))dt-\bbI_T^{\xi_i}}{\log T}=0 \ \text{ for $\mu_{\xi_i}$-a.e. $\tomega_i$, whence } 
$$
$$
\lim_{T\to\infty}\frac{\log\int_0^T h(\phi^t(\tomega_2))dt-\log\int_0^T h(\phi^t(\tomega_1))dt}{\log T}+\frac{\bbI_T^{\xi_1}-\bbI_T^{\xi_2}}{\log T}=0 \ \text{ for $(\mu_{\xi_1}\times\mu_{\xi_2})$-a.e. $(\tomega_1,\tomega_2)$.} 
$$
By Theorem \ref{c.diff-xi-diff-asymp}, there is a positive constant $r>0$ such that 
$$
\liminf_{T\to\infty}\frac{\log\int_0^T h(\phi^t(\tomega_2))dt-\log\int_0^T h(\phi^t(\tomega_1))dt}{\log T}>r,  \ \text{ for $(\mu_{\xi_1}\times\mu_{\xi_2})$-a.e. $(\tomega_1,\tomega_2)$.} 
$$
Necessarily, 
$\DS
\liminf_{T\to\infty}\frac{\bbI_T^{\xi_2}-\bbI_T^{\xi_1}}{\log T}>r
$, and we proved part (2).

Finally, we  prove (3). The upper  bound $\bbI_T^\xi\leq\log T$ is trivial. To prove the lower bound,   we bound the right-hand-side of \eqref{e.wanderer} from below.
 Fix some continuous function $h\geq 0$, with compact support and positive integral. 
By Theorem \ref{c.diff-xi-diff-asymp}, there exists $0<r<1$ such that for $\mu_\xi$-a.e. $\tomega$, and $\mu_{\xi/2}$-a.e. $\tomega'$, for all $T$ large enough, 
$\DS
\int_0^T h(\phi^t(\tomega))dt>T^r\int_0^T h(\phi^t(\tomega'))dt.
$
Necessarily, for all $T$ large enough, for $\mu_\xi$-a.e. generic point $\tomega$, 
$\DS
\int_0^T h(\phi^t(\tomega))dt>T^r.
$

Together with Proposition \ref{p.Bikhoff-Sums-Estimate}, this implies that for $\mu_\xi$-a.e. $\tomega$, for all $T$ large enough, 
$
\DS T^{r-1}<\const.\frac{1}{\sqrt{\log T}}\exp\left(-V_n \HH_n\left(\frac{z_{nN_0}(\tomega)}{V_n}\right)\right)
$, 
whence 
$$
V_n\HH_n\left(\frac{z_{nN_0}(\tomega)}{V_n}\right)
<\log\left(\const \frac{T^{1-r}}{\sqrt{\log T}}\right)<(1-r)\log T+o(\log T). 
$$
By  \eqref{e.Rod}, 
$\DS
V_n\HH_n\left(\frac{z_{n}^\xi}{V_n}\right)
<(1-r)\log T+o(\log T)+o(V_n) 
$. Since $V_n\asymp n\asymp\log T$, 
$$
V_n\HH_n\left(\frac{z_{n}^\xi}{V_n}\right)<(1-\frac{r}{2})\log T\text{ for all $T$ large}. 
$$
Plugging this at \eqref{e.wanderer}, 
and recalling from the proof of  (1) that 
$\int_{\Rect_0} |\eps_T(\tomega)|d\mu_\xi\!\!=\!\!o(\log T)$, we find that for all $T$ large enough, 
$
\log T\geq \bbI_T^\xi\geq \frac{r}{3}\log T.
$ 
\qed

%
%
%

\smallskip
\noindent
{\bf Proof of Theorem \ref{PrBalance-second half}.} 
The theorem states that for Lebesgue a.e. $\tomega$, for  every $\eps>0$, for all $T$ large enough, $\int_0^T 1_{\St_{-}}(\phi^t(\tomega))dt>\frac{T}{(\log\log T)^{1+\eps}}$, and there is a subsequence $T_j\to\infty$ so that $\int_0^{T_j} 1_{\St_{-}}(\phi^t(\tomega))dt<\frac{T_j}{(\log\log T_j)^{1-\eps}}$.

We begin with two preparatory sublemmas.
Let $V_n:=\Var(S_n)$. 
\setcounter{sublem}{0}
\begin{sublem}\label{l.LIL-for-z}
For $\mu_0$-a.e. $\tomega\in\St$, \index{Law of the iterated logarithm}
$$
\limsup_{n\to\infty}\frac{z_{nN_0}(\tomega)}{\sqrt{2V_n\log\log V_n}}=+1\ , \ \liminf_{n\to\infty}\frac{z_{nN_0}(\tomega)}{\sqrt{2V_n\log\log V_n}}=-1.
$$
\end{sublem}
\begin{proof}
Recall that $D:\St\to\St$ is the automorphism which maps $\Rect_k$ to $\Rect_{k+1}$ by translation. Then $z_{nN_0}\circ D=z_{nN_0}+1$, and it is therefore  sufficient to prove the statement on  $\Rect_0$. 
On $\Rect_0$, $z_{nN_0}(\tomega)=\Phi_n(\tomega)$, and by Proposition \ref{Prop-Markov-Representation-of-Phi}, there is a uniformly bounded sequence of functions $\bar{\eps}_n(\tomega)$ such that  $\Phi_n-\bar{\eps}_n$ is equal in distribution to $S_n$.  
Thus the lemma reduces to the assertion that 
$$
\limsup_{n\to\infty}\frac{S_n}{\sqrt{2V_n\log\log V_n}}=+1\ , \ \liminf_{n\to\infty}\frac{S_n}{\sqrt{2V_n\log\log V_n}}=-1\text{ almost surely}.
$$
Recall that  $\{X_n\}$ is uniformly elliptic, $f_n(X_n)$ are uniformly bounded,  $\E(S_n)=O(1)$, and $V_n\to\infty$. Applying the law of the iterated logarithm for additive functionals of inhomogeneous Markov chains (Theorem \ref{t.LIL}), we obtain the result.
\end{proof}

\begin{sublem}\label{Lem-LIL-Asymp} 
For every $\eps>0$, for all $n$ sufficiently large, for any sequence $a_n\in\mathfrak P_n$, 
$$
\frac{1}{(\log V_n)^{1+\eps}}\leq 
\Prob\left(S_n\geq \sqrt{(2+o(1))V_n\log\log V_n}\, \bigg|\,  X_n=a_n\right)\leq \frac{1}{(\log V_n)^{1-\eps}}
$$
\end{sublem}
\begin{proof}
First we prove the sublemma without the conditioning on $X_n$.

Fix $\eps>0$. 
Let $\mathcal F_n(\xi):=\frac{1}{V_n}\log \E(e^{tS_n})$ and let $\mathcal I_n(\eta)$ be the Legendre transform of $\mathcal F_n(\xi)$. 
By the local limit theorem for large deviations (Theorem \ref{ThLD} in Appendix \ref{AppMC}), there are constants $C>1$, $\delta>0$ and $n_0>1$ so that for all  $n>n_0$, 
$$
\DS \Prob[S_n=k]=C^{\pm 1}(2\pi V_n)^{-\frac{1}{2}}\exp\left[-V_n\mathcal I_n\left(\frac{k}{V_n}\right)\right]\text{ for all $k\in\Z$ s.t. $\frac{k-\E(S_n)}{V_n}\in (-2\delta,2\delta)$}.
$$
{Otherwise  Theorem \ref{ThLD}  would fail along a sequence $k_n\in\Z$ such that $\frac{k_n-\E(S_n)}{V_n}\to 0$. }

By the uniform estimates of $\mathcal I_n$ in Theorem \ref{PrRateFn}, we may choose $\delta$ sufficiently small  so that for every $k\in\Z$ such that $\frac{k-\E(S_n)}{V_n}\in (-2\delta,2\delta)$,
$$
V_n\mathcal I_n\left(\frac{k}{V_n}\right)=\frac{1\pm(\eps/3)}{2}\left(
\frac{k-\E(S_n)}{\sqrt{V_n}}
\right)^2.
$$ 

Recall that $\E(S_n)=O(1)$ and $V_n\to\infty$. Therefore,  if $n_0$ is large enough, then  for all $n>n_0$, and $\sqrt{V_n\log\log V_n}<k<\delta V_n$, $\left|\frac{k-\E(S_n)}{V_n}\right|<\delta\ ,\ 
\frac{1\pm(\eps/2)}{2}\left(
\frac{k-\E(S_n)}{\sqrt{V_n}}
\right)^2=\frac{1\pm(\eps/2)}{2}\frac{k^2}{V_n}.
$
It follows that for all $k\in\Z$ such that $\sqrt{V_n\log\log V_n}<k<\delta V_n$, 
\begin{equation}\label{e.k-prob}
 \Prob[S_n=k]=C^{\pm 1}(2\pi V_n)^{-\frac{1}{2}}\exp\left[-\frac{(1\mp(\eps/2))k^2}{2V_n}\right].
\end{equation}

Let $A_n:=\lceil\sqrt{(2+o(1))V_n\log\log V_n}\rceil$ and $B_n:=\lfloor\delta V_n\rfloor$. Then, 
\begin{align*}
&\Prob[A_n\leq S_n\leq B_n]=\sum_{k\in [A_n,B_n]\cap\Z} \Prob[S_n=k]\asymp \frac{1}{\sqrt{V_n}}\sum_{k\in [A_n,B_n]\cap\Z} \exp\left[-\frac{(1\mp\eps)k^2}{2V_n}\right].
\end{align*}
The right-hand-side can be bounded above and below by expressions of the form 
$$
\mathsf I_n:=\frac{1}{\sqrt{V_n}}\int_{A_n'}^{B_n'} e^{-\frac{ct^2}{2V_n}}dt, \text{ where } \begin{array}{l}
A_n'=A_n+O(1)=\sqrt{(2+o(1))V_n\log\log V_n}\\
B_n':=B_n+O(1)\sim \delta V_n\\
c=1-\frac{\eps}{2}\text{ or }c=1+\frac{\eps}{2}. 
\end{array}
$$
After the substitution $t=x\sqrt{2V_n/c}$, the integral $\mathsf I_n$ becomes
\begin{align*}
&\mathsf I_n=\sqrt{2/c}\int_{A_n'\sqrt{c/2V_n}}^{B_n'\sqrt{c/2V_n}} e^{-x^2}dx
=\sqrt{\frac{\pi}{2c}}\left[\mathrm{erfc}(A_n'\sqrt{c/2V_n})-\mathrm{erfc}(B_n'\sqrt{c/2V_n})\right],
\end{align*}
where $\mathrm{erfc}(x):=\frac{2}{\sqrt{\pi}}\int_x^\infty e^{-t^2}dt$. 
It is  well-known that $\mathrm{erfc}(x)\sim \frac{e^{-x^2}}{x\sqrt{\pi}}$ as $x\to\infty$, therefore 
\begin{align*}
\mathsf I_n&\sim\sqrt{\frac{1}{2c}}\left[
\frac{e^{-(c+o(1))\log\log V_n}}{\sqrt{(c+o(1))\log\log V_n}}-\frac{e^{-(\frac{1}{2}+o(1))c\delta^2 V_n }}{\delta\sqrt{cV_n/2} }
\right]\sim \frac{1}{\sqrt{2c^2}} \left(\frac{1}{\log V_n}\right)^{c+o(1)}\frac{1}{\sqrt{\log\log V_n}}.
\end{align*} 
Clearly, $\frac{1}{\sqrt{2c^2}}\frac{1}{\sqrt{\log\log V_n}}=\left(\frac{1}{\log V_n}\right)^{o(1)}$, therefore for all $n$ large enough, $\mathsf I_n=\left(\frac{1}{\log V_n}\right)^{c+o(1)}$. This shows that 
$\DS
\left(\frac{1}{\log V_n}\right)^{1+\eps}<\Prob[A_n\leq S_n\leq B_n]<\left(\frac{1}{\log V_n}\right)^{1-\eps}.
$

Next we consider $\Prob[S_n>B_n]$. By the uniform estimates for $\mathcal F_n(\xi)$ in Theorem \ref{PrFreeEn}, there  are positive constants $c_\delta,R_\delta$ such that $\frac{B_n}{V_n}\in [\mathcal F_n'(c_\delta),\mathcal F_n'(R_\delta)]$ for all $n$. It now follows from  Theorem 7.9 in \cite{Dolgopyat-Sarig-Book}, that 
$$
\Prob[S_n>B_n]\asymp \frac{1}{\sqrt{V_n}}e^{-V_n\mathcal I_n(\frac{B_n}{V_n})}=O(V_n^{-1/2})=o((\log V_n)^{1+\eps}). 
$$

Thus, for all $\eps>0$, for all $n$ large enough, 
$\DS
\left(\frac{1}{\log V_n}\right)^{1+\eps}<\Prob[S_n\geq A_n]<\left(\frac{1}{\log V_n}\right)^{1-\eps}.
$
This is the statement of the  sublemma, but without conditioning. 
Now we consider  the effect of conditioning on $X_n=a_n$. 

\medskip
\noindent
{\bf Upper Bound:} 
 $\Prob[X_n=a_n]=\frac{1}{2}\times$area of parallelogram $a_n\in\mathfrak P_n$. Such parallelograms have bounded geometry, and  their areas are bounded away from zero. So 
\begin{align*}
\Prob[S_n\geq A_n|X_n=a_n]= \frac{\Prob[S_n\geq A_n, X_n=a_n]}{\Prob[X_n=a_n]}\leq \const. \Prob[S_n\geq A_n]\leq \left(\frac{1}{\log V_n}\right)^{1-\eps},
\end{align*}
for all $n$ large enough. 

\medskip
\noindent
{\bf Lower Bound:} \footnote{The argument relying on $\phi$ mixing produces a stronger result than 
what is needed for the lower bound. The stronger bound will be used later 
in the proof of Theorem \ref{ThLT-ZD}
}
By Corollary \ref{Cor-phi-mixing},  $\exists C_{mix}>1, \rho\in (0,1)$  such that 
$$
\|\E(1_{[X_n=a_n]}|X_{n-k})-\Prob[X_n=a_n]\|_\infty<C_{mix}\rho^k.
$$
Recall that the jump functions are uniformly bounded, say by a constant $M$. Then
\begin{align*}
&\Prob[S_n\geq A_n|X_n=a_n]\geq \Prob[S_{n-k}\geq A_n'|X_n=a_n],\text{ where }A_n':=A_n+Mk\\
&=\frac{\E\left(1_{[S_{n-k}>A_n']} 1_{[X_n=a_n]}\right)}{\Prob[X_n=a_n]}
=\frac{\E\left(1_{[S_{n-k}>A_n']} \E(1_{[X_n=a_n]}|X_1,\ldots,X_{n-k})\right)}{\Prob[X_n=a_n]}\\
&=\frac{\E\left[1_{[S_{n-k}>A_n']} \E(1_{[X_n=a_n]}|X_{n-k})\right]}{\Prob[X_n=a_n]},\text{ by the Markov property}\\
&\geq \Prob[S_{n-k}>A_n']\times \frac{\Prob[X_n=a_n]-C_{mix}\rho^k}{\Prob[X_n=a_n]}.
\end{align*}
Since $\Prob[X_n=a_n]$ is uniformly bounded from below, we can choose $k$ so that for $n>k$, 
$$
\Prob[S_n\geq A_n|X_n=a_n]\geq \const \Prob[S_{n-k}>A_n']
$$
{with a constant independent of $n.$}

Since $V_{n}\asymp n$, $V_{n-k}\asymp V_n$, and therefore  $A_{n}'=\sqrt{(2+o(1)V_{n-k}\log\log V_{n-k}}$. So for every $\eps>0$, for all $n$ large enough,  $\Prob[S_{n}>A_n]>\const.\left(\frac{1}{\log V_n}\right)^{1+\eps/2}$. Slightly increasing $\eps$, we can get rid of the constant.
\end{proof}

We are now ready to prove the theorem. 
As in the proofs of Propositions \ref{p.Bikhoff-Sums-Estimate} and \ref{p.Large-Deviations}, $\frac{1}{T}\int_0^T 1_{\St^-}(\phi^t(\tomega))dt$ can be sandwiched between quantities of the form 
$$
\const. \sum_{i=1}^{m} |W_{n,i}|\Prob_{W_{n,i}}\bigl[\wt{\zc}_0 \circ\tpsi_0^{-1}\circ\cdots\circ\tpsi_{n-1}^{-1}-\wt{\zc}_n<z_{n,i}\bigr],
$$ where $m=m(T)$ is uniformly bounded;  $n=n(T)\asymp\log T$; $W_{n,i}$ are pure $E^s_{n}$-segments of constant depth $d$ and height $z_{n,i}+O(1)$, $z_{n,i}=z_{nN_0}(\tomega)+O(1),$ and 
$\DS \sum_{i=1}^m |W_{n,i}|\asymp 1$. 
We saw in the proof of Proposition \ref{p.Large-Deviations} that there are $a_{n-d,i}\in\mathfrak P_{n-d}$ and bounded constants $C_{n,i}$ such that 
$$
\Prob_{W_{n,i}}\bigl[\wt{\zc}_0 \circ\tpsi_0^{-1}\circ\cdots\circ\tpsi_{n-1}^{-1}-\wt{\zc}_n<z_{n,i}\bigr]=\Prob[S_{n-d}>z_{n,i}-C_{n,i}|X_{n-d}=a_{n-d,i}].
$$
It follows that $\frac{1}{T}\int_0^T 1_{\St^-}(\phi^t(\tomega))dt$ can be sandwiched between quantities of the form 
\begin{equation}\label{e.sandwich-bounds}
\const. \sum_{i=1}^{m} |W_{n,i}|\Prob[S_{n-d}>z_{nN_0}(\tomega)+C|X_{n-d}=a_{n-d,i}]
\end{equation}
(To get $C$ independent of $n,i$ we replace $C_{n,i}$ by the maximum or minimum of these bounded constants.)

By Sublemma \ref{l.LIL-for-z},   
$z_{nN_0}(\tomega)<\sqrt{(2+o(1))V_n\log\log V_n}$ for $\mu_0$-a.e. $\tomega$. For such $\tomega$, Sublemma \ref{Lem-LIL-Asymp} and \eqref{e.sandwich-bounds} say that for all $T$ sufficiently large, 
$$
\frac{1}{T}\int_0^T 1_{\St^{-}}(\phi^t(\tomega))dt>\const.\left(\frac{1}{\log V_n}\right)^{1+\frac{\eps}{2}}>\const.\left(\frac{1}{\log \log T}\right)^{1+\frac{\eps}{2}}>\left(\frac{1}{\log \log T}\right)^{1+\eps}.
$$
(The penultimate inequality uses the bounds $V_n\asymp n \asymp\log T$.)

Sublemma \ref{l.LIL-for-z} also says that for a.e. $\tomega$ there exists a sequence $n_j\to\infty$ such that $z_{n_j N_0}(\tomega)>\sqrt{(2+o(1))V_n\log\log V_n}$ (the little-oh is negative). 
Recall that $n(T)$ satisfies that $|\tpsi_{n-1}\circ\cdots\tpsi_0(A_T)|\in [1,L]$. Clearly, there are $T_j\to\infty$ such that $n(T_j)=n_j$. Along this sequence, we have  
$$
\frac{1}{T_j}\int_0^{T_j} 1_{\St^{-}}(\phi^t(\tomega))dt<\const.\left(\frac{1}{\log V_{n_j}}\right)^{1-\frac{\eps}{2}}\!\!\!\!<\const.\left(\frac{1}{\log \log T_j}\right)^{1-\frac{\eps}{2}}\!\!\!\!<\left(\frac{1}{\log \log T}\right)^{1-\eps}.
$$
{The last two displayed equations prove the claim} of 
Theorem \ref{PrBalance-second half} on $\frac{1}{T}\int_0^T 1_{\St^{-}}(\phi^t(\tomega))dt$. 

The assertion on $\frac{1}{T}\int_0^T 1_{\St^{+}}(\phi^t(\tomega))dt$ can be proved in a similar way, by analyzing the asymptotic behavior of 
$\Prob\left(S_n\leq -\sqrt{(2+o(1))V_n\log\log V_n}\, \bigg|\,  X_n=a_n\right)$.
\qed

\begin{proof}[Proof of Theorem \ref{ThLT-ZD}]
Let $\mathcal N$ be the Gaussian distribution with mean zero and standard deviation one. 
The first part of the theorem says that there exists a function $a_T\asymp T/\sqrt{\log T}$ such that for every smooth compactly supported probability density $\vf$ on $\St$, if $\tomega$ is sampled randomly uniformly on $\St$ according to $\vf d\mu_0$, then for each non-negative $h\in L^1(\mu_0)$ such that $\int h d\mu_0=1$, 
$$
\frac{1}{a_T}\int_0^T h(\phi^t(\tomega))dt\xrightarrow[T\to\infty]{\text{dist}}\exp\left(-\frac{1}{2}\mathcal N^2\right).
$$
We begin with a few reductions. 

Recall that $\mu_0$ is ergodic. By the ratio ergodic theorem, $$\DS\int_0^T h(\phi^t(\tomega))dt=[1+o(1)]\int_0^T 1_{\Rect_0}(\phi^t(\tomega))dt, \text{ for a.e. $\tomega\in\St$.}$$
 Thus, if the assertion holds for $1_{\Rect_0}$, then it holds for all $h\in L^1(\mu_0)$ with integral equal to one. Henceforth, we only consider $h=1_{\Rect_0}$.

It is also easy to see, by sandwiching $\vf$ between two non-negative linear combinations of indicators of sets $C_i=D^j (\pi|_{\wt{F}_0})^{-1}(\<\dot{a}_0,\ldots,a_n\>)$, that it is sufficient to prove the theorem for $\vf$ proportional to $1_{C}$ with  $C$ as above. We will analyze  the special case $C=(\pi|_{\wt{F}_0})^{-1}\<\dot{a}_0\>$. The  modifications needed for the general case are routine.
 
Fix $\eps>0$. By Lemma \ref{l.approx-by-pure-segments} and as in the proofs of Propositions \ref{p.Bikhoff-Sums-Estimate} and \ref{p.Large-Deviations}, $\frac{1}{T}\int_0^T 1_{\Rect_0}(\phi^t(\tomega))dt$ can be sandwiched between quantities of the form 
$$
\frac{1}{|A_T^\ast|}\sum_{i=1}^{m} |W_{n,i}|\Prob[S_{n-d}=z_{n,i}|X_{n-d}=a_{n-d,i}],
$$ where $m=m(T)$ is uniformly bounded;  $n=n(T)\asymp\log T$;  $|A_T^\ast|\asymp 1$; and where $W_{n,i}$ are pure $E^s_{n}$-segments of constant depth $d$ and height $z_{n,i}+O(1)$. Moreover, we can choose $W_{n,i}$ so that 
$z_{n,i}=z_{nN_0}(\tomega)+O(1)\text{ and } 
e^{-\eps}|A_T^\ast|<\sum_{i=1}^m |W_{n,i}|<e^{\eps}|A_T^\ast|.
$

Fix some $\eps>0$. 
By Proposition \ref{Prop-Markov-Representation-of-Phi} and Dobrushin's central limit theorem 
(Theorem~\ref{ThMCCLT}), for the Markov chain $\{X_n\}$ conditioned on $X_0=a_0$, there is an $R>0$ so that 
$$
\Prob(E_n)>1-\eps,\text{ where }E_n:=\{\tomega: |z_{nN_0}(\tomega)|<R\sqrt{V_n}\},
$$
and where $\Prob=\Prob(\,\cdot\,|X_0=a_0)=\const.\varphi d\mu_0$.

By the mixing local limit theorem (Theorem \ref{ThLLTArrays}),  for all $n$ large enough, for all $z\in\Z$ such that $|z|\leq R\sqrt{V_n}$, 
$$
\Prob[S_{n-d}=z|X_{n-d}=a_{n-d,i}]=e^{\pm\eps} \left(\frac{1}{\sqrt{2\pi V_n}}e^{-\frac{z^2}{2V_n}}\right).
$$
(To see  uniformity in $z$, assume by contradiction that there is $|z_n|\leq R\sqrt{V_n}$ so that the estimates break for $z_n$. Passing to a subsequence, we may assume that $\frac{z_n}{\sqrt{V_n}}\to z_0$. But then $\Prob[S_{n-d}=z_n|X_{n-d}=a_{n-d,i}]$ violates the mixing local limit theorem.) 

Let $\DS a_T:=\frac{T}{\sqrt{2\pi V_n}}\asymp  \frac{T}{\sqrt{\log T}}.$ Then  for every $\tomega\in E_n$, for all $T$ sufficiently large,
$$
\frac{1}{a_T}\int_0^T 1_{\Rect_0}(\phi^t(\tomega))dt=e^{\pm 3\eps}\exp\left(-\frac{z_{nN_0}(\tomega)^2}{2V_n}\right).
$$
Since  $\Prob(E_n)>1-\eps$, this implies that for all $x$, for all $T$ large enough,\\
$\DS
\Prob\left[\frac{1}{a_T}\int_0^T 1_{\Rect_0}(\phi^t(\tomega))dt<x\right]
$
is sandwiched between expressions of the form 
$$
\Prob\left[\exp\left(-\frac{z_{nN_0}(\tomega)^2}{2V_n}\right)<e^{\pm 3\eps} x\right]\pm 2\eps.
$$
By Proposition \ref{Prop-Markov-Representation-of-Phi} and the CLT for $\{X_n\}$ conditioned on $X_0=a_0$,  
$
\DS \frac{z_{nN_0}(\tomega)}{\sqrt{V_n}}\xrightarrow[n\to\infty]{\text{dist}}\mathcal N$ when $\tomega\sim \vf d\mu_0.
$
Thus $
\Prob\left[\exp\left(-\frac{z_{nN_0}(\tomega)^2}{2V_n}\right)<e^{\pm 3\eps} x\right]\pm 2\eps
$
is eventually sandwiched between expressions of the form 
$\Prob\left[e^{ -N^2/2}<e^{\pm \eps}x\pm\eps\right]\pm 3\eps$, where $N$ is a random variable with distribution $\mathcal N$. Since $\eps$ was arbitrary, 
$\DS\frac{1}{a_T}\int_0^T 1_{\Rect_0}(\phi^t(\tomega))dt\xrightarrow[n\to\infty]{\text{dist}}e^{-\frac{1}{2}\mathcal N^2}$, as required.

\medskip
Let $\mathcal U$ denote the uniform distribution on $[0,1]$. 
The second part of the theorem says that if we sample $\tomega$ randomly  from $\vf d\mu_0$, then 
$
\DS \frac{1}{T}\int_0^T  1_{\St_-}(\phi^t(\tomega))dt\xrightarrow[T\to\infty]{\text{dist}}\mathcal U
$.
Again, it is sufficient to consider the case $\vf=\const.1_C$, where $C:=D^j (\pi|_{\wt{F}_0})^{-1}\<\dot{a}_0,\ldots,a_n\>$, and again we restrict our attention to the special case $j=n=0$. 

As in the first part of the theorem, $\int_0^T 1_{\St_-}(\phi^t(\tomega))dt$ can be sandwiched between expressions of the form 
$$
\frac{1}{|A_T^\ast|}\sum_{i=1}^m |W_{n,i}|\Prob[S_{n-d}>z_{n,i}|X_{n-d}=a_{n-d,i}],
$$
with $m,n,W_{n,i}$ and $z_{n,i}$ as above.  Observe that 
$$
\Prob[S_{n-d}>z_{n,i}|X_{n-d}=a_{n-d,i}]=
\Prob\left[\frac{S_{n-d}-\E(S_{n-d})}{\sqrt{V_{n-d}}}>\frac{z_{n,i}-\E(S_{n-d})}{\sqrt{V_{n-d}}}\bigg|X_{n-d}=a_{n-d,i}\right].
$$

 $V_{n-d}\asymp (n-d)\asymp n\asymp V_n$, and an  
exponential mixing argument (as in Sublemma~\ref{Lem-LIL-Asymp} in the proof of Theorem \ref{PrBalance-second half}) allows to get rid of the conditioning and shows that 
$
\Prob\left[\frac{S_{n-d}-\E(S_{n-d})}{\sqrt{V_{n-d}}}>\frac{z_{n,i}-\E(S_{n-d})}{\sqrt{V_{n-d}}}\bigg|X_{n-d}=a_{n-d,i}\right]
$
can be sandwiched between expressions of the form 
$
e^{o(1)}\Prob\left[\frac{S_{n-d}-\E(S_{n-d})}{\sqrt{V_{n-d}}}>\frac{z_{nN_0}(\tomega)-\E(S_{n-d})}{\sqrt{V_{n-d}}}+o(1)\right]+o(1)
$.

We see that  $\frac{1}{a_T}\int_0^T1_{\St_{-}}(\phi^t(\tomega))dt$ is sandwiched between expressions of the form \\
$
e^{\pm\eps}\Prob\left[\frac{S_{n-d}-\E(S_{n-d})}{\sqrt{V_{n-d}}}>\frac{z_{nN_0}(\tomega)-\E(S_{n-d})}{\sqrt{V_{n-d}}}+o(1)\right]+o(1). 
$
Thus to prove the second part of the theorem, it is sufficient to show that 
for all $0<t<1$, 
\begin{equation}\label{e.aimmm}
\Prob\{\tomega: \Prob\left[\frac{S_{n-d}-\E(S_{n-d})}{\sqrt{V_{n-d}}}>\frac{z_{nN_0}(\tomega)-\E(S_{n-d})}{\sqrt{V_{n-d}}}+o(1)\right]<t\}\xrightarrow[n\to\infty]{}t,
\end{equation}
when $\tomega$ is sampled uniformly from $C$. 

Let 
 $Y_n:=\frac{S_{n-d}-\E(S_{n-d})}{\sqrt{V_{n-d}}}$ and $Z_n:= \frac{z_{nN_0}(\tomega)-\E(S_{n-d})}{\sqrt{V_{n-d}}}+o(1)$. 
We saw above that $Y_n,Z_n$  converge in distribution to $\mathcal N$. Then \eqref{e.aimmm} follows from the following general fact:

\medskip
\noindent
{\em Claim.\/} Suppose $Y_n$ and $Z_n$ are random variables (perhaps on different probability spaces), with associated probability measures $\Prob_{Y_n}$ and $\Prob_{Z_n}$. If $Y_n,Z_n\xrightarrow[n\to\infty]{\text{dist}}\mathcal N$, then $\Prob_{Y_n}[Y_n>Z_n]\xrightarrow[n\to\infty]{dist}\mathcal U$.  

\medskip
\noindent
{\em Proof of the Claim.\/}
Fix $\eps>0$ and 
let $\Phi(x):=\frac{1}{\sqrt{2\pi}}\int_{-\infty}^x e^{-t^2/2}dt$, then 
$$
\Prob_{Y_n}[Y_n>z]=1-\Phi(z)+o(1)\text{  uniformly on compact sets of $z$. }
$$
Since $Z_n\xrightarrow[n\to\infty]{dist}\mathcal N$, there is $R>0$ such that $\Prob_{Z_n}[|Z_n|<R]>1-\eps$ for all $n$. 
For all $0<t<1$:
\begin{align*}
&\Prob_{Z_n}\left[\Prob_{Y_n}[Y_n>Z_n]<t\right]=
\Prob_{Z_n}\left[|Z_n|<R\ , \Prob_{Y_n}[Y_n>Z_n]<t\right]+O(\eps)\\
&=\Prob_{Z_n}[|Z_n|<R\ , \ 1-\Phi(Z_n)+o(1)<t]{ +O(\eps)}
=\Prob_{Z_n}[\Phi(Z_n)>1-t+o(1)]{+O(\eps)}\\
&=
\Prob_{Z_n}[Z_n>\Phi^{-1}(1-t+o(1))]{+O(\eps)}\\
&\overset{!}{=}[1-\Phi(\Phi^{-1}(1-t+o(1))]+O(\eps)=[t+o(1)]+O(\eps),\text{ as $n\to\infty$}, 
\end{align*} 
(identity $\overset{!}{=}$ uses  $Z_n\xrightarrow[n\to\infty]{dist}\cN$).
Taking $\eps\to 0$, we obtain  $\Prob_{Y_n}[Y_n>Z_n]\xrightarrow[n\to\infty]{dist}\mathcal U$. 
\end{proof}

\begin{proof}[Proof of Theorem \ref{ThLT-NZD}]
Let $\xi>0$. The theorem states that there are $a_T(\xi)\asymp \sqrt{\log T}$ 
and  $b_T$ satisfying
$(1-c)\log T\leq b_T(\xi)\leq \log T$ with $0<c<1$ as follows:
Suppose $h\in L^1(\mu_\xi)$ is non-negative and let $\vf:\St\to\R$ be a  smooth,  compactly supported,  non-negative function such that  $\int\vf d\mu_\xi=1$. If  $\tomega$ is sampled randomly from the probability measure $\vf d\mu_\xi$,    then 
$$
\frac{\log \int_0^T h(\phi^t(\tomega))dt-b_T(\xi)}{a_T(\xi)}\xrightarrow[T\to\infty]{\text{dist}}\mathcal N. 
$$
This includes  the case $h=1_{\St_{-}}$, because 
$\DS
\|1_{\St_{-}}\|_{L^1(\mu_\xi)}
=\sum_{k=1}^\infty e^{-\xi k}<\infty
$
{ (recall that $\mu_\xi\circ D=e^{\xi}\mu_\xi$, and $\mu_\xi(\Rect_0)=1$,
see Corollary \ref{c.HHW}). }

As in the proof of Theorem \ref{ThLT-ZD}, it is sufficient to consider $h=1_{\Rect_0}$ and $\varphi$ proportional to the indicator of a set of the form $D^j (\pi|_{\wt{F}_0})^{-1}\<\dot{a}_0,\ldots,a_n\>$. We consider the case $j=n=0$ and leave the (routine) modification needed in the general case to the reader.

{It is sufficient to consider $\tomega$ which are $\mu_\xi$-generic, because  $\mu_\xi$-a.e. $\tomega$ are like that. By Theorem \ref{t.Char-mu-xi-Gen}, $\xi_n(\tomega)\xrightarrow[n\to\infty]{}\xi$ $\mu_\xi$-almost everywhere. 
By Proposition \ref{p.Bikhoff-Sums-Estimate}, there are  integers $n=n(T)\asymp \log T$ such that for $\mu_\xi$-a.e. $\tomega$,
\begin{equation}\label{e.pre-Taylor}
\frac{\log \int_0^T 1_{\Rect_0}(\phi^t(\tomega))dt-\log T}{\sqrt{V_n}}=-\sqrt{V_n} \; 
\HH_n\left(\frac{z_{nN_0} ({\tomega})}{V_n}\right)
+O\left(\frac{\log n}{\sqrt{n}}\right).
\end{equation}
The error estimate uses the bound $V_n\asymp n$ from Proposition \ref{Prop-Var-Growth}.

Let  $\DS z_n^\xi:=\wt{\E}^{\xi}(\wt{S}^\xi_n),$ the expectation of $S_n$ with respect to the change of measure $\Prob^{\xi}$. By \eqref{Pres1stDer}, $z_n^\xi=\PP_n'(\xi)+O(1)$, and by \eqref{Pres2dDer} and the estimate $V_n\asymp n$, there exists $R>0$ such that $\frac{z_n^\xi}{V_n}\asymp 1$, and that $\frac{z_n^\xi}{V_n}$ belong to the intervals $[\bar{a}^R_n,\bar{b}^R_n]$ from Lemma \ref{l.rate-function}. It follows that 
$\DS
\sqrt{V_n}\HH_n\left(\frac{z_n^\xi}{V_n}\right)\asymp\sqrt{V_n}\left(\frac{z_n^\xi-\E(S_n)}{{V_n}}\right)^2\asymp \sqrt{V_n}\asymp\sqrt{\log T}.
$

In addition, since  $z_n^\xi=\PP_n'(\xi)+O(1)$ and  $\HH_n'$ is uniformly bounded on $[\bar{a}^R_n,\bar{b}^R_n]$, 
\begin{align*}
&\HH_n'\left(\frac{z_n^\xi}{V_n}\right)=\HH_n'\left(\frac{\PP_n'(\xi)}{V_n}+o(1)\right)=\HH_n'\left(\frac{\PP_n'(\xi)}{V_n}\right)+o(1)=\xi+o(1), 
\end{align*}
where the last identity is because  $\HH_n$ is the Legendre transform of $\PP_n/V_n$.

Finally, by Lemma \ref{l.rate-function}, $\sup |\HH_n''|=O(1)$ on $[\bar{a}^R_n,\bar{b}^R_n]$, and $\HH_n$ is real-analytic on a neighborhood of $[\bar{a}^R_n,\bar{b}^R_n]$.
Expanding the right-hand-side of \eqref{e.pre-Taylor} into the Taylor expansion, and rearranging terms, we obtain 
\begin{align*}
&\frac{\log \int_0^T 1_{\Rect_0}(\phi^t(\tomega))dt-\log T}{\sqrt{V_n}}+\sqrt{V_n}\HH_n\left(\frac{z_n^\xi}{V_n}\right)=\\
&
=
-\HH_n'\left(\frac{z_n^\xi}{V_n}\right)\frac{z_{nN_0}({ \tomega})-z_n^\xi}{\sqrt{V_n}}+
O\left(\frac{(z_n^\xi-z_{nN_0}({ \tomega}))^2}{\sqrt{V_n^3}} \right)+O\left(\frac{\log\log T}{\sqrt{\log T}}\right)\\
&=-(\xi+o(1))\sqrt{\frac{\wt{V}^\xi_n}{V_n}}\left[\frac{z_{nN_0}(\tomega)-z_n^\xi}{\sqrt{\wt{V}^\xi_n}}+
O\left(\frac{\wt{V}^\xi_n}{V_n^{3/2}}\left(\frac{z_n^\xi-z_{nN_0}(\tomega)}{
\sqrt{\wt{V}^\xi_n}} \right)^2\right)\right]+O\left(\frac{\log\log T}{\sqrt{\log T}}\right),
\end{align*}
where $\wt{V}^\xi_n:=\Var(\wt{S}^\xi_n)$. 
By Lemma \ref{l.V-n-xi}, $\wt{V}^\xi_n\asymp V_n$, so $\wt{V}^\xi_n\asymp n\asymp\log T$

By Lemma \ref{l.uniform-Doeblin}, the Markov chain $\{\wt{X}_n^\xi\}$ is uniformly elliptic, and since $z_{nN_0}(\tomega)=\wt{S}^\xi_n+O(1)$ in distribution (Lemma \ref{l.Phi-Decomp}), 
{ $\frac{z_{nN_0}(\tomega)-z_n^\xi}{\sqrt{\wt{V}^\xi_n}}$ converges in distribution to $\mathcal N$.
Therefore 
$\frac{\wt{V}^\xi_n}{V_n^{3/2}}\left(\frac{z_n^\xi-z_{nN_0}(\tomega)}{
\sqrt{\wt{V}^\xi_n}} \right)^2$ converges to $0$ in probability. Hence}
the term in the square brackets converges in distribution to $\mathcal N$. Letting 
\begin{align*}
& a_T(\xi):={\sqrt{V_n}/\xi}, \text{ and } 
b_T(\xi):=\log T-V_{n}\HH_{n}\left(\frac{z_n^\xi}{V_n}\right)
\end{align*}
we obtain  $\DS \frac{\log \int_0^T 1_{\Rect_0}(\phi^t(\tomega))dt-b_T(\xi)}{a_T(\xi)}\xrightarrow[T\to\infty]{\text{dist}}\mathcal N$.

The denominator satisfies the bounds 
$ a_T(\xi)=\sqrt{V_n}/\xi\asymp \sqrt{n}/\xi\asymp \sqrt{\log T}$. 
In particular, $\DS \frac{\log \int_0^T 1_{\Rect_0}(\phi^t(\tomega))dt-b_T(\xi)}{\log T}$ converges in probability to zero. 

Recall that  $\DS \mu_\xi[\zc\leq 1]=\sum_{k\geq -1} e^{-k\xi}<\infty$, therefore $1_{[\zc\leq 1]}\in L^1(\mu_\xi)$, and by the   ratio ergodic theorem
$\DS
\log \int_0^T 1_{\Rect_0}(\phi^t(\tomega))dt=\bbI_T(\tomega)+\const+o(1)\ \ \text{$\mu_\xi$-a.e.}
$
with $\bbI_T(\tomega)$ as in Theorem \ref{CrLogOTDetXi}. It follows that 
$\DS \frac{\bbI_T(\tomega)-b_T(\xi)}{\log T}$ converges in probability to zero. However, by Theorem \ref{CrLogOTDetXi}, $\DS \frac{\bbI_T(\tomega)-\bbI_T^\xi}{\log T}$ converges in probability (indeed a.e.) to zero. Necessarily, $b_T(\xi)=\bbI_T^\xi+o(\log T)$. It now follows from Theorem \ref{CrLogOTDetXi} that for some $0<c(\xi)<1$, for all $T$ large enough, $(1-c(\xi))\log T\leq b_T(\xi)\leq \log T$. }
\end{proof}

\section{Exponentially Tilted $D_N(\alpha)$ }
\label{SSZeroTilted}
In this section we prove Theorem \ref{t.exp-tilted-exist} on the existence, for each $\eta\in\R$, of irrationals $\alpha$ of bounded type for which  $D_N(\alpha)$ is exponentially tilted with tilt parameter $\eta$. To do this we first prove Theorem \ref{CrWeakBal}, 
 which claims
that for each $\eta\in\R$ there exists $\beta\in(-1,1)$ of
bounded type, such that $\tomega_0$, the midpoint of the top side of $\Rect_0$, is $\mu_\eta$-generic for the linear flow  with velocity vector ${\beta\choose 1}$. 

We will look for  $\beta$ inside the set $
\Lambda_6:=\{\beta\in (-1,1)\setminus \Q: M(\beta)\leq 6\}, \text{ see \eqref{e.M-beta}}
$ (why $6$ will be explained later). 
We view $\beta$  as a parameter. For $\beta\in\Lambda_6$, we may take the constant $N_0$ in \eqref{e.mat-prod-form} to be independent of $\beta$. But the other objects in the construction  depend on $\beta$, and we will signify this  by including $\beta$ in the subscript. For example,
\begin{itemize}
\item $\phi_\beta:=$ the linear flow with velocity vector ${\beta\choose 1}$;
\item $\tpsi_{n,\beta}:=$  the renormalizing automorphisms from Proposition \ref{Prop-Contracting-Automorphisms}; 
\item $\PP_{n,\beta}(\xi)=p_{1,\beta}(\xi)+\cdots+p_{n,\beta}(\xi)$,  the pressure functions from Section \ref{ScPressure};
\item $z_{n,\beta}:=$ $\Z$-coordinate of the horizontal rectangle which contains the beginning of the ray $(\tpsi_{n-1,\beta}\circ\cdots\circ\tpsi_{0,\beta})(\{\phi_\beta^t(\tomega_0):t>0\})$; 
\item $\eta_{n,\beta}:=$ solution to the equation $\PP_{n,\beta}'(\eta_{n,\beta})=z_{n,\beta}$ (if it exists). 
\end{itemize}
We will construct $\beta\in\Lambda_6$ such that $\eta_{n,\beta}\to\eta$. It will then follow from 
Proposition~\ref{p.Pressure-and-Geometric-Pressure}, and Theorem \ref{t.Char-mu-xi-Gen} that $\tomega_0$ is $\mu_\eta$-generic, whence by \eqref{e.Lin-Flow-1},   $D_N(\frac{\beta+1}{2})$ is exponentially tilted with tilt parameter $\eta$.

\subsection{Preparations} 
For every irrational $\beta\in (-1,1)$, let $m_i(\beta)$ be the unique non-zero integers such that 
$
\beta=[0;2m_1(\beta),2m_2(\beta),\ldots].
$

\begin{lem}
\label{LmBetaCont}
There exists constants $\sK(R)>1$ such that 
if $\beta',\beta''\in\Lambda_6$ satisfy $m_i(\beta')=m_{i+i_0N_0}(\beta'')$ for all  
$2n_1 N_0\leq i \leq  2 n_2 N_0$,   then
$$
\left|\wt{\E}^\xi_{\beta'}\left(\sum_{n=n_1}^{n_2-1} f_{n,\beta'}(\wt{X}^\xi_{n,\beta'})\right)-\wt{\E}^\xi_{\beta''}\left(\sum_{n=n_1+i_0}^{n_2+i_0-1} f_{n,\beta''}(\wt{X}^\xi_{n,\beta''})\right)\right|\leq \sK\text{ for all }\xi\in [-R,R].
$$
\end{lem}

\begin{proof}
Let $m_i$ be the common value of $m_i(\beta')$ and $m_{i+i_0 N_0}(\beta'')$, for $2n_1 N_0\leq i\leq 2n_2 N_0$. 
By the proof of Proposition \ref{Prop-Contracting-Automorphisms}, for every $n_1\leq n\leq n_2-1$, 
$
\tpsi_{n,\beta'}=\tpsi_{n+i_0,\beta''}=
$ the unique homogeneous automorphism with zero drift and derivative 
\begin{align*}
&A_{N_0}([0;2m_{2nN_0+1}, 2m_{2nN_0+2},\ldots])
=\\
&\hspace{1cm}{\tiny{\left(
       \begin{array}{cc}
         1 & -2m_{2(n+1)N_0} \\
         0 & 1 \\
       \end{array}
     \right)
     \left(
       \begin{array}{cc}
         1 & 0 \\
         -2m_{2(n+1)N_0-1} & 1 \\
       \end{array}
     \right)
     \cdots
      \left(
       \begin{array}{cc}
         1 & -2m_{2nN_0+2} \\
         0 & 1 \\
       \end{array}
     \right)
     \left(
       \begin{array}{cc}
         1 & 0 \\
         -2m_{2nN_0+1} & 1 \\
       \end{array}
     \right)}}.
\end{align*}
Equations   \eqref{e.E-k} show that  the slopes of $E^t_{n,\beta'},E^u_{n+i_0,\beta''}$ $(t=u,s)$ are numbers whose even continued fraction expansion agrees on the first $2N_0 d(n)$ coordinates, where 
$$
 d(n):=
 \min(n-n_1, n_2-n).
$$
Since $\beta',\beta''\in\Lambda_6$, there is a global constant $\theta_1\in (0,1)$ such that 
$$
|\measuredangle(E^t_{n,\beta'},E^t_{n+i_0,\beta''})|\leq C\theta_1^{d(n)} \text{ for all $n_1\leq n\leq n_2$ and  
$t=u,s$}$$
(see Lemmas \ref{l.v-contracts} and \ref{l.limit-cfe-exists}). 

Let $\Pi_{n,\beta}:=\Pi(E^s_{n,\beta},E^u_{n,\beta})$ be the polygons from Figure  \ref{Figure-Pi(n)} on page \pageref{Figure-Pi(n)}. The angles and sides of    $\Pi_{n,\beta'}$ and $\Pi_{n+i_0,\beta''}$ are the same up to errors of size $O(\theta_1^{d(n)})$, $n_1\leq n\leq n_2$.  

Recall that  $\Pi_{n,\beta}$ is partitioned into two parallelograms $Q_{1,\beta}(n), Q_{2,\beta}(n)$ as Figure \ref{Figure-Pi(n)}. 
Since $\beta',\beta''\in\Lambda_6$, the lengths of the $u$-sides of $Q_{i,\beta}(n)$ are bounded below. Necessarily, for some global constant $\sk\in\N$ and  for all $n_1+\sk\leq n\leq n_2-\sk$, $\tpsi_{n,\beta'}(Q_{i,\beta'}(n))$ and $\tpsi_{n+i_0,\beta''}(Q_{i,\beta''}(n+i_0))$ wrap around $Q_{j,\beta'}(n+1)$ and $Q_{j,\beta''}(n+i_0+1)$ the same number of times. Thus the Markov partitions $\mathfrak P_{n,\beta'}, \mathfrak P_{n+i_0,\beta''}$ can be enumerated by
\begin{equation}
\begin{aligned}\label{e.P-enumerate}
&Q_{\beta',1,1}^n, \ldots,Q_{\beta',1,N_{1}(n)}^n;
Q_{\beta',2,1}^n, \ldots,Q_{\beta',2,N_{2}(n)}^n \\
&Q_{\beta'',1,1}^{n+i_0}, \ldots,Q_{\beta'',1,N_{1}(n)}^{n+i_0};
Q_{\beta'',2,1}^{n+i_0}, \ldots,Q_{\beta'',2,N_{2}(n)}^{n+i_0} \ \ \end{aligned}
\end{equation}
with the same $N_i(n)$. 
Moreover, by making $\sk$ sufficiently large we can  guarantee that the incidence matrices  $(t^{(n)}_{\beta',a,b}), (t^{(n+i_0)}_{\beta'',a,b})$ in \eqref{e.incidence-matrix} are the same, i.e.
$$
\psi_{n,\beta'}(\mathrm{int}(Q_{\beta',i,j}^n))\cap \mathrm{int}(Q_{\beta',k,\ell}^{n+1})\neq \emptyset\Leftrightarrow \psi_{n+i_0,\beta''}(\mathrm{int}(Q_{\beta'',i,j}^{n+i_0}))\cap \mathrm{int}(Q_{\beta'',k,\ell}^{n+i_0+1})\neq \emptyset.
$$
To see this recall that $\psi_{n,\beta'}=\psi_{n+i_0,\beta''}$. Next, note that the infimum of the stable widths of $\psi_n(\mathrm{int}(Q_{\beta,i,j}^n))$ over all $\beta\in\Lambda_6$, and $i,j,n$  is a positive number $\delta$ (bounded geometry). Therefore,  
if $\psi_n(\mathrm{int}(Q_{\beta,i,j}^n))\cap \mathrm{int}(Q_{\beta,k,\ell}^{n+1})\neq \emptyset$ then this intersection contains $u$-fibres within distance at least $\delta/2$ from the boundary of $Q_{\beta,k,\ell}^n$. By passing from $(n,\beta')$ to $(n+i_0,\beta'')$ we are perturbing $\psi_{n,\beta'}(\mathrm{int}(Q_{\beta,i,j}^n))$ and $\mathrm{int}(Q_{\beta',k,\ell}^{n+1})$ by a an amount uniformly bounded by $O(\theta_1^{d(n)})=O(\theta_1^{\mathsf k})$ so the non-emptiness of the intersection is preserved. 

For similar reasons, if $\sk$ is sufficiently large and $n_1+\sk\leq n\leq n_2-\sk$, then 
$$
f_{n,\beta'}=f_{n+i_0,\beta''}.
$$

Finally, we claim that if $\sk$ is sufficiently large and $n_1+\sk\leq n\leq n_2-\sk$, then the transition kernels from Proposition \ref{Prop-Markov-Kernel} satisfy
\begin{equation}
\label{TrPrClose}
\pi_{\beta',n,n+1}-\pi_{\beta'',n+i_0,n+i_0+1}=O\left(\theta_1^{d(n)}\right) \text{ for all }n_1+\sk\leq n\leq n_2-\sk.
\end{equation}
Indeed, $\DS\pi_{\beta',n,n+1}(Q_{\beta',i,j}^n,Q_{\beta',k,\ell}^{n+1})=\frac{\mu_0[\psi_{n,\beta'}(Q_{\beta',i,j}^n)\cap Q_{\beta',k,\ell}^{n+1}]}{\mu_0[Q_{\beta',k,\ell}^{n+1}]}$,  and  $$\DS\pi_{\beta'',n+i_0,n+i_0+1}(Q_{\beta'',i,j}^{n+i_0},Q_{\beta'',k,\ell}^{n+i_0+1})=\frac{\mu_0[\psi_{n,\beta''}(Q_{\beta'',i,j}^{n+i_0})\cap Q_{\beta'',k,\ell}^{n+i_0+1}]}{\mu_0[Q_{\beta'',k,\ell}^{n+i_0+1}]}.$$ By the previous discussion, the numerators and the denominators  are  equal  up to error $O(\theta_1^{d(n)})$, and by Lemma \ref{Lemma-Bounded-Geometry-Gamma(M)}, the denominator is bounded away from zero for $\beta\in\Lambda_6$.  

Now we consider the change of measure in \S \ref{ss.change-of-measure}, starting with the functions $h_{n,\beta}(\cdot,\xi)$ from Lemma \ref{LmAddPres}. These functions are constructed in the proof of Lemma~7.3 in \cite{Dolgopyat-Sarig-Book}, as follows. Using the enumeration \eqref{e.P-enumerate}, we identify $\mathfrak P_{n,\beta}$ with a set $\mathfrak S_n$ independent of $\beta$, and $\pi_{\beta,n,n+1}$ with Markov kernels on $\mathfrak S_n\times\mathfrak S_{n+1}$. Let 
\begin{itemize}
\item $V_{n}:=\ell^\infty(\mathfrak S_n)$;
\item $\cC_{n}:=\{h\in V_{n}: h\geq 0\}$;
\item $L_{n,\beta}^\xi:V_{n+1}\to V_{n}$, $(L_{n,\beta}^\xi h)(a)=\sum_{b\in\mathfrak S_{n+1}}e^{\xi f_{n+1}(b)}\pi_{\beta,n,n+1}(a,b)h(b)$. 
\end{itemize}
It is shown in the proof of Lemma~7.3 in \cite{Dolgopyat-Sarig-Book} that the limits 
$\DS \lim_{\ell\to\infty} L_{n+1,\beta}^\xi\cdots L_{n+\ell,\beta}^\xi 1$ exist, and that $h_n(\cdot,\xi)$ can be taken to be of the form $\DS e^{c_n\xi}\times \lim_{\ell\to\infty} L_{n+1,\beta}^\xi\cdots L_{n+\ell,\beta}^\xi 1$ for suitable $c_n$. The choice of $c_n$ does not affect $\pi_{\beta,n,n+1}^\xi$.

The existence of the limit is proved in \cite{Dolgopyat-Sarig-Book} by first showing that $L_{n,\beta}^\xi$ maps the closure of $\cC_{n+1}$ into the interior of $\cC_n$, and that $L_{n,
\beta}^\xi L_{n+1,\beta}^\xi: \cC_{n+2}\to \cC_n$ are uniform contractions in the Hilbert-Birkhoff metric. 
The contraction rate $\theta_2\in (0,1)$ is uniform in $\beta$ and $n$, provided  $\beta\in \Lambda_6$. Thus,
$$
h_{n,\beta}^\xi=C_n[(L_{n+2,\beta}^\xi \cdots L_{n_2-\sk,\beta}^\xi 1+O(\theta_2^{d(n)})],
$$
where $C_n$ are positive proportionality constants. Comparing these expressions for $\beta=\beta',\beta''$ using \eqref{TrPrClose}, we conclude that for some global constants $\theta_3\in (0,1)$, $\sK'(R)$,  for all $|\xi|\leq R$, for all $n_1+\sk\leq n\leq n_2-\sk$, 
$
\|h_{\beta',n}^\xi-h_{\beta'',n+i_0}^\xi\|_\infty=\sK'(R)\theta_3^{d(n)}
$, 
and so
$ 
\|\pi_{\beta',n,n+1}^\xi-\pi_{\beta'',n+i_0,n+i_0+1}^\xi\|\leq O\left(\theta_3^{d(n)}\right) $
uniformly for $|\xi|\leq R.$ 

Using the identity  $f_{n,\beta'}=f_{n+i_0,\beta''}$, we see by direct calculations that for some constant $\sK(R)$, if $|\xi|\leq R$, then (after identifying the state spaces $\mathfrak P_{n,\beta'}$ and $\mathfrak P_{n+i_0,\beta''}$ as in \eqref{e.P-enumerate})
$$
\left|
\wt{\E}^{\xi}_{\beta'}\left(\sum_{n=n_1+\sk}^{n_2-\sk}f_{n,\beta'}(\wt{X}^{\xi}_{n,\beta'})\bigg| \wt{X}^{\xi}_{n_1+\sk,\beta'}\right)-
\wt{\E}^{\xi}_{\beta''}\left(\sum_{n=n_1+i_0+\sk}^{n_2+i_0-\sk}f_{n,\beta''}(\wt{X}^{\xi}_{n,\beta''})\bigg| \wt{X}^{\xi}_{n_1+i_0+\sk,\beta''}\right)
\right|\leq \sK(R).
$$
All Markov chains associated to $\beta',\beta''\in\Lambda_6$ are uniformly elliptic with the same ellipticity constant $\eps_0$ (Proposition \ref{Prop-Doeblin}).  Therefore by Lemma \ref{l.uniform-Doeblin}, the Markov chains 
 $\wt{X}^{\xi}_{n,\beta'}, \wt{X}^{\xi}_{n,\beta''}$ are uniformly elliptic with the same ellipticity constant $\eps_0$ for all $\xi\in [-R,R]$ and $\beta',\beta''\in\Lambda_6$. It now follows from Corollary \ref{Cor-phi-mixing} that for some constant $\sK'=\sK'(R)$, for all $\beta',\beta''$ as above, 
$$
\left|
\wt{\E}^{\xi}_{\beta'}\left(\sum_{n=n_1+\sk}^{n_2-\sk}f_{n,\beta'}(\wt{X}^{\xi}_{n,\beta'})\right)-
\wt{\E}^{\xi}_{\beta''}\left(\sum_{n=n_1+i_0+\sk}^{n_2+i_0-\sk}f_{n,\beta''}
(\wt{X}^{\xi}_{n,\beta''})\right)
\right|\leq \sK(R).
$$
By Lemma \ref{Lemma-f(X)} there is a constant $\sK''$ such that for all $\beta\in\Lambda_6$, $\sup_{i\geq 1}|f_{i,\beta}|\leq \sK''$. So
$$
\left|
\wt{\E}^{\xi}_{\beta'}\left(\sum_{n=n_1}^{n_2}f_{n,\beta'}(\wt{X}^{\xi}_{n,\beta'})\right)-
\wt{\E}^{\xi}_{\beta''}\left(\sum_{n=n_1+i_0}^{n_2+i_0}f_{n,\beta''}(\wt{X}^{\xi}_{n,\beta''})\right)
\right|\leq \sK(R)+4\sK'',
$$
and the lemma follows.
\end{proof}

Henceforth we set \index{Even continued fraction expansion!examples}
$$
\beta_2:=2\sqrt{2}-3=[0;-6,6,-6,6,\ldots]\ \ \ ,\ \ \  \beta_3:=2\sqrt{3}-3=[0;2,6,2,6,\ldots].
$$
Note that $\beta_2,\beta_3\in\Lambda_6$ (this is why we chose six). 
\begin{lem}
\label{LmPres23} \ 
\begin{enumerate}[(1)]
\item $\exists\eta_{n,\beta_3}\xrightarrow[n\to\infty]{}0$ such that 
$\PP_{n,{\beta_3}}'(\eta_{n,{\beta_3}})=z_{n,\beta_3}$.
\item For every $\xi\in\R$, for all $n$ large enough, $\PP_{n,\beta_2}'(\xi)<z_{n,\beta_2}$.
\end{enumerate}
\end{lem}

\begin{proof}
Let  $\beta=\beta_3$. By Corollary \ref{CrOrigin},   $\tomega_0$ is $\mu_0$-generic. By Theorem \ref{t.Char-mu-xi-Gen}, $\xi_{n,\beta_3}(\tomega)\to 0$. 
By Proposition \ref{p.Pressure-and-Geometric-Pressure}, $\eta_{n,\beta_3}$ eventually exists, and  $\eta_{n,\beta_3}\to 0$.

Now consider   $\beta=\beta_2$.  By Theorem \ref{ThDriftOfOrigin}   
$
z_{n,\beta_2}=nN_0
$. Suppose first that for all $n$ large enough, there are no $\xi_n$ such that $\PP_{n,\beta_2}'(\xi_n)\geq z_{n,\beta_2}$, then part (2) of the lemma holds automatically. Otherwise, there are $n_k\uparrow\infty$ and $\xi_{n_k}\in\R$ such that $\PP_{n_k,\beta_2}'(\xi_{n_k})\geq z_{n_k,\beta_2}$.
 Since $\PP_{n_k,\beta_2}'(0)=\E(S_{n,\beta_2})=O(1)$ and $\PP_{n_k,\beta_2}'(\xi_{n_k})\geq z_{n_k,\beta_2}\to\infty$, there must exist $\eta_{n_k}\in\R$ such that $\PP_{n_k,\beta_2}'(\eta_{n_k})=z_{n_k,\beta_2}$. We will show that $\eta_{n_k}\to\infty$, and deduce part~(2) from the strict convexity of  $\PP_{n_k,\beta_2}$ for large $k$ (Lemma \ref{LmPressure}).
 
Write ${\beta_2}=2\alpha_2-1$, where 
$
\alpha_2:=\sqrt{2}-1$. Since $\alpha_2=\sqrt{2}\mod 1
$,   $D_n(\alpha_2)=D_n(\sqrt{2})$. It is shown in  \cite{Boshernitzan-Ralston} that $D_n(\sqrt{2})$ is uniformly bounded below. By \eqref{e.Lin-Flow-1}, $\zc[\phi_{\beta_2}^t(\tomega_0)]$ is  bounded below, and  none of the $\mu_\xi$ can be  historical for $\tomega_0$. 
By Theorem \ref{CrHist} and  Proposition \ref{p.Pressure-and-Geometric-Pressure}, $\eta_{n_k}$ cannot have finite limit points. It follows that $|\eta_{n_k}|\xrightarrow[n\to\infty]{}\infty.$
Since 
 $\PP_{n_k,\beta_2}'(\eta_{n_k})=z_{n_k,\beta_2}\to+\infty$ and 
 $\PP_{n_k,\beta_2}'(0)=O(1)$,  $\eta_{n_k}>0$, whence $\eta_{n_k}\to\infty.$
%
%
%
%
%
\end{proof}

\begin{lem}\label{l.14.3} For all $\beta\in\Lambda_6$
there are  constants $c_\beta(R),\mathsf C_\beta(R)>0$ such that for all $\xi\in [-R,R]$, for all $L$ sufficiently large, for all $n$
$$
\left|\PP_{L,\beta}'(\xi)-\wt{\E}^{\xi}_{\beta}\left(\sum_{k=1}^{L}f_{k,\beta}(\wt{X}^{\xi}_{k,\beta})\right)\right|\leq \mathsf C_\beta(R)\ \text{ and }\ \PP_{L,\beta}''(\xi)\geq c_\beta(R)L.
$$ 
\end{lem}
\noindent
This is a direct consequence of Lemma \ref{LmPressure} and Proposition \ref{Prop-Var-Growth}.  The methods of \cite{Dolgopyat-Sarig-Book} do not guarantee that  the constants can be chosen uniformly in $\Lambda_6$.

\begin{lem}\label{l.Dima-Trick}
For every $C,r,R>0$,  for all  $L\in\N$ large enough and  $\beta\in \Lambda_6$, 
\begin{enumerate}[(1)]
\item  If $(2m_{2nN_0}(\beta),\ldots,2m_{2(n+L)N_0}(\beta))=(-6,6,\ldots,-6,6)$ and $\xi\leq R$, then 
\begin{equation}\label{e.Dima-Trick-1}
\wt{\E}^{\xi}_{\beta}\left(\sum_{k=n+1}^{n+L}f_{k,\beta}(\wt{X}^{\xi}_{k,\beta})\right)< z_{n+L,\beta}-z_{n,\beta}-C.
\end{equation}
\item  If $(2m_{2nN_0}(\beta),\ldots,2m_{2(n+L)N_0}(\beta))=(2,6,\ldots,2,6)$ and $\xi>r$, then 
\begin{equation}\label{e.Dima-Trick-2}
\wt{\E}^{\xi}_{\beta}\left(\sum_{k=n+1}^{n+L}f_{k,\beta}(\wt{X}^{\xi}_{k,\beta})\right)> z_{n+L,\beta}-z_{n,\beta}+C.
\end{equation}
\end{enumerate}
\end{lem}
\begin{proof}
Choose constants $\sK:=\sK(R+1)$, $c_{\beta}=c_{\beta}(R+1)$ and $\mathsf C_\beta(R+1)$ as in the previous lemmas. By Lemma \ref{LmPres23} there exists  $L_1$ so that if $\ell>L_1$, then 
$$\PP_{\ell,\beta_2}'(R+1)<z_{\ell,\beta_2},\text{ and }\PP_{\ell,\beta_3}'(\eta_{\ell,\beta_3})=0\text{ for some }\eta_{\ell,\beta_3}\in (-r/2,r/2).$$ 
If $(2m_{2nN_0}(\beta),\ldots,2m_{2(n+L)N_0}(\beta))=(-6,6,\ldots,-6,6)$ and $\xi\leq R$, then for all $L$
\begin{align*}
&\wt{\E}^{\xi}_{\beta}\left(\sum_{k=n+1}^{n+L}f_{k,\beta}(\wt{X}^{\xi}_{k,\beta})\right)\leq \sK+
\wt{\E}^{\xi}_{\beta_2}\left(\sum_{k=1}^{L}f_{k,\beta_2}(\wt{X}^{\xi}_{k,\beta_2})\right),\text{ by Lemma \ref{LmBetaCont}}\\
&\leq \sK+2\mathsf C_{\beta_2}(R)+
\PP_{L,\beta_2}'(\xi),\text{ by Lemma \ref{l.14.3}}\\
&\leq \sK+2\mathsf C_{\beta_2}(R)+\PP_{L,\beta_2}'(R)
\leq \sK+2\mathsf C_{\beta_2}(R)+\PP_{L,\beta_2}'(R+1)-c_{\beta_2}(R+1)L\\
&< \PP_{L,\beta_2}'(R+1)-C, \text{ provided }L>(C+\sK+2\mathsf C_{\beta_2}(R))/c_{\beta_2}(R+1)\\
&\leq z_{L,\beta_2}-C,\text{ provided $L>L_1$} \\
&=\frac{1}{2}\sum_{j=1}^{LN_0}[\sgn m_{2j}(\beta_2)-\sgn m_{2j-1}(\beta_2)]-C \text{ by Theorem \ref{ThDriftOfOrigin}}\\
&=\frac{1}{2}\sum_{j=nN_0}^{(n+L)N_0}[\sgn m_{2j}(\beta)-\sgn m_{2j-1}(\beta)]-C=z_{n+L,\beta}-z_{n,\beta}-C. 
\end{align*}
Thus, part (1) holds whenever  $L>\max\{L_1, (\sK+2\mathsf C_{\beta_2}(R))/c_{\beta_2}(R)\}$. 

 If $(2m_{2nN_0}(\beta),\ldots,2m_{2(n+L)N_0}(\beta))=(2,6,\ldots,2,6)$ and $\xi>r$,   then 
\begin{align*}
&\hskip-25mm
\wt{\E}^{\xi}_{\beta}\left(\sum_{k=n+1}^{n+L}f_{k,\beta}(\wt{X}^{\xi}_{k,\beta})\right)\geq -\sK+
\wt{\E}^{\xi}_{\beta_3}\left(\sum_{k=1}^{L}f_{k,\beta_3}(\wt{X}^{\xi}_{k,\beta_3})\right)\\
&\hskip-25mm
\geq -\sK-2\mathsf C_{\beta_3}(R)+\PP_{L,\beta_3}'(\xi)
\geq -\sK-2\mathsf C_{\beta_3}(R)+\PP_{L,\beta_3}'(r)\\
&\hskip-25mm
\geq -\sK-2\mathsf C_{\beta_3}(R)+\PP_{L,\beta_3}'(r/2)+c_{\beta_3}(R)L(r/2)\\
&\hskip-25mm
\geq \PP_{L,\beta_3}'(r/2)+C,\text{ provided $L>(2/r)(C+\sK+2\mathsf C_{\beta_3}(R))/c_{\beta_3}(R)$}\\
&\hskip-25mm
\geq z_{L,\beta_3}+C,\text{ provided $L>L_1$}
\end{align*}
\hskip7mm$=z_{n+L,\beta}-z_{n,\beta}+C$ by the assumption on the even CFE of $\beta$ and Theorem \ref{ThDriftOfOrigin}.
\end{proof}

We will use the following abbreviations:
\begin{align*}
S_{n,\beta}&:=\sum_{k=1}^n f_{n,\beta}(X_{n,\beta}) \ , \ \wt{S}_{n,\beta}^\xi:=\sum_{k=1}^n f_{n,\beta}(\wt{X}_{n,\beta}^\xi),\\
V_{n,\beta}&:=\mathrm{Var}(S_{n,\beta})\ , \  \wt{V}_{n,\beta}^\xi:=\mathrm{Var}(\wt{S}_{n,\beta}^\xi).
\end{align*}
Note that 
 $S_{n,\beta}$ and $\wt{S}^\xi_{n,\beta}$ are the same functions, viewed on different probability spaces. Thus,  for every Borel function $\varphi:\R\to\R$,  
 $$
\wt{\E}^\xi_\beta(\vf(\wt{S}_{n,\beta}^\xi))=\wt{\E}_\beta^\xi\left(\vf(S_{n,\beta})\right)
=\E_\beta\left(
\frac{e^{\xi S_{n,\beta}}h_{n,\beta}(X_{n,\beta},\xi)}{e^{\PP_{n,\beta}(\xi)} h_{1,\beta}(X_{1,\beta},\xi)}\cdot \vf(S_{n,\beta})
\right).
$$

\begin{lem}\label{l.V-beta-xi-growth}
For every $R>0$ there exist $C_1,C_2>1$ such that for all $\beta\in\Lambda_6$, $|\xi|\leq R$, and  $n$, $C_1^{-1} n-C_2\leq \wt{V}_{n,\beta}^\xi\leq C_1 n +C_2.$
\end{lem}
\begin{proof}
By Lemma \ref{Lemma-f(X)},  $\DS K:=\sup_{\beta\in\Lambda_6}\sup_n \|f_{n,\beta}\|<\infty$. 

By Proposition \ref{Prop-Doeblin}, there exists $\eps_0'>0$ independent of $\beta$ such that  $\{X_{n,\beta}\}$ are uniformly elliptic with the same ellipticity constant $\eps_0'$ for every $\beta\in\Lambda_6$. 
By  Lemma \ref{LmAddPres},\index{Change of measure!density estimates} there are constants $C(R)$ so that for all $\beta\in\Lambda_6$, $|\xi|\leq R$, and $n\in\N$,  
\begin{equation}\label{e.RN-est}
C_1(R)^{-1}\leq \frac{d\pi^\xi_{\beta,n,n+1}}{d\pi_{\beta,n,n+1}}\leq C_1(R). 
\end{equation}
It follows that for some constant $\eps_0>0$, for all $\beta\in\Lambda_6$ and $|\xi|\leq R$, $\{\wt{X}_{n,\beta}\}$ are uniformly elliptic with the same ellipticity constant $\eps_0$. 

Let  $u_{n,\beta}(\xi)$ denote the structure constants  associated to $\{\wt{X}_{n,\beta}^\xi\}$  in  \eqref{e.u}. By Theorem~\ref{LmVarCycles}, there are $C_2,C_3$ depending only on $\eps_0$ and $K$ such that 
\begin{equation}\label{e.v-beta-n-xi-ineq}
C_2^{-1}\sum_{k=3}^n u_{k,\beta}^2(\xi)-C_3\leq \wt{V}_{n,\beta}^\xi\leq C_2 \sum_{k=3}^n u_{k,\beta}^2(\xi)+C_3. 
\end{equation}

The same argument leading to the  equation array \eqref{e.state-lower-bound-xi} gives a constant $\eps_1>0$ such that $\Prob[\wt{X}^\xi_{n,\beta}=a]$ is either zero or bigger than $\eps_1$,  for all $\beta\in\Lambda_6$, $|\xi|\leq R$, and $n\in\N$. This and  \eqref{e.RN-est} imply that the hexagon measures of $\{\wt{X}_{n,\beta}^\xi\}$ and $\{{X}_{n,\beta}\}$ are equal up to a density which is uniformly bounded away from zero and infinity by constants independent of $\beta$ and $\xi$. It follows that there is a positive constant $C_4$ such that 
$$
C_4^{-1}\leq u_{k,\beta}(\xi)/u_{k,\beta}(0)\leq C_4\text{ for all $\beta\in\Lambda_6$, $|\xi|\leq R$, and $n\in\N$.}
$$

Looking at \eqref{e.v-beta-n-xi-ineq}, we see that for some $C_5,C_6$, for all $\beta\in\Lambda_6$, $|\xi|\leq R$, and $n\in\N$, 
$$
C_5^{-1}\sum_{k=3}^n u_{k,\beta}^2-C_6\leq \wt{V}_{n,\beta}^\xi\leq C_5 \sum_{k=3}^n u_{k,\beta}^2+C_6.
$$
Here $u_{k,\beta}:=u_{k,\beta}(0)$, the structure constants of $\{X_{n,\beta}\}$. 
By Theorem \ref{Theorem-D_N-Growth}, there are constants $C_7,C_8$ such that for all $\beta\in\Lambda_6$, 
$\DS C_7^{-1}n-C_8\leq \sum_{k=3}^n u_{n,\beta}^2\leq C_7 n+C_8$, and the lemma follows.
\end{proof}

\begin{lem}\label{l.nu(R)}
For every $R>0$, there exist constants $\nu_1(R)$ and $v(R)>0$  such that for all $n>\nu_1(R)$, $\beta\in\Lambda_6$, and $0\!<\!\xi_1\!<\!\xi_2\!<\!R$,\;\;
$\DS 
\wt{\E}^{\xi_2}_{\beta}(\wt{S}_{n,\beta}^{\xi_2})-\wt{\E}^{\xi_1}_{\beta}(\wt{S}_{n,\beta}^{\xi_1})\geq v(R)(\xi_2-\xi_1)n.
$
\end{lem}
\begin{proof}
Fix $R>0$. 
By Lemma \ref{LmAddPres}, there is constant $\mathsf C$ such that  for every $\beta\in\Lambda_6$, $|\xi|\leq R$, and $n\in\N$,   
$\mathsf C(R)^{-1}\leq h_{n,\beta}(\cdot,\xi)\leq \mathsf C(R)$. If  $|\xi_1|,|\xi_2|\leq R$, then 
\begin{align*}
1&=\wt{\E}^{\xi_2}_\beta(1)=\E_\beta\left(\frac{e^{\xi_2 S_{n,\beta}} h_{n,\beta}(X_{n,\beta},\xi_2)}{e^{\PP_{n,\beta}(\xi_2)}h_{1,\beta}(X_{1,\beta},\xi_2)}\right)\\
&=\wt{\E}_\beta\left(\frac{e^{\xi_1 S_{n,\beta}}h_{n,\beta}(X_{n,\beta},\xi_1)}{e^{\PP_{n,\beta}(\xi_1)}h_{1,\beta}(X_{1,\beta},\xi_1)}\cdot \frac{e^{\PP_{n,\beta}(\xi_1)}h_{1,\beta}(X_{1,\beta},\xi_1)}{h_{n,\beta}(X_{n,\beta},\xi_1)}\frac{h_{n,\beta}(X_{n,\beta},\xi_2)}{e^{\PP_{n,\beta}(\xi_2)}h_{1,\beta}(X_{1,\beta},\xi_2)}\cdot  e^{(\xi_2-\xi_1)S_{n,\beta}}\right)\\
&= \mathsf C(R)^{\pm 4} e^{\PP_{n,\beta}(\xi_1)-\PP_{n,\beta}(\xi_2)}\wt{\E}^{\xi_1}_\beta(e^{(\xi_2-\xi_1)S_{n,\beta}}), \text{ whence}
\end{align*}
\begin{equation}\label{e.delta-P-id}
e^{\PP_{n,\beta}(\xi_2)-\PP_{n,\beta}(\xi_1)}=\mathsf C(R)^{\pm 4} \wt{\E}^{\xi_1}_\beta(e^{(\xi_2-\xi_1)S_{n,\beta}}).
\end{equation}

Similarly, writing $\eta:=\xi_2-\xi_1$ and $\wh{S}_{n,\beta}^{\xi_1}:={S}_{n,\beta}-\wt{\E}_\beta^{\xi_1}(\wt{S}_{n,\beta}^{\xi_1})$ we have 
\begin{align}
&\wt{\E}^{\xi_2}_\beta (\wt{S}_{n,\beta}^{\xi_2})-\wt{\E}^{\xi_1}_\beta (\wt{S}_{n,\beta}^{\xi_1})=
\wt{\E}_\beta^{\xi_2} \left(\wt{S}_{n,\beta}^{\xi_2}-\wt{\E}^{\xi_1}_\beta (\wt{S}_{n,\beta}^{\xi_1})\right)
={\E}_\beta \left(\wh{S}_{n,\beta}^{\xi_1}\cdot \frac{e^{\xi_2 S_{n,\beta}} h_{n,\beta}(X_{n,\beta},\xi_2)}{e^{\PP_{n,\beta}(\xi_2)}h_{1,\beta}(X_{1,\beta},\xi_2)}\right)\notag\\
&={\E}_\beta \left(\wh{S}_{n,\beta}^{\xi_1} e^{\eta S_{n,\beta}}e^{\PP_{n,\beta}(\xi_1)-\PP_{n,\beta}(\xi_2)}\cdot \frac{e^{\xi_1 S_{n,\beta}} h_{n,\beta}(X_{n,\beta},\xi_1)}{e^{\PP_{n,\beta}(\xi_1)}h_{1,\beta}(X_{1,\beta},\xi_1)}\cdot \frac{h_{n,\beta}(X_{n,\beta},\xi_2)h_{1,\beta}(X_{1,\beta},\xi_1)}{h_{n,\beta}(X_{n,\beta},\xi_1)h_{1,\beta}(X_{1,\beta},\xi_2)}\right) \notag\\
&=\wt{\E}_\beta^{\xi_1} \left(\wh{S}_{n,\beta}^{\xi_1} e^{\eta S_{n,\beta}}\cdot\frac{e^{\PP_{n,\beta}(\xi_1)}}{e^{\PP_{n,\beta}(\xi_2)}}\cdot \frac{h_{n,\beta}(X_{n,\beta},\xi_2)h_{1,\beta}(X_{1,\beta},\xi_1)}{h_{n,\beta}(X_{n,\beta},\xi_1)h_{1,\beta}(X_{1,\beta},\xi_2)}\right)
\geq \frac{\wt{\E}_\beta^{\xi_1} \left(\wh{S}_{n,\beta}^{\xi_1} e^{\eta S_{n,\beta}}\right)}{\mathsf C(R)^{8} \wt{\E}_\beta^{\xi_1} \left(e^{\eta S_{n,\beta}}\right)},\text{ by \eqref{e.delta-P-id}} \notag
\end{align}
\begin{align}
&=\mathsf C(R)^{-8} \frac{\wt{\E}_\beta^{\xi_1} \left(\wh{S}_{n,\beta}^{\xi_1} e^{\eta \wh{S}_{n,\beta}^{\xi_1}}\right)}{ \wt{\E}_\beta^{\xi_1} \left(e^{\eta \wh{S}_{n,\beta}^{\xi_1}}\right)}=\mathsf C(R)^{-8}\int_0^\eta \frac{d}{d\theta}\left[\frac{\wt{\E}_\beta^{\xi_1} \left(\wh{S}_{n,\beta}^{\xi_1} e^{\theta \wh{S}_{n,\beta}^{\xi_1}}\right)}{ \wt{\E}_\beta^{\xi_1} \left(e^{\theta \wh{S}_{n,\beta}^{\xi_1}}\right)}\right]d\theta \notag\\
&=\mathsf C(R)^{-8}\int_0^\eta \frac{d}{d\theta}\left[\frac{\wt{\E}_\beta^{\xi_1} \left({S}_{n,\beta} e^{\theta {S}_{n,\beta}}\right)}{ \wt{\E}_\beta^{\xi_1} \left(e^{\theta {S}_{n,\beta}}\right)}\right]d\theta.\label{e.Newton-Leibnitz-Trick}
\end{align}

To  calculate the integrand, it is convenient to use the following elementary fact: 
\begin{sublem*}
(\cite[Lemma 7.22]{Dolgopyat-Sarig-Book}).
Let $X,Y$ be bounded random variables on the same probability space $\Omega$,  such that $\mathrm{Var}(X)<\infty$ and $C^{-1}y\leq Y\leq Cy$ for some constants $C>1$, $y>0$. Then the variance of $X$ with respect to the measure $\frac{Yd\Prob}{\E(Y)}$, given by \index{Change of measure!variance estimates}
$$
\Var^Y(X):=\frac{\E(X^2 Y)}{\E(Y)}-\left(\frac{\E(XY)}{\E(Y)}\right)^2, 
$$ 
satisfies the bounds $C^{-4}\Var(X)\leq \Var^Y(X)\leq C^{4}\Var(X)$. 
\end{sublem*}
\noindent
{\em Proof }
Let $(X_1,Y_1)$, $(X_2,Y_2)$ be two independent copies of $(X,Y)$, then \\
$\Var^Y(X)=\frac{1}{2}\frac{\E[(X_1-X_2)^2 Y_1 Y_2]}{\E(Y_1 Y_2)}=C^{\pm 4} \frac{1}{2}\E[(X_1-X_2)^2]=C^{\pm 4}\Var(X)$.
\qed

\smallskip 
Let $Y_\theta\!\! :=\!\!\frac{e^{\PP_{n,\beta}(\xi_1)}}{e^{\PP_{n,\beta}(\xi_1+\theta)}}\cdot \frac{h_{n,\beta}(X_{n,\beta},\xi_1)h_{1,\beta}(X_{1,\beta},\xi_1+\theta)}{h_{n,\beta}(X_{n,\beta},\xi_1+\theta)h_{1,\beta}(X_{1,\beta},\xi_1)}
\!\!=\!\!\mathsf C(R)^{\pm 4}\frac{e^{\PP_{n,\beta}(\xi_1)}}{e^{\PP_{n,\beta}(\xi_1+\theta)}}$.
By direct differentiation
\begin{align*}
&\frac{d}{d\theta}\left[\frac{\wt{\E}_\beta^{\xi_1} \left({S}_{n,\beta} e^{\theta {S}_{n,\beta}}\right)}{ \wt{\E}_\beta \left(e^{\theta {S}_{n,\beta}}\right)}\right]=
\frac{\wt{\E}_\beta \left(({S}_{n,\beta})^2 e^{\theta {S}_{n,\beta}}\right)}{ \wt{\E}_\beta^{\xi_1} \left(e^{\theta {S}_{n,\beta}}\right)}
-\left[\frac{\wt{\E}_\beta^{\xi_1} \left({S}_{n,\beta} e^{\theta {S}_{n,\beta}}\right)}{ \wt{\E}_\beta^{\xi_1} \left(e^{\theta {S}_{n,\beta}}\right)}\right]^2\\
&\!\!=\!\!\frac{\wt{\E}_\beta^{\xi_1+\theta} \left(({S}_{n,\beta})^2 Y_\theta\right)}{ \wt{\E}_\beta^{\xi_1+\theta} \left(Y_\theta\right)}
\!-\!\left[\frac{\wt{\E}_\beta^{\xi_1+\theta} \left({S}_{n,\beta}Y_\theta\right)}{ \wt{\E}_\beta^{\xi_1+\theta} \left(Y_\theta\right)}\right]^2
\!\!\!\!=\!\!\Var^{Y_\theta}(\wt{S}_{n,\beta}^{\xi_1+\theta})\geq 
\frac{\Var(\wt{S}_{n,\beta}^{\xi_1+\theta})}{\mathsf C(R)^4}
\!\!=\!\!\frac{\wt{V}_{n,\beta}^{\xi_1+\theta}} {\mathsf C(R)^4}
.
\end{align*}
If $0\leq \theta\leq \eta$, then $\xi_1+\theta\in [\xi_1,\xi_2]\subset [0,R]$, and by Lemma \ref{l.V-beta-xi-growth}, there are constants $C_1(R), C_2(R)>0$  such that for all $\beta\in\Lambda_6$ and $0\leq \theta\leq \eta$, 
$
\wt{V}_{n,\beta}^{\xi_1+\theta}\geq C_1(R)^{-1}n-C_2(R).
$
Looking at \eqref{e.Newton-Leibnitz-Trick}, we deduce that if $\eta:=\xi_2-\xi_1>0$, and $n>2C_1(R)C_2(R)$,  then
$\DS
\wt{\E}^{\xi_2}_\beta (\wt{S}_{n,\beta}^{\xi_2})-\wt{\E}^{\xi_1}_\beta (\wt{S}_{n,\beta}^{\xi_1})\geq 
\mathsf C(R)^{-12}\int_0^\eta \wt{V}_{n,\beta}^{\xi_1+\theta}d\theta
\geq 
\frac{n(\xi_2-\xi_1)}{2\mathsf C(R)^{12}C_1(R)}>0. 
$.
\end{proof}

\subsection{Proof of Theorem \ref{CrWeakBal}} 
The theorem states that for every $\xi\in\R$ there is an irrational $\beta$ of bounded type such that $\tomega_0$ (the midpoint of the top horizontal side of $\Rect_0$) is $\mu_\xi$-generic for the linear flow $\phi_\beta:\St\to\St$.

Recall that every irrational $\beta\in (-1,1)$ can be uniquely expanded into a continued fraction 
$\beta\!\!=\!\! [0;2m_1(\beta),2m_2(\beta),\ldots]$  where $m_i(\beta)\!\!\in\!\! \Z\setminus\{0\}$.  
Given $b_1,\ldots,b_n\!\!\in\!\! 2\Z\setminus\{0\}$, let  
$$
[0;b_1,\ldots,b_n]:=\{\beta\in (-1,1)\setminus\Q: 2m_i(\beta)=b_i\ (1\leq i\leq n)\}.
$$
The key to the proof are the following two lemmas:

\begin{lem}\label{l.P1-P3}
Suppose $0<\xi<\xi+1<R$. 
 There exist positive integers $n_0, N_0, L,c_0$ and  even integers $b_1,b_2,\ldots\in \{2,6,-6\}$ such that for every $n> n_0$ and $\beta\in [0;b_1,\ldots,b_{nLN_0}]$, the following  holds:
\begin{enumerate}[(P1)]
\item $\exists!\xi_{nL,\beta}\in [0,R]$ such that  $\wt{\E}^{\xi_{nL,\beta}}_\beta(\wt{S}_{nL,\beta}^{\xi_{nL,\beta}})=z_{nL,\beta}$;
\item Suppose $n> n_0+1$, then 
\begin{enumerate}[--]
\item If $\xi_{(n-1)L,\beta}\leq \xi$, then $\xi_{nL,\beta}>\xi_{(n-1)L,\beta}$
\item If $\xi_{(n-1)L,\beta}\geq  \xi+\frac{c_0}{n-1}$, then $\xi_{nL,\beta}<\xi_{(n-1)L,\beta}$
\item If $\xi_{(n-1)L,\beta}\in (\xi, \xi+\frac{c_0}{n-1})$, then $|\xi_{nL,\beta}-\xi|\leq \frac{c_0}{n-1}$.
\end{enumerate}
\item  Suppose $n> n_0+1$, then $|\xi_{nL,\beta}-\xi_{(n-1)L,\beta}|<c_0/(n-1)$.
\end{enumerate}
\end{lem}
\begin{lem}\label{l.beta-is-nice}
Let $b_i$ be as in the previous lemma, and set $\beta=[0;b_1,b_2,b_3,\ldots]$. Then $\xi_{nL,\beta}\xrightarrow[n\to\infty]{}\xi$. 
\end{lem}
\begin{proof}[Proof of Lemma \ref{l.P1-P3}]
We fix once and for all the following constants:
\begin{itemize}
\item $N_0$, a positive  even integer independent of $\beta\in\Lambda_6$ as in \eqref{e.mat-prod-form}. 
\item $K$, a constant such that $|f_{n,\beta}|\leq K$ for all $n\in\N$ and $\beta\in\Lambda_6$ (Lemma \ref{Lemma-f(X)}).
\item $\mathsf K$, the constant $\mathsf K(R)$ from Lemma \ref{LmBetaCont}.
\item $E$, a positive integer such that $|\E_\beta(S_{n,\beta})|\leq E$ for all $n\geq 1$ and $\beta\in\Lambda_6$. Such a constant exists because of Proposition \ref{Prop-Drift}. 
\item $\nu(R)\in\N$ and $v(R)>0$ such that for all $n\geq \nu(R)$ and $\beta\in\Lambda_6$, 
\begin{enumerate}[(i)]
\item $\forall 0<\xi_1<\xi_2<R+1$, $\wt{\E}_\beta^{\xi_2}(\wt{S}^{\xi_2}_{n,\beta})>\wt{\E}_\beta^{\xi_1}(\wt{S}^{\xi_1}_{n,\beta})+v(R)(\xi_2-\xi_1)n$ (Lemma \ref{l.nu(R)});
\item Thus, $\eta\mapsto \wt{\E}_\beta^\eta(\wt{S}^\eta_{n,\beta})$ is increasing on $[0,R]$, and 
\item  $\forall \eta\in [\xi,R]$,  $\wt{\E}^\eta_\beta(\wt{S}^\eta_{n,\beta})>v(R)\xi n-E$. 
\end{enumerate}
\item $L$, a positive  integer such that \eqref{e.Dima-Trick-1} and \eqref{e.Dima-Trick-2} hold with $C:=\max\{E,2\mathsf K\}$, and so that $L>\max\{1,\nu(R),E/N_0,1/v(R)\}$.
\item $M:=4(N_0+K+\mathsf K)$.
\item  $c_0:=M/v(R)$.
\item $n_0$, a positive integer such that 
$\DS n_0>\max\left\{\frac{2E+LN_0}{v(R)\xi}, \frac{M}{v(R)}\right\}$.
\end{itemize}

Define the strings
$\DS
\un{B}_{\sqrt{2}}=(\,\underset{2LN_0}{\underbrace{-6,6,\ldots,-6,6}}\,)\ ,\ \un{B}_{\sqrt{3}}=(\,\underset{2LN_0}{\underbrace{2,6,\ldots,2,6}}\,).
$ 
The sequence $(b_1,b_2,b_3,\ldots)$ will be a concatenation $(\un{B}_1,\un{B}_2,\ldots)$ where $\un{B}_i=\un{B}_{\sqrt{2}}$ or $\un{B}_{\sqrt{3}}$. 

First we declare that $\un{B}_1,\cdots,\un{B}_{n_0}:=\un{B}_{\sqrt{3}}$, and $\un{B}_{n_0+1}:=\un{B}_{\sqrt{2}}$. 
Let us see  (P1). Suppose $\beta\in [0;\un{B}_1,\ldots,\un{B}_{n_0+1}]$, then 
\begin{align*}
z_{(n_0+1) L,\beta}&=\frac{1}{2}\sum_{k=1}^{(n_0+1) LN_0} [\sgn(b_{2k}(\beta))-\sgn(b_{2k-1}(\beta))]\\
&=\frac{1}{2}\sum_{k=1}^{n_0 LN_0} [\sgn(6)-\sgn(2)]+\frac{1}{2}\sum_{k=1}^{LN_0} [\sgn(6)-\sgn(-6)]=LN_0;\\
\wt{\E}^0_{\beta}(\wt{S}^0_{(n_0+1) L})&<z_{(n_0+1)L,\beta},\text{ because }\wt{\E}^0_{\beta}(\wt{S}^0_{(n_0+1) L,\beta})\leq E<LN_0;\\
\wt{\E}^{\xi}_{\beta}(\wt{S}^\xi_{(n_0+1) L,\beta})&>z_{(n_0+1)L,\beta},\text{ because  } \wt{\E}^{\xi}_{\beta}(\wt{S}^\xi_{(n_0+1) L,\beta})>v(R)\xi n_0-E>LN_0.
\end{align*}
The function $\eta\mapsto \wt{\E}^{\eta}_{\beta}(\wt{S}^\eta_{(n_0+1) L,\beta})$ is continuous (Lemma \ref{LmAddPres}), and strictly increasing on $[0,R]$ (since $(n_0+1)L>\nu(R)$). Therefore there exists a unique $\xi_{(n_0+1)L,\beta}\in [0,R]$ such that $\wt{\E}^{\xi_{(n_0+1)L,\beta}}_{\beta}(\wt{S}^\xi_{(n_0+1) L,\beta})=z_{(n_0+1)L,\beta}$. This shows (P1) for $n=n_0+1$. (P2) and (P3) hold vacuously for $n=n_0+1$. 

\medskip
Assume by induction that we have already constructed  $\un{B}_1,\ldots,\un{B}_n\in \{\un{B}_{\sqrt{2}},\un{B}_{\sqrt{3}}\}$ with $n\geq n_0$, so that   properties (P1)--(P3) hold. We define 
$$
\un{B}_{n+1}:=
\begin{cases}
\un{B}_{\sqrt{2}} & \exists\gamma\in [0;\un{B}_1,\ldots,\un{B}_n]\cap \Lambda_6$ s.t.  $\xi_{nL,\gamma}\leq \xi;\\
\un{B}_{\sqrt{3}} & \forall\gamma\in [0;\un{B}_1,\ldots,\un{B}_n]\cap \Lambda_6$,   $\xi_{nL,\gamma}>\xi.
\end{cases}
$$
We claim that (P1)--(P3) hold for $n+1$.

\medskip
\noindent
{\bf Case 1.\/} Suppose  $\exists\gamma\in [0;\un{B}_1,\ldots,\un{B}_n]\cap \Lambda_6$ s.t.  $\xi_{nL,\gamma}\leq \xi$. 
Then  $\un{B}_{n+1}=\un{B}_{\sqrt{2}}$. 

 Fix $\beta\in [0;\un{B}_1,\ldots,\un{B}_{n+1}]\cap \Lambda_6$. By Theorem \ref{ThDriftOfOrigin}, $z_{nL,\beta}=z_{nL,\gamma}$, and
 \[z_{(n+1)L,\beta}=z_{nL,\beta}+\frac{1}{2}\sum_{k=1}^{LN_0} [\sgn(6)-\sgn(-6)]=z_{nL,\beta}+LN_0=z_{nL,\gamma}+LN_0.
\]
By Lemma \ref{LmBetaCont}, $\wt{\E}^{\xi_{nL,\gamma}}_\beta(\wt{S}^{\xi_{nL,\gamma}}_{nL,\beta})=\wt{\E}^{\xi_{nL,\gamma}}_{\gamma}(\wt{S}^{\xi_{nL,\gamma}}_{nL,\gamma})\pm \mathsf K=z_{nL,\gamma}\pm\mathsf K=z_{nL,\beta}\pm\mathsf K$.\\
By Lemma \ref{l.Dima-Trick}(1), 
$\wt{\E}^{\xi_{nL,\gamma}}_\beta(\wt{S}^{\xi_{nL,\gamma}}_{(n+1)L,\beta}-\wt{S}^{\xi_{nL,\gamma}}_{nL,\beta})
\!\!<\!\! z_{(n+1)L,\beta}-z_{ nL,\beta}-C\leq z_{(n+1)L,\beta}-z_{ nL,\beta}
-2\mathsf K.$\\
Adding the two inequalities, we obtain  
$$\wt{\E}^{\xi_{nL,\gamma}}_\beta(\wt{S}^{\xi_{nL,\gamma}}_{(n+1)L,\beta})<z_{(n+1)L,\beta}.$$

Let $M:=4(N_0+K+\mathsf K)$ and  $\eta:=\xi_{nL,\gamma}+\frac{M}{2v(R) n}$. Note that $\eta\leq \xi+{\frac{1}{2}}<R$. Therefore, 
\begin{align*}
\wt{\E}^{\eta}_\beta(\wt{S}^{\eta}_{(n+1)L,\beta})&\geq \wt{\E}^{\eta}_\beta(\wt{S}_{nL,\beta}^\eta)-KL
\geq \wt{\E}^{\xi_{nL,\gamma}}_\beta(\wt{S}_{nL,\beta}^\eta)+{v(R)(\eta-\xi_{nL,\gamma})nL}-KL\\
&=z_{nL,\beta}\pm\mathsf K+ML-KL=(z_{(n+1)L,\beta}-LN_0)\pm\mathsf K+ML-KL\\
&>z_{(n+1)L,\beta}+L(M-N_0-K-\mathsf K)>z_{(n+1)L,\beta}.
\end{align*}
By the intermediate value theorem, 
\begin{equation}\label{e.sol-interval}
\exists!\xi_{(n+1)L,\beta}\in \left(\xi_{nL,\gamma}\ ,\ \xi_{nL,\gamma}+\frac{M}{2v(R) n}\right)\;\; \text{ s.t. }
\wt{\E}^{\xi_{nL,\beta}}_\beta(\wt{S}^{\xi_{nL,\gamma}}_{(n+1)L,\beta})=z_{(n+1)L,\beta}.
\end{equation}
Note that $0<\xi_{(n+1)L,\beta}<\xi+\frac{M}{v(R)n_0}\leq\xi+1$, therefore $\xi_{(n+1)L,\beta}\in [0,R]$. 

If 
$\xi_{nL,\beta}\leq \xi$, then we can repeat the previous argument with $\gamma:=\beta$,  and obtain 
$$
\xi_{(n+1)L,\beta}>\xi_{nL,\beta}
\ , \  |\xi_{(n+1)L,\beta}-\xi_{nL,\beta}|\leq \frac{M}{v(R) n}.
$$

If $\xi_{nL,\beta}\geq \xi+\frac{M}{v(R)n}$, then $\xi_{(n+1)L,\beta}<\xi_{nL,\gamma}+\frac{M}{v(R)n}<\xi+\frac{M}{v(R)n}\leq \xi_{nL,\beta}$, and 
$$
\xi_{(n+1)L,\beta}<\xi_{nL,\beta}.
$$

If $\xi_{nL,\beta}\in(\xi,\xi+\frac{M}{v(R)n})$, then  we cannot determine whether $\xi_{(n+1)L,\beta}<\xi_{nL,\beta}$. But we can estimate 
$|\xi_{(n+1)L,\beta}-\xi|$, as follows. We saw above that 
$$
\wt{\E}^{\xi_{nL,\gamma}}_\beta(\wt{S}^{\xi_{nL,\gamma}}_{nL,\beta})=\wt{\E}^{\xi_{nL,\gamma}}_{\gamma}(\wt{S}^{\xi_{nL,\gamma}}_{nL,\gamma})\pm \mathsf K=z_{nL,\gamma}\pm\mathsf K=z_{nL,\beta}\pm\mathsf K.
$$
Since  $z_{nL,\beta}=\wt{\E}^{\xi_{nL,\beta}}_\beta(\wt{S}^{\xi_{nL,\beta}}_{nL,\beta})$,  this leads to 
 $
 |\wt{\E}^{\xi_{nL,\gamma}}_\beta(\wt{S}^{\xi_{nL,\gamma}}_{nL,\beta})-\wt{\E}^{\xi_{nL,\beta}}_\beta(\wt{S}^{\xi_{nL,\beta}}_{nL,\beta})|\leq \mathsf K.
 $
However   $ |\wt{\E}^{\xi_{nL,\gamma}}_\beta(\wt{S}^{\xi_{nL,\gamma}}_{nL,\beta})-\wt{\E}^{\xi_{nL,\beta}}_\beta(\wt{S}^{\xi_{nL,\beta}}_{nL,\beta})|\geq v(R)|\xi_{nL,\gamma}-\xi_{nL,\beta}|nL$. By our assumptions,  $\xi_{nL,\gamma}$ and $\xi_{nL,\beta}$ are on opposite sides of $\xi$, so  
$
|\xi_{nL,\gamma}-\xi|<\frac{\mathsf K}{v(R)nL}<\frac{M}{4v(R)nL}.  
$ By  \eqref{e.sol-interval},
\[
|\xi_{(n+1)L,\beta}-\xi|<|\xi_{(n+1)L,\beta}-\xi_{nL,\gamma}|+|\xi_{nL,\gamma}-\xi|<\frac{3M}{4v(R) n}<\frac{M}{v(R)n}\equiv \frac{c_0}{n}
\]
We obtained (P1) and (P2) for $n+1$. 

Next we check (P3). By definition 
\begin{align*}
&\wt{\E}^{\xi_{(n+1)L,\beta}}_\beta(\wt{S}^{\xi_{(n+1)L,\beta}}_{nL,\beta})
=\wt{\E}^{\xi_{(n+1)L,\beta}}_\beta(\wt{S}^{\xi_{(n+1)L,\beta}}_{(n+1)L,\beta})\pm KL=z_{(n+1)L,\beta}\pm KL
=z_{nL,\beta}\pm 2KL\\
&=\wt{\E}^{\xi_{nL,\beta}}_\beta(\wt{S}^{\xi_{nL,\beta}}_{nL,\beta})\pm 2KL.
\end{align*}
Therefore $v(R) n{L}|\xi_{(n+1)L,\beta}-\xi_{L,\beta}|<|\wt{\E}^{\xi_{(n+1)L,\beta}}_\beta(\wt{S}^{\xi_{(n+1)L,\beta}}_{nL,\beta})-\wt{\E}^{\xi_{nL,\beta}}_\beta(\wt{S}^{\xi_{nL,\beta}}_{nL,\beta})|<2KL$, whence 
$|\xi_{(n+1)L,\beta}-\xi_{L,\beta}|<\frac{2KL}{v(R) nL}<c_0/n$. 

\medskip
\noindent
{\bf Case 2.}  Suppose $\forall\gamma\in [\un{B}_1,\ldots,\un{B}_n]$, $\xi_{nL,\gamma}>\xi$. Then $\un{B}_{n+1}=\un{B}_{\sqrt{3}}$. 

Fix $\beta\in [0;\un{B}_1,\ldots,\un{B}_{n+1}]$.
By Lemma \ref{l.Dima-Trick}(2), 
\begin{align*}
\wt{\E}^{\xi_{nL,\beta}}_\beta(\wt{S}^{\xi_{nL,\beta}}_{(n+1)L,\beta})&=
\wt{\E}^{\xi_{nL,\beta}}_\beta(\wt{S}^{\xi_{nL,\beta}}_{L,\beta})+\wt{\E}^{\xi_{nL,\beta}}_\beta(\wt{S}^{\xi_{nL,\beta}}_{(n+1)L,\beta}-\wt{S}^{\xi_{nL,\beta}}_{nL,\beta})\\
&>z_{nL,\beta}+ [z_{(n+1)L,\beta}-z_{nL,\beta}+C]>z_{(n+1)L,\beta}+C>z_{(n+1)L,\beta}.
\end{align*}
Also,  
 $\DS
z_{(n+1)L,\beta}=z_{nL,\beta}+\frac{1}{2}\sum_{k=1}^{LN_0}[\sgn(6)-\sgn(2)]=z_{nL,\beta}
$, so by the choice of $n_0$,
\begin{align*}
z_{(n+1)L,\beta}&=z_{nL,\beta}=
\wt{\E}^{\xi_{nL,\beta}}_\beta(\wt{S}^{\xi_{nL,\beta}}_{nL,\beta})>v(R)\xi_{nL,\beta}nL-E>v(R)\xi n_0 L
-E>E>\wt{\E}^{0}_\beta(\wt{S}^{0}_{nL,\beta}).
\end{align*}
So $\wt{\E}^{0}_\beta(\wt{S}^{0}_{nL,\beta})<z_{(n+1)L,\beta}
<\wt{\E}^{\xi_{nL,\beta}}_\beta(\wt{S}^{\xi_{nL,\beta}}_{(n+1)L,\beta})
$. Recalling that $\eta\mapsto \wt{\E}^{\eta}_\beta(\wt{S}^{\eta}_{(n+1)L,\beta})$ is continuous and strictly increasing on $[0,R]$, we find that 
$$
\exists!\xi_{(n+1)L,\beta}\in (0,\xi_{nL,\beta})\text{ s.t. }
\wt{\E}^{\xi_{(n+1)L,\beta}}_\beta(\wt{S}^{\xi_{(n+1)L,\beta}}_{(n+1)L,\beta})=z_{(n+1)L,\beta}.
$$
Since $\xi_{nL,\beta}\in [0,R]$ by the induction hypothesis, 
$\xi_{(n+1)L,\beta}\in [0,R]$ and we have established (P1) for $n+1$.

(P3) for $n+1$ can be proved as in Case 1. 

We check (P2). 
In Case 2, $\xi_{nL,\beta}>\xi$, and we saw  above that 
$\xi_{(n+1)L,\beta}<\xi_{nL,\beta}$. 
It remains to check that $\xi_{nL,\beta}\in (\xi,\xi+\frac{c_0}{n})\Rightarrow |\xi_{(n+1)L,\beta}-\xi|<c_0/n$.

On the one hand, 
$
\xi_{(n+1)L,\beta}<\xi_{nL,\beta}<\xi+\frac{c_0}{n}$. On the other hand,  
by (P3) for $n+1$, 
$
\xi_{(n+1)L,\beta}=\xi_{nL,\beta}-|\xi_{nL,\beta}-\xi_{(n+1)L,\beta}|
>\xi-\frac{c_0}{n}
$. 
So  $|\xi_{(n+1)L,\beta}-\xi|<\frac{c_0}{n}$ as required. 

\noindent
This completes the inductive construction of $\un{B}_i$. Set $(b_1,b_2,\ldots):=(\un{B}_1,\un{B}_2,\ldots)$.  
\end{proof}

\begin{proof}[Proof of Lemma \ref{l.beta-is-nice}] Let $\un{B}_i$ be the sequence of blocks constructed in the previous proof, and set
$
\beta:=[0;\un{B}_1,\un{B}_2,\ldots].
$

We claim that $\un{B}_i$ is not  equal to $\un{B}_{\sqrt{2}}$ for all $i$ large enough. Otherwise, $\un{B}_{n}=(-6,6,\ldots,-6,6)$ for all $n>n_1$. Fix $C>|z_{n_1,\beta}|+Kn_1$, then by Lemma \ref{l.Dima-Trick}(1), for all $n$ large enough, for all $\eta\leq R$, 
$$
\wt{\E}^\eta_\beta(\wt{S}_{n,\beta}^\eta)<z_{n,\beta}-z_{n_1,\beta}-C+Kn_1<z_{n,\beta}.
$$  
Necessarily, $\xi_{nL,\beta}>R$ for all $n$ large enough, in contradiction to (P1). 

Similarly, $\un{B}_i$ is not  equal to $\un{B}_{\sqrt{3}}$ for all $i$ large enough.  Otherwise, $\un{B}_{n}=(2,6,\ldots,2,6)$ for all $n>n_1$. Fix $C>|z_{n_1,\beta}|+Kn_1$ and $r>0$ arbitrarily small, then by Lemma \ref{l.Dima-Trick}(2), for all $n$ large enough, for all $r<\eta\leq R$, 
$$
\wt{\E}^\eta_\beta(\wt{S}_{n,\beta}^\eta)>z_{n,\beta}-z_{n_1,\beta}+C-Kn_1>z_{n,\beta}.
$$  
Necessarily, $\xi_{nL,\beta}<r$ for all $n$ large enough, and $\xi_{nL,\beta}\xrightarrow[n\to\infty]{}0$. But this is impossible, because by (P2), this sequence is increases when its values drop below $\xi$.

We are now ready to show that $\xi_{nL,\beta}\xrightarrow[n\to\infty]{}\xi$.
Consider the following  three cases:

(A)  $\forall\gamma\in [0;\un{B}_1,\ldots,\un{B}_n]\cap \Lambda_6$, $\xi_{nL,\gamma}\leq \xi$.   This case cannot happen  for all $n$ large enough, because $\un{B}_n$ is not eventually constant equal to $\un{B}_{\sqrt{2}}$.

(B) $\forall\gamma\in [0;\un{B}_1,\ldots,\un{B}_n]\cap\Lambda_6$, $\xi_{nL,\gamma}>\xi$. This case cannot happen  for all $n$ large enough, because  $\un{B}_n$ is not   eventually constant equal to $\un{B}_{\sqrt{3}}$. 

(C) $\exists\gamma_1,\gamma_2\in [0;\un{B}_1,\ldots,\un{B}_n]\cap\Lambda_6$ such that  $\xi_{nL,\gamma_1}<\xi<\xi_{nL,\gamma_2}$.

Arguing as in Case 1 in the proof of the previous lemma, it is possible to show that in Case (C), if $n>n_0$, then   all the $\xi_{nL,\gamma}$ with $\gamma\in  [0;\un{B}_1,\ldots,\un{B}_n]\cap\Lambda_6$ are within  distance $c_0/n$ from each other. In particular $\xi_{nL,\gamma_1}, \xi_{nL,\gamma_2}$ and $\xi_{nL,\beta}$ are all within distance $c_0/n$ from each other. Since $\xi_{nL,\gamma_i}$ are on opposite sides of $\xi$, 
$|\xi_{nL,\gamma_i}-\xi|<c_0/n$, whence also 
$
|\xi_{nL,\beta}-\xi|<c_0/n.$

Similarly, if at time $n$ we are in case (A) and in time $n+1$ we are in case (B) (or the opposite), then $\xi$ is between $\xi_{nL,\beta}$ and $\xi_{(n+1)L,\beta}$, whence  $|\xi_{nL,\beta}-\xi|<c_0/n$ by (P3).

 By the foregoing discussion 
there is an infinite sequence $n_k$ such that
$|\xi_{n_k L,\beta}-\xi|\leq \frac{c_0}{n_k}.$ 
Therefore to complete the proof
it suffices to show that  if for some sufficiently
large $\brn$, 
$|\xi_{\brn L,\beta}-\xi|\leq c_0/\brn$, then 
\begin{equation}
\label{Trap}
|\xi_{n L,\beta}-\xi|\leq c_1/(\brn-1)\text{ for all $n>\brn$, where $c_1:=2c_0+1$.}
\end{equation}

We prove this by induction. The base $n=\brn$ holds by assumption, so we just need to establish the step
of induction. So suppose that the claim holds for $n-1$.

If $\xi_{(n-1)L,\beta}\in (\xi, \xi+\frac{c_0}{n-1})$
 then \eqref{Trap} holds by the third
condition in (P2).

If $\xi-c_0/(\brn-1)\leq \xi_{(n-1)L, \beta}<\xi$, then \eqref{Trap} holds 
since
$$\DS \xi-c_0/(\brn-1)\leq \xi_{(n-1)L, \beta}\overset{(P2)}{<}\xi_{nL, \beta}\overset{(P3)}{<}\xi_{(n-1)L, \beta}+c_1/(n-1)<\xi+c_0/(\brn-1) .$$
The case where
$\xi+\frac{c_0}{n-1}<\xi_{(n-1)L,\beta}<\xi+c_1/(\brn-1)$ 
can be analyzed similarly. This completes the proof of \eqref{Trap} and finishes the proof of the lemma.
\end{proof}






\medskip
\noindent
{\bf Proof of Theorem \ref{CrWeakBal}.} Let $\beta:=[0;\un{B}_1,\un{B}_2,\ldots]$ as above. Henceforth we work with this $\beta$ only, and drop the subscripts $\beta$. 

Let  $\tomega_0$ be the midpoint of the top of rectangle $\Rect_0$.  
Note that $z_{NN_0}(\tomega_0)=z_{N,\beta}$ (see \eqref{e.z-n}).
We saw above that there exist $\xi_{nL}\xrightarrow[n\to\infty]{}\xi$ such that $\wt{\E}^{\xi_{nL}}(\wt{S}^{\xi_{nL}}_{nL})=z_{nL,\beta}$. 
 For integers $N\in [nL,(n+1)L)$, let  $\eta_N:=\xi_{nL}$ and $z_{N,\beta}:=z_{NN_0}(\tomega_0)$.  Then it is easy to see using the uniform boundedness of the jump functions that 
$$
\wt{\E}^{\eta_{N}}(\wt{S}^{\eta_{N}}_{N})=z_{N,\beta}+O(1)\text{ for all }N.
$$  
By Lemma \ref{l.14.3}, 
$
\PP_N'(\eta_N)=z_{N,\beta}+O(1), \text{ and }\eta_N\to\xi. 
$
By Lemma \ref{l.solution-xi}, there are $\eta_N'\to\xi$ such that 
$
\PP_N'(\eta_N')=z_{NN_0}(\tomega_0).
$
Looking at the proof of Proposition \ref{PrBulk}, we see that $\tomega_0$ must be  $\mu_\xi$-generic for the linear flow with velocity ${\beta\choose 1}$. \qed

\subsection{Proof of Theorem \ref{t.exp-tilted-exist}.} The theorem says that for every $\xi/2$, there is an irrational of bounded type $\alpha$ such that $D_N(\alpha)$ is exponentially tilted with parameter $\xi/2$.

By Theorem \ref{CrWeakBal}, there is an irrational $\beta$ of bounded type such that $\tomega_0$, the midpoint of the top side of $\Rect_0$, is $\mu_{\xi/2}$-generic for the linear flow with velocity ${\beta\choose 1}$. 
Let $\alpha:=\frac{\beta+1}{2}$. Since $\beta$ is badly approximable by rationals, $\alpha$ is badly approximable by rationals. In particular, $\alpha$ is an irrational of bounded type.

By Lemma \ref{l.Generic-Good-For-Rectangles} and the $\mu_{\xi/2}$-genericity of $\tomega_0$,  for all $k,\ell\in\Z$, 
$$
\lim_{T\to\infty}\frac{\int_0^T 1_{\Rect_k}(\phi^t(\tomega_0))dt}{\int_0^T 1_{\Rect_\ell}(\phi^t(\tomega_0))dt}=e^{\xi(k-\ell)/2}, 
$$
and $\int_0^T 1_{\Rect_\ell}(\phi^t(\tomega_0))dt\to\infty$ for all $\ell\in\Z$. It now follows from  \eqref{e.Lin-Flow-1}, that for all $k,\ell\in\Z$, 
$$
\lim_{n\to\infty}\frac{\#\{1\leq N\leq n: D_N(\alpha)=k/2\}}{\#\{1\leq N\leq n: D_N(\alpha)=\ell/2\}}=e^{\xi(k-\ell)/2}. 
$$
In other words, $D_N(\alpha)$ is exponentially tilted with parameter $\xi$.\qed

\begin{remark}
\label{RemZeroNG}
\normalfont
A slight modification of the construction   allows to build an irrational of bounded type $\alpha$, for which $D_N(\alpha)$  is not $\rho$-equidistributed for any $\rho$. To do so fix two thresholds $\hxi_1$ and $\hxi_2$ so that $\hxi_1<\hxi_2$, and build  a sequence of blocks equal to $\un{B}_{\sqrt{2}}$ or $B_{\sqrt{3}}$ as follows. We start with
block $\un{B}_{\sqrt{3}}$. Then  we will keep adding blocks $\un{B}_{\sqrt{2}}$
until the time $k_1$ when $\xi_{k_1L,\gamma}>\hxi_2$ for all $\gamma\in [\un{B}_1,\ldots,\un{B}_{k_1}]$. Then we  keep adding  blocks $\un{B}_{\sqrt{3}}$ until the time $k_2$ when $\xi_{k_2 L,\gamma}<\hxi_1$ for all $\gamma\in [\un{B}_1,\ldots,\un{B}_{k_2}]$. Then we switch back again to blocks $\un{B}_{\sqrt{2}}$ until $\xi_{k_3 L,\gamma}>\hxi_2$ for all $\gamma\in [\un{B}_1,\ldots,\un{B}_{k_1}]$, and so on. 
A modification of the  proof of Theorem \ref{t.exp-tilted-exist} shows if  $\beta$ has even CFE $[0;\un{B}_1,\un{B}_2,\cdots]$, then $\tomega_0$  has a historical measure $\mu_\xi$ for $\phi_\beta$,   for every $\xi\in [\hxi_1, \hxi_2].$ It follows that $D_N(\alpha)$ is not $\rho$-equidistributed for any $\rho$.
We leave the details to the reader.
\end{remark}

\section{Hausdorff Dimension and Baire Category of Generic Points}
\label{ScNonGen}
In this section we prove Theorems \ref{CrHDGen} and \ref{CrHDExcept} on the Hausdorff dimensions and Baire categories of the sets of  generic and non-generic points. 

Theorem \ref{CrHDGen} says that for some $c_1,c_2,\xi_0>0$, the Hausdorff dimension of   $\mathsf{Gen}(\xi):=\{\tomega:\tomega\text{ is $\mu_\xi$-generic}\}$ satisfies the bounds 
$$
2-c_1\xi^2\leq \mathrm{HD}(\mathsf{Gen}(\xi))\leq 2-c_2\xi^2\ \ \ (|\xi|\leq \xi_0). 
$$

We divide the proof into several steps. 
Given $\brxi$, let 
$$
S_\brxi:=\{\tomega : \exists \xi\geq \brxi\text{  s.t. $\mu_\xi$ is a historical measure of $\tomega$}\}.
$$ 

\begin{lem}[Upper Bound]
\label{LmHDUpper} 
For some $c_1>0$, for all $\brxi>0$ sufficiently close to zero, for every  $E^u_0$-segment $W$,
$\DS \HD(S_\brxi\cap W)\leq 1-c_1 \brxi^2. $
\end{lem}

\begin{proof}
 It is sufficient to prove the lemma for  pure $E^u_0$-segments with depth zero and height zero (in particular 
 $\wt{\zc}_0\!\!\equiv\!\! 0$ on $W$, and $\pi(W)$ is a maximal $u$-fibre in some $\<\dot{a}_0\>\!\!\in\! \mathfrak P_0$). 
Fix some constant $C$ such that for all $n\!\!\geq\!\! 0$ and for all $\eta\!\!\in\!\! [0,R]$,
$\sup|\wt{\zc}_n-\zc|\!\!<\!\!C$ and 
$\left|\PP_n'(\eta)-\left.\frac{d}{d\xi}\right|_{\xi=\eta}\log\E(e^{\xi S_n})\right|<C
$, 
see Lemmas \ref{Lemma-Modified-Approx-Canonical} and \ref{LmPressure}.

By Theorem \ref{t.historical-crit} and  
Proposition \ref{p.Pressure-and-Geometric-Pressure},  if $\tomega\in S_{\brxi}\cap W$, then for some  $n_k\uparrow\infty$, there are   $\eta_{n_k}\to\brxi$ such that  $\PP_{n_k}'(\eta_{n_k})=z_{n_kN_0}(\tomega)$.
 
 Since $\brxi>0$ and $\PP_n$ are strictly convex with second derivative bounded below by $\const n$ on $[\brxi/2,\brxi]$, there exists $k_0$ such that  for  all $k\geq k_0$,  $z_{n_kN_0}(\tomega)=\PP_{n_k}'(\brxi)>\PP_{n_k}'(\brxi/2)+10C$, 
whence $\wt{\zc}_{n_k}(\tomega)-\wt{\zc}_0(\tomega)\geq \PP_{n_k}'(\brxi/2)+8C$. By the definition of the jump functions, 
$
\DS \ov{S}_{n_k}(\tomega)\equiv\sum_{j=1}^{n_k} f_{j}\big(X^{(0)}_j(\pi(\tomega))\big)=\wt{\zc}_{n_k}(\tomega)-\wt{\zc}_0(\tomega)\geq \PP_{n_k}'(\brxi/2)+8C.
$
Thus, $$\DS S_\brxi\subset\bigcap_{n_0} \bigcup_{n\geq n_0} K_n,\text{ where }
K_n:=\{\tomega\in W: \ov{S}_n(\tomega)>\PP_{n}'(\brxi/2)+8C\}.
$$

Each $K_n$ is a union of non-overlapping $E^s_0$-segments $I_{j,n}$, which project isometrically to  maximal unstable segments in cylinders of the form 
$
A_{j,n}=\<\dot{a}_0^{j,n},\ldots,a_{n}^{j,n}\>
$, $a_0^{j,n}=a_0$. In particular,  $|I_{j,n}|=\ell^u(A_{j,n})\geq\const e^{-\kappa n}$, where $\kappa>0$ comes from a uniform upper bound on the expansion of the automorphisms $\psi_i$ in the unstable direction, and the constant comes from a uniform lower bound on the unstable side of the Markov rectangles in $\mathfrak P_n$ for all $n$. It follows that for all $j,n$, and $0<s<1$,
$$
|I_{j,n}|^{1-s}\leq \const |I_{j,n}|e^{\kappa n s}. 
$$

By construction, $S_{\brxi}$ is covered by $\{I_{j,n}\}$, and 
\begin{align*}
&\sum_{j,n}|I_{j,n}|^{1-s}\leq  \const\sum_n e^{\kappa n s}\sum_j |I_{j,n}|\leq
\const\sum_n  e^{\kappa n s} \Prob_W[\ov{S}_n(\tomega)>\PP_n'(\brxi/2)+8C]  \\
&\leq
\const\sum_n  e^{\kappa n s} \Prob_W\biggl[\ov{S}_n(\tomega)>\frac{d}{d\xi}\bigg|_{\xi=\brxi/2}\log\E(e^{\xi S_n})+7C\biggr]  \\
&\equiv \const \sum_n e^{ n \kappa s} \Prob \biggl[{S}_n>\frac{d}{d\xi}
\bigg|_{\xi=\brxi/2}\log\E(e^{\xi S_n})+7C\bigg| X_0=a_0\biggr]\\
&\leq \const \sum_n e^{ n \kappa s} \Prob \biggl[{S}_n>\frac{d}{d\xi}\bigg|_{\xi=\brxi/2}\log\E(e^{\xi S_n})+7C\biggr],\text{because $a_0$ is constant}\\
&\leq \const \sum_n e^{ n \kappa s} \Prob\left[S_n>V_n\mathcal F_n'(\brxi/2)\right],
\text{ where }\mathcal F_n(\eta):=\frac{1}{V_n}\log \E(e^{\eta S_n}).
\end{align*}
Let $\cI_n$ denote  the {\em rate functions,}\index{Rate functions} equal by definition to the   Legendre transforms of $\mathcal F_n$. By the bounds on $\cI_n$ in  Proposition \ref{PrRateFn}, and  Chernoff's bounds,
 there is a positive constant $c_1$ such that for all  $\brxi\in [-R,R]$, 
$\Prob[S_n>V_n\mathcal F_n'(\brxi/2)]=O(e^{-c_1 n\brxi^2+o(n)}).
$
So
$\DS    \sum_{j,n} |I_{j,n}|^{1-s}\!\!\leq\!\! \const\sum_n e^{n(\kappa s-c_1 \brxi^2)+o(n)}<\infty$ for  $\DS s\!\!<\!\!\frac{c_1 \brxi^2}{\kappa}.$ Thus 
 $\mathrm{HD}(S_{\brxi})\leq 1-c_1\kappa^{-1}\brxi^2$.
\end{proof}

The lower bound requires some preparation. Let $X_n$ be the Markov chain from Section \ref{s.Markov-Chain}, and let $\PP_n(\xi)$ be the pressure functions from \S\ref{ss.change-of-measure}.

\begin{lem}\label{l.delta-pressure-is-uni-conv}
For every $R>0$ there are constants $\mathsf{C},n_0>0$ such that for all $k,n>n_0$ and $|\xi|\leq R$, 
$\mathsf C^{-1} n\leq \mathrm{Var}(\wt{S}_{k+n}^\xi-\wt{S}_k^\xi)\leq \mathsf C n$ and 
 $\mathsf{C}^{-1}n\leq \PP''_{k+n}(\xi)-\PP''_{k}(\xi)\leq \mathsf{C}n$.
\end{lem}
\begin{proof}
By Lemma \ref{l.uniform-Doeblin}, the Markov chains $\{\wt{X}_i^\xi\}_{i\geq k}$ are  uniformly elliptic with the same ellipticity constant for all $|\xi|\leq R$. 
By Theorem \ref{LmVarCycles}, 
there are  constants $C_1,C_2>0$  such that for all $k,n\geq n_0$, and for every $|\xi|\leq R$,  
$$
C_1 \sum_{j=k+3}^k u_j^2(\xi)-C_1< \mathrm{Var}(\wt{S}_{k+n}^\xi-\wt{S}_k^\xi)<C_1 \sum_{j=k+3}^k u_j^2(\xi)+C_2,
$$
where $u_j(\xi)$ are the structure constants of $\{\wt{X}_i^\xi\}_{i\geq 1}$. 

For $|\xi|\leq R$, the hexagon measures of  $\{\wt{X}_i^\xi\}_{i\geq 1}$ and the hexagon measures of $\{{X}_i\}_{i\geq 1}$ are related by Radon-Nikodym derivatives uniformly bounded away from zero and infinity by constants independent of $j$ and $\xi$. Therefore $u_j^2(\xi)\asymp u_j^2$, where $u_j$ are the structure constants of $\{X_n\}_{n\geq 1}$, and  by Claim 2 on page \pageref{page-hexagon-claim}, $\sum_{j=k+3}^{k+n} u_j^2(\xi)\asymp n$ uniformly in $k$ and $\xi\in [-R,R]$. Necessarily, 
$$
\exists c_1,n_1\text{ s.t. }\forall k,n\geq n_1 \forall|\xi|\leq R,\  c_1^{-1} n\leq \mathrm{Var}(\wt{S}_{k+n}^\xi-\wt{S}_k^\xi)\leq c_1 n.
$$

Now we appeal to Lemma 7.17 in \cite{Dolgopyat-Sarig-Book} and its proof. Fixing $k$, and applying Equation (7.21) in the proof of this lemma to the truncated chain $\{{X}_i\}_{i\geq k}$, we see that 
\begin{equation}\label{e.big-Ohs}
\PP_{k+n}''(\xi)-\PP_k''(\xi)=\mathrm{Var}(\wt{S}_{k+n}^\xi-\wt{S}_k^\xi)
+O(1)\sqrt{\mathrm{Var}(\wt{S}_{k+n}^\xi-\wt{S}_k^\xi)}+O(1), 
\end{equation}
where the big-Ohs are uniform in $n$ and $\xi\in [-R,R]$. Combining this with the estimates for $\mathrm{Var}(\wt{S}_{k+n}^\xi-\wt{S}_k^\xi)$ above, we obtain that for all $k>n_1$, for some $n_2(k)$ and $\mathsf C(k)$, 
$$
\forall n\geq n_2(k) \forall \xi\in[-R,R],\  \mathsf C(k)^{-1} n\leq \PP_{k+n}''(\xi)-\PP_k''(\xi)\leq \mathsf C(k)n
$$

It remains to show that $\mathsf C(k)$ can be replaced by constants independent of $k$. It is sufficient to verify that the big-Ohs in \eqref{e.big-Ohs} are uniform in $k$. 
A close look at the proof of Lemma 7.17 in \cite{Dolgopyat-Sarig-Book} shows that  the big Ohs depend only on the $L^\infty$ norms of   
$$
\epsilon_{k,n}(x,y,\xi):=\frac{h_k(x,\xi)}{h_{k+n}(y,\xi)}\frac{d}{d\xi}\left(\frac{h_{k+n}(y,\xi)}{h_k(x,\xi)}\right)\text{ and }\frac{d}{d\xi}\epsilon_{k,n}(x,y).
$$
By Lemmas 7.14 and 7.16 in \cite{Dolgopyat-Sarig-Book},  $h_j(x,\xi), 1/h_j(x,\xi), \partial_\xi h_j(x,\xi)$ and $\partial_\xi^2 h_j(x,\xi)$ are uniformly bounded by constants independent of $j$ and $\xi$, provided $|\xi|\leq R$. Thus $\epsilon_{k,n}$ and $\partial_\xi\epsilon_{k,n}$ are uniformly bounded, and $\mathsf C(k)$ can be taken independent of $k$.  
\end{proof}

\begin{lem}\label{l.LLT-LD-Blocks}
There exists $c>0$ such that the following holds. Suppose  $0<|\xi|<1$ and $z_n=z_n(\xi):=\lfloor\PP_n'(\xi)\rfloor$. Then for all sufficiently large integers $L$,
$$
\Prob\left(S_{n+L}-S_n=z_{n+L}-z_{n}\, \big|\, X_n=x_n\right)>e^{-c\xi^2 L} \text{ for all $n\in\N$ and  $x_n\in\mathfrak P_n$}. 
$$
\end{lem}
\begin{proof}
Assume by way of contradiction that the lemma is false. Then there are  $c_i\to\infty$, $0<|\xi_i|<1$,  $L_i\uparrow\infty$,  $n_i\in\N$,  and $x_{n_i}\in\mathfrak P_{n_i}$ such that 
\begin{equation}\label{e.ass-by-cont}
\Prob\left(S_{n_i+L_i}-S_{n_i}=z_{n_i+L_i}(\xi_i)-z_{n_i}(\xi_i)\, \big|\, X_{n_i}=x_{n_i}\right)<e^{-c_i \xi_i^2 L_i}. 
\end{equation}

Let $\wh{S}_i:=S_{n_i+L_i}-S_{n_i}$, $\wh{\PP}_i(t)=\PP_{n_i+L_i}(t)-\PP_{n_i}(t)$, and $\wh{z}_i:=z_{n_i+L_i}(\xi_i)-z_{n_i}(\xi_i)$. 
Let $\wt{X}_n^\xi$ denote the change of measure built in \S\ref{ss.change-of-measure}, and let 
$\DS \wh{S}_i^\xi:=\sum_{j=n_i+1}^{n_i+L_i}f_i(\wt{X}_i^\xi)$.

\medskip
\noindent
{\em Step 1.\/} Uniform  estimates of the expectation and variance of $\wh{S}_i^{\xi_i}$, given $\wt{X}_{n_i}^{\xi_i}=x_{n_i}$.

\medskip

By Lemma \ref{l.uniform-Doeblin}, $\wt{X}^\xi_n$ are uniformly elliptic with the same ellipticity constant for all $|\xi|<1$. An exponential mixing argument based on Corollary \ref{Cor-phi-mixing} shows that
\begin{align}
&\wt{\E}^{\xi}(\wh{S}_i^{\xi_i}|\wt{X}_{n_i}^{\xi_i}=x_{n_i})=\wt{\E}^{\xi_i}(\wh{S}_i^{\xi_i})+O(1)=
\wt{\E}^{\xi_i}(\wt{S}_{n_i+L}^{\xi_i})-\wt{\E}^{\xi}(\wt{S}_{n_i}^{\xi_i})+O(1)\\
&=\PP_{n_i+L_i}'(\xi_i)-\PP_{n_i}'(\xi_i)+O(1),\text{ see \eqref{Pres1stDer}. Thus,  } \notag\\
&\wt{\E}^{\xi_i}(\wh{S}_i^{\xi_i}|\wt{X}_{n_i}^{\xi_i}=x_{n_i})=\wh{z}_i+O(1).\label{e.com-works}
\end{align}

Conditioning on the first state does not affect the ellipticity constant, therefore the  Markov chains $(\wt{X}_n^\eta)_{n\geq n_i}$ conditioned on $\wt{X}_{n_i}^{\eta}=x_{n_i}$ are also all uniformly elliptic with the same ellipticity constant $\eps_0$ . By Theorem \ref{LmVarCycles},  there are constants $C_1, C_2$ such that for all $|\eta|\leq 1$ and $i\geq 1$, 
$$
 C_1^{-1}\sum_{j=n_i+3}^{n_i+L_i} u_n^2(\eta)-C_2\leq \Var(\wh{S}^\eta_i|\wt{X}_{n_i}^{\eta}=x_{n_i})\leq C_1\sum_{j=n_i+3}^{n_i+L_i} u_n^2(\eta)+C_2,
$$
where $u_n(\eta)$ are the structure constants of $(\wt{X}_n^\eta)_{n\geq 1}$, see \eqref{e.u}. It is easy to see that the hexagon measures of  
$(\wt{X}_n^\eta)_{n\geq 1}$ differ from those of $({X}_n)_{n\geq 1}$ by densities which are uniformly bounded away from zero and infinity. 
Hence $u_n(\eta)\asymp u_n$, where $u_n=u_n(0)$ are the structure constants of $(X_n)_{n\geq 1}$. Looking at Theorem \ref{Theorem-D_N-Growth}, and specializing to $\eta=\xi_i$, we deduce that there is a constant $C_3$ such that for all  $i\geq 1$, 
\begin{equation}\label{e.var-xi-block-cond}
C_3^{-1} L_i\leq \Var(\wh{S}_i^{\xi_i}|\wt{X}_{n_i}^{\xi_i}=x_{n_i})\leq C_3 L_i. 
\end{equation}

\medskip
\noindent
{\em Step 2.\/} Uniform estimates for $\wh{\PP}_i'(0)$ and $\wh{\PP}_i''(\eta)$, $|\eta|\leq 1$.

%
%

$L_i\to\infty$, so by Lemma \ref{l.delta-pressure-is-uni-conv}, there exists  $C_5>0$ such that  for all $i$  large enough, 
$$
C_5^{-1} L_i\leq \wh{\PP}_i''(\cdot)\leq C_5 L_i\text{ on }[-1,1].
$$
Next by   \eqref{Pres1stDer},  $\PP_i'(0)=\PP_{n_i+L_i}'(0)-\PP_{n_i}'(0)=\E(S_{n_i+L_i})-\E(S_{n_i})$, so by Proposition \ref{Prop-Drift}, 
\begin{equation}\label{e.P-i-O(1)}
\wh{\PP}_i'(0)=O(1). 
\end{equation}

\medskip
\noindent
{\em Step 3.\/} Uniform estimate of $\wh{\PP}_i(\xi_i)-\xi_i\wh{\PP}_i'(\xi_i)$. 

\medskip
Let $\varphi_i(\eta):=\wh{\PP}_i(\eta)-\eta\wh{\PP}_i'(\eta)$, then $\varphi_i'(\eta)=-\eta\wh{\PP}_i''(\eta)$, and 
$\varphi_i(0)=0$, see \eqref{Pres0thDer}. Therefore, 
\begin{align*}
\varphi_i(\xi_i)&=\int_0^{\xi_i} \varphi_i'(\eta)d\eta\geq -C_5L_i\int_0^{\xi_i}\eta d\eta=-\frac{1}{2}C_5 L_i{\xi_i}^2. 
\end{align*}
It follows that for some {\em positive} constant $C_6$ (which depends on $\xi$), for all $i$ large enough, 
\begin{equation}\label{e.Legendre-Block-estimate}
\wh{\PP}_i(\xi)-\xi\wh{\PP}_i'(\xi)\geq -C_6 L_i\xi_i^2. 
\end{equation}

\medskip
\noindent
{\em Step 4.\/} Estimate of $\Prob[\wh{S}_i=\wh{z}_i|X_{n_i}=x_{n_i}]$.

\medskip
Recall that $C_4^{-1}\leq h_k(\cdot,\xi)\leq C_4$ for all $k$. Therefore, 
\begin{align}
&\Prob[\wh{S}_i=\wh{z}_i|X_{n_i}=x_{n_i}]\geq C_4^{-2}e^{\wh{\PP}_i(\xi_i)-\xi_i\wh{z}_i}\E\left(1_{[\wh{S}_i=\wh{z}_i]}\cdot e^{\xi_i\wh{S}_i}\frac{h(X_{n_i+L_i},\xi_i)}{e^{\wh{\PP}_i(\xi_i)} h_{n_i}(X_{n_i},\xi_i)}\, \bigg|\, X_{n_i}=x_{n_i}\right)\notag\\
&\equiv C_4^{-2} e^{\wh{\PP}_i(\xi_i)-\xi(\lfloor{\PP}_{n_i+L_i}'(\xi_i)\rfloor-\lfloor{\PP}_{n_i}'(\xi_i)\rfloor)}\wt{\E}^{\xi_i}(1_{[\wh{S}_i^{\xi_i}=\wh{z}_i]}
|\wt{X}^{\xi_i}_{n_i}=x_{n_i})\ \ (\wt{X}^{\xi_i}_{n_i}=x_{n_i}\Leftrightarrow X_{n_i}=x_{n_i})\notag\\
&\geq C_4^{-2}e^{-2} e^{\wh{\PP}_i(\xi_i)-\xi_i\wh{\PP}_{i}'(\xi_i)}\wt{\E}^{\xi_i}(1_{[\wh{S}_i^{\xi_i}=\wh{z}_i]}
|\wt{X}^{\xi_i}_{n_i}=x_{n_i})\notag\\
&\geq C_7 e^{-C_6 L_i\xi_i^2}\Prob[\wh{S}_i^\xi=\wh{z}_i|\wt{X}^\xi_{n_i}=x_{n_i}],\text{ where $C_7:=C_4^{-2} e^{-2}$, see  \eqref{e.Legendre-Block-estimate}}\label{e.scottish}
\end{align}
Note that $C_7$ is independent of $i$.

To estimate $\Prob[\wh{S}_i^{\xi_i}=\wh{z}_i|\wt{X}^{\xi_i}_{n_i}=x_{n_i}]$, we consider the Markov array
$$
\begin{array}{cccclll}
  \wt{X}^{\xi_1}_{n_1+1} & \cdots &  \wt{X}^{\xi_1}_{n_1+L_1} &                     & &\text{conditioned on }\wt{X}^{\xi_1}_{n_1}=x_{n_1}&\\
  \wt{X}^{\xi_2}_{n_2+1} & \cdots &       \cdots          & \wt{X}^{\xi_2}_{n_2+L_2} & &\text{conditioned on }\wt{X}^{\xi_2}_{n_2}=x_{n_2}&\\
  \wt{X}^{\xi_3}_{n_3+1} & \cdots &       \cdots          & \cdots               &  \wt{X}^{\xi_3}_{n_2+L_3} & \text{conditioned on }\wt{X}^{\xi_3}_{n_3}=x_{n_3}  & \\
        \vdots       &        &                       &                      &                 & \ddots&\\                 
\end{array}
$$
together with additive functional generated by the jump functions $f_k(\wt{X}_k^{\xi_i})$. The rows have lengths $L_i\uparrow\infty$, and  the row sum are equal in distribution to $\wh{S}_i^{\xi_i}$, conditioned on $\wh{X}^{\xi_i}_{n_i}=x_{n_i}$. Moreover,

(1)  This array is uniformly elliptic, because $(\wt{X}^\eta_n)_{n\geq 1}$ is uniformly elliptic with the same ellipticity constant for all $|\eta|\leq 1$ (Lemma \ref{l.uniform-Doeblin}). 

(2) The variance of the row sums tends to infinity by  \eqref{e.var-xi-block-cond} and since  $L_i\to\infty$.

(3) Our additive functional is  hereditary\index{Hereditary property} in the sense of  Definition \ref{d.hereditary} in Appendix~\ref{AppMC}, by    Claim 2 on page \pageref{page-hexagon-claim} 
and Lemma \ref{Lemma-D_N-Growth}.

(4)  Irreducibility: The previous argument also shows that the co-range of our array is ${2\pi}\Z$. By Theorem \ref{t.co-range},  the essential range is $\Z$. Since the jump functions are integer-valued, this is also the algebraic range, and we obtain irreducibility. 

(5) Finally, and crucially,  
by  \eqref{e.com-works} and \eqref{e.var-xi-block-cond}, 
 $$
 \frac{\wt{\E}^{\xi_i}(\wh{S}_i^{\xi_i}|\wt{X}^{\xi_i}_{n_i}=x_{n_i})-\wh{z}_i}{{\sqrt{\Var(\wh{S}_i^\xi|\wt{X}^{\xi_i}_{n_i}=x_{n_i})}}}
 =O\left(\frac{1}{\sqrt{L_i}}\right)\xrightarrow[i\to\infty]{}0.
 $$
By the local limit theorem for arrays, Theorem \ref{ThLLTArrays}, $$\Prob[\wh{S}_i^\xi=\wh{z}_i|\wt{X}^\xi_{n_i}=x_{n_i}]\sim 1/\sqrt{2\pi\Var(\wh{S}_i^\xi|\wt{X}_{n_i}^\xi=x_{n_i})}.$$ By  \eqref{e.var-xi-block-cond}, there is a constant $C_8$ such that for all $i$ large enough, 
$$
\Prob[\wh{S}_i^\xi=\wh{z}_i|\wt{X}^\xi_{n_i}=x_{n_i}]\geq \frac{C_8}{\sqrt{L_i}}. 
$$
Plugging this into \eqref{e.scottish}, we obtain that there is a positive constant $C_9$ such that  for all $i$ large enough, 
$$
\Prob[\wh{S}_i=\wh{z}_i|X_{n_i}=x_{n_i}]\geq \frac{C_7 C_8 e^{-C_6 L_i\xi_i^2}}{\sqrt{L_i}}>e^{-C_9 L_i \xi_i^2}.
$$
But this contradicts \eqref{e.ass-by-cont}, because $c_i\to\infty$. 
\end{proof}

\begin{lem}[Lower Bound]
\label{LmHDLower}
For some $c_2>0$, for all $|\xi|\leq 1$, for every $E^u_0$-segment $W$, 
$\DS \HD(\mathsf{Gen}(\xi)\cap W)\geq 1-c_3 \xi^2. $
\end{lem}

\begin{proof}
It is sufficient to consider $W$ such that $\wt{\zc}_0\equiv 0$ on $W$, and so that $\pi(W)$ is a maximal $E^u_0$-segment in some cylinder $\<\dot{a}_0,\ldots,a_k\>$. For simplicity, we take $k=0$ (the modifications needed for the general case are routine). Let 
$$
\ov{S}_n(\wt{\omega}):=\sum_{j=1}^n f_j\big(X^{(0)}_k(\pi(\wt{\omega}))\bigr)=\wt{\zc}_n\bigl((\tpsi_{n-1}\circ\cdots\circ\tpsi_0)(\wt{\omega})\bigr)-\wt{\zc}_0(\wt{\omega}).
$$
By  Lemma \ref{Lemma-Modified-Approx-Canonical}, 
$z_{nN_0}(\tomega)=\ov{S}_n(\wt{\omega})+O(1)\text{ on }W.$

Fix $|\xi|\leq 1$, let $z_n=\lfloor\PP_n'(\xi)\rfloor$, and let  $L\!\!=\!\!L(\xi)$ be a constant,  to be specified later.
Let
$$
K_n(\xi):=\{\wt{\omega}\in W: \ov{S}_{nL}(\wt{\omega})=z_{nL}\}, \ K(\xi):=\bigcap_{n} K_n(\xi).
$$ 
If  $\tomega\in K(\xi)$, then $z_{nLN_0}(\tomega)=\ov{S}_{nL}(\tomega)+O(1)=\PP_{nL}'(\xi)+O(1)$.  By Lemma \ref{l.solution-xi}, $\exists \xi_n'\to\xi$ such that $\PP_{nL}'(\xi_n')=z_{nLN_0}(\tomega)$. Also, $\mathcal P_{nLN_0}'(\xi_{nLN_0}(\tomega))=z_{LN_0}(\tomega)$, by definition. Looking at  Proposition \ref{p.Pressure-and-Geometric-Pressure}, we deduce that  $\xi_{nLN_0}(\tomega)\to\xi$, and by  Theorem \ref{t.historical-crit}, $\mu_\xi$ is a historic measure for $\tomega$. This shows that
$$
S_{\brxi}\cap W\supset K(\xi)\text{ for all }\xi\geq \brxi.
$$

By Frostman's lemma,\index{Frostman's lemma} \index{Hausdorff dimension} to prove that $\HD(K(\xi))>1-s$, it is sufficient to construct  a probability measure $\nu$ on $K(\xi)$ such that 
\begin{equation}\label{e.mass-dist}
\exists\delta_0,c_0>0\text{ s.t. for all intervals $I$, }|I|<\delta_0\Rightarrow \nu(I)\leq c_0|I|^{1-s}.
\end{equation}
$K_n(\xi)$ is a union of non overlapping linear segments $I_{n,j}$ so that for each $j$, $\pi(I_{n,j})$ is a  maximal $E^u_{0}$-segment in some cylinder $\<\dot{x}_0^{n,j},\ldots,x_{nL}^{n,j}\>$, $x_0^{n,j}=a_0$. It follows that  for all $j,j'$, either $I_{n,j}$ contains the interior of $I_{n+1,j'}$,  or the interiors of $I_{n,j}$ and $I_{n+1,j'}$ are disjoint.  Thus the interiors of $I_{n,j}$ form a semi-algebra $\mathcal S$ of subsets of $W\setminus\{\text{endpoints of $I_{n,j}$}\}$.

We define a set function $\nu:\mathcal S\to [0,1]$ as follows. If $n=0$, there is just one interval $I_{0,1}=W$ and we define $\nu(I_{0,1})=1$. Suppose by induction we defined $\nu(I_{n,j})$ for all $j$. Let 
$\{I_{n,j, \ell}\}$ be the enumeration of the $I_{n+1,\ell}$ with interiors inside $I_{n,j}$, then we set 
$$ 
\nu(I_{n,j,\ell}):=\nu(I_{n,j})\times\frac{|I_{n, j, \ell}|}{\sum_t |I_{n, j, t}|}. 
$$
A standard argument shows that  $\nu$ is $\sigma$-additive on $\mathcal S$. By Carath\'eodory's extension theorem, $\nu$   extends to a Borel probability measure on $W$. 

The key to the structure of $\nu$ is the following estimate:

\medskip
\noindent
{\em Step 1.\/}  There is a constant $c>0$ such that for every $\xi$, 
if $L=L(\xi)$ is large enough,  then for some  $n_0>0$ and every $n\geq n_0$, 
$
\DS \sum_\ell |I_{n, j, \ell}|\geq |I_{n, j}| e^{-c \xi^2 L} .
$

\medskip
\noindent
{\em Proof of Step 1.\/} 
Set $\Delta S_{n,L}:=S_{nL}-S_{(n-1)L}, $ $\Delta V_{n,L}:= \mathrm{Var}(\Delta S_{n,L}),$
$\Delta z_{n,L}=z_{nL}-z_{(n-1)L}$,  $\Delta \PP_{n,L}=\PP_{nL}(\xi)-\PP_{(n-1)L} (\xi)$,  $\Delta\mathcal F_{n_i,L_i}(t):=(1/\Delta V_{n_i,L_i})\log\E(e^{t\Delta S_{n_i,L_i}})$. 

Suppose $\pi(I_{n,j})$ is a maximal $E^u_0$-segment in $\<\dot{x}_0,\ldots,x_{nL}\>$, then 
$\DS \sum_{j=1}^{nL}f_j(x_j)=z_{nL}$, and  each 
 $\pi(I_{n,j,\ell})$ is a maximal $E^u_{0}$-segment in 
some   cylinder  $\<\dot{x}_0,\ldots,x_{(n+1)L}\>$ such that 
$\DS \sum_{j=1}^{(n+1)L}f_j(x_j)=z_{(n+1)L}$. By the Markov property, 
$$
\frac{1}{|I_{n, j}|}\sum_\ell |I_{n, j, \ell}|=\Prob[\Delta S_{n,L}=\Delta z_{n,L}|X_{nL}=x_{nL}].
$$
By Lemma \ref{l.LLT-LD-Blocks}, for some $c>0$, and $L=L(\xi)$ large,   for all $n$ sufficiently large, 
\begin{equation}\label{e.snl}
\Prob[\Delta S_{n,L}=\Delta z_{n,L}|X_{nL}=x_{nL}]>\exp[-cL\xi^2].
\end{equation} 
Step 1 follows.
Henceforth we fix $L=L(\xi)$ as in Step 1.

\medskip
There is a cylinder $\<\dot{x}_0,\ldots,x_{(n+1)L}\>$ such that  $|I_{n,j,\ell}|=\ell^u(\<\dot{x}_0,\ldots,x_{(n+1)L}\>)$ and $|I_{n,j}|=
\ell^u(\<\dot{x}_0,\ldots,x_{nL}\>)$. Looking at Lemma \ref{Lemma-Cylinder-Area}, we see that if $L$ is sufficiently large, then  there are constants $\rho =\rho_L$ and $\eps>0$ such that 
\begin{equation}\label{e.rho-eps}
\DS \rho\leq \frac{|I_{n,j,\ell}|}{|I_{n,j}|}\leq e^{-\eps L}.
\end{equation}

\medskip
\noindent
{\em Step 2.\/} For some constant $C_0>0$, for every $n\geq n_0$,  $I_{n,j}\in\mathcal S$, and $s>c\xi^2/\eps$, 
$$
\nu(I_{n,j})\leq C_0 |I_{n,j}|^{1-s}.
$$

\medskip
\noindent
{\em Proof of Step 2.\/}  Fix $n_0$ as in Step 1.
Since there are only finitely many intervals $I_{n_0,j}$, there is a constant $C_0$ such that $\nu(I_{n_0,j})\leq C_0 |I_{n_0,j}|^{1-s}$ for all $j$ and $s>c\xi^2/\eps$. 
We claim that $\nu(I_{n,j})\leq C_0 |I_{n,j}|^{1-s}$ for all $n\geq n_0$, $j$ and $s>c\xi^2/\eps$. 

For $n=n_0$ this follows from the definition of $C_0$. Assume by induction that  the statement holds for $n$. Every $I_{n+1,i}$ equals $I_{n,j,\ell}$ for some $j,\ell$. 
By Step 1 and the definition of $\nu$,
$\DS
\frac{\nu(I_{n,j,\ell})}{\nu(I_{n,j})}\leq \frac{|I_{n,j,\ell}|}{|I_{n,j}|} e^{c\xi^2 L}.
$
If $s\geq c\xi^2/\eps$,  then $\DS e^{c\xi^2 L}\leq e^{s\eps L}\leq \biggl(\frac{|I_{n,j,\ell}|}{|I_{n,j}|}\biggr)^{-s}$, whence 
$\DS
\frac{\nu(I_{n,j,\ell})}{\nu(I_{n,j})}\leq \frac{|I_{n,j,\ell}|^{1-s}}{|I_{n,j}|^{1-s}}
$, and therefore 
$$
\nu(I_{n+1,i})=\nu(I_{n,j,\ell})\leq \frac{\nu(I_{n,j})}{|I_{n,j}|^{1-s}}\times |I_{n,j,\ell}|^{1-s}\leq C_0|I_{n,j,\ell}|^{1-s}=C_0|I_{n+1,i}|^{1-s}.
$$

\medskip
\noindent
{\em Step 3.\/} For some $\delta_0>0$, for every linear segment $I\subset W$, if $|I|<\delta_0$ and $s>c\xi^2/\eps$, then $\nu(I)\leq C_0 |I|^{1-s}$.

\medskip
\noindent
{\em Proof of Step 3.\/} First note that $\nu$ is non-atomic, because by Step 2, the measure of a set with diameter $r$ tends to zero as $r\to 0$.

Let $\delta_0:=\min\limits_{j} |I_{n_0,j}|$, and suppose $|I|<\delta_0$. It is convenient to think of $I$ and $W$ as of intervals. 
Let $J$ be a closed sub-interval of $\ov{I}$ such that both endpoints of $J$ are inside $\ov{K(\eps)}$, and $\ov{I}\cap \ov{K(\eps)}=\ov{J}\cap \ov{K(\eps)}$. 

Let $I_{n_0,j_1},I_{n_0,j_2}$ be the one or two intervals which contain the endpoints of $J$ in their closure.  The gap between $I_{n_0,j_1}$ and $I_{n_0,j_2}$ is less than $|J|$, whence less than $\delta_0$. Necessarily,  the interior of the gap  does not intersect any other $I_{n_0,i}$, and is  disjoint from $K_{n_0}(\xi)$. Thus 
$\ov{J}\cap \ov{K(\xi)}\subset \ov{I_{n_0,j_1}}\cup \ov{I_{n_0,j_2}}$. 

\medskip
\noindent
{\em Case 1:\/} Neither of $\ov{I_{n_0,j_1}}$, $\ov{I_{n_0,j_2}}$ covers $J$. 
In this case, each of $\ov{I_{n_0,j_1}}$, $\ov{I_{n_0,j_2}}$ has an endpoint inside the interior of $J$. 

Let $N\geq n_0$ be the maximal integer for which there exist intervals $I_{N,j_1},I_{N,j_2}$ such that $\ov{J}\cap \ov{I_{N,j_i}}=\ov{J}\cap \ov{I_{n_0,j_i}}$ ($i=1,2$). Clearly,  $I_{N,j_1},I_{N,j_2}$ both contain some endpoint of $J$ in their closures, and at least one of their endpoints is inside the interior of $J$. Also  
$
\ov{J}\cap \ov{K(\xi)}\subset  \ov{I_{N,j_1}}\cup \ov{I_{N,j_2}}.
$
Without loss of generality, $\nu(\ov{J}\cap \ov{I_{N,j_1}})\geq \nu(\ov{J}\cap\ov{I_{N,j_2}})$. Then, 
$
\nu(I)=\nu(\ov{J}\cap \ov{K(\xi)})\leq 2\nu(\ov{J}\cap \ov{I_{N,j_1}}). 
$

Let $N'$ be the maximal integer $N'\geq N$ for which there exists  $I_{N',j}$ such that $\ov{I_{N',j}}$ contains an endpoint of $J$, has an endpoint in the interior of $J$, and 
$
\nu(\ov{J})\leq 2\nu(\ov{J}\cap \ov{I_{N',j}}).
$
By construction, $\mathrm{int}(I_{N',j})\cap K(\eps)$ is a finite disjoint union of  intervals $\mathrm{int}(I_{N'+1,j'})$.  By the  maximality of $N'$, the one whose closure contains an endpoint of $J$ does not cover $J\cap K(\eps)$. Necessarily, another one, say $I_{N'+1,j'}$ is completely contained in $\ov{J}$. 
By \eqref{e.rho-eps}, 
 $|I_{N'+1,j'}|\geq \rho |I_{N',j}|$, whence $|J|\geq \rho |I_{N',j}|$. 

In summary, 
$
\nu(I)=\nu(\ov{J}\cap \ov{K(\xi)})\leq 2\nu(\ov{I_{N',j}})\leq 2C_0|I_{N',j}|^{1-s}
\leq 2C_0(\rho^{-1}|J|)^{1-s}
$. Since $J\subset \ov{I}$, $\nu(I)\leq 2C_0\rho^{-(1-s)}|I|^{1-s}$. 

\medskip
\noindent
{\em Case 2.\/} $J\subset \ov{I_{n_0,j}}$.

Let $n\geq n_0$ be the maximal integer such that  $J\subset \ov{I_{n,j}}$ for some $j$. Then there are $\ell_1,\ell_2$ such that $\ov{I_{n,j,\ell_i}}$ cover some endpoint of $J$, and their other endpoint is inside the interior of $J$. 

If $\ov{J}\subset \ov{I_{n,j,\ell_1}}\cup  \ov{I_{n,j,\ell_2}}$, then we continue as in Case 1.
Otherwise, there is some $\ell_3$ such that $I_{n,j,\ell_3}\subset J$, and $\rho|I_{n,j}|\leq |I_{n,j,\ell_3}|\leq |J|$  by \eqref{e.rho-eps}. Recalling that $J\subset \ov{I}_{n,j}$ we obtain that for all $s>c\xi^2/\eps$
$$
\nu(I)=\nu(J)\leq \nu(I_{n,j})\leq C_0|I_{n,j}|^{1-s}\leq C_0 \rho^{-(1-s)}|J|^{1-s}\leq C_0 \rho^{-(1-s)}|I|^{1-s}.
$$
This completes the proof that \eqref{e.mass-dist} holds for each $s>c\xi^2/\eps$. 

By Frostman's lemma, $\HD(W\cap K(\xi))\geq 1-c\xi^2/\eps$.
\end{proof}

\smallskip
\noindent
{\bf Proof of Theorem \ref{CrHDGen}.}
First we cover $\St^\ast$ be a countable collection of  $(E^s_0,E^u_0)$-parallelograms $Q_i$.  Then $\HD(\mathsf{Gen}(\xi))=\sup_i \HD(\mathsf{Gen}(\xi)\cap Q_i)$, and it is sufficient to show that there are constants $c_1,c_2, \xi_0$ such that for every $i$ and $|\xi|<\xi_0$,  
$$
1-c_1\xi^2\leq \HD(\mathsf{Gen}(\xi)\cap Q_i)\leq 1-c_2\xi^2.
$$

We fix $i$, set $Q:=Q_i$, and denote the $E^u_0$-side of $Q$ by $W$. 

Let $\vartheta_{Q}:Q\to W$ be the projection along the $E^s_0$-direction.
If $\vartheta_{Q}(\tomega_1)=\vartheta_Q(\tomega_2)$, 
then the linear segment $J_{\tomega_1,\tomega_2}$ from $\tomega_1$ to $\tomega_2$  is parallel to $E^s_0$, therefore 
$$
\dist((\tpsi_{n-1}\circ\cdots\circ\tpsi_0)(\tomega_1),(\tpsi_{n-1}\circ\cdots\circ\tpsi_0)(\tomega_2))\leq \diam[(\tpsi_{n-1}\circ\cdots\circ\tpsi_0)(J_{\tomega_1,\tomega_2})]\xrightarrow[n\to\infty]{}0.
$$
It follows that $|z_{nN_0}(\tomega_1)-z_{nN_0}(\tomega_2)|=O(1)$. 
By Lemma \ref{l.solution-xi} and Proposition \ref{p.Pressure-and-Geometric-Pressure}, for any sequence of integers $N_k\uparrow\infty$, if 
$
\exists \eta_{N_k}\to\eta$ s.t. $\mathcal P_{N_k}'(\eta_{n_k})=z_{N_k}(\tomega_1)$, then $\exists \eta_{N_k}'\to\eta$ s.t.  $\mathcal P_{N_k}'(\eta_{N_k}')=z_{N_k}(\tomega_2)$, and vice versa.
It follows that  $\tomega_1$  is $\mu_\xi$-generic iff $\tomega_2$ is $\mu_\xi$-generic.
Thus, $\mathsf{Gen}(\xi)\cap Q=\vartheta_Q^{-1}(W\cap\mathsf{Gen}(\xi))$.  

Since $\vartheta_Q^{-1}(W\cap\mathsf{Gen}(\xi))$  is a bi-Lipschitz image of $(W\cap\mathsf{Gen}(\xi))\times [0,1]$
$$
\HD\big[\mathsf{Gen}(\xi)\cap Q\big]=\HD\big[(W\cap\mathsf{Gen}(\xi))\times [0,1]\big]=\HD\big[W\cap\mathsf{Gen}(\xi)\big]+1.
$$
Now the theorem follows directly from Lemmas \ref{LmHDUpper} and \ref{LmHDLower}.
\qed

\medskip
\noindent
{\bf Proof of Theorem \ref{CrHDExcept}.} Let $\mathsf{Gen}$ denote the set of $\tomega$ which are generic for some locally finite ergodic invariant measure. The theorem states that $\mathsf{NonGen}:=\St^\ast\setminus\mathsf{Gen}$ is a residual set, with  zero Lebesgue measure and full Hausdorff dimension. 

\medskip
\noindent
{\bf Residual:} Recall the definition of $\ov{S}_n(\tomega)$ from the proof of Lemma \ref{LmHDUpper}.  
Fix $\xi>0$ and let
$
E_{n_0}(\xi)=\{\tomega: \exists n\geq n_0 \text{ such that }\ov{S}_n(\tomega)=[\PP_n'(\xi)]\}.
$

 $\DS \mathsf{NonGen}\supset \bigcap_{n_0} \mathrm{int}\left[ E_{n_0}(\xi/2) \cap E_{n_0}(\xi) \right]$, because as shown in the proof of Lemma \ref{LmHDUpper},  if $\tomega$ is in the right-hand-side, then $\exists n_k, n_k'\uparrow\infty$ and $\xi_{n_k}\to\xi$, $\xi_{n_k'}\to\xi/2$ such that $\PP_{n_k}'(\xi_{n_k})=z_{n_kN_0}(\tomega)+O(1)$ and $\PP_{n_k'}'(\xi_{n_k}')=z_{n_k'N_0}(\tomega)+O(1)$, and this implies that $\tomega$ has at least two different historical measures, $\mu_{\xi}$ and $\mu_{\xi/2}$.

The set $\mathrm{int}\left[ E_{n_0}(\xi/2) \cap E_{n_0}(\xi) \right]$ is clearly open. 

We claim that it is also dense. It  
suffices to show that every cylinder $\<\dot{a}_0,\ldots,a_n\>$ has an extension  $\<\dot{a}_0,\ldots,a_{n+m+\ell}\>$ where 
$$\DS \sum_{j=1}^{n+m} f_j(a_i)=[\PP_{n+m}'(\xi/2)]\text{ and }
\DS \sum_{j=1}^{n+m+\ell} f_j(a_i)=[\PP_{n+m+\ell}'(\xi)].$$

Let $S_i^j:=f_i(X_i)+\cdots+ f_j(X_j)$. By the local limit theorem {for large deviations (Theorem \ref{ThLD})}, 
$$
\Prob\left[ S_n^{n+m}=[\PP_{n+m}'(\xi/2)]-\sum_{j=1}^{n-1} f_j(a_j) \bigg|X_n=a_n\right]\neq 0\text{ for all $m$ large enough}. 
$$
Fix $m$ like that. The event on the right is a union of non-empty cylinders $\<\dot{a}_0,\ldots,a_{n+m}\>$ such that 
$\DS \sum_{j=1}^{n+m} f_j(a_j)=[\PP_{n+m}'(\xi/2)]$. 
Fix one of these cylinders, then again, 
$$
\Prob\left[ S_{n+m}^{n+m+\ell}=[\PP_{n+m+\ell}'(\xi)]-\sum_{j=1}^{n+m-1} f_j(a_j) \bigg|X_{n+m}=a_{n+m}\right]\neq 0\text{ for all $\ell$ large enough}.  
$$
Fix $m$ like that, then the event on the right is a union of cylinders  $\<\dot{a}_0,\ldots,a_{n+m+\ell}\>\neq \emptyset$ such that  $\DS \sum_{j=1}^{n+m+\ell} f_j(a_i)=[\PP_{n+m}'(\xi)]$ and 
$\DS \sum_{j=1}^{n+m} f_j(a_i)=[\PP_{n+m}'(\xi/2)]$.

In summary, $\mathsf{NonGen}$ contains a dense $G_\delta$ set. It follows that it is residual.

\medskip
\noindent
{\bf Zero Lebesgue Measure:} This is obvious since $\mu_0$ is proportional to the Lebesgue measure, and $\mu_0$-a.e. $\tomega$ is $\mu_0$-generic, by the ratio ergodic theorem.

\medskip
\noindent
{\bf Full Hausdorff Dimension:} 
Fix $C, \hat n_0>0$ such that 
$$|\PP_n'(0)|=|\E(S_n)|\leq C\text{ and  }C^{-1}n\leq \PP_n''\leq Cn\text{ on $[-1,1]$ for all $n>\hat n_0$,}
$$ see \eqref{Pres1stDer}, \eqref{Pres2dDer} and Proposition \ref{Prop-Var-Growth}. 

By Lemma \ref{l.delta-pressure-is-uni-conv}, there are $c, \hat n_1>0$  such that for all $k,n\geq \hat n_1$ and $|\xi|\leq 1$,
$$
c^{-1}n\leq \mathrm{Var}(S_{k+n}-S_k)\leq c n\text{ and } c^{-1}n\leq\PP_{k+n}''(\xi)-\PP_k''(\xi)\leq c n.
$$


\noindent
Let $K:=(C+c+1)(c+1)$. Fix  $0<\xi<1$ and choose an integer 
$$
L>\max\left\{\hat n_0, \hat n_1, \frac{2(C+1)K}{\xi}, \frac{2Cc}{1-\xi} \right\}.
$$
Let $n_i:=4^i$, and
define $\xi_{n_i L}:=\frac{\xi}{2K}$ for $i$ odd, $\xi_{n_i L}:=\frac{\xi}{K}$ for $i$ even, and extend by linear interpolation
$$
\DS \xi_{n}:=\frac{n-n_iL}{n_{i+1}L-n_iL}\xi_{n_{i+1}L}+\frac{n_{i+1}L-n}{n_{i+1}L-n_iL}\xi_{n_{i}L}\ \ (n_i L\leq n\leq n_{i+1}L).
$$
Obviously, $|\xi_n|\leq \frac{1}{K}<1$ for all $n$. We claim that if $n>4$ then 
$$
|\xi_{nL}-\xi_{(n-1)L}|<\frac{\xi}{Kn}.
$$
To see this, let  $i$ be the integer such that $n-1\in [n_i,n_{i+1}-1]$. Then  $(n-1)L, nL\in [n_iL,n_{i+1}L]$, and therefore   
$
|\xi_{nL}-\xi_{(n-1)L}|=\frac{|\xi_{n_{i+1}L}-\xi_{n_{i}L}|}{n_{i+1}-n_i}=\frac{\xi/(2K)}{(3/4)n_{i+1}}<\xi/(Kn).
$

Let $z_n:=[\PP_n'(\xi_n)]$. Then 
\begin{enumerate}[(a)]
\item  $z_{n_i L}:=[\PP_{n_i L}'(\xi)]$ for $i$ even;
\item  $z_{n_i L}:=[\PP_{n_i L}'(\xi/2)]$ for $i$ odd.
\end{enumerate}
In addition, we claim that 
\begin{enumerate}[(a)]
\item[(c)]   $ |z_{nL}-z_{(n-1)L}|\leq \xi \mathrm{Var}(S_{nL}-S_{(n-1)L})$ for all $n$;
\item[(d)]  $\exists\eta_{nL}\in [-1,1]$ such that  $z_{nL}-z_{(n-1)L}=\PP_{nL}'(\eta_{nL})-\PP_{(n-1)L}'(\eta_{nL})$. 
\end{enumerate}
Proof of  (c): $|z_n-\PP_n'(\xi_n)|\leq 1$, therefore
\begin{align*}
&|z_{nL}-z_{(n-1)L}|\leq \bigl|\PP_{nL}'(\xi_{nL})-\PP_{(n-1)L}'(\xi_{(n-1)L})\bigr|+2\\
&\leq \bigl|\PP_{nL}'(\xi_{nL})-\PP_{(n-1)L}'(\xi_{nL})\bigr|+\bigl|\PP_{(n-1)L}'(\xi_{nL})-\PP_{(n-1)L}'(\xi_{(n-1)L})\bigr|+2\\
&\leq \left(\bigl|\PP_{nL}'(\eta)-\PP_{(n-1)L}'(\eta)\bigr]^{\eta=\xi_{nL}}_{\eta=0}\bigr|\!+\!|\PP_{nL}'(0)-\PP_{(n-1)L}'(0)|\right)\\
&\hspace{7cm}+\!\bigl|\PP_{(n-1)L}'(\xi_{nL})\!-
\!\PP_{(n-1)L}'(\xi_{(n-1)L})\bigr|\!+\!2\\
&\leq \max_{[-1,1]}(\PP_{nL}''-\PP_{(n-1)L}'')\cdot |\xi_{nL}|+2C+(\max_{[-1,1]}\PP_{(n-1)L}'')\cdot|\xi_{nL}-\xi_{(n-1)L}|+2\\
&\leq c L(\xi/K)+2C+C(n-1)L\cdot \xi/(Kn)+2<L\xi\left(\frac{c+C}{K}+\frac{2(C+1)K}{K L\xi}\right)\\
&<c^{-1} L\xi \leq \xi\Var(S_{nL}-S_{(n-1)L})\ \ \text{ by choice of $L$, $K$ and  $c$}.
\end{align*}
Proof of (d): By the previous calculations, $|z_{nL}-z_{(n-1)L}|<c^{-1}\xi L$. At the same time, 
\begin{align*}
\PP_{nL}'(1)-\PP_{(n-1)L}'(1)&\geq \PP_{nL}'(0)-\PP_{(n-1)L}'(0)+\min_{[-1,1]}(\PP_{nL}''-\PP_{(n-1)L}'')\cdot(1-0)\\
&\geq -2C+c^{-1}L=c^{-1}L\xi +(c^{-1}L(1-\xi)-2C)\geq c^{-1}L\xi.\\
\PP_{nL}'(-1)-\PP_{(n-1)L}'(-1)&\leq \PP_{nL}'(0)-\PP_{(n-1)L}'(0)+\min_{[-1,1]}(\PP_{nL}''-\PP_{(n-1)L}'')\cdot(-1-0)\\
&\leq 2C-c^{-1}L=-c^{-1}L\xi -(c^{-1}L(1-\xi)-2C)\leq -c^{-1}L\xi.
\end{align*}
So $z_{nL}-z_{(n-1)L}$ lies between $\PP_{nL}'(-1)-\PP_{(n-1)L}'(-1)$ and $\PP_{nL}'(1)-\PP_{(n-1)L}'(1)$. Part (d) follows from the intermediate value theorem.

Let $K(\xi)=\{\tomega: \ov{S}_{nL}(\tomega)=z_{nL}$ for all $n\}$. Every $\tomega$ in $K$ has at least two different historical measures, and is therefore in $\mathsf{NonGen}$. 
Property (d) allows us to push through the  proof of Lemmas \ref{l.LLT-LD-Blocks} and \ref{LmHDLower}, and  deduce  that $\HD(K(\xi))\geq 1-O(\xi^2)$.
Thus $\HD(\mathsf{NonGen})\geq \HD(K(\xi))+1\geq 
2-O(\xi^2)\xrightarrow[\xi\to 0]{}2$.\qed

\part{Appendices}

\appendix

\section{Even Continued Fraction Expansions}
\label{AppEven}
Let $[x_0;x_1,\ldots,x_k]:=x_0+1/(x_1+1/(x_2+1/(\cdots+1/x_k)))$, whenever the expression makes sense.  
Given $m_0\in\Z$ and  $m_i\in\Z\setminus\{0\}$, let \index{Even continued fraction expansion}\index{Even CFE|see {Even continued fraction expansion}}
$$
[2m_0;2m_1,2m_2,\ldots]:=\lim_{k\to\infty} [2m_0;2m_1,\ldots,2m_k]
$$
In this appendix we will show that the limit exists, and belongs to $2m_0+[-1,1]$, and that every irrational $\alpha\in (-1,1)$ admits a unique  expansion
\begin{equation}\label{e.ecfe}
\alpha=[0;2m_1,2m_2,\ldots ],\text{ where }m_i\in\Z\setminus\{0\}.
\end{equation}
We call \eqref{e.ecfe} the {\em even continued fraction expansion (even CFE)} of $\alpha$. \index{Even continued fraction expansion}
We call  $2m_i$ the  partial quotients of the even CFE of $\alpha$. Note that all the partial quotients of the  even CFE must be even, but   some may be negative.  The even CFE has been considered before by Smith \cite{Smith},  Schweiger \cite{Schweiger}, and Kraaikamp \& Lopes \cite{Kraaikamp}.

\smallskip

We begin with some preparations. {\em M\"obius transformations} shall be denoted by \index{M\"obius transformations}  ${\tiny \left(
\begin{array}{cc}
a & b \\
c & d \\
\end{array}
\right)}\cdot z=\frac{az+b}{cz+d}$. Recall that the composition of two M\"obius transformations transforms to usual matrix product. Given $a\in 2\Z\setminus\{0\}$, define  
$$ v_a:[-1,1]\to \left[\frac{1}{a+1},\frac{1}{a-1}\right]\ , v_a(t)={\tiny \left(
\begin{array}{cc}
0 & 1 \\
1 & a \\
\end{array}
\right)}\cdot t=\frac{1}{a+t}. $$ 
\begin{lem}[Contraction]\label{l.v-contracts}
For every $a,b\in 2\Z\setminus\{0\}$, $|(v_a\circ v_b)'|\leq 1$ on $[-1,1]$. If $(a,b)\neq (\pm 2,\mp2)$, then $|(v_a\circ v_b)'|\leq \frac{1}{9}$ on $[-1,1]$.
\end{lem}
\begin{proof}
$(v_a\circ v_b)(t)={\tiny \left(
\begin{array}{cc}
0 & 1 \\
1 & a
\end{array}
\right)\left(\begin{array}{cc}
0 & 1 \\
1 & b
\end{array}
\right)\cdot t=\left(\begin{array}{cc}
1 & b \\
a & ab+1
\end{array}
\right)}\cdot t=\frac{t+b}{at+ab+1}$, so
$$
(v_a\circ v_b)'(t)=\frac{1}{(at+ab+1)^2}.
$$
Suppose $|t|\leq 1$.

\noindent
$\bullet$ 
If $(a,b)=(2,-2)$, then
$
\DS |(v_a\circ v_b)'(t)|{=}\frac{1}{(2t-3)^2}\leq 1
$.
 
\noindent
$\bullet$  
 If $(a,b)=(-2,2)$, then
$
\DS |(v_a\circ v_b)'(t)|{=} \frac{1}{(-2t-3)^2}\leq 1
$.

\noindent
$\bullet$ 
  If $(a,b)=\pm (2,2)$ then
$
\DS |(v_a\circ v_b)'(t)|\leq \frac{1}{(5-2|t|)^2}\leq \frac{1}{9}.$

\noindent
$\bullet$ 
In all other cases $|ab|\geq 8$ and $|ab|>|at|+1$, so \\
$\DS
|(v_a\circ v_b)'(t)|\leq \frac{1}{(|ab|-|a|-1)^2}\leq \frac{1}{(|ab|-\frac{1}{2}|ab|-1)^2}\leq \frac{1}{(4-1)^2}=\frac{1}{9}.
$
\end{proof}

\begin{lem}[Convergence]\label{l.limit-cfe-exists}
The limit \eqref{e.ecfe} exists for all $m_i\in 2\Z\setminus\{0\}$, and belongs to $[-1,1]$.
\end{lem}
\begin{proof}
Suppose the sequence of $2m_i$ does not terminate in an infinite sequence of the form $(2,-2,2,-2,\cdots)$, then
the intervals $(v_{2m_1}\circ\cdots\circ v_{2m_N})[-1,1]$ form a decreasing sequence of intervals with diameter tending to zero. So
$$
\lim_{n\to\infty}[0;2m_1,2m_2,\ldots,2m_n]=\lim_{n\to\infty}(v_{2m_1}\circ\cdots\circ v_{2m_n})(0)\text{ exists, and belongs to }[-1,1].
$$

Next, suppose  $2m_i$ does terminate in an infinite sequence of the form $(2,-2,2,-2,\ldots)$.
Assume first that  $(2m_1,2m_2,\ldots)$ equals $(2,-2,2,-2,\ldots)$. For $i$ odd,
$$
(v_{2m_i}\circ v_{2m_{i+1}})(x)=\frac{x-2}{2x-3}=\biggl(\begin{array}{ll}
1 & -2\\
2 & -3
\end{array}\biggr)\cdot x.
$$
The is a parabolic M\"obius transformation, therefore for all $\xi\in\R$, $(v_{2m_i}\circ v_{2m_{i+1}})^N(\xi)$ converges to the unique fixed point in $\R\cup\{\pm\infty\}$, equal to $1$, uniformly on compact subsets of 
$[-1,1]$. Let $\xi_n:=0$ for $n$ even and $\xi_n:=v_2(0)=\frac{1}{2}$ for $n$ odd, then 
$$
\lim_{n\to\infty}[0;2m_1,2m_2,\ldots,2m_{2n}]=\lim_{n\to\infty}(v_{2m_1}\circ\cdots\circ v_{2m_n})(0)=\lim_{n\to\infty}\biggl(\begin{array}{ll}
1 & -2\\
2 & -3
\end{array}\biggr)^n\cdot \xi_n=1.
$$
The case when $(m_1,m_2,\ldots)$ terminates with $(2,-2,2,-2,\cdots)$ easily follows.
\end{proof}

\begin{lem}[Existence and Uniqueness of the Even CFE]\label{l.f-expansion}
Every $\alpha\in (-1,1)\setminus\Q$ admits a unique sequence of $m_i\in \Z\setminus\{0\}$ such that
$\alpha=[0;2m_1,2m_2,\ldots]$.\index{Even continued fraction expansion!existence and uniqueness}
\end{lem}
\begin{proof}[Proof (Schweiger \cite{Schweiger-Book})]
By the previous lemma,  for every sequence $n_i\in 2\Z\setminus\{0\}$,
$
-1\leq [0;n_1,n_2,\cdots]\leq 1.
$
The maximum is obtained at $t:=[0;2,-2,2,-2,\ldots]$ and the minimum is obtained at $s:=[0;-2,2,-2,\ldots]$. Since $t=[0;2,-2,t]$, $t=1$. Since $s=[0;-2,2,s]$, $s=-1$.

Recall the  {\em even Gauss map}  $E:(-1,1)\to (-1,1)$ from \eqref{EvenGauss}.
Suppose $\alpha\in (-1,1)\setminus\Q$. Set $r_0:=\alpha$, $m_0=0$,
$r_n:=E^n(\alpha)$, and let $m_n\in\Z\setminus\{0\}$ be the unique integer such that $\frac{1}{r_n}\in 2m_{n+1}+(-1,1).$
Starting from the identity $\DS r_{n+1}=\frac{1}{r_n}-2m_{n+1}$, we find that
\begin{align*}
r_n&=\frac{1}{2m_{n+1}+r_{n+1}}=\left[0;2m_{n+1},\frac{1}{r_{n+1}}\right]\\
r_{n-1}&=\frac{1}{2m_{n}+r_{n}}=\cfrac{1}{2m_{n}+\cfrac{1}{2m_{n+1}+r_{n+1}}}=\left[0;2m_{n},2m_{n+1},
\frac{1}{r_{n+1}}\right]\\
&\cdots\cdots\cdots\cdots\\
r_1&=\left[0;2m_2,2m_3,\cdots,2m_{n+1},\frac{1}{r_{n+1}}\right]\\
\alpha&=r_0=\frac{1}{2m_1+r_1}=[0;2m_1,2m_2,\ldots,2m_{n+1},\frac{1}{r_{n+1}}].
\end{align*}
 Looking at Lemma \ref{l.limit-cfe-exists} and its proof, we deduce that the last expression tends to $[0;2m_1,2m_2,\ldots]$. Therefore $\alpha=[0;2m_1,2m_2,\cdots]$.
This is the even CFE of $\alpha$.

Next we show that the even CFE is unique. Suppose $\alpha=[0;2m_1,2m_2,\ldots]$ and $\alpha$ is irrational. Then
$
\frac{1}{\alpha}=2m_1+[0;2m_2,2m_3,\ldots]\in 2m_1+(-1,1)
$
This determines $m_1$ uniquely, and we see that $E(\alpha)=\frac{1}{\alpha}-2m_1=[0;2m_2,2m_3,\ldots]$, so 
$
\frac{1}{E(\alpha)}=2m_2+[0;2m_3,2m_4,\ldots]\in 2m_2+(-1,1).
$
This determines $m_2$ uniquely, and shows that $E^2(\alpha)=[0;2m_3,2m_4,\ldots]$.   Continuing in this way we see that each $m_n$ is uniquely determined by the condition $\frac{1}{E^{n-1}(\alpha)}\in 2m_n+(-1,1)$.
\end{proof}

\begin{lem}[Bounded Type]\label{l.BT-bounded-pm-two}
An irrational number\index{Even continued fraction expansion!of numbers of bounded type}
 $\alpha$ is of bounded type iff the even continued fraction expansion of $\alpha$ is uniformly bounded, 
 and  the lengths of the sunwards $(2,-2,\ldots,\pm 2)$ appearing in the even CFE of $\alpha$ are bounded.  Moreover, if 
\begin{align}
B(\alpha):=& \sup\{\text{partial quotients of the regular CFE of $\alpha$}\}\notag\\
M_1(\alpha):=&\sup\{\text{absolute values of partial quotients of the even CFE of $\alpha$}\}\notag\\
M_2(\beta):=&\sup\{\text{length of the sub-words $(2,-2,\ldots,\pm 2)$ or $(-2,2,\ldots,\pm 2)$ }\notag\\
&\hspace{6cm}\text{ which appear in the even CFE of $\alpha$}\}\notag\\
M(\alpha):=&\max\{M_1(\alpha),M_2(\alpha)\}+1.
\end{align}
 Then $M(\alpha)\geq 3$ and 
 $B(\alpha)-3\leq M(\alpha)\leq B(\alpha)+3$.
 \end{lem}

\begin{proof} $M(\alpha)\geq 3$ because $M_1(\alpha)\geq 2$.

The qualitative part of the lemma asserting that $M(\alpha)$ is finite iff $B(\alpha)$ is finite  has short geometric proof, relying on the well-known relation between the regular (resp. even) CFEs of $\alpha$ and the dynamical properties of projection of a geodesic ray in $\H$ with endpoint $\alpha$ to $\mathrm{SL}(2,\Z)\setminus\H$ (resp. $\Gamma(2)\setminus\H$):
Let $\alpha:=[2m_{0};2m_{1},2m_{2},\ldots]$ be the even CFE of $\alpha$.  $M(\alpha)=\infty$ iff   any geodesic $\gamma$ on $\H$ which terminates at $\alpha$ projects to a geodesic $\Gamma(2)\gamma$ on $\Gamma(2)\setminus\H$ which enters at least one of the cusps of $\Gamma(2)\setminus\H$ arbitrarily deeply.
Since $\Gamma(2)\lhd \mathrm{SL}(2,\Z)$, the surface $\Gamma(2)\setminus\H$ is a  regular cover of $\mathrm{SL}(2,\Z)\setminus\H$. Necessarily, $M(\alpha)=\infty$ iff  $\mathrm{SL}(2,\Z)\gamma$ enters the cusp of the modular surface arbitrarily deeply. This happens iff  $B(\alpha)=\infty$.\index{Kraaikamp-Lopes theorem}

Next we give another proof, which yields the stronger quantitative assertion that $B(\alpha)-3\leq M(\alpha)\leq B(\alpha)+3$.  The proof relies on  a procedure due to Kraaikamp \cite{Kraaikamp}, as presented in \cite{Kraaikamp-Lopes}, which  transforms the regular CFE to the even CFE.\footnote{\cite{Kraaikamp-Lopes} use a different definition for the even CFE which results in a different (though equivalent) expansion. What follows is  translation to our language.}  

Kraaikamp's procedure is based on the following algebraic identities \cite{Kraaikamp,Kraaikamp-Lopes}:
	\begin{enumerate}[(S1)]
		\item $[A,1,B,z]=[A+1,-(B+1), -z]$
		\item $[A,B,z]=[A+1,-1,-(B-1),-z]$
		\item $[A,B,z]=[A-1,1,-(B+1),-z]$
	\end{enumerate}	
To verify these identities, we first check by direct calculation that
	$$
	\left(\begin{array}{cc}
		1 & a \\
		0 & 1
	\end{array}
	\right)\left(\begin{array}{cc}
		0 & \eps \\
		1 & 1
	\end{array}
	\right)\left(\begin{array}{cc}
		0 & 1 \\
		1 & b
	\end{array}
	\right)=\left(\begin{array}{cc}
		1 & a+\eps \\
		0 & 1
	\end{array}
	\right)\left(\begin{array}{cc}
		0 & -\eps \\
		1 & b+1
	\end{array}
	\right).
	$$
	Now suppose  $\eps=\pm 1$, and let us replace each matrix by the M\"obius transformation it defines. Applying the M\"obius transformations on the two sides of the equation to $1/z$, and noting that $\eps^{-1}=\eps$, we obtain
$		[a;\eps,\eps b,\eps z]=[a+\eps;-\eps(b+1),-\eps z].$

	Substituting $\eps=1$, $a=A$, $b=B$, gives (S1).
	Substituting $\eps=-1$, $a=A+1$, $b=B-1$, gives (S2).
	Substituting $\eps=1$, $a=A-1$, $b=-(B+1)$, $z=-w$ gives (S3).

	Suppose now that $\alpha$ has regular CFE
	$
	\alpha=[B_0;B_1,B_2,\ldots]
	$, i.e. $B_0\in\Z$ and $B_i\in\N$. If all partial quotients $B_i$ are even, this is also the even CFE, and the lemma is trivial. If not, choose the minimal $n\geq 0$ such that   $B_0,\ldots,B_{n-1}$ are even and $B_n$ is odd.

\smallskip
\noindent
{\em Case 1.\/}  $B_{n+1}=1$. Write using (S1)
	\begin{align*}
		\alpha&=[\underset{\text{even}}{\underbrace{B_0;B_1,\ldots,B_{n-1}}},\underset{\text{odd}}{\underbrace{B_n}},1,B_{n+2},z],\text{ where }z:=[B_{n+3};B_{n+4},\ldots]\\
		&=[\underset{\text{even}}{\underbrace{B_0;B_1,\ldots,B_{n-1}}},\underset{\text{even}}{\underbrace{B_n+1}},-(B_{n+2}+1),-z],
	\end{align*}
and we have gained one additional even digit.

\smallskip
\noindent
{\em Case 2.\/}  $B_{n+1}>1$. In this case we write
	\begin{align*}
		\alpha&=\big[B_0;B_1,\ldots,B_{n-1},
		\big[B_n;B_{n+1},z\big]\big],\text{ where }z:=[B_{n+2};B_{n+3},\ldots]\\
		&=\big[B_0;B_1,\ldots,B_{n-1},
		[B_n+1;-1,-(B_{n+1}-1),-z]\big], \text{ by (S2)}\\
		&=\big[B_0;B_1,\ldots,B_{n-1},
		B_n+1,[-1;-(B_{n+1}-1),-z]\big] \text{ trivially}\\
		&=\big[B_0;B_1,\ldots,B_{n-1},
		B_n+1,[-2;1,B_{n+1}-2,z]\big],\text{ by (S3)}\\
		&=\big[B_0;B_1,\ldots,B_{n-1},
		B_n+1,-2,[1,B_{n+1}-2,z]\big],\text{ trivially}\\
				&=\big[B_0;B_1,\ldots,B_{n-1},
		B_n+1,-2,[2,-1,-(B_{n+1}-3),-z]\big],\text{ by (S2)}\\
				&= \big[B_0;B_1,\ldots,B_{n-1},
		B_n+1,-2,2,-[1;(B_{n+1}-3),z]\big]\\
		& \hspace{5cm}\cdots\cdots\cdots  \text{ (repeated applications of last two steps)}\\
		&= \big[B_0;B_1,\ldots,B_{n-1},B_n+1,
		\underset{B_{n+1}-2}{\underbrace{-2,2,-2,\cdots,(-1)^{B_{n+1}}2}},(-1)^{B_{n+1}+1}[1;1,z]\big].
	\end{align*}
	Write $[1;1,z]=[1;1,B_{n+2},w]$ with $w=[B_{n+3};B_{n+4},\ldots]$. By (S1),  
	$$[1;1,B_{n+2},w]=[2;-(B_{n+2}+1),-w]=[2;-[B_{n+2}+1,w]].$$
	Substituting this in the above leads to
	\begin{align*}
		\alpha&= \big[B_0;B_1,\ldots,B_{n-1},B_n+1,
		\underset{B_{n+1}-1}{\underbrace{-2,2,-2,\cdots,(-1)^{B_{n+1}+1}2}},\\
		&\hspace{4.9cm}(-1)^{B_{n+1}}[B_{n+2}+1;
		B_{n+3},B_{n+4},\ldots]\big].
	\end{align*}
	This representation starts with at least $n+B_{n+1}>n$ even digits.
	
\smallskip
\noindent
	In both cases the initial CFE starts with $n$ even digits,  and the final CFE starts with at least $n+1$ even digits. The new expansion may contain a stretch of $(2,-2,\ldots,\pm 2)$ of  length bounded by
	$
	M:=\max\{B_i:i\geq 0\}.
	$	
	Finally, the two expansions agree with each other on the  first $n$ digits.
	
	Proceeding by induction, we obtain an infinite sequences $C_i$ and $w_i\in (-1,1)$ s.t.
	\begin{enumerate}[(1)]
		\item $\alpha=[C_0;C_1,\ldots,C_n,w_n]$,  and $|w_n|$ is uniformly bounded;
		\item  $C_i\leq M+1$, $C_i$ are all even, and if $i\neq 0$ , then $C_i\neq 0$;
		\item All sub-words  $(-2,2,-2,\ldots,\pm 2)$ of  $(C_0,C_1,C_2\ldots)$ have length less than $M$.
	\end{enumerate}
Looking at the proof of Lemma \ref{l.limit-cfe-exists}, it is easy to verify that
	$
	\alpha=[C_0;C_0,C_1,\ldots]
	$.  

Looking closer 	at the previous procedure we see that it either modifies a digit by  $+1$ or a sign or both at most  twice, or it replaces a digit $B$ by a word of the form $\pm (-2,2,\ldots,\pm 2)$ of length $\leq B
+2$. It follows that $B(\alpha)-3\leq M(\alpha)\leq B(\alpha)+3$.
\end{proof}

\begin{lem}[Compactness]\label{l.BT-compact}
For every $M>1$ there is a compact subset $K\subset (-1,1)$ as follows. Suppose $M(\beta)\leq M$ and the even CFE of
$\beta$ is $\beta=[0;2m_1,2m_2,\ldots]$. Then
$$ r_k:=[0;2m_k;2m_{k+1},\ldots]\in K\text{ for all }k.$$
\end{lem}
\begin{proof}
$\DS
\pm \frac{1}{r_k}\in \bigcup_{n=1}^{M} \underset{m\neq -1}{\bigcup_{0<|m|\leq M}}(\underset{n\text{ even (incl. zero)} }{\underbrace{v_{2}v_{-2}\cdots v_{-2}}})\circ (v_2 v_{2m})([-1,1]) $
$$\hspace{-1cm}\bigcup \bigcup_{n=1}^M \underset{m\neq 1}{\bigcup_{0<|m|\leq M}}(\underset{n\text{ odd } }{\underbrace{v_{2}v_{-2}\cdots v_{2}}})\circ (v_{-2} v_{2m})([-1,1])\bigcup \bigcup_{2\leq m\leq M} v_{2m}[-1,1].
$$
The first and second  sets take care of $1/r_k=[2;\ast,\ast,\ldots]$, the third set takes care of $1/r_k=[2m;\ast,\ast,\ldots]$ where $2m>2$, and the plus/minus in front of $1/r_k$ takes care of the remaining (symmetric) options.

Call the set on the right-hand-side $K'$.
We claim that $K'$ is a finite union of compact subsets of $(-1,1)$. In what follows, $\mathrm{Int}[c,d]$ is the closed interval with endpoints $c,d$ (even if $c>d$).
\begin{enumerate}[$\bullet$]
\item $v_2v_{2m}[-1,1]=v_2[\frac{1}{2m+1},\frac{1}{2m-1}]=\mathrm{Int}[\frac{2m+1}{4m+3},\frac{2m-1}{4m-1}]$. We know that $v_a:[-1,1]\hookrightarrow[-1,1]$ are continuous homeomorphisms onto their images.  If $m\neq -1$, then $\pm 1$ is not in the interval $\mathrm{Int}[\frac{2m+1}{4m+3},\frac{2m-1}{4m-3}]$. So
    $\DS
    v_2v_{2m}[-1,1]$ is a compact subset of $(-1,1).
    $
    
    Similarly, $v_2v_{-2}\cdots v_{-2}:[-1,1]\hookrightarrow[-1,1]$ is a homeomorphisms onto its image. As such, it must map a compact subset of $(-1,1)$ onto a compact subset of $(-1,1)$. This takes care of the first union of sets.
\item The second union of sets can be handled in the same way, using the observation that $(v_{-2}v_{2m})[-1,1]=\mathrm{Int}[\frac{2m+1}{-4m-1},\frac{2m-1}{-4m+3}]$ does not contain $\pm 1$, when $m\neq 1$.
\item The third union of sets obviously does not contain $\pm 1$, and is compact.
\end{enumerate}
Since $K'$ is a finite union of compact subsets of $(-1,1)$, it too is a compact subset of $(-1,1)$. The lemma follows with $K:= K'\cup(-K')$.
\end{proof}

\section{Denseness Bounds: Proof of Lemma \ref{LmDenseWs}}\label{appendix-denseness}
Let $\BT(M):=\{[a_0; a_1,a_2,\cdots]: a_0\in\Z, a_1,a_2,\ldots\in\N, \sup|a_i|< M\}$.  We say that a line on $\mathsf{St}_{0}$ has slope $r$, if it has slope $r$ in the horizontal rectangle model for $\mathsf{St}_{0}$ in Figure \ref{Figure-Staircase}(b). In this appendix we prove Lemma \ref{LmDenseWs}:\index{Denseness constant}

\medskip
\noindent
{\bf Lemma \ref{LmDenseWs}.} 
{\em  For every $M>1$, $0<\epsilon<1$ there exists $T$ as follows: If  $r\in\BT(M)$,  
   then every line with slope $r$ and length $T$ on $\St_0$ is $\epsilon$-dense in $\St_0$.}

\noindent
The proof relies of the following fact:
\begin{lem}\label{l.rho-BT}
 For each $M$ there is a constant $A=A(M)$ as follows: If $r>0$ and $r\in \BT(M)$, then  $\frac{1+r}{1-r}\in \BT(A)$. 
\end{lem}
\begin{proof}
Fix $r>0$ in $\BT(M)$, $0<c(M)<(M+2)^{-1}(M+1)^{-2},$ and  let  $\rho:=\frac{1+r}{1-r}$.

Since $r\in\BT(M)$,  for any 
 fraction $\frac{p}{q}$, $$\left|r-\frac{p}{q}\right|> \frac{c(M)}{q^2}$$ (see e.g.  the proof of  \cite[Thm 23]{Khintchine-Continued-Fractions}). In particular
$
|r-1|>c(M).
$

Let $\varphi(t):=\frac{1+t}{1-t}$. The inverse function of $\varphi(t)$ is  $\psi(s)=\frac{s-1}{s+1}$. Clearly, $\psi(t)$ is bounded away from $s=-1$.  Therefore, there exists $K=K(M)$ such that 
$$
|\psi(s)|>K\Rightarrow |s+1|<c(M)/2.
$$

\medskip
\noindent
{\em Step 1.\/} There exists a positive constant $c^\#(M)$ which only depends on $M$ as follows. 
Suppose $\frac{p}{q}$ is a 
fraction such that  (a) $|\frac{p}{q}+1|<\frac{M+3}{c(M)}$; (b) $|\psi(\frac{p}{q})|\leq K$; and \\(c) the interval with endpoints $r, \psi(\frac{p}{q})$ does not contain $1$.  
Then $|\rho-\frac{p}{q}|>c^\#(M)/q^2$. 

\noindent
\medskip
{\em Proof of Step 1.\/} Note that $\psi(p/q)=\frac{p-q}{p+q}$. 
By (c), $\varphi$ is smooth on the interval $I$ with endpoints $r,\psi(p/q)$. By the  mean value theorem,
\begin{align*}
&\left|\rho-\frac{p}{q}\right|=\left|\varphi(r)
-\varphi\left(\frac{p-q}{p+q}\right)\right|\geq \inf_{\xi\in I}{|\varphi'(\xi)|}\left|r-\frac{p-q}{p+q}\right|.
\end{align*}

By our assumptions on $r$,  $0<r<M+1$. By (b), $I\subseteq \{t:|t|<\max\{M+1, K\}\}$. Let 
$
\mu(M):=\inf\{|\varphi'(\xi)|:\; |\xi|\leq \max\{M+1, K\},\; \xi\neq 1\}
$ 
(this is a positive number because $\DS |\varphi'(\xi)|=2/(\xi-1)^2$). 
Thus, by the choice of $c(M)$ and the assumption on $r$,
\begin{align*}
& 
\left|\rho-\frac{p}{q}\right|\geq \mu(M)\cdot \left|r-\frac{p-q}{p+q}\right|
\geq \frac{\mu(M)c(M)}{(p+q)^2}=\frac{\mu(M)c(M)}{(1+\frac{p}{q})^2}\cdot \frac{1}{q^2}\geq \frac{\mu(M)c(M)}{c(M)^{-2}(M+3)^2}\cdot \frac{1}{q^2}
\end{align*}
 by (a).
 Now take $c^\#(M):=\mu(M)c(M)^3(M+3)^{-2}$.

\medskip
\noindent
{\em Step 2.\/} There exists a constant $c^\flat(M)$ which only depends on $M$ as follows.  If a 
fraction $p/q$ does not satisfy some of the conditions (a)--(c) above, then $|\rho-\frac{p}{q}|>\frac{c^\flat(M)}{q^2}$.

\medskip
\noindent
{\em Proof of Step 2.\/} 
If (a) fails, then $|\frac{p}{q}+1|\geq \frac{M+3}{c(M)}$. Since  $|\rho+1|=\frac{2}{|r-1|}<2/c(M)$, we have
$\left|
\rho-\frac{p}{q}\right|\geq \left||\rho+1|-|\tfrac{p}{q}+1|\right|>\frac{M+1}{c(M)}>1\geq \frac{1}{q^2}.
$

If (b) fails, then $|\psi(p/q)|>K$, and by the choice of $K$, $|\frac{p}{q}+1|<c(M)/2$. On the other hand, $|\rho+1|=\frac{2}{|r-1|}>\frac{2}{M+1}>2(M+1)(M+2)c(M)>c(M)$, whence
$\left|
\rho-\frac{p}{q}\right|=|(\rho+1)-(\frac{p}{q}+1)|>c(M)/2>(c(M)/2)/q^2
$.

Finally, if (c) fails and $1$ lies between $r$ and $\psi(p/q)$, then either $0<r<1<\psi(p/q)$ or $0<\psi(p/q)<1<r$ or $\psi(p/q)\leq 0<1<r$. In the first two cases, $\rho=\varphi(r)$ and $p/q=\varphi(\psi(p/q))$ have different signs, whence 
$$
\left|
\rho-\frac{p}{q}
\right|>|\rho|=\frac{1+r}{|1-r|}>\frac{1}{|r-1|}>\frac{1}{c(M)}>1\geq \frac{1}{q^2}.
$$
In the third case, since $\varphi(t)$ is increasing on $(-\infty,1)$ and on $(1,\infty)$, \\ $\rho=\varphi(r)<\varphi(M+1)=-\frac{M+2}{M}$ and $\frac{p}{q}=\varphi(\psi(p/q))>\lim\limits_{s\to-\infty}\varphi(s)=-1$, so 
$$
\left|
\rho-\frac{p}{q}
\right|>\left|1-\frac{M+2}{M}\right|=\frac{2}{M}\geq \frac{2M^{-1}}{q^2}.
$$
Step 2 follows with $c^\flat(M):=\min\{c(M)/2,1,2/M\}=c(M)/2$. 

\medskip
\noindent
{\em Proof of the Lemma.\/} Let $b=b(M):=\min\{c^\sharp(M), c^\flat(M)\}$, then $|\rho-\frac{p}{q}|\geq b/q^2$ for every 
fraction $p/q$. By the theory of continued fractions (see e.g. equation (34) in \cite{Khintchine-Continued-Fractions}), all the partial quotients with positive index of $\rho$ are bounded above by $1/b$. 
The zeroth partial quotient is bounded above by 
$$
|\rho|+1\leq  |\rho+1|+2\leq \frac{2}{|r-1|}+2\leq \frac{2}{c(M)}+2. 
$$
The lemma follows with $A:=\max\{b(M),2c(M)^{-1}+2\}$. 
\end{proof}

\begin{lem}
 For every $\epsilon, M$ there exists $T$ as follows. Suppose $|\rho|\in\BT(M)$. Then any line segment in $\R^2$ with slope $\rho$ and length $L>T$ projects to an  $\epsilon$-dense curve in $\T^2:=\R^2/\Z^2$.
\end{lem} 
\begin{proof}[Sketch of Proof] We use the fundamental domain $\T^2:=[0,1)^2$.  Applying one of the symmetries of $\T^2$, we may reduce the claim to the case when $0<\rho<1$. Then any line of length $L$ contains at least 
$
n\geq L(1+\rho^2)^{-1/2}-2
$
parallel line segments of length $\sqrt{1+\rho^2}$ and slope $\rho$ extending across $[0,1)^2$ from left to right, and starting at heights
\begin{equation}\label{e.set-of-heights}
\{x_0\}, \{x_0+\rho\},\ldots,\{x_0+n\rho\}.
\end{equation}

Clearly, if these points are $\epsilon$ dense in $\T$, then our line is $\epsilon$-dense in $\T^2$.
Thus the claim reduces to the following well-known fact in diophantine geometry: For all $\epsilon,M$ there exists $n=n(\epsilon,M)$ so that the set \eqref{e.set-of-heights} is $\epsilon$-dense in $\T$ whenever all the partial quotients of $\rho$ are less than $M$. (See  \cite[Prop 5.6]{Adiceam-Solomon-Weiss} for a deeper statement.)
\end{proof}

\smallskip
\noindent
{\bf Proof of Lemma \ref{LmDenseWs}.}
Reflecting the fundamental domain $[0,2]\times [0,1]$ across the axis $x=1$ we may assume without loss of generality that $r>0$. 
 Recall that $\mathsf{St}_0$ is isometric to torus $\R^2/\sqrt{2}\Z^2$. In this identification, lines with slope $r$ on $\mathsf{St}_0$ are mapped to lines with slope $\rho:=\frac{1+r}{1-r}$ on $\R^2/\sqrt{2}\Z^2$ (see Figure \ref{Figure-Staircase} (c)). By Lemma  \ref{l.rho-BT},  $\rho\in\BT(A)$ for some constant $A$ which only depends on $M$. Now apply the previous lemma. \qed


\section{Twists: Proof of Theorem \ref{ThDriftOfOrigin}}\label{Appendix-Twist}
Recall the number $N_0$ from \eqref{e.mat-prod-form}. In this appendix we provide the proof of: \index{Automorphisms of the infinite staircase!twist}\index{Twist}

\medskip
\noindent
{\bf Theorem \ref{ThDriftOfOrigin}.} {\em 
Let $\beta\in (-1,1)$ be an irrational with even CFE  $\beta:=[0;2m_1,2m_2,\ldots]$, and let $L_0$ denote the ray in direction $\beta\choose 1$,  which begins at  the singularity $\tomega_0$, in the middle of the top side of rectangle $\mathsf R_0$. Let $\frak{z}_n$ denote the $\Z$-coordinate of the horizontal rectangle which contains the beginning of  $(\wt{\psi}_{n-1}\circ\cdots\circ\wt{\psi}_0)(L_0)$. Then
$$
\frak{z}_n=\frac{1}{2}\sum_{k=1}^{nN_0} \bigl[\sgn(m_{2k})-\sgn(m_{2k-1})\bigr].
$$
}

\medskip
\noindent
{\bf Horizonal and Vertical Coordinates.}
Recall that the infinite staircase can be partitioned
into {\em horizontal rectangles} $\mathsf R_{\frak{z}}$ $(\frak{z}\in\Z)$, equal up to translation to $(0,2]\times (0,1]$. 
There is a {\em deck transformation} $D:\St\to\St$ which maps $\mathsf R_{\frak{z}}$ onto $\mathsf R_{\frak{z}+1}$, by translation. \index{Infinite staircase!horizontal and vertical coordinates}\index{Horizontal and vertical coordinates}

 Suppose $\tomega\in \mathsf R_\zc$, and $\tomega$ is not a singularity. Then $\tomega$ is uniquely determined by the vector
$
{\tiny \left(
                                         \begin{array}{c}
                                           x \\
                                           y \\
                                           \zc \\
                                         \end{array}
                                       \right)_H}
$
where $x$ is the distance of $\tomega$ from the left vertical side, $y$ is the distance of $\tomega$ from the bottom horizontal side, and $\zc$ is the horizontal $\Z$-coordinate. We call this vector the {\em horizontal coordinates} of $\tomega$.

The infinite staircase can also be divide into {\em vertical rectangles}, equal up to translation to $(0,1]\times(0,2]$. Exactly one of these rectangles has a bottom half inside horizontal rectangle zero, and we will call it {\em vertical rectangle zero}. Every other vertical rectangle is the image of this rectangle by $D^\zc$ for some unique $\zc$, which we call the {\em vertical $\Z$-coordinate} of this rectangle. Suppose $\tomega$ is a point in vertical rectangle $\xi$, and $\tomega$ is not a singularity. Then $\tomega$ is uniquely determined by the vector
$
{\tiny \left(
                                         \begin{array}{c}
                                           x \\
                                           y \\
                                           \zc \\
                                         \end{array}
                                       \right)_V}
$
where $x$ is the distance of $\tomega$ from the left vertical side, $y$ is the distance of $\tomega$ from the bottom horizontal side, and $\zc$ is the vertical $\Z$-coordinate. We call this vector the {\em vertical coordinates} of $\tomega$.

With these conventions, each vertical rectangle $\zc$ is covered by horizontal rectangles $\zc$ and $\zc+1$, and
$$
\left(
\begin{array}{c}
x \\
y \\
\zc \\
\end{array}
\right)_V=
\begin{cases}
{\tiny    \left(
\begin{array}{c}
x+1 \\
y \\
\zc \\
\end{array}
\right)_H}\vspace{0.1cm} & 0<y\leq 1,\\
{\tiny    \left(
\begin{array}{c}
x \\
y-1 \\
\zc+1 \\
\end{array}
\right)_H} & 1<y\leq 2.
\end{cases}
$$

\medskip
\noindent
{\bf Building Blocks.}
Given $\ell,m\in\Z$, let  $H_{-2m},V_{-2\ell}:\St\to\St$ be the unique the homogeneous automorphisms with zero drift and derivatives ${\tiny \left(
                                     \begin{array}{cc}
                                       1 & -2m \\
                                       0 & 1 \\
                                     \end{array}
                                   \right)
}
$ and ${\tiny \left(
                                     \begin{array}{cc}
                                       1 & 0 \\
                                       -2\ell & 1 \\
                                     \end{array}
                                   \right)
}
$. 
Write  $\{t\}_{(0,2]}:=\text{  unique $s\in (0,2]$ s.t. }t\in s+2\Z.$ Then these automorphisms are given explicitly by the formulas \index{Building blocks}
$$
H_{-2m}:\left(
    \begin{array}{c}
      x \\
      y \\
      \zc \\
    \end{array}
  \right)_H\mapsto \left(
    \begin{array}{c}
      \{x-2my\}_{(0,2]} \\
      y \\
      \zc \\
    \end{array}
  \right)_H \ ,\quad  V_{-2 \ell}:\left(
    \begin{array}{c}
      x \\
      y \\
      \zc \\
    \end{array}
  \right)_V\mapsto \left(
    \begin{array}{c}
      x \\
      \{y-{2 \ell x}\}_{(0,2]} \\
      \zc \\
    \end{array}
  \right)_V\!\!.
$$

\begin{lem}\label{l.bb-id}
$\wt{\psi}_{n-1}\circ\cdots\circ\wt{\psi}_0=(H_{-2m_{2nN_0}}\circ V_{-2m_{2nN_0-1}})\circ \cdots\circ (H_{-2m_2}\circ V_{-2m_1})$.
\end{lem}
\begin{proof}
If $\delta(\wt{\psi})$ denotes the drift defined in \eqref{Eq-Drift-Defn}, then $\delta(\wt{\psi}_1\circ \wt{\psi}_2)=\delta(\wt{\psi}_1)+\delta(\wt{\psi}_2)$, see \cite{Avila-Dolgopyat-Duryev-Sarig}.  Therefore, on both sides of the equation we have homogeneous automorphisms with zero drift. 
The derivatives  are the same: By \eqref{e.mat-prod-form}, the derivative of  $\wt{\psi}_{n-1}\circ\cdots\circ\wt{\psi}_0$ equals 
$
A_{nN_0}(\beta)
$, and by \eqref{e.A-k}, this is also the derivative of the right-hand-side.
By   Proposition \ref{Prop-Automorphisms}, a homogeneous automorphism with zero drift is  uniquely determined by its derivative. Therefore the two sides of the equation are equal.
\end{proof}

We call  $H_{-2m_{i+1}}\circ V_{-2m_{i}}$ {\em building blocks}. To calculate them explicitly, it useful to partition $\St$ into the sets
\begin{align*}
&H^i_{00}(\zc):=\left\{{\tiny
\left(
  \begin{array}{c}
    x \\
    y \\
    \zc \\
  \end{array}
\right)_H
}: 0<x\leq 1\ , \ \{y-2m_i x\}_{(0,2]}\in (0,1]
\right\}\\
&H^i_{01}(\zc):=\left\{{\tiny
\left(
  \begin{array}{c}
    x \\
    y \\
    \zc \\
  \end{array}
\right)_H
}:  0<x\leq 1\ , \  \{y-2m_i x\}_{(0,2]}\in (1,2]\right\}\\
&H^i_{10}(\zc):=\left\{{\tiny
\left(
  \begin{array}{c}
    x \\
    y \\
    \zc \\
  \end{array}
\right)_H
}: 1<x\leq 2\ , \ \{y-2m_i x\}_{(0,2]}\in (0,1]
\right\}\\
&H^i_{11}(\zc):=\left\{{\tiny
\left(
  \begin{array}{c}
    x \\
    y \\
    \zc \\
  \end{array}
\right)_H
}:  1<x\leq 2\ , \  \{y-2m_i x\}_{(0,2]}\in (1,2]\right\}\\
\end{align*}
Let $\DS H_{01}^i:=\biguplus_{\zc\in\Z}H^i_{01}(\zc)$, 
$\DS H_{11}^i:=\biguplus_{\zc\in\Z}H^i_{11}(\zc)$, and define the function 
$Z_i:\St^\ast\to\Z$,
\begin{equation}\label{e.Z-def}
Z_i:=1_{H^i_{01}}-1_{H_{11}^i}.
\end{equation}
\begin{lem}\label{l.composition-formula}
$H_{-2m_{i+1}}\circ V_{-2m_i}$ is given in coordinates by 
$$
\left(
\begin{array}{c}
x \\
y \\
\zc
\end{array}
\right)_{\!\!\!H}\!\!\!\!\mapsto
\left(
  \begin{array}{ccc}
    4m_i m_{i+1}+1 & -2m_{i+1} & 0 \\
    -2m_i & 1  & 0 \\
    0 & 0 & 1 \\
  \end{array}
\right)
\left(
\begin{array}{c}
x \\
y \\
\zc
\end{array}
\right)_{\!\!\!H}+\left(
\begin{array}{c}
Z_i \\
0 \\
-Z_i
\end{array}
\right)
\begin{array}{l}
\mod 2\Z\\
\mod\Z\\
\ \\
\end{array}
$$
\end{lem}
\begin{proof}
 $H^i_{00}(\zc)$ is inside the top half of vertical rectangle $\zc-1$, so on $H^i_{00}(\zc)$ we have
\begin{align*}
\left(
  \begin{array}{c}
    x \\
    y \\
    \zc \\
  \end{array}
\right)_H&\!\!\!=
\left(
  \begin{array}{c}
    x \\
    y+1 \\
    \zc-1 \\
  \end{array}
\right)_V \xmapsto[]{V_{-2m_i}}
\left(
  \begin{array}{c}
    x \\
    \{y+1-2m_i x\}_{(0,2]} \\
    \zc-1 \\
  \end{array}
\right)_V\\
&=\left(
  \begin{array}{c}
     x\\
     \{y-2m_i x\}_{(0,2]}\\
     \zc\\
  \end{array}
\right)_H,\text{ because }\{y+1-2m_i x\}_{(0,2]}\in (1,2]\\
&\xmapsto[]{H_{-2m_{i+1}}}
\left(
  \begin{array}{c}
     \left\{x-2m_{i+1}\{y-2m_i x\}_{(0,2]}\right\}_{(0,2]}\\
     \{y-2m_i x\}_{(0,2]}\\
     \zc\\
  \end{array}
\right)_H.
\end{align*}
Modulo $2\Z$, the $x$-coordinate equals to  $(4m_{i+1}m_i+1)x-2m_{i+1} y$. Thus on $H^i_{00}(\zc)$
\begin{equation}\label{e.H00}
(H_{-2m_{i+1}}\circ V_{-2m_i}):
\left(
  \begin{array}{c}
    x \\
    y \\
    \zc \\
  \end{array}
\right)_H\mapsto  \left(
  \begin{array}{rll}
    (4m_{i+1}m_i+1)x&-2m_{i+1} y &\mod 2\Z\\
    -2m_i x&+y &\mod \Z \\
    \zc &&\\
  \end{array}
\right)_H .
\end{equation}

Similarly, on $H_{01}^i(\zc)$
\begin{align}
\left(
  \begin{array}{c}
     x \\
     y\\
    \zc \\
  \end{array}
\right)_H&\!\!\!=
\left(
  \begin{array}{c}
    x \\
    y+1 \\
    \zc-1 \\
  \end{array}
\right)_V \xmapsto[]{V_{-2m_i}}
\left(
  \begin{array}{c}
    x \\
    \{y+1-2m_i x\}_{(0,2]} \\
    \zc-1 \\
  \end{array}
\right)_V\notag\\
&=\left(
  \begin{array}{c}
     x+1 \\
     \{y+1-2m_i x\}_{(0,2]}\\
     \zc-1\\
  \end{array}
\right)_H,\text{ because }\{y+1-2m_i x\}_{(0,2]}\in (0,1]\notag\\
&\xmapsto[]{H_{-2m_{i+1}}}
\left(
  \begin{array}{c}
     \left\{x+1-2m_{i+1}\{y+1-2m_i x\}_{(0,2]}\right\}_{(0,2]}\\
     \{y+1-2m_i x\}_{(0,2]}\notag\\
     \zc-1\\
  \end{array}
\right)_H.
\end{align}
Modulo $2\Z$, the $x$-coordinate equals $(4m_{i+1}m_i+1)x-2m_{i+1}y+1$. Modulo $\Z$, the $y$-coordinate equals $-2m_i x+y$. Thus,  on $H^i_{01}(\zc)$, 
\begin{equation}\label{e.H01}
(H_{-2m_{i+1}}\circ V_{-2m_i}):
\left(
  \begin{array}{c}
    x \\
    y \\
    \zc \\
  \end{array}
\right)_H\!\!\!\mapsto\!  \left(
  \begin{array}{rlll}
    (4m_{i+1}m_i+1)x&-2m_{i+1} y &+1 &\mod 2\Z\\
    -2m_i x &+y & &\mod \Z \\
    \zc & &-1&\\
  \end{array}
\right)_H
\end{equation}

$H_{10}^i(\zc)$ is inside the bottom half of vertical rectangle $\zc$, so on $H_{10}^i(\zc)$,
\begin{align*}
\left(
  \begin{array}{c}
     x \\
     y\\
    \zc \\
  \end{array}
\right)_H&\!\!\!=
\left(
  \begin{array}{c}
    x-1 \\
    y \\
    \zc \\
  \end{array}
\right)_V \xmapsto[]{V_{-2m_i}}
\left(
  \begin{array}{c}
    x-1 \\
    \{y-2m_i x\}_{(0,2]} \\
    \zc \\
  \end{array}
\right)_V\\
&=\left(
  \begin{array}{c}
  x\\
  \{y-2m_i x\}_{(0,2]}\\
  \zc\\
    \end{array}
\right)_H,\text{ because }\{y-2m_i x\}_{(0,2]}\in (0,1] \\
&\xmapsto[]{H_{-2m_{i+1}}}
\left(
  \begin{array}{c}
    \{x-2m_{i+1}\{y-2m_i x\}_{(0,2]}\}_{(0,2]}\\
    \{y-2m_i x\}_{(0,2]}\\
    \zc\\
  \end{array}
\right)_{H}.
\end{align*}
The $x$-coordinate equals $(4m_i m_{i+1}+1)x-2m_{i+1}y\mod 2\Z$, and the $y$-coordinate equals $y-2m_i x\mod\Z$. So on $H^i_{10}(\zc)$
\begin{equation}\label{e.H10}
(H_{-2m_{i+1}}\circ V_{-2m_i}):
\left(
  \begin{array}{c}
    x \\
    y \\
    \zc \\
  \end{array}
\right)_H\mapsto  \left(
  \begin{array}{rll}
    (4m_{i+1}m_i+1)x&-2m_{i+1} y &\mod 2\Z\\
    -2m_i x &+y &\mod \Z \\
    \zc & &\\
  \end{array}
\right)_H
\end{equation}

Similarly, on $H_{11}^i(\zc)$ we have
\begin{align*}
\left(
  \begin{array}{c}
     x \\
     y\\
    \zc \\
  \end{array}
\right)_H&\!\!\!=
\left(
  \begin{array}{c}
     x-1\\
     y\\
     \zc \\
  \end{array}
\right)_V \xmapsto[]{V_{-2m_i}}
\left(
  \begin{array}{c}
    x-1 \\
     \{y-2m_i(x-1)\}_{(0,2]}\\
     \zc\\
  \end{array}
\right)_V=\left(
  \begin{array}{c}
    x-1 \\
     \{y-2m_ix\}_{(0,2]}\\
     \zc\\
  \end{array}
\right)_V\\
&=\left(
  \begin{array}{c}
  x-1\\
  \{y-2m_i x\}_{(0,2]}-1\\
  \zc+1\\
    \end{array}
\right)_H,\text{ because }\{y-2m_ix\}_{(0,2]}\in (1,2]  \\
&\xmapsto[]{H_{-2m_{i+1}}}
\left(
  \begin{array}{c}
\{x-1-2m_{i+1}\{y-2m_ix\}_{(0,2]}\}_{(0,2]}\\
\{y-2m_ix\}_{(0,2]}-1\\
\zc+1\\
  \end{array}
\right)_{H}.
\end{align*}
The $x$-coordinate equals $(4m_i m_{i+1}+1)x-2m_{i+1}y -1\mod 2\Z$, and the $y$-coordinate equals $y-2m_i x \mod\Z$. So on $H^i_{11}(\xi)$
\begin{equation}\label{e.H11}
(H_{-2m_{i+1}}\circ V_{-2m_i}):
\left(
  \begin{array}{c}
    x \\
    y \\
    \zc \\
  \end{array}
\right)_H\!\!\!\!\!\mapsto\!\!  \left(
  \begin{array}{rlll}
    (4m_{i+1}m_i+1)x&-2m_{i+1} y &-1&\mod 2\Z\\
    -2m_i x &+y & &\mod \Z \\
    \zc & & +1&\\
  \end{array}
\right)_H
\end{equation}

Together, \eqref{e.H00}--\eqref{e.H11} imply the lemma.
\end{proof}

\medskip
\noindent
{\bf Twist Calculations for the Building Blocks.}
\begin{lem}\label{l.sign}
Suppose $n_k\in\Z\setminus\{0\}$ for all $k$, and $[0;2n_1,2n_2,\ldots]\not\in\Q$. Then
$$
2n_1\cdot [0;2n_1,2n_2,\ldots]\text{ is in }\begin{cases}
(0,1) & \text{when }n_1 n_2>0\\
(1,2) & \text{when }n_1 n_2<0.
\end{cases}
$$
\end{lem}
\begin{proof}
We saw in Lemma \ref{l.limit-cfe-exists} that 
$
\DS[0;2k_1,2k_2,\ldots]=\lim_{N\to\infty}(v_{2k_1}\circ\cdots\circ v_{2k_N})(0)$,  where $\DS v_{2m}(t)=\frac{1}{2m+t}$. Since $v_{2m}[-1,1]\subset [-1,1]$ for  $m\in 2\Z\setminus\{0\}$, 
$
[0;2k_1,2k_2,\ldots]\in [-1,1]$. Since $-v_{2m}(t)=v_{-2m}(-t)$, 
$
-[0;2k_1,2k_2,\ldots]=[0;-2k_1,-2k_2,\ldots].
$
We claim that
\begin{equation}\label{e.sign-of-ecfe}
\mathrm{sgn}[0;2n_1,2n_2,\ldots]=\mathrm{sgn}(n_1).
\end{equation}
Indeed, $[0;2n_1,2n_2,\ldots]=\mathrm{sgn}(n_1)\cdot [0;|2n_1|,\pm 2n_2,\pm 2n_3,\ldots]\in\{\frac{\mathrm{sgn}(n_1)}{|2n_1|+t}: t\in [-1,1]\}$, and the denominator is always positive.

Let $s_i:=\mathrm{sgn}(n_i)$, then
\begin{align*}
&2n_1\cdot [0;2n_1,2n_2,\ldots]=|2n_1|\cdot [0;|2n_1|,2n_2s_1,2n_3s_1,2n_4s_1,\ldots]\\
&=|2n_1|\cdot [0;|2n_1|,|2n_2|s_1 s_2,2n_3s_1,2n_4s_1,\ldots]=\frac{|2n_1|}{|2n_1|+[0;|2n_2|s_1 s_2,2n_3 s_1,\ldots]}
\end{align*}
By \eqref{e.sign-of-ecfe}, the even CFE in the denominator is in $[0,1]$ when  $s_1 s_2=1$, and in $[-1,0]$,  when $s_1 s_2=-1$. In the first case the fraction is in $(0,1]$. In the second case it is inside
$
\{\frac{|2n_1|}{|2n_1|-t}: t\in [0,1]\}\subset [1,2]
$. The possibilities $0,1,2$ are ruled out by irrationality.
\end{proof}

Suppose $o$ is a singularity, $\vec{v}\in\R^2\setminus\{0\}$, and  $\zc\in\Z$. Let $L(o,\vec{v},\zc)$ denote the ray emanating from $o$,  in the direction of the vector  $\vec{v}$, and initially lying in horizontal rectangle $\zc$. The length of $\vec{v}$ does not matter; it just needs to be non-zero. For each $o$, only some of the $\zc$'s are possible (odd or even, depending on 
$o$ and $\vec{v}$).

\begin{lem}\label{l.twist-positive-1}
Suppose $\beta_1=[0;2m_1,2m_2,\ldots]$ is irrational, and $o$ is the singularity at the bottom left side of horizontal rectangle $\zc$, and $\ho$ the singularity at the middle top side of horizontal rectangle $\zc\pm 1$ ($\ho$ is congruent to $o$).
\begin{enumerate}[(1)]
\item If $m_1 m_2>0$, then
$\DS
(H_{-2m_2}\circ V_{-2m_1})\left[L(o,{\beta_1\choose 1},\zc)\right]=L\left(o, {E^2(\beta_1)\choose 1},\zc\right).
$
\item If $m_1 m_2<0$, then
$\DS
(H_{-2m_2}\circ V_{-2m_1})\left[L(o,{\beta_1\choose 1},\zc)\right]=
L\left(\ho, -{E^2(\beta_1)\choose 1},\zc-\sgn(m_1)\right).
$
\end{enumerate}
Here $E$ is the even Gauss map, i.e. $E^2(\beta_1)=[0;2m_3,2m_4,\ldots]$.
\end{lem}

%

\begin{proof}
\begin{enumerate}[$\bullet$]
\item If $m_1>0$, then $\beta_1>0$, and $L(o,{\beta_1\choose 1},\zc)$ begins with the segment
$\left(
   \begin{array}{c}
     t\beta_1 \\
     t \\
     \zc \\
   \end{array}
 \right)_H
 $,  $0<t\ll 1$. In particular, the ray starts at $ H_{00}^1\cup H_{01}^1$. To see which of these possibilities happens, we use Lemma \ref{l.sign} to see that  for small positive $t$,
 $$
 \{y-2m_1 x\}_{(0,2]}=\{t(1-2m_1[0;2m_1,2m_2,\ldots])\}_{(0,2]}\in \begin{cases}
 (0,1) & m_1 m_2>0\\
 (1,2) & m_1 m_2<0.
 \end{cases}
 $$
Thus, if $m_1>0$ and $m_1m_2>0$ then $L(o,{\beta_1\choose 1},0)$ starts at $H_{00}^1$, and if $m_1>0$ and $m_1 m_2<0$, then $L(o,{\beta_1\choose 1},0)$ starts at $H_{01}^1$.

\item
If $m_1<0$, then $\beta_1<0$, and
 $L(o,{\beta_1\choose 1},\zc)$ begins with
$\left(
   \begin{array}{c}
     2+t\beta_1 \\
     t \\
     \zc \\
   \end{array}
 \right)_H\!\!\! ,0<t\ll 1$, so the ray starts at $H_{10}^i$ or $H_{11}^i$. By Lemma \ref{l.sign}:
  \begin{align*}
\{y-2m_1 x\}_{(0,2]}=\{t-2m_1(2+t\beta_1)\}_{(0,2]}=\{t(1-2m_1\beta_1)\}_{(0,2]}\in
\begin{cases}
(0,1) & m_1 m_2>0\\
(1,2) & m_1 m_2<0.
\end{cases}
\end{align*}
Thus, if $m_1<0$ and $m_1 m_2>0$ then $L(o,{\beta_1\choose 1},\zc)$ starts at $H_{10}^1$, and if $m_1<0$ and  $m_1 m_2<0$, then 
$L(o,{\beta_1\choose 1},\zc)$ starts at $H_{11}^1$.
 \end{enumerate}

We now apply $H_{-2m_2}\circ V_{-2m_1}$ to $L(o,{\beta_1\choose 1},\zc)$. By Lemma \ref{l.composition-formula}, the result is a linear ray which begins with the segment
\begin{align*}
&\left(
  \begin{array}{c}
    \!\!\!\left\{(4m_1 m_{2}+1) \beta_1 t -2m_2t+Z_1\right\}_{(0,2]}\!\!\! \\
    \left\{-2m_1 \beta_1 t+t\right\}_{(0,1]} \\
    \zc-Z_1 \\
  \end{array}
\right)_H (0<t\ll 1),
\text{ where }
Z_1\!\!=\!\!\begin{cases}
+1 & m_1\!\!>\!\!0, \, m_1 m_2\!\!<\!\!0\\
-1 & m_1\!\!<\!\!0,\, m_1 m_2\!\!<\!\!0\\
0& \text{otherwise.}
\end{cases}
\end{align*}
Equivalently,  $Z_1=\mathrm{sgn}(m_1)\times \frac{1-\mathrm{sgn}(m_1 m_2)}{2}.$

Let us simplify the horizontal coordinates.

Suppose first that $m_1 m_2>0$. Then $Z_1=0$ and $2m_1\beta_1\in (0,1)$. Therefore for all $t$ small, $-2m_1\beta_1 t+t\in (0,1)$. Changing coordinates $\tau=-2m_1\beta_1t+t$, and noting that $\tau>0$,  we see that our ray starts
with
$$
\left(
  \begin{array}{c}
    \left\{\tau \left(\frac{(4m_1 m_{2}+1) \beta_1  -2m_2}{-2m_1 \beta_1 +1}\right)\right\}_{(0,2]} \\
    \tau \\
    \zc \\
  \end{array}
\right)_H\ \ (0<\tau\ll 1).
$$
Next, ${\tiny
\left(
  \begin{array}{cc}
    4m_1 m_2+1 & -2m_2 \\
    -2m_1 & 1 \\
  \end{array}
\right)=
\left(
  \begin{array}{cc}
    1 & -2m_2 \\
    0 & 1 \\
  \end{array}
\right)
\left(
  \begin{array}{cc}
    1 & 0 \\
    -2m_1 & 1 \\
  \end{array}
\right)}
$, therefore
$$
\frac{(4m_1 m_{2}+1) \beta_1  -2m_2}{-2m_1 \beta_1 +1}=-2m_2+\cfrac{1}{-2m_1+\cfrac{1}{\beta_1}}=[0;2m_3,2m_4,\ldots]=E^2(\beta).
$$
This is positive when $m_3>0$ and negative when $m_3<0$. So our ray starts with
$$
\left(
  \begin{array}{c}
    \tau E^2(\beta)+(1-\sgn(m_3)) \\
    \tau \\
    \zc \\
  \end{array}
\right)_H\ \ (0<\tau\ll 1).
$$
Either way, if $m_1 m_2>0$, then 
$(H_{-2m_2}\circ V_{-2m_1})[L(o,{\beta_1\choose 1},\zc)]
=L(o,{\beta_2\choose 1},\zc)$, where $$\beta_2:=E^2(\beta_1)=[0;2m_3,2m_4,\ldots].$$

Next suppose $m_1 m_2<0$ and $m_1>0$. In this case $Z_1=+1$ and $2m_1\beta_1\in (1,2)$. Therefore for all $t>0$ small, we have the following identities for the $x$ and $y$ coordinates:
$$
\{(4m_1 m_2+1)\beta_1t-2m_2 t+Z_1\}_{(0,2]}=
(4m_1 m_2+1)\beta_1t-2m_2 t+1,
$$
$$
\{-2m_1\beta_1 t+t\}_{(0,1]}=t(1-2m_1\beta_1)+1.
$$
It follows that our ray begins with
\begin{align*}
\left(
  \begin{array}{c}
    (4m_1 m_2+1)\beta_1t-2m_2 t+1 \\
    t(1-2m_1\beta_1)+1 \\
    \zc-1 \\
  \end{array}
\right)_H\ \ (0<t\ll 1)
\end{align*}
Changing coordinates into $\tau:=-t(1-2m_1\beta_1)$ (we need the minus to keep $\tau$ positive), we can rewrite this as
\begin{align*}
\left(
  \begin{array}{c}
    -\tau\left(\frac{(4m_1 m_2+1)\beta_1-2m_2}{-2m_1\beta_1+1}\right) +1 \\
    -\tau+1 \\
    \zc-1 \\
  \end{array}
\right)_H\ \ (0<\tau\ll 1)
\end{align*}
(Note that  as  $\tau\!\!\to \!\!0^{+}$,  we obtain the  singularity $\wh{o}
=(1,1,\zc-1)_H$, which is congruent to $o=(0,0,\zc)_H$.)
The $x$-coordinate simplifies as before to $-\tau\beta_2+1$, where $\beta_2=[0;2m_3,2m_4,\ldots]$. In summary, if $m_1 m_2<0$ and $m_1>0$, then
$$
(H_{-2m_2}\circ V_{-2m_1})[L(o,{\beta_1\choose 1},\zc)]
=L(\wh{o},-{\beta_2\choose 1},\zc-1)\text{,  } \beta_2:=E^2(\beta_1)=[0;2m_3,2m_4,\ldots].
$$

Finally, suppose $m_1 m_2<0$ and $m_1<0$. In this case, $Z_1=-1$ and $2m_1\beta_1\in (1,2)$. Therefore, for all $t>0$ small, the $x$ and $y$ coordinates are given for $t>0$ small by
$$
\left\{(4m_1 m_{2}+1) \beta_1 t -2m_2t-1\right\}_{(0,2]}=(4m_1 m_{2}+1) \beta_1 t -2m_2t+1
$$
$$
 \left\{-2m_1 \beta_1 t+t\right\}_{(0,1]}=t(1-2m_1\beta)+1
$$
So our ray begins with
$$
\left(
  \begin{array}{c}
    (4m_1 m_{2}+1) \beta_1 t -2m_2t+1 \\
     t(1-2m_1\beta)+1 \\
    \zc+1 \\
  \end{array}
\right)_H \ (0<t\ll 1).
$$
As before, we can rewrite this in the form
$
\left(
  \begin{array}{c}
    -\tau\beta_2+1 \\
     -\tau+1 \\
    \zc+1 \\
  \end{array}
\right)_H \ (0<\tau\ll 1),
$
whence $(H_{-2m_2}\circ V_{-2m_1})[L(o,{\beta_1\choose 1},\zc)]=
L(\wh{o},-{\beta_2\choose 1},\zc+1)$.

Looking at all these cases, we obtain the lemma.
\end{proof}

\begin{lem}\label{l.twist-negative}
Suppose $\beta_1=[0;2m_1,2m_2,\ldots]$ is irrational, and $\ho$ is the singularity at the middle top  side of horizontal rectangle $\eta$. Let $o$ be the singularity at the bottom left corner of horizontal rectangle $\eta\pm 1$ (which is congruent to $\ho$).
\begin{enumerate}[(1)]
\item If $m_1 m_2>0$, then
$$
(H_{-2m_2}\circ V_{-2m_1})\left[L\left(\ho,-{\beta_1\choose 1},\eta\right)\right]=
L\left( \ho,-{\beta_2\choose 1},\eta\right),\text{ where }\beta_2=E^2(\beta).
$$
\item If $m_1 m_2<0$, then
$$
(H_{-2m_2}\circ V_{-2m_1})\left[L\left(\ho,-{\beta_1\choose 1},\eta\right)\right]=L\left(o,{\beta_2\choose 1},\eta-\sgn(m_1)\right),
$$
\end{enumerate}
\end{lem}
\begin{proof}

$\bullet$ If $m_1>0$, then $\beta_1>0$, and $L(\wh{o},-{\beta_1\choose 1},\eta)$ begins with the segment
$\left(
   \begin{array}{c}
     1-t\beta_1 \\
     1-t \\
     \eta \\
   \end{array}
 \right)_{\!\!\!H}
 $, $0\!\!<\!\!t\!\ll\! 1$. In particular, the ray starts at $H_{00}^1\!\cup\! H_{01}^1$. By  Lemma \ref{l.sign} for small positive $t$
 $$
 \{y-2m_1 x\}_{(0,2]}=\{1-t-2m_1+2m_1\beta_1t\}_{(0,2]}\in \begin{cases}
 (0,1) & m_1 m_2>0\\
  (1,2) & m_1 m_2<0.
 \end{cases}
 $$
 So if $m_1>0$ and $m_1m_2>0$ then $L(\wh{o},{ -{\beta_1\choose 1}},0)$ starts at $H_{00}^1$, and if $m_1>0$ and $m_1 m_2<0$, then $L(\wh{o},{- {\beta_1\choose 1}},\xi)$ starts at $H_{01}^1$.

$\bullet$
If $m_1<0$, then $\beta_1<0$, and
 $L(\wh{o},-{\beta_1\choose 1},\eta)$ begins with
$\left(
   \begin{array}{c}
     1-t\beta_1 \\
     1-t \\
     \eta \\
   \end{array}
 \right)_{\!\!\!H},$ $0<t\ll 1$, so the ray starts at $H_{10}^1$ or $H_{11}^1$. As before, by  Lemma \ref{l.sign},
  \begin{align*}
\{y-2m_1 x\}_{(0,2]}=\{1-t-2m_1+2m_1\beta_1 t\}_{(0,2]}\in
\begin{cases}
(0,1) & m_1 m_2>0\\
(1,2) & m_1 m_2<0.
\end{cases}
\end{align*}
So if $m_1<0$ and $m_1 m_2>0$ then $L(\wh{o},-{\beta_1\choose 1},\eta)$ starts at $H_{10}^i$, and if $m_1<0$ and  $m_1 m_2<0$, then $L(\wh{o},-{\beta_1\choose 1},\eta)$ starts at $H_{11}^i$.

By Lemma \ref{l.composition-formula},  $(H_{-2m_2}\circ V_{-2m_1})
[L(\wh{o},-{\beta_1\choose 1},\eta)]$ begins with the segment
\begin{align*}
&\left(
  \begin{array}{c}
    \left\{(4m_1 m_{2}+1) (1-\beta_1 t) -2m_2(1-t)+Z_1\right\}_{(0,2]} \\
    \left\{-2m_1 (1-\beta_1  t)+(1-t)\right\}_{(0,1]} \\
    \eta-Z_1 \\
  \end{array}
\right)_H\ \ (0<t\ll 1),
\text{ where }\\
&Z_1=\begin{cases}
+1 & m_1>0, m_1 m_2<0\\
-1 & m_1<0, m_1 m_2<0\\
0& \text{otherwise.}
\end{cases}
\end{align*}
$Z_1$ is exactly the same as we had before. But the horizontal coordinates are now
\begin{align*}
&\left(
  \begin{array}{c}
    \left\{-(4m_1 m_{2}+1)\beta_1 t +2m_2t+Z_1+1\right\}_{(0,2]} \\
    \left\{2m_1 \beta_1t-t\right\}_{(0,1]} \\
    \eta-Z_1 \\
  \end{array}
\right)_H\ \ (0<t\ll 1),
\end{align*}

Suppose first that $m_1 m_2>0$. Then $Z_1=0$ and $2m_1\beta_1\in (0,1)$. Therefore for all $t>0$ small, $2m_1\beta_1 t-t<0$ and $\{2m_1\beta_1 t-t\}_{(0,1]}=1+2m_1\beta_1 t-t$. Changing coordinates $\tau=-2m_1\beta_1t+t$, we see that our ray starts with
$$
\left(
  \begin{array}{c}
    \{-\tau \left(\frac{(4m_1 m_{2}+1) \beta_1  -2m_2}{-2m_1 \beta_1 +1}\right)+1\}_{(0,2]} \\
    \{1-\tau\}_{(0,1]} \\
    \eta \\
  \end{array}
\right)_H
=\left(
  \begin{array}{c}
    1-\tau \beta_2 \\
    1-\tau \\
    \eta \\
  \end{array}
\right)_H
\ \ (0<\tau\ll 1).
$$
Here, as before, $\beta_2=[0;2m_3,2m_4,\ldots]=E^2(\beta_1)$.
In summary, if $m_1 m_2>0$, then $$(H_{-2m_2}\circ V_{-2m_1})\left[L\left(\wh{o},-{\beta_1\choose 1},\eta\right)\right]
=L\left(\wh{o},-{\beta_2\choose 1},\eta\right).$$

Next suppose $m_1 m_2<0$ and $m_1>0$. In this case $Z_1=+1$ and $2m_1\beta_1\in (1,2)$. For $t>0$ small, the $y$-coordinate  equals
$\DS
\{2m_1\beta t-t\}_{(0,1]}=(2m_1\beta-1) t=:\tau>0.
$
This allows us to simplify the $x$-coordinate into
\begin{align*}
&\{-(4m_1 m_2+1)\beta_1 t+2m_2t\}_{(0,2]}=\left\{\tau\cdot \frac{-(4m_1 m_2+1)\beta_1 +2m_2}{2m_1\beta-1} \right\}_{(0,2]}\\
&=
\left\{\tau\cdot \frac{(4m_1 m_2+1)\beta_1 -2m_2}{-2m_1\beta_1+1} \right\}_{(0,2]}=\left\{\tau\left(-2m_2+\cfrac{1}{-2m_1+\cfrac{1}{\beta_1}}\right) \right\}_{(0,2]}\\
&=\{\tau \beta_2\}_{(0,2]},\text{ where }\beta_2=[0;2m_3,2m_3,\ldots]=E^2(\beta_1)\\
&=(1-\mathrm{sgn}(m_3))-\tau\beta_2\text{ (because $m_3<0\Leftrightarrow \tau\beta_2<0$ and then we must add $2$).}
\end{align*}
It follows that our ray begins with
$\DS
\left(
  \begin{array}{c}
    \tau\beta_2+(1-\mathrm{sgn}(m_3)) \\
    \tau \\
    \eta-1 \\
  \end{array}
\right)_H\ \ (0<\tau\ll 1).
$
So if $m_1 m_2<0$ and $m_1>0$, then
$$
(H_{-2m_2}\circ V_{-2m_1})\left[L\left(\wh{o},-{\beta_1\choose 1},\eta\right)\right]=L\left(o,{\beta_2\choose 1},\eta-1\right), \text{ where }\beta_2:=E^2(\beta_1).
$$

Finally, suppose $m_1 m_2<0$ and $m_1<0$. In this case, $Z_1=-1$ and $2m_1\beta_1\in (1,2)$. Therefore,  the $y$ and $x$ coordinates are given for $t>0$ small by
$$
 \left\{2m_1 \beta_1t-t\right\}_{(0,1]}=(2m_1\beta_1-1)t=:\tau>0
$$
\begin{align*}
&\left\{-(4m_1 m_{2}+1)\beta_1 t +2m_2t\right\}_{(0,2]}=\{\tau\left(\tfrac{-(4m_1 m_1+1)\beta_1+2m_2}{2m_1\beta_1-1}\right)\}_{(0,2]}\\
&=\{\tau\beta_2\}_{(0,2]},\text{ where }\beta_2=E^2(\beta)=[0;2m_3,2m_3,\ldots]\\
&=(1-\mathrm{sgn}(m_3))+\tau\beta_2.
\end{align*}
So our ray begins with
$
\left(
  \begin{array}{c}
    \tau\beta_2+ (1-\mathrm{sgn}(m_3))\\
     \tau \\
    \eta+1 \\
  \end{array}
\right)_H \ (0<\tau\ll 1).
$
In summary,  if $m_1 m_2<0$ and $m_1<0$, then
$
(H_{-2m_2}\circ V_{-2m_1})\left[L\left(\wh{o},-{\beta_1\choose 1},\eta\right)\right]=
L\left(o,{\beta_2\choose 1},\eta+1\right).
$

The lemma follows.
\end{proof}

\smallskip
\noindent
{\bf Proof of Theorem \ref{ThDriftOfOrigin}.}
By Lemma \ref{l.bb-id}, to prove the theorem we need to find the  $\Z$-coordinate of the horizontal rectangle which contains the beginning of
\begin{align*}
L_n&=\left(H_{-2m_{2nN_0}}\circ V_{-2m_{2nN_0-1}})\circ \cdots\circ (H_{-2m_2}\circ V_{-2m_1})\right)\left[L\left(\widehat{o},{\beta\choose 1},(-1)\right)\right]
\end{align*}
(the minus one is needed  because $\widehat{o}$ is on the top side of the horizontal rectangle, and we want the initial ray to start at rectangle zero). Since  the  midpoint of the top side of horizontal rectangle (-1) is the bottom left corner of horizontal rectangle zero,  
\begin{align*}
L_n&=\left(H_{-2m_{2nN_0}}\circ V_{-2m_{2nN_0-1}})\circ \cdots\circ (H_{-2m_2}\circ V_{-2m_1})\right)\left[L\left({o},{\beta\choose 1},0\right)\right]
\end{align*}

$L_n$ is given recursively  by $L_0=L\left({o},{\beta\choose 1},0\right)$ and $L_{k+1}=(H_{-2m_{2k+2}}\circ V_{-2m_{2k+1}})(L_k)$. 
An inductive argument using Lemmas \ref{l.twist-positive-1} and \ref{l.twist-negative} shows that either   $L_k=L(o,{E^{2k}(\beta)\choose 1},\cdot)$ (``case A") or $L_k=L(\widehat{o},-{E^{2k}(\beta)\choose 1},\cdot)$ (``case B"). Let $\Delta_{k+1}$ denote the horizontal $\Z$ coordinate of the beginning of $L_{k+1}$ minus  that of $L_k$. $\Delta_{k+1}$ is determined by   Lemma~\ref{l.twist-positive-1} in case A, and by  Lemma~\ref{l.twist-negative} in case B. In both cases we have
\begin{align*}
\Delta_{k+1}&=\begin{cases}
0 & m_{2k+1}m_{2k+2}>0,\\
-\sgn(m_{2k+1}) & m_{2k+1}m_{2k+2}<0.
\end{cases}\\
&=-\sgn(m_{2k+1})\times \frac{1}{2}(1-\sgn(m_{2k+1})\sgn(m_{2k+2}))=
\frac{1}{2} \left[ \sgn(m_{2k+2})-\sgn(m_{2k+1})\right].
\end{align*}
 Summing $\Delta_{k+1}$ over $k=0,\ldots,nN_0-1$, we obtain the theorem.
 \qed

\section{Limit Theorems for Markov Chains and Markov Arrays}
\label{AppMC}
This appendix includes a summary of the  results of \cite{Dolgopyat-Sarig-Book} used in this work. 

\medskip
\noindent
{\bf The Basic Setup:} An {\em inhomogeneous Markov chain}  is a stochastic process $\mathsf X:=\{X_n\}_{n\geq 0}$ such that \index{Markov chain}
$
\Prob [X_0=a_0, \ldots, X_n=a_n]:=\pi_0(a_0) \pi_{0,1}(a_0, a_1)
\cdots \pi_{n-1,n}(a_{n-1}, a_n)
$, 
where for each $i$, 
\begin{itemize}
\item $a_i\in\mathfrak S_i$, where $\mathfrak S_i$ is a finite set called the {\em state space} at time $i$; \index{State spaces} 
\index{Markov chain!state spaces}
\item $\DS \pi_0:=\sum_{a\in\mathfrak S_0}\pi_0(a) \delta_{a}$ is a probability measure on $\mathfrak S_0$ called the {\em initial distribution}\index{Initial distribution}
\item $\DS \pi_{i,i+1}(a_i, E):=\sum_{a_{i+1}\in E}\pi_{i,i+1}(a_i, a_{i+1})$ $(a_i\in\frak S_i)$ is a family of  probability measures on $\mathfrak S_{i+1}$, called the {\em transition kernel} at time $i$.\index{Transition kernel}\index{Markov chain!transition kernels}
\end{itemize}
An {\em additive functional} \index{Additive functional} on $\mathsf X$ is a sequence of functions $\mathsf f=\{f_i\}_{i\geq 0}$,  $f_i:\mathfrak S_i\times\mathfrak S_{i+1}\to\R$. We say that $\mathsf f$ is {\em uniformly bounded} \index{Additive functional!uniform boundedness} (on $\mathsf X$), if 
$\DS K:=\sup_n \esssup  f_{n}(X_n, X_{n+1})<\infty$.  
Let
$$
S_n=f_0(X_0,X_1)+\cdots+f_{n}(X_n,X_{n+1})\ ,\ V_n:=\Var(S_n).
$$

More generally, a  (finite state) {\em inhomogeneous Markov array} is an array $\mathsf X$ \index{Markov array}
$$
\begin{array}{cccccc}
X_0^{(0)} &                 &                 &  &                        &  \\
X_0^{(1)} & X_1^{(1)} &                 &            && \\
X_0^{(2)} & X_1^{(2)} &   X_2^{(2)}              &            & &\\
\vdots      &  \vdots    &                   & \ddots &  &\\
X_0^{(n)} & X_1^{(n)} & X^{(n)}_2 & \cdots & X^{(n)}_{k_N+1}&\\
\vdots      & \vdots       &  \vdots    &             &    &\ddots
\end{array}
$$
where $k_N+1$ is a strictly increasing, and where each row has joint distribution 
$$
\Prob [X_0^{(N)}=a_0, \ldots, X_{k_N+1}^{(N)}=a_{k_N+1}]:=\pi_{0}^{(N)}(a_0)\pi^{(N)}_{0,1}(a_0,a_1)\cdots 
\pi_{k_N,k_N+1}^{(N)}(a_{k_N},a_{k_N+1}),
 $$ 
 such that  for each $i$, 
\begin{itemize}
\item $a_i\in\mathfrak S_i^{(N)}$, where $\mathfrak S_i^{(N)}$ is a finite set called the {\em state space} at time $i$ for row $n$; 
\item $\DS\pi_0^{(N)}:=\sum_{a\in\mathfrak S_0^{(N)}}\pi_0^{(N)}(a) \delta_{a}$ are probability measures  called the {\em initial distributions};
\item $\DS\pi_{i,i+1}^{(N)}(a_i, E):=\sum_{  a_{i+1}\in E}
\pi_{i,i+1}^{(N)}(a_i, a_{i+1})$ $(a_i\in\frak S_i^{(n)})$ is  a family of  probability measures on $\mathfrak S_{i+1}^{(N)}$ called  the {\em transition kernel} at  time $i$.
\end{itemize}
The {\em $N^{\text{th}}$ row} is  $(X_0^{(N)}, X_1^{(N)}, \cdots , X^{(N)}_{k_N+1})$.  
No joint distribution is specified for elements of different rows. 

An {\em additive functional}\index{Additive functional} on $\mathsf X$ is an array of functions 
$
\mathsf f=\{f_i^{(N)}\!:\! i=0,\ldots, k_N; N\geq 0\}$, $f_i^{(N)}:\mathfrak S_i^{(N)}\times \mathfrak S_{i+1}^{(N)}\to\R. 
$ We say that $\mathsf f$ is {\em uniformly bounded} (on $\mathsf X$), if there is a constant $K$ such that $|f_{i}^{(N)}(X_i^{(N)}, X_{i+1}^{(N)})|<K$ a.s. for all $i,N$.
Let 
$$
\DS S_N:=f_0(X_0^{(N)},X_1^{(N)})+\cdots+f_{k_N}^{(N)}(X_{k_N}^{(N)},X_{k_N+1}^{(N)}), 
\quad V_N:=\mathrm{Var}(S_N).
$$

The Markov arrays in this work will mostly arise  from the ``change of measure'' construction in \S\ref{ss.change-of-measure}. In Lemma  \ref{l.LLT-LD-Blocks}, we will study Markov arrays with rows equal to blocks of increasing size of a given  Markov chain. 

Results on Markov arrays yield results on Markov chains $\{X_n\}$, when applied to the ``trivial  array" with rows $(X_0,\ldots,X_{N+1})$. 
Array results imply various uniformities in Markov chain results, see for instance the proof of Lemma \ref{l.LLT-LD-Blocks}.

In this work,  
all additive functionals will be integer valued, and they will all satisfy  $f_0\equiv 0$, $f_{i}=f_{i}(X_{i}^{(N)})$, $f_0^{(N)}\equiv 0$ and $f_{i}^{(N)}=f_{i}^{(N)}(X_{i}^{(N)})$.    But the results in the appendix work in the greater generality described above. There are also analogous results for continuous state spaces and real-valued additive functionals, see \cite{Dolgopyat-Sarig-Book}. 

\medskip\noindent
{\bf Uniform Ellipticity and Exponential Mixing.} Uniform ellipticity of Markov chains was defined in 
\S\ref{ss.uniform-ellipticity}. The definition for Markov arrays is similar:\index{Uniform ellipticity}

\begin{defi}
A Markov array $\mathsf X$ is called  {\em uniformly elliptic},\index{Uniform ellipticity} if there is a  constant $\epsilon_0$, called an {\em ellipticity constant}, and a representation $\pi_{n,n+1}^{(N)}(y,y')=p_n^{(N)}(y,y')\mu_{n+1}^{(N)}(y')$, where for all $N$ and $0\leq n\leq k_N$, 
\begin{enumerate}[(1)]
\item $\mu_n^{(N)}$ are probability measures on $\mathfrak S_n$;
\item $0\leq p_n^{(N)}(\cdot,\cdot)\leq 1/\epsilon_0$;
\item $\DS \sum_{y\in\mathfrak S_{n+1}} p_n^{(N)}(x,y)p_{n+1}^{(N)}(y,z)>\epsilon_0$ for all $x\in\mathfrak S_n^{(N)}, z\in\mathfrak S_{n+2}^{(N)}$. 
\end{enumerate} 
\end{defi} 
Subject to this assumption we have the following analogue of Corollaries \ref{Cor-rho-mixing} and \ref{Cor-phi-mixing} (\cite{Sethuraman-Varadhan},\cite[Prop. 2.13]{Dolgopyat-Sarig-Book}): \index{Mixing}\index{Exponential mixing}
\begin{prop}\label{p.mixing-for-arrays}
Let $\mathsf X$ be a uniformly elliptic Markov array with row lengths $k_N+1$. There exist $\theta\in(0,1)$ and $C_{mix}>1$, which only depend on the ellipticity constant, such that
for any functions 
$g_m^{(N)}:\fS_m^{(N)}\to\reals, g_n^{(N)}:\fS_n^{(N)}\to\reals$ and  $0\leq m<n<k_N+1$,  
\begin{align}
& \left\|\EXP(g_n^{(N)}(X_n^{(N)})|X_m^{(n)})-\E(g_n^{(N)})\right\|_{L^\infty}<C_{mix} \theta^{n-m} \|g_n^{(N)}\|_{L^\infty}\label{EM-1}\\ 
& \left|\Cov(g_m^{(N)}(X_m), g_n^{(N)}(X_n))\right|<C_{mix} \theta^{n-m}\sqrt{\Var(g_m^{(N)}(X^{(n)}_m))\Var(g_n^{(N)}(X^{(N)}_n))} 
\label{EM-2} 
\end{align}
\end{prop}
\noindent
Note that $C_{mix}$ and $\theta$ are independent of $N$. 

\medskip
\noindent
{\bf Structure Constants.} In \S\ref{ss.structure-constants}, we defined  the structure constants \index{Structure constants}
$u_n, d_n(\xi)$ $(n\!\!\geq\!\!1)$ for additive functionals on  Markov chains.
The structure constants
 for Markov arrays
are defined  for each row in exactly the same way.  We denote them by 
$$
u_n^{(N)}\ , \ d_n^{(N)}(\xi)\ \ (N\geq 1\ , \ 3\leq n\leq k_N, \xi\in\R)$$
$$U_N:=\sum_{n=3}^{k_N} (u_n^{(N)})^2\ ,\  
 D_N(\xi):=\sum_{n=3}^{k_N} d_n^{(N)}(\xi)^2.
$$
It is easy to check that  $d_n^{(N)}(\xi)\leq |\xi| u_n^{(N)}$, see \cite[Lem 2.16]{Dolgopyat-Sarig-Book}.

\medskip
\noindent
{\bf Variance Growth and Center-Tightness.} 
The following is  \cite[Thm 3.6]{Dolgopyat-Sarig-Book}:
\begin{thm}
\label{LmVarCycles}
Let $\mathsf X$ be a uniformly\index{Variance of $S_N$!in terms of structure constants} elliptic Markov array, and let $\mathsf f$ be a uniformly bounded additive functional on $\mathsf X$. Then there are constants $C_1, C_2$ 
which depend only $\DS \sup_n\|f_n\|_\infty$ and the ellipticity constant,  such that
$$ \frac{U_N}{C_1}-C_2\leq V_N \leq C_1 U_N+C_2. $$
In particular $V_N\to\infty$ iff $U_N\to\infty$, and if $\xi\neq 0$, then  
$V_N\geq C_1^{-1}\xi^{-2}D_N(\xi)-C_2.$
\end{thm}
\begin{remark}\normalfont
In particular, if $\mathsf X$ is a uniformly elliptic Markov \underline{chain}, and $\mathsf f$ is a uniformly bounded additive functional on $\mathsf X$, then either $V_N\to\infty$ or $V_N$ is bounded. The reason is that in this case $u_n^{(N)}=u_n$, and therefore $U_N\xrightarrow[N\to\infty]{}\sum u_n^2$, which is either finite or infinite. The remark is false for arrays, for obvious reasons. 
\end{remark}

Next we study  the case when $V_N$ is bounded. We need the following definition. \index{Variance of $S_N$!bounded variance}\index{Center-tightness!and variance growth}
\begin{defi}
An additive functional $\mathsf f$ on a Markov array $\mathsf X$ is called {\em center-tight}, if there are constants $m_N\in\R$ such that for all $\eps>0$, there is some $M>0$ such that 
$$
\Prob[|S_N-m_N|>M]<\eps\text{ for all }N. 
$$
\end{defi}

The following result is shown in \cite[Chap. 3]{Dolgopyat-Sarig-Book}: 
\begin{thm}\label{t.V-n-bounded}
Let $\mathsf X$ be a uniformly elliptic Markov array, and let $\mathsf f$ be a uniformly bounded additive functional on $\mathsf X$. Then the following are equivalent:
\begin{enumerate}[(1)]
\item $V_N$ is bounded;
\item $\mathsf f$ is center-tight;
\item There exist uniformly bounded additive functionals $\mathsf a$ and $\mathsf h$ on $\mathsf X$ such that \\
$a_n^{(N)}=a_n^{(N)}(X_n^{(N)})$, $h_n^{(N)}=h_n^{(N)}(X_n^{(N)},X_{n+1}^{(N)})$, $\DS\sup_N\sum_{n=1}^{k_N} \Var [h_n^{(N)}(X_n^{(N)},X_{n+1}^{(N)})]<\infty$, and
\begin{equation}\label{e.coboundary}
f_n^{(N)}=a_{n+1}^{(N)}-a_n^{(N)}+h^{(N)}\text{ a.s.}
\end{equation}
\end{enumerate}
In the special case of Markov chains, (3) holds with $a_n^{(N)},h_n^{(N)}$ independent of $N$, and 
$\DS \sum_{n\geq 1} \left[ h_n-\E(h_n)\right]$ converges a.s.
\end{thm}
\noindent
The right hand side of \eqref{e.coboundary} should be viewed as the inhomogeneous analogue of a ``coboundary."

\medskip
\noindent
{\bf Central Limit Theorem.}  
The following result, up to small differences in the setup (discussed below), is due to Dobrushin \cite{Dobrushin}. For a modern proof, see Sethuraman \& Varadhan \cite{Sethuraman-Varadhan}:\index{Dobrushin's theorem}\index{Central limit theorem!for Markov arrays}\index{Central limit theorem!Dobrushin's theorem}
\begin{thm}[Dobrushin]
\label{ThMCCLT}
Suppose $\mathsf X$ is a uniformly elliptic Markov array, with a uniformly bounded additive functional $\mathsf f$. If $V_N\to\infty$, then for all $a<b$, 
$$
\Prob\left[\frac{S_N-\E(S_N)}{\sqrt{V_N}}\in (a,b)\right]
\xrightarrow[N\to\infty]{}\frac{1}{\sqrt{2\pi}}\int_a^b e^{-t^2/2} dt.
$$
\end{thm}
The papers of Dobrushin and of Sethuraman \& Varadhan allow  more general state spaces and (certain) unbounded additive functionals,  but they use a stronger ellipticity condition, which enables them to replace the assumption $V_N\to\infty$ by a more transparent  assumption on the individual variances of $f_n$. The uniform ellipticity we use was tailored to fit the Markov chain arising from the  staircase. 
 The proof in  \cite{Sethuraman-Varadhan} works without major changes, under the assumptions stated above, 
see \cite[\S 3.2.4]{Dolgopyat-Sarig-Book}. 

\medskip
\noindent
{\bf Invariance Principle and Law of the Iterated Logarithm.} 
The following results are  Theorem 3.24 and Corollary 3.25 in \cite{Dolgopyat-Sarig-Book}:

\begin{thm}[Almost Sure Invariance Principle]
\label{ThASIP}
Suppose $\mathsf X$ is a uniformly elliptic Markov \underline{chain}, and $\mathsf f$ is a uniformly bounded additive functional on $\mathsf X$, such that $V_N\to\infty$. Then there exist $\delta>0$,\index{Almost sure invariance principle}
a probability space $(\tilde\Omega, \tilde\Prob)$ and  measurable functions
$\tS_N,W(t):\tOmega\to\reals$ $(N\in \naturals, t\in \reals^+)$ 
so that
\begin{enumerate}[(1)]
\item  the stochastic processes $S_N-\E(S_N)$ and $\tS_N$ are equal in distribution;

\item  $W(t)$ is equal in distribution to the  standard Brownian motion;

\item  For a.e. $\omega\in\tOmega$, for all $N$ large enough, $|\tS_N-W(V_N)|\leq V_N^{1/2-\delta}$. 
\end{enumerate}
\end{thm}

The almost sure invariance principle  and standard properties of Brownian motion imply  the following:

\begin{thm}[Law of Iterated Logarithm]\label{t.LIL}
Suppose $\mathsf X$ is a uniformly elliptic Markov \underline{chain}, and $\mathsf f$ is a uniformly bounded additive functional on $\mathsf X$, such that $V_N\to\infty$. Then the following holds with probability one:\index{Law of the iterated logarithm}
\begin{equation}
\label{LIL}
\limsup_{N\to\infty} \frac{S_N-\E(S_N)}{\sqrt{V_N\log \log V_N}}=\sqrt{2},
\quad
\liminf_{N\to\infty} \frac{S_N-\E(S_N)}{\sqrt{V_N\log \log V_N}}=-\sqrt{2}.
\end{equation}
\end{thm}
\noindent
These results  do not make sense for Markov arrays. The rows of a Markov array are not endowed with a joint distribution, therefore the random variables $S_1, S_2, S_3,\ldots$ do not form  a stochastic process. 

\medskip
\noindent
{\bf Algebraic Range, Essential Range, and Irreducibility.} 
The {\em central} limit theorem describes  $\Prob[S_N\in (\E(S_N)-a\sqrt{V_N}, \E(S_N)+b\sqrt{V_N})]$. The  {\em local} limit theorem describes  $\Prob[S_N=z_N]$, with constants $z_N$ not too far from $\E(S_N)$.  

The LLT requires additional assumptions on the arithmetic properties of $S_N$ and $z_N$. For example, if for each $n$, $f_n(X_n,X_{n+1})= c_n\, \mathrm{mod }\, \ell$ a.s., then $\Prob[S_N=z_N]$ can only be non-zero if 
$\DS z_N=\sum_{i=0}^{k_N} c_i\mod\ell$ for all $n$. 

For real-valued $\mathsf f$ we have a more subtle pathology, when the distribution of $S_n$ concentrates ``near but not exactly on''  cosets of $\ell\Z$, due to some center-tight term which takes values outside $\ell\Z$. This leads to the following definitions.

Suppose $\mathsf f$ is an additive functional on a Markov array. 

\begin{defi}\label{d.alg-rng}
The {\em algebraic range}\index{Algebraic range} of $(\mathsf X,\mathsf f)$ is the intersection $G_{alg}(\mathsf X,\mathsf f)$ of all closed subgroups $G$ of $(\R,+)$ for which there are constants $c^{(N)}_n$ such that 
\begin{equation}\label{e.alg-range}
\Prob[f_n^{(N)}(X_n^{(N)},X_{n+1}^{(N)})-c^{(N)}_n\in G]=1\text{ for all }0\leq n\leq k_N, N\in\N.
\end{equation}
\end{defi}
\noindent
It can be shown that $G_{alg}(\mathsf X,\mathsf f)$ itself satisfies \eqref{e.alg-range}, therefore $G_{alg}(\mathsf X,\mathsf f)$ is the {\em minimal} closed group with property \eqref{e.alg-range}, see 
 \cite[Lem 4.15]{Dolgopyat-Sarig-Book}.

\begin{defi}\label{d.ess-rng}
The {\em essential range}\index{Essential range} of $(\mathsf X,\mathsf f)$ is the intersection 
$$
G_{ess}(\mathsf X,\mathsf f):=\bigcap \left\{G_{alg}(\mathsf X, \mathsf g):\mathsf g-\mathsf f\text{ is center-tight} \right\}. 
$$
\end{defi}

\begin{defi}
An additive functional $\mathsf f$ such that $G_{ess}(\mathsf X,\mathsf f)=G_{alg}(\mathsf X,\mathsf f)$, is called {\em irreducible}\index{Irreducibility} (on $\mathsf X$).
\end{defi}
\noindent
The reduction lemma of \cite{Dolgopyat-Sarig-Book} says that there exists a center-tight additive functional $\mathsf h$ such that $G_{alg}(\mathsf X,\mathsf f-\mathsf h)=G_{ess}(\mathsf X,\mathsf f-\mathsf h)=G_{ess}(\mathsf X,\mathsf f)$. Thus, if $\mathsf f$ is not irreducible, then we can reduce its algebraic range  to the minimum possible by subtracting an  ``inhomogeneous coboundary" of the form \eqref{e.coboundary}. After this modification, the algebraic range equals the essential range, and it cannot be reduced any further.

The following result \cite[Chap 4]{Dolgopyat-Sarig-Book} calculates the essential range from the structure constants:
\begin{defi}
The {\em co-range} of $(\mathsf X,\mathsf f)$ is \index{Co-range}
$
\DS H(\mathsf X,\mathsf f):=\left\{\xi\in\R: 
\sup_N D_N(\xi)<\infty 
\right\}
$ 
\end{defi}

\begin{thm}\label{t.co-range}
Let  $\mathsf X$ be a uniformly elliptic Markov array, with a uniformly bounded additive functional $\mathsf f$. Then  $H(\mathsf X,\mathsf f)$ is a closed subgroup of $(\R,+)$, and:\index{Essential range!and the co-range}
\begin{enumerate}[(1)]
\item If $H(\mathsf X,\mathsf f)=\{0\}$, then $G_{ess}(\mathsf X,\mathsf f)=\R$;
\item If $H(\mathsf X,\mathsf f)=\R$, then $G_{ess}(\mathsf X,\mathsf f)=\{0\}$, and $\mathsf f$ is center-tight;
\item If $H(\mathsf X,\mathsf f)=t\Z$ with $t\neq 0$, then $G_{ess}(\mathsf X,\mathsf f)=\frac{2\pi}{t}\Z$.
\end{enumerate}
\end{thm}

\begin{cor}
Suppose $\mathsf X$ is uniformly elliptic and $\mathsf f$ is a uniformly bounded  \underline{$\Z$-valued} additive functional. Then $(\mathsf X,\mathsf f)$ is irreducible with essential range $\Z$, iff for each integer  $\ell>1$, \index{Irreducibility}
$\DS \sup_N D_N(2\pi/\ell)=\infty$ (or for Markov chains, equivalently,  
$\DS \sum d_n^2\left(\frac{2\pi}{\ell}\right)=\infty$).
\end{cor}
\begin{proof}
Since $\mathsf f$ is integer valued,  $d_n^{(N)}(2\pi k)=0$ for all $k\in\Z$, and  $H(\mathsf X,\mathsf f)\supset 2\pi\Z$. Necessarily, $H(\mathsf X,\mathsf f)=\R$, or $H(\mathsf X,\mathsf f)=\frac{2\pi}{\ell}\Z$ for some $\ell\in\N$. 

Suppose $\DS (\forall\ell\in\N\setminus\{1\},\ \sup_N D_N(2\pi/\ell)=\infty)$, then  $H(\mathsf X,\mathsf f)\neq \R$,  and $H(\mathsf X,\mathsf f)\neq \frac{2\pi}{\ell}\Z$ with $\ell\in\N\setminus\{1\}$. Thus $H(\mathsf X,\mathsf f)=2\pi\Z$, and $G_{ess}(\mathsf X,\mathsf f)=\Z$. So
$$\Z=G_{ess}(\mathsf X,\mathsf f)\subset G_{alg}(\mathsf X,\mathsf f)\subset\Z,
$$
and $(\mathsf X,\mathsf f)$ is irreducible with range $\Z$. 

Conversely, if 
$G_{ess}(\mathsf X,\mathsf f)\!\!=\!\!G_{alg}(\mathsf X,\mathsf f)\!\!=\!\!\Z$, then 
$H(\mathsf X,\mathsf f)\!\!=\!\!2\pi\Z$, whence 
$\DS \sup_N D_N\left(\frac{2\pi}{\ell}\right)\!\!=\!\!\infty$ for all $\ell\in\N\setminus\{1\}$.
\end{proof}

\begin{lem}\label{LmChange}
Let $\mathsf f$ be a uniformly bounded additive functional on a Markov chain $\mathsf X$. Suppose $\xi_n$ is a bounded sequence of real numbers, and let $\wt{\mathsf X}$ denote the Markov arrays with rows $(\wt{X}^{\xi_n}_0,\ldots,\wt{X}_{n+1}^{\xi_n})$ $(n\geq 1)$, see Section \ref{ss.change-of-measure}. Then \index{Change of measure!structure constants}\index{Change of measure!essential range and irreducibility} \index{Change of measure!hereditary property}\index{Irreducibility!preserved under change of measure}
\begin{enumerate}[(1)]
\item The structure constants of $(\wt{\mathsf X},\mathsf f)$ are equal up to bounded multiplicative error to the  structure constants of   $(\mathsf X,\mathsf f)$. 
\item $(\wt{\mathsf X},\mathsf f)$ and $(\mathsf X,\mathsf f)$ have the same essential ranges.
\item $\mathsf f$ is irreducible on $\wt{\mathsf X}$ iff $\mathsf f$ is irreducible on ${\mathsf X}$.
\end{enumerate}
\end{lem}
\noindent
(1) is checked in  \cite[Example 4.12]{Dolgopyat-Sarig-Book}, and (2)--(3) follow from Theorem \ref{t.co-range}.

\medskip
\noindent
{\bf Hereditary and Stably Hereditary Arrays.} A Markov array $\mathsf X'$ is called a  {\em sub-array} of a Markov array $\mathsf X$, if every row of $\mathsf X'$ is a row of $\mathsf X$.  In this case, every additive functional on $\mathsf X$ restricts in the natural way to an additive functional $\mathsf f|_{\mathsf X'}$. \index{Hereditary property}

The reader will likely not be surprised to learn that  the local limit theorem for Markov arrays must assume that $(\mathsf X,\mathsf f)$ satisfies the following property:
\begin{defi}\label{d.hereditary}
Suppose $\mathsf f$ is an additive functional on a Markov array $\mathsf X$. We say that $(\mathsf X,\mathsf f)$ is {\em hereditary}, if $G_{ess}(\mathsf X', \mathsf f|_{\mathsf X'})=G_{ess}(\mathsf X, \mathsf f)$ for all sub-arrays $\mathsf X'$ of $\mathsf X$. 
\end{defi}
\noindent
The hereditary property holds automatically for the arrays in this work, because of the following lemma:
\begin{lem}\label{l.hereditary}
Suppose $\{X_n\}_{n\geq 1}$ is a uniformly elliptic Markov \underline{chain}, and $\{f_n\}_{n\geq 1}$ is a uniformly bounded  additive functional on $\{X_n\}_{n\geq 1}$. Let  $f^{(N)}_n:=f_n$. 
\begin{enumerate}[(1)]
\item Suppose $N_k\uparrow\infty$. The  Markov array with rows $(X_1,\ldots,X_{N_k+1})$,  is  hereditary.
\item  The Markov array with rows $(\wt{X}^{\xi_N}_1,\ldots,\wt{X
}_{N+1}^{\xi_N})$ obtained from the change of measure in \S\ref{ss.change-of-measure} is hereditary, whenever $\xi_N$ is bounded.
\item Suppose $f_n$ are integer-valued, and for every  $\eta\in \R\setminus 2\pi\Z$, 
there are constants $c, M>0$ with the following property:
$$
\forall n\in\N\, \exists k\in (n,n+M)\cap\N\text{ s.t. } d_k(\eta)>c.
$$
Suppose  $\xi_N\in\R$ are bounded, $n_N\in\N$, and  $L_N:=n_{N+1}-n_N\to\infty$. Then the array with rows $(\wt{X}^{\xi_N}_{n_N+1},\wt{X}_{n_N+2}^{\xi_N},\ldots,\wt{X}_{n_{N}+L_N}^{\xi_N})$ is hereditary.
\end{enumerate} 
\end{lem}
\begin{proof} We prove part (3). 
Parts (1) and (2) have similar proofs (see  \cite[\S 4.2.3]{Dolgopyat-Sarig-Book} for details). Let   $H(\mathsf X,\mathsf f)$ denote  the co-range of our array. To prove the hereditary property, it is sufficient to show that 
$
D_N(\eta)\xrightarrow[N\to\infty]{}\infty\text{ on $\R\setminus H(\mathsf X,\mathsf f)$}
$, because this implies that $H(\mathsf X', \mathsf f|_{\mathsf X'})=H(\mathsf X, \mathsf f)$, whence 
$G_{ess}(\mathsf X', \mathsf f|_{\mathsf X'})=G_{ess}(\mathsf X, \mathsf f)$,
 for every sub-array $\mathsf X'$. 

We denote the structure constants of $\{X_n\}$ by $d_n(\eta)$, and the structure constants of the array in (3) by $d^{(N)}_n(\eta)$. 
By Lemma \ref{LmChange} and the boundedness of $\xi_N$,  $d_j^{(N)}(\eta)\!\!\asymp \!\! d_{n_N+j}(\eta)$. By  assumption, if $\eta\in\R\setminus 2\pi\Z$, then 
$
\DS \sum_{j=n_N+3}^{n_N+L_N-3}\!\!\!\!\! d_j(\eta)^2\!\!\gtrsim \!\! c^2 L_N/M
$, and therefore 
 $\DS D_N(\eta)=\sum_{j=3}^{L_N} d_j^{(N)}(\eta)^2\xrightarrow[i\to\infty]{}\infty$ on $\R\setminus 2\pi\Z$. 
 In particular, $H(\mathsf X,\mathsf f)\subset 2\pi\Z$. 
 
 At the same time $H(\mathsf X,\mathsf f)\supset 2\pi\Z$, because $\mathsf f$ is integer-valued. So  
 $H(\mathsf X,\mathsf f)=2\pi\Z$ and $D_N(\eta) \xrightarrow[i\to\infty]{}\infty$ on $\R\setminus H(\mathsf X,\mathsf f)$. 
 \end{proof}

\begin{remark}\normalfont
The local limit theorem for irreducible {\em real-valued} additive functionals requires the stronger {\em stable hereditary property}: $D_N(\eta)\to\infty$ {\em uniformly} on compacts in $\R\setminus H(\mathsf X,\mathsf f)$. But for $\Z$-valued functions, such as those used in this work,  the stable hereditary property and the hereditary property are equivalent \cite[Thm 4.14]{Dolgopyat-Sarig-Book}.
\end{remark}

\medskip
\noindent
{\bf Local Limit Theorem.} The following  theorems are called, respectively, the {\em lattice local limit theorem} and the {\em mixing lattice local limit theorem}.  \index{Local limit theorem} \index{Local limit theorem!mixing LLT for local deviations} \index{Local limit theorem!for local deviations}
\begin{thm}[LLT]\label{t.ThLLT}
Let $\mathsf f$ be an integer-valued,  uniformly bounded, irreducible, and hereditary additive functional on a uniformly elliptic Markov array $\mathsf X$, with algebraic range $\Z$. Then $V_N\to\infty$, and for every sequence of   $z_N\in\Z$ such that $\frac{z_N-\E(S_N)}{\sqrt{V_N}}\xrightarrow[N\to\infty]{}s$,
$$  
\Prob[S_N=z_N]=[1+o(1)]\frac{e^{-s^2/2}}{\sqrt{2\pi V_N}},\text{ as $N\to\infty$}. 
$$
\end{thm}
\noindent
This is a special case of \cite[Thm 5.2']{Dolgopyat-Sarig-Book}. 

The next result  adds  conditioning on the first and  last entries of the $N$-th row. Such conditioning is the symbolic dynamical counterpart of conditioning on  $E^s_n$-pure segments, see \S\ref{ss.pure-segments}. Recall that the last entry in the $N$-th row has index $k_N+1$. 
\begin{thm}[Mixing LLT]\label{ThLLTArrays}
Under the assumptions of  Theorem \ref{t.ThLLT}, if
 $x_1^{(N)}\in\mathfrak S^{(N)}_1$,  $x_{k_N+1}^{(N)}\!\!\in\!\!\mathfrak S^{(N)}_{k_N+1}$, and $\Prob[X^{(N)}_{k_N+1}=x_{k_N+1}^{(N)}]$ is bounded away from zero, then
$$  
\Prob\biggl[S_N=z_N\big|X^{(N)}_{1}=x_{1}^{(N)}\ , \  X^{(N)}_{k_N+1}=x_{k_N+1}^{(N)}\biggr]=[1+o(1)]\frac{e^{-s^2/2}}{\sqrt{2\pi V_N}},\text{ as $N\to\infty$}. 
$$
\end{thm}
\noindent
This is a special case of \cite[Thm 5.4']{Dolgopyat-Sarig-Book}. 

Chapter 5 of \cite{Dolgopyat-Sarig-Book}  includes  non-lattice  versions of these results, for irreducible real-valued additive functionals with essential range $\R$, but we will omit the statements, since they are not needed for this work. Local limit theorems for additive functionals which are not irreducible, are given in \cite[Chapt. 6]{Dolgopyat-Sarig-Book}. The asymptotic behavior in this case is different, and non-universal.

\begin{remark}\normalfont
The local limit theorem in the inhomogeneous case was first worked out 
 for sums of bounded  {\em independent} (non-identically distributed) integer-valued random variables, by Prohorov
\cite{Prohorov}. Rozanov extended this work in \cite{Rozanov} to the  unbounded independent case, subject to appropriate tightness conditions. 
The  non-integer  independent case was done first in
in \cite{Dolgopyat-LLT-2016}.

The local limit theorem for homogeneous Markov chains is due to Nagaev \cite{Nagaev}. Theorems \ref{t.ThLLT} and \ref{ThLLTArrays} on the inhomogeneous Markov case above are from \cite{Dolgopyat-Sarig-Book}. There is a related result by Merlev\'ede, M. Peligrad and C. Peligrad in \cite{Peligrad-LLT}. For a fuller historical account and further references, see \cite[\S 1.7]{Dolgopyat-Sarig-Book}. 

 The rate of convergence in the local limit theorem is discussed in
 \cite{Dolgopyat-Hafouta-iid, Dolgopyat-Hafouta-mc}.
\end{remark}

\medskip
\noindent
{\bf Log Moment Generating Functions, Rate Functions, and Large Deviations.}
In this section, we only consider Markov {\em chains}. 

Suppose $\mathsf f$ is a uniformly bounded additive  functional on a Markov chain $\mathsf X$.
As always, $S_N=f_0(X_0,X_1)+\cdots+f_N(X_N,X_{N+1})$, and $V_N=\Var(S_N)$. 
 \begin{defi} \label{DefFreeEn}
 Suppose $V_N\neq 0$ for all $N\geq \nu_0$. \index{Log-moment generating function}
 The {\em log moment generating function} with parameter $N\geq \nu_0$ is the function   $\cF_N:\R\to\R$ ($N\geq\nu_0$) given by 
$$ \cF_N(\xi)=\frac{1}{V_N} \log\EXP\left(e^{\xi S_N}\right).$$
\end{defi}

The following result \cite[Thm 7.3]{Dolgopyat-Sarig-Book} controls the shape of $\cF_N$ as $N\to\infty$:\index{Log-moment generating function!properties}
\begin{prop}
\label{PrFreeEn}
Suppose $\mathsf f$ is a uniformly bounded additive functional on a uniformly elliptic Markov chain $\mathsf X$, such that $V_N>0$ for all $N\geq \nu_0$. Then for all $N\geq\nu_0$, 
\begin{enumerate}[(1)]
\item  $\cF_N(0)=0$,    $\DS \cF_N'(0)=\frac{\EXP(S_N)}{V_N}$ and $\cF_N''(0)=1.$
\item  $\cF_N$ is real-analytic and  strictly convex on $\R$. Moreover,  given $R>0$ there is a constant
$c_R\!\!>\!\!0$ such that  $\DS  c_R\leq \cF_N''\leq c_R^{-1}\text{ on $[-R, R]$ for all $N\geq \nu_0$.}$
\item Suppose $V_N\to\infty$. For every $\eps>0$ there are $\delta, N_\eps>0$ such that for all $|\xi|\leq \delta$ and $N>N_\eps$, we have 
$e^{-\eps}\cdot \frac{1}{2}\xi^2\leq \cF_N(\xi)-\frac{\E(S_N)}{V_N}\xi\leq e^{+\eps}\cdot \frac{1}{2}\xi^2$.
\end{enumerate}
\end{prop}

Define 
the intervals $(a_N^\infty, b_N^\infty)$ and $(a_N^R, b_N^R)$ $(R>0)$ with endpoints 
 \begin{align*}
a_N^\infty:=\cF'_N(-\infty):=\lim_{\xi\to -\infty}\cF_N'(\xi) & &  a_N^R:=\cF'_N(-R)\\
b_N^\infty:=\cF'_N(+\infty):=\lim_{\xi\to +\infty}\cF_N'(\xi) & & b_N^R:=\cF_N'(R).
\end{align*}
If $V_N\neq 0$, then  these intervals are  non-empty, 
by the strict convexity of $\cF_N$. 
\begin{defi}
Suppose $\mathsf f$ is a uniformly bounded additive functional\index{Log-moment generating function!Legendre transform}\index{Rate function} on a uniformly elliptic Markov chain $\mathsf X$, and $V_N>0$ for all $N\geq \nu_0$. 
The {\em rate function} with parameter $N\geq \nu_0$ is the Legendre transform of $\cF_N$, namely $\cI_N:(a_N^\infty,b_N^\infty)\to\R$ given by 
$$
\cI_N(\eta):=\xi\eta-\cF_N(\xi)\text{ for the unique $\xi\in\R$ such that $\cF_N'(\xi)=\eta$}.
$$
\end{defi}
\noindent
The definition is proper, by the real-analyticity and strict convexity of $\cF_N$.

The following result \cite[Thm 7.4]{Dolgopyat-Sarig-Book} controls the shape of $\cI_N$ as $N\to\infty$:
\begin{prop}\label{PrRateFn}
Suppose $\mathsf f$ is a uniformly bounded additive functional on a uniformly elliptic Markov chain $\mathsf X$, such that $V_N>0$ for all $N\geq \nu_0$. Then
\begin{enumerate}[(1)]
\item $\exists c,\nu_1,R>0$ such that for all $N>\nu_1$,
$$
\dom(\cI_N)\supset [a_N^R,b_N^R]\supset\left[\frac{\E(S_N)}{V_N}-c,\frac{\E(S_N)}{V_N}+c \right].
$$
\item For each $R$ there exists $\rho=\rho(R)>0$ such that for all $N>\nu_1$, 
$$
\rho^{-1}\leq \cI_N''\leq\rho\text{ on }[a_N^R,b_N^R].
$$
\item Suppose $V_N\to\infty$. For every $\eps>0$, there exists $\delta>0$ and $N_\eps$ such that for all $\eta\in \left[\frac{\E(S_N)}{V_N}-\delta, \frac{\E(S_N)}{V_N}+\delta\right]$ and $N>N_\eps$, 
    $$
    e^{-\eps}\cdot \frac{1}{2}\left(\eta-\frac{\E(S_N)}{V_N}\right)^2\leq \cI_N(\eta)\leq     e^{\eps}\cdot \frac{1}{2}\left(\eta-\frac{\E(S_N)}{V_N}\right)^2.
    $$
\item Suppose $V_N\to\infty$ and $\frac{z_N-\E(S_N)}{V_N}\xrightarrow[N\to\infty]{}0$, then 
$$
V_N\cI_N\left(\frac{z_N}{V_N}\right)=\frac{1+o(1)}{2}
\left(\frac{z_N-\E(S_N)}{\sqrt{V_N}}\right)^2\text{, as $N\to\infty$.}
$$  
\end{enumerate}
\end{prop}

The following  result on large deviations is Theorem 7.9 in \cite{Dolgopyat-Sarig-Book}.
Note that there are no irreducibility assumptions on the additive functional.\index{Large deviations bounds}

\begin{thm}\label{t.LDP}
Suppose $\mathsf f$ is a uniformly bounded additive functional on a uniformly elliptic Markov chain such that  $V_N\to\infty$. For every $\eps,R>0$, there are constants $D, \nu_2>0$   such that for all $N>\nu_2$:
\begin{enumerate}[(1)]
\item  If $\frac{z_N}{V_N}\in [\cF_N'(\eps), \cF_N'(R)]$, then 
$$
\frac{D^{-1}}{\sqrt{V_N}}e^{-V_N\cI_N\left(\frac{z_N}{V_N}\right)}\leq \Prob(S_N\geq z_N)\leq 
\frac{D}{\sqrt{V_N}}e^{-V_N\cI_N\left(\frac{z_N}{V_N}\right)};
$$

\item  If $\frac{z_N}{V_N}\in [\cF_N'(-R),\cF_N'(-\eps)]$, then 
$$
\frac{D^{-1}}{\sqrt{V_N}}e^{-V_N\cI_N\left(\frac{z_N}{V_N}\right)}\leq \Prob(S_N\leq z_N)\leq 
\frac{D}{\sqrt{V_N}}e^{-V_N\cI_N\left(\frac{z_N}{V_N}\right)}.
$$
\end{enumerate}
Moreover,  $D$ can be chosen to depend only on $\eps$, $R$, and 
$\DS K:=\sup_n\esssup|f_n(X_n,X_{n+1})|$.
\end{thm}

\medskip
\noindent
{\bf Local Limit Theorems for Moderate Deviations and for Large Deviations.} When studying the asymptotic behavior of $\Prob[S_N=z_N]$, it is convenient to distinguish between the following three regimes, listed in increasing level of generality:
\begin{enumerate}[$\bullet$]
\item {\em Local Deviations:} $\frac{z_N-\E(S_N)}{\sqrt{V_N}}$ is  bounded;\index{Local deviations regime}
\item {\em Moderate Deviations:} $\frac{z_N-\E(S_N)}{{V_N}}\xrightarrow[N\to\infty]{}0$;\index{Moderate deviations regime}
\item {\em Admissible Large Deviations:} $\frac{z_N-\E(S_N)}{{V_N}}$ is bounded but ``not too large" (see below). \index{Admissible large deviations regime}\index{Large deviations bounds!LLT for moderate deviations}\index{Local limit theorem!for moderate deviations}
\end{enumerate}

The regime of local deviations is covered by Theorem \ref{t.ThLLT}. We now discuss the other two regimes. The following result is a special case of Theorem 7.6 in \cite{Dolgopyat-Sarig-Book}:

\begin{thm}[LLT for Moderate Deviations]
Let $\mathsf f$ be an irreducible, uniformly bounded,  additive functional with algebraic range $\Z$, on a uniformly elliptic Markov chain $\mathsf X$. For every sequence of integers $z_N$ such that $\frac{z_N-\E(S_N)}{V_N}\to 0$, 
\begin{align*}
\Prob[S_N=z_N]&=\frac{1+o(1)}{\sqrt{2\pi V_N}}\exp\left[-V_N\cI_N\left(\frac{z_N}{V_N}\right)\right]\\
&=\frac{1+o(1)}{\sqrt{2\pi V_N}}\exp\left[-\frac{1+o(1)}{2}\left(\frac{z_N-\E(S_N)}{\sqrt{V_N}}\right)^2\right].
\end{align*}
\end{thm}

Next we discuss the regime of large deviations. The first step is to determine the requirement that  $\frac{z_N-\E(S_N)}{{V_N}}$ be ``not too large."
\begin{defi}
Suppose $\mathsf f$ is a uniformly bounded additive functional on a uniformly elliptic Markov chain, such that $V_N\to\infty$.  
A sequence of numbers $\{z_N\}$  is called  {\em admissible} for $(\mathsf X, \mathsf f)$,  if there exists $R$ such that 
$\frac{z_N}{V_N} \in [a_N^R, b_N^R]$ for all $N$ large.\index{Admissible large deviations regime}
\end{defi}
\noindent
Any sequence $z_N$ such that $|\frac{z_N-\E(S_N)}{V_N}|<c$, with $c$ as in Proposition \ref{PrRateFn}, is admissible. 
The following result is Theorem 7.8(3) in  \cite{Dolgopyat-Sarig-Book}:
\begin{thm}[LLT for Large Deviations]
\label{ThLD}
Let $\mathsf f$ be an integer-valued uniformly bounded additive functional on a uniformly elliptic Markov chain $\mathsf X$. Suppose $(\mathsf X,\mathsf f)$ is irreducible  with algebraic range $\Z$.  Then there are positive functions $\rho_N(\cdot)$ as follows:\index{Large deviations bounds!LLT for admissible large deviations}\index{Local limit theorem!for admissible large deviations}
\begin{enumerate}[(1)]
\item  For every admissible  sequence of integers $z_N$, 
$$
\Prob[S_N=z_N]=[1+o(1)]\frac{e^{-V_N\cI_N\left(\frac{z_N}{V_N}\right)}}{\sqrt{2\pi V_N}}\times\rho_N\left(\frac{z_N-\E(S_N)}{V_N}\right), \text{ as $N\to\infty$}.
$$
\item $\rho_N(\eta)\xrightarrow[\eta\to 0]{}1$ uniformly in $N$, and  $\forall R>0$, $\exists c_R>1$ such that for $N$ and $z_N$ such that $\frac{z_N}{V_N}\in [a_N^R,b_N^R]$, 
$
c_R^{-1}\leq \rho_N\left(\frac{z-\E(S_N)}{V_N}\right)\leq c_R.
$
\item In particular, if $z_N$ is admissible, then $\rho_N\left(\frac{z_N-\E(S_N)}{V_N}\right)\asymp 1$, and if $\frac{z-\E(S_N)}{V_N}\to 0$, then $\rho_N\left(\frac{z_N-\E(S_N)}{V_N}\right)\to 1$.
\end{enumerate}
\end{thm}
\noindent
Analogous results for irreducible real-valued additive functionals are given in \cite{Dolgopyat-Sarig-Book}.

\bibliographystyle{alpha}
\bibliography{Bounded-Type-Bib}{}

\printindex

\newpage
\section*{Proof Locator}
\label{page-proof-locator}
The following table  summarizes the connections between our number-theoretic results and our dynamical results, together with information on the location of proofs.
\begin{table}[h!]
\centering
\begin{tabular}{|c|c|c|}
\hline 
{\sc Number Theoretic} & {\sc  Dynamical Results} &  {\sc Location of }\\ 
{\sc Results}          & {\sc  Which Imply Them} &  {\sc Proofs}     \\
\hline
                                 & Thm \ref{t.Char-mu-xi-Gen} & \S\ref{ScGen}\\\cline{2-3}
Thm \ref{t.char-of-equidist-seq} & Thm \ref{CrLebGen} &  \S\ref{ScGen} \\ \cline{2-3}
                                 & Thm \ref{ThDriftOfOrigin} & Appendix \ref{Appendix-Twist}\\\hline
Thm \ref{t.tilt-and-exp-tilt}    & Thm \ref{t.historical-crit} &  \S\ref{ScGen} \\\hline

Thm \ref{t.exp-tilted-exist}    & Thm \ref{CrWeakBal} & \S\ref{SSZeroTilted}\\\hline

Thm \ref{t.exponential-tilt-HD} & Thm \ref{CrHDGen} & \\\cline{1-2}

Thm \ref{t.no-equidist-no-tilt} & Thm \ref{CrHDExcept} & \S\ref{ScNonGen}\\\hline

Thm \ref{t.A_n-B_n-estimates}   & Thm \ref{ThStairTLT} & \S\ref{ScCircle}\\\hline

Thm \ref{t.A_n}                 & Thm \ref{ThStairTLT} & \S\ref{ScCircle}\\\cline{2-3}

                                & Thm \ref{ThDriftOfOrigin} & Appendix \ref{Appendix-Twist}\\\hline
                                
Thm \ref{t.LLT-for-D_N}         & Thm \ref{t.Local-TCLT}   & \S\ref{ScCircle}\\\hline\hline
                                & {\sc Other Dynamical} & {\sc Location of} \\
                                & {\sc Results}         & {\sc Proofs}\\\cline{2-3}
                                & Thm \ref{CrLebGen2}   &  \S\ref{ScBirkhoff}\\\cline{2-3}
                                & Thm \ref{t.Char-mu-xi-Gen} & \S\ref{ScGen}\\\cline{2-3}
                                & Thm \ref{CrDriftSeq} &  \S\ref{ScGen} \\\cline{2-3}
                               & Thm \ref{CrLogOTDetXi} & \S \ref{ScBirkhoff} \\\cline{2-3}
                                & Thm \ref{CrHist}  & \S\ref{ScGen}\\\cline{2-3}
                                & Prop \ref{PrBalance-second half} & \\\cline{2-2}
                                & Thm \ref{ThLT-ZD} & \S\ref{ScBirkhoff} \\\cline{2-2}
                                & Thm \ref{ThLT-NZD} & \\\hline
\end{tabular}
\end{table}

\newpage
\section*{Notation Index}

\noindent
\begin{tabular}{llll}
{\bf A:} 
&$a_T(\tomega)$ & & the centering term in the temporal CLT, see Theorem \ref{ThStairTLT}\\
&$[a_N^R,b_N^R]$ & & the interval $[\mathcal F_N'(-R),\mathcal F_N'(R)]$, see Proposition \ref{PrRateFn}\\
&$[\ov{a}_N^R,\ov{b}_N^R]$ & & the interval $[\PP_N'(-R)/V_N,\PP_N'(R)/V_N]$, see Lemma \ref{l.rate-function}\\ 
&$A_n(\beta)$ && the $n$-th canonical matrix product, see \eqref{e.A-k}--\eqref{e.A-k-pq}\\
&$A_n$ && $A_{nN_0}(\beta)$, the derivative of the $n$-th renormalizing automorphism,\\ 
&&&  see Proposition \ref{Prop-Contracting-Automorphisms} and Remark \ref{r.mat-prod-form}\\  
&$A_T$, $A_T(\tomega)$  && the line segment $\{\phi^t(\tomega):0<t<T\}$\\
&$A_T^\ast$, $A_T^\ast(\tomega)$ && $(\tpsi_{n-1}\circ\cdots\circ\tpsi_0)(A_T)$ for the smallest $n=n(T)$ such that\\ 
&&& $|A_T^\ast|\in [1,L]$, see \eqref{DefBrN}\\
&
$\alpha$ & & an irrational number\\
{\bf B:}
& $\mathsf{BT}(M)$ && $\{[a_0;a_1,a_2,\cdots]: a_0\in\Z, a_1,a_2,\ldots\in\N, \sup|a_i|<M\}$\\ 
 & $\beta$ && $2\alpha-1$\\
{\bf C:} &
$C_c(X)$ && the space of real-valued continuous functions on $X$, with compact\\
&&&  support inside $X\setminus\{\text{singularities}\}$\\
&$C_{mix}$ && the mixing constant in Corollaries \ref{Cor-rho-mixing} and \ref{Cor-phi-mixing}\\
&$\Cov$ && the covariance:  $\Cov(X,Y)=\E[(X-\E(X))(Y-\E(Y))]$\\
&$\wt{C}_\ell$ && the lift of $C$ to $D^\ell(\wt{F}_0)$, see \eqref{e.curly-C-k}\\
&$\mathfs C_k$ && $\{\wt{C}_\ell: C=\<\dot{a}_0,\ldots,a_{k-1}\>,\ell\in\Z,\mathrm{int}(\wt{C}_\ell)\cap\mathrm{int}(\mathsf{R}_0)\neq \emptyset\}$\\
{\bf D:} &$d_n(\xi)$ && a structure constant, defined in \eqref{e.d}\\
&$d\tpsi_{\tomega}$ && the differential of an automorphism  $\tpsi$ at $\tomega$, viewed as a matrix\\
&$\dist$ && the metric on $\St$ associated to the locally Euclidean structure\\
&$D$ && the staircase automorphism which maps $\mathsf R_k$ to $\mathsf R_{k+1}$ by translation\\
&$D_N(\alpha,I)$ && $\#\{0\leq j\leq N-1: \{j\alpha\}\in I\}-N|I|$\\
&$D_N(\alpha)$ && $\#\{0\leq j\leq N-1: \{j\alpha\}\in [0,\frac{1}{2})\}-\frac{1}{2}N$\\
{\bf$\boldsymbol{\Delta}$:} & $\delta_x$ && the point mass measure at $x$\\
&$\delta(\tpsi)$ & & the average drift of an automorphism $\tpsi$, see \eqref{Eq-Drift-Defn}\\
&$\Delta\zc_0$, $\Delta\zc_{n,\ell}$ && edge effect terms in the decomposition $\Upsilon_N=\ov{S}_n-\Delta\zc_0+\Delta\zc_{n,\ell}$,\\
&&& defined in  Lemma \ref{l.Phi-Decomp} \\
{\bf E:} 
&$E^s_n,E^u_n$ && the stable and unstable spaces of $\psi_n$, $\tpsi_n$, see Proposition \ref{Prop-Contracting-Automorphisms}\\
&$\mathrm{erfc}(x)$ && the complementary error function, $\mathrm{erfc}(x)=\frac{2}{\sqrt{\pi}}\int_x^\infty e^{-t^2}dt$\\
&$\mathrm{ess\,sup}$ && the essential supremum\\
&$E(x)$ && the even Gauss map, see \eqref{EvenGauss}\\
&$\E$ && the expectation (of a random variable)\\
&$\wt{\E}^\xi$ &&  the expectation with respect to the change of measure, see \S\ref{ss.change-of-measure}\\
&
$\epsilon_0$ && the ellipticity constant, see Definition \ref{d.uniform-ellipticity} and Proposition \ref{Prop-Doeblin}\\
&$\epsilon(\alpha)$ && a uniform bound on the geometry of $E^s_n, E^u_n$, see Proposition \ref{Prop-Contracting-Automorphisms}\\
&$\epsilon_n, \ov{\epsilon}_n$ && error terms in Proposition \ref{Prop-Markov-Representation-of-Phi} and Lemma \ref{Lemma-Xi-and-S_n}\\
{\bf F:} &$f_k, \wt{f}_k^{(n)}$ && the jump functions, see \eqref{Eq-Modified-Frobenius} and Lemma \ref{Lemma-f(X)}\\
&$F_n$ && $\Pi(n)$ modulo side pairings (a model for $\St_0$), see \S\ref{Section-Polygons}\\
&$\wt{F}_n$ && the lift of $F_n$ to $\St$ described in \S\ref{Section-Modified-Z-coordinates}. Equivalently, $\{\tomega:\wt{\zc}_n(\tomega)=0\}$\\
&$\mathcal F_N(\xi)$ && the log moment generating function $(1/V_n)\log\E(e^{tS_N})$,\\
&&& see Definition \ref{DefFreeEn} and Proposition \ref{PrFreeEn}\\
&$\mathcal F_0^k, \mathcal F_\ell^\infty$ && the $\sigma$-algebras defined on page \pageref{page-sigma-fields}\\
& $\mathsf F_i$ && usually a finite union of horizontal rectangles $\mathsf R_j$\\
& $\mathsf{Frob}_{\tpsi}(\tomega)$ && the Frobenius function, see page \pageref{page-Frobenius}\\
\end{tabular}

\noindent
\begin{tabular}{llll}
{$\boldsymbol{\Phi}$:} &
$\varphi^t_{\vec{v}}$ && the linear flow on the infinite staircase, with constant velocity  $\vec{v}$\\
&$\phi^t$ && the linear flow on the infinite staircase, with constant velocity  ${\beta\choose 1}$\\
&$\Phi_n$ && the forward random walk, see \S\ref{s.Forward-RW}\\
{\bf G:} &$\mathsf{Gen}(\xi)$ && the set of generic points of $\mu_\xi$\\
&$\gcd$ && greatest common divisor\\
&$G$ && the group generated by $T(x,y)=(x+1,y+1)$ and  $S(x,y)=(x+2,y)$\\
&$G_{alg}, G_{ess}$ && the algebraic range and the essential range, see Definitions \ref{d.alg-rng}, \\
&&& \ref{d.ess-rng}, and  Theorem \ref{t.co-range}\\
&$\mathfs G_\ell, \mathfs G$ && see \S\ref{ScBulk}\\
{\bf$\boldsymbol{\Gamma}$:} & $\Gamma(2)$ && the group $\left\{{\tiny \biggl(\begin{array}{ll}
a & b\\
c & d
\end{array}\biggr)}:a,b,c,d\in\Z,\ ad-bc=1,\  {\tiny \biggl(\begin{array}{ll}
a & b\\
c & d
\end{array}\biggr)}=
{\tiny \biggl(\begin{array}{ll}
1 & 0\\
0 & 1
\end{array}\biggr)}\mod 2
\right\}$\\
&$\Gamma(M)$ && the ``geometric constant" in Lemma \ref{Lemma-Bounded-Geometry-Gamma(M)} (a bound on the geometry of $\Pi(n)$)
\end{tabular}

\noindent
\begin{tabular}{lll}
&&$\Gamma\left({
x_{n-2} \begin{array}{l}
x_{n-1}\\
y_{n-1}
\end{array}
\begin{array}{l}
x_{n}\\
y_{n}
\end{array} y_{n+1}}
\right)$ the balance of a hexagon, see \S\ref{ss.structure-constants}\\
\end{tabular}

\noindent
\begin{tabular}{llll}
{\bf H:} &
$h_n(\cdot,\xi)$ && the functions appearing in the change of measure, see Lemma \ref{LmAddPres}\\
&$\mathrm{Hex}(n)$ && the space of hexagons at position $n$, see \S\ref{ss.structure-constants}\\
&$\mathrm{HD}$ && the Hausdorff dimension (of a set)\\
&$H(\mathsf X,\mathsf f)$ && the co-range, see Theorem \ref{t.co-range}\\
&$\mathcal H(\tomega)$ && the collection of historical measures of $\tomega$, see \S\ref{SSNonGen}\\
&$\mathsf H_n(\eta)$ && the Legendre transform of $\mathsf P_n(\xi)/V_n$,  see Lemma \ref{l.rate-function}\\
{\bf I:} &
$\mathrm{int}(\cdot)$ && the interior (of a set)\\
&$\mathcal I_n(\eta)$ && the rate function, the Legendre transform of $\mathcal F_n(\xi)$, see Prop. \ref{PrRateFn}\\
&$\bbI_T(\tomega)$ && $\log \int_0^T 1_{[\zc\leq 1]}(\phi^t(\tomega))dt$, see Theorem \ref{CrLogOTDetXi}\\
&$\bbI_T^\xi$ &&  $\int_{\Rect_0}\bbI_T(\tomega)d\mu_\xi(\tomega)$, see Theorem \ref{CrLogOTDetXi}\\
&$\frak{I}_T(\tomega,h)$ && $\int_0^T h(\phi^t(\tomega))dt$\\
{\bf$\boldsymbol{\L}$:} &
$\lambda^s_n, \lambda^u_n$ && $\lambda^s_n:=\|A_n|_{E^s_{n}}\|\text{ and }\lambda^u_n:=\|A_n|_{E^u_{n}}\|$\\
&$\lambda$ && a constant bigger than one, chosen in  \eqref{e.choice-of-lambda}, such that  $\lambda^s_n\leq \lambda^{-1}$ \\
&&& and $\lambda^u_n\geq \lambda$ for all $n$, see Proposition \ref{Prop-Contracting-Automorphisms}\\
& $\Lambda_6$ && $\{\beta\in (-1,1)\setminus\Q: M(\beta)\leq 6\}$\\
{\bf L:} &$\log$ && the natural logarithm (the same as $\ln$)\\
&$\ell^s, \ell^u$ && the lengths of the stable and unstable sides\\
&$L(o,\vec{v},k)$ && the ray  from $o$, in direction $\vec{v}$,  initially lying in $\mathsf R_k$\\
{\bf M:} &
$\mu_0$ &&  $\frac{1}{2}\times$ area measure on $\St$\\
&$\mu_\xi$ && the unique ergodic $\phi^t$-invariant locally finite measure such that\\
&&&  $\mu_\xi\circ D=e^\xi\mu_\xi$ and  $\mu_\xi(\mathsf R_0)=1$, see Corollary \ref{c.HHW}\\
&$\mu_{\Rect_0}$ && $\frac{1}{2}\times$ area measure on $\Rect_0$\\
&$\mathrm{mes}$ && Lebesgue's measure on $\R$\\
&$m_{\mathrm{Hex}}^n$ && the hexagon measure, see \S\ref{ss.structure-constants}\\
&$M(\beta)$ && a characteristic of the  even CFE of $\beta$, see \eqref{e.M-beta}\\
&$M_\xi(\wt{C})$ && see \eqref{DefMM}
\end{tabular}

\noindent
\begin{tabular}{llll}
{\bf N:} &
$n(T)$ && the minimal number $n$ such that $|A_T^\ast|\equiv |(\tpsi_{n-1}\circ\cdots\circ\psi_0)(A_T)|$\\
&&& lies in $[1,L]$, see \eqref{DefBrN}. Satisfies $n(T)\asymp\log T$.\\
&$N_0$ && the integer such that $A_{nN_0}(\beta)=A_{n-1}\cdots A_1 A_0$ and\\
&&& $\wt{\Psi}_{nN_0}=\tpsi_{n-1}\circ\cdots\circ\tpsi_0$ for all $n$. See  Remark \ref{r.mat-prod-form}\\
&$N_i(n)$ && $(i=1,2)$ see Lemma \ref{Lemma-Partition} and Proposition \ref{Prop-Markov}\\
&$\mathcal N$ && depending on the context,  the normal distribution with zero \\
&&& mean and unit variance, or a   random variable with this\\
&&& distribution\\
{\bf $\boldsymbol{\Omega}$:} & $\omega$ && (usually) a point in $\St_0$\\
& $\tomega$ && (usually) a point in $\St$\\
& $\omega(x,k)$ && the point on the top side of $\mathsf R_k$, $y$ units of distance from\\
&&& right of the middle, where $y\in [-1,1)$, $y=2x\mod 2$. See Theorem \ref{t.HHW} \\  
&$\tomega_0$ && $\omega(0,0)$, see Theorem \ref{t.HHW} and \eqref{e.Lin-Flow-1}\\
&$\tomega_{-x}$ && $\omega(-x,0)$, see Theorem \ref{t.HHW} and \eqref{e.Lin-Flow-2}\\
{\bf P:} 
&$p_n(\xi)$ && parameters appearing in the change of measure, see Lemma \ref{LmAddPres}\\
&$\mathsf P_n(\xi)$ && the pressure function, equal to $p_1(\xi)+\cdots+p_n(\xi)$, see Lemma \ref{LmPressure}\\
&&& and Proposition \ref{p.Pressure-and-Geometric-Pressure}\\
&$\mathcal P_n(\xi)$ && the geometric pressure function, equal to $\log \E_{\mu_{\Rect_0}}(e^{\xi\Upsilon_n})$. See \S\ref{ss.geom-press}\\
&$\mathcal P_n'(\pm\infty)$ && $\lim_{t\to\pm \infty}\mathcal P_n'(t)$\\
&$\mathfrak P_n$ && the $n$-th Markov partition, see Lemma \ref{Lemma-Partition} and  Proposition \ref{Prop-Markov}\\
&$\Prob$ && Probability (of an event)\\
&$\Prob_W$ && normalized length measure on a linear segment $W$\\
&$\Prob_Q$ && normalized area measure on a parallelogram $Q$\\
$\boldsymbol{\Pi:}$&$\pi$ && the natural covering map  $\pi:\St^\ast\to\St_0$, see \S\ref{ss.Z-d-cover}\\
&$\pi_{i,i+1}$ && the transition kernel of $\{X_n\}_{n\geq 1}$, see Proposition \ref{Prop-Markov-Kernel}\\
&${\pi}_{i,i+1}^\xi$ && the transition kernel of the change of measure, see \eqref{ProbXi}\\
&$\Pi(E^u,E^s)$ && the polygon in  Figure \ref{Figure-Pi(n)}\\
&$\Pi(n)$ && $\Pi(E^u_n,E^s_n)$\\ 
$\boldsymbol{\Psi:}$ &$\wt{\Psi}_n$ && the homogeneous automorphism with zero drift and derivative $A_n(\beta)$\\
&$\tpsi_i$ && the renormalizing automorphisms  from Proposition \ref{Prop-Contracting-Automorphisms}.\\
&&& Note that $\tpsi_{n-1}\circ\cdots\circ\tpsi_0=\wt{\Psi}_{nN_0}$\\
&$\psi_i$ && the projections of $\tpsi_i$ to automorphisms of $\St_0$, see \eqref{Eq-Projected-Auto}\\
{\bf Q:} & $\wt{Q}_1, \wt{Q}_2$ && the two parallelograms which combine to  $\Pi$, see Lemma \ref{Lemma-Q}\\
& $Q_i(n)$ && $(i=1,2)$ the $\wt{Q}_i$  which combine to  $\Pi(n)$, see Lemma \ref{Lemma-Q}\\
& $Q_{i,j}^n$ && the elements of the Markov partition, see Lemma \ref{Lemma-Partition}\\
{\bf R:} & $\mathsf R_k$ && the $k$-th horizontal rectangle, minus the corners and the\\
&&&   bottom side. See Figure \ref{Figure-Staircase}\\
{\bf S:} & 
 $S_N$ && $f_1(X_1)+\cdots+f_N(X_N)$, except for Appendix \ref{AppMC}, where more\\
&&& general sums are considered\\
& $\wt{S}_N^\xi$ && $f_1(\wt{X}_1^\xi)+\cdots+f_N(\wt{X}_N^\xi)$. Equivalently, $S_N$ but with the joint \\
&&& distribution of the change of measure from  \S\ref{ss.change-of-measure} \\
& $\ov{S}_N$ &&  $\sum_{k=1}^N f_k(\ov{X}_k(\tomega))$, where 
$\ov{X}_k(\tomega)$ is the unique element of $\mathfrak P_k$ \\
&&& containing $(\psi_{k-1}\circ\cdots\circ\psi_0)(\pi(\tomega))$. See Lemma \ref{l.Phi-Decomp}.\\
&$\sgn(x)$ && the sign of $x$: (+1) for $x>0$, $(-1)$ for $x<0$, and $0$ for $x=0$.
\end{tabular}

\noindent
\begin{tabular}{llll}
& $\mathrm{SL}(2,\R)$ && the group of $2\times 2$ real-valued matrices with determinant one\\
& $\mathrm{SL}(2,\Z)$ && the group of $2\times 2$ integer-valued matrices with determinant one\\
& $\St$ && the infinite staircase, see \S\ref{s.infinite-staircase} and Figure \ref{Figure-Staircase}\\
& $\St^\ast$ && the infinite staircase, minus the singularities (the corners)\\
& $\St_0$ && the finite area translation surface in Figure \ref{Figure-Staircase}(b) together with its singularities.\\ 
&&&  $\St_0$ is isometric to $\R^2/\sqrt{2}\Z^2$. \\
& $\St_0^\ast$ && $\St_0$ minus the congruence classes of the corners.\\
&&& $\St^\ast$ is a $\Z$-cover of $\St_0^\ast$\\
& $\St_\pm$ && $\St_+=\bigcup_{k\geq 0}\mathsf R_k$, $\St_-=\bigcup_{k< 0}\mathsf R_k$\\
$\boldsymbol{\Sigma:}$ & $\Sigma$ && The union of the top sides of $\mathsf R_k$. A Poincar\'e section for $\phi^t:\St\to\St$ with\\ 
&&&  section map canonically isomorphic to the cylinder map, see Theorem \ref{t.HHW}\\
& $\sigma(\vec{v})$ && the slope: $\sigma{x\choose y}=y/x$\\
& $\hsigma(\vec{v})$ && the co-slope: $\hsigma{x\choose y}=x/y$\\
& $\sigma_T$ && the ``standard deviation" in the temporal CLT and LLT  (Theorem \ref{ThStairTLT}\\
&&& and Theorem \ref{t.Local-TCLT})\\
{\bf T:} & $T(\alpha)$ && $\{x: D_N(\alpha,x+I)\text{ is tilted}\}$, see \S\ref{ss.tilt}\\
& $T_t(\alpha)$ && $\{x: D_N(\alpha,x+I)\text{ is exponentially tilted with parameter $t$}\}$, see \S\ref{ss.tilt}\\
& $T_0(M)$ && the ``denseness constant" from the end of  \S\ref{ss.density-const}\\
& $\T$ && $\R/\Z$\\
& $t_{a,b}^{(i)}$ && the entries of the $i$-th transition matrix, see \eqref{e.incidence-matrix}\\
& $\tau(x)$ && the cocycle defining the cylinder map, $\tau(x):=2\cdot 1_{[0,\frac{1}{2})}(\{x\})-1$, see \S\ref{ss.Cylinder-Map}\\
& $\tau_N(x)$ && $\sum_{j=0}^{N-1}\tau(x+j\alpha)$\\
{\bf U:} & $\mathcal U$ && depending on the context, the uniform distribution on $[0,1]$\\
&&& or a random variable with this distribution\\
{\bf $\boldsymbol{\Upsilon:}$} & $\Upsilon_N$ && the forward random walk driven by renormalizations, see \S\ref{s.RW-Driven-By-Renormalizations}\\
{\bf V:} & $v_a(z)$ && the M\"obius transformation $\frac{1}{a+z}$ \\
&$\mathrm{Var}$ && the variance (of a random variable), $\mathrm{Var}(X)=\E[(X-\E(X))^2]$\\
& $V_N$ && $\mathrm{Var}(S_N)$\\
& $\wt{V}_N^\xi$ && the variance of $S_N$ with respect to the change of measure in \S\ref{ss.change-of-measure}.\\ 
&&&Equivalently, $\mathrm{Var}(\wt{S}_N^\xi)$\\
{\bf W:} & $w$ && (unless stated otherwise) a uniform geometric bound on the elements\\
&&&  of the Markov partitions. Every $Q_{ij}^n\in\mathfrak P_n$ contains a ball with radius\\
&&& $w$. See Lemma \ref{Lemma-Partition}\\
& $W^u(\omega,F_n)$ && the open $u$-fiber in $F_n$ through $\omega$\\
& $W^s(\omega,F_n)$ && the open $s$-fiber in $F_n$ through $\omega$\\
& $W^u(\omega,Q)$ && $W^u(\omega,F_n)\cap Q$ ($Q$ an $(E^u_n,E^s_n)$-parallelogram)\\
& $W^s(\omega,Q)$ && $W^s(\omega,F_n)\cap Q$ ($Q$ an $(E^u_n,E^s_n)$-parallelogram)\\
& $W_{n,i}$ && (usually) pure $E^s_n$-segments which appear in approximations of \\
&&& non-pure  $E^s_n$-segments as in Lemma \ref{l.approx-by-pure-segments}\\
{\bf X:} & $X_i$ && the elements of a Markov chain (unless stated otherwise,  the one  in \S\ref{ss.IMC})\\
& $X_i^{(n)}(\omega)$ && the element of $\mathfrak P_i$ which contains $(\psi_{i}^{-1}\circ\cdots\circ\psi_{n-1}^{-1})(\omega)$ when  $i<n$;\\
&&& $(\psi_{i-1}\circ\cdots\circ \psi_n)(\omega)$ when $i>n$; and 
 $\omega$ when $i=n$.  See \eqref{e.X-new}\\
& $\ov{X}_i(\tomega)$ && the element of $\mathfrak P_k$ containing $(\psi_{k-1}\circ\cdots\circ \psi_0)(\pi(\tomega))$, see \eqref{e.X-bar}\\
& $\mathsf X$ && a Markov chain or a Markov array\\
$\boldsymbol{\Xi:}$ & $\Xi_n^{W_n}(\tomega)$ && the backward finite random walk, see \S\ref{s.RW-Driven-By-Renormalizations}\\
& $\xi_n(\tomega)$ & & the unique $\xi$ s.t.  $\mathcal P_n'(\xi)=z_n(\tomega)$, when $z_n(\tomega)\in(\mathcal P_n'(-\infty),\mathcal P_n'(\infty))$.\\
&&& Otherwise, if $z_n(\tomega)\leq \inf \mathcal P_n'$, we set $\xi_n(\tomega):=-\infty$; and\\
&&& if $z_n(\tomega)\geq \inf \mathcal P_n'$, we set $\xi_n(\tomega):=+\infty$. See \S\ref{SSGP-Tlit}.
\end{tabular}

\noindent
\begin{tabular}{llll}
{\bf Z:} & $z_n(\tomega)$ && For $\tomega\in\St$ not equal to  a singularity, this is the canonical $\Z$-coordinate of\\
&&& $(\tpsi_{n-1}\circ\cdots\circ\tpsi_0)(\tomega)$.
For singular $\tomega$,  this is the canonical $\Z$-coordinate of\\
&&& $(\tpsi_{n-1}\circ\cdots\circ\tpsi_0)(\phi^t(\tomega))$ for all $t>0$ small enough.\\
& $\zc$ && canonical $\Z$-coordinate (equal to $k$ on $\mathsf R_k$)\\
& $\wt{\zc}_n$ && modified $\Z$-coordinate (equal to $k$ on $D^k(\wt{F}_k)$)\\ 
\end{tabular}

\noindent
\begin{tabular}{ll}
&\\
$1_E$ & the indicator function, equal to $1$ on $E$ and $0$ outside $E$\\
$[a_0;a_1,\ldots]$ & the limit $\lim\limits_{n\to\infty}a_0+\cfrac{1}{a_1+\cfrac{1}{a_2+\cdots}}$ (when it exists)\\
$\<a_{n-k},\ldots,\dot{a}_n,\ldots,a_{n+m}\>$ & a cylinder with beginning $n-k$, center $n$, ending $n+m$,\\
& and coordinates $a_i$, see \S\ref{ss.cylinders}\\ 
$\{x\}$ & unique number in $[0,1)$ such that $x\in\{x\}+\Z$\\
$\{x\}_{[0,2)}$ & unique number in $[0,2)$ such that $x\in\{x\}+2\Z$\\
$\lfloor x\rfloor$, $\lceil x\rceil$ & $\lfloor x\rfloor:=\max\{n\in\Z:n\leq x\}$, $\lceil x\rceil:=\min\{n\in\Z:n\geq x\}$\\
$\|x\|$ & $\dist(x,\Z)$\\
$|I|$ & length or area, depending on the context;\\
$\#A$ & the number of elements in a set $A$\\
$x+I$ & $\{\{x+y\}:y\in I\}$\\
$:=$ & ``is defined to be "\\
$\equiv$ & ``is equal by definition or by some trivial reason"\\
$\overset{!}{=}$ & an equality which shall be justified below\\
$\overset{?}{=}$ & an equality which could be false\\ 
$a_n\asymp b_n$ & $\exists M,N>1$ s.t. $\forall n>N$, $M^{-1}\leq a_n/b_n\leq M$\\
$a_n\sim b_n$ & $a_n/b_n\xrightarrow[n\to\infty]{}1$\\
$b\pm C$ & a quantity in $[b-C,b+C]$ \\
$C^{\pm 1}b$ & a quantity between $C^{-1}b, Cb$\\
$({\tiny \begin{array}{cc} a & b \\ c & d \end{array}})\cdot z$ & $\frac{az+b}{cz+d}$\\
$X_n\xrightarrow[n\to\infty]{dist}Y$ & the random variables $X_n$ converge to $Y$ in distribution
\end{tabular}

\medskip
\noindent
\begin{tabular}{llll}
a.e. & almost everywhere\\
a.s. & almost surely\\
CFE & continued fraction expansion\\
CLT & central limit theorem\\
IID & independent and identically distributed (random variables)\\
s.t. & such that
\end{tabular}

\end{document}